\documentclass[graybox,envcountchap,sectrefs]{svmono}

\usepackage{type1cm}         

\usepackage{makeidx}         
\usepackage{graphicx}        
\usepackage{multicol}        
\usepackage[bottom]{footmisc}

\usepackage{newtxtext}       %
\usepackage[varvw]{newtxmath}       

\makeindex             

\usepackage{chapterbib}
\usepackage{hyperref,url}
\usepackage{mhchem}
\usepackage{physics}
\usepackage{verbatim}
\usepackage{xcolor,soul}
\usepackage{float} 
\usepackage{tabulary}
\newcolumntype{C}[1]{>{\centering\arraybackslash}p{#1}}

\theoremstyle{definition}
\newtheorem{thm}{Theorem}[section]
\newtheorem{dfn}[thm]{Definition}
\newtheorem{pro}[thm]{Problem}
\newtheorem{cor}[thm]{Corollary}
\newtheorem{lem}[thm]{Lemma}
\newtheorem{rem}[thm]{Remark}
\newtheorem{exa}[thm]{Example}

\newtheorem{conj}[thm]{Conjecture}

\newtheorem{prot}[thm]{Protocol}
\newtheorem{prop}[thm]{Proposition}

\newcommand{\C}{\mathbb{C}}
\newcommand{\R}{\mathbb{R}}
\newcommand{\Q}{\mathbb{Q}}
\newcommand{\Z}{\mathbb{Z}}
\newcommand{\Or}{\mathrm{O}}
\newcommand{\GL}{\mathrm{GL}}
\newcommand{\SL}{\mathrm{SL}}
\newcommand{\SO}{\mathrm{SO}}
\newcommand{\Eu}{\mathrm{E}}
\newcommand{\SE}{\mathrm{SE}}
\newcommand{\BD}{\mathrm{BD}}
\newcommand{\HD}{\mathrm{HD}}
\newcommand{\CIS}{\mathrm{CIS}}
\newcommand{\CDS}{\mathrm{CDS}}

\newcommand{\CRS}{\mathrm{CRS}}
\newcommand{\CRIS}{\mathrm{CRIS}}
\newcommand{\CIMS}{\mathrm{CIMS}}
\newcommand{\CRDS}{\mathrm{CRDS}}
\newcommand{\CRHS}{\mathrm{CRHS}}

\newcommand{\bri}{\mathrm{BRI}}
\newcommand{\bris}{\mathrm{BRIS}}
\newcommand{\trin}{\mathrm{TRIN}}
\newcommand{\brain}{\mathrm{Brain}}

\newcommand{\bid}{\mathrm{BID}}
\newcommand{\bib}{\mathrm{BIB}}
\newcommand{\RMS}{\mathrm{RMS}}

\newcommand{\SM}{\mathrm{SM}}
\newcommand{\LAC}{\mathrm{LAC}}
\newcommand{\PCM}{\mathrm{PCM}}
\newcommand{\PCI}{\mathrm{PCI}}
\newcommand{\WMI}{\mathrm{WMI}}

\newcommand{\cov}{\mathrm{Cov}}

\newcommand{\SRD}{\mathrm{SRD}}
\newcommand{\SPD}{\mathrm{SPD}}
\newcommand{\PDD}{\mathrm{PDD}}
\newcommand{\SCD}{\mathrm{SCD}}
\newcommand{\SDD}{\mathrm{SDD}}
\newcommand{\RDD}{\mathrm{RDD}}
\newcommand{\wRDD}{\widetilde{\RDD}}
\newcommand{\ORD}{\mathrm{ORD}}
\newcommand{\mORD}{\overline{\ORD}}
\newcommand{\OSD}{\mathrm{OSD}}
\newcommand{\mOSD}{\overline{\OSD}}
\newcommand{\OCD}{\mathrm{OCD}}

\newcommand{\mSCD}{\overline{\SCD}}

\newcommand{\ADD}{\mathrm{ADD}}
\newcommand{\ASD}{\mathrm{ASD}}

\newcommand{\MSD}{\mathrm{MSD}}
\newcommand{\WSD}{\mathrm{WSD}}
\newcommand{\WDD}{\mathrm{WDD}}
\newcommand{\wSDD}{\widetilde{\SDD}}
\newcommand{\VID}{\mathrm{VID}}
\newcommand{\VSM}{\mathrm{VSM}}

\newcommand{\DM}{\mathrm{DM}}

\newcommand{\CDM}{\mathrm{CDM}}

\newcommand{\MCD}{\mathrm{MCD}}
\newcommand{\MCS}{\mathrm{MCS}}

\newcommand{\TS}{\mathrm{TS}}
\newcommand{\CT}{\mathrm{CT}}

\newcommand{\CI}{\mathrm{CI}}
\newcommand{\DI}{\mathrm{DI}}
\newcommand{\CR}{\mathrm{CR}}
\newcommand{\DR}{\mathrm{DR}}
\newcommand{\CIM}{\mathrm{CIM}}
\newcommand{\DIM}{\mathrm{DIM}}
\newcommand{\CRM}{\mathrm{CRM}}
\newcommand{\DRM}{\mathrm{DRM}}

\newcommand{\lcm}{\mathrm{lcm}}

\newcommand{\DC}{\mathrm{DC}}
\newcommand{\DT}{\mathrm{DT}}
\newcommand{\QS}{\mathrm{QS}}
\newcommand{\VF}{\mathrm{VF}}
\newcommand{\CF}{\mathrm{CF}}
\newcommand{\RC}{\mathrm{RC}}
\newcommand{\PC}{\mathrm{PC}}
\newcommand{\RI}{\mathrm{RI}}
\newcommand{\PRF}{\mathrm{\overline{\RI}}}
\newcommand{\PI}{\mathrm{PI}}

\newcommand{\Oct}{\mathrm{Oct}}
\newcommand{\TC}{\mathrm{TC}}
\newcommand{\TP}{\mathrm{TP}}
\newcommand{\QT}{\mathrm{QT}}
\newcommand{\MS}{\mathrm{MS}}

\newcommand{\RM}{\mathrm{RM}}
\newcommand{\PM}{\mathrm{PM}}
\newcommand{\SIM}{\mathrm{SIM}}
\newcommand{\SRM}{\mathrm{SRM}}
\newcommand{\SDM}{\mathrm{SDM}}
\newcommand{\SHM}{\mathrm{SHM}}

\newcommand{\RIS}{\mathrm{RIS}}
\newcommand{\PIN}{\mathrm{PIN}}
\newcommand{\SBI}{\mathrm{SBI}}
\newcommand{\SBR}{\mathrm{SBR}}
\newcommand{\SBD}{\mathrm{SBD}}
\newcommand{\SBH}{\mathrm{SBH}}
\newcommand{\SLM}{\mathrm{SLM}}
\newcommand{\HS}{\mathrm{HS}}
\newcommand{\EP}{\mathrm{EP}}

\newcommand{\LIS}{\mathrm{LIS}}
\newcommand{\LRS}{\mathrm{LRS}}
\newcommand{\LDS}{\mathrm{LDS}}
\newcommand{\LHS}{\mathrm{LHS}}

\newcommand{\RB}{\mathrm{Red}}

\newcommand{\BT}{\mathrm{BT}}
\newcommand{\GC}{\mathrm{GC}}

\newcommand{\AMD}{\mathrm{AMD}}
\newcommand{\ADA}{\mathrm{ADA}}
\newcommand{\PDA}{\mathrm{PDA}}
\newcommand{\AND}{\mathrm{AND}}
\newcommand{\PND}{\mathrm{PND}}

\newcommand{\LND}{\mathrm{LND}}

\newcommand{\EMD}{\mathrm{EMD}}

\newcommand{\IT}{\mathrm{IT}}

\newcommand{\PPC}{\mathrm{PPC}}

\newcommand{\sign}{\mathrm{sign}}

\newcommand{\Ga}{\Gamma}
\newcommand{\de}{\delta}

\newcommand{\ep}{\varepsilon}
\newcommand{\es}{\emptyset}
\newcommand{\Ra}{\Rightarrow}
\newcommand{\al}{\alpha}
\newcommand{\la}{\lambda}
\newcommand{\La}{\Lambda}
\newcommand{\ga}{\gamma}
\newcommand{\be}{\beta}
\newcommand{\si}{\sigma}
\newcommand{\Si}{\Sigma}
\newcommand{\ph}{\varphi}
\newcommand{\om}{\omega}
\newcommand{\ti}{\tilde}
\newcommand{\vl}{\, | \,}
\newcommand{\bd}{\partial}

\newcommand{\sym}{\mathrm{Sym}}
\newcommand{\iso}{\mathrm{Iso}}
\newcommand{\vol}{\mathrm{vol}}

\newcommand{\aff}{\mathrm{aff}}
\newcommand{\lra}{\leftrightarrow}

\newcommand{\angstrom}{\textup{\AA}}

\newcommand{\ar}[1]{\langle #1 \rangle}
\newcommand{\cbox}[2]{\colorbox{#1}{$\displaystyle #2$}}
\newcommand{\ve}[1]{\overrightarrow{#1}}
\newcommand{\ov}[1]{\overline{#1}}
\newcommand{\pd}[2]{\frac{\partial #1}{\partial #2}}

\newcommand{\mc}[2]{\colorbox{#1}{$\displaystyle #2$}}
\newcommand{\vect}[2]{ \left( \begin{array}{c} 
 #1 \\ #2 \end{array} \right)}
\newcommand{\colthree}[3]{ \left( \begin{array}{c} 
 #1 \\ #2 \\ #3 \end{array} \right)}

\newcommand{\mat}[4]{ \left( \begin{array}{cc} 
 #1 & #2 \\ #3 & #4 \end{array} \right)}
 \newcommand{\matfour}[4]{\left(\begin{array}{ccc}
#1 & #2 \\ #3 & #4 \end{array}\right)}
\newcommand{\matv}[4]{ \left( \begin{array}{cc} 
 #1 & #3 \\ #2 & #4 \end{array} \right)}

\newcommand{\mathree}[6]{ \left( \begin{array}{cc} 
 #1 & #4 \\ #2 & #5 \\ #3 & #6 \end{array} \right)}

\newcommand{\bs}{\hfill $\blacksquare$}
\newcommand{\bt}{\hfill $\blacktriangle$}
\newcommand{\eexa}{\hfill $\Diamondblack$}
\newcommand{\elem}{\hfill $\blacksquare$}
\newcommand{\erem}{\hfill $\Diamondblack$}
\newcommand{\ecor}{\hfill $\blacksquare$}
\newcommand{\ethm}{\hfill $\blacksquare$}
\newcommand{\epro}{\hfill $\bigstar$}
\newcommand{\edfn}{\hfill $\blacktriangle$}

\newcommand{\tb}[1]{\textbf{#1}}
\newcommand{\nt}{\noindent}
\newcommand{\sskip}{\smallskip}
\newcommand{\myskip}{\medskip}
\newcommand{\bskip}{\bigskip}

\begin{document}

\author{Olga D. Anosova, Vitaliy A. Kurlin}
\title{Geometric Data Science}
\subtitle{Moduli spaces of real data under equivalences} 
\maketitle

\frontmatter

%
%

\preface




\hfill
\emph{Where there is Matter, there is Geometry}
\myskip

\hfill
\emph{--- Johannes Kepler (1571- 1630)}
\myskip

\hfill
\emph{a key figure in the 17th-century Scientific Revolution}.
\bskip

This book introduces the new research area of \emph{Geometric Data Science}, where data can represent any real objects through geometric measurements. 
Some of the simplest inputs of real data objects are finite and periodic sets of unordered points.
\myskip

For example, a molecule can be fully described by the positions of its atoms in a 3-dimensional space.
However, many descriptions are highly ambiguous, especially to a computer, which operates only with numbers.
For example, a photograph is ambiguous, because any object can have an astronomically large number of photographs.
\myskip

All attempts to standardise photographs, as in passports, have shifted towards more reliable biometric data.
Indeed, the identification of living organisms was dramatically improved due to the discovery of a DNA structure.  
However, \emph{geometric structures} remained ambiguous for many objects, including proteins and materials, which are still represented by photograph-style inputs depending on arbitrary coordinate systems.
\myskip

The major obstacle to progress from trial-and-error in chemistry and biology to a justified design of materials and drugs was the absence of rigorous definitions and problem statements.
Geometric Data Science fills this gap by developing foundations based on equivalences, invariants, distance metrics, and polynomial-time algorithms.
\myskip

The main \emph{geo-mapping problem}\index{Geo-Mapping Problem}
 is to analytically describe moduli spaces of \emph{geometric structures} that are classes of data objects modulo an equivalence relation.
These moduli spaces are prototypes of `treasure maps' containing all known objects of a certain type as well as all not yet discovered ones.
A discrete example is Mendeleev's table of chemical elements, which was initially half-empty, but importantly guided an efficient search for new elements.  
A continuous example is a geographic map of the Earth, where any location is unambiguously identified by the latitude and longitude. 
\myskip

Geometric Data Science aims to develop universal geographic-style coordinates for all real data objects under practically important equivalences, such as rigid motion. 
\myskip

The first part of the book focuses on finite point sets.
The most important result is a complete and continuous classification of all finite clouds of unordered points under rigid motion in any Euclidean space. 
The key challenge was to avoid the exponential complexity arising from permutations of the given unordered points.
For a fixed dimension of the ambient Euclidean space, the times of all algorithms for the resulting invariants and distance metrics depend polynomially on the number of points.
\myskip

The second part of the book advances a similar classification in the much more difficult case of periodic point sets, which model all periodic crystals at the atomic scale.
The most significant result is the hierarchy of invariants from the ultra-fast to complete ones.
The key challenge was to resolve the discontinuity of crystal representations that break down under almost any noise.
Experimental validation on all major materials databases confirmed the \emph{Crystal Isometry Principle}: any real periodic crystal has a unique location in a common moduli space of all periodic structures under rigid motion. 
The resulting moduli space contains all known and not yet discovered periodic crystals and hence continuously extends Mendeleev's table to the full crystal universe.
\myskip
 
The book was written for research students and professionals who work in mathematics and need rigorously justified and computationally efficient methods for real data. such as crystalline materials and molecules, including proteins.
The pre-requisite knowledge is linear algebra, metric geometry, and calculus at the undergraduate level. 
\myskip

We finish by extending Johannes Kepler's quote from the 17th century to inspire a transformation from brute-force computations, which currently `burn' our planet, to a 21st-century \emph{Maths for Science} revolution: \emph{where there is Data, there is Geometry}.
\myskip
 
\nt
\textbf{Acknowledgments}.
We thank Tatiana Kurlina (University College London) for helpful comments on the initial draft, and the Data Science Theory and Applications group at the Materials Innovation Factory (Liverpool) for developing Geometric Data Science.
\myskip

Most significantly, Dr Daniel Widdowson implemented the invariant-based \emph{Crystal Geomap} to visualise in real time the crystal universe of all known materials.
Dr Yury Elkin substantially contributed by correcting past claims on computational complexities of nearest neighbour search and by parallelising computations for molecular structures.
Dr Matthew Bright produced the first geographic-style maps of 2.6+ million 2D lattices extracted from real periodic crystals.
Dr Philip Smith's software for density functions demonstrated the first impact by detecting a missed crystal, which was confused with a different one by the authors of the original work.
The implementations of Mr Ziqiu Jiang and Mr William Jeffcott exposed thousands of duplicate chains in the Protein Data Bank.  
Dr Jonathan Balasingham adapted ultra-fast crystal invariants for predicting material properties.
Dr Miloslav Torda, Dr Jonathan McManus, Dr Marjan Safi-Samghabadi, and Mr Surya Majumder validated Geometric Data Science methods on real data.
\myskip

We are grateful to Prof Marjorie Senechal, Prof Nikolai Dolbilin, Prof Andy Cooper FRS, Prof Sally Price FRS, Prof Graeme Day, Prof Simon Billinge, Prof Ram Seshadri, Prof John Helliwell, and Prof Yulia Gel for support and insightful discussions.
The latest version is
at \url{http://kurlin.org/Geometric-Data-Science-book.pdf}.

\begin{flushright}\noindent
United Kingdom,
\hfill {\it Dr Olga D. Anosova}\\
November 2025
\hfill {\it Prof Vitaliy A. Kurlin}\\
\end{flushright}

\graphicspath{{images/}}

\tableofcontents

\mainmatter

%
%
%

\chapter{Introduction: from practical challenges to fundamental problems}
\label{chap:intro} 


\abstract{
This chapter discusses how practical challenges in object recognition and data comparison can be converted into formally stated mathematical problems.
After introducing the necessary concepts of equivalences, invariants, and metrics, we state the general \emph{geo-mapping problem} to continuously parametrise moduli spaces for any data under a given equivalence.
With this foundation in place, further chapters examine specific types of data objects that allow recently developed solutions in this book.
}

\section{What questions should we ask about real data objects?}
\label{sec:questions}

The initial question that can be asked about any real data object is \emph{what is it?} or (more formally) \emph{how is it defined?} or (more deeply) \emph{how can we make sense of this data?}
\myskip

The first obstacle in achieving these goals is to embrace differences between real objects and their digital representations.
For example, a car is a physical object that is very different from a pixel-based image of this car, which is only a matrix of integers.
\myskip

The second obstacle is the ambiguity of digital representations in the sense that any real object can have many representations that look very different to a computer.
\myskip

If measurements have continuous real values, the resulting space of representations is infinite.
Even if we fix a finite resolution of physical measurements, all potential data values still live in a huge space.
For example, all images of size $2\times 2$ pixels and greyscale intensities $0,\dots,255$ form a huge collection of $256^4>4$ billion images.
This \emph{combinatorial explosion}
(or the \emph{curse of dimensionality}) has blocked many brute-force attempts to make sense of the data so a different scientific approach is called for \cite{berisha2021digital}.
\myskip
  
All concepts and results in this book are introduced for very general data, such as discrete sets of points, and hence are relevant to many applied areas, e.g. point clouds in Computer Vision and Graphics.
However, since our latest work is joint with chemists and biologists, our motivations and examples will include data at the atomic scale, including molecules, atomic clouds, and solid crystalline materials (periodic crystals).
\myskip
  
Some physical objects can be exactly represented in a digital form, for example, by listing the coordinates of atoms.
This atomistic representation much better describes a real molecule than pixel-based images of a car.
Though we can exactly describe a molecule by the positions of its atoms, is this description unambiguous? 
To better understand the underlying obstacles, we split this question into more questions below.
\myskip

\nt
\textbf{First question}: which objects are the \emph{same or different?}
Indeed, if we shift all atomic positions by a fixed vector, the digital representation changes, but does the underlying object remain \emph{the same}?
The question is deceptively simple, but because the definition of a crystal structure was incomplete in practice \cite{brock2021change}, this problem has attracted considerable attention, even appearing in the titles of papers \cite{sacchi2020same}.
The missing ingredient was the concept of an equivalence, which should accompany all newly defined objects.
\myskip

\nt
\textbf{Second question}: \emph{if different, by how much}?
Indeed, all real data is uncertain at least due to measurement noise.
Moreover, all atoms vibrate so their relative positions are always uncertain.
This basic fact in Richard Feynman's first lecture on physics \cite[Chapter 1 ``Atoms in motion'']{feynman2011lectures} called for a continuous quantification of similarities.
The resulting problem is algorithmically difficult even for macroscopic objects.
Indeed, when walking or driving, our brains (but not computers) easily recognise obstacles whose visual representations change in our moving coordinate system.
If a car moved or the wind slightly deformed a bush, humans can still identify them as perturbations of the original objects, while a computer program needs an exact formula for a distance. 
\myskip

\nt
\textbf{Third question}: \emph{where can we find new objects?}
Discovery sciences, such as molecular or materials synthesis, struggle to find or even recognise new objects in the vast chemical space. 
Indeed, all known molecules, say for a fixed number $m$ of atoms, live in a common space of $m$-atom configurations.
This space potentially contains unknown molecules, which have not yet been discovered, but where should we look for them?
Humans faced similar challenges in their early exploration of our planet to discover new places to live and thrive.
The slightly rephrased question \emph{where do all real objects live?} requires us to build a geographic-style map on a space of all potential objects of a given type.
\myskip

It took cartographers over two centuries (1400--1600), during the Age of Geographic Discoveries \cite{arnold2013age}, to build a map of the Earth based on latitude and longitude coordinates. 
Scientists can build geographic-style maps of continuous spaces for other real objects.  
\myskip

\begin{figure}[H]
\centering
\includegraphics[width=\textwidth]{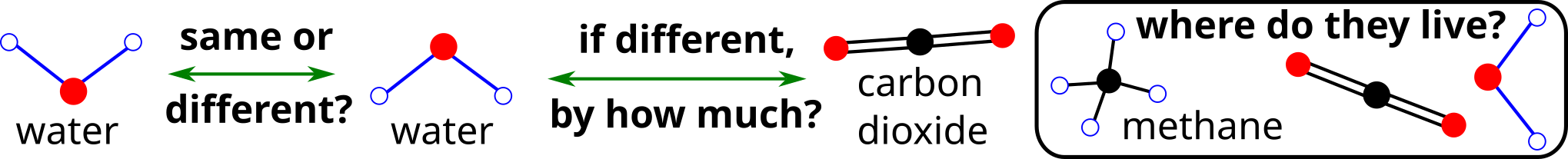}
\caption{The main questions of Geometric Data Science are illustrated for molecules: \ce{H2O}, \ce{CO2}, \ce{CH4}.}
\label{fig:same-different-molecules}
\end{figure}

In summary, \emph{Geometric Data Science} (GDS) aims to mathematically formalise and answer the following questions for real data under practical equivalences, see Fig.~\ref{fig:same-different-molecules}.
\sskip

\nt
\tb{The first question}:
Same or different?
\sskip

\nt
\tb{The second question}:
If different, by how much?
\sskip

\nt
\tb{The third question}:
Where do all (known and new) real objects live?

We finish this section by describing a few conventions used in the book.
All spelling is British as we are based in the UK.
All acronyms that are harder to guess than the UK are listed at the end of the book before the index.
All environments are highlighted in the bold font and numbered according to sections, e.g. Definition~\ref{dfn:equivalence} is followed by Example~\ref{exa:non-equivalences} in section~\ref{sec:equivalences}.
All figures and tables are numbered consecutively within each chapter, as Fig.~\ref{fig:same-different-molecules}.
All environments have the following end symbols:
\sskip

\nt
$\square$ for proofs;
\sskip

\nt
$\blacktriangle$ for definitions;
\sskip

\nt
$\Diamondblack$ for examples and remarks;
\sskip

\nt
$\bigstar$ for problems and conjectures;
\sskip

\nt
$\blacksquare$ for theorems, corollaries, propositions, and lemmas.
\myskip

All new concepts in definitions are highlighted in the \emph{italic} font, which is also used for emphasising keywords.
$\R^n$ denotes
the Euclidean $n$-dimensional space
with a fixed coordinate system of the standard orthonormal basis and origin $0$.
Any vector $\vec p$ with real coordinates $p_1,\dots,p_n$can be positioned with the tail at $0\in\R^n$ and the head at
the point $p=(p_1,\dots,p_n)\in\R^n$, i.e. $p$ and $\vec p$ are often used interchangeably.

\section{Abstract and practical equivalence relations on data objects}
\label{sec:equivalences}

This section recalls an equivalence relation and various types of invariant under a given equivalence. 
These concepts will help formalise the first question: \emph{same or different?} 
\myskip

For any type of object, such as real numbers or all finite sets $A$ of unordered points in $\R^n$, a binary relation describes ordered pairs $(A,B)$ that satisfy this relation.
If $x,y$ are real numbers, one simple relation is the strict inequality $x<y$.

\index{equivalence relation}
\index{transitivity axiom}
 
\begin{dfn}[equivalence relation]
\label{dfn:equivalence}
A binary relation between objects of a given type is called an \emph{equivalence} and denoted by $\sim$ if the following axioms hold:
\sskip

\nt
\tb{(a)}
\emph{reflexivity:} 
any object is equivalent to itself, i.e. $A\sim A$;
\sskip

\nt
\tb{(b)}
\emph{symmetry:} 
for any objects $A,B$, if $A\sim B$, then $B\sim A$;
\sskip

\nt
\tb{(c)}
\emph{transitivity:} 
for any objects $A,B,C$, if $A\sim B$ and $B\sim C$, then $A\sim C$.
\myskip

\nt
Any object $A$ generates its \emph{equivalence class} $[A]=\{\text{all objects equivalent to } A\}$.
\edfn
\end{dfn}

One widely used equivalence on real numbers is the usual equality ($=$), which can be extended to vectors (points in $\R^n$), matrices, and multisets of elements with multiplicities or weights.
The axioms in Definition~\ref{dfn:equivalence} justify a classification under a given equivalence as a splitting or \emph{partition} into disjoint equivalence classes such that every objects belongs to exactly one class.
Indeed, any such classes, say $[A]$ and $[C]$ share a common object $B$, i.e. $A\sim B$ and $B\sim C$, then $A\sim C$ and hence $[A]=[C]$ due to the transitivity axiom.
If a classification is based on a labelled dataset, often with labels produced by humans or a computer program, 
this finite classification is hard to extend to many other real objects, which are called ``out-of-distribution''. 
\myskip

\begin{exa}[non-equivalences]
\label{exa:non-equivalences}
\tb{(a)}
The relation $x<y$ (strict inequality) on real numbers is not an equivalence, because the reflexivity axiom fails: $x<x$ is false. 
\myskip

\nt
\tb{(b)} 
The relation $x\leq y$ (non-strict inequality) on real numbers is not an equivalence.
Though the reflexivity holds, the symmetry is also expected to hold for all real $x,y$ but fails for any non-equal numbers: $x\leq y$ does not imply that $y\leq x$ for $x\neq y$.
\myskip

\nt
\tb{(c)} 
For any fixed real $\ep>0$, the relation $|x-y|\leq\ep$ ($\ep$-closeness) satisfies the reflexivity and symmetry but fails the transitivity axiom.
For instance, if $x=-\ep$, $y=0$, and $z=\ep$, then $|x-y|=\ep=|y-z|$, but $|x-z|=2\ep$.
\eexa
\end{exa}

Example~\ref{exa:non-equivalences}(c) illustrates the \emph{sorites} paradox \cite{hyde2011sorites}\index{sorites paradox}, which has been discussed since ancient times: ``does a heap of sand remain a heap if grains of sand are removed one by one?''
Removing one grain of sand plays the role of an $\ep$-perturbation applied to a data object, such as a heap of millions of grains.
Such a single grain can be considered similar to an outlier in data.
If we are allowed to remove a point from a given set, such as an outlier, without noticing any difference, then all point sets can be made equivalent. 
\smallskip

Similarly, if we assume that a given object remains the same (equivalent to the original one) under all perturbations up to any tiny threshold $\ep>0$, the transitivity axiom will imply that sufficiently many perturbations can make all objects equivalent.  
For instance, if we are ignore slight deviations of vertices in a triangle, the resulting classification of triangles becomes trivial, consisting of a single class of all triangles.
Hence, noise in real data cannot be ignored but should be properly measured.
\myskip

A simple example of an equivalence (not restricted to a fixed dataset) is an equality for a specific property.
For instance, two finite sets $A,B$ can be called equivalent if they have the same size: $|A|=|B|$.
However, many objects that share one property might differ in other properties,
For finite objects, the equivalence relation defined by their size is \emph{weak} in the sense that many substantially different objects have the same size, e.g. molecules of the same number of atoms, belong to the same equivalence class.
\myskip

We will look for a \emph{stronger} equivalence that better separates given objects.
For objects that are more complicated than points in $\R^n$, such as sets of points, the identity relation is overkill in practice (too strong), because shifting all points of a cloud changes only its coordinate representation rather than physical properties.
\myskip
 
Though many equivalence relations make sense for real objects, such as molecules or materials, one equivalence relation stands out in our world: a rigid motion usually preserves all meaningful properties and hence is the strongest relation for most applications.  
When comparing physical objects, the first thing people try to do is to superimpose them by rigid motion. 
Even if a given object, such as a human hand or a molecule, is intrinsically flexible (non-rigid), its different rigid conformations (classes under rigid motion) often have different properties and hence should be reliably distinguished.
\myskip
 
Recall that a \emph{basis} of $\R^n$ consists of $n$ vectors $\vec v_1,\dots,\vec v_n$ such that any $\vec v\in\R^n$ can be written as a linear combination $\vec v=\sum\limits_{i=1}^n t_i\vec v_i$ for some $t_1,\dots,t_n\in\R$.
The basis vectors are \emph{linearly independent} in the sense that if $\sum\limits_{i=1}^n t_i\vec v_i=0$, then $t_1=\cdots=t_n=0$.
\myskip

Later, we will discuss more technical concepts, such as the determinant, which has an algebraic definition.
In the geometric spirit of this book, we only mention here that the \emph{determinant} $\det(\vec v_1,\dots,\vec v_n)$ of the $n\times n$ matrix with columns $\vec v_1,\dots,\vec v_n$ is the signed volume of the parallelepiped on the edge vectors $\vec v_1,\dots,\vec v_n$.
In particular, any vectors $\vec v_1,\dots,\vec v_n\in\R^n$ form a linear basis of $\R^n$ if and only if $\det(\vec v_1,\dots,\vec v_n)\neq 0$.
Our default distance between any points $\vec a=(a_1,\dots,a_n)$ and $\vec b=(b_1,\dots,b_n)$ in $\R^n$ is Euclidean, denoted as $|\vec a-\vec b|=\sqrt{\sum\limits_{i=1}^n (a_i-b_i)^2}$.

\index{rigid motion}
\index{isometry}
\index{dilation}
\index{homothety}

\begin{exa}[rigid motion, isometry, dilation, and homothety]
\label{exa:isometry}
\tb{(a)}
The \emph{translation} along a fixed vector $\vec v\in\R^n$ is the map $T[\vec v]:\R^n\to\R^n$ such that $T[\vec v](\vec p)=\vec p+\vec v$.
\myskip

\nt
\tb{(b)}
A rotation in $\R^n$ is a linear map $R[Q]:\R^n\to\R^n$, $R[Q](\vec p)=Q\vec p$ represented by a \emph{special orthogonal} $n\times n$ matrix $Q$ that has the determinant $\det(Q)=1$ and satisfies $Q^T Q=I_n=Q^T Q$ is the identity matrix, where $Q^T$ is the transpose of $Q$.
All such matrices $Q$ form the \emph{special orthogonal} group $\SO(\R^n)$.
\myskip

\nt
\tb{(c)}
A \emph{rigid motion} $f:\R^n\to\R^n$ is a composition of translations and rotations in $\R^n$ and can be written as $f(\vec p)=Q\vec p+\vec v$ for any $\vec p\in\R^n$, a fixed $\vec v\in\R^n$ and $Q\in\SO(\R^n)$.  
Any sets $A,B\subset\R^n$ that are related by rigid motion are called \emph{rigidly equivalent} (denoted by $A\cong B$).
All rigid motions in $\R^n$ form the \emph{Special Euclidean} group $\SE(\R^n)$.
\myskip

\nt
\tb{(d)}
The \emph{mirror reflection} relative to an $(n-1)$-dimensional hyperspace $H\subset\R^n$ with a normal vector $\vec v$ is the map defined by $f[H](p)=p$ (any point $p\in H$ is fixed) and $f[H](\vec u)=-\vec u$ for any vector $\vec u$ parallel to $\vec v$.
\myskip

\nt
\tb{(e)}
A \emph{Euclidean isometry} is any map $f:\R^n\to\R^n$ preserving Euclidean distance, i.e. $|f(\vec a)-f(\vec b)|=|\vec a-\vec b|$ for any vectors $\vec a,\vec b\in\R^n$.
Alternatively, any \emph{Euclidean isometry} is a composition of a rigid motion and a mirror reflection, and can be written as $f(\vec a)=Q\vec a+\vec v$ for any $\vec a\in\R^n$, a fixed vector $\vec v\in\R^n$, and an \emph{orthogonal} matrix $Q$ satisfying $Q^T Q=I_n=Q^T Q$.
Any subsets $A,B\subset\R^n$ are related by isometry are called \emph{isometric} (denoted by $A\simeq B$).
All isometries in $\R^n$ form the \emph{Euclidean} group $\Eu(\R^n)$.
\myskip

\nt
\tb{(f)}
For a fixed factor $s>0$, the \emph{uniform scaling} is the map $u:\R^n\to\R^n$, $u(\vec a)=s\vec a$ for any vector $\vec a\in\R^n$.
A \emph{dilation} is a composition of a rigid motion and a uniform scaling.
A \emph{homothety} is a composition of an isometry and a uniform scaling.
\eexa
\end{exa}

Any rigid motion $f$ preserves \emph{orientation} of $\R^n$, which can be defined as the sign of the determinant of the $n\times n$ matrix consisting of the columns $f(\vec v_1),\dots,f(\vec v_n)$, where $\vec v_1,\dots,\vec v_n$ is a basis of $\R^n$.
A mirror reflection, for example, changing the sign of the first coordinate (the reflection relative to the hyperspace $a_1=0$) is not a rigid motion, because the orientation is changed.
Hence, compositions or rigid motion with mirror reflections form a slightly wider collection of equivalences, which do not distinguish mirror images. 
The identities $Q^T Q=I_n=Q Q^T$ imply that $\det(Q)=\pm 1$.
All such orthogonal matrices form the \emph{orthogonal} group $O(\R^n)$.
We avoid the notation $O(n)$, which will be later used to denote a linear-time complexity of algorithms.

\index{equivalence relation}
\begin{dfn}[weaker vs stronger equivalences]
\label{dfn:stronger-weaker}
For a fixed collection of objects, one equivalence relation $\sim_1$ is (non-strictly) \emph{weaker} than another $\sim_2$ (then $\sim_2$ is called \emph{stronger} than $\sim_1$) if any objects equivalent under the stronger relation $\sim_2$ are equivalent under the weaker relation, i.e. $A\sim_2 B$ always implies that $A\sim_1 B$.
\edfn
\end{dfn}

If one equivalence $\sim_1$ is weaker than $\sim_2$, then the stronger equivalence $\sim_2$ refines the partition into equivalence classes defined by the weaker equivalence $\sim_1$.
\smallskip

Rigid motions, isometries, dilations, and homotheties define equivalence relations in the sense of Definition~\ref{dfn:equivalence}.
Among these four equivalences, rigid motion is the strongest ($\cong$).
Isometry ($\simeq$) is slightly weaker because any pair of mirror images is in the same isometry class, not necessarily in the same class under rigid motion.
\myskip

Dilation is weaker than rigid motion because all uniformly scaled objects belong to the same class.
Homothety is the weakest of the four so that any objects related by isometry or dilation are homothetic.
A substantially weaker equivalence is defined by bijection, which is a 1-1 map between all points of two objects.
\myskip

After an equivalence is fixed, the next challenge is to classify all given objects under this equivalence.
Such a classification should answer the first question (\emph{same or different?}) by a practical algorithm that determines whether given objects are equivalent or not.
A mathematically justified approach to any classification is to develop invariant descriptors that can reliably distinguish objects under a given equivalence, as defined below.  
Invariant values can be numbers, vectors, matrices, or more complicated objects in a metric space that should still be easier to compare than the original ones.  

\index{invariant}
\index{complete invariant}
\begin{dfn}[invariants and complete invariants] 
\label{dfn:invariants}
\tb{(a)}
Fix an equivalence on some objects.
An \emph{invariant} $I$ is a function that takes the same value on all equivalent objects, i.e. $A\sim B$ implies that $I(A)=I(B)$.
Alternatively, if $I(A)\neq I(B)$, then $A\not\sim B$. 
In other words, $I$ is a descriptor with \tb{no} \emph{false negatives} defined as pairs $A,B$ that represent equivalent objects $A\sim B$ but have different values of this descriptor.
\myskip

\nt
\tb{(b)}
An invariant $I$ is called \emph{complete} if $I$ distinguishes all non-equivalent objects, i.e. if $A\not\sim B$, then $I(A)\neq I(B)$.
Alternatively, if $I(A)=I(B)$, then $A\sim B$, i.e. $I$ takes the same value only on equivalent objects.
In other words, $I$ has \tb{no} \emph{false positives} defined as pairs of non-equivalent $A\not\sim B$ that are indistinguishable by $I$, i.e. $I(A)=I(B)$.
\edfn
\end{dfn}

A constant function $I$ taking the same value on all objects satisfies Definition~\ref{dfn:invariants} but does not distinguish any objects.
We will always assume that an invariant is not the same for all equivalence classes.
Then the implication $[I(A)\neq I(B)] \Ra [A\not\sim B]$ allows us to distinguish complicated) objects by using simpler invariants. 

\begin{exa}[invariants vs non-invariants]
\label{exa:invariants}
\tb{(a)}
A simple invariant of a finite set $A$ under bijection or any stronger equivalence is the size of $A$, which we denote by $|A|$.
\myskip

\nt
\tb{(b)}
For a finite set $A\subset\R^n$, 
the \emph{centre of mass} $\bar A=\dfrac{1}{|A|}\sum\limits_{p\in A}p$ is not invariant of $A$ even under translations and rotations in $\R^n$, and hence under all weaker equivalences, including rigid motion, isometry, and bijection.
\myskip

\nt
\tb{(c)}
For sets of two ordered points $x,y\in\R$, the difference $x-y$ is invariant under rigid motion (only translations in $\R$), but not under isometry that can swap the order of $x,y$.
The Euclidean distance $|x-y|$ is a complete invariant of two ordered point sets under isometry in $\R$, but not under rigid motion.
Indeed, the ordered pairs $(0,1)$ and $(1,0)$ of numbers are not rigidly equivalent but have the same inter-point distance 1.
\eexa
\end{exa}

A complete invariant $I$ fully answers the first main question (\emph{same or different?}) by checking if $I(A)=I(B)$, which is equivalent to $A\sim B$ by~Definition~\ref{dfn:invariants}.
\myskip

Any function $h$ generates its equivalence relation on the domain where $h$ is defined: $A\sim_h B$ if and only if $h(A)=h(B)$.
Then $h$ is a complete invariant under its equivalence $\sim_h$.
For a fixed collection of objects, invariants can be compared by strength similar to equivalence relations in Definition~\ref{dfn:stronger-weaker}.
For ordered pairs $(x,y)$, the difference $x-y$ is a stronger invariant than the distance $|x-y|$.
For a fixed equivalence relation, a complete invariant is the strongest one among all invariants under this equivalence.

\section{Distance metrics on invariant values and equivalence classes}
\label{sec:metrics}

To rephrase the second main question (\emph{if different, by how much?}) in mathematical terms, this section introduces a distance metric between arbitrary objects, which can be equivalence classes or values of an invariant under a given equivalence. 

\index{metric}
\index{Hausdorff distance}
\index{bottleneck distance}
\begin{dfn}[metrics and pseudo-metrics]
\label{dfn:metrics}
\tb{(a)}
A real-valued function $d$ on pairs of objects under an equivalence relation $\sim$ is a \emph{metric} if these axioms hold:
\sskip

\noindent
(1) \emph{coincidence:} $d(A,B)=0$ if and only if $A\sim B$;
\sskip

\noindent
(2) \emph{symmetry:} $d(A,B)=d(B,A)$ for any objects $A,B$; 
\sskip

\noindent
(3) $\triangle$ 
\emph{triangle inequality:} 
 $d(A,B)+d(B,C)\geq d(A,C)$ for any objects $A,B,C$.
\myskip
 
\nt 
\tb{(b)}
If axiom (1) is replaced with the weaker version 
$(1')$ $d(A,A)=0$ for any $A$, 
then non-equivalent objects $A\not\sim B$ can have $d(A,B)=0$, and $d$ is called a \emph{pseudo-metric}.
\edfn
\end{dfn}

The axioms in Definition~\ref{dfn:metrics}(a) imply the non-negativity of a metric as follows: $2d(A,B)=d(A,B)+d(B,A)\geq d(A,A)=0$.
The word ``metric'' is often used in applications and evaluation functions that depend on a single object.
We emphasise that all metrics measure a distance between two objects.
The concept of a distance becomes more general if some of the three metric axioms are weakened \cite{deza2009encyclopedia}.

\index{metric space}
\index{cloud}
\begin{dfn}[metric spaces and clouds]
\label{dfn:clouds}
\tb{(a)}
Any set $M$ of objects with a metric $d:M\times M\to\R$ is called a \emph{metric space}.
\myskip
 
\nt
\tb{(b)}
A \emph{cloud} is any finite set of unordered points in a metric space.
A \emph{Euclidean} cloud is any finite set $A\subset\R^n$ of unordered points with the Euclidean distance.
\edfn
\end{dfn}

We will usually consider metrics on invariant values rather than on original objects.
If $I$ is a complete invariant under a given equivalence, then any metric on invariant values is the metric on equivalence classes of original objects so that $d(I(A),I(B))=0$ is equivalent to $A\sim B$.
However, if $I$ is incomplete, then $d(I(A),I(B))=0$ only guarantees that $I(A)=I(B)$, not necessarily $A\sim B$.
Then any metric on (values of) an incomplete invariant defines only a pseudo-metric on equivalence classes of original objects in the sense of Definition~\ref{dfn:metrics}(b). 
Example~\ref{exa:metrics} introduces well-known metrics on vectors and arbitrary subsets in a metric space.

\index{Minkowski metric}
\index{Chebyshev metric}
\index{Hausdorff distance}
\index{bottleneck distance}

\begin{exa}[Minkowski metrics, Hausdorff and bottleneck distances]
\label{exa:metrics}
\tb{(a)}
Fix a real parameter $q\in[1,+\infty)$.
For any points $\vec a=(a_1,\dots,a_n)$ and $\vec b=(b_1,\dots,b_n)$ in $\R^n$,
the \emph{Minkowski metric} is $L_q(a,b)=\left(\sum\limits_{i=1}^n |a_i-b_i|^q\right)^{1/q}$.
In the limit case $q=+\infty$, the  metric is defined as $L_\infty(a,b)=\max\limits_{i=1,\dots,n} |a_i-b_i|$, also called the \emph{Chebyshev} metric.
\myskip

\nt
\tb{(b)}
Let $A,B$ be subsets of a space $X$ with a metric $d$.
The distance from $a\in A$ to $B$ is $d(a,B)=\inf\limits_{b\in B} d_X(a,b)$.
The \emph{directed} distance is $d_H(A,B)=\sup\limits_{a\in A} d(a,B)$.
The \emph{Hausdorff} distance is $\HD(A,B)=\max\{d_H(A,B),d_H(B,A)\}$.
The \emph{bottleneck distance} 
$\BD(A,B)=\inf\limits_{g:A\to B} \sup\limits_{p\in A}d(g(p),p)$ is minimised for all bijections $g:A\to B$.
\eexa
\end{exa}

If there are no bijections $A\to B$, one can set $\BD(A,B)=+\infty$, so $\BD$ is a well-defined metric only on subsets that allow bijections. 
In Example~\ref{exa:metrics}(a), the parameters $q=1,2,+\infty$ define the metrics that are also called Manhattan (sum metric), Euclidean, and Chebyshev (max metric), respectively.
We will often consider Minkowski metrics for all parameters $q\in[1,+\infty]$, including the limit case $q=+\infty$.
\myskip

For any $\ep\geq 0$ and a subset $C$ of a metric space $M$, the $\ep$-offset of $C$ consists of all points $q\in M$ at a maximum distance $\ep$ from $C$, i.e. $d(q,C)\leq\ep$.
\myskip

In other words, the $\ep$-offset of $C$ is the union of closed balls with the radius $\ep$ and centres at all points $p\in C$.  Then the Hausdorff distance $\HD(A,B)$ can be defined as the minimal $\ep\geq 0$ such that the $\ep$-offset of $A$ covers $B$ and the $\ep$-offset of $B$ covers $A$.
\myskip

\begin{figure}[h!]
\centering
\includegraphics[width=\textwidth]{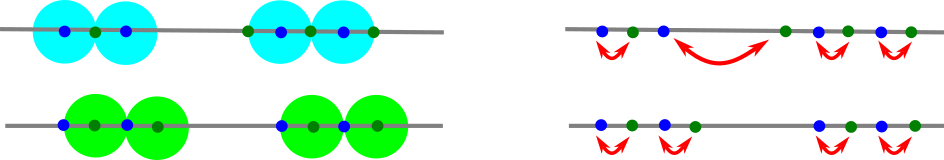}
\caption{
\textbf{Left}: 
in the Euclidean line $\R$. the clouds $A$ of 4 green points and $B$ of 4 blue points have a small Hausdorff distance $\HD$.
\textbf{Right}: the same clouds $A,B\subset\R$ have a large bottleneck distance $\BD$ based on a bijection $g:A\to B$ (shown by red arrows), which minimises the maximum deviation of points, see Example~\ref{exa:metrics}(b). }
\label{fig:Hausdorff-vs-bottleneck}
\end{figure}

The Hausdorff distance $\HD$ is illustrated in terms of $\ep$-offsets in Figure~\ref{fig:Hausdorff-vs-bottleneck}, where a single ball around a blue or green point can cover a cluster of several points from a different subset of another colour.
The bottleneck distance $\BD$ is stricter by measuring a minimum required deviation for a bijective matching of points, as in Figure~\ref{fig:Hausdorff-vs-bottleneck}~(right).

Definition~\ref{dfn:metrics} allows a \emph{discrete} metric that takes a constant value on all non-equivalent objects , e.g. $d(A,B)=1$ for all $A\not\sim B$, and $d(A,B)=0$ for all $A\sim B$.
However, this discrete metric is purely theoretic because all real objects slightly differ due to noise, so this metric would almost always have the same value 1.
Definition~\ref{dfn:Lipschitz} formalises the practically useful continuity under perturbations.
\myskip

For simplicity, we will consider a collection $X$ of objects that allow bijections between each other.
For finite sets, this restriction means that all clouds from $X$ have the same size.
Bijections always exist between any infinite \emph{discrete subsets} of $\R^n$ (that have a positive minimum inter-point distances), because such subsets are countable.  

\index{Lipschitz continuity}
\begin{dfn}[Lipschitz continuity]
\label{dfn:Lipschitz}
Let $X$ be a space of objects with a distance metric $d_X$.
Let $I:X\to Y$ be a function to a space $Y$ with a metric $d_Y$, e.g. an invariant under a given equivalence.
Then $I$ is \emph{Lipschitz continuous} with a Lipschitz constant $\la$ if, for any $\ep\geq 0$, the $\ep$-closeness $d_X(A,B)\leq\ep$ implies that $d_Y(I(A),I(B))\leq\la\ep$.
\edfn  
\end{dfn}

For any discrete subsets in a metric space, the condition $\BD(A,B)\leq\ep$ can be replaced with the following: $B$ is obtained from $A$ by perturbing every point of $A$ within its $\ep$-neighbourhood.
For molecules and materials consisting of atoms, such $\ep$-perturbations and hence
the bottleneck distance are motivated by thermal vibrations and experimental noise under which atoms can slightly change their positions but cannot disappear.
For other applications, e.g. to point clouds in Computer Vision, Definition~\ref{dfn:Lipschitz} can be made stricter by requiring the Lipschitz continuity under less restrictive perturbations of input data in the Hausdorff distance.
\myskip

The classical continuity in terms of $\ep,\de$ is much weaker than the more practical Definition~\ref{dfn:Lipschitz}.
For instance, the function $y=\dfrac{1}{x}$ is continuous for all $x\neq 0$ but has no Lipschitz constant because $\dfrac{1}{x}\to\infty$ as $x\to 0$.
Hence, almost any function in practice can be called continuous in the weak sense of $\ep,\de$ away from singular points.
\myskip

The Lipschitz continuity brings a physical meaning due to an explicit constant $\la$.
For example, if any atom is perturbed up to $\ep$, any inter-atomic distance changes up to $2\ep$ due to the triangle inequality of Euclidean distance.
Hence, inter-atomic distances in physically meaningful units have Lipschitz constant $\la=2$. 

\section{The geo-mapping problem for objects under equivalence}
\label{sec:problem}

This section introduces auxiliary concepts of a metric moduli space and a computational complexity before stating the main geo-mapping problem.
Several non-trivial cases of this problem will be solved in later chapters for point clouds and periodic point sets.  

\index{moduli space}
\begin{dfn}[metric moduli spaces] 
\label{dfn:moduli}
Let $X$ be a collection of objects with the bottleneck distance $\BD$ and an equivalence relation $\sim$.
The \emph{metric moduli space} is the set $Y=X/\sim$ of all equivalence classes $[A]$ for $A\in X$, equipped with a distance metric $d$ satisfying all metric axioms in Definition~\ref{dfn:metrics}(b) so that the \emph{class map} $X\to Y$ defined by $A\mapsto [A]$ is Lipschitz continuous in the sense of Definition~\ref{dfn:Lipschitz}.
If an equivalence is defined by an action of a group $G$ on $X$, the moduli space is also denoted by $X/G$.
\edfn
\end{dfn}

Moduli spaces can also be called quotient spaces or spaces of orbits.
The adjective \emph{metric} means that a moduli space is equipped with a metric satisfying Definition~\ref{dfn:metrics}(a).
\myskip

Algebraic geometry studied moduli spaces in more general settings \cite{tschinkel2006geometry}, usually for varieties defined by polynomial equations and considered 
under actions of linear groups, not involving permutations and metrics.
However, the recently emerged area of Metric Algebraic Geometry \cite{breiding2024metric} started to explore metrics on moduli spaces.
Geometric Data Science (GDS) goes further by requiring polynomial-time algorithms for complete invariants with continuous metrics, and adds the realisability and Euclidean embeddability to parametrise moduli spaces similar to geographic maps of Earth.

\index{cloud space}
\index{Cloud Rigid Space}
\index{Cloud Isomery Space}
\begin{dfn}[cloud spaces $\CIS(\R^n;m)$ and $\CRS(\R^n;m)$] 
\label{dfn:cloud_spaces}
\tb{(a)}
For the collection $S(\R^n;m)$ of all $m$-point sequences $p_1,\dots,p_m\in\R^n$, the moduli space $S(\R^n;m)/\SE(\R^n)$ was previously called a \emph{shape space} $\Sigma_n^m$ \cite{kendall2009shape}. 
Under the extra action of the permutation group, the moduli space $S(\R^n;m)/(\SE(\R^n)\times S_m)=\Sigma_n^m/S_m$ will be called the \emph{Cloud Rigid Space} and denoted by $\CRS(\R^n;m)$. 
\myskip

\nt
\tb{(b)}
Under isometry, not distinguishing mirror images, the space $S(\R^n;m)/(\Eu(\R^n)\times S_m)$ will be called the \emph{Cloud Isometry Space} and denoted by $\CIS(\R^n;m)$.
\edfn
\end{dfn}

Definition~\ref{dfn:cloud_spaces} is motivated by the fact that points are unordered (unlabelled) in most practical scenarios.
Then the equivalence relation is defined by the actions of $\Eu(\R^n)$ and the permutation group $S_m$.
Though atoms in molecules are labelled by chemical elements and sometimes electric charges, many simple molecules such as benzene \ce{C6H6} consist of many indistinguishable atoms, whose permutation group $S_6\times S_6$ consists of $(6!)^2$ (more than half a million) permutations. 
These challenges motivated geo-mapping problems for finite and periodic sets of unordered points. 
\myskip

For all 2-point clouds (sets of two distinct unordered points) in $\R^n$, their inter-point distance $d$ is a complete invariant under isometry.
In this case, the Cloud Isometry Space $\CIS(\R^n;2)$ is the interval $(0,+\infty)$ parametrised by the distance $d$.
\myskip

For any Euclidean cloud $A$, its \emph{input size} $|A|$ is the number of points, because the required computer memory is proportional to $|A|$, for a fixed dimension $n$.
If $A$ is a subset of a metric space, then $A$ can be given by a distance matrix of size $|A|^2$.  
\myskip

If $A\subset\R^n$ is a periodic set of points, its input size $|A|$ can be defined as the number of points in a minimal cell whose periodic translations define the infinite set $A$.
\myskip

All computational complexities will be considered in the Random Access Memory (RAM) model, where any numerical value can be accessed in a constant time.
 
\index{polynomial-time complexity} 
\begin{dfn}[the big $O$ notation for computational complexities] 
\label{dfn:complexity}
Let an algorithm have an input size $m$.
For a function $f(m)$, an algorithm has the \emph{computational complexity} $O(f(m))$ if the total number of required operations, including additions, multiplications, and evaluations of elementary functions 
has an upper bound $c f(m)$ for a constant $c$ and all sufficiently large $m$.
If $f$ is a linear or polynomial function of $m$, the resulting algorithms have a \emph{linear} or \emph{polynomial} time, respectively. 
\edfn
\end{dfn}

For ordered points $p_1,\dots,p_n\in\R^n$, their distance matrix can be computed in time $O(n^2)$, because we need only $\dfrac{n(n-1)}{2}$ distances $|p_i-p_j|$ for $1\leq i<j\leq n$.

\index{embedding}
\begin{dfn}[homeomorpism and embedding]
\label{dfn:embedding}
A \emph{homeomorphism} $f$ is a bi-continuous bijection, i.e. both $f$ and $f^{-1}$ are continuous.
Then an \emph{embedding} $f:X\to Y$ is a homeomorphism on image, i.e. $X\to f(X)$ is a homeomorphism.
\edfn
\end{dfn}

Since all our spaces have metrics, the goal is to guarantee the Lipschitz continuity in the sense of Definition~\ref{dfn:Lipschitz} so that all 
embeddings have Lipschitz constants.  
\myskip

While we state Problem~\ref{pro:geocodes} in full generality below, it can help to keep in mind the partial case of 3-point clouds (triangles) under Euclidean isometry in the plane $\R^2$.

\index{geocode}
\index{invariant}
\index{complete invariant}
\index{Geo-Mapping Problem}
\begin{pro}[Geo-Mapping Problem]
\label{pro:geocodes}
For any space $X$ of objects with a metric $d_X$ and an equivalence relation $\sim$, design a \emph{geocode} defined as an invariant $I:X\to M$ with values in a metric space $M$ satisfying the following conditions.
\myskip

\nt
\tb{(a)} 
\emph{Completeness:} objects $A, B\in X$ are equivalent ($A\sim B$) if and only if $I(A)=I(B)$.
\myskip

\nt
\tb{(b)} 
\emph{Reconstruction:} any object $A\in X$ can be reconstructed from its invariant value $I(A)\in M$, uniquely under the given equivalence.
\myskip

\nt
\tb{(c)}
\emph{Metric:} 
there is a metric $d_M$ in the \emph{invariant space} $I\{X\}=\{I(A) \vl A\in X\}\subset M$, satisfying all metric axioms in Definition~\ref{dfn:metrics}(a).
\myskip

\nt
\tb{(d)}
\emph{Continuity:} 
$I$ is Lipschitz continuous in the sense of Definition~\ref{dfn:Lipschitz}:
there is a constant $\la>0$ such that, for any $\ep\geq 0$, if $d_X(A,B)\leq\ep$, then $d_M(I(A),I(B))\leq\la\ep$.
\myskip

\nt
\tb{(e)}
\emph{Inverse continuity:} 
there is a constant $\mu>0$ such that, for any $\de\geq 0$, if $d_M(I(A),I(B))\leq\de$, there is an equivalence $f$ satisfying 
$d_X(f(A),B)\leq\mu\de$.
\myskip

\nt
\tb{(f)}
\emph{Realisability:} 
the invariant space $I\{X\}=\{I(A) \vl A\in X\}$ can be parametrised so that we can generate any value $I(A)\in I\{X\}$ realisable by some object $A\in X$.
\myskip

\nt
\tb{(g)}
\emph{Euclidean embedding:} 
the invariant space $I\{X\}=\{I(A) \vl A\in X\}$ with the metric $d$ allows a bi-Lipschitz embedding into a Euclidean space $\R^N$ for some $N$.
\myskip

\nt
\tb{(h)}
\emph{Computability:} 
fix a metric space containing all objects of $X$, then the invariant $I(A)$ in (a), a reconstruction of $A$ from $I(A)$ in (b), the metric $d_M(I(A),I(B))$ in (c), an equivalence $f$ in (e), the generation of a new value in $I\{X\}$, and an embedding $I\{X\}\subset\R^N$ can be algorithmically computed in polynomial times of the input size.
\epro
\end{pro}

Completeness in~\ref{pro:geocodes}(a) formalises the first main question (\emph{same or different?}) in Geometric Data Science by requiring that a geocode is a complete invariant code unambiguously representing any given object.
Completeness alone is impractical because one can define a complete invariant as the entire collection of images $I(A)=\{f(A) \text{ for all equivalences }f\}$, which is infinite for most equivalence relations. 
\myskip

The reconstruction in~\ref{pro:geocodes}(b) is stronger than the completeness, because a complete invariant can be abstract or too complicated without an algorithmic reconstruction.
For instance, a human fingerprint and a genetic code are practically used for identifying humans but are insufficient (yet) to grow a genetic replica of a living person.
\myskip
  
The metric requirements in~\ref{pro:geocodes}(c) are justified by recently designed distances $d$ on point clouds \cite{rass2024metricizing}, which guarantee pre-determined outputs of several clustering algorithms, such as $k$-means and DBSCAN, if a distance $d$ between points is allowed to fail the triangle inequality, even with any small additive error.  
\myskip
 
The Lipschitz continuity in~\ref{pro:geocodes}(d) in the bottleneck distance is motivated for atomic-scale objects by the fact that atoms vibrate around their average positions.
\myskip

The inverse continuity in~\ref{pro:geocodes}(e) 
allows us to deform a geocode $I(A)$ and continuously trace the evolution of reconstructed objects $A$.
Exact values of Lipschitz constants $\la,\mu$ are less important than their existence, because one can always scale down a metric $d$ to make $\la$ in \ref{pro:geocodes}(d) smaller, then the constant $\mu$ in \ref{pro:geocodes}(e) will be larger.
\myskip

Conditions~\ref{pro:geocodes}(c,d,e) formalise the second main question (\emph{if different, by how much?}) by requiring that a geocode $I$ is a bi-continuous invariant. 
\myskip

For most objects in this book, the initial metric $d_X$ will be the bottleneck distance $\BD$.
However, the continuity in conditions \ref{pro:geocodes}(d,e) for $\BD$ might be unrealistic for some infinite objects, such as periodic lattices.
In this case, an initial space $X$ will consist of finite inputs with bottleneck-type metrics, e.g. lattice bases under isometry.
\myskip

The realisability in~\ref{pro:geocodes}(f)
 justifies the name \emph{geocode} as an analogue of geographic coordinates and requires an explicit description of all realisable values in the invariant space $I\{X\}$ similar to all hospitable places on Earth.
This realisability in~\ref{pro:geocodes}(f) is stated in purely mathematical terms and can be extended for practical applications by additionally requiring that $I(A)$ is realised by a physical object $A$.
For instance, a distance $d$ between atoms cannot be any positive number, so 2-atom molecules have this distance $d$ in a small range within the full moduli space $\CIS(\R^3;2)=(0,+\infty)$.
\myskip

The Euclidean embeddability in~\ref{pro:geocodes}(f) converts $I(A)$ into a vector in some in $\R^N$ with usual Euclidean distance, which can be used as an input of machine learning algorithms.
However, the space exploration, such as a deformation or sampling of invariant values, should be performed in the invariant space $I\{X\}$, because the complement 
$\R^N\setminus I\{X\}$ consists of artificial values that are not realisable by any objects.
\myskip
 
The polynomial-time computability in~\ref{pro:geocodes}(h) glues all previous conditions and makes Geo-Mapping Problem~\ref{pro:geocodes} notoriously hard even for finite sets of $m\geq 4$ unordered points under isometry in $\R^2$, which were classified only into discrete types such as
squares and parallelograms, through they live in a continuous 5D space.
\myskip

A full solution to Problem~\ref{pro:geocodes} enables a continuous exploration of a complicated moduli space $X/\sim$ via polynomial-time geocodes on the invariant space $I\{X\}$.
Machine learning often relies on \emph{latent} spaces of descriptor values, which can be ambiguous due to non-invariance or incompleteness, or discontinuous under noise. 
\myskip

Geometric Data Science aims to replace all latent spaces with \emph{invariant spaces} of real data objects, which should be continuously parametrised by fast geocodes.
\myskip

In the geographic analogy, a example geocode of any position on Earth (considered as a round sphere) consists of the latitude and longitude coordinates in the realisable ranges $[-90^\circ,90^\circ]$ and $(-180^\circ,+180^\circ]$.
More exactly, the (interior of the) rectangle $R=[-90^\circ,90^\circ]\times [-180^\circ,+180^\circ]$ continuously maps to sphere $S^2$.
This parametrisation assumes that the horizontal edges of $R$ (all geocodes with a latitude $90^\circ$ and $-90^\circ$) map to the north and south pole, respectively.
We should also glue the vertical edges of $R$ (all geocodes with longitudes $\pm 180^\circ$ and a fixed latitude) to a single meridian of $S^2$.
\myskip

Apart from these boundary identifications, geocodes have real values in known ranges and have enabled navigation on Earth.
Indeed, the shortest way from the US to Japan is to cross the International Date Line over the Pacific Ocean, where the longitude changes from $-180^\circ$ to $+180^\circ$.
Hence, complete invariants of data objects become much more valuable with a continuous metric to find shortest paths in a moduli space.
\myskip

The vision of Geometric Data Science is to develop such geographic-style maps (briefly, \emph{geomaps}) for moduli spaces of all real objects under practical equivalences.
\myskip

These geomaps have analytically defined invariant coordinates and substantially differ from outputs of dimensionality reduction algorithms for the following reasons.
\myskip

Firstly, many such algorithms are stochastic in the sense that their outputs for the same input data can differ on runs with random seeds or on different machines.
\myskip

Secondly, even if a dimensionality reduction is deterministic, such as Principal Component Analysis, the underlying algorithm is data-driven in the sense that adding new data changes the output projection of all data.
Moreover, the coordinates of the resulting projections are so complicated that it is impractical to write them down. 
\myskip

Thirdly, any dimensionality reduction as a function $h:\R^{m}\to\R^n$ for $m>n\geq 1$ is either discontinuous, i.e. makes close points distant, or collapses an unbounded region of $\R^m$ to a single point, i.e. loses an infinite amount of data \cite{landweber2016fiber}.
Hence, dimensionality reductions can produce nice pictures, but a justified analysis of similarities and differences should use invariants and distances in the original high-dimensional space.
\myskip

When we choose 2 or 3 invariants for a low-dimensional projection of a geomap, we know all other skipped invariants and hence can expand any cluster or hot spot from the first projection in other coordinates. 
Most importantly, adding new data to geomaps keeps the locations of all past data similar to mapping new places on a geographic map, because the invariant coordinates are defined in a data-independent way. 
\myskip

We considered the name \emph{Metric Data Science} since it is similar to \emph{Metric Algebraic Geometry} \cite{breiding2024metric}.
The progress beyond metrics towards \emph{geo}graphic-style maps \cite{bright2023geographic} of moduli spaces motivated the extra prefix in the name of \emph{Geo}metric Data Science. 

\section{Solutions to the geo-mapping problem in the simplest cases}
\label{sec:solutions}

This section discusses Examples~\ref{exa:geocodes_n=1}-\ref{exa:cyclic_polygons}, which solve Problem~\ref{pro:geocodes} for finite sets of unordered points under rigid motion and isometry in $\R^n$ for the simplest known cases of dimension $n=1$, up to $m=3$ points in $\R^n$, and for cyclic polygons in $\R^2$.

\index{geocode}
\begin{exa}[geocodes of finite point sets in $\R$]  
\label{exa:geocodes_n=1}
In dimension $n=1$, any finite set $A\subset\R$ consists of naturally ordered points $p_1<p_2<\cdots<p_m$.
\myskip

\nt
\tb{(a)}
Since any rigid motion in $\R$ is a translation, the first point $p_1$ can be fixed at the origin $0\in\R$.
Then the sequence $p_1<\cdots<p_m$ is uniquely determined by the geocode $I(A)=(d_1,\dots,d_{m-1})$ of the $m-1$ distances $d_i=p_{i+1}-p_{i}$, where the only realisability condition is $d_i>0$, $i=1,\dots,m-1$.
The Cloud Rigid Space is $\CRS(\R;m)=\R_+^{m-1}$.
\myskip

\nt
\tb{(b)}
Any isometry in $\R$ is a translation or its composition with the reflection $x\mapsto -x$, which reverses the order of all points of $A$ and the order of the distances $d_i$, so that $(d_1,\dots,d_{m-1})\mapsto(d_{m-1},\dots,d_1)$.
Under isometry in $\R$, the geocode $I(A)$ is the unordered pair of these distance vectors in $\R_+^{m-1}$
The Cloud Isometry Space is $\CIS(\R;m)=\R_+^{m-1}/\sim$, where the equivalence relation $\sim$ reverses the order of all coordinates.
If $m=2$, then $\CIS(\R;2)=(0,+\infty)=\CRS(\R;2)$.
If $m=3$, then $\CIS(\R;2)=\{(x,y)\in\R^2 \vl 0<x<y\}$ for $x=\min\{d_1,d_2\}$, $y=\max\{d_1,d_2\}$.  
\eexa
\end{exa}

\index{geocode}
\begin{exa}[geocodes for $m=2$ points in $\R^n$]  
\label{exa:geocodes_m=2}
\tb{(a)}
For pairs of unordered points $p,q$ in $\R^n$, one complete invariant under isometry is the inter-point distance $d=|p-q|$, because we can fix $p$ at the origin $0\in\R^n$ by translation and then apply rotation from $\Or(\R^n)$ to the point $q$ at the distance $d$ in the positive 1st coordinate axis of $\R^n$.
\sskip

The distance $d$ has Lipschitz constant $\la=2$, because perturbing each of the points $p,q$ up to $\ep$ changes their distance $d=|p-q|$ up to $2\ep$ due to the triangle inequality.
\sskip

To check the inverse continuity, let 2-point clouds $A=\{p,q\}$ and $B=\{u,v\}$ in $\R^n$  have $\de$-close distances $d(A)=|\vec p-\vec q|$ and $d(B)=|\vec u-\vec v|$ so that $|d(A)-d(B)|=\de$.
Let $f$ be the isometry that translates the point $p$ to $u$ and then rotates the vector $\vec q-\vec p$ around the point $u$ with an orthogonal matrix from $\Or(\R^n)$ to make $f(\vec q-\vec p)$ parallel to the fixed vector $\vec v-\vec u$.
Since $f(p)=u$, the difference of parallel vectors can be estimated by the difference of their lengths: $|f(\vec q)-\vec v|=\big|\, |f(\vec q-\vec p)|-(\vec v-\vec u)| \,\big|=|d(A)-d(B)|=\de$, so the image $f(A)$ is $\de$-close to $B$.
We can even additionally shift the 2-point cloud $f(A)$ along the straight line through the points of $B$ to put each point of (the image of) $A$ at a distance of $\dfrac{\de}{2}$ from its closest point of $B$.
Hence, the Lipschitz constant in condition~\ref{pro:geocodes}(e) is $\mu=\dfrac{1}{2}$ so that $\la\mu=1$.
The realisability condition for an inter-point distance is $d>0$.
The moduli space $\CIS(\R^n;2)=(0,+\infty)$ is embedded in $\R$. 
\myskip

\nt
\tb{(b)}
If the given points $p,q$ are ordered, all conclusions in part (a) remain valid in dimensions $n\geq 2$, also under rigid motion instead of isometry, because any vectors in $\R^n$ can be made parallel by rigid motion.
In $\R$, the complete invariant of two ordered points $p,q$ under rigid motion (translation) is the difference $p-q$. 
The Cloud Rigid Space $\CRS(\R;2)=\R\setminus\{0\}$ excludes the degenerate case of identical points $p=q$.
\eexa
\end{exa}

\begin{figure}[h!]
\centering
\includegraphics[height=27mm]{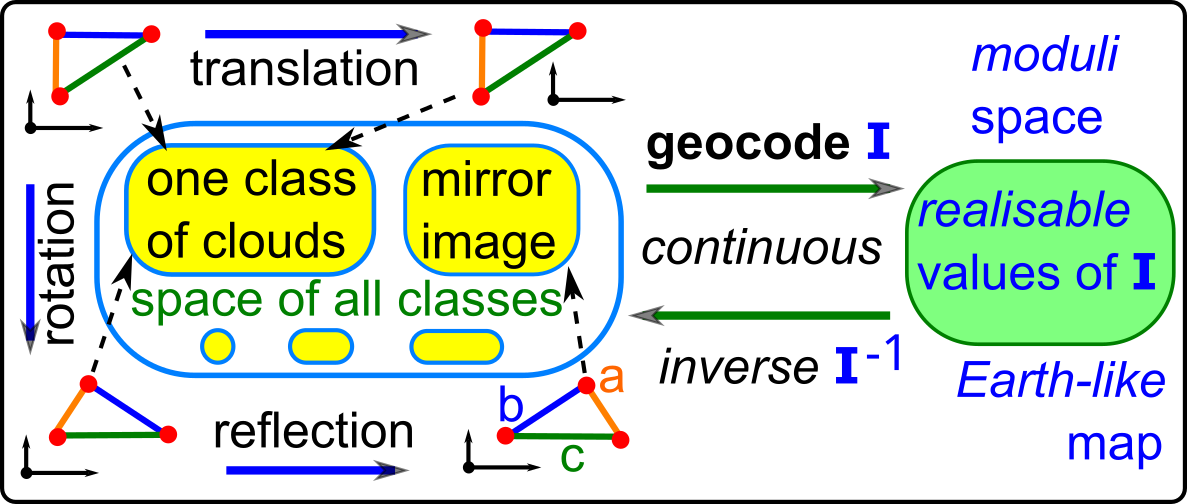}
\hspace*{2mm}
\includegraphics[height=27mm]{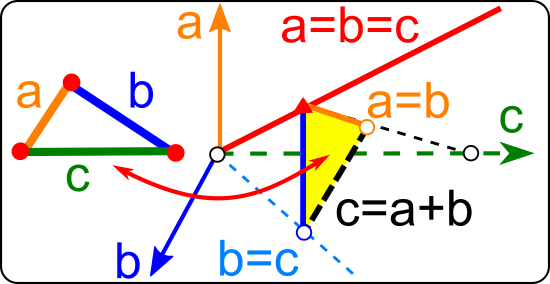}
\caption{
\textbf{Left}: a geocode $I$ from Problem~\ref{pro:geocodes} is illustrated for triangles (3-point clouds) whose isometry classes form a moduli space, which can be mapped like Earth.
\textbf{Right}: the Cloud Isometry Space $\CIS(\R^n;3)$ is 
continuously parametrised by triples of inter-point distances $0<a\leq b\leq c\leq a+b$.}
\label{fig:geomaps-triangles}
\end{figure}

Fig.~\ref{fig:geomaps-triangles}~(left) illustrates Geo-Mapping Problem~\ref{pro:geocodes} by  geocodes parametrising geographic-style maps for moduli spaces of 3-point clouds (triangles) under isometry. 

\index{geocode}
\begin{exa}[geocodes for $m=3$ points in $\R^n$]  
\label{exa:geocodes_m=3}
\tb{(a)}
The side-side-side theorem in Euclidean geometry says that any triangles are congruent (isometric, in our language) if and only if they have the same triple of side lengths, under permutations.
\sskip

An isometry in the plane 
can reverse orientation, so the vertices of a triangle are unordered.
Since a triangle is considered a cloud of 3 unordered points, its three inter-point distances can be written in increasing order, say $0<a\leq b\leq c$.
\sskip

The side-side-side theorem implies that the ordered triple $(a,b,c)$ is a complete invariant of 3 unordered points under isometry in the plane and hence in any $\R^n$.
Similar to Example~\ref{exa:geocodes_m=2},
the distances $a,b,c$ are continuous under perturbations with Lipschitz constant $\la=2$.
The inverse continuity in \ref{pro:geocodes}(e) is harder and will be tackled in forthcoming work.
The only realisability condition is the single triangle inequality $c\leq a+b$, which is the upper bound for the largest distance.
Then the Cloud Isometry Space $\CIS(\R^n;3)$ is the triangular cone $\{(a,b,c)\in\R^3 \vl 0<a\leq b\leq c\leq b+c\}$.
\sskip

Under homothety (isometry composed with uniform scaling), this cone projects to the smaller moduli space represented by the yellow triangle in Fig.~\ref{fig:geomaps-triangles}~(right).
\sskip

The red diagonal $\{a=b=c\}$ represents all equilateral triangles.
The boundary planes $\{a=b\leq c\}$ and $\{a\leq b=c\}$ represent two types of isosceles triangles: ``more horizontal'' and ``more vertical'', respectively.
The third boundary plane $\{c=a+b\}$ represents degenerate triangles of three points in a straight line.
\myskip

\nt
\tb{(b)}
Under rigid motion in any $\R^n$ for $n\geq 3$, all conclusions remain valid,
because a mirror reflection in $\R^2$ can be realised by a rotation in $\R^n$ with a matrix from $\SO(\R^n)$.
Then, for $n\geq 3$, the Cloud Rigid Space $\CRS(\R^n;3)$ is the same cone as in part (a).
\sskip

In the plane $\R^2$, any rigid motion preserves the cyclic order of 3 points (vertices of a triangle).
Then the complete invariant is a triple $(a,b,c)$ of inter-point distances satisfying $\max\{a,b\}\leq c\leq a+b$, so the shorter distances $a,b$ can be in any order.
The resulting Cloud Rigid Space $\CRS(\R^2;3)$ is obtained by gluing two copies of triangular cones $\CIS(\R^2;3)$ along their boundaries, where 3-point clouds are mirror-symmetric.
\myskip

While $\CIS(\R^n;3)$ is embedded in $\R^3$ for any $n\geq 2$, Euclidean embeddability of $\CRS(\R^2;3)$ needs a higher-dimensional space $\R^N$ in condition~\ref{pro:geocodes}(g). 
\eexa
\end{exa}

\index{geocode}
\begin{exa}[geocodes of cyclic polygons in $\R^2$]  
\label{exa:cyclic_polygons}
A polygon polygon in $\R^2$ is \emph{cyclic} if its vertex set $A$ is a subset of a circle of a radius (say) $r>0$.
Then all points of $A$ are cyclically ordered along this circle, say as $p_1,\dots,p_m$.
Since the centre of the circle can be fixed at $0\in\R^2$,
the set $A$ is determined, uniquely under rotation from $\SO(\R^2)$, by the $m$ cyclically ordered inter-point distances $b_i=|p_{i+1}-p_i|$, $i=1,\dots,m$, where $p_{m+1}=p_1$. 
The geocode of $A$ under rigid motion in  $\R^2$ is the sequence $I(A)=(b_1,\dots,b_m)$ under cyclic permutations.
\cite[Theorem~1.8]{penner2012decorated} provides the realisability condition: $b_j<\sum\limits_{i\neq j} b_i$ for $j=1,\dots,m$.
Under isometry in $\R^2$, the sequence $(b_1,\dots,b_m)$ should also be considered under reversing the order:
$(b_1,\dots,b_m)\mapsto(b_m,\dots,b_1)$.
\eexa
\end{exa}

Euclid might have drawn a geomap of triangles from Example~\ref{exa:geocodes_m=3}(a) on sand more than 2000 years ago. 
Hence, it was surprising that even the case of $m=4$ unordered points in $\R^2$ remain opened until complete, Lipschitz continuous, and polynomial-time isometry invariants were developed in 2023 for $m$ unordered points in any $\R^n$ \cite{widdowson2023recognizing}.  


\section{Related areas and connections of Geometric Data Sciences}
\label{sec:related_areas}

This section briefly relates Geometric Data Science to other areas in mathematics, data science, and computer science.
Later chapters will review past work on specific data.
\myskip

Our classifications of geometric objects by invariants were inspired by the famous question \emph{Can we hear the shape of a drum?} \cite{kac1966can}, which has the negative answer in terms of 2D polygons indistinguishable by spectral invariants \cite{gordon1992isospectral,gordon1992one}. 
Problem~\ref{pro:geocodes} went beyond complete classifications for even better invariants that satisfy extra conditions in \ref{pro:geocodes}(b-h) to allow us not only `hear' but more fully `sense' geometric shapes, e.g. equivalence classes under rigid motion in any $\R^n$.
Problem~\ref{pro:geocodes} can be informally rephrased as a short question: \emph{can we sense the shape of a real object?}
\myskip

\index{Geo-Mapping Problem}
Though Geo-Mapping Problem~\ref{pro:geocodes} will be re-phrased for various discrete objects, the original statement covers all possible data under arbitrary equivalence relations.
In this book, our data objects will be finite and periodic sets of unordered points, which represent atoms in molecules and materials.
Other important objects include embedded  graphs, polygonal surface meshes, and simplicial complexes.
Our standard equivalences are rigid motion, isometry, and their compositions with uniform scaling.
Weaker but still practical equivalences are defined by affine, projective, and conformal maps, or actions of specific linear groups on subsets of $\R^n$ as in classical algebraic geometry.
\myskip

\index{Geo-Mapping Problem}
The generality of Geo-Mapping Problem~\ref{pro:geocodes} and the recent progress in the practical cases of finite and periodic point sets justified the birth of Geometric Data Science \cite{widdowson2023recognizing,kurlin2024mathematics,anosova2025recognition}
 as a new area on the interface between metric geometry and data science.
\myskip

Statistics, data analysis, and shape analysis considered a similar object-oriented approach \cite{marron2021object}, often for continuous shapes, such as curves and surfaces, under more complicated equivalences, including re-parametrisations or diffeomorphisms.
\myskip

The carefully written book ``Object-oriented data analysis'' discussed the concepts of equivalence relations and classes (fibres or orbits) in \cite[section~1.2.2]{marron2021object} in the case of triangles under congruence as in Example~\ref{exa:geocodes_m=3}(a), though without mapping the Cloud Isometry Space $\CIS(\R^n;3)$ as in Fig.~\ref{fig:geomaps-triangles}(b).
Since the keyword \emph{invariant} appeared once 
in \cite[p.184]{marron2021object}, section~\ref{sec:equivalences} added more motivations and examples of invariants.
\myskip

Several books on classical invariant theory \cite{olver1999classical,kraft2000classical,dolgachev2003lectures} discuss invariants in the context of algebraic geometry for solutions of polynomial equations, usually under actions of linear groups with elements in $\Q$ or $\C$, not involving translations or permutations, as discussed with
B.Hasset \cite{hassett2018stable} and  F.Kirwan \cite{mumford1994geometric} in private communications.  
\myskip

In real algebraic geometry, the invariants are also studied for (semi-)algebraic sets \cite{bochnak2013real,scheiderer2024course}, including important computational aspects \cite{derksen2001computational,theobald2024real}. The recent book ``Metric Algebraic Geometry'' \cite{breiding2024metric} moved beyond invariants towards distance metrics after the earlier workshop on 
``Emerging applications of algebraic geometry'' \cite{putinar2008emerging}.
\myskip

The key difference of this book from algebraic geometry is the focus on discrete sets of unordered points coming from real data, such as atomic configurations.
\myskip

Discrete sets of points have been studied in discrete geometry \cite{liberti2017euclidean} and rigidity theory \cite{alfakih2018euclidean}, often for ordered points.
This case also has practical applications in proteins, though unordered (unlabelled) sets are more common in real data.
\myskip

Topological Data Analysis \cite{edelsbrunner2010computational,dey2022computational} developed persistent homology summarising the evolution of complexes built on discrete data \cite{carlsson2021topological,joharinad2023mathematical}.
For point clouds, the resulting persistence diagrams are invariants under isometry \cite{smith2024generic}, which are usually computed in dimensions 0 and 1 due to the high complexity and turned out to be weaker than previously anticipated.
\myskip

In Computer Science, the interest in the geometry of data has risen due to the influential area of ``Geometric Deep Learning'' \cite{bronstein2017geometric,bronstein2021geometric} 
advocating for the invariance of inputs or outputs under actions of $\SE(\R^3)$ or $\Eu(\R^3)$ in machine learning algorithms.
\myskip

\begin{figure}[h!]
\centering
\includegraphics[height=38mm]{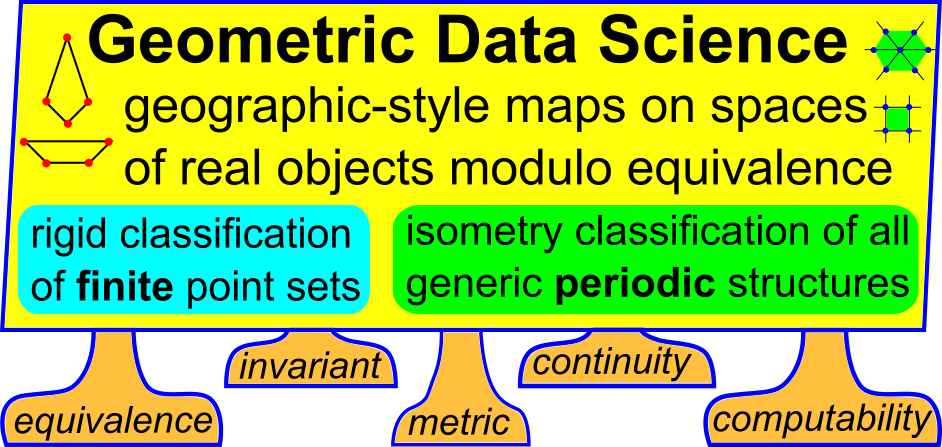}
\hspace*{1mm}
\includegraphics[height=38mm]{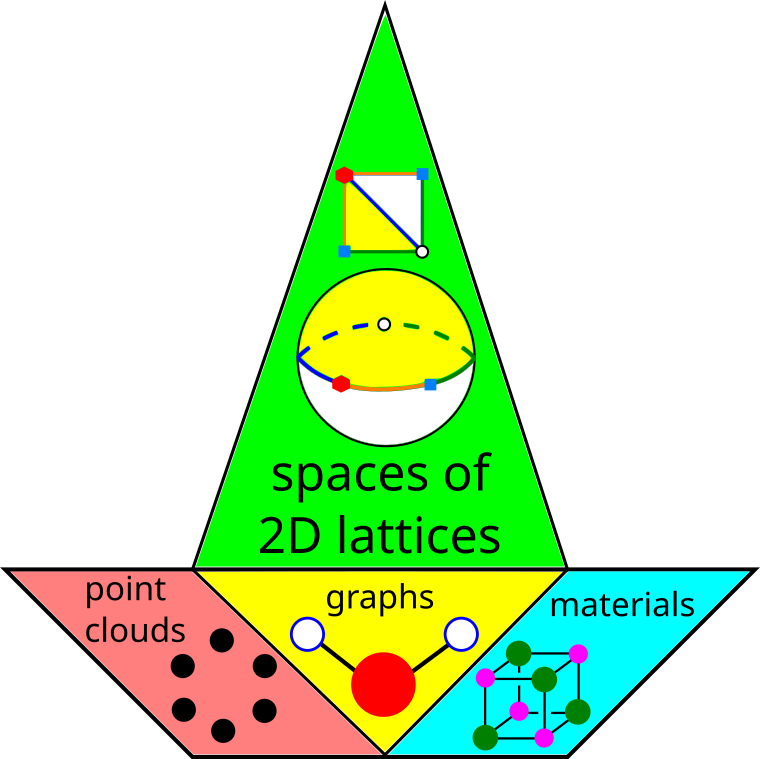}
\caption{
\textbf{Left}: the key concepts 
are introduced in Definitions~\ref{dfn:equivalence}, 
\ref{dfn:invariants}, \ref{dfn:metrics}, \ref{dfn:Lipschitz}, and \ref{dfn:complexity}, all linked in Problem~\ref{pro:geocodes}.
\textbf{Right}: the main objects are finite and periodic sets of unordered points, including lattices in $\R^2$ whose space under rigid motion was the first solution to Problem~\ref{pro:geocodes}.
}
\label{fig:GDS-foundations}
\end{figure}

All these related developments essentially inspired the new area of Geometric Data Science, whose foundational concepts are highlighted in Fig.~\ref{fig:GDS-foundations}~(left).
\myskip

The logo-style image in 
Fig.~\ref{fig:GDS-foundations}~(right) 
shows two quadrilaterals (a vertical kite in green and yellow, and a horizontal trapezium in red, yellow, and blue) whose vertex sets (clouds of 4 unordered points) are indistinguishable by 6 pairwise distances.
\myskip

\section{The chapter plan: from easier to more challenging data}
\label{sec:plan}

This section outlines the plan of all further chapters, which are split into two big parts: finite point sets and periodic point sets, which will be mostly studied under 
isometry.
In the first part, Chapters 2-6 solve partial cases of Problem~\ref{pro:geocodes} for finite point sets.
\myskip

\nt
\tb{Chapter~\ref{chap:proteins}} 
discusses complete, bi-continuous, and linear-time invariants \cite{anosova2025complete} for  
finite sets of ordered points under rigid motion in $\R^3$, which exposed thousands of duplicate chains in the Protein Data Bank within a few hours on a modest desktop computer.
\myskip

\nt
\tb{Chapter~\ref{chap:directional}} leverages Principal Component Analysis to get complete and polynomial-time invariants for finite clouds of unordered points under rigid motion in $\R^n$.
\myskip

\nt
\tb{Chapter~\ref{chap:PDD-finite}} 
introduces Pointwise Distance Distribution.
The $\PDD$ is a fast and generically complete isometry invariant of finite and periodic sets of unordered points.
This chapter proves 
that the PDD is complete for all 4-point clouds under isometry in $\R^n$.  
\myskip

\nt
\tb{Chapter~\ref{chap:SDD}} 
refines the isometry invariant $\PDD$ to a stronger Simplexwise Distance Distribution ($\SDD$) for a finite cloud of unordered points in any metric space.
Despite having a higher (polynomial-time) complexity, a simple definition of $\SDD$ allows us to distinguish all known non-isometric clouds in $\R^3$ that have identical $\PDD$s.
\myskip

\nt
\tb{Chapter~\ref{chap:SCD}} 
improves $\SDD$ to a Simplexwise Centred Distribution, which is a complete, Lipschitz continuous and polynomial time invariant of all point clouds under rigid motion in $\R^n$.
The hardest obstacle in the proof of Lipschitz continuity was resolved by a \emph{strength} of a simplex, which is a linear-growth analogue of the simplex volume.
\bskip

\nt
In the second part, Chapters 7-11 solve partial cases of Problem~\ref{pro:geocodes} for periodic sets.
\myskip

\nt
\tb{Chapter~\ref{chap:1-periodic}} 
defines complete invariants for ordered sequences of points (under several versions of isometries in $\R\times\R^{n-1}$) that are periodic along the first coordinate axis.
\myskip

\nt
\tb{Chapter~\ref{chap:lattices2D}} 
expands the classical approaches of Gauss, Lagrange \cite{lagrange1773recherches}, and Delone \cite{delone1934mathematical}, who studied lattices via quadratic forms, and the more recent work of Conway and Sloane \cite{conway1992low} to solve Problem~\ref{pro:geocodes} for all periodic lattices under rigid motion in $\R^2$.
\myskip

\nt
\tb{Chapter~\ref{chap:densities}} 
discusses density functions, which extend the point density of periodic point sets to generically complete invariants under isometry in $\R^3$.
These density functions will be analytically described for all periodic sequences of intervals within $\R$. 
\myskip

\nt
\tb{Chapter~\ref{chap:PDD-periodic}} 
extends $\PDD$s from finite to periodic point sets, proves their generic completeness under isometry in $\R^n$, and describes their asymptotic behaviour. 
\myskip

\nt
\tb{Chapter~\ref{chap:isosets}}
refines the seminal work of Dolbilin, Lagarias, and Senechal \cite{dolbilin1998multiregular} to build a complete invariant \emph{isoset} with Lipschitz continuous metrics, which can be approximated by polynomial-time algorithms for all periodic point sets in $\R^n$. 
\myskip   

\nt
\tb{Chapter~\ref{chap:conclusions}} 
summarises the most significant results, the already verified principles, and includes several open problems, some of which can be accessible to school students, partially inspired by V.I.Arnold's ``Problem for School Pupils'' in \cite[Chapter~6]{arnold2013real}.

\bibliographystyle{plain}
\bibliography{Geometric-Data-Science-book}

%
%
%

\begin{partbacktext}
\part{Geometric Data Science of finite point sets}
\end{partbacktext}

%
%
%

\chapter{Sequences of ordered points under rigid motion in Euclidean spaces}
\label{chap:proteins} 

\abstract{
This chapter discusses complete invariants of finite sets of ordered points under rigid motion in Euclidean space $\R^n$.
After discussing distance matrices, we adapt Geo-Mapping Problem~\ref{pro:geocodes} for protein backbones, which are non-degenerate polygonal chains in $\R^3$.
Further sections describe the linear-time Backbone Rigid Invariant ($\bri$) and all-vs-all comparisons of chains in the Protein Data Bank.
}

\section{Classical invariants and shape spaces of ordered points}
\label{sec:Euclidean_ordered}

This section reviews past approaches to
classify sequences of ordered points $p_1,\dots,p_m$ under rigid motion or isometry in $\R^n$.
In the case of isometry, a complete invariant of the sequence $p_1,\dots,p_m$, known at least since 1935 \cite{schoenberg1935remarks}, is the matrix of pairwise distances.
An alternative (complete) isometry invariant is the Gram matrix of scalar products $\vec p_i\cdot\vec  p_j$ \cite[chapter 2.9]{weyl1946classical}, which can be expressed in terms of the distance matrix and vice versa. 
Since these matrices do not distinguish mirror images, we state the Euclidean version of Problem~\ref{pro:geocodes} for any finite sets of ordered points below.

\begin{pro}[partial case of Problem~\ref{pro:geocodes} for sequences under rigid motion in $\R^n$]
\label{pro:Euclidean_ordered}
Design a map $I$ on finite sets of ordered points in $\R^n$ satisfying the conditions below.
\myskip

\noindent
\tb{(a)} 
\emph{Completeness:} 
any sequences $A,B\subset\R^n$ are related by rigid motion ($A\cong B$)  in $\R^n$ if and only if $I(A)=I(B)$.
\myskip

\noindent
\tb{(b)} 
\emph{Reconstruction:} 
any sequence $A\subset\R^n$ of ordered points can be reconstructed from its invariant value $I(A)$, uniquely under rigid motion.
\myskip

\noindent
\tb{(c)} \emph{Metric:} 
there is a distance $d$ on the invariant space $\{I(A) \vl \text{sequences } A\subset\R^n\}$ satisfying all metric axioms in Definition~\ref{dfn:metrics}(a). 
\myskip

\nt
\tb{(d)} 
\emph{Continuity:} 
there is a constant $\la$ such that, for any $\ep>0$, if $B$ is obtained from $A$ by perturbing every point of $A$ up to Euclidean distance $\ep$, then $d(I(A),I(B))\leq \la\ep$.  
\myskip

\nt
\tb{(e)} 
\emph{Computability:} 
for a fixed dimension $n$, the invariant $I(A)$, and the metric $d(I(A),I(B))$ can be computed in times that depend polynomially on the maximum size $\max\{|A|,|B|\}$ of sets $A,B\subset\R^n$.
\epro 
\end{pro} 

Problem~\ref{pro:Euclidean_ordered} will be solved in section~\ref{sec:Euclidean_finite}.
Further sections of this chapter will solve a non-degenerate version of Problem~\ref{pro:Euclidean_ordered} for protein backbones in $\R^3$, see Problem~\ref{pro:proteins}.
The following concepts will be essentially used later.

\index{affine dimension}
\begin{dfn}[affine dimension]
\label{dfn:aff}
The \emph{affine dimension} $\aff(A)$ of a cloud $A\subset\R^n$ of points $p_1,\dots,p_m$ is the maximum dimension of the vector space generated by all inter-point vectors $\vec p_i-\vec p_j$ for $i,j\in\{1,\dots,m\}$.
\edfn
\end{dfn}

The affine dimension $\aff(A)$ is an isometry invariant and is independent of the order of points of $A$.
Any cloud $A$ of 2 distinct points has $\aff(A)=1$.
Any cloud $A$ of 3 points that are not in the same straight line has $\aff(A)=2$.
Lemma~\ref{lem:distance_realisation} provides a criterion for a matrix to be realisable by squared distances of a cloud in $\R^n$.
\smallskip

\begin{lem}[realisation of distances]
\label{lem:distance_realisation}
\textbf{(a)}
A symmetric $m\times m$ matrix of $s_{ij}\geq 0$ with $s_{ii}=0$ is realisable as a matrix of squared distances between $p_0=0,p_1,\dots,p_{m-1}\in\R^n$ for some $n$ \emph{if and only if} the $(m-1)\times(m-1)$ matrix $G$ of $g_{ij}=\dfrac{s_{0i}+s_{0j}-s_{ij}}{2}$ has only non-negative eigenvalues.
\smallskip

\noindent
\textbf{(b)}
If $G$ has only non-negative eigenvalues, then $\aff(0,p_1,\dots,p_{m-1})$ equals the number $k\leq m-1\leq n$ of positive eigenvalues of $G$. 
In this case, $g_{ij}=p_i\cdot p_j$ define the \emph{Gram matrix} of the vectors $p_1,\dots,p_{m-1}\in\R^n$, which are uniquely determined in time $O(m^3)$ under a map from $\Or(n)$.
\end{lem}
\newcommand{\lemdistancerealisation}{
\begin{proof}[\textbf{Proof of Lemma}~\ref{lem:distance_realisation}]
\textbf{(a,b)}
We extend \cite[Theorem~1]{dekster1987edge} to the case $m<n+1$ and find $p_1,\dots,p_{m-1}\in\R^n$ in time $O(m^3)$, uniquely under an orthogonal map from 
$\Or(n)$.
\myskip

The part \emph{only if} $\Rightarrow$.
Let a symmetric matrix $S$ consist of squared distances between points the $p_0=0,p_1,\dots,p_{m-1}\in\R^n$.
For $i,j=1,\dots,m-1$, the matrix of
$$g_{ij}=\dfrac{s_{0i}+s_{0j}-s_{ij}}{2}=\frac{|p_i|^2+|p_{j}|^2-|p_i-p_j|^2}{2}=p_i\cdot p_j$$ is the Gram matrix, which can be written as $G=P^T P$, where the columns of the $n\times(m-1)$ matrix $P$ are the vectors $p_1,\dots,p_{m-1}$.
For any vector $v\in\R^{m-1}$, we have 
$$0\leq |P v|^2=(Pv)^T(Pv)=v^T(P^T P)v=v^T G v.$$
Since the quadratic form $v^T G v\geq 0$ for any $v\in\R^{m-1}$, the matrix $G$ is positive semi-definite, i.e. $G$ has only non-negative eigenvalues, see \cite[Theorem~7.2.7]{horn2012matrix}.  
\myskip

The part \emph{if} $\Leftarrow$.
For any positive semi-definite matrix $G$, there is an orthogonal matrix $Q$ such that $Q^T G Q=D$ is the diagonal matrix, whose $m-1$ diagonal elements are non-negative eigenvalues of $G$.
The diagonal matrix $\sqrt{D}$ consists of the square roots of the eigenvalues of $G$.
The number of positive eigenvalues of $G$ equals the dimension $k=\aff(\{0,p_1,\dots,p_{m-1}\})$ of the subspace that is linearly spanned by $p_1,\dots,p_{m-1}$. 
\myskip

We may assume that all $k\leq n$ positive eigenvalues of $G$ correspond to the first $k$ coordinates of $\R^n$.
Since $Q^T=Q^{-1}$, the given matrix $G=Q D Q^T=(Q\sqrt{D})(Q\sqrt{D})^T$ becomes the Gram matrix of the columns of $Q\sqrt{D}$.
These columns become the reconstructed vectors $p_1,\dots,p_{m-1}\in\R^n$.
\myskip

If there is another diagonalisation $\ti Q^T G \ti Q=\ti D$ for $\ti Q\in \Or(n)$, then $\ti D$ differs from $D$ by a permutation of eigenvalues, which is realized by an orthogonal map, so we set $\ti D=D$.
Then $G=\ti Q D \ti Q^T=(\ti Q\sqrt{D})(\ti Q\sqrt{D})^T$ is the Gram matrix of the columns of $\ti Q\sqrt{D}$.
The new columns differ from the previously reconstructed vectors $p_1,\dots,p_{m-1}\in\R^n$ by the orthogonal map $Q\ti Q^T$.
Hence the reconstruction is unique under $\Or(n)$-transformations.
Computing eigenvectors $p_1,\dots,p_{m-1}$ needs a diagonalisation of $G$ in time $O(m^3)$, see section~11.5 in \cite{press2007numerical}.
\end{proof}
}\lemdistancerealisation

Chapter~3 in \cite{liberti2017euclidean} discusses realisations of a complete graph given by a distance matrix in $\R^n$.
Lemma~\ref{lem:seq_reconstruction}(a) holds for all clouds, including degenerate ones, e.g. for 3 points in a straight line.
Any points $p_1,\dots,p_{n-1}\in A$ have $\aff(p_1,\dots,p_{n-1})\leq n-2$.
For example, any two distinct points in $A\subset\R^3$ generate a straight line. 

\begin{lem}[sequence reconstruction]
\label{lem:seq_reconstruction}
\textbf{(a)}
Any sequence of ordered points $p_1,\dots,p_m$ in $\R^n$ can be reconstructed (uniquely under isometry) from the matrix of the Euclidean distances $|p_i-p_j|$ in time $O(m^3)$.
If all distances are divided by $R=\max\limits_{i=1,\dots,m}|p_i|$, the reconstruction of $p_1,\dots,p_m$ is unique under isometry and uniform scaling in $\R^n$.
\myskip

\noindent
\textbf{(b)}
If $m\leq n$, the uniqueness of reconstructions in part (a) remains true if we replace isometry with rigid motion in $\R^n$.
\elem
\end{lem}
\newcommand{\lemreconstruction}{
\begin{proof}[\textbf{Proof of Lemma}~\ref{lem:seq_reconstruction}]
\textbf{(a)}
By translation, we can fix $p_1$ at the origin.
Let $G$ be the $(m-1)\times (m-1)$ matrix $g_{ij}=\dfrac{|p_i|^2+|p_j|^2-|p_i-p_j|^2}{2}=p_i\cdot p_j$, where $i,j=2,\dots,m$, which is obtained from the squared distances between the points $p_1=0,p_2,\dots,p_m$.
\sskip

By Lemma~\ref{lem:distance_realisation} if $G$ has $k\leq n$ positive eigenvalues, then $p_1=0,\dots,p_m$ can be uniquely determined under isometry in $\R^k\subset\R^n$ in time $O(m^3)$.
If all distances are divided by the same radius $R=\max\limits_{i=1,\dots,m}|p_i|$, the above construction guarantees uniqueness under isometry and uniform scaling. 
\myskip

\noindent
\textbf{(b)}
If $m\leq n$, any mirror images of $p_1,\dots,p_m\in\R^n$, after a suitable rigid motion, can be assumed to belong to an $(n-1)$-dimensional hyperspace $H\subset\R^n$, where they are matched by a mirror reflection $H\to H$ with respect to an $(n-2)$-dimensional subspace $S\subset H$, which is realized by the $180^\circ$ orientation-preserving rotation around $S$.
\end{proof}
}\lemreconstruction

Lemma~\ref{lem:seq_reconstruction}(b) for $m=n=3$ implies that any triangle is determined by its sides, uniquely under rigid motion in $\R^3$.
For example, sides $3,4,5$ define a right-angled triangle whose mirror images are not related by rigid motion within a plane $H\subset\R^3$, but are matched by a rigid motion in $H$ and a $180^\circ$ rotation of $\R^3$ around a line in $H$.   
\myskip

The difference between the matrices of distances or scalar products can be converted into a continuous metric by taking a matrix norm.
These matrices are preserved under any mirror reflections.
Hence, these invariants are incomplete under rigid motion.
\myskip

One can define the sign of orientation on some (or all subsets of) $n+1$ points from a given sequence.
This extra sign is discrete and vanishes for degenerate configurations for $n+1$ points that affinely span a $k$-dimensional subspace in $\R^n$ for $k\leq n-1$.
\myskip 

Another attempt to satisfy the Lipschitz continuity in Problem~\ref{pro:geocodes}(d) is to multiply the sign of orientation by the volume of the simplex spanned by $n+1$ points, say $p_1,\dots,p_{n+1}$.
In other words, one can take the signed volume of the parallelepiped spanned by the vectors $\vec p_i-\vec p_1$, $i=2,\dots,n+1$.
If the first $n+1$ points are degenerate, then the zero volume of their spanned parallelepiped does not give any extra information to distinguish mirror images of the full sequence.
\myskip

More importantly, the resulting signed volume is not Lipschitz continuous already in dimension $n=2$.
Indeed, let us consider the triangle $A(t)$ on the vertices $(\pm l,0)$ and $(0,t \ep)$, where $l,\ep>0$ are fixed constants ($l$ is large, $\ep>0$ is small), and $t\in[-1,1]$ is a time parameter.  
The signed area of $A(t)$ is $tl\ep$, changing from $-l\ep$ at $t=-1$ to $l\ep$ at $t=1$.
Then the Lipschitz constant cannot be smaller than $\la=\dfrac{2l\ep}{2}=l\ep$.
Then the signed area can have a fixed Lipschitz constant only for bounded triangles, not for all triangles, because we can choose $\ep=\dfrac{1}{\sqrt{l}}$ to make $\la=l\ep=\sqrt{l}$ unbounded.
\myskip

We will resolve this obstacle to Lipschitz continuity by a different function (the strength of a simplex) in a later chapter.
This section finishes by noting that a naive extension of the complete isometry invariant (matrix of distances or scalar products) from the ordered to $m$ unordered points requires $m!$ permutations.
This exponential complexity is ruled out by the polynomial-time requirement in \ref{pro:geocodes}(h).
Indeed, $m!$ becomes too large already for $m=4$ points: $4!=24$ matrices of size $4\times 4$.
\myskip

Each distance matrix is symmetric and has zeros on the diagonal, and hence can be represented by only 6 distances.
However, the total number of $24\times 6=144$ distances seems overwhelmingly unnecessary to unambiguously and continuously encode 4 unordered points under isometry in $\R^2$.
Chapter~4 will prove that a smaller $4\times 3$ matrix invariant is complete for any 4 unordered points under isometry in $\R^n$.
\myskip

In 1977, Kendall \cite{kendall1977diffusion} started to study configuration spaces of ordered points modulo rigid motion in $\R^n$ under the name of \emph{size-and-shape spaces} \cite{kendall2009shape}.
If we consider sequences equivalent also under uniform scaling, the smaller \emph{shape space} $\Si_2^m$ of $m$ ordered points in $\R^2$ can be described as a complex projective space $\C P^{m-1}$ due to the group $\SO(\R^2)$ being identified with the unit circle in the complex space $\C^1=\R^2$. 
However, there is no easy description of the moduli space $\Si_3^m$ of $m$-point sequences in $\R^3$, which has no multiplicative group structure similar to $\R^2=\C^1$.
\myskip

In a general metric space, let a sequence $A$ of $m$ ordered points be given by their $m\times m$ distance matrix $D$.
Multidimensional scaling 
 \cite{kruskal1978multidimensional} finds an embedding $A\subset\R^k$ (if it exists) preserving all distances of $M$ for a minimum dimension $k\leq m$.
The underlying computation of $m$ eigenvalues of the Gram matrix expressed via $D$ needs $O(m^3)$ time.
The resulting representation of $A\subset\R^k$ uses orthonormal eigenvectors whose ambiguity up to signs for potential comparisons leads to the time factor $2^k$, which can be close to $2^m$ and hence exponential in the number $m$ of points.
\myskip

Further sections in this chapter follow papers \cite{anosova2025complete,wlodawer2025duplicate}.

\section{The geo-mapping problem for protein backbones in $\R^3$}
\label{sec:backbones}

\begin{figure}[h!]
\includegraphics[height=23mm]{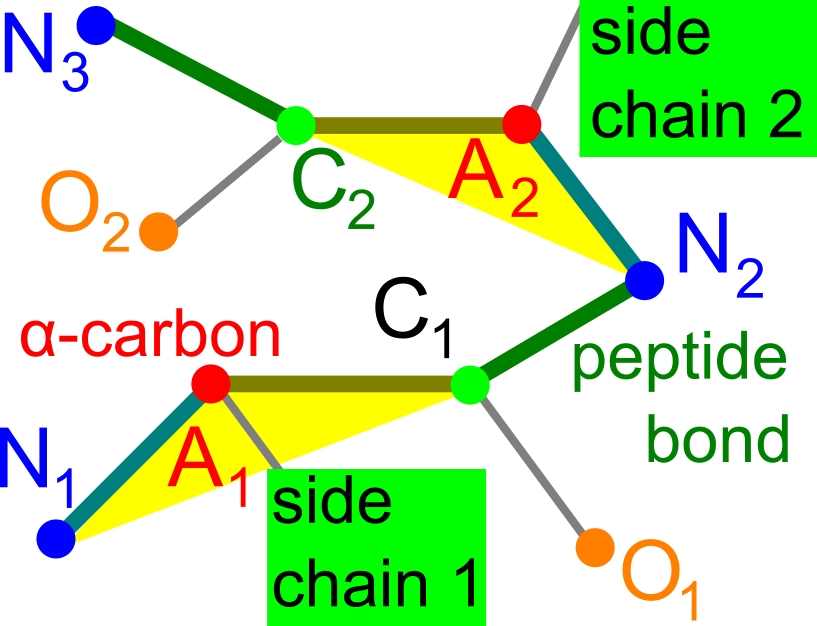}
\hspace*{1mm}
\includegraphics[height=23mm]{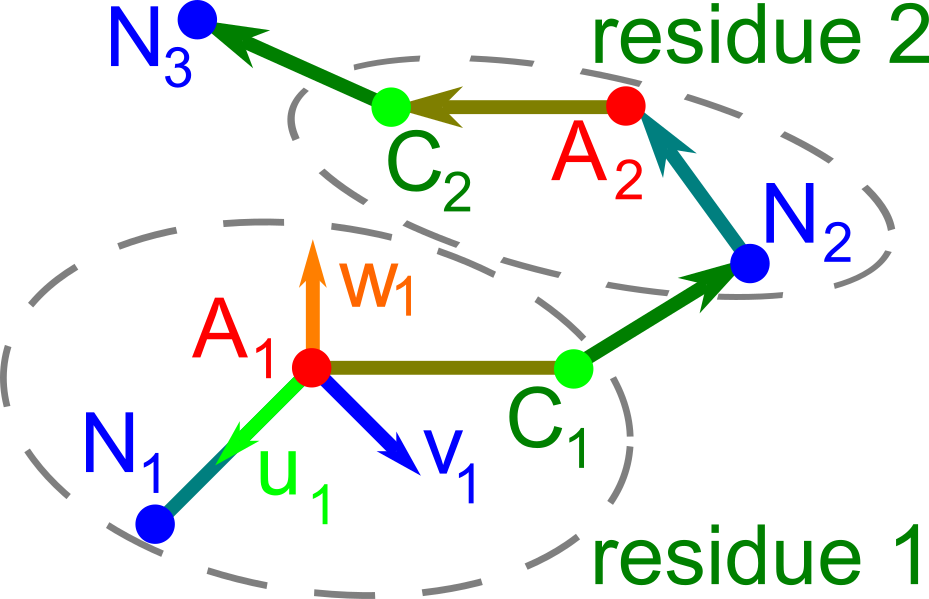}
\hspace*{1mm}
\includegraphics[height=23mm]{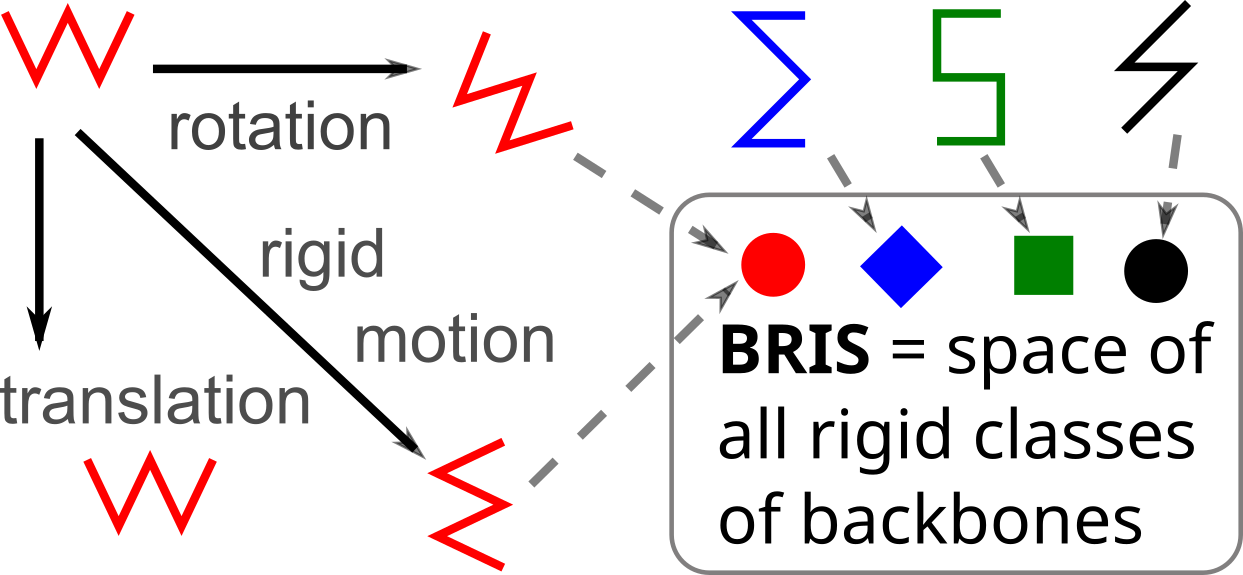}
\caption{\textbf{Left}: 
all main atoms $N_i$, $A_i$, $C_i$  of a protein chain form a \emph{backbone} embedded in $\R^3$.
\textbf{Middle}: each triangle $\triangle N_i A_i C_i$ defines an orthonormal basis $\vb*{u}_i,\vb*{v}_i,\vb*{w}_i$.
The coordinates of the bonds $\ve{C_i N_{i+1}}$, $\ve{N_{i+1}A_{i+1}}$,  $\ve{A_{i+1}C_{i+1}}$ in this basis form the complete Backbone Rigid Invariant $\bri$.
\textbf{Right}: 
All rigidly equivalent backbones form a single \emph{rigid class}.
All rigid classes form the \emph{Backbone Rigid Space}. 
The image schematically illustrates four different classes of simple polygonal chains in $\R^3$.
}
\label{fig:backbones}
\end{figure}

A \emph{protein} is a large biomolecule consisting of one or several chains of amino acid residues.
The \emph{primary structure} (\emph{sequence}) of a protein chain is a string of residue labels (represented by one or three letters), each denoting one of (usually) 20 standard amino acids \cite{scott2017mathematical}. 
The \emph{secondary} structure consists of frequent semi-rigid subchains such as $\al$-helices and $\be$-strands \cite{linderstrom1952lane}. 
A sequence of a protein is relatively easy to experimentally determine but important functional properties, such as interactions with drug molecules, depend on a 3-dimensional geometric \emph{fold} (a \emph{tertiary structure}) represented by an embedding of all its atoms in $\R^3$ , 
see Fig.~\ref{fig:backbones}~(left). 
\medskip

In 1973, Nobel laureate Anfinsen conjectured that the sequence of any protein chain determines its 3D geometric fold \cite{anfinsen1973principles}.
Following this conjecture, neural networks such as AlphaFold2 and RosettaFold \cite{jumper2021highly,baek2021accurate} 
 optimise millions of parameters to predict a protein fold from its sequence, but need re-training \cite{jones2022impact} on the growing Protein Data Bank (PDB), which is considered a `gold standard' for experimental structures \cite{burley2017protein}.
The reported accuracies of prediction are often based on the LDDT (Local Distance Difference Test) \cite[p.~2728]{mariani2013lddt} and TM-score \cite{zhang2004scoring}, which fail the metric axioms.
Then clustering 
can produce pre-determined clusters and may not be trustworthy \cite{rass2024metricizing}.
\medskip

Backbones of the same length (number of residues) can be optimally aligned to minimise the Root Mean Square Deviation (RMSD) between corresponding atoms \cite{holm2024dali}.
This RMSD is slow to compute for all pairs of proteins and gives only distances without mapping the protein universe (moduli space of proteins under rigid motion). 
\medskip

We develop a different approach by mapping the space of protein backbones in analytically defined coordinates similar to geographic-style maps of a new planet.  
\medskip

Any embedded protein in $\R^3$ can be rigidly moved, which changes all atomic coordinates.
However, the underlying structure remains the same in the sense that different images of a protein under rigid motion have the same properties in a fixed environment. 
Though proteins are flexible, it is important to distinguish their rigid structures that can interact differently \cite{heifetz2003effect} with other molecules, including medical drugs.
\myskip


\index{Backbone Rigid Space}
\index{backbone structure}
\begin{dfn}[Backbone Rigid Space $\bris_m$]
\label{dfn:bris}
A protein \emph{backbone} is a sequence of $m$ ordered triplets of main chain atoms (nitrogen $N_i$, $\al$-carbon $A_i$, and carbonyl carbon $C_i$) given by their positions in $\R^3$.
The \emph{structure} of a backbone, or a protein chain (with all side chains), or a biomolecule consisting of several chains is the equivalence class of this geometric object under rigid motion in $\R^3$. 
For any $m\geq 1$, the classes of all backbones of $m$ triplets under rigid motion form the \emph{Backbone Rigid Space} $\bris_m$.
\edfn
\end{dfn}

Backbones were studied by incomplete invariants such as torsion angles, which allow \emph{false positive} pairs of non-equivalent backbones $S\not\cong Q$ with $I(S)=I(Q)$.
Because all atoms in a backbone $S$ are ordered, their distance matrix determines $S\subset\R^3$, uniquely under isometry, but has a large quadratic size in the number $m$ of residues and fails to distinguish mirror images.
Adding a sign of orientation creates discontinuity for polygonal chains that are almost mirror-symmetric. 
\medskip

Problem~\ref{pro:proteins} adapts Geo-Mapping Problem~\ref{pro:geocodes} to  protein backbones.
The completeness in \ref{pro:geocodes}(a) is restricted to polygonal chains, where each triplet of atoms $N_i,A_i,C_i$ is not in a straight line, as we have checked for all experimental structures in the PDB.
The polynomial-time condition in \ref{pro:geocodes}(h) is strengthened to linear time.

\begin{pro}[geo-mapping for protein backbones]
\label{pro:proteins}
For any $m\geq 1$, design a map $I:\bris_m\to\R^N$ for some $N$ satisfying the following conditions.
\medskip

\noindent
\tb{(a)} 
\emph{Completeness:} 
any backbones $S,Q\subset\R^3$ are rigidly equivalent if and only if $I(S)=I(Q)$, i.e. 
$I$ has \emph{no false negatives} and \emph{no false positives}.
\medskip

\noindent
\tb{(b)} 
\emph{Reconstruction:}
any protein backbone $S\subset\R^3$ can be reconstructed from its invariant value $I(S)$ uniquely under rigid motion.
\medskip

\noindent
\tb{(c)} \emph{Metric:} 
there is a distance $d$ on invariant values satisfying all metric axioms in Definition~\ref{dfn:metrics}(a).
\medskip

\noindent
\tb{(d)} \emph{Continuity:} 
there is a constant $\la$ such that, for any $\ep>0$, if $Q$ is obtained from $S$ by perturbing every atom up to Euclidean distance $\ep$, then $d(I(S),I(Q))\leq\la\ep$.
\medskip

\noindent
\tb{(e)} 
\emph{Atom matching:} 
there is a constant $\mu$ such that, for any backbones $S,Q$ 
with $\de=d(I(S),I(Q))$, all their atoms can be matched up to a distance $\mu\de$ by a rigid motion. 
\medskip

\noindent
\tb{(f)} 
\emph{Realisability:} 
the invariant space $I\{X\}=\{I(A) \vl A\in X\}$ can be parametrised so that we can generate any value $I(A)\in I\{X\}$ realisable by some object $A\in X$.
\medskip

\noindent
\tb{(g)} 
\emph{Respecting subchains:} 
for any subchain of residues $R_{i}\cup\dots\cup R_{i+j}$ in a backbone $S$,
the invariant $I(R_{i}\cup\dots\cup R_{i+j})$ can be obtained from $I(S)$ in linear time $O(j)$ with respect to the length of the subchain.
\medskip

\noindent
\tb{(h)} 
\emph{Linear-time computability:} the invariant $I$, the metric $d$, a reconstruction in (b), and a rigid motion in (e) can be computed in time $O(m)$ for any backbone of $m$ residues.
\epro 
\end{pro}

The completeness in \ref{pro:proteins}(a) means that $I$ is the strongest invariant and hence distinguishes all protein backbones that cannot be exactly matched by rigid motion.
The reconstruction in \ref{pro:proteins}(b) is more practical because $I$ may not allow an efficiently computable inverse map $I^{-1}$ from an invariant value $I(S)$ to a backbone $S\subset\R^3$.
\medskip

The continuity in \ref{pro:proteins}(d) fails for invariants based on principal directions that can discontinuously change in degenerate cases when eigenvalues become equal.  
The atom matching in \ref{pro:proteins}(e) says that, after finding a rigid motion $f$ in $\R^3$, any atom $p\in S$ 
has Euclidean distance at most $\mu\de$ to the corresponding atom $q\in f(Q)$. 
\medskip

Conditions~\ref{pro:proteins}(d,e) guarantee the Lipschitz continuity of $I$ and its inverse on the image $I(\bris_m)\subset\R^N$.
Though Lemma~\ref{lem:distance_realisation} gives a two-sided criterion for the realisability of distances by ordered points $p_1,\dots,p_m\in\R^n$, the space of distance matrices is highly singular and cannot be easily sampled.
Since a random matrix of potential distances for $m>n+1$ is unlikely to be realisable by $m$ ordered points in $\R^n$, the realisability condition in \ref{pro:proteins}(g) is non-trivial for the distance matrix.
\myskip

Since Problem~\ref{pro:proteins} asked for an invariant $I:\bris_m\to\R^N$, the Euclidean embeddability in~\ref{pro:geocodes}(g) hold automatically and has been replaced with condition~\ref{pro:proteins}(g), motivated by secondary structures, which are subchains in full backbones.
\medskip

The linear time in \ref{pro:proteins}(h) makes all previous conditions practically useful because even the distance matrix needs $O(m^2)$ time and space, substantially slower than linear time $O(m)$ for thousands of residues.
\myskip

Past work on similarities of proteins is reviewed in \cite[section~2]{anosova2025complete}.
Section~\ref{sec:BRI} introduces the \emph{Backbone Rigid Invariant} $\bri:\bris_m\to\R^{9m-6}$ to 
solve Problem~\ref{pro:proteins} by Theorems~\ref{thm:BRI_complete}, \ref{thm:BRI_continuous}, \ref{thm:BRI_inverse}.
The numerical components of $\bri$ play the role of geocodes, which are geographic-style coordinates on the space $\bris_m$, where any protein backbone has a uniquely defined location.
Section~\ref{sec:duplicates} describes how $\bri$ detected thousands of geometric duplicates in the PDB, some of which need updates.

\section{Complete and bi-continuous Backbone Rigid Invariant}
\label{sec:BRI}

We start with the simpler \emph{triangular invariant} that describes the rigid class of each residue triangle $\triangle N_i A_i C_i$ on three main atoms: nitrogen $N_i$, $\alpha$-carbon $A_i$, and carbonyl carbon $C_i$, for $i=1,\dots,m$, see Fig.~\ref{fig:backbones}~(middle).
For any points $A,B\in\R^3$, let $|\ve{AB}|$ be the Euclidean length of the vector $\ve{AB}$ from $A$ to $B$.
The \emph{scalar} and \emph{vector} products of vectors $\vb*{u},\vec v\in\R^3$ are denoted by $\vb*{u}\cdot\vb*{v}$ and $\vb*{u}\times\vb*{v}$, respectively. 

\index{triangular invariant}
\begin{dfn}[triangular invariant $\trin$]
\label{dfn:trin}
Let a backbone $S\subset\R^3$ have $3m$ ordered atoms $N_i$, $A_i$, $C_i$, $i=1,\dots,m$.
In the plane of $\triangle N_i A_i C_i$, for the 2D basis obtained by Gaussian orthogonalisation of $\ve{A_i N_i},\ve{A_i C_i}$, the vector $\ve{A_i N_i}$ has the coordinates $x(A_i N_i)=|\ve{A_i N_i}|$ and $y(A_i N_i)=0$.
Let $\vec x=\dfrac{\ve{A_i N_i}}{|\ve{A_i N_i}|}$ be the unit vector. 
Then $\ve{A_i C_i}$ has the coordinates $x(A_i C_i)=\ve{A_i C_i}\cdot\vec x$ and
$y(A_i C_i)=
|\ve{A_i C_i} - x(A_i C_i)\vec x|$ 
in the direction orthogonal to $\vec x$.
The \emph{triangular invariant} $\trin(S)$ is the $m\times 3$ matrix whose $i$-th row consists of the coordinates $x(A_i N_i)$, $x(A_i C_i)$, and $y(A_i C_i)$ for $i=1,\dots,m$.
\edfn
\end{dfn}

The $i$-th row of $\trin(S)$ uniquely determines the rigid class of $\triangle N_i A_i C_i$.
\myskip

On May 4, 2024, the PDB had 213,191 entries with 1,091,420 chains. 
Protocol~\ref{prot:cleanPDB} below produced $104,688\approx 49\%$ entries with $707410\approx 65\%$ chains in 4 hours 48 min 11 sec.
All experiments were run on 
CPU Core i7-11700 @2.50GHz RAM 32Gb.

\begin{prot}[selecting a subset of 707K+ chains in the PDB]
\label{prot:cleanPDB}
The PDB was filtered by removing the following entries and individual chains. \\
(1) 4513 non-proteins (the entity is labeled as `not a protein'). \\
(2) 178153 disordered chains, where some atoms have occupancies $<1$. \\
(3) 201648 chains with residues having non-consecutive indices. \\
(4) 9941 incomplete chains missing one of the main atoms $N_i,A_i,C_i$. \\
(5) 4364 chains with non-standard amino acids.
\end{prot}
 
To guarantee new condition~\ref{pro:proteins}(e) respecting subchains,
Definition~\ref{dfn:bri} will represent atoms $N_{i+1},A_{i+1},C_{i+1}$ in a basis of the  previous $i$-th residue.
The first residue needs only three invariants from Definition~\ref{dfn:trin} to determine the rigid class of $\triangle N_1 A_1 C_1$ in $\R^3$. 
Due to cleaning in Protocol~\ref{prot:cleanPDB}, all consecutive atoms along any backbone have distances $d\geq 0.01\angstrom$ and all angles in any residue triangle $\triangle N_i A_i C_i$ are at least $3^\circ$, which makes the bases of all residue triangles well-defined in Definition~\ref{dfn:bri} below.

\index{Backbone Rigid Invariant}
\begin{dfn}[backbone rigid invariant $\bri(S)$ of a protein backbone $S$]
\label{dfn:bri}
In the notations of Definition~\ref{dfn:trin}, 
define the orthonormal basis vectors 
$\vb*{u}_i=\dfrac{\ve{A_i N_i}}{|\ve{A_i N_i}|}$, 
$\vb*{v}_i=\dfrac{\vb*{h}_i}{|\vb*{h}_i|}$ for $\vb*{h}_i=\ve{A_i C_i}-b_i\ve{A_i N_i}$, 
$b_i=\dfrac{\ve{A_i C_i}\cdot\ve{A_i N_i}}{|\ve{A_i N_i}|^2}$, and $\vb*{w}_i=\vb*{u}_i\times\vb*{v}_i$.
The \emph{backbone rigid invariant} $\bri(S)$ is the $m\times 9$ matrix whose $i$-th row for $i=2,\dots,m$ contains the coefficients $x,y,z$ of the vectors $\ve{C_{i-1}N_i}$, $\ve{N_i A_i}$, $\ve{A_i C_i}$ in the basis $\vb*{u}_{i-1},\vb*{v}_{i-1},\vb*{w}_{i-1}$.
So, for $i=2,\dots,m$, the nine columns of $\bri(S)$ contain the coordinates $x(C_{i-1} N_i),y(C_{i-1} N_i),z(C_{i-1} N_i)$ of 
$\ve{C_{i-1}N_i}$, 
followed by the three coordinates $x(N_i A_i)$, $y(N_i A_i)$, $z(N_i A_i)$ of $\ve{N_i A_i}$ and three coordinates $x(A_i C_i)$, $y(A_i C_i)$, $z(A_i C_i)$ of $\ve{A_i C_i}$. 
In the exceptional case $i=1$, the first row of $\bri(S)$ has only three non-zero coordinates $x(N_1 A_1)$, $x(A_1 C_1)$ and $y(C_1)=y(A_1 C_1)$ from the first row of the invariant $\trin(S)$ in Definition~\ref{dfn:trin}. 
\edfn
\end{dfn}

For a backbone of $m$ residues, the first row of the $m\times 9$ matrix $\bri(S)$ contains only three non-zero coordinates.
Hence the matrix $\bri(S)$ can be considered a vector of length $9(m-1)+3=9m-6$.
The simplest metric on $\bri$s as vectors in $\R^{9m-6}$ is $L_\infty$ equal to the maximum absolute difference between all corresponding coordinates.
\medskip

A small value $\de$ of $L_\infty(\bri(S),\bri(Q))$ guarantees by Theorem~\ref{thm:BRI_inverse} that backbones $S,Q$ are closely matched by rigid motion.
Another metric, such as Euclidean distance or its normalisation by the chain length, has no such guarantees and can be small even for a few outliers that can affect the rigid structure and hence functional properties of a protein. 
Theorem~\ref{thm:BRI_complete} proves conditions~\ref{pro:proteins}(a,b,c,e,h) in Problem~\ref{pro:proteins}. 
\myskip

All stated results below have references to the original papers with detailed proofs.

\index{Backbone Rigid Invariant}
\begin{thm}[completeness, reconstruction, and subchains {\cite[Theorem~3.5]{anosova2025complete}}]
\label{thm:BRI_complete}
\textbf{(a)}
Under any rigid motion in $\R^3$, the matrix $\trin(S)$ in Definition~\ref{dfn:trin} is invariant, while $\bri(S)$ in Definition~\ref{dfn:bri} is a complete invariant, so any backbones $S,Q\subset\R^3$ are matched by rigid motion if and only if $\bri(S)=\bri(Q)$.
\medskip

\noindent
\textbf{(b)}
For any backbone $S$ of $m$ residues, 
 the invariant $\bri(S)$, metric $L_\infty$ on $\bri$s, and a reconstruction of $S\subset\R^3$ from $\bri(S)$ can be computed in time $O(m)$. 
\medskip

\noindent
\textbf{(c)}
Let $Q$ be a subchain of $j$ consecutive residues in a backbone $S\subset\R^3$.
If $Q$ includes the first residue of $S$, then $\bri(Q)$ consists of the first $j$ rows of $\bri(S)$.
If $Q$ starts from the $i$-th residue of $S$ for $i>1$, the rows $2,\dots,j$ of $\bri(Q)$ coincide with the rows $i+1,\dots,i+j-1$ of $\bri(S)$.
The 1st row of $\bri(Q)$ is computed from the $i$-th row of $\bri(S)$ in a constant time, so $\bri(Q)$ is computed from $\bri(S)$ in time $O(j)$. 
\ethm
\end{thm}

\begin{cor}[completeness under isometry {\cite[Corollary~3.6]{anosova2025complete}}]
\label{cor:isometry}
Any mirror image $\bar S$ of a backbone $S\subset\R^3$ has the invariant $\ov{\bri}(S):=\bri(\bar S)$ obtained by reversing the signs in all $z$-columns of $\bri(S)$.
The unordered pair of $\bri(S)$ and $\ov{\bri}(S)$ is complete under isometry.
\ethm
\end{cor}

Since the realisability in condition~\ref{pro:proteins}(f) did not appear in \cite[Problem~1.2]{anosova2025complete}, new Lemma~\ref{lem:BRI_realisability} describes the geometric realisability of non-degenerate polygonal lines in $\R^3$.
The physical realisability of protein backbones will be tackled in future work.
 
\begin{lem}[realisability of $\bri$]
\label{lem:BRI_realisability}
A sequence of $m$ ordered triplets of points $N_i,A_i,C_i$, $i=1,\dots,m$, is called \emph{non-degenerate} if the vectors $\ve{N_i A_i}$ and $\ve{A_i C_i}$ are not parallel for $i=1,\dots,m-1$. 
The invariant space $I(\bri)$, i.e. the collection of $\bri(S)\in\R^{9m-6}$ for all non-degenerate sequences $S$ of $m$ ordered triplets $N_i,A_i,C_i$, consists of any sequence of
numbers $l>0$, $x$, $y\neq 0$, followed by $m-1$ triples of vectors $\vec a_i,\vec b_i,\vec c_i$, $i=2,\dots,m$, such that $\vec a_i$ and $\vec b_i$ are not parallel for $i=2,\dots,m-1$.
\elem
\end{lem}
\begin{proof}
The first three numbers $l=|N_i A_i|$, $x=x(A_i C_i)$, $y=y(A_i C_i)$ form the triangular invariant $\trin(S)$ from Definition~\ref{dfn:trin}.
The realisability conditions $l>0$ and $y\neq 0$ mean that the vectors $\ve{A_i N_i}$ and $\ve{A_i C_i}$ are not parallel and hence define the orthonormal basis $\vec u_1,\vec v_1, \vec w_1$ associated with the first residue triangle $\triangle N_1 A_1 C_1$.
Similarly, every next pair of vectors $\vec a_i=\ve{N_i A_i}$ and $\vec b_i=\ve{A_i C_i}$ should not be parallel so that we can define an orthonormal basis of the $(i+1)$-st residue for $i=2,\dots,m-1$.
\end{proof}

Theorem~\ref{thm:BRI_continuous} will prove the Lipschitz continuity of $\bri$ in condition~\ref{pro:proteins}(c).
For a given backbone $S$ and its perturbation $Q$, let $l_{N,A}$ and $L_{N,A}$ denote the minimum and maximum bond length between any $\al$-carbon $A_i$ and nitrogen $N_i$ in $S,Q$, respectively.
The maximum bond lengths $L_{A,C},L_{C,N}$ are similarly defined for other types of bonds.

\index{Backbone Rigid Invariant}
\begin{thm}[Lipschitz continuity of $\bri$, {\cite[Theorem~4.1]{anosova2025complete}}]
\label{thm:BRI_continuous}
For any $\ep>0$, let $Q$ be obtained from a backbone $S\subset\R^3$ by perturbing every atom of $S$ up to Euclidean distance $\ep$.
Let $h=\min_i|y(A_i C_i)|$ be the minimum height in triangles $\triangle N_i A_i C_i$ at $C_{i}$ for all residues in the backbones $S,Q$. 
Set $L=\max\{L_{C,N},L_{N,A},L_{A,C}\}$, $K=\dfrac{1}{l_{N,A}}+\dfrac{2}{h}\Big(1+2\dfrac{L_{A,C}}{l_{N,A}}\Big)$, and $\la=2(1+2LK)$.
Then $L_\infty(\bri(S),\bri(Q))\leq \la\ep$. 
\ethm
\end{thm} 

\begin{exa}[continuity of $\bri$]
\label{exa:continuity}
For all 707K+ cleaned chains, the median upper bound for $\la$ is about 34.5, but the real values are smaller as in the example below.
Consider the backbone $S$ of the chain A (141 residues) from the standard hemoglobin 2hhb in the PDB.
We perturb $S$ to $Q$ by adding to each coordinate $x,y,z$ of all atoms in $S$ some uniform noise up to various thresholds $\ep=0.01,0.02,\dots,0.1\angstrom$.
Fig.~\ref{fig:hemoglobins}~(top left) shows how the distance $L_\infty(\bri(S),\bri(Q))$ averaged over 20 perturbations depends on $\ep$
As expected by Theorem~\ref{thm:BRI_continuous}, the metric $L_\infty$ is perturbed linearly up to $\la\ep$ for $\la\approx 4$.
\eexa
\end{exa}

Since the metric $L_\infty$ between invariants $\bri$ ($m\times 9$ matrices) can be computed in linear time $O(m)$, Theorem~\ref{thm:BRI_continuous} also completes condition (\ref{pro:proteins}f) in Problem~\ref{pro:proteins}.
Theorem~\ref{thm:BRI_inverse} will prove 
condition in \ref{pro:proteins}(d).

\index{Backbone Rigid Invariant}
\begin{thm}[inverse continuity of $\bri$, {\cite[Theorem~4.8]{anosova2025complete}}]
\label{thm:BRI_inverse}
For any $\de>0$ and backbones $S,Q\subset\R^3$ with $L_\infty(\bri(S),\bri(Q))<\de$, there is a rigid motion $f$ of $\R^3$ such that any atom of $S$ is $\mu\de$-close to the corresponding atom of $f(Q)$ for $\mu=\sqrt{3}\dfrac{(8LK)^{m-1}-1}{8LK-1}$.
Let $\widehat{\bri}(S)$ be $\bri(S)$ after multiplying the $i$-th row by $\dfrac{(8LK)^{i-1}-1}{8LK-1}$ for $i=2,\dots,m$.
Then $L_\infty(\widehat{\bri}(S),\widehat{\bri}(Q))<\de$ guarantees a rigid motion $f$ of $\R^3$ such that any atom of the backbone $S$ is $\sqrt{3}\de$-close to the corresponding atom of $f(Q)$.
\ethm
\end{thm}

\section{Average invariant, diagrams, and barcodes of backbones}
\label{sec:barcode}

This section simplifies the complete invariant $\bri$ to its average vector in $\R^9$ and introduces the diagram and barcode that visually represent the high-dimensional $\bri$. 

\index{backbone rigid average invariant}
\index{backbone invariant diagram}
\index{backbone invariant barcode}
\begin{dfn}[average invariant $\brain$, diagram $\bid$, and barcode $\bib$]
\label{dfn:diagrams}
\tb{(a)}
For any protein backbone $S$ of $m$ residues, the \emph{backbone rigid average invariant} $\brain(S)\in\R^9$ is the vector of nine column averages in $\bri(S)$ excluding the first row. 
\myskip

\nt
\tb{(b)}
The \emph{backbone invariant diagram} $\bid(S)$ consists of nine polygonal curves going through the points $(i,c(i))$, $i=2,\dots,m$, where $c$ is one of the coordinates (columns) of $\bri(S)$, see Fig.~\ref{fig:hemoglobins}~(middle).
\myskip

\nt
\tb{(c)}
For each atom type such as $N$, the coordinates $(x(N_i),y(N_i),z(N_i))$ are linearly converted into the RGB color value for $i=1,\dots,m$.
The resulting color bars for the ordered atoms $N,A,C$ form the \emph{backbone invariant barcode} $\bib(S)$, see Fig.~\ref{fig:hemoglobins}~(bottom).
\edfn
\end{dfn}

\index{backbone rigid average invariant}
\index{backbone invariant diagram}
\index{backbone invariant barcode}
\begin{exa}[hemoglobins]
\label{exa:hemoglobins}
The PDB contains thousands of hemoglobin structures.
We consider here the structure 2hhb as a standard, and compare it with oxygenated 1hho, which contains an extra oxygen whose transport is facilitated by hemoglobin.
In both cases, we considered the main chains (entity 1, model 1, chain A) of 141 residues. 
\smallskip

\begin{figure}[ht!]
\includegraphics[height=19mm]{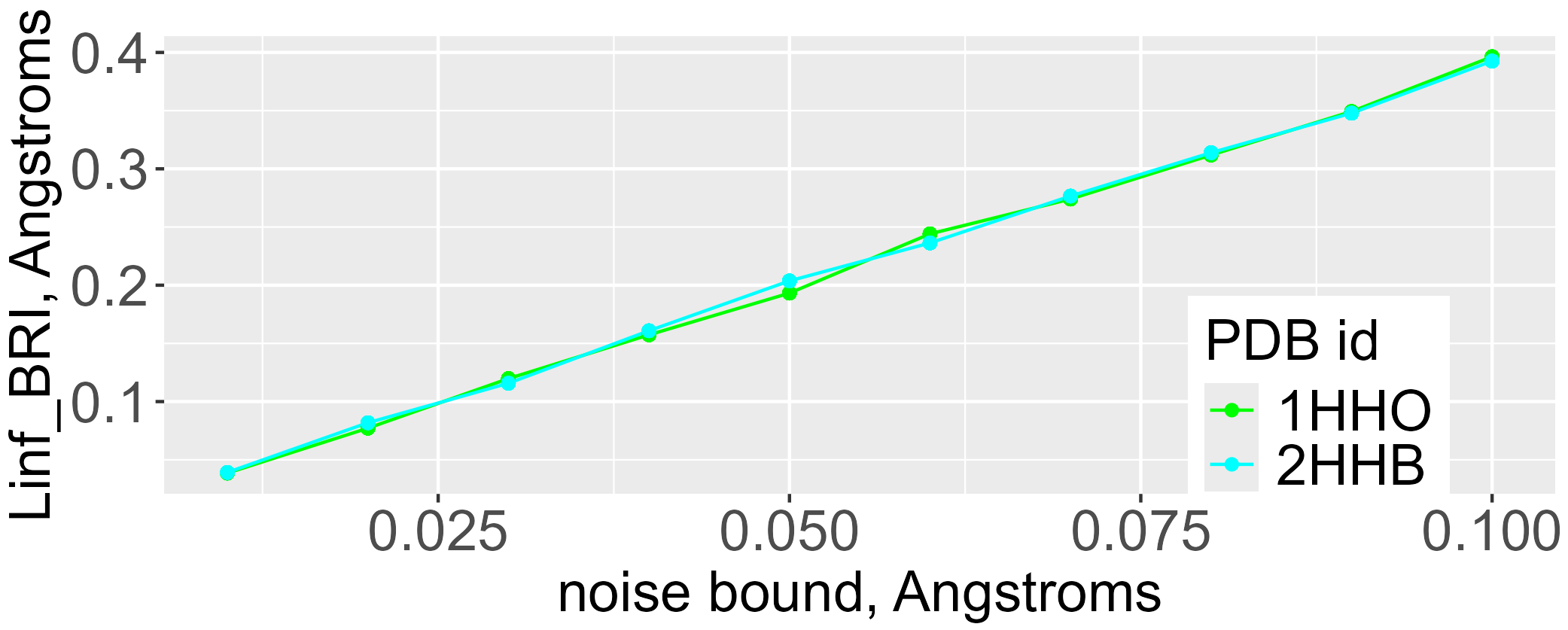}
\includegraphics[height=19mm]{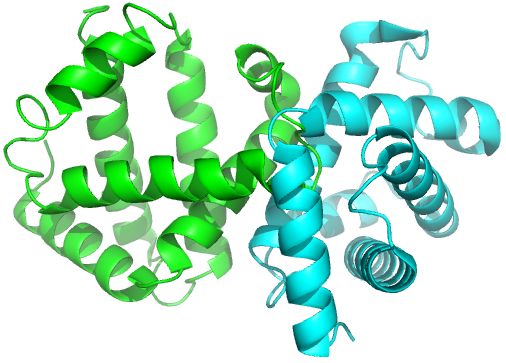}
\includegraphics[height=19mm]{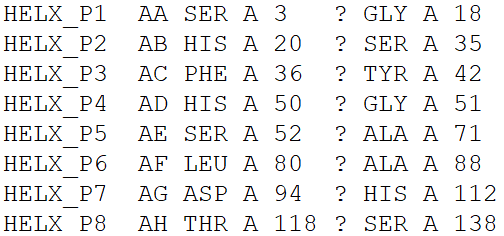}
\medskip

\includegraphics[height=38mm]{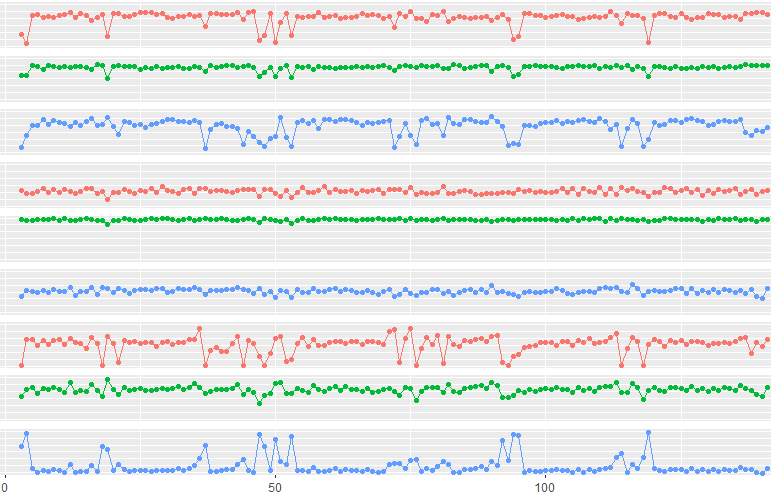}
\includegraphics[height=38mm]{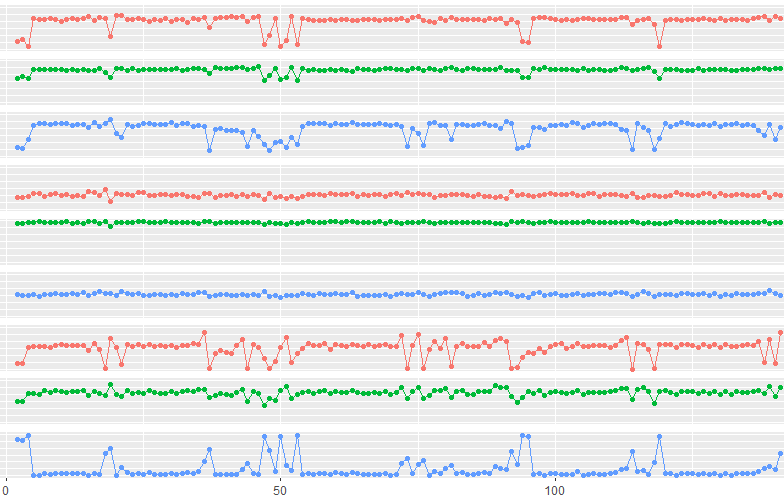}
\medskip

\includegraphics[height=15mm]{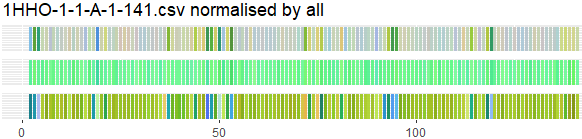}
\includegraphics[height=15mm]{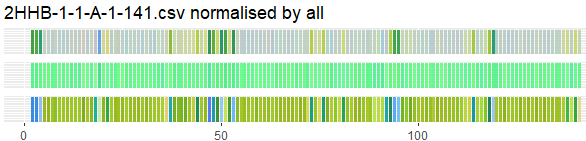}

\caption{
\textbf{Row 1}: the Lipschitz continuity of $\bri$ from Theorem~\ref{thm:BRI_continuous} is illustrated on the left by perturbing hemoglobins in
Example~\ref{exa:continuity}, whose main chains A of 141 residues are shown in the middle (oxygenated 1hho in green, standard 2hhb in cyan) and eight $\al$-helices found by \cite{kabsch1983dictionary} 
on the right.
\textbf{Row 2}: the Backbone Invariant Diagram ($\bid$) of the hemoglobins 1hho vs 2hhb in the PDB, see Definition~\ref{dfn:diagrams}.
\textbf{Row 3}: the Backbone Invariant Barcode ($\bib$), see Example~\ref{exa:hemoglobins}. }
\label{fig:hemoglobins}
\end{figure}

The top left image in Fig.~\ref{fig:hemoglobins}~(top) shows that the Lipschitz constant from Theorem~\ref{thm:BRI_continuous} is $\la\approx 4$ for both hemoglobins.
Fig.~\ref{fig:hemoglobins}~(middle) illustrates the complexity of identifying similar proteins with distant coordinates.
The similarity under rigid motion becomes clear by comparing their diagrams and barcodes in Fig.~\ref{fig:hemoglobins}~(rows 2, 3).
\smallskip

More importantly, a rigidly repeated pattern such as $\al$-helix or $\be$-strand has constant invariants over several residue indices, which are easily detectable in $\bid$ and visible in $\bib$ as intervals of uniform color.
The PDB uses the baseline algorithm DSSP (Define Secondary Structure of Proteins) \cite{kabsch1983dictionary}, which depends on several manual parameters and sometimes outputs $\al$-helices of only two residues.
\myskip

For instance, the PDB entries 1hho and 2hhb in Fig.~\ref{fig:hemoglobins}~(right) include HELX\_P4  consisting of only residues 50 and 51, and HELX\_P5  of length 20 over residue indices $i=52,\dots,71$.
Fig.~\ref{fig:hemoglobins} shows that a `constant' interval of little noise appears only for $i=54,\dots,70$.
Hence new invariants allow a more objective detection of secondary structures,  which will be explored in future work.
\eexa
\end{exa}

\section{A fast detection of duplicate chains in the Protein Data Bank}
\label{sec:duplicates}

The linear time of the complete invariant $\bri(S)$ has enabled all-vs-all comparisons for all tertiary structures in the PDB, which was additionally cleaned by Protocol~\ref{prot:cleanPDB}.
To speed up comparisons, Lemma~\ref{lem:brain} proves that the metric  $L_\infty(\bri(S),\bri(Q))$ between complete invariants is not smaller than 
the faster distance $L_\infty(\brain(S),\brain(Q))$ between the averaged invariants (vectors of 9 coordinates) from Definition~\ref{dfn:diagrams}.

\begin{lem}[metrics on $\bri$ and $\brain$, {\cite[Lemma~6.1]{anosova2025complete}}]
\label{lem:brain}
Any protein backbones $S,Q$ of the same number of residues
satisfy the inequality $L_\infty(\brain(S),\brain(Q))\leq L_\infty(\bri(S),\bri(Q))$.
\elem
\end{lem}

The complete invariants and their statistical summaries were computed in 3 hours 18 min 21 sec.
After comparing all (888+ million) pairs of same-length backbones within 1 hour, we found 13907 pairs $S,Q$ with the \emph{exact zero-distance} $L_\infty(\bri(S),\bri(Q))=0$ between complete invariants meaning that all these backbones $S,Q$ are related by rigid motion, but they may not be geometrically identical.
\medskip

However, 9366 of these pairs turned out to have $x,y,z$ coordinates of all main atoms \emph{identical to the last digit} despite many of them (763) coming from \emph{different PDB entries}. 
Table~\ref{tab:PDB9pairs} lists nine pairs whose geometrically identical chains unexpectedly differ in the sequences of amino acids.
The duplicates from Table~\ref{tab:PDB9pairs} were shown to the PDB validation team, who did not know about the found coincidences (in coordinates) and differences (in amino acids), because the PDB validation is currently done only for an individual protein  (checking atom clashes, outliers etc).

\begin{table}
\caption{Chains with identical backbones but different sequences of amino acid residues.}
\label{tab:PDB9pairs}
\centering
\begin{tabular}{l|l|l|l|l}
PDB id1 & method and & PDB id2  & 
all atoms have  & different\\
\& chain & resolutions, $\angstrom$ & \& chain & 
identical $x,y,z$ & residues \\
\hline
1a0t-B  & X-ray, 2.4, 2.4
 & 1oh2-B 
 & all $3\times 413$ & 9 \\
1ce7-A & X-ray, 2.7, 2.7
 & 2mll-A 
 & all $3\times 241$ & 1, GLY$\neq$HIS  \\
1ruj-A & X-ray, 3, 3
 & 4rhv-A 
 & all $3\times 237$ & 1, GLY$\neq$SER \\
1gli-B/D & X-ray, 2.5, 1.7
 & 3hhb-B/D 
 & all $3\times 146$ & 1, MET$\neq$VAL \\
2hqe-A & X-ray, 2, 2
 & 2o4x-A 
 & all $3\times 217$ & 1, GLN$\neq$GLU \\
5adx-T &  EM, 4, 8.2
 & 5afu-Z 
 & all $3\times 165$ & 1, ILE$\neq$VAL \\
5lj3-O & EM, 3.8, 10
 & 5lj5-P 
 & all $3\times 252$ & 1, ALA$\neq$VAL \\
8fdz-A & X-ray, 2.5, 2.2
 & 8fe0-A 
 & all $3\times 200$ & 1, THR$\neq$SER \\
\end{tabular}
\end{table}

In the row starting with 2hqe in Table~\ref{tab:PDB9pairs}, the chain IDs A, B refer to two pairs of duplicates: chain A of 2hqe is identical to chain A of 2o4x, similarly for B.
The notation $\{$B,D$\}$ in the row starting with 1gli means 4 duplicates: each of the chains B,D in 1gli is identical to each of the chains B,D in 3hhb. 
\myskip

The histograms in \cite[Fig.~5]{anosova2025complete} reveal about 220K pairs of near-duplicates among 707K+ cleaned chains up to $L_\infty\leq 0.01\angstrom$. 
The bound of $0.01\angstrom$ is considered noise because the smallest inter-atomic distance is about 100 times larger at $1\angstrom=10^{-10}$ m.
\medskip

The physical meaning of distances follows from the bi-continuity conditions (c,d) in Problem~\ref{pro:proteins}.
If every atom of a backbone $S$ is shifted up to Euclidean distance $\ep$, then $\bri(S)$ changes up to $\la\ep$ in $L_\infty$.
The Lipschitz constant $\la$ was expressed in Theorem~\ref{thm:BRI_continuous} and estimated as $\la\approx 4$ for the hemoglobin chains in Example~\ref{exa:hemoglobins}.
So any small perturbation of atoms yields a small value of $L_\infty$ in Angstroms. 
\myskip

The inverse Lipschitz continuity in \ref{pro:proteins}(d) implies that a small Chebyshev distance $L_\infty(\bri(S),\bri(Q))=\de$ guarantees that all atoms of the backbones $S,Q$ can be matched (under a suitable rigid motion) up to Euclidean distance $\mu\de$ in Theorem~\ref{thm:BRI_inverse}.
\myskip

One potential explanation of identical coordinates is the molecular replacement method \cite{rossmann1990molecular}, which uses an existing protein structure, often a previous PDB deposit or part thereof, to solve a new structure. 
If the newly calculated electron density map does not allow for further refinement then the coordinates may (reasonably) remain unchanged.
The same coincidences can happen with lower-quality cryo-EM maps in which an existing PDB structure may be placed but where the resolution may not allow for further refinement of atomic coordinates \cite{murshudov2011refmac5,hekkelman2024pdbredo}.
\medskip

We have checked that the found duplicate backbones also have identical distance matrices on $3m$ ordered atoms, which were slower to compute in time $O(m^2)$ over two days on a similar machine.
The widely used DALI server \cite{holm2024dali} also confirmed the found duplicates by the traditional Root Mean Square Deviation (RMSD) through optimal alignment.
The DALI took about 30 min on average to find a short list of nearest neighbors of one chain in the whole PDB.
Extrapolating this time to all pairwise comparisons for 707K+ cleaned chains yields 40+ years, slower by orders of magnitude than 6 hours needed for all comparisons of $\bri$s on the same desktop computer.
\medskip

The ultra-fast speed of all-vs-all comparisons by $\bri$ is explained by the hierarchical nature of this complete invariant.
To find near-duplicates in the PDB, we first compared only average invariants $\brain(S)\in\R^9$.
By Lemma~\ref{lem:brain} the full comparisons by $\bri$ are needed only for a tiny proportion of backbones with the closest vectors $\brain(S)$.
This hierarchical speed-up is unavailable for any distance without underlying invariants.

\bibliographystyle{plain}
\bibliography{Geometric-Data-Science-book}

%
%
%

\chapter{Complete and polynomial-time invariants of unordered points in $\R^n$}
\label{chap:directional} 

\abstract{
This chapter is the first in the book to focus on Euclidean clouds of unordered points under rigid motion in $\R^n$. 
We leverage Principal Component Analysis to construct a direction-based invariant of point clouds, whose continuity and completeness under isometry proved for principally generic clouds.
This invariant is extended to a larger distribution that is complete for all clouds under rigid motion. 
The main novelty are polynomial-time algorithms for these invariants and distance metrics. 
}

\section{Towards complete and polynomial-time invariants for clouds}
\label{sec:Euclidean}

\begin{figure}[h]
\includegraphics[height=19mm]{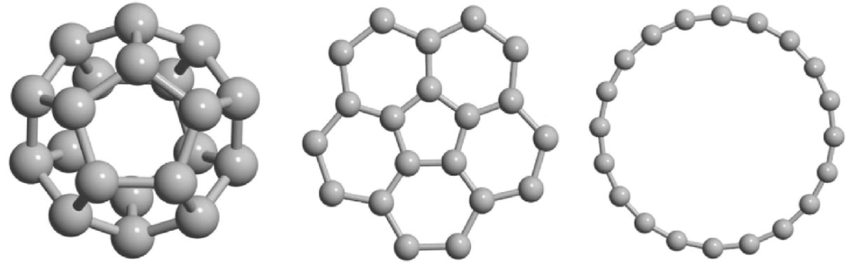}
\hspace*{1mm}
\includegraphics[height=19mm]{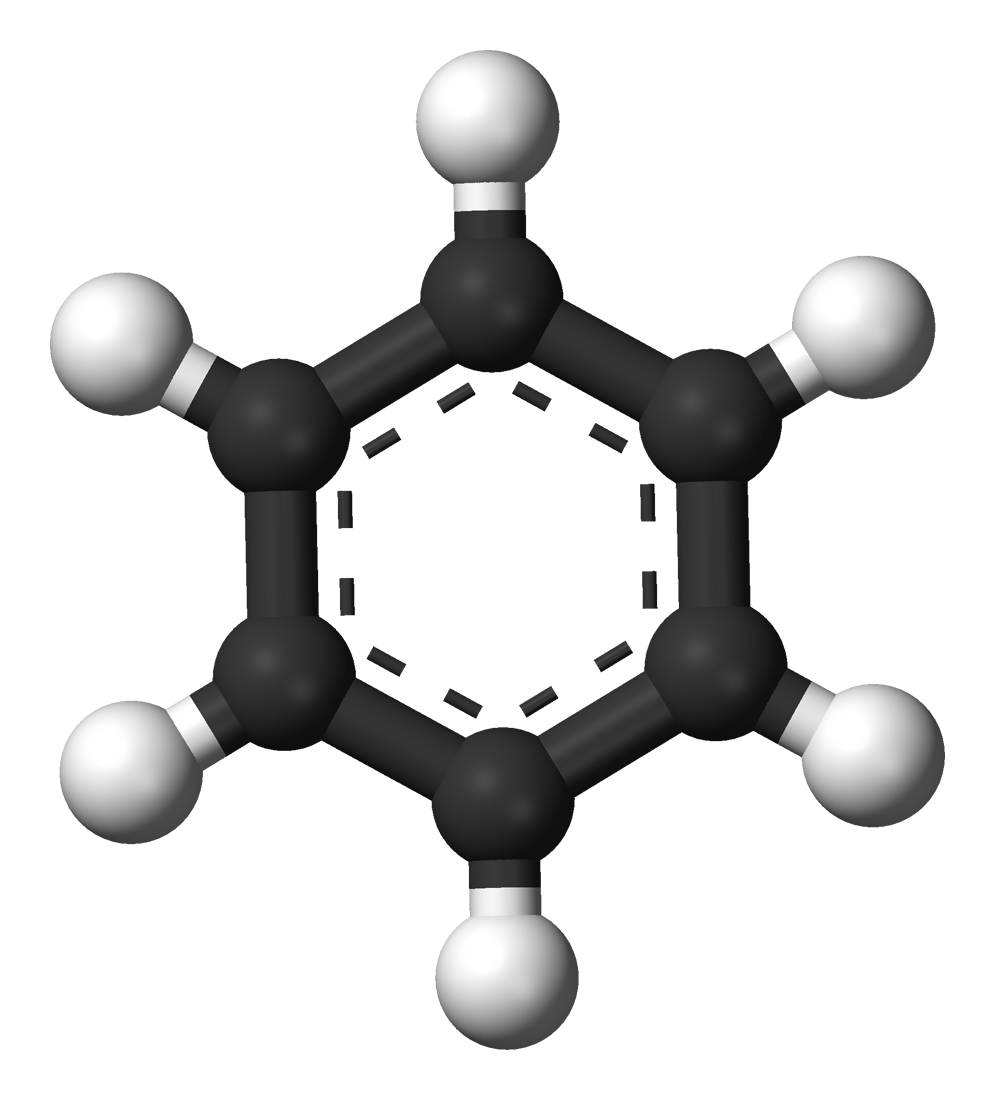}
\hspace*{0.5mm}
\includegraphics[height=19mm]{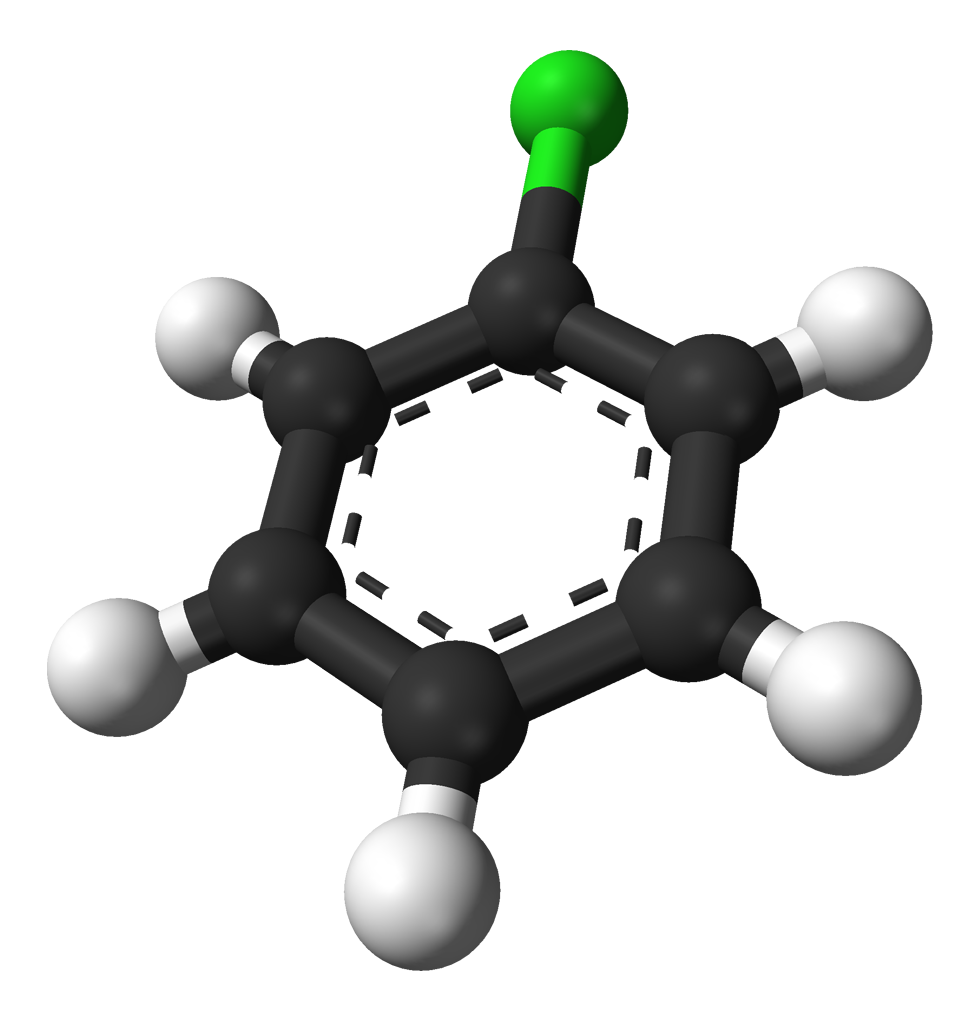}
\hspace*{0.5mm}
\includegraphics[height=19mm]{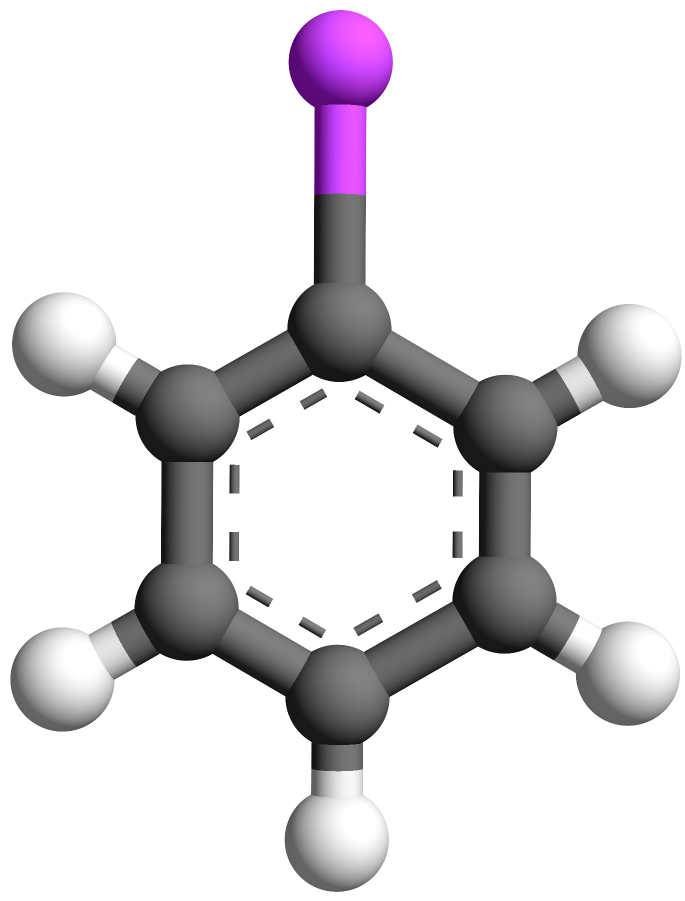}
\caption{
Isomers of \ce{C20}, benzene \ce{C6 H6}, phenyllithium \ce{C6 H5 Li}, chlorobenzene \ce{C6 H_5 Cl} have many indistinguishable atoms.
}
\label{fig:molecules}
\end{figure}

All sections in this chapter follow paper \cite{kurlin2024polynomial} with minor updates.
Any finite chemical system, such as a molecule, can be represented as a cloud of atoms whose nuclei are real physical objects \cite{widdowson2022average}, while chemical bonds are not real sticks and only abstractly represent inter-atomic interactions. 
In the hardest scenario, all atoms are modelled as zero-sized points at all atomic centres without any labels such as chemical elements.
For example, the \ce{C60} molecule \cite{kroto1985c60} consists of 60  unordered carbons.
Allowing different compositions enables a quantitative comparison of isomers, 
see Fig.~\ref{fig:molecules}.
 
This chapter studies finite clouds of unordered points in $\R^n$ for a fixed dimension $n$.
Problem~\ref{pro:directional} adjusts Geo-Mapping Problem~\ref{pro:geocodes} to unordered clouds under rigid motion.
The stronger problem with Lipschitz continuity will be solved in Chapter~6.

\begin{pro}[complete and polynomial-time invariants of clouds in $\R^n$]
\label{pro:directional}
Design an invariant $I$ of all clouds of unordered points in $\R^n$ 
satisfying the conditions below.
\medskip

\noindent
\tb{(a)} 
\emph{Completeness:}
any finie clouds $A,B\subset\R^n$ of unordered points are related by rigid motion ($A\cong B$) if and only if $I(A)=I(B)$.
\medskip

\noindent
\tb{(b)} 
\emph{Reconstruction:}
any cloud $A\subset\R^n$ of unordered points can be reconstructed from its invariant value $I(A)$, uniquely under rigid motion in $\R^n$. 
\medskip

\noindent
\tb{(c)} \emph{Metric:} 
there is a distance $d$ on the space $\{I(A) \vl \text{unordered clouds }A\subset\R^n\}$ satisfying all metric axioms in Definition~\ref{dfn:metrics}(a). 
\medskip

\noindent
\tb{(d)} 
\emph{Computability:} 
for a fixed dimension $n$, the invariant $I(A)$, a reconstruction of $A\subset\R^n$ from $I(A)$, and the metric $d(I(A),I(B))$ are computable in times that depend polynomially on the maximum size $\max\{|A|,|B|\}$ of any clouds $A,B\subset\R^n$.
\epro
\end{pro}

Based on Principal Component Analysis, section~\ref{sec:PCI} introduces the Principal Coordinates Invariant ($\PCI$) to uniquely identify under isometry in $\R^n$ all point clouds that allow a unique alignment by principal directions.
Section~\ref{sec:SM} defines a symmetrised metric on $\PCI$s, which is continuous under perturbations in general position and can be computed (for a fixed dimension $n$) in a subquadratic time in the number of unordered points.
Section~\ref{sec:WMI} extends the $\PCI$ to the Weighted Matrices Invariant ($\WMI$), which is complete for all point clouds under isometry in $\R^n$.
Section~\ref{sec:LAC+EMD} applies the Linear Assignment Cost and Earth Mover's Distance to define metrics on $\WMI$s.
\myskip

For a fixed dimension $n$ of the ambient space $\R^n$, all these invariants and metrics have polynomial-time algorithms in the number $m$ of the given points.
For $n=2$, the time $O(m^{3.5}\log m)$ improves the time $O(m^5\log m)$ of the only previous exact algorithm \cite{chew1997geometric} for the Hausdorff distance on isometry classes of clouds.
\myskip

\begin{figure}[h]
\includegraphics[height=18mm]{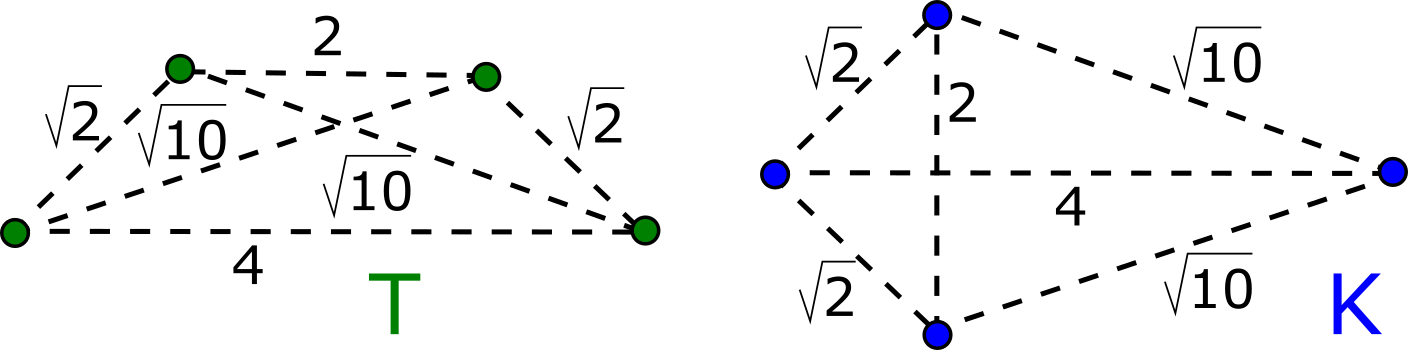}
\hspace*{1mm}
\includegraphics[height=18mm]{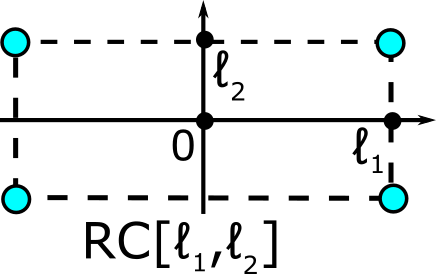}
\hspace*{1mm}
\includegraphics[height=18mm]{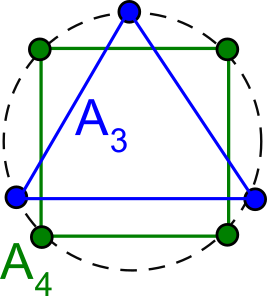}
\caption{
\textbf{First and second}: non-isometric sets $T\not\cong K$ of 4 points have the same 6 pairwise distances.
\textbf{Third}: the vertex set $\RC[l_1,l_2]$ of a $2l_1\times 2l_2$ 
 rectangle. 
\textbf{Fourth}: what is the distance between an equilateral 
triangle $A_3$ and a square $A_4$?
See new invariants and metrics in Examples~\ref{exa:PCI},~\ref{exa:SM},~\ref{exa:triangle-vs-square}.}
\label{fig:4-point_sets}
\end{figure}

As a potential extension of the side-side-side theorem to $m$ unordered points in $\R^n$, the seminal work \cite{boutin2004reconstructing} in 2004 proved that the total distribution of pairwise distances is a complete invariant under isometry in $\R^n$ for generic  clouds whose point coordinates are not solutions of a complicated polynomial equation.
However, infinitely many counter-examples to the full completeness of this invariant were constructed even for $m=4$ points in $\R^2$ \cite{caelli1979generating}.
The first two pictures of Fig.~\ref{fig:4-point_sets} show the simplest non-isometric clouds $T\not\simeq K$ of 4 points in $\R^2$.
Other past work was reviewed in \cite[section~2]{kurlin2024polynomial}.

\section{The Principal Coordinates Invariant of unordered clouds in $\R^n$}
\label{sec:PCI}

This section recalls Principal Component Analysis (PCA), whose principal directions \cite{abdi2010principal} will be used to introduce the Principal Coordinates Invariant (PCI) in Definition~\ref{dfn:PCI}.
We assume that all coordinates in $\R^n$ have the same units.
In practice, we should first normalise all features given in different units.
\myskip

Any cloud $A\subset\R^n$ of $m$ unordered points has the \emph{centre of mass} $O(A)=\dfrac{1}{m}\sum\limits_{p\in A} p$.
Shifting $A$ by the vector $-O(A)$ allows us to always assume that $O(A)$ is the origin $0$. 
Then Problem~\ref{pro:directional} reduces to invariants only under orthogonal maps from the orthogonal group $O(\R^n)$ instead of the full Euclidean group $\Eu(\R^n)$.

\index{covariance matrix}
\begin{dfn}[\emph{covariance} matrix $\cov(A)$ of a point cloud $A$]
\label{dfn:covariance}
If we arbitrarily order points $p_1,\dots,p_m$ of a cloud $A\subset\R^n$, we get the sample $n\times m$ matrix (or data table) $P(A)$, whose $i$-th column consists of $n$ coordinates of the point $p_i\in A$, $i=1,\dots,m$.
The \emph{covariance} $n\times n$ matrix $\cov(A)=\dfrac{P(A) P(A)^T}{n-1}$ is symmetric and positive semi-definite meaning that $v^T\cov(A)v\geq 0$ for any vector $v\in\R^n$.
Hence the matrix $\cov(A)$ has real \emph{eigenvalues} $\la_1\geq\dots\geq\la_n\geq 0$ satisfying $\cov(A)\vec v_j=\la_j \vec v_j$ for an \emph{eigenvector} $\vec v_j\in\R^n$, which can be scaled by any real $s\neq 0$.
\edfn
\end{dfn}

If all eigenvalues of the covariance matrix $\cov(A)$ are distinct and positive, there is an orthonormal basis of eigenvectors $\vec v_1,\dots,\vec v_n$ ordered according to the decreasing eigenvalues $\la_1>\cdots>\la_n>0$.
This \emph{eigenbasis} is unique under reflection $\vec v_j\lra -\vec v_j$ of each eigenvector, $j=1,\dots,n$.

\index{principally generic cloud}
\begin{dfn}[\emph{principally generic} cloud]
\label{dfn:principally_generic_cloud}
A point cloud $A\subset\R^n$ is \emph{principally generic} if, after shifting 
$O(A)$ to the origin, the covariance matrix $\cov(A)$ has distinct eigenvalues $\la_1>\cdots>\la_n>0$.
The $j$-th eigenvalue $\la_j$ defines the $j$-th \emph{principal direction} parallel to an eigenvector $\vec v_j$, which is uniquely determined under scaling.  
\edfn
\end{dfn}

The vertex set of any rectangle in $\R^2$, but not a square, is principally generic.

\index{Principal Coordinates Invariant}
\index{Principal Coordinates Matrix} 
\begin{dfn}[matrix $\PCM$ and invariant $\PCI$]
\label{dfn:PCI}
For $n\geq 1$, let $A\subset\R^n$ be a principally generic cloud of points $p_1,\dots,p_m$ with the centre of mass $O(A)$ at $0\in\R^n$.
Then $A$ has principal directions along unit eigenvectors $\vec v_1,\dots,\vec v_n$, defined up to a sign.
In the orthonormal basis $V=(\vec v_1,\dots,\vec v_n)^T$, any point $p_i\in A$ has the \emph{principal coordinates} $p_i\cdot \vec v_1,\dots,p_i\cdot \vec v_n$, which can be written as a vertical column $n\times 1$ denoted by $Vp_i$.
The \emph{Principal Coordinates Matrix} is the $n\times m$ matrix $\PCM(A)$ whose $m$ columns are the \emph{coordinate sequences} $Vp_1,\dots,Vp_m$.
Two such matrices are \emph{equivalent} under changing signs of rows due to the ambiguity $\vec v_j\lra -\vec v_j$ of unit length eigenvectors in the basis $V$.
The \emph{Principal Coordinates Invariant} $\PCI(A)$ is an equivalence class of $\PCM(A)$. 
\edfn
\end{dfn}

For simplicity, we skip the dependence on a basis $V$ in the notation $\PCM(A)$. 
The columns of $\PCM(A)$ are unordered, though we can write them according to any order of points in the cloud $A$ considered as the vector $(p_1,\dots,p_m)$.
Then $\PCM(A)$ can be viewed as the matrix product $VA$ consisting of the $m$ columns $Vp_1,\dots,Vp_m$.
One can minimise the ambiguity under re-ordering of columns and switching signs $\vec v_j\lra -\vec v_j$ as follows.
For each $j=1,\dots,n$, choose a sign so that a coordinate with a largest value $|p_i\cdot \vec v_j|$ is positive.
Then write all columns in the lexicographically decreasing order:  $(c_{1j},\dots,c_{nj})>(c_{1k},\dots,c_{nk})$ if a few first values (possibly none) coincide $c_{ij}=c_{ik}$ and then $c_{ij}>c_{ik}$ for the next index $i$.

\begin{exa}[computing PCI]
\label{exa:PCI}
\textbf{(a)}
For any $l_1>l_2>0$, let the \emph{rectangular cloud} $\RC[l_1,l_2]$ consist of the four vertices $(\pm l_1,\pm l_2)$ of the rectangle $[-l_1,l_1]\times[-l_2,l_2]$.
Then $\RC[l_1,l_2]$ has the centre at $0\in\R^2$ and the sample $2\times 4$ matrix $P=\left( \begin{array}{cccc} 
l_1 & l_1 & -l_1 & -l_1\\
l_2 & -l_2 & l_2 & -l_2
\end{array} \right)$ whose columns are in a 1-1 correspondence with (arbitrarily) ordered points $(l_1,l_2)$, $(l_1,-l_2)$, $(-l_1,l_2)$, $(-l_1,-l_2)$.
The covariance matrix 
$\cov(\RC[l_1,l_2])=\mat{4l_1^2}{0}{0}{4l_2^2}$ has eigenvalues $\la_1=4l_1^2>\la_2=4l_2^2$. 
If we choose unit length eigenvectors $\vec v_1=(1,0)$ and $\vec v_2=(0,1)$, then $\PCM(\RC[l_1,l_2])$ coincides with the matrix $P$ above.
The invariant $\PCI(\RC[l_1,l_2])$ is the equivalence class of all matrices obtained from $P$ by changing signs of rows and re-ordering columns.
\medskip

\noindent
\textbf{(b)}
The vertex set $T$ of the trapezium in the first picture of Fig.~\ref{fig:4-point_sets} has four points written in the columns of
the sample matrix $P(T)=\left( \begin{array}{cccc} 
2 & 1 & -1 & -2 \\
-1/2 & 1/2 & 1/2 & -1/2
\end{array} \right)$ so that the centre of mass $O(T)$ is the origin $0$.
Then $\cov(T)=\mat{10}{0}{0}{1}$ has eigenvalues 10, 1 with orthonormal eigenvectors $(1,0)$, $(0,1)$, respectively.
The invariant 
$\PCI(T)$ is the equivalence class of the matrix $P(T)$ above.
The vertex set $K$ of the kite in the second picture of Fig.~\ref{fig:4-point_sets} consists of four points written in the columns of
the sample matrix $P(K)=\left( \begin{array}{cccc} 
5/2 & -1/2 & -1/2 & -3/2 \\
0 & 1 & -1 & 0
\end{array} \right)$ so that the centre of mass $O(K)$ is the origin $0$.
Then $\cov(K)=\mat{9}{0}{0}{2}$ has eigenvalues 9, 2 with orthonormal eigenvectors $(1,0),(0,1)$, respectively.
The invariant 
$\PCI(K)$ is the equivalence class of the matrix $P(K)$.
\eexa
\end{exa}

All results in this chapter have details proofs in the original paper \cite{kurlin2024polynomial}.

\index{Principal Coordinates Invariant}
\begin{thm}[generic completeness of $\PCI$, {\cite[Theorem~3.5]{kurlin2024polynomial}}]
\label{thm:PCI_completeness}
Any principally generic clouds $A,B\subset\R^n$ of $m$ unordered points 
are isometric if and only if their PCI invariants 
coincide as equivalence classes of matrices: $\PCI(A)=\PCI(B)$.
\ethm
\end{thm}

\begin{lem}[time complexity of $\PCI$, {\cite[Lemma~3.6]{kurlin2024polynomial}}]
\label{lem:PCI_comp}
For a principally generic cloud $A\subset\R^n$ of $m$ points, a matrix $\PCM(A)$ from the invariant $\PCI(A)$ in Definition~\ref{dfn:PCI} can be computed in time  $O(n^2m+n^3)$.
\elem
\end{lem}

Theorem~\ref{thm:PCI_completeness} requires that clouds $A,B$ are principally generic, which holds with 100\% probability due to noise.
If real clouds are close to symmetric configurations with equal eigenvalues, to avoid numerical instability, we should use the slower but always complete invariants from section~\ref{sec:WMI}.

\section{A symmetrised metric on principally generic clouds in $\R^n$}
\label{sec:SM}

This section defines a metric on $\PCI$ invariants, whose polynomial-time computation and continuity will be proved in Theorems~\ref{thm:SM_complexity} and~\ref{thm:PCI_continuity}.
For any $v=(x_1,\dots,x_n)\in \R^n$, the \emph{maximum norm} is $||v||_{\infty}=\max\limits_{i=1,\dots,n}|x_i|$.
Below we use the \emph{Chebyshev} distance $L_\infty(u,v)=||\vec u-\vec v||_{\infty}$ between points $u,v\in\R^n$ and
the bottleneck distance $\BD$ from Example~\ref{dfn:metrics}(b) on matrices $P$ interpreted as clouds $[P]$ of column-vectors in $\R^n$. 

\begin{dfn}[$m$-point cloud ${[P]}\subset\R^n$ of an $n\times m$ matrix $P$]
\label{dfn:matrix-cloud}
For any $n\times m$ matrix $P$, let $[P]$ denote the unordered set of its $m$ columns considered as vectors in $\R^n$.
The set $[P]$ of $m$ columns 
can be interpreted as a cloud of $m$ unordered points in $\R^n$.
\edfn
\end{dfn}

For any $n\times m$ matrices $P,Q$, let $g:[P]\to[Q]$ be a bijection of columns
indexed by $1,2,\dots,m$.
Then the Chebyshev distance $L_{\infty}(v,g(v))$ between columns $v\in[P]$ and $g(v)\in[Q]$ is the maximum absolute difference of corresponding coordinates in $\R^n$.
The minimisation over all column bijections $g:[P]\to[Q]$ gives the bottleneck distance $\BD([P],[Q])=\min\limits_{g:[P]\to[Q]}\max\limits_{v\in[P]} L_{\infty}(v,g(v))$ between the sets $[P]$, $[Q]$ considered as clouds of unordered points in $\R^n$. 
\medskip

An algorithm for detecting a potential isometry $A\cong B$ will check if $\SM(A,B)=0$ for the metric $\SM$ defined via changes of signs.
A change of signs in $n$ rows can be represented by a binary string $\si$ in the product group $\Z_2^n$, where $\Z_2=\{\pm 1\}$, 1 means no change, $-1$ means a change.
\medskip

For instance, the binary string $\si=(1,-1)\in\Z_2^2$ acts on the matrix $P=\PCM(\RC[l_1,l_2])$ from Example~\ref{exa:PCI} as follows:
$$\si\left( \begin{array}{cccc} 
l_1 & l_1 & -l_1 & -l_1\\
l_2 & -l_2 & l_2 & -l_2
\end{array} \right)
=\left( \begin{array}{cccc} 
l_1 & l_1 & -l_1 & -l_1\\
-l_2 & l_2 & -l_2 & l_2
\end{array} \right).$$

\index{symmetrised metric}
\begin{dfn}[symmetrised metric $\SM$ on matrices and clouds]
\label{dfn:SM}
For any $n\times m$ matrices $P,Q$, the minimisation for $2^n$ changes of signs represented by strings $\si\in\Z_2^n$ acting on rows gives the \emph{symmetrised metric} $\SM([P],[Q])
=\min\limits_{\si\in\Z_2^n} \BD([\si(P)],[Q])$.
For any principally generic clouds $A,B\subset\R^n$, the \emph{symmetrised metric} is
$\SM(A,B)=\SM([\PCM(A)],[\PCM(B)])$ for matrices $\PCM(A),\PCM(B)$ in Definition~\ref{dfn:PCI}.
\edfn
\end{dfn}

If we denote the action of a column permutation $g$ on a matrix $P$ as $g(P)$, the matrix difference $g(P)-Q$ has the maximum norm $\max\limits_{v\in[P]} L_{\infty}(v,g(v))$. 
Then $\BD([P],[Q])$ will be computed by an efficient algorithm for bottleneck matching in Theorem~\ref{thm:SM_complexity}.

\begin{lem}[metric axioms for the symmetrised metric $\SM$, {\cite[Lemma~4.4]{kurlin2024polynomial}}]
\label{lem:SM_axioms}
\textbf{(a)}
The metric $\SM(P,Q)$ from Definition~\ref{dfn:SM} is well-defined on
 equivalence classes of $n\times m$ matrices $P,Q$ considered under changes of signs of rows and permutations of columns, and satisfies all metric axioms.
\medskip

\noindent
\textbf{(b)}
The metric $\SM(A,B)$ from Definition~\ref{dfn:SM} is well-defined on isometry classes of principally generic clouds $A,B$ and satisfies all axioms.
\elem
\end{lem}

\index{symmetrised metric}
\begin{exa}[computing the symmetrised metric $\SM$]
\label{exa:SM}
\textbf{(a)}
By Example~\ref{exa:PCI}(a), the vertex set $\RC[l_1,l_2]$ of any rectangle with sides $2l_1>2l_2$ in the plane has $\PCI$ represented by the matrix 
$\PCM(\RC[l_1,l_2])=\left( \begin{array}{cccc} 
l_1 & l_1 & -l_1 & -l_1\\
l_2 & -l_2 & l_2 & -l_2
\end{array} \right)$.
The vertex set $\RC[l_1',l_2']$ of any other rectangle has a similar matrix whose element-wise subtraction from $\PCM(\RC[l_1,l_2])$ consists of $\pm l_1\pm l_1'$ and $\pm l_2\pm l_2'$.
Re-ordering columns and changing signs of rows minimises the maximum absolute value of these elements to $\max\{|l_1-l_1'|,|l_2-l_2'|\}$, which should equal $\SM(\RC[l_1,l_2],\RC[l'_1,l'_2])$.
\medskip

\noindent
\textbf{(b)}
The invariants $\PCI$ of the vertex sets $T$ and $K$ in Fig.~\ref{fig:4-point_sets} were computed in Example~\ref{exa:PCI}(b) and represented by these matrices in Definition~\ref{dfn:PCI}:
$$\PCM(T) =
\left( \begin{array}{cccc} 
2 & 1 & -1 & -2 \\
-1/2 & 1/2 & 1/2 & -1/2
\end{array} \right), \qquad
\PCM(K) =
\left( \begin{array}{cccc} 
5/2 & -1/2 & -1/2 & -3/2 \\
0 & 1 & -1 & 0
\end{array} \right).$$

The maximum absolute value of the element-wise difference of these matrices is $|1-(-\frac{1}{2})|=\frac{3}{2}$, which cannot be smaller after permuting columns and changing signs of rows.
The symmetrised metric equals $\SM(T,K)=\BD(\PCM(T),\PCM(K))=\frac{3}{2}$.
\eexa
\end{exa}

\index{symmetrised metric}
\begin{thm}[time of the metric $\SM$, {\cite[Theorem~4.6]{kurlin2024polynomial}}]
\label{thm:SM_complexity}
\textbf{(a)}
Given any $n\times m$ matrices $P,Q$, the symmetrised metric $\SM(P,Q)$ in Definition~\ref{dfn:SM} is computable in time $O(m^{1.5}2^n\log^{n}m)$.
If $n=2$, the time is $O(m^{1.5}\log m)$.
\medskip

\noindent
\textbf{(b)}
The above conclusions hold for $\SM(A,B)$ 
of any principally generic $m$-point clouds $A,B\subset\R^n$ represented by $n\times m$ matrices $\PCM(A),\PCM(B)$. 
\ethm 
\end{thm}

Theorem~\ref{thm:PCI_continuity} proves the continuity in the weaker sense of $\ep,\de$ because $\PCI$ is defined only for generic clouds anyway.
Explicit estimates in the proof from \cite[section~4]{kurlin2024polynomial} are based on recent bounds for perturbations of eigenvectors from  of \cite[Theorem~3]{fan2018eigenvector}.

\index{Principal Coordinates Invariant}
\begin{thm}[continuity of $\SM$, {\cite[Theorem~4.9]{kurlin2024polynomial}}]
\label{thm:PCI_continuity}
For any principally generic cloud $A\subset\R^n$ 
and any $\ep>0$, there is $\de>0$ (depending on $A$ and $\ep$) such that
if any principally generic cloud $B\subset\R^n$ has $\BD(A,B)<\de$, then $\SM(A,B)<\ep$.
\ethm
\end{thm}

\section{A complete invariant for all clouds of unordered points in $\R^n$}
\label{sec:WMI}

This section extends the $\PCI$ from Definition~\ref{dfn:PCI} to a complete $\WMI$ (Weighted Matrices Invariant) of all possible clouds, which may not be principally generic.
\smallskip

If a cloud $A\subset\R^n$ is not principally generic, some of the eigenvalues $\la_1\geq\dots\geq\la_n\geq 0$ of the covariance matrix $\cov(A)$ coincide or vanish.
Let us start with the most singular case when all eigenvalues are equal to $\la>0$.
The case $\la=0$ means that $A$ is a single point. 
Though $A$ has no preferred (principal) directions, $A$ still has the well-defined centre of mass $O(A)=\dfrac{1}{m}\sum\limits_{p\in A}p$, which is at the origin $0\in\R^n$ as always.
For $n=2$, we consider $m$ possible vectors from the origin $0$ to every point of $A\setminus\{0\}$.

\index{Weighted Matrices Invariant}
\begin{dfn}[Weighted Matrices Invariant $\WMI(A)$ for clouds  $A\subset\R^2$]
\label{dfn:WMI2}
Let a cloud $A$ of $m$ points $p_1,\dots,p_m$ in $\R^2$ have the centre of mass at the origin $0$.
For any point $p_i\in A\setminus\{0\}$, let $\vec v_1$ be the unit length vector parallel to $p_i\neq 0$.
Let $\vec v_2$ be the unit length vector orthogonal to $\vec v_1$ whose anti-clockwise angle from $\vec v_1$ to $\vec v_2$ is $+\dfrac{\pi}{2}$.
The $2\times m$ matrix $M(p_i)$ consists of the $m$ pairs of coordinates of all points $p\in A$ written in the orthonormal basis $\vec v_1,\vec v_2$, for example, $\vec p_i=\vect{|\vec p_i|}{0}$. 
Each matrix $M(p_i)$ is considered under re-ordering of columns.
If one point $p$ of $A$ is the origin $0$, there is no basis defined by $p=0$, let $M(p)$ be the zero matrix in this centred case.
If $k>1$ of the matrices $M(p_i)$ are \emph{equivalent} under re-ordering of columns, we collapse them into one matrix with the weight $\dfrac{k}{m}$.
The unordered collection of the equivalence classes of $M(p)$
with weights for all $p\in A$ is called the \emph{Weighted Matrices Invariant} $\WMI(A)$.
\edfn
\end{dfn}

In comparison with the generic case in Definition~\ref{dfn:PCI}, for any fixed $i=1,\dots,m$, if $p_i\neq 0$, then the orthonormal basis $\vec v_1,\vec v_2$ is uniquely defined without the ambiguity of signs, which will re-emerge for higher dimensions $n>2$ in Definition~\ref{dfn:WMI} later.
The vertex sets of regular polygons $A_m$ have $\WMI$ consisting of a single matrix due to extra symmetries as shown below.

\begin{exa}[regular clouds $A_m\subset\R^2$]
\label{exa:reg_clouds}
Let $A_m$ be the vertex set of a regular $m$-sided polygon inscribed into a circle of a radius $r$, see the last picture in Fig.~\ref{fig:4-point_sets}.
Due to the $m$-fold rotational symmetry of $A_m$, the invariant $\WMI(A_m)$ consists of a single matrix (with weight 1) whose columns
 are the vectors $\vect{r\cos\frac{2\pi i}{m}}{r\sin\frac{2\pi i}{m}}$, $i=1,\dots,m$.
For instance, the vertex set $A_3$ of the equilateral triangle has $\WMI(A_3)=
\left\{\left( \begin{array}{ccc} 
r & -r/2 & -r/2 \\
0 & r\sqrt{3}/2 & -r\sqrt{3}/2
\end{array} \right)\right\}$.
The vertex set $A_4$ of the square has $\WMI(A_4)=
\left\{\left( \begin{array}{cccc} 
r & 0 & 0 & -r \\
0 & r & -r & 0
\end{array} \right)\right\}$.
\myskip

Let  $B_m$ be obtained from $A_m$ by adding the origin $0\in\R^2$.
Then $\WMI(B_m)$ has the matrix from $\WMI(A_m)$ with the weight $\dfrac{m}{m+1}$ and the zero $2\times 4$ matrix with the weight $\dfrac{1}{m+1}$ representing the added origin $0$.
\eexa
\end{exa}

Definition~\ref{dfn:WMI} applies to all point clouds $A\subset\R^n$ including the most singular case when all eigenvalues of the covariance matrix $\cov(A)$ are equal, so we have no preferred directions at all.

\index{Weighted Matrices Invariant}
\begin{dfn}[Weighted Matrices Invariant $\WMI$ for any cloud $A\subset\R^n$]
\label{dfn:WMI}
Let a cloud $A\subset\R^n$ of $m$ points $p_1,\dots,p_m$ have the centre of mass at the origin $0$.
For any ordered sequence of points $p_1,\dots,p_{n-1}\in A$, build an orthonormal basis $\vec v_1,\dots,\vec v_n$ as follows.
The first unit length vector $\vec v_1$ is $p_1$ normalised by its length.
For $j=2,\dots,n-1$, the unit length vector 
$\vec v_j$ is $p_j-\sum\limits_{k=1}^{j-1}(p_j\cdot \vec v_k)\vec v_k$ normalised by its length.
\sskip

Then every $\vec v_j$ is orthogonal to all previous vectors $\vec v_1,\dots,\vec v_{j-1}$ and belongs to the $j$-dimensional subspace spanned by $p_1,\dots,p_j$.
Define the last unit length vector $\vec v_n$ by its orthogonality to $\vec v_1,\dots,\vec v_{n-1}$ and the positive sign of the determinant $\det(\vec v_1,\dots,\vec v_n)$ of the matrix with the columns $\vec v_1,\dots,\vec v_n$.
\smallskip

The $n\times m$ matrix $M(p_1,\dots,p_{n-1})$ consists of column vectors of all points $p\in A$ in the basis $\vec v_1,\dots,\vec v_n$, for example, $p_1=(||p_1||_2,0,\dots,0)^T$. 
If $p_1,\dots,p_{n-1}\in A$ are affinely dependent, let $M(p_1,\dots,p_{n-1})$ be the $n\times m$ matrix of zeros in this centred case.
If $k>1$ matrices are \emph{equivalent} under re-ordering of columns, we collapse them into a single matrix with the weight $\dfrac{k}{N}$, where $N=m(m-1)\dots(m-n+1)$.
\sskip

The \emph{Weighted Matrices Invariant} $\WMI(A)$ is the unordered set of equivalence classes of matrices $M(p_1,\dots,p_{n-1})$ with weights for all sequences of $p_1,\dots,p_{n-1}\in A$.
\edfn
\end{dfn}

If $\cov(A)$ has some equal eigenvalues, $\WMI(A)$ can be made smaller by choosing bases only for subspaces of eigenvectors with the same eigenvalue.

\index{Weighted Matrix Invariant}
\begin{thm}[completeness of $\WMI$ under rigid motion in $\R^n$, {\cite[Theorem~5.4]{kurlin2024polynomial}}]
\label{thm:WMI_completeness}
\textbf{(a)}
Any clouds $A,B\subset\R^n$ are related by rigid motion (orientation-preserving isometry) if and only if 
there is a bijection $\WMI(A)\to\WMI(B)$ preserving all weights or, equivalently,
some matrices $P\in\WMI(A)$, $Q\in\WMI(B)$ are related by re-ordering of columns.
So $\WMI(A)$ is a complete invariant of $A$ under rigid motion.
\medskip

\noindent
\textbf{(b)}
Any mirror reflection $f:A\to B$ induces a bijection $\WMI(A)\to\WMI(B)$ respecting their weights and changing the sign of the last row of every matrix.
This pair of $\WMI$s is a complete invariant of $A$ under isometry including reflections.
\ethm
\end{thm}

It suffices to store in computer memory only one matrix $M(p_1,\dots,p_{n-1})$ from the full $\WMI(A)$.
Any such matrix suffices to reconstruct a point cloud $A$, uniquely under rigid motion in $\R^n$ by Theorem~\ref{thm:PCI_completeness}, as required in Problem~\ref{pro:directional}(b).
The full invariant $\WMI(A)$ can be computed from the reconstructed cloud $A\subset\R^n$. 

\begin{lem}[time of $\WMI$, {\cite[Lemma~5.5]{kurlin2024polynomial}}]
\label{lem:WMI_complexity}
For any cloud $A\subset\R^n$ of $m$ points and a fixed sequence of points $p_1,\dots,p_{n-1}\in A$, the matrix $M(p_1,\dots,p_{n-1})$ from Definition~\ref{dfn:WMI} can be computed in time $O(nm+n^3)$.
All $N=m(m-1)\dots(m-n+1)=O(m^{n-1})$ matrices in the Weighted Matrices Invariant $\WMI(A)$ can be computed in time $O((nm+n^3)N)=O(nm^n+n^3m^{n-1})$.  
\elem
\end{lem}

\section{Polynomial-time metrics for all clouds of unordered points in $\R^n$}
\label{sec:LAC+EMD}

This section introduces two metrics on Weighted Matrices Invariants ($\WMI$s), which are computable in polynomial time by Theorems~\ref{thm:LAC_time} and~\ref{thm:EMD_time}.
Since any rigid motion $f:A\to B$ induces a bijection $\WMI(A)\to\WMI(B)$, we will use a linear assignment cost \cite{jonker1987shortest} based on permutations of matrices.

\index{Linear Assignment Cost}
\begin{dfn}[Linear Assignment Cost LAC]
\label{dfn:LAC_clouds}
Recall that Definition~\ref{dfn:SM} introduced the bottleneck distance $\BD$ on matrices considered under re-ordering of columns.
For any clouds $A,B\subset\R^n$ of $m$ points, consider the \emph{Linear Assignment Cost}  
$\LAC(A,B)=\min\limits_{g}\sum\limits_{P\in\WMI(A)} \BD(P,g(P))$ minimised \cite{jonker1987shortest} over all bijections $g:\WMI(A)\to\WMI(B)$ of full Weighted Matrices Invariants consisting of all $N=m(m-1)\dots(m-n+1)$ equivalence classes of matrices.
\edfn
\end{dfn}

\index{Linear Assignment Cost}
\begin{lem}[$\LAC$ on clouds, {\cite[Lemma~6.2]{kurlin2024polynomial}}]
\label{lem:LAC}
\textbf{(a)}
The Linear Assignment Cost from Definition~\ref{dfn:LAC_clouds} satisfies all metric axioms on 
clouds under rigid motion. 
\medskip

\noindent
\textbf{(b)}
Let $O(A)$ be any mirror image of a cloud $A\subset\R^n$.
Then $\min\{\LAC(A,B),\LAC(O(A),B)\}$ is a metric on classes of clouds under  isometry.
\elem
\end{lem}

\index{Linear Assignment Cost}
\begin{thm}[time complexity of $\LAC$ on $\WMI$s, {\cite[Theorem~6.3]{kurlin2024polynomial}}]
\label{thm:LAC_time}
For any clouds $A,B\subset\R^n$ of $m$ points, the invariants $\WMI(A),\WMI(B)$ consists of at most $N=m(m-1)\dots(m-n+1)=O(m^{n-1})$ matrices.
Then the metric $\LAC(A,B)$ from Definition~\ref{dfn:LAC_clouds} 
can be computed in time $O(m^{1.5}(\log^n m)N^2+N^3)=O(m^{2n-0.5}\log^n m+m^{3n-3})$.
If $n=2$, the time is $O(m^{3.5}\log m)$.
\ethm
\end{thm}

The worst-case estimate $N=O(m^{n-1})$ of the size (number of matrices in) $\WMI(A)$ is very rough.
If the covariance matrix $\cov(A)$ has equal eigenvalues, $\WMI(A)$ is often smaller due to extra symmetries of $A$.
\medskip

However, for $n=2$, even the rough estimate of the LAC time $O(m^{3.5}\log m)$ improves the time $O(m^5\log m)$ for computing the exact Hausdorff distance between $m$-point clouds under Euclidean motion in $\R^2$.
\medskip

Since real noise may include erroneous points, it is practically important to continuously quantify the similarity between close clouds consisting of different numbers of points.
The weights of matrices allow us to match them more flexibly via the Earth Mover's Distance \cite{rubner2000earth} than via strict bijections $\WMI(A)\to\WMI(B)$. 
The Weighted Matrices Invariant $\WMI(A)$ can be considered as a finite distribution $C=\{C_1,\dots,C_k\}$ of matrices (equivalent up to re-ordering columns) with weights.

\index{Earth Mover's Distance}
\begin{dfn}[Earth Mover's Distance on weighted distributions]
\label{dfn:EMD}
Let $C=\{C_1,\dots,C_k\}$ and $D=\{D_1,\dots,D_l\}$ be finite unordered set of objects with weights $w(C_i)$, $i=1,\dots,k$, and $w(D_j)$, $j=1,\dots,l$, respectively such that $\sum\limits_{i=1}^k w(C_i)=1=\sum\limits_{j=1}^l w(D_j)$.
Let $d$ be a \emph{ground metric} between any $C_i$ and $D_j$.
A \emph{flow} from $C$ to $D$ is a $k\times l$ matrix  whose entry $f_{ij}$ represents a \emph{flow} from $C_i$ to $D_j$.
The \emph{Earth Mover's Distance} is  
$\EMD(C,D)=\sum\limits_{i=1}^{k} \sum\limits_{j=1}^{l} f_{ij} d(C_i,D_j)$ minimised for $f_{ij}\in[0,1]$ subject to 
$\sum\limits_{j=1}^{l} f_{ij}\leq w(C_i)$, $i=1,\dots,k$, 
$\sum\limits_{i=1}^{k} f_{ij}\leq w(D_j)$, $j=1,\dots,l$, and
$\sum\limits_{i=1}^{k}\sum\limits_{j=1}^{l} f_{ij}=1$.
\edfn
\end{dfn}

The first condition $\sum\limits_{j=1}^{l} f_{ij}\leq w(C_i)$ means that not more than the weight $w(C_i)$`flows' into all objects $D_j$ via $f_{ij}$, $j=1,\dots,l$. 
Similarly, the second condition $\sum\limits_{i=1}^{k} f_{ij}\leq w(D_j)$ means that all $f_{ij}$ `flow' from $C_i$, $i=1,\dots,k$ into $D_j$ up to its weight $w(D_j)$.
\medskip

The last condition
$\sum\limits_{i=1}^{k}\sum\limits_{j=1}^{l} f_{ij}=1$
 forces to `flow' all $C_i$ to all $D_j$.  
The EMD is a partial case of more general Wasserstein metrics \cite{vaserstein1969markov} in transportation theory \cite{kantorovich1960mathematical}.
For finite distributions as in Definition~\ref{dfn:EMD}, the metric axioms for $\EMD$ were proved in \cite[appendix]{rubner2000earth}. 
$\EMD$ can compare any weighted distributions of different sizes.
Instead of the bottleneck distance $\BD$ on columns on $\PCM$ matrices, one can consider $\EMD$ on the distributions of columns (with equal weights) in these matrices.

\begin{lem}[time complexity of $\EMD$, {\cite[Lemma~6.5]{kurlin2024polynomial}}]
\label{lem:EMD_matrices}
Any matrix $P$ of a size $n\times m(P)$ can be considered as a distribution of $m(P)$ columns with equal weights $\frac{1}{m(P)}$.
For two such matrices $P,Q$ having the same number $n$ of rows but potentially different numbers $m(P),m(Q)$ of columns, measure the distance between any columns by the Chebyshev metric $L_{\infty}$ in $\R^n$.
For the matrices $P,Q$ considered as weighted distributions of columns, the Earth Mover's Distance $\EMD(P,Q)$ can be computed in time $O(m^3\log m)$, where $m=\max\{m(P),m(Q)\}$.
\elem
\end{lem}

\begin{thm}[time of $\EMD$ on clouds, {\cite[Theorem~6.6]{kurlin2024polynomial}}]
\label{thm:EMD_time}
Let clouds $A,B\subset\R^n$ of up to $m$ points have pre-computed invariants $\WMI(A),\WMI(B)$ of sizes at most $N\leq m(m-1)\dots(m-n+1)=O(m^{n-1})$.
Measure the distance between any matrices $P\in\WMI(A)$ and $Q\in\WMI(B)$ as 
$\EMD(P,Q)$ from Lemma~\ref{lem:EMD_matrices}.
Then  the Earth Mover's Distance $\EMD(\WMI(A),\WMI(B))$ 
from Definition~\ref{dfn:EMD} 
can be computed in time $O(m^3(\log m) N^2+N^3\log N)=O((m^{2n+1} +nm^{3n-3})\log m)$.
\ethm
\end{thm}

\begin{exa}[$\EMD$ for a square and an equilateral triangle]
\label{exa:triangle-vs-square}
Let $A_4$ and $A_3$ be the vertex sets of a square and equilateral triangle inscribed into the circle of a radius $r$ in Example~\ref{exa:reg_clouds}.
$\PCM(A_3)=
\left( \begin{array}{ccc} 
r & -r/2 & -r/2 \\
0 & r\sqrt{3}/2 & -r\sqrt{3}/2
\end{array} \right)$ and
$\PCM(A_4)=
\left( \begin{array}{cccc} 
r & 0 & 0 & -r \\
0 & r & -r & 0
\end{array} \right)$.
Notice that switching the signs of the 2nd row keeps the PCI matrices the same up to permutation of columns.
The weights of the three columns in $\PCM(A_3)$ are $\dfrac{1}{3}$.
The weights of the four columns in $\PCM(A_4)$ are $\dfrac{1}{4}$.
The EMD optimally matches the identical first columns of $\PCM(A_3)$ and $\PCM(A_4)$ with weight $\dfrac{1}{4}$ contributing the cost $0$.
The remaining weight $\dfrac{1}{3}-\dfrac{1}{4}=\dfrac{1}{12}$ of the first column $\vect{r}{0}$ in $\PCM(A_3)$ can be equally distributed between the closest (in the $L_{\infty}$ distance) columns $\vect{0}{\pm r}$ contributing the cost $\dfrac{r}{12}$.
The column $\vect{-r}{0}$ in $\PCM(A_4)$ has equal distances $L_{\infty}=\dfrac{r}{2}$ to the last columns $\vect{-r/2}{\pm r\sqrt{3}/2}$ in $\PCM(A_3)$ contributing the cost $\dfrac{r}{8}$.
Finally, the distance $L_{\infty}=\dfrac{r}{2}$ between the columns $\vect{0}{\pm r}$ and $\vect{-r/2}{\pm r\sqrt{3}/2}$ with the common signs is counted with the weight $\dfrac{5}{24}$ and contributes the cost $\dfrac{5r}{48}$.
The final optimal flow $(f_{jk})$ matrix 
$\left( \begin{array}{cccc} 
1/4 & 1/24 & 1/24 & 0 \\
0 & 5/24 & 0 & 1/8 \\
0 & 0 & 5/24 & 1/8
\end{array} \right)$
gives $\EMD(\PCM(A_3),\PCM(A_4))=\dfrac{r}{12}+\dfrac{r}{8}+\dfrac{5r}{48}=\dfrac{5r}{16}$.
\eexa
\end{exa}

\bibliographystyle{plain}
\bibliography{Geometric-Data-Science-book}

%
%
%

\chapter{Fast isometry invariants of finite clouds of unordered points}
\label{chap:PDD-finite} 

\abstract{This chapter adapts the general Geo-Mapping Problem to finite clouds under isometry in a metric space.
We start by discussing Sorted Pairwise Distances, which distinguish all generic clouds under Euclidean isometry.
Then we introduce the stronger invariant $\PDD$ (Pointwise Distance Distribution), which was recently proved to be complete for all 4-point clouds under isometry in any $\R^n$.
}

\section{Geo-mapping problem for generic clouds in a metric space}
\label{sec:distances}

This chapter studies finite clouds of unordered points under isometry in any metric space, though $\R^n$ remains an important partial case.
Problem~\ref{pro:metric_space} adjusts Geo-Mapping Problem~\ref{pro:geocodes} to include only realistically achievable conditions in an arbitrary metric space.
One realistic condition is a general position, as formally defined below.

\index{general position}
\begin{dfn}[general position in a metric space]
\label{dfn:general_position}
Let $A$ be a set of $m$ points in a metric space $M$ with a metric $d_M$.
A cloud $A$ is called \emph{generic} (or in a \emph{general position} in $M$) if all $\dfrac{m(m-1)}{2}$ inter-point distances $d_M(p,q)$ for $p,q\in A$ are not solutions of a certain polynomial equation.
\edfn
\end{dfn}

To give an example of a general position in 
 Definition~\ref{dfn:general_position}, recall the the \emph{lexicographic order} on ordered pairs: $(i,j)<(k,l)$ if $i<j$ or $i=j$ and $k<l$.

\begin{exa}[a polynomial for a general position]
\label{exa:general_position} 
Let $p_1,\dots,p_m$ be all points of a set $A$ in a space $M$ with a metric $d_M$.
Set $P(A)=\prod\limits_{(i,j)\neq (k,l)} (d_{i,j}-d_{k,l})$, where $d_{i,j}=d_M(p_i,p_j)$ for $1\leq i<j\leq m$.
If $m=3$, then $P(A)=(d_{1,2}-d_{1,3})(d_{1,2}-d_{2,3})(d_{1,3}-d_{2,3})$.
By Definition~\ref{dfn:general_position}, the polynomial condition $P(A)\neq 0$ describes the general position for all clouds $A\subset M$ that have distinct distances.
\eexa
\end{exa} 

The concept of a general position allows us to weaken the concept of a complete invariant from Definition~\ref{dfn:invariants}(b) to a generically complete invariant below.

\index{generically complete invariant}
\begin{dfn}[generically complete invariant]
\label{dfn:generically_complete}
Let $\sim$ be an equivalence relation on finite sets in a metric space $M$.
An invariant $I$ under this equivalence is  \emph{generically complete} if the implication $I(A)=I(B) \Rightarrow A\sim B$ holds for all sets $A,B\subset M$ in a general position for a certain polynomial on inter-point distances in Definition~\ref{dfn:general_position}.
\edfn
\end{dfn}

Instead of the full completeness in Problem~\ref{pro:geocodes}, Problem~\ref{pro:metric_space} asks for a more realistic generic completeness under isometry in the sense of Definition~\ref{dfn:generically_complete}.

\begin{pro}[geo-mapping for generic clouds under isometry in a metric space]
\label{pro:metric_space}
For a space $M$ with a metric $d_M$, find an isometry invariant $I$ of generic clouds of unordered points in $M$ with values in a metric space satisfying the following conditions.
\myskip

\noindent
\tb{(a)} 
\emph{Generic completeness:} 
any generic clouds $A,B\subset M$ are isometric  in $M$, i.e. $A\simeq B$, if and only if $I(A)=I(B)$.
\myskip

\noindent
\tb{(b)} 
\emph{Reconstruction:} 
any generic cloud $A\subset M$ can be reconstructed from its invariant value $I(A)$, uniquely under isometry in $M$.
\myskip

\noindent
\tb{(c)} \emph{Metric:} 
there is a distance $d$ on the \emph{invariant space} $\{I(A) \vl A\subset M\}$ satisfying all metric axioms in Definition~\ref{dfn:metrics}(a). 
\myskip

\nt
\tb{(d)} 
\emph{Continuity:} 
there is a constant $\la$ such that, for any $\ep>0$, if $B$ is obtained from $A$ by perturbing every point of $A$ up to $\ep$ in the metric $d_M$, then $d(I(A),I(B))\leq \la\ep$.  
\myskip

\nt
\tb{(e)} 
\emph{Computability:} 
for a fixed metric space $M$, the invariant $I(A)$, the metric $d(I(A),I(B))$, and a reconstruction of $A\subset M$ from $I(A)$ can be computed in a time that depends polynomially on the maximum size $\max\{|A|,|B|\}$ of clouds $A,B$.
\epro 
\end{pro} 

Definition~\ref{dfn:SPD} introduces the invariant that nearly solved Problem~\ref{pro:metric_space} in 2004 by \cite[Theorem~2.6]{boutin2004reconstructing} under Euclidean isometry in $\R^n$.
Though this seminal work \cite{boutin2004reconstructing} talks about reconstructing point configurations, the main result actually proves the generic completeness of the following invariant under isometry in $\R^n$ as stated in \ref{pro:metric_space}(a).
 
\index{sorted pairwise distances} 
\begin{dfn}[Sorted Pairwise Distances $\SPD$]
\label{dfn:SPD}
For any finite cloud $A$ of unordered points in a metric space $M$, the vector $\SPD(A)$ of \emph{Sorted Pairwise Distances} consists of all $\dfrac{m(m-1)}{2}$ distances between all points of $A$, written in increasing order.
\edfn
\end{dfn}
 
Any isometry in a metric space $M$ preserves distances and hence $\SPD(A)$ for any cloud $A\subset M$.
In the Euclidean case $M=\R^n$, if $A$ consists f $m=3$ points, $\SPD(A)$ coincides with the geocode of three inter-point distances, which classified all triangles under isometry in Example~\ref{exa:geocodes_m=3}(a).
For a cloud $A\subset\R^n$ of any $m$ unordered points, 
\cite[Theorem~2.6]{boutin2004reconstructing} proved that $\SPD(A)$ is generically complete  under isometry in $\R^n$.  
We leave as an exercise that $\SPD(A)$ is Lipschitz continuous, for example, in any Minkowski metric $L_q$, because this Lipschitz continuity will be proved for stronger invariants in the next section.
Since $\SPD(A)$ needs a quadratic time of the size $|A|=m$, this invariant solves Problem~\ref{pro:metric_space} for generic clouds in $\R^n$.
\myskip

\begin{figure}[h!]
\centering
\includegraphics[width=\textwidth]{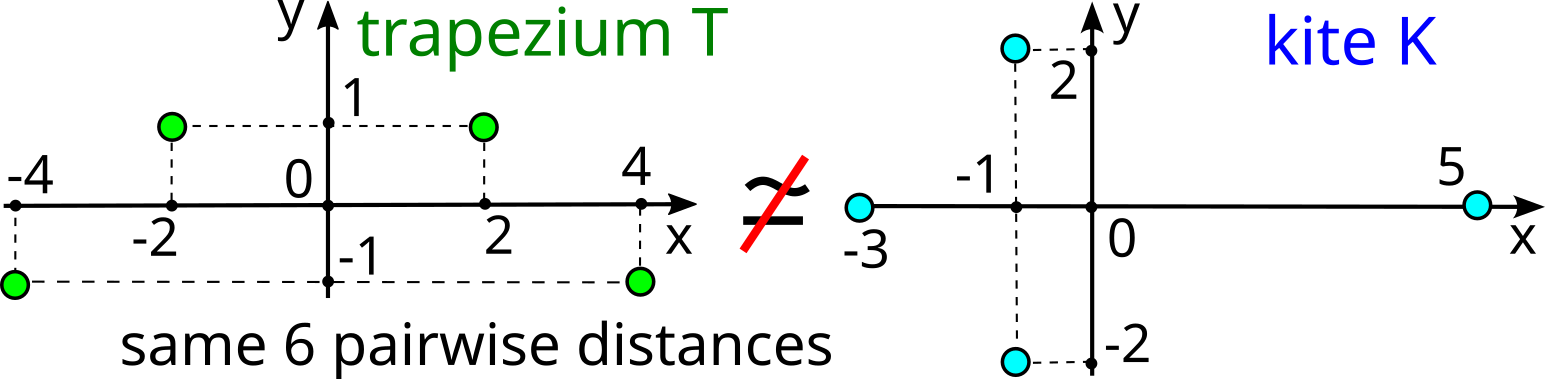}
\caption{
Non-isometric clouds of 4 points with the same 6 pairwise distances.
\textbf{Left}: the trapezium $T$ has the vertices $(\pm 2,1)$, $(\pm 4,-1)$.
\textbf{Right}: the kite $K$ has the vertices $(5,0)$, $(-3,0)$, $(-1,\pm 2)$.}
\label{fig:4-point_clouds_origin}
\end{figure}

\begin{figure}[h!]
\centering
\includegraphics[width=\textwidth]{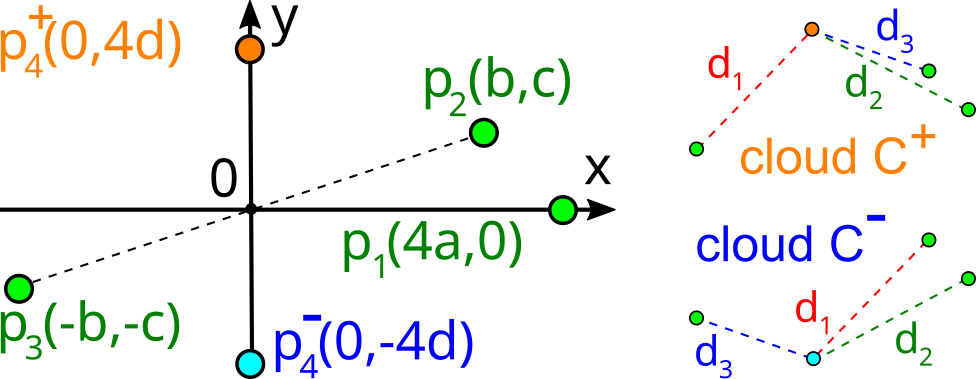}
\caption{Infinitely many non-isometric clouds $C^+\not\simeq C^-$ 
depending on free parameters $a,b,c,d>0$ \cite{caelli1979generating}.}
\label{fig:4-point_clouds_family}
\end{figure}

However, infinitely many counterexamples to the completeness of $\SPD$ have been known at least since 1979 \cite{caelli1979generating} even for $m=4$ points in $\R^2$.  
Fig.~\ref{fig:4-point_clouds_origin} shows the most famous pair of a trapezium and a kite, which inspired the flagship image of Geometric Data Science in Fig.~\ref{fig:GDS-foundations}~(right).
Fig.~\ref{fig:4-point_clouds_family} illustrates infinitely many non-isometric 4-point clouds in $\R^2$, which share three points $p_1,p_2,p_3$ in green and differ only in points $p_4^\pm$, but share 3 distances $d_1,d_2,d_3$ from $p_4^\pm$ to three others. 
\myskip

We can leverage the Euclidean structure of $\R^n$ to introduce a simpler invariant below.

\begin{dfn}[Sorted Radial Distances $\SRD$]
\label{dfn:SRD}
For any finite cloud $A$ of unordered points in $\R^n$, a translation can fix the centre of mass $\bar A$ of $A$ at the origin $0\in\R^n$.
The vector $\SRD(A)$ of \emph{Sorted Radial Distances} consists of all $m$ Euclidean distances from $\bar A=0$ to all points of $A$, written in decreasing order.
\edfn
\end{dfn}

The orders of distances in Definitions~\ref{dfn:SPD} and~\ref{dfn:SRD} are motivated by applications to molecules.
The most characteristic inter-atomic distances are the smallest ones between bonded atoms at the beginning of the $\SPD$, written in increasing order.
\myskip

On another hand, the simpler invariant $\SRD(A)$ describes the global shape of a molecule by the largest distances to atoms from the centre of mass. 

\begin{exa}[$\SPD$ and $\SRD$ for 4-point clouds in {Fig.~\ref{fig:4-point_clouds_origin}}]
\label{exa:SPD4-points}
The vertex sets $T,K$ of the trapezium and kite in
Fig.~\ref{fig:4-point_clouds_origin} have $\SPD=(2\sqrt{2},2\sqrt{2},4,2\sqrt{10},2\sqrt{10},8)$, but are distinguished by 
$\SRD(T)=(\sqrt{17},\sqrt{17},\sqrt{5},\sqrt{5})$ and $\SRD(K)=(5,3,\sqrt{5},\sqrt{5})$. 
\eexa
\end{exa}

The next two sections follow the finite (non-periodic) case of paper \cite{widdowson2025pointwise}.

\section{Pointwise Distance Distributions of unordered points}
\label{sec:PDD}

\index{sorted pairwise distances} 
This section defines our main isometry invariant, which we first introduced for periodic point sets \cite{widdowson2022resolving} in 2022 and only after that understood its importance in the finite case, where it was previously studied under the name of a \emph{local distribution of distances} \cite{memoli2011gromov}. 

\index{Pointwise Distance Distribution} 
\begin{dfn}[Pointwise Distance Distribution $\PDD(A;k)$ for a finite cloud $A$]
\label{dfn:PDD_finite}
Let $A=\{p_1,\dots,p_m\}$ be a finite cloud of unordered points in a metric space $M$.
\sskip

Fix an integer $k\geq 1$.
For every point $p\in A$, let $d_1(p)\leq\dots\leq d_k(p)$ be the distances from $p$ to its $k$ nearest neighbours in $A$.
The matrix $D(A;k)$ has $m$ rows consisting of the distances $d_1(p_i),\dots,d_k(p_i)$ for $i=1,\dots,m$.
If any $l\geq 2$ rows coincide, we collapse them into a single row with the weight $\dfrac{l}{m}$.
\sskip

The resulting matrix of maximum $m$ rows and $k+1$ columns, including the extra column of weights, is called the \emph{Pointwise Distance Distribution} $\PDD(A;k)$. 
\edfn
\end{dfn}

The rows of $\PDD(A;k)$ are unordered, though we might write them in a lexicographic order only for convenience.
Hence $\PDD(A;k)$ can be considered a weighted distribution of (say) $m$ rows of $k$ distances.
Each row can also be interpreted as a point in $\R^k$.
\myskip

Then $\PDD(A;k)$ can be viewed as a cloud of $m$ unordered points in $\R^k$.
The crucial difference with the original cloud $A$ under isometry in $\R^n$ is the fixed coordinate system for $\PDD(A;k)\subset\R^k$, not under any equivalence.   

\begin{exa}[$\PDD$ for 4-point clouds $T,K$ in {Fig.~\ref{fig:4-point_clouds_origin}}]
\label{exa:PDD4-points}
Table~\ref{tab:ordered_distances_T+K} shows the $4\times 3$ matrices $D(S;3)$ from Definition~\ref{dfn:PDD_finite}.
The matrix $D(T;3)$ in Table~\ref{tab:ordered_distances_T+K} has two pairs of identical rows, so the matrix $\PDD(T;3)$ consists of two rows of weight $\frac{1}{2}$ below.
\myskip

The matrix $D(K;3)$ in Table~\ref{tab:ordered_distances_T+K} has only one pair of identical rows, so $\PDD(K;3)$ has three rows of weights $\frac{1}{2}$, $\frac{1}{4}$, $\frac{1}{4}$.
Then $T,K$ are distinguished by $\PDD$s even for $k=1$.

\begin{table}[h!]
\caption{Each point of $T,K\subset\R^2$ in Figure~\ref{fig:4-point_clouds_origin} has three distances to other points in increasing order.
After keeping only distances (not neighbours), the resulting $\PDD$s distinguish $T\not\simeq K$.} 
\label{tab:ordered_distances_T+K}
\centering
\begin{tabular}{l|r|r|r}
points of $T$ & distance to neighbour 1 & distance to neighbour 2 & distance to neighbour 3 \\
\hline
$(-4,-1)$ & $2\sqrt{2}$ to $(-2,+1)$ & $2\sqrt{10}$ to $(+2,+1)$  & $8$ to $(+4,-1)$ \\
$(+4,-1)$ & $2\sqrt{2}$ to $(+2,+1)$ & $2\sqrt{10}$ to $(-2,+1)$ & $8$ to $(-4,-1)$ \\
$(-2,+1)$ & $2\sqrt{2}$ to $(-4,-1)$ & $8$ to $(+2,+1)$ & $2\sqrt{10}$ to $(+4,-1)$ \\
$(+2,+1)$ & $2\sqrt{2}$ to $(+4,-1)$ & $4$ to $(-2,+1)$ & $2\sqrt{10}$ to $(-4,-1)$  
\end{tabular}
\smallskip

\begin{tabular}{l|r|r|r}
points of $K$ & distance to neighbour 1 & distance to neighbour 2 & distance to neighbour 3  \\
\hline
$(-3,0)$ & $2\sqrt{2}$ to $(-1,-2)$ & $2\sqrt{2}$  to $(-1,+2)$ & $8$ to $(5,0)$ \\
$(+5,0)$ & $2\sqrt{10}$ to $(-1,-2)$ & $2\sqrt{10}$ to $(-1,+2)$ & $8$ to $(-3,0)$ \\
$(-1,-2)$ & $2\sqrt{2}$ to $(-3,0)$ & $4$ to $(-1,+2)$ & $2\sqrt{10}$ to $(5,0)$\\
$(-1,+2)$ & $2\sqrt{2}$ to $(-3,0)$ & $4$ to $(-1,-2)$ & $2\sqrt{10}$ to $(5,0)$
\end{tabular}
\end{table}

\noindent
$\PDD(T)=\left(\begin{array}{c|ccc}
1/2 & 2\sqrt{2} & 4 & 2\sqrt{10} \\
1/2 & 2\sqrt{2} & 2\sqrt{10} & 8
\end{array}\right)\neq 
\PDD(K)=\left(\begin{array}{c|ccc}
1/4 & 2\sqrt{2} & 2\sqrt{2} & 8 \\
1/2 & 2\sqrt{2} & 4 & 2\sqrt{10} \\
1/4 & 2\sqrt{10} & 2\sqrt{10} & 8
\end{array}\right)$.
\eexa
\end{exa}

Since any isometry preserves distances, $\PDD(A;k)$ is an isometry invariant of $A$.
The brute-force algorithm for $\PDD(A;k)$ needs only a quadratic time in the size $|A|=m$.
In a general metric space with certain expansion constants, we found counter-examples \cite{elkin2022counterexamples} to past estimates for a parametrised complexity a nearest neighbour search and proved new linear-time complexities  \cite{elkin2023new} with extra parameters depending, for example, on a dimension $n$.  
Hence, the invariant $\PDD(A;k)$ satisfies the computability in \ref{pro:metric_space}(e).
\myskip
 
Interpreting $\PDD(A;k)$ as a discrete distribution of rows (or points in $\R^k$) with weights as probabilities allows us to compare $\PDD$s by many metrics on probability distributions.
If we use the Earth Mover's Distance from Definition~\ref{dfn:EMD} with a ground metric $L_q$ on rows of $\PDD$, we denote the resulting metric by $\EMD_q$ for all parameters $q\in[1,+\infty]$.
For $\PDD(A;k)$, the notation $\EMD$ without any subscript means that the ground metric is the Root Mean Square $\RMS=\frac{L_2}{\sqrt{k}}$.
\myskip

The EMD satisfies all metric axioms \cite[appendix]{rubner2000earth}, needs $O(m^3\log m)$ time for distributions of a maximum size $m$, and can be approximated in $O(m)$ time \cite{shirdhonkar2008approximate}.

\index{Pointwise Distance Distribution} 
\begin{thm}[Lipschitz continuity of $\PDD$ for a finite cloud, {\cite[Theorem~4.2(a)]{widdowson2025pointwise}}]
\label{thm:PDD_continuous_finite} 
Let $A$ be a finite cloud in a space $M$ with a metric $d_M$.
For any $\ep>0$, let $B$ be obtained from $A$ by perturbing every point of $A$ up to $\ep$ in the metric $d_M$.
Fix any real $q\in[1,+\infty]$ and an integer $k\geq 1$.
Interpret $\sqrt[q]{k}$ as 1 in the limit case $q=+\infty$.
Then $\EMD_q(\PDD(A;k),\PDD(B;k))\leq 2\ep\sqrt[q]{k}$.
\ethm
\end{thm}

For any cloud $A\subset\R^n$ of $m$ unordered points, the vector $\SPD(A)$ of Sorted Pairwise Distances obtained from $\PDD(A;m-1)$ by writing all distances in a single distribution and collapsing each pair of equal distances into one.
Indeed, any distance $|p_i-p_j|$ appears in both rows $i,j$ of $\PDD(A;m-1)$.
Due to Example~\ref{exa:PDD4-points}, $\PDD(A;m-1)$ is strictly stronger than $\SPD(A)$.
Due to this strength, the generic completeness of $\PDD(A;m-1)$ under isometry in $\R^n$ is much easier to prove than for $\SPD(A)$.

\index{Pointwise Distance Distribution} 
\begin{thm}[generic completeness of $\PDD$ for a finite cloud, {\cite[Theorem~5.1]{widdowson2025pointwise}}]
\label{thm:PDD_finite_gen_complete}
Any cloud $A\subset\R^n$ of $m$ unordered points with distinct inter-point distances can be reconstructed from $\PDD(A;m-1)$, uniquely under isometry. 
\end{thm}
\begin{proof}
Since all inter-point distances are distinct, every such distance $|p-q|$ between points $p,q\in A$ appears twice in $\PDD(A;m-1)$: once in the row of $p$ and once in the row of $q$.
Hence, after choosing an arbitrary order of points, we can use $\PDD(A;m-1)$ to reconstruct the classical distance matrix on ordered points.
This distance matrix determines $A\subset\R^n$ uniquely under isometry \cite{schoenberg1935remarks}. 
\end{proof}

The following open conjecture should be understandable to schools students. 

\begin{conj}[completeness of $\PDD$ in $\R^2$]
\label{conj:PDD_complete_n=2}
Any cloud $C\subset\R^2$ of $m$ unordered points can be reconstructed from $\PDD(C;m-1)$, uniquely under isometry in $\R^2$. 
\epro
\end{conj}

In $\R^3$, the known non-isometric clouds with the same $\PDD$ inspired the stronger invariant, which will distinguish all these examples in the next chapter.
In a general metric space, Problem~\ref{pro:metric_space} is notoriously hard, but provides targets for further research.   

\section{Extending the side-side-side theorem from 3 to 4 points in $\R^n$}
\label{sec:4-points}

Many authors considered criteria of congruence for plane quadrilaterals
\cite{vance1982minimum}, whose vertices are ordered.
The $m\times m$ matrix of pairwise distances \cite{schoenberg1935remarks} and the Gram matrix of scalar products \cite{weyl1946classical} are complete and continuous invariants of $m$ ordered points under isometry in $\R^n$, known at least since 1935.
The extension of this approach to $m$ unordered points leads to the exponential complexity because of $m!$ permutations. 
\myskip

For $m=4$ unordered points, Theorem~\ref{thm:PDD_complete_m<5} will prove the completeness of $\PDD(C;m-1)$ under isometry in any $\R^n$.
For any $m$, the invariant $\PDD(C;m-1)$ can be computed in quadratic time $O(m^2)$.
For $m=4$, $\PDD(C;3)$ contains only 12 numbers (6 pairs of distances between 4 points), while $4!=24$ distance matrices on 4 points contain at least 144 numbers if we take only distances above the diagonal. 

\index{Pointwise Distance Distribution} 
\begin{thm}[completeness of $\PDD$ for $m\leq 4$ points, {\cite[Theorem~5.3]{widdowson2025pointwise}}]
\label{thm:PDD_complete_m<5}
The Pointwise Distance Distribution $\PDD(C;m-1)$ from Definition~\ref{dfn:PDD_finite} is a complete isometry invariant of all clouds $C\subset\R^n$ of any $m\leq 4$ unordered points. 
\ethm
\end{thm}

Since Theorem~\ref{thm:PDD_complete_m<5} finally extends the side-side-side criterion of congruence to $m=4$ unordered points, without relying on brute-force permutations, we include the detailed proof, which previously appeared only in supplementary materials of \cite{widdowson2025pointwise}. 
\myskip

If a cloud $A$ of $m$ points has a line or plane of symmetry $L$ in $\R^2$ or $\R^3$, then all points $A\setminus L$ split into pairs of points that are symmetric in $L$ and hence have equal rows in $\PDD(A;m-1)$.
Lemma~\ref{lem:PDD_symmetry} shows that the converse holds for $m=4$.

\begin{figure}[h!]
\centering
\includegraphics[width=\textwidth]{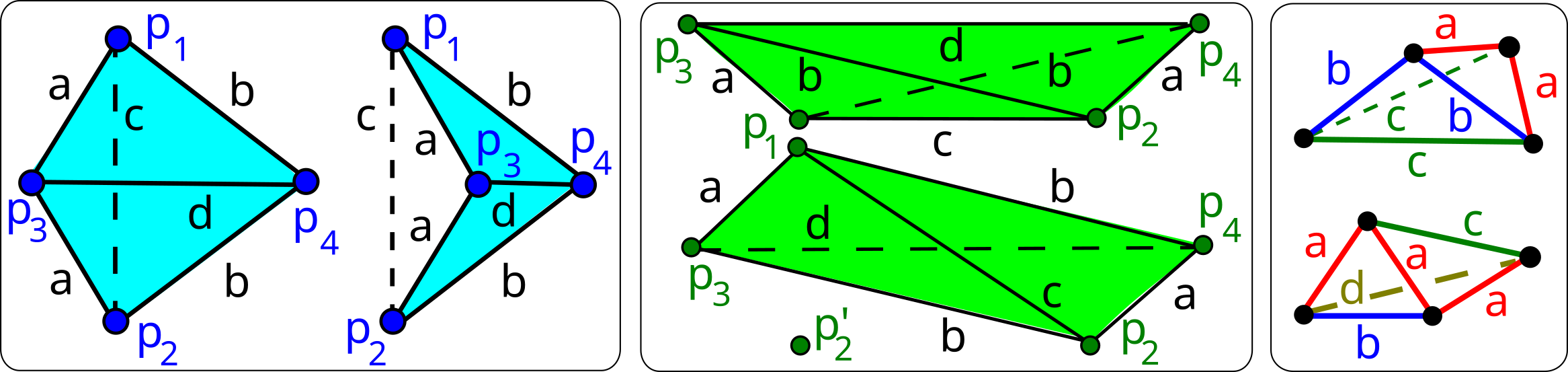}
\caption{\textbf{Left}: 
in $\R^2$, the convex and non-convex kites have two equal rows $\{a,b,c\}$ in $\PDD(A;3)$ and are distinguished by $d=|p_3-p_4|$, see Lemma~\ref{lem:PDD_symmetry}.
\textbf{Middle}: an isosceles trapezoid and parallelogram have $\PDD(A;3)$ with two pairs of equal rows $\{a,b,c\}$ and $\{a,b,d\}$, e.g. a rectangle for $c=d$. 
\textbf{Top right}: a \emph{trisosceles} cloud.
\textbf{Bottom right}: a \emph{3-chain-equal} cloud, see Example~\ref{exa:trisosceles+trichainequal}. 
}
\label{fig:symmetric+special_clouds}
\end{figure}

\begin{lem}[$\PDD$ detects symmetry of $m=4$ points, {\cite[Lemma~SM3.5]{widdowson2025pointwise}}]
\label{lem:PDD_symmetry}
For any cloud $A\subset\R^n$ of $m=4$ points for $n=2,3$, if $\PDD(A;3)$ has two equal rows, then $A$ is either (1) mirror symmetric in the plane passing through two points of $A$ orthogonally to the line segment joining the other points of $A$, or (2) symmetric by the $180^\circ$ degree rotation around the line through the mid-points of two pairs of points of $A$.   
If $n=2$, then $A$ defines a kite, or a parallelogram or an isosceles trapezoid; see Fig.~\ref{fig:symmetric+special_clouds}. 
\elem
\end{lem}

\begin{proof}[Proof of Lemma~\ref{lem:PDD_symmetry}]
Let points $p_1,p_2\in A$ have the same row $a\leq b\leq c$ in $\PDD(C;3)$.
One of the distances $a,b,c$ is between the points $p_1,p_2$.
Without loss of generality, assume that $|p_1-p_2|=c$.
Then $p_1,p_2$ have distances $a,b$ to the points $p_3,p_4\in A\setminus\{p_1,p_2\}$, but it is unknown which distance corresponds to which point.
\smallskip

\emph{Isosceles case}.
Let $|p_1 - p_3| = a = |p_2 - p_3|$ and $|p_1 - p_4| = b = |p_2 - p_4|$, see Fig.~\ref{fig:symmetric+special_clouds}~(left).
Then $A$ has two equal triangles $\triangle p_1 p_3 p_4=\triangle p_2 p_3 p_4$
 and two isosceles triangles $\triangle p_3 p_1 p_2$ and $\triangle p_4 p_1 p_2$ with equal sides at $p_3,p_4$, respectively.
Let $L$ be the plane that passes through $p_3,p_4$ and is orthogonal to the line segment $[p_1,p_2]$.
Then the mirror reflection in $L$ swaps $p_1,p_2$. 
If $n=2$, $A$ defines a (non-)convex kite. 
\smallskip
 
\emph{Non-isosceles case}.
Then $|p_1-p_3|=a=|p_2-p_4|$ and $|p_2-p_3|=b=|p_1-p_4|$, see Fig.~\ref{fig:symmetric+special_clouds}~(middle).
Let $L$ be the perpendicular bisector of the line segment $[p_3,p_4]$.
The mirror reflection in $L$ swaps $p_3\lra p_4$ and either swaps $p_1\lra p_2$ 
(then $A$ defines an isosceles trapezoid in $\R^2$) or maps $p_2$ to $p'_2$, so that $p_1,p'_2,p_3,p_4$ satisfy the previous case.
In the latter case, the composition with the reflection in the plane through $p_3,p_4$ orthogonal to $[p_1,p'_2]$ is the $180^\circ$ degree rotation that
swaps the points as $p_1\lra p_2$ and $p_3\lra p_4$.
If $n=2$, then $A$ defines a parallelogram, see  Fig.~\ref{fig:symmetric+special_clouds}~(bottom middle).
\end{proof}

\begin{exa}[trisosceles and 3-chain-equal clouds in $\R^3$]
\label{exa:trisosceles+trichainequal}
Fig.~\ref{fig:symmetric+special_clouds}~(right) shows \emph{trisosceles} and \emph{3-chain-equal} clouds that have 3 pairs of equal distances and a chain of 3 equal distances, their $\PDD$s are
$\left(\begin{array}{ccc} a & a & c \\ a & b & b \\ a & b & c \\ b & c & c \end{array} \right)$ and
$\left(\begin{array}{ccc} a & a & b \\ a & a & c \\ a & b & d \\ a & c & d \end{array} \right)$, 
respectively.
\eexa
\end{exa}

\begin{proof}[Proof of Theorem~\ref{thm:PDD_complete_m<5}]
\emph{Case $m=2$}.
Any cloud $A\subset\R^n$ of $m=2$ unordered points $p_1,p_2$ (labelled only for convenience) has $\PDD(A;1)$ consisting of the single distance $|p_1-p_2|$, which uniquely determines $A$ under isometry in any $\R^n$.  
\smallskip

\emph{Case $m=3$}.
Any cloud $A\subset\R^n$ of $m=3$ unordered points with pairwise distances $a\leq b\leq c$ has $\PDD(A;2)
=\left(\begin{array}{cc} a & b \\ a & c \\ b & c \end{array} \right)$.
The (lexicographically) first row of $\PDD(A;2)$ gives us $a\leq b$.
Each of the remaining two rows of $\PDD(A;2)$ should contain at least one value of $a$ or $b$, including in all degenerate cases such as $a=b$.
Removing these repeated values from the other two rows gives us $c$, also in the case $b=c$. 
So $\PDD(A;2)$ identifies $a\leq b\leq c$ and hence $A$, uniquely under isometry in any $\R^n$.  
\smallskip

\emph{Case $m=4$}, then $n\leq 3$. 
For a cloud $A\subset\R^3$ of $m=4$ unordered points, $\PDD(A;3)$ is a $4\times 3$ matrix.
Assume that $\PDD(A;3)$ has two equal rows $a\leq b\leq c$. 
\smallskip

\emph{Isosceles case}.
In the first case of Lemma~\ref{lem:PDD_symmetry} in Fig.~\ref{fig:symmetric+special_clouds}~(left), 
$\PDD(A;3)$ has two more rows $\{a,a,d\}$ and $\{b,b,d\}$ including two repeated distances (say, $a,b$) among $a,b,c$.
We can form two isosceles triangles with sides $a,a,c$ and $b,b,c$, which can be rotated in $\R^3$ around their common side of the length $c$, but their positions are fixed under isometry in $\R^3$ by the distance $d$ between their non-shared vertices.
\smallskip

\emph{Non-isosceles case}.
In the second case of Lemma~\ref{lem:PDD_symmetry} in Fig.~\ref{fig:symmetric+special_clouds}~(middle), $\PDD(A;3)$ has two pairs of equal rows of (unordered) distances $\{a,b,c\}$ and $\{a,b,d\}$.
Each of these triples uniquely determines a pair of equal triangles with a common side that are symmetric in the perpendicular bisector to this side.
For example, if we start with a fixed position of $[p_3,p_4]$ in $\R^3$, the union of equal triangles $\triangle p_1 p_3 p_4=\triangle p_2 p_3 p_4$ in Fig.~\ref{fig:symmetric+special_clouds}~(middle) is uniquely determined under isometry by the length $d$ of $[p_1,p_2]$.   
In $\R^2$, the parallelogram and isosceles trapezoid are distinguished by this distance 
$d$. 
\smallskip

Now we can assume that all rows of $\PDD(A;3)$ are different.
Then all points can be uniquely labelled as $p_1,p_2,p_3,p_4$ according to the lexicographic order of rows. 
Our aim is to get $\PDD(\{p_2,p_3,p_4\};2)$, reconstruct $\triangle p_2 p_3 p_4$, and then uniquely add $p_1$.
\smallskip

\emph{Case of a row with 3 equal distances}.
Let $\PDD(A;3)$ have a row of (say) $p_1$ with 3 equal distances $a$.
After removing the row of $p_1$, the distance $a$ from the rows of $p_2,p_3,p_4$, we get $\PDD(\{p_2,p_3,p_4\};2)$.
This smaller $3\times 2$ matrix determines $\triangle p_2 p_3 p_4$, uniquely under isometry in $\R^3$.
For a fixed $\triangle p_2 p_3 p_4$, the position of $p_1$ in $\R^3$ is determined by its distance $a$ to $p_2,p_3,p_4$, uniquely under the mirror reflection relative to the plane of $\triangle p_2 p_3 p_4$.   
If $n=2$, then $p_1$ is the unique circumcenter of $\triangle p_2 p_3 p_4$.
\smallskip

\emph{Case of a row with 3 unique distances}.
Let $\PDD(C;3)$ have a row of (say) $p_1$, where each of the distances $a,b,c$ (say, to $p_2$, $p_3$, $p_4$) appears in at most one other row (then $a,b,c$ are distinct).
After removing the row of $p_1$, the distance $a$ from the row $p_2$, the distance $b$ from the row of $p_3$, and the distance $c$ from the row of $p_4$, we get $\PDD(\{p_2,p_3,p_4\};2)$.
This $3\times 2$ matrix determines $\triangle p_2 p_3 p_4$, uniquely under isometry in $\R^3$.
Then the position of $p_1$ in $\R^3$ is determined by its distances $a,b,c$ to $p_2,p_3,p_4$, respectively, under a mirror reflection relative to the plane of the triangle $\triangle p_2 p_3 p_4$.
\smallskip

\emph{Case of one distance in 4 rows}.
Then two pairs of points have disjoint edges of the same length, e.g. $|p_1-p_2|=a=|p_3-p_4|$, so 
$\PDD(A;3)=\left(\begin{array}{ccc} a & b & c \\ a & d & e \\ a & b & d \\ a & c & e \end{array} \right)$ for $b=|p_1-p_3|$, $c=|p_1-p_4|$, $d=|p_2-p_3|$, $e=|p_2-p_4|$.
Then $c\neq d$ and $b\neq e$, else $\PDD(A;3)$ has two equal rows (considered above), similarly when $b=c$ and $d=e$.
\sskip

If $a$ equals one of $b,c,d,e$ (say, $e$), then $A$ is a 3-chain-equal cloud in Fig.~\ref{fig:symmetric+special_clouds}~(bottom right) and the argument below still works.
If $b\neq c$, we remove the row of $p_1$, the distance $b$ from the only row of $p_3$ containing $b$, the distance $c$ from the only row of $p_4$ containing $c$, and then remove $a$ from the remaining row of $p_2$.
This reduction to $\PDD(\{p_2,p_3,p_4\};2)$ allows us to reconstruct $A$, uniquely under isometry in $\R^3$ as in the case of a row with 3 unique distances.
If $b=c$ but $d\neq e$, we remove the row of $p_2$, the distance $d$ from the only row of $p_3$ containing $d$, the distance $e$ from the only row of $p_4$ containing $e$, and then remove $a$ from the remaining row of $p_1$, which allows us to uniquely reconstruct $A$ as in the case of a row with 3 unique distances above.
\smallskip

\emph{The final case}: no distance appears in all 4 distinct rows, but every row has a distance appearing in 3 rows, hence at least four times, including two times in the same row.
Then $A$ is a trisosceles cloud in Fig.~\ref{fig:symmetric+special_clouds}~(top right).
If any of the remaining distances $a,b,c$ are equal, $\PDD(A;3)$ has two equal rows (the case considered above).
Then we remove any row (say $a,b,b$) with two repeated distances, the distance $b$ from the only two rows containing $b$, and the distance $a$ from the remaining row.
\sskip

This reduction to $\PDD(\{p_2,p_3,p_4\};2)$, allows us to reconstruct $\triangle p_2 p_3 p_4$, uniquely under isometry in $\R^3$.
Though $p_1$ has equal distances to two of the vertices (say $p_2,p_3$), the ambiguity of reconstructing $p_1$ in $\R^3$ by its distances to $p_2,p_3,p_4$, is only under the mirror reflections relative to the bisector plane of $[p_2,p_3]$ and the plane of $\triangle p_2 p_3 p_4$.
\end{proof}

\begin{figure}[h!]
\centering
\includegraphics[width=\textwidth]{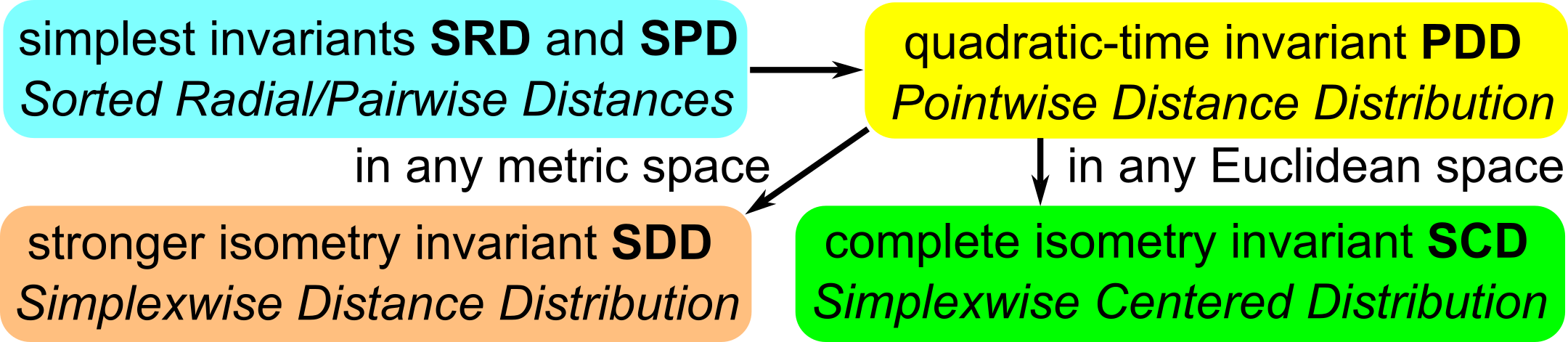}
\caption{ 
A hierarchy of invariants from the fastest (linear-time $\SRD$ and quadratic-time $\SPD$) to the stronger $\PDD$ and $\SDD$ in Chapter~\ref{chap:SDD} up to the Simplexwise Centered Distribution ($\SCD$)  in Chapter~\ref{chap:SCD}, which will satisfy the majority of conditions in Geo-Mapping Problem~\ref{pro:geocodes} for any finite $n$-dimensional clouds of unordered points under rigid motion in Theorems~\ref{thm:SCD_complete}, \ref{thm:SCD_metrics}, and \ref{thm:OSD+SCD_continuous}. 
}
\label{fig:SRD+SPD+PDD+SDD+SCD}
\end{figure}

Chapters~\ref{chap:SDD} and~\ref{chap:SCD} will extend the $\PDD$ to the stronger invariants in a metric space ($\SDD$) and complete invariant ($\SCD$) under rigid motion in any $\R^n$, as shown in Fig.~\ref{fig:SRD+SPD+PDD+SDD+SCD}.

\bibliographystyle{plain}
\bibliography{Geometric-Data-Science-book}

%
%
%

\chapter{Higher order distance distributions of unordered points in a metric space}
\label{chap:SDD} 

\abstract{
This chapter presents further advances towards a solution of the geo-mapping problem under isometry in any metric space, as stated in the previous chapter.
The Pointwise Distance Distribution (PDD) will be extended to stronger isometry invariant by collection distances to $h$-point subsets.
The resulting Simplexwise Distance Distribution ($\SDD$) is Lipschitz continuous and computable in a polynomial time of the number of points, for a fixed order $h$.
For $h=2$, the $\SDD$ distinguishes all (infinitely many) known counter-examples to the completeness of the PDD under isometry in $\R^3$.
}

\section{Simplexwise Distance Distributions of a cloud in a metric space}
\label{sec:SDD}

This chapter follows paper \cite{widdowson2023recognizing} and its extension \cite{kurlin2023simplexwise} to metric spaces with measures.
\myskip

We continue solving Problem~\ref{pro:metric_space} to find geocodes of finite clouds in any metric space.
The first section extends the Pointwise Distance Distribution (PDD) from Definition~\ref{dfn:PDD_finite} to a stronger invariant, which requires a few auxiliary definitions.
\myskip

The key idea of a stronger invariant is to use a base sequence of $h>1$ ordered points instead of $h=1$ point in the $\PDD$. 
\myskip

The \emph{lexicographic} order $\vec u<\vec v$ on vectors $\vec u=(u_1,\dots,u_h)$ and $\vec v=(v_1,\dots,v_h)$ means that if the first $i$ coordinates (where $i$ might be $0$) of $u,v$ coincide, then $u_{i+1}<v_{i+1}$.
Let $S_h$ denote the permutation group on indices $1,\dots,h$. 

\index{Relative Distance Distribution}
\begin{dfn}[Relative Distance Distribution $\RDD(C;A)$]
\label{dfn:RDD}
Let $C$ be a cloud of $m$ unlabelled points in a space with a metric $d$.
A \emph{base} sequence $A=(p_1,\dots,p_h)\in C^h$ consists of $1\leq h<m$ distinct points.
Let $D(A)$ be the \emph{triangular distance} matrix whose entry $D(A)_{i,j-1}$ is $d(p_i,p_j)$ for $1\leq i<j\leq h$, all other entries are zeros.
\sskip

Any permutation $\xi\in S_h$ acts on $D(A)$ by mapping $D(A)_{ij}$ to $D(A)_{kl}$, where $k\leq l$ is the pair of indices $\xi(i),\xi(j)-1$ written in increasing order.
For any other point $q\in C-A$, write distances from $q$ to $p_1,\dots,p_h$ as a column.
The $h\times (m-h)$-matrix $R(C;A)$ is formed by these $m-h$ lexicographically ordered columns.
The action of $\xi$ on $R(C;A)$ maps any $i$-th row to the $\xi(i)$-th row, after which all columns can be written in the lexicographic order.
The \emph{Relative Distance Distribution} $\RDD(C;A)$ is the equivalence class of the pair $[D(A),R(C;A)]$ of matrices under  permutations $\xi\in S_h$.
\edfn
\end{dfn}

For $h=1$ and a base sequence $A=(p_1)$, the matrix $D(A)$ is empty and $R(C;A)$ is a single row of distances (in the increasing order) from $p_1$ to all other points $q\in C$.
For $h=2$ and a base sequence $A=(p_1,p_2)$, the matrix $D(A)$ is the single number $d(p_1,p_2)$ and $R(C;A)$ consists of two rows of distances from $p_1,p_2$ to all other $q\in C$.

\begin{figure}[h!]
\centering
\includegraphics[width=\linewidth]{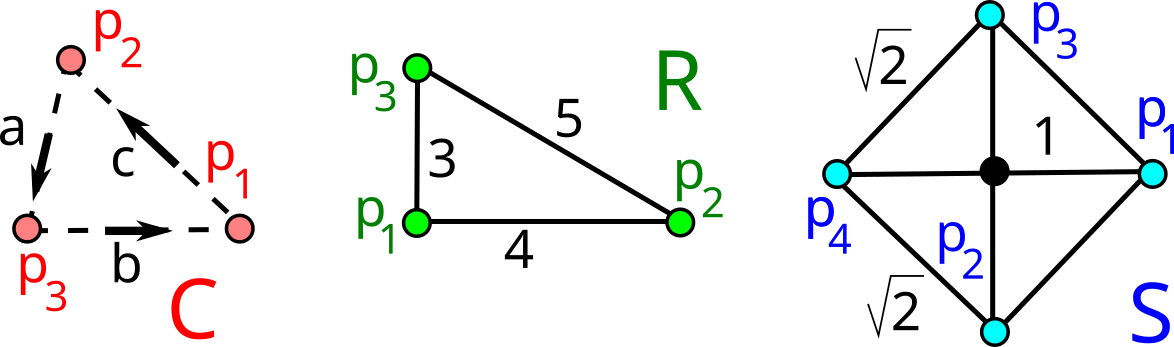}
\caption{\textbf{Left}: triangular cloud $C$ of points $p_1,p_2,p_3$ with inter-point distances $a\leq b\leq c$.
\textbf{Middle}: right-angled cloud $R$ of points $(0,0)$, $(4,0)$, $(0,3)$. 
\textbf{Right}: square cloud $S$ of points $(1,0)$, $(-1,0)$, $(0,1)$, $(-1,0)$. 
}
\label{fig:triangular_clouds}
\end{figure}

\begin{exa}[$\RDD$ for a 3-point cloud $C$]
\label{exa:RDD}
Let $C\subset\R^2$ consist of $p_1,p_2,p_3$ with inter-point distances $a\leq b\leq c$ ordered counter-clockwise as in Fig.~\ref{fig:triangular_clouds}~(left).  
Then
$$\RDD(C;p_1)=[\es;(b,c)], \qquad
\RDD(C;\vect{p_2}{p_3})=[a;\vect{c}{b}],$$
$$\RDD(C;p_2)=[\es;(a,c)], \qquad
\RDD(C;\vect{p_3}{p_1})=[b;\vect{a}{c}],$$
$$\RDD(C;p_3)=[\es;(a,b)], \qquad
\RDD(C;\vect{p_1}{p_2})=[c;\vect{b}{a}].$$
We have written $\RDD(C;A)$ for a base sequence $A=(p_i,p_j)$ of ordered points represented by a column.
Swapping the points $p_1\lra p_2$ makes the last $\RDD$ above equivalent to 
$\RDD\big(C;\vect{p_2}{p_1}\big)=\big[c;\vect{a}{b}\big]$.
\eexa
\end{exa}

Though $\RDD(C;A)$ is defined up to a permutation $\xi\in S_h$ of $h$ points in $A\in C^h$, comparisons of $\RDD$s will be practical for $h=2,3$ with metrics independent of $\xi$. 

\index{Simplexwise Distance Distribution}
\begin{dfn}[Simplexwise Distance Distribution $\SDD(C;h)$]
\label{dfn:SDD}
Let $C$ be a cloud of $m$ unlabelled points in a metric space.
For an integer $1\leq h<m$, the \emph{Simplexwise Distance Distribution} $\SDD(C;h)$ of order $h$ is the unordered set of $\RDD(C;A)$ for all unordered $h$-point subsets $A\subset C$.
\edfn
\end{dfn}

For order $h=1$ and any $m$-point cloud $C$, the distribution $\SDD(C;1)$ can be considered as a matrix of $m$ rows of ordered distances from every point $p\in C$ to all other $m-1$ points.
If we lexicographically order these $m$ rows and collapse any $l>1$ identical rows into a single one with the weight $l/m$, then we get the
Pointwise Distance Distribution $\PDD(C;m-1)$ introduced in Definition~\ref{dfn:PDD_finite}.

\index{moments of a weighted distribution}
\begin{dfn}[moments of a weighted distribution]
\label{dfn:moments}
Let $A$ be any unordered set of real numbers $a_1,\dots,a_m$ with weights $w_1,\dots,w_m$, respectively, such that $\sum\limits_{i=1}^m w_i=1$.
The 1st moment (\emph{average}) is the $\mu_1(A)=\sum\limits_{i=1}^m w_i a_i$.
The 2nd moment 
is $\mu_2(A)=\sqrt{\dfrac{1}{m}\sum\limits_{i=1}^m w_i a_i^2}$.
For $t\geq 3$, the $t$-th \emph{moment}  is $\sqrt[t]{m^{1-t}\sum\limits_{i=1}^m w_i a_i^t}$, see \cite[section~2.7]{keeping1995introduction}.
\edfn
\end{dfn}

The vector $\SPD(A)$ of Sorted Pairwise Distances was introduced in Definition~\ref{dfn:SPD}.  

\index{Simplexwise Distance Moments}
\begin{dfn}[Simplexwise Distance Moments $\SDM$]
\label{dfn:SDM}
For any $m$-point cloud $C$ in a metric space,
 let $A\subset C$ be a subset of $h$ unordered points.
The vector $\vec R(C;A)\in\R^{m-h}$ is obtained from the $h\times(m-h)$ matrix $R(C;A)$ in Definition~\ref{dfn:RDD} by writing the vector of $m-h$ column averages in increasing order.
\sskip

The pair $[\SPD(A);\vec R(C;A)]$ is the \emph{Average Distance Distribution} $\ADD(C;A)$ considered a vector of length $\frac{h(h-3)}{2}+m$.
The unordered collection of $\ADD(C;A)$ for all $\binom{m}{h}$ unordered subsets $A\subset C$ is the \emph{Average Simplexwise Distribution} $\ASD(C;h)$.
\sskip

The \emph{Simplexwise Distance Moment} $\SDM(C;h,t)$ is the $t$-th moment of $\ASD(C;h)$ considered a probability distribution of $\binom{m}{h}$ vectors, separately for each coordinate.
\edfn
\end{dfn}

\begin{figure}[h!]
\centering
\includegraphics[width=\linewidth]{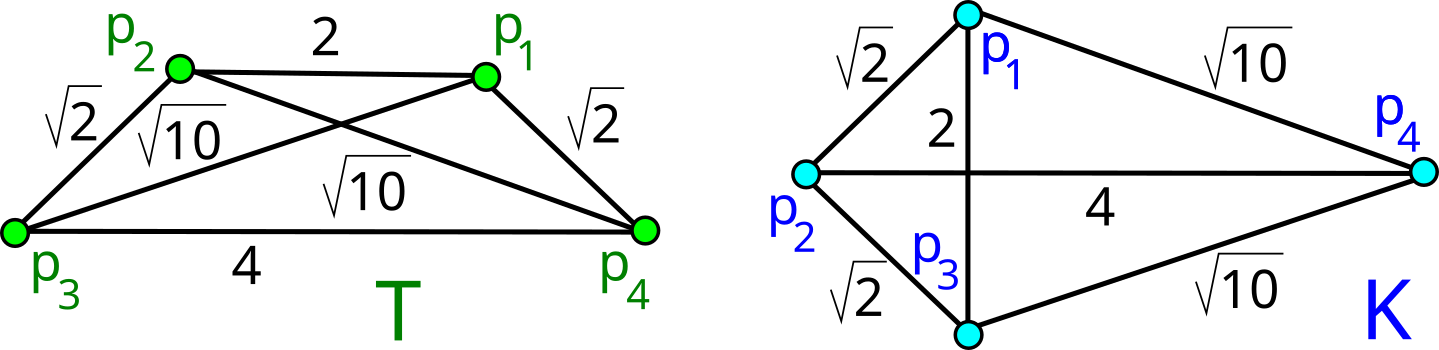}
\caption{\textbf{Left}: trapezium cloud $T$ of points $(1,1)$, $(-1,1)$, $(-2,0)$, $(2,0)$. 
\textbf{Right}: kite cloud $K$ of points $(0,1)$, $(-1,0)$, $(0,-1)$, $(3,0)$. 
}
\label{fig:quad_clouds}
\end{figure}

\begin{exa}[$\SDD$ and $\SDM$ for the 4-point clouds $T,K$]
\label{exa:SDD}
Fig.~\ref{fig:quad_clouds} shows the non-isometric 4-point clouds $T,K$ with the same Ordered Pairwise Distances: $\SPD=\{\sqrt{2},\sqrt{2},2,\sqrt{10},\sqrt{10},4\}$, see infinitely many examples in \cite{boutin2004reconstructing}.
The arrows on the edges of $T,K$ show orders of points in each pair of vertices for $\RDD$s.
Then $T,K$ are distinguished under isometry by $\SDD(T;2)\neq \SDD(K;2)$ in Table~\ref{tab:SDD+TK}.
The 1st coordinate of $\SDM(C;2,1)\in\R^3$ is the average of the six distances from $\SPD$ (the same for $T,K$) but the other two coordinates (column averages from $R(C;A)$ matrices) differ.
\eexa
\end{exa}

\begin{table}
  \centering
  \begin{tabular}{@{}l|l@{}}
   	$\RDD(T;A)$ in $\SDD(T;2)$ & $\RDD(K;A)$ in $\SDD(K;2)$ \\
    
    	$[\sqrt{2},\matv{2}{\cbox{yellow}{\sqrt{10}}}{\sqrt{10}}{4}]\times 2$ &
	$[\sqrt{2},\matv{2}{\cbox{yellow}{\sqrt{2}}}{\sqrt{10}}{4}]\times 2$ 
	\smallskip \\
	
$[2,\matv{\sqrt{2}}{\cbox{yellow}{\sqrt{10}}}{\sqrt{10}}{\cbox{yellow}{\sqrt{2}}}]$ & 
$[2,\matv{\sqrt{2}}{\cbox{yellow}{\sqrt{2}}}{\sqrt{10}}{\cbox{yellow}{\sqrt{10}}}]$ 
\smallskip \\
    
	$[\sqrt{10},\matv{\sqrt{2}}{\cbox{yellow}{2}}{\cbox{yellow}{4}}{\cbox{yellow}{\sqrt{2}}}]\times 2$ &
	$[\sqrt{10},\matv{\sqrt{2}}{\cbox{yellow}{4}}{\cbox{yellow}{2}}{\cbox{yellow}{\sqrt{10}}}]\times 2$ \smallskip \\
	
	    $[4,\matv{\sqrt{2}}{\sqrt{10}}{\cbox{yellow}{\sqrt{10}}}{\cbox{yellow}{\sqrt{2}}}]$ &
	$[4,\matv{\sqrt{2}}{\sqrt{10}}{\cbox{yellow}{\sqrt{2}}}{\cbox{yellow}{\sqrt{10}}}]$  \\

    \hline
  	$\ADD(T;A)$ in $\ASD(T;2)$ & $\ADD(K;A)$ in $\ASD(K;2)$ \\
    $[\sqrt{2},(\frac{2+\cbox{yellow}{\sqrt{10}}}{2},\frac{4+\sqrt{10}}{2})]\times 2$ &
    $[\sqrt{2},(\frac{2+\cbox{yellow}{\sqrt{2}}}{2},\frac{4+\sqrt{10}}{2})]\times 2$ \\

    $[2,(\cbox{yellow}{\frac{\sqrt{2}+\sqrt{10}}{2},\frac{\sqrt{2}+\sqrt{10}}{2}})]$ &
    $[2,(\cbox{yellow}{\sqrt{2},\sqrt{10}})]$ \\

    $[\sqrt{10},(\frac{2+\cbox{yellow}{\sqrt{2}}}{2},\frac{4+\sqrt{2}}{2})]\times 2$ &
    $[\sqrt{10},(\frac{2+\cbox{yellow}{\sqrt{10}}}{2},\frac{4+\sqrt{2}}{2})]\times 2$ \\
    
       $[4,(\frac{\sqrt{2}+\sqrt{10}}{2},\frac{\sqrt{2}+\sqrt{10}}{2})]$ & $[4,(\frac{\sqrt{2}+\sqrt{10}}{2},\frac{\sqrt{2}+\sqrt{10}}{2})]$     \\

    \hline
   	$\SDM_1=\dfrac{3+\sqrt{2}+\sqrt{10}}{3}$ & 
   	$\SDM_1=\dfrac{3+\sqrt{2}+\sqrt{10}}{3}$ \\

   	$\SDM_2=\dfrac{\cbox{yellow}{6+2\sqrt{2}+4\sqrt{10}}}{12}$ & 
   	$\SDM_2=\dfrac{\cbox{yellow}{8+5\sqrt{2}+3\sqrt{10}}}{12}$ \\

   	$\SDM_3=\frac{16+\cbox{yellow}{4\sqrt{2}+4\sqrt{10}}}{12}$ & 
   	$\SDM_3=\frac{16+\cbox{yellow}{3\sqrt{2}+5\sqrt{10}}}{12}$
  \end{tabular}
  \caption{\textbf{Top}: Relative Distance Distributions from Definition~\ref{dfn:RDD} for all 6 base sequences $A$ in the 4-point clouds $T,K$ in Fig.~\ref{fig:quad_clouds}. 
The symbol $\times 2$ indicates a doubled $\RDD$.
The three bottom rows show coordinates of $\SDM(C;2,1)\in\R^3$ from Definition~\ref{dfn:SDM} for $h=2$, $t=1$ and $C=T,K$.
Different elements are \hl{highlighted} and imply that all invariants $\SDD,\ADD,\SDM$ distinguish $T\not\simeq K$.
}
\label{tab:SDD+TK}
\end{table}

Some of the $\binom{m}{h}$ $\RDD$s in $\SDD(C;h)$ can concide as in Example~\ref{exa:SDD}.
If we collapse any $l>1$ identical $\RDD$s into a single $\RDD$ with the \emph{weight} $l/\binom{m}{h}$, $\SDD$ can be considered as a weighted probability distribution 
of $\RDD$s.
\smallskip

In a general metric space, a point cloud $C$ is usually given by a distance matrix on (arbitrarily ordered) points of $C$. 
Hence, we assume that the distance between any points of $C$ is accessible in a constant time.

\index{polynomial-time complexity}
\index{Simplexwise Distance Distribution}
\begin{thm}[invariance and time of $\SDD$, {\cite[Theorem~3.6]{kurlin2023simplexwise}}]
\label{thm:SDD_time}
For any order $h\geq 1$ and any cloud $C$ of $m$ unlabelled points in a metric space, $\SDD(C;h)$ is an isometry invariant, which can be computed in time $O(m^{h+1}/(h-1)!)$.
For any $t\geq 1$, the invariant $\SDM(C;h,t)\in\R^{m+\frac{h(h-3)}{2}}$ has the same asymptotic time.
\ethm
\end{thm}

\section{The expressiveness of Simplexwise Distance Distributions}
\label{sec:SDD-examples}


This section shows that $\SDD(C;2)$ distinguishes all infinitely many known pairs \cite[Fig.~S4]{pozdnyakov2020incompleteness} of non-isometric clouds $S,Q\subset\R^3$ that have equal $\PDD(S)=\PDD(Q)$ 

Examples~\ref{exa:5-point_sets} and~\ref{exa:7-point_sets} distinguish clouds of 5 points and 7 points, respectively, in $\R^3$ by comparing their $\SDD$s of order 2.
In Example~\ref{exa:6-point_sets}, the invariant $\SDD(C;2)$ distinguishes 6-point clouds in a family of pairs depending on three parameters.  


\begin{exa}[5-point clouds]
\label{exa:5-point_sets}
Fig.~\ref{fig:5-point_sets} shows the 5-point clouds $S_{\pm}\subset\R^3$ taken from \cite[Figure S4(A)]{pozdnyakov2020incompleteness}.
The clouds $S_{\pm}$ are not isometric, because $S_+$ has the triple of points $ B_+, G_+, R_+$ with pairwise distances $\sqrt{2},\sqrt{6},\sqrt{6}$, but $S_-$ has no such a triple.

\begin{figure}[h!]
\centering
\includegraphics[width=\linewidth]{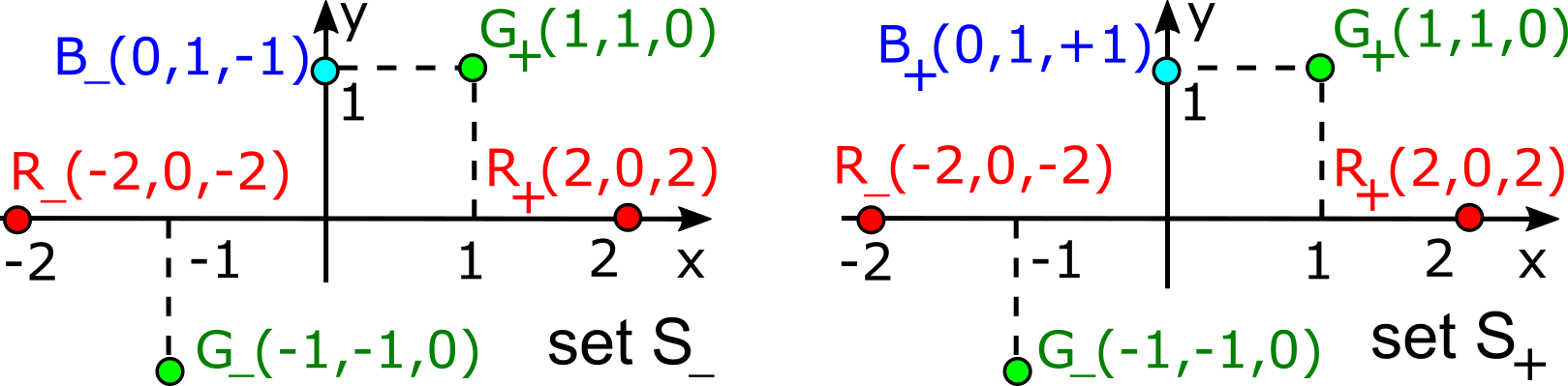}
\caption{
See Example~\ref{exa:5-point_sets}
\textbf{Left}: $(x,y)$-projection of the 5-point cloud $S_-\subset\R^3$ consisting of the green points $G_-=(-1,-1,0)$ and $G_{+}=(1,1,0)$, the red points $R_-=(-2,0,-2)$ and $R_+=(2,0,2)$, and the blue point $B_-=(0,1,-1)$.
\textbf{Right}: to get $S_+\subset\R^3$ from the cloud $S_-$, replace the point $B_-$ with another point $B_+=(0,1,1)$.}
\label{fig:5-point_sets}
\end{figure}

Table~\ref{tab:distances_S-+} \hl{highlights} differences between distance matrices.
If we order distances to neighbours, the matrices in Table~\ref{tab:ordered_distances_S-+} differ only in one pair.

\begin{table}[h!]
\begin{tabular}{@{}l|ccccc@{}}
 \hline
distances of \hl{$S_-$} & $R_-$ & $R_+$ & $G_-$ & $G_+$ & \hl{$B_-$} \\
\hline
$R_-(-2,0,-2)$ & $0$ & $\sqrt{32}$ & $\sqrt{6}$ & $\sqrt{14}$ & \hl{$\sqrt{6}$} \smallskip \\
$R_+(+2,0,+2)$ & $\sqrt{32}$ & 0 & $\sqrt{14}$ & $\sqrt{6}$ & \hl{$\sqrt{14}$} \smallskip \\
$G_-(-1,-1,0)$ & $\sqrt{6}$ & $\sqrt{14}$ & 0 & $\sqrt{8}$ & $\sqrt{6}$ \smallskip \\
$G_+(+1,+1,0)$ & $\sqrt{14}$ & $\sqrt{6}$ & $\sqrt{8}$ & 0 & $\sqrt{2}$ \smallskip \\
\hl{$B_-(0,+1,-1)$} & \hl{$\sqrt{6}$} & \hl{$\sqrt{14}$} & $\sqrt{6}$ & $\sqrt{2}$ & 0 \\
\hline
\end{tabular}
\hspace*{5mm}
\begin{tabular}{@{}l|ccccc@{}}
 \hline
distances of \hl{$S_+$} & $R_-$ & $R_+$ & $G_-$ & $G_+$ & \hl{$B_+$} \\
\hline
$R_-(-2,0,-2)$ & $0$ & $\sqrt{32}$ & $\sqrt{6}$ & $\sqrt{14}$ & \hl{$\sqrt{14}$} \smallskip \\
$R_+(+2,0,+2)$ & $\sqrt{32}$ & 0 & $\sqrt{14}$ & $\sqrt{6}$ & \hl{$\sqrt{6}$} \smallskip \\
$G_-(-1,-1,0)$ & $\sqrt{6}$ & $\sqrt{14}$ & 0 & $\sqrt{8}$ & $\sqrt{6}$ \smallskip \\
$G_+(+1,+1,0)$ & $\sqrt{14}$ & $\sqrt{6}$ & $\sqrt{8}$ & 0 & $\sqrt{2}$ \smallskip \\
\hl{$B_+(0,+1,+1)$} & \hl{$\sqrt{14}$} & \hl{$\sqrt{6}$} & $\sqrt{6}$ & $\sqrt{2}$ & 0 \\
\hline
\end{tabular}
\caption{Distance matrices of the clouds $S_{\mp}\subset\R^3$ in Fig.~\ref{fig:5-point_sets}.}
\label{tab:distances_S-+}
\end{table}

\begin{table}[h]
\centering
\begin{tabular}{@{}l|cccc@{}}
\hline
\hl{$S_-$} distances to & 1st neighbour & 2nd neighbour & 3rd neighbour & 4th neighbour \\
\hline
$R_-=(-2,0,-2)$ & $\sqrt{6}$ & \hl{$\sqrt{6}$}  & $\sqrt{14}$ & $\sqrt{32}$ \smallskip \\
$R_+=(+2,0,+2)$ & $\sqrt{6}$ & \hl{$\sqrt{14}$}  & $\sqrt{14}$ & $\sqrt{32}$ \smallskip \\
$G_-=(-1,-1,0)$  & $\sqrt{6}$ & $\sqrt{6}$ & $\sqrt{8}$ & $\sqrt{14}$  \smallskip \\
$G_+=(+1,+1,0)$ & $\sqrt{2}$ & $\sqrt{6}$ & $\sqrt{8}$ & $\sqrt{14}$ \smallskip \\
$B_-=(0,+1,-1)$ & $\sqrt{2}$ & $\sqrt{6}$ & $\sqrt{6}$ & $\sqrt{14}$ \\
\hline
\end{tabular}
\myskip

\begin{tabular}{@{}l|cccc@{}}
\hline
\hl{$S_+$} distances to & 1st neighbour & 2nd neighbour & 3rd neighbour & 4th neighbour \\
\hline
$R_-=(-2,0,-2)$ & $\sqrt{6}$ & \hl{$\sqrt{14}$}  & $\sqrt{14}$ & $\sqrt{32}$ \smallskip \\
$R_+=(+2,0,+2)$ & $\sqrt{6}$ & \hl{$\sqrt{6}$}  & $\sqrt{14}$ & $\sqrt{32}$ \smallskip \\
$G_-=(-1,-1,0)$  & $\sqrt{6}$ & $\sqrt{6}$ & $\sqrt{8}$ & $\sqrt{14}$  \smallskip \\
$G_+=(+1,+1,0)$ & $\sqrt{2}$ & $\sqrt{6}$ & $\sqrt{8}$ & $\sqrt{14}$ \smallskip \\
$B_+=(0,+1,-1)$ & $\sqrt{2}$ & $\sqrt{6}$ & $\sqrt{6}$ & $\sqrt{14}$ \\
\hline
\end{tabular}
\caption{For each point from the 5-point cloud $S_{+}$ in Fig.~\ref{fig:5-point_sets}, the distances to neighbours from Table~\ref{tab:distances_S-+} are ordered in each row.}
\label{tab:ordered_distances_S-+}
\end{table}

If we ignore the labels of all points in columns,
Table~\ref{tab:ordered_distances_S-+} implies that $S_{\pm}$ have identical Pointwise Distance Distribution ($\PDD$). 
For easier visualisation, the matrix below is obtained by lexicographically ordering the rows in Table~\ref{tab:ordered_distances_S-+}:
$$\PDD(S_{\pm})=\SDD(S_{\pm};1)=\left( \begin{array}{ccccc}
\sqrt{2} & \sqrt{6} & \sqrt{6} & \sqrt{14} \\
\sqrt{2} & \sqrt{6} & \sqrt{8} & \sqrt{14} \\
\sqrt{6} & \sqrt{6} & \sqrt{8} & \sqrt{14} \\
\sqrt{6} & \sqrt{6}  & \sqrt{14} & \sqrt{32}\\
\sqrt{6} & \sqrt{14}  & \sqrt{14} & \sqrt{32} 
\end{array} \right).$$
Now we show that $\SDD(S_-;2)\neq \SDD(S_+;2)$.
For $h=2$, the Simplexwise Distance Distribution $\SDD(C;h)$ consists of $\RDD(C;A)$ for 2-point subsets $A\subset C$.
Both sets $S_{\pm}$ have a single pair of points $(G_+,B-)$ and $(G+,B_+)$ at distance $\sqrt{2}$. 
Hence it suffices to show that the Relative Distance Distributions differ for this pair:
$$\RDD\left(S_-,\vect{G_+}{B_-}\right)=\left[\sqrt{2},
\left( \begin{array}{ccc}
\sqrt{8} & \sqrt{14} & \sqrt{6} \\
\sqrt{6} & \cbox{yellow}{\sqrt{6}} & \cbox{yellow}{\sqrt{14}} \\
G_- & R_- & R_+ \\
\end{array} \right) \right],$$

$$\RDD\left(S_+,\vect{G_+}{B_+}\right)=\left[\sqrt{2},
\left( \begin{array}{ccc}
\sqrt{8} & \sqrt{14} & \sqrt{6} \\
\sqrt{6} & \cbox{yellow}{\sqrt{14}} & \cbox{yellow}{\sqrt{6}} \\
G_- & R_- & R_+ \\
\end{array} \right) \right].$$

The last rows in the above $3\times 3$ matrices indicate a complementary point $q\in C-A$ for indexing columns of the $2\times 3$ matrices $R(C;A)$ in Definition~\ref{dfn:RDD}. 
The resulting $\RDD$s differ because any permutation of rows or columns of $R(S_+;\{G_+,B_+\})$ keeps the pair $\sqrt{6},\sqrt{6}$ in the same column but $R(S_-;\{G_+,B_+\})$ has no pair $\sqrt{6},\sqrt{6}$ in one column.
Hence $\SDD(S_-;2)\neq\SDD(S_+;2)$. 
\eexa
\end{exa}


\begin{figure}[h!]
\centering
\includegraphics[width=\linewidth]{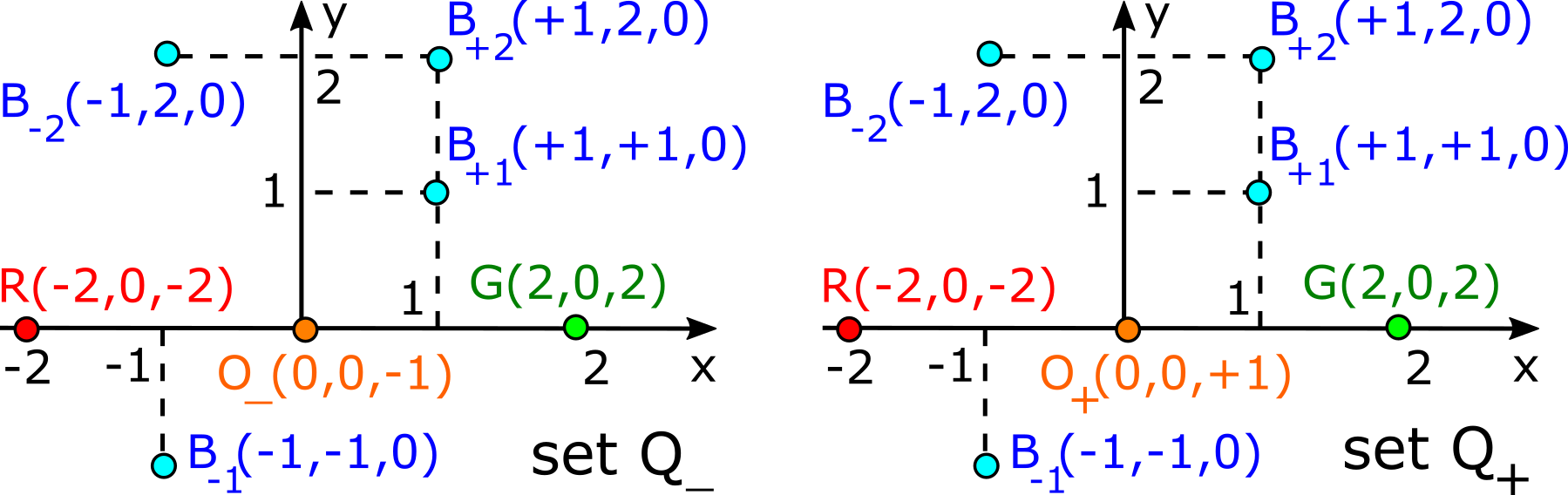}
\caption{
See Example~\ref{exa:7-point_sets}.
\textbf{Left}: $(x,y)$-projection of the 7-point cloud $Q_-\subset\R^3$, which consists of the red point $R=(-2,0,-2)$, green point $G=(2,0,2)$, four blue points $B_{\pm 1}=(\pm 1,\pm 1,0)$, $B_{\pm 2}=(\pm 1,2,0)$, orange point $O_-=(0,0,-1)$.
\textbf{Right}: to get the cloud $Q_+$ from the cloud $Q_-\subset\R^3$, replace the point $O_-$ with $O_+=(0,0,+1)$.}
\label{fig:7-point_sets}
\end{figure}

\begin{exa}[7-point clouds]
\label{exa:7-point_sets}
The clouds $Q_{\pm}$ in Fig.~\ref{fig:7-point_sets} taken from \cite[Figure S4(B)]{pozdnyakov2020incompleteness} have distances in Table~\ref{tab:distances_Q-+}.
Both sets have only two pairs of points at distance $\sqrt{6}$.
Hence it suffices to compare $\RDD$s for these pairs below.

\begin{table}[h!]
\centering
\begin{tabular}{@{}l|ccccccc@{}}
\hline
distances of \hl{ $Q_-$ } & $R$ & $G$ & $B_{-1}$ & $B_{+1}$ & $B_{-2}$ & $B_{+2}$ & \hl{ $O_{-}$ } \\
\hline
$R=(-2,0,-2)$ & $0$ & $\sqrt{32}$ & $\sqrt{6}$ & $\sqrt{14}$ & $3$ & $\sqrt{17}$ & \hl{ $\sqrt{5}$ } \smallskip \\
$G=(+2,0,+2)$ & $\sqrt{32}$ & 0 & $\sqrt{14}$ & $\sqrt{6}$ & $\sqrt{17}$ & $3$ & \hl{ $\sqrt{13}$ } \smallskip \\
$B_{-1}=(-1,-1,0)$ & $\sqrt{6}$ & $\sqrt{14}$ & 0 & $\sqrt{8}$ & $3$ & $\sqrt{13}$ & $\sqrt{3}$ \smallskip \\
$B_{+1}=(+1,+1,0)$ & $\sqrt{14}$ & $\sqrt{6}$ & $\sqrt{8}$ & 0 & $\sqrt{5}$ & $1$ & $\sqrt{3}$ \smallskip \\
$B_{-2}=(-1,2,0)$ & $3$ & $\sqrt{17}$ & $3$ & $\sqrt{5}$ & 0 & 2 & $\sqrt{6}$ \smallskip \\
$B_{+2}=(+1,2,0)$ & $\sqrt{17}$ & $3$ & $\sqrt{13}$ & $1$ & 2 & 0 & $\sqrt{6}$ \smallskip \\
\hl{ $O_{-}=(0,0,-1)$ } & \hl{ $\sqrt{5}$ } & \hl{ $\sqrt{13}$ } & $\sqrt{3}$ & $\sqrt{3}$ & $\sqrt{6}$ & $\sqrt{6}$ & 0 \smallskip \\
\hline
\end{tabular}
\smallskip

\begin{tabular}{@{}l|ccccccc@{}}
\hline
distances of \hl{ $Q_+$ } & $R$ & $G$ & $B_{-1}$ & $B_{+1}$ & $B_{-2}$ & $B_{+2}$ & \hl{ $O_{+}$ } \\
\hline
$R=(-2,0,-2)$ & $0$ & $\sqrt{32}$ & $\sqrt{6}$ & $\sqrt{14}$ & $3$ & $\sqrt{17}$ & \hl{ $\sqrt{13}$ } \smallskip \\
$G=(+2,0,+2)$ & $\sqrt{32}$ & 0 & $\sqrt{14}$ & $\sqrt{6}$ & $\sqrt{17}$ & $3$ & \hl{ $\sqrt{5}$ } \smallskip \\
$B_{-1}=(-1,-1,0)$ & $\sqrt{6}$ & $\sqrt{14}$ & 0 & $\sqrt{8}$ & $3$ & $\sqrt{13}$ & $\sqrt{3}$ \smallskip \\
$B_{+1}=(+1,+1,0)$ & $\sqrt{14}$ & $\sqrt{6}$ & $\sqrt{8}$ & 0 & $\sqrt{5}$ & $1$ & $\sqrt{3}$ \smallskip \\
$B_{-2}=(-1,2,0)$ & $3$ & $\sqrt{17}$ & $3$ & $\sqrt{5}$ & 0 & 2 & $\sqrt{6}$ \smallskip \\
$B_{+2}=(+1,2,0)$ & $\sqrt{17}$ & $3$ & $\sqrt{13}$ & $1$ & 2 & 0 & $\sqrt{6}$ \smallskip \\
\hl{ $O_{+}=(0,0,+1)$ } & \hl{ $\sqrt{13}$ } & \hl{ $\sqrt{5}$ } & $\sqrt{3}$ & $\sqrt{3}$ & $\sqrt{6}$ & $\sqrt{6}$ & 0 \\
\hline
\end{tabular}
\caption{The distance matrices of the 7-point clouds $Q_{\mp}$ in Fig.~\ref{fig:7-point_sets} taken from \cite[Figure S4(B)]{pozdnyakov2020incompleteness}.}
\label{tab:distances_Q-+}
\end{table}

$$R\left(Q_-;\vect{G}{B_{+1}}\right)=
\left( \begin{array}{ccccc}
\sqrt{32} & \sqrt{14} & \sqrt{17} & 3 & \cbox{yellow}{\sqrt{13}}  \\
\sqrt{6} & \sqrt{8} & \sqrt{5} & 1 & \sqrt{3} \\
R & B_{-1} & B_{-2} & B_{+2} & O_-
\end{array} \right),$$
$$R\left(Q_-;\vect{R}{B_{-1}}\right)=
\left( \begin{array}{ccccc}
\sqrt{32} & \sqrt{14} & 3 & \sqrt{17} & \cbox{yellow}{\sqrt{5}}  \\
\sqrt{14} & \sqrt{8} & 3 & \sqrt{13} & \sqrt{3} \\
G & B_{+1} & B_{-2} & B_{+2} & O_-
\end{array} \right).$$

The pair above has submatrices $\mat{3}{\sqrt{13}}{1}{\sqrt{3}}$ and $\mat{3}{\sqrt{5}}{3}{\sqrt{3}}$ but the pair below has no such submatrices.

$$R\left(Q_+;\vect{G}{B_{+1}}\right)=
\left( \begin{array}{ccccc}
\sqrt{32} & \sqrt{14} & \sqrt{17} & 3 & \cbox{yellow}{\sqrt{5}}  \\
\sqrt{6} & \sqrt{8} & \sqrt{5} & 1 & \sqrt{3} \\
R & B_{-1} & B_{-2} & B_{+2} & O_+
\end{array} \right),$$
$$R\left(Q_+;\vect{R}{B_{-1}}\right)=
\left( \begin{array}{ccccc}
\sqrt{32} & \sqrt{14} & 3 & \sqrt{17} & \cbox{yellow}{\sqrt{13}}  \\
\sqrt{14} & \sqrt{8} & 3 & \sqrt{13} & \sqrt{3} \\
G & B_{+1} & B_{-2} & B_{+2} & O_+
\end{array} \right)$$

The pair of $\RDD\left(Q_-;\vect{G}{B_{+1}}\right)$ and $\RDD\left(Q_-;\vect{R}{B_{-1}}\right)$ differs from the pair $\RDD\left(Q_+;\vect{G}{B_{+1}}\right)$ and $\RDD\left(Q_+;\vect{R}{B_{-1}}\right)$.
Hence $\SDD(Q_-;2)\neq\SDD(Q_+;2)$. 
\eexa
\end{exa}


\begin{figure}[h!]
\centering
\includegraphics[width=\linewidth]{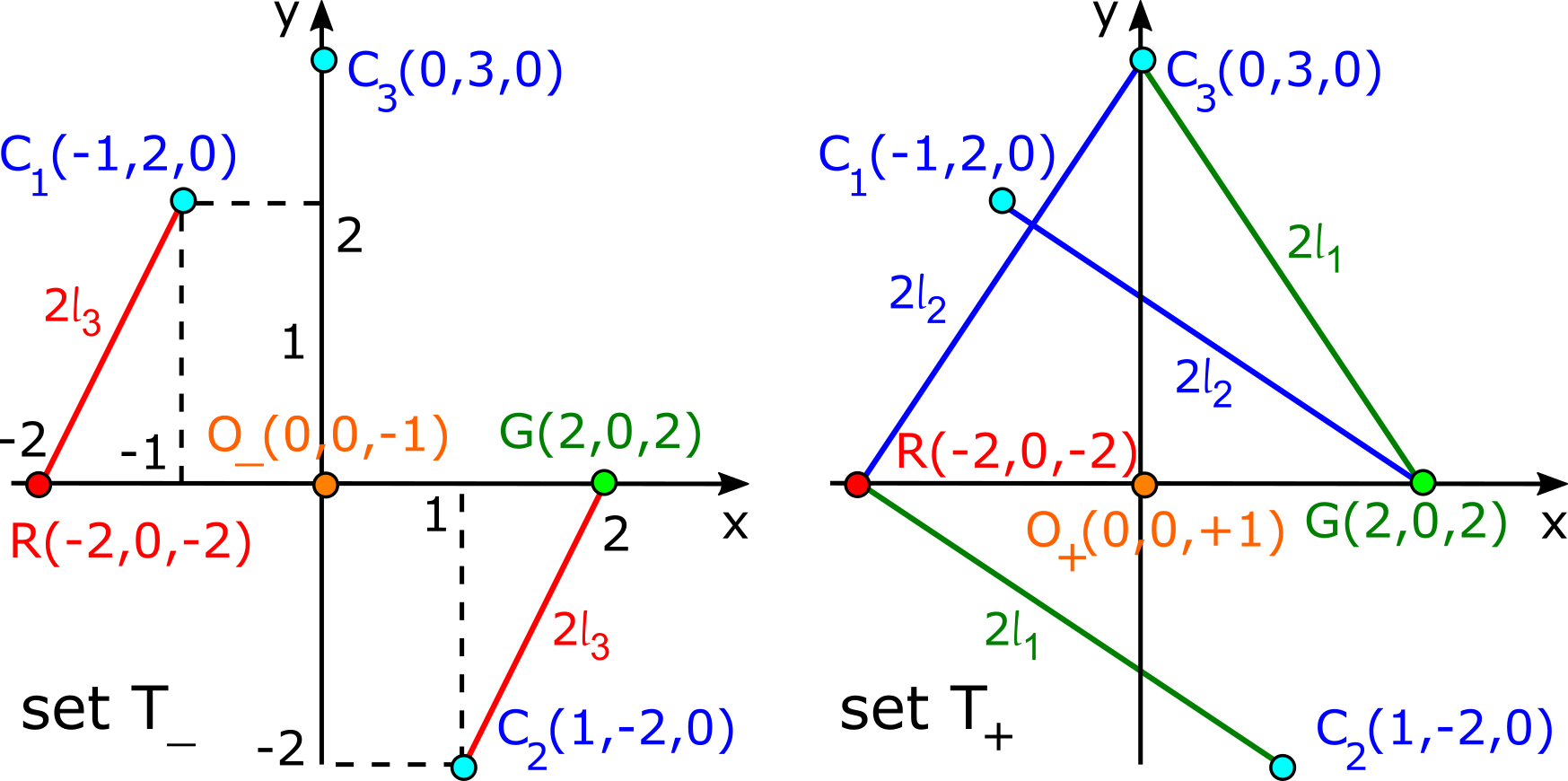}
\caption{
See Example~\ref{exa:6-point_sets}.
\textbf{Left}: $(x,y)$-projection of the 6-point cloud $T_-\subset\R^3$ consisting of the red point $R=(-2,0,-2)$, green point $G=(2,0,2)$, three blue points $C_1=(x_1,y_1,0)$, $C_2=(x_2,y_2,0)$,  $C_3=(x_3,y_3,0)$, and orange point $O_-=(0,0,-1)$ so that $|RC_1|=2l_3=|GC_2|$, $|RC_2|=2l_1=|GC_3|$, $|RC_3|=2l_2=|GC_1|$.
\textbf{Right}: to get the cloud $T_+\subset\R^3$ from the cloud $T_-$, replace the point $O_-$ with $O_+=(0,0,+1)$.}
\label{fig:6-point_sets}
\end{figure}

\begin{exa}[6-point clouds]
\label{exa:6-point_sets}
The clouds $T_{\pm}$ in Fig.~\ref{fig:6-point_sets}, which was motivated by \cite[Figure S4(C)]{pozdnyakov2020incompleteness}, have the points $R,G,O_{\pm}$ from the clouds $Q_{\pm}$ in Example~\ref{exa:7-point_sets} and three new points $C_1(x_1,y_1,0)$, $C_2(x_2,y_2,0)$, $C_3(x_3,y_3,0)$ such that $|RC_1|=|GC_2|$, $|RC_2|=|GC_3|$, $|RC_3|=|GC_1|$.
Denote by $2l_1,2l_2,2l_3$ the lengths of these three pairs of line segments after their projection to the $xy$-plane so that
$$\leqno{(\ref{exa:6-point_sets}.1)}
\left\{\begin{array}{l}
(x_2+2)^2+y_2^2=|RC_2|^2-4=(2l_1)^2, \\
(x_3-2)^2+y_3^2=|GC_3|^2-4=(2l_1)^2;
\end{array} \right.$$ 
$$\leqno{(\ref{exa:6-point_sets}.2)}
\left\{\begin{array}{l}
(x_3+2)^2+y_3^2=|RC_3|^2-4=(2l_2)^2, \\
(x_1-2)^2+y_1^2=|GC_1|^2-4=(2l_2)^2;
\end{array} \right.$$ 
$$\leqno{(\ref{exa:6-point_sets}.3)}
\left\{\begin{array}{l}
(x_1+2)^2+y_1^2=|RC_1|^2-4=(2l_3)^2, \\
(x_2-2)^2+y_2^2=|GC_2|^2-4=(2l_3)^2.
\end{array} \right.$$ 

\begin{table}[h!]
\centering
\begin{tabular}{@{}l|cccccc@{}}
\hline
distances of \hl{ $T_-$ } & $R$ & $G$ & $C_{1}$ & $C_{2}$ & $C_{3}$ & \hl{ $O_{-}$ } \\
\hline
$R=(-2,0,-2)$ & $0$ & $\sqrt{32}$ & $2\sqrt{l_3^2+1}$ & $2\sqrt{l_1^2+1}$ & $2\sqrt{l_2^2+1}$ & \hl{ $\sqrt{5}$ } \smallskip \\
$G=(+2,0,+2)$ & $\sqrt{32}$ & 0 & $2\sqrt{l_2^2+1}$ & $2\sqrt{l_3^2+1}$ & $2\sqrt{l_1^2+1}$  & \hl{ $\sqrt{13}$ } \smallskip \\
$C_{1}=(x_1,y_1,0)$ & $2\sqrt{l_3^2+1}$ & $2\sqrt{l_2^2+1}$ & 0 & $|C_1C_2|$ & $|C_3C_1|$ & $\sqrt{2l_2^2+2l_3^2-3}$ \smallskip \\
$C_{2}=(x_2,y_2,0)$ & $2\sqrt{l_1^2+1}$ & $2\sqrt{l_3^2+1}$ & $|C_1C_2|$ & 0 & $|C_2C_3|$ & $\sqrt{2l_3^2+2l_1^2-3}$ \smallskip \\
$C_{3}=(x_3,y_3,0)$ & $2\sqrt{l_2^2+1}$ & $2\sqrt{l_1^2+1}$ & $|C_3C_1|$ & $|C_2C_3|$ & 0 & $\sqrt{2l_1^2+2l_2^2-3}$ \smallskip  \\
\hl{ $O_{-}=(0,0,-1)$ } & \hl{ $\sqrt{5}$ } & \hl{ $\sqrt{13}$ } & $\sqrt{2l_2^2+2l_3^2-3}$ & $\sqrt{2l_3^2+2l_1^2-3}$ & $\sqrt{2l_1^2+2l_2^2-3}$ & 0 \smallskip \\
\hline
\end{tabular}
\smallskip

\begin{tabular}{@{}l|cccccc@{}}
\hline
distances of \hl{ $T_+$ } & $R$ & $G$ & $C_{1}$ & $C_{2}$ & $C_{3}$ & \hl{ $O_{+}$ } \\
\hline
$R=(-2,0,-2)$ & $0$ & $\sqrt{32}$ & $2\sqrt{l_3^2+1}$ & $2\sqrt{l_1^2+1}$ & $2\sqrt{l_2^2+1}$ & \hl{ $\sqrt{13}$ } \smallskip \\
$G=(+2,0,+2)$ & $\sqrt{32}$ & 0 & $2\sqrt{l_2^2+1}$ & $2\sqrt{l_3^2+1}$ & $2\sqrt{l_1^2+1}$  & \hl{ $\sqrt{5}$ } \smallskip \\
$C_{1}=(x_1,y_1,0)$ & $2\sqrt{l_3^2+1}$ & $2\sqrt{l_2^2+1}$ & 0 & $|C_1C_2|$ & $|C_3C_1|$ & $\sqrt{2l_2^2+2l_3^2-3}$ \smallskip \\
$C_{2}=(x_2,y_2,0)$ & $2\sqrt{l_1^2+1}$ & $2\sqrt{l_3^2+1}$ & $|C_1C_2|$ & 0 & $|C_2C_3|$ & $\sqrt{2l_3^2+2l_1^2-3}$ \smallskip \\
$C_{3}=(x_3,y_3,0)$ & $2\sqrt{l_2^2+1}$ & $2\sqrt{l_1^2+1}$ & $|C_3C_1|$ & $|C_2C_3|$ & 0 & $\sqrt{2l_1^2+2l_2^2-3}$  \smallskip \\
\hl{ $O_{+}=(0,0,+1)$ } & \hl{ $\sqrt{13}$ } & \hl{ $\sqrt{5}$ } & $\sqrt{2l_2^2+2l_3^2-3}$ & $\sqrt{2l_3^2+2l_1^2-3}$ & $\sqrt{2l_1^2+2l_2^2-3}$ & 0 \smallskip \\
\hline
\end{tabular}
\caption{The distance matrices of the 6-point clouds $T_{\mp}$ in Fig.~\ref{fig:6-point_sets} motivated by \cite[Figure S4(C)]{pozdnyakov2020incompleteness}.}
\label{tab:distances_T-+}
\end{table}

Comparing the 1st part of $(\ref{exa:6-point_sets}.1)$ with the 2nd part of $(\ref{exa:6-point_sets}.3)$, we get $(2l_1)^2-4x_2=(2l_3)^2+4x_2$, so $x_2=\dfrac{l_1^2-l_3^2}{2}$.
Similarly, 
$x_3=\dfrac{l_2^2-l_1^2}{2}$, $x_1=\dfrac{l_3^2-l_2^2}{2}$ so that $x_1+x_2+x_3=0$.
From the second part of $(\ref{exa:6-point_sets}.2)$, we get
$x_1^2-4x_1+4+y_1^2=4l_2^2$, so
$$|O_{\pm} C_1|^2=x_1^2+y_1^2+1=4l_2^2+4x_1-3=2l_2^2+2l_3^2-3,$$
$$\text{similarly }|O_{\pm} C_2|^2=2l_3^2+2l_1^2-3,\;
|O_{\pm} C_3|^2=2l_1^2+2l_2^2-3.$$
Then $|C_1C_2|^2=(x_1-x_2)^2+(y_1-y_2)^2=x_1^2+y_1^2$.
\myskip

The last columns in Tables~\ref{tab:distance_pairs_T-} and \ref{tab:distance_pairs_T+} show the pairs of distances that distinguish $T_{+}\not\simeq T_-$.
The distributions $\wSDD(T_{\pm};2)$ can differ only by $\wRDD$s of the pairs $\{R,O_{\pm}\}, \{G,O_{\pm}\}, \{R,C_i\}, \{G,C_i\}$, where $i\in\{1,2,3\}$ is considered modulo 3 so that $1-1\equiv 3\pmod{3}$.
In rows of corresponding pairs of points, some pairs of distances are the same in both $\wSDD(T_{\pm};2)$, but other pairs differ.
If $l_1,l_2,l_3$ are pairwise distinct, the rows $\{R,O_-\},\{G,O_+\}$ include three different pairs of distances, so $\wSDD(T_{-};2)\neq \wSDD(T_{+};2)$.

\begin{table}[h!]
\centering
\begin{tabular}{@{}l|c|c|c@{}}
\hl{$T_-$} pair & distance & common pairs in $\wSDD(T_{\pm};2)$
& pairs that differ in $\wSDD($\hl{$T_{-}$}$;2)$ \smallskip \\
\hline
$\{R,O_{-}\}$ & $\sqrt{5}$ & $(\sqrt{13},\sqrt{32})$ to $G$ &
\begin{tabular}{l}
$(2$\hl{$\sqrt{l_3^2+1}$}$,\sqrt{2l_2^2+2l_3^2-3}) \text{ to } C_1$, \smallskip \\
$(2$\hl{$\sqrt{l_1^2+1}$}$,\sqrt{2l_3^2+2l_1^2-3}) \text{ to } C_2$, \smallskip \\
$(2$\hl{$\sqrt{l_2^2+1}$}$,\sqrt{2l_1^2+2l_2^2-3}) \text{ to } C_3$  \smallskip
\end{tabular} \\
\hline
$\{G,O_{-}\}$ & $\sqrt{13}$ & $(\sqrt{5},\sqrt{32})$ to $R$ &
\begin{tabular}{c}
\myskip
$(2$\hl{$\sqrt{l_2^2+1}$}$,\sqrt{2l_2^2+2l_3^2-3}) \text{ to } C_1$, \smallskip \\
$(2$\hl{$\sqrt{l_3^2+1}$}$,\sqrt{2l_3^2+2l_1^2-3}) \text{ to } C_2$, \smallskip\\
$(2$\hl{$\sqrt{l_1^2+1}$}$,\sqrt{2l_1^2+2l_2^2-3}) \text{ to } C_3$ 
\myskip 
\end{tabular} \\
\hline
$\{R,C_{i+1}\}$ & $2\sqrt{l_{i}^2+1}$ &
$\begin{array}{l}
(2\sqrt{l_{i-1}^2+1},\sqrt{32}) \text{ to } G,\smallskip \\ 
(2\sqrt{l_{i+1}^2+1},|C_{i+1} C_{i-1}|) \text{ to } C_{i-1},\smallskip \\
(2\sqrt{l_{i-1}^2+1},|C_i C_{i+1}|) \text{ to } C_i \smallskip 
\end{array}$ &
(\hl{$\sqrt{5}$}$,\sqrt{2l_{i-1}^2+2l_i^2-3})$ to $O_-$ \\
\hline
$\{G,C_{i-1}\}$ & $2\sqrt{l_i^2+1}$ &
$\begin{array}{l}
(2\sqrt{l_{i+1}^2+1},\sqrt{32}) \text{ to } R, \smallskip \\ 
(2\sqrt{l_{i-1}^2+1},|C_{i+1} C_{i-1}|) \text{ to } C_{i+1}, \smallskip \\
(2\sqrt{l_{i+1}^2+1},|C_{i-1} C_i|) \text{ to } C_i \smallskip 
\end{array}$ &
(\hl{$\sqrt{13}$}$,\sqrt{2l_i^2+2l_{i+1}^2-3})$ to $O_-$ 
\end{tabular}
\caption{Pairs of distances in the simplified invariant $\wSDD(T_{-};2)$.
For comparison with $\wSDD(T_{+};2)$, see Table~\ref{tab:distance_pairs_T+} .
The highlighted differences imply that $\wSDD(T_{-};2)\neq \wSDD(T_{+};2)$, so $T_-\not\simeq T_+$.}
\label{tab:distance_pairs_T-}
\end{table}

\begin{table}
\begin{tabular}{@{}l|c|c|c@{}}
\hl{ $T_+$ } pair & distance & common pairs in $\wSDD(T_{\pm};2)$
& pairs that differ in $\wSDD($\hl{$T_{+}$}$;2)$  \\
\hline
$\{G,O_{+}\}$ & $\sqrt{5}$ & $(\sqrt{13},\sqrt{32})$ to $R$ &
\begin{tabular}{c}
\bskip
$(2$\hl{$\sqrt{l_2^2+1}$}$,\sqrt{2l_2^2+2l_3^2-3}) \text{ to } C_1$, \smallskip \\ 
$(2$\hl{$\sqrt{l_3^2+1}$}$,\sqrt{2l_3^2+2l_1^2-3}) \text{ to } C_2$, \smallskip \\ 
$(2$\hl{$\sqrt{l_1^2+1}$}$,\sqrt{2l_1^2+2l_2^2-3}) \text{ to } C_3$ \smallskip 
\end{tabular} \\
\hline
$\{R,O_{+}\}$ & $\sqrt{13}$ & $(\sqrt{5},\sqrt{32})$ to $G$ &
\begin{tabular}{c}
\myskip
$(2$\hl{$\sqrt{l_3^2+1}$}$,\sqrt{2l_2^2+2l_3^2-3}) \text{ to } C_1$, \smallskip \\ 
$(2$\hl{$\sqrt{l_1^2+1}$}$,\sqrt{2l_3^2+2l_1^2-3}) \text{ to } C_2$, \smallskip \\ 
$(2$\hl{$\sqrt{l_2^2+1}$}$,\sqrt{2l_1^2+2l_2^2-3}) \text{ to } C_3$ 
\end{tabular} 
\myskip \\
\hline
$\{R,C_{i+1}\}$ & $2\sqrt{l_i^2+1}$ &
$\begin{array}{l}
(2\sqrt{l_{i-1}^2+1},\sqrt{32}) \text{ to } G, \smallskip  \\ 
(2\sqrt{l_{i+1}^2+1},|C_{i+1} C_{i-1}|) \text{ to } C_{i-1}, \smallskip \\
(2\sqrt{l_{i-1}^2+1},|C_i C_{i+1}|) \text{ to } C_i \smallskip 
\end{array}$ &
(\hl{$\sqrt{13}$}$,\sqrt{2l_{i-1}^2+2l_i^2-3})$ to $O_+$ \\
\hline
$\{G,C_{i-1}\}$ & $2\sqrt{l_i^2+1}$ &
$\begin{array}{l}
(2\sqrt{l_{i+1}^2+1},\sqrt{32}) \text{ to } R, \smallskip \\ 
(2\sqrt{l_{i-1}^2+1},|C_{i+1} C_{i-1}|) \text{ to } C_{i+1}, \smallskip \\
(2\sqrt{l_{i+1}^2+1},|C_{i-1} C_i|) \text{ to } C_i \smallskip 
\end{array}$ &
(\hl{$\sqrt{5}$}$,\sqrt{2l_i^2+2l_{i+1}^2-3})$ to $O_+$ 
\end{tabular}
\caption{Pairs of distances in the simplified invariant $\wSDD(T_{+};2)$.
For comparison with $\wSDD(T_{-};2)$, see Table~\ref{tab:distance_pairs_T-}.
The highlighted differences imply that $\wSDD(T_{-};2)\neq \wSDD(T_{+};2)$, so $T_-\not\simeq T_+$.} 
\label{tab:distance_pairs_T+}
\end{table}

Table~\ref{tab:distances_T-+} contains all pairwise distances between the points of $T_{\mp}$.
We show that $T_{\pm}$ differ by the simplified invariants $\wSDD(T_{\pm};2)$ below.
In each column of $R(C;A)$, we additionally allow any permutation of elements independent of other columns, so we could order each column (a pair of distances) lexicographically.
Denote the resulting simplification of $\RDD$ by $\wRDD$.  
Then $\wSDD(T_{\pm};2)$ have identical $\wRDD$s for the 2-point subsets $A$ from the list $\{R,G\}, \{O_{\pm},C_i\}, \{C_i,C_j\}$ for distinct $i,j=1,2,3$.
\myskip

For example, both $\wRDD(T_{\pm};\{R,G\})$ start with the distance $|R-G|=\sqrt{32}$ and then include the same four pairs $(\sqrt{5},\sqrt{13})$, $(2\sqrt{l_i^2+1},2\sqrt{l_{i-1}^2+1})$ for $i\in\{1,2,3\}$ modulo $3$, which should be ordered and written lexicographically.
\myskip

Hence, it makes sense to compare $\wSDD(T_{\pm};2)$ only by the remaining $\wRDD(T_{\pm};A)$ for $A$ from the list $\{R,O_{\pm}\}, \{G,O_{\pm}\}, \{R,C_i\}, \{G,C_j\}$ in Tables~\ref{tab:distance_pairs_T-} and~\ref{tab:distance_pairs_T+}.
\myskip

Without loss of generality, assume that $l_1\geq l_2\geq l_3$.
If all the lengths are distinct, then $l_1>l_2>l_3$.
Then the rows for $\{R,O_{-}\}$ and $\{G,O_{+}\}$ differ in Tables~\ref{tab:distance_pairs_T-} and~\ref{tab:distance_pairs_T+} even after ordering each pair so that a smaller distance precedes a larger one, and after writing all pairs lexicographically. 
So $\wSDD(T_-;2)\neq \wSDD(T_+;2)$ unless two of $l_i$ are equal.
If (say) $l_1=l_2$, the lexicographically ordered rows of $\{R,O_{-}\}$ and $\{G,O_{+}\}$ coincide in $\wSDD(T_{\pm};2)$, similarly for the rows of $\{G,O_{-}\}$ and $\{R,O_{+}\}$.
Hence, it suffices to compare only the six rows for the remaining pairs $\{R,C_i\},\{G,C_j\}$ in $\wSDD(T_{\pm};2)$.
\medskip

For $l_1=l_2$, we get $x_3=\dfrac{l_2^2-l_1^2}{2}=0$ and 
$x_1=-x_2=\dfrac{l_3^2-l_2^2}{2}$.
In equation $(\ref{exa:6-point_sets}.3)$ the equality 
$(x_1+2)^2+y_1^2=(x_2-2)^2+y_2^2$ with $x_1=-x_2$ implies that $y_1^2=y_2^2$.
The more degenerate case $l_1=l_2=l_3$, means that $x_1=x_2=x_3=0$ and $y_1^2=y_2^2=y_3^2$, hence at least two of $C_1,C_2,C_3$ should coincide.
The above contradiction means that it remains to consider the case $l_1=l_2>l_3$ when $x_1=-x_2\neq 0=x_3$ and $y_1=\pm y_2$, see Fig.~\ref{fig:6-point_sets}.
\medskip

If $y_1=y_2$, the clouds $T_{\pm}$ are isometric by $(x,y,z)\mapsto(-x,y,-z)$. 
If $y_1=-y_2$ and $y_3=0$, the clouds $T_{\pm}$ are isometric by the isometry $(x,y,z)\mapsto(-x,-y,-z)$.
If $y_1=-y_2$ and $y_3\neq 0$, then $C_1=(x_1,y_1,0)$, $C_2=(-x_1,-y_1,0)$, $C_3\neq(0,0,0)$.
Then among the six remaining rows, only the rows of $\{R,C_1\}$, $\{G,C_2\}$ have points at the distance $2\sqrt{l_3^2+1}$, see  Tables~\ref{tab:distance_pairs_T-} and~\ref{tab:distance_pairs_T+} for $i=3$ considered modulo 3.
Then $i+1\equiv 1\pmod{3}$, $i-1\equiv 2\pmod{3}$, so $l_{i+1}=l_1=l_2=l_{i-1}$.
\medskip

Looking at the rows of $\{R,C_1\}$, $\{G,C_2\}$, the three common pairs in each of $\SDD(T_{\pm};2)$ include the same distance $2\sqrt{l_1^2+1}=2\sqrt{l_2^2+1}$ but differ by $|C_{i-1} C_i|=|C_2 C_3|\neq |C_3C_1|=|C_i C_{i+1}|$ as $C_1=\pm C_2$, $C_3\neq(0,0,0)$.
\medskip

This couple of different rows implies that $\SDD(T_{-};2)\neq \SDD(T_{+};2)$ due to the swapped distances $\sqrt{5},\sqrt{13}$ in the remaining pairs, see Tables~\ref{tab:distance_pairs_T-+ex} for the clouds $T_{\pm}$ in Fig.~\ref{fig:6-point_sets} with
$l_1=l_2=\frac{\sqrt{13}}{2}$, $l_3=\frac{\sqrt{5}}{2}$. 
\eexa
\end{exa}

\begin{table}[h!]
\centering
\begin{tabular}{@{}l|c|c|c|c|c|@{}}
\hline
\hl{ $T_-$ } pair & distance &  dist. to neighb. 1 & dist. to neighb. 2 & dist. to neighb. 3 & dist. to neighb. 4 \\
\hline
$\{R,C_{1}\}$ & $3$ 
& $(\sqrt{2},\sqrt{17})$ to $C_3$ 
& \hl{ $(\sqrt{5},\sqrt{6})$ } to $O_{-}$ 
& $(\sqrt{17},\sqrt{20})$ to $C_2$ 
& $(\sqrt{17},\sqrt{32})$ to $G$ \smallskip \\
\hline
$\{G,C_{2}\}$ & $3$ 
& \hl{ $(\sqrt{6},\sqrt{13})$ } to $O_{-}$ 
& $(\sqrt{17},\sqrt{20})$ to $C_1$ 
& $(\sqrt{17},\sqrt{26})$ to $C_3$ 
& $(\sqrt{17},\sqrt{32})$ to $R$  \smallskip \\
\hline
\end{tabular}
\smallskip

\begin{tabular}{@{}l|c|c|c|c|c|@{}}
\hl{ $T_+$ } pair & distance &  dist. to neighb. 1 & dist. to neighb. 2 & dist. to neighb. 3 & dist. to neighb. 4 \\
\hline
$\{R,C_{1}\}$ & $3$ 
& $(\sqrt{2},\sqrt{17})$ to $C_3$ 
& \hl{ $(\sqrt{6},\sqrt{13})$ } to $O_{+}$ 
& $(\sqrt{17},\sqrt{20})$ to $C_2$ 
& $(\sqrt{17},\sqrt{32})$ to $G$  \smallskip\\
\hline
$\{G,C_{2}\}$ & $3$ 
& \hl{ $(\sqrt{5},\sqrt{6})$ } to $O_{+}$ 
& $(\sqrt{17},\sqrt{20})$ to $C_1$ 
& $(\sqrt{17},\sqrt{26})$ to $C_3$ 
& $(\sqrt{17},\sqrt{32})$ to $R$ \smallskip \\
\hline
\end{tabular}
\caption{The above rows show that $\SDD(T_-;2)\neq\SDD(T_+;2)$ for the clouds $T_{\pm}$ with $C_1=(-1,2,0)$, $C_2=(1,-2,0)$, $C_3=(0,3,0)$ so that $l_1=l_2=\frac{\sqrt{5}}{2}$, $l_3=\frac{\sqrt{13}}{2}$ in Tables~\ref{tab:distance_pairs_T-} and~\ref{tab:distance_pairs_T+}.}
\label{tab:distance_pairs_T-+ex}
\end{table}

Examples~\ref{exa:5-point_sets}, \ref{exa:7-point_sets}, and \ref{exa:6-point_sets} motivate the following conjecture.

\begin{conj}[completeness of $\SDD(C;h)$ in $\R^n$]
For any $n\geq 2$, there is some $2\leq h\leq n$ such that the Simplexwise Distance Distribution $\SDD(C;h)$ from Definition~\ref{dfn:SDD} is a complete isometry invariant of all clouds $C\subset\R^n$.
\epro
\end{conj}

\section{Continuous metrics on Simplexwise Distance Distributions}
\label{sec:SDD-metrics}

This section defines Lipschitz continuous metrics on $\SDD$s, which can be computable in a polynomial time of the number $m$ of points, for a fixed order $h$.
The $m-h$ permutable columns of the matrix $R(C;A)$ in $\RDD$ from Definition~\ref{dfn:RDD} can be interpreted as $m-h$ unlabelled points in $\R^h$.
Since any isometry is bijective, the simplest metric respecting bijections is the bottleneck distance $\BD$ from Example~\ref{exa:metrics}(b). 

\index{max metric}
\index{Relative Distance Distribution}
\begin{dfn}[the max metric $M_{\infty}$ on $\RDD$s]
\label{dfn:RDD_max-metric}
For any $m$-point clouds and ordered $h$-point base sequences $A\subset C$ and $A'\subset C'$, set 
$$d(\xi)=\max\{L_{\infty}(\xi(D(A)),D(A')),\BD(\xi(R(C;A)),R(C';A'))\}$$ for a permutation $\xi\in S_h$ on $h$ points. 
Then the \emph{max metric} on Relative Distance Distributions is defined as $M_{\infty}(\RDD(C;A),\RDD(C';A'))=\min\limits_{\xi\in S_h}d(\xi)$.
\edfn
\end{dfn}

We will use only $h=n$ for Euclidean space $\R^n$, so the factor $h!$ in Definition~\ref{dfn:RDD_max-metric} is practically small for $n=2,3$.
For $h=1$ and a 1-point sequence $A\subset C$, the matrix $D(A)$ is empty, so $d(\xi)=\BD(\xi(R(C;A)),R(C';A'))$.
The metric $M_{\infty}$ on $\RDD$s will be used for intermediate costs to get metrics on unordered collections of $\RDD$s ($\SDD$s) by using the standard tools in Definitions~\ref{dfn:EMD} and~\ref{dfn:LAC_matrix} below.  

\begin{dfn}[Linear Assignment Cost of a matrix {\cite{fredman1987fibonacci}}]
\label{dfn:LAC_matrix}
For any $k\times k$ matrix of costs $c(i,j)\geq 0$, $i,j\in\{1,\dots,k\}$, the \emph{Linear Assignment Cost}   
$\LAC=\frac{1}{k}\min\limits_{g}\sum\limits_{i=1}^k c(i,g(i))$ is minimized for all bijections $g$ on the indices $1,\dots,k$.
\edfn
\end{dfn}

\index{Earth Mover's Distance}
\index{Linear Assignment Cost}

The normalisation factor $\frac{1}{k}$ in $\LAC$ makes this metric better comparable with $\EMD$ in Definition~\ref{dfn:EMD} whose weights sum up to 1.
For both $\LAC$ and $\EMD$, the matrix of initial costs will consist of max metrics between all $\RDD$s in two given $\SDD$s. 
\myskip
  
Theorem~\ref{thm:SDD_metrics}(b) extends the $O(m^{1.5}\log^n m)$ algorithm for fixed clouds of $m$ unlabelled points in \cite[Theorem~6.5]{efrat2001geometry} to the harder case of isometry classes but keeps the polynomial time in $m$ for a fixed dimension $n$.

\index{polynomial-time complexity}
\index{Simplexwise Distance Distribution}
\begin{thm}[time of metrics on $\SDD$s, {\cite[Theorem~5.5]{kurlin2023simplexwise}}]
\label{thm:SDD_metrics}
\tb{(a)}
For any $m$-point clouds $C,C'$ in their own metric spaces and $h\geq 1$, let the Simplexwise Distance Distributions $\SDD(C;h)$ and $\SDD(C';h)$ consist of $k=\binom{m}{h}$ $\RDD$s with equal weights $\frac{1}{k}$ without collapsing identical $\RDD$s.
\myskip

\noindent
\textbf{(b)}
Using the $k\times k$ matrix of costs computed by the max metric $M_{\infty}$ between $\RDD$s from $\SDD(C;h)$ and $\SDD(C';h)$, the Linear Assignment Cost $\LAC$ from Definition~\ref{dfn:LAC_matrix} satisfies all metric axioms on $\SDD$s and can be computed in time $O(h!(h^2 +m^{1.5}\log^h m)k^2 + k^3\log k)$.
\smallskip

\noindent
\textbf{(b)}
Let $\SDD(C;h)$ and $\SDD(C';h)$ have a maximum size $l\leq k$ after collapsing identical $\RDD$s. 
Using the same matrix of max metrics as in part (b), the $\EMD$ from Definition~\ref{dfn:EMD} satisfies all metric axioms  on $\SDD$s and can be computed in time $O\big(h!(h^2 +m^{1.5}\log^h m) l^2 +l^3\log l\big)$.
\ethm
\end{thm}

Theorem~\ref{thm:SDD_continuity} substantially generalizes the fact that perturbing two points in their $\ep$-neighbourhoods changes the distance between these points by at most $2\ep$. 

\index{Lipschitz continuity}
\index{Earth Mover's Distance}
\index{Simplexwise Distance Distribution}

\begin{thm}[Lipschitz continuity of $\SDD$s, {\cite[Theorem~5.8]{kurlin2023simplexwise}}]
\label{thm:SDD_continuity}
In any metric space, let $C'$ be obtained from a cloud $C$ by perturbing every point of $C$ within its $\ep$-neighbourhood.
For any order $h\geq 1$, $\SDD(C;h)$ changes by at most $2\ep$ in the $\LAC$ and $\EMD$ metrics.
The lower bound holds: $\EMD\big(\SDD(C;h),\SDD(C';h)\big)\geq|\SDM(C;h,1)-\SDM(C';h,1)|_\infty$.
\ethm
\end{thm}

\section{Measured Simplexwise Distributions for metric-measure spaces}
\label{sec:MSD}

This section adapts  Simplexwise Distance Distributions $\SDD$ to metric-measure spaces.

\index{metric-measure space}
\begin{dfn}[metric-measure space]
\label{dfn:metric-measure space}
A \emph{metric-measure space} $(X,d_X,\mu_X)$ is a compact space $X$ with a metric $d_X$ and a Borel measure $\mu_X$ such that $\mu_X(X)<+\infty$.
An \emph{isomorphism} between metric-measure spaces is an isometry $f:X\to Y$ that respects the measures in the sense that $\mu_Y(U)=\mu_X(f^{-1}(X))$ for any subset $U\subset Y$.
\edfn
\end{dfn}

Dividing $\mu_X(U)$ by the measure $\mu_X(X)<+\infty$ for any $U\subset X$, we can assume that $\mu_X(X)=1$, so $\mu_X$ is a probability measure.
Any metric space $X$ of $m$ points can be considered a metric-measure space with the uniform measure $\mu_X(p)=\dfrac{1}{m}$ for all $p\in X$.
On two points $0,1$ in $\R$, the metric-measure spaces $X=(\{0,1\},1,\{\frac{1}{2},\frac{1}{2}\})$ and $Y=(\{0,1\},1,\{\frac{1}{3},\frac{2}{3}\})$ are isometric but not isomorphic because of different weights.
\myskip

Problem~\ref{pro:metric_space} becomes much harder if we replace isometries between metric spaces with isomorphisms between metric-measure spaces because all known isometry invariants should be further refined to distinguish under isomorphism.
Definition~\ref{dfn:MSD} extends the local distribution of distances from \cite[Definition~5.5]{memoli2011gromov} to orders $h>1$.

\index{metric-measure space}
\index{measured simplexwise distribution}
\begin{dfn}[Measured Simplexwise Distribution $\MSD$]
\label{dfn:MSD}
Let $(X,d_X,\mu_X)$ be any metric-measure space.
For any base sequence $A=(p_1,\dots,p_h)$ of $h\geq 1$ ordered points of $X$, write the triangular distance matrix $D(A)$ from Definition~\ref{dfn:RDD} row-by-row as the vector $\vec v(A)\in\R_+^{h(h-1)/2}$ so that $\vec v_k=d_X(p_i,p_j)$ for $k=h(i-1)+j-1$, $1\leq i<j\leq h$.
For a vector $\vec d=(d_1,\dots,d_h)\in\R_+^h$ of distance thresholds, the vector $\vec m(A;\vec d)\in\R_+^h$ consists of $h$ values $\mu_X(\{q\in X\vl d_X(q,p_i)\leq d_i\})$ for $i=1,\dots,h$.
\sskip
  
The \emph{Measured Simplexwise Distribution} of order $h\geq 1$ is the function $\MSD[X;h]: X^h\times\R_+^h\to\R_+^{h(h+1)/2}$ mapping any $A\in X^h$ and $\vec d\in\R_+^h$ to the pair $[\vec v(A),\vec m(A;\vec d)]$ considered as a concatenated vector in $\R_+^{h(h+1)/2}$.
\edfn
\end{dfn}

For $h=1$, the vector $\vec v(A)$ is empty and the Measured Simplexwise Distribution of order $h=1$ coincides with the local distribution of distances \cite[Definition~5.5]{memoli2011gromov} $\MSD[X;1]: X\times\R_+\to\R_+$ mapping any point $p\in X$ and a threshold $d\in\R_+$ to the measure value 
$\mu_X(\{q\in X \vl d_X(q,p)\leq d\})$.
\smallskip

Any permutation $\xi$ on indices $1,\dots,h$ naturally permutes the components of $\MSD[X;h]$.
If $X$ consists of $m$ points, $\MSD[X;h]$ reduces to the finite collection
 of $\binom{m}{h}$ vectors $\VID(A)$ paired with fields $\VSM(A;\vec d):\R_+^h\to\R_+^{h}$ only for unordered $h$-point subsets $A\subset X$, which can be refined to a stronger analogue of $\SDD$ below.  

\begin{figure}[h!]
\centering
\includegraphics[width=\linewidth]{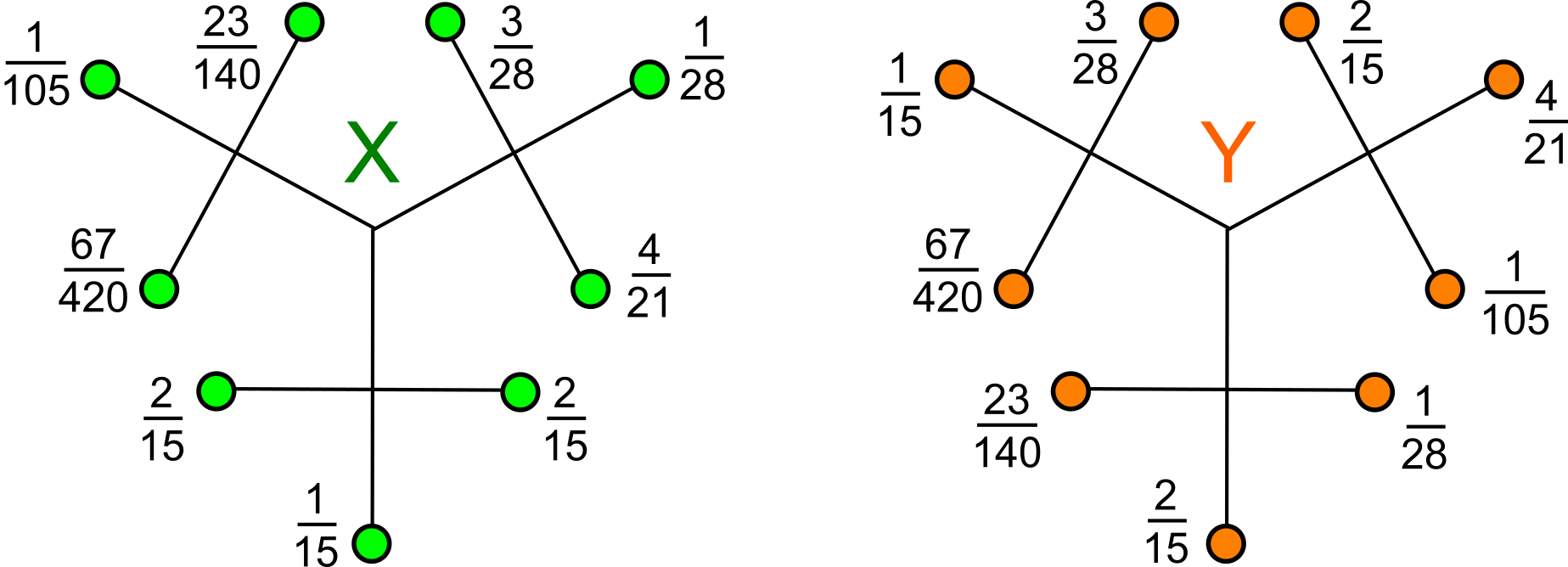}
\caption{Non-isomorphic metric-measure spaces $X,Y$ from \cite[Fig.~8]{memoli2011gromov} have equal local distributions of distances but are distinguished by the new Weighted Simplexwise Distribution of order 1 and the Measured Simplexwise Distributions of order 2, see details in Example~\ref{exa:9-point_trees}.
All edges have length $\frac{1}{2}$. 
}
\label{fig:9-point_trees}
\end{figure}

\index{metric-measure space}
\index{weighted simplexwise distribution}
\begin{dfn}[Weighted Simplexwise Distribution $\WSD$]
\label{dfn:WSD}
Let $X$ be a finite metric-measure space whose any point $p$ has a weight $w_(p)$.
For $h\geq 1$ and a base sequence $A=(p_1,\dots,p_h)$ of $h$ ordered points of $X$ in Definition~\ref{dfn:RDD}, endow any distance $d(p,q)$ in $D(A)$ with the unordered pair $w(p),w(q)$ of weights.
For every point $q\in X-A$, put the weight $w(q)$ in the extra $(h+1)$-st row of the matrix $M(X;A)$ whose columns are indexed by unordered $q\in X-A$.
If $h=1$ and $A=p_1$, set $D(A)=w(p_1)$.
\sskip

The \emph{Weighted Distance Distribution} $\WDD(X;A)$ is the equivalence class of the pair $[D(A);M(X;A)]$ under permutations $\xi\in S_h$ acting on $A$. 
The \emph{Weighted Simplexwise Distribution} $\WSD(X)$ is the unordered collection of $\WDD(X;A)$ for all subsets $A\subset X$ of unordered $h$ points.
\edfn
\end{dfn}

\index{Earth Mover's Distance}

For finite metric-measure spaces, a metric on $\WDD$s can be defined similar to $M_\infty$ from Definition~\ref{dfn:RDD_max-metric} by combining the weights and distances.
Then $\LAC$ and $\EMD$ from Definitions~\ref{dfn:LAC_matrix} and~\ref{dfn:EMD} can be computed as in Theorem~\ref{thm:SDD_metrics}.

\begin{exa}[the strength of $\WSD$ for 9-point trees]
\label{exa:9-point_trees}
Fig.~\ref{fig:9-point_trees} shows metric-measure spaces $X,Y$ on 9 points visualised as trees \cite[Fig.~8]{memoli2011gromov}.
All edges have length $\frac{1}{2}$ and induce the shortest-path metrics $d_X,d_Y$.
The sum of weights in every small branch of 3 nodes is $\frac{1}{3}$.
These metric-measure spaces $X,Y$ have all inter-point distances only 1 and 2, and equal local distributions of distances $\MSD[X;1]=\MSD[Y;1]$ by \cite[Example~5.6]{memoli2011gromov}.
\myskip

Indeed, both $\MSD$s can be considered the same set of 9 piecewise constant functions $\mu(p)$ taking values $w(p)$, $\frac{1}{3}$, and $1$ on the intervals $[0,1)$, $[1,2)$, $[2,+\infty)$, respectively.
\myskip

However, $\WSD$s have more pointwise data: $\WSD[X;1]$ has $A(D)=w(p)=\frac{23}{140}$ and the following $2\times 8$ matrix 
$$M(X;p)=\left(\begin{array}{cccccccc}
1 & 1 & 2 & 2 & 2  & 2 & 2 & 2 \\
\dfrac{1}{105} & \dfrac{67}{240} & \dfrac{2}{15} & \dfrac{1}{15} & 
\dfrac{2}{15} & \dfrac{4}{21} & \dfrac{1}{28} &  \dfrac{3}{28}
\end{array}\right),$$ but 
$\WSD[Y;1]$ has another matrix for $w(p)=\frac{23}{140}$.
$$M(Y;p)=\left(\begin{array}{cccccccc}
1 & 1 & 2 & 2 & 2  & 2 & 2 & 2 \\
\dfrac{2}{15} & \dfrac{1}{28} & \dfrac{1}{105} & \dfrac{4}{21} & 
\dfrac{2}{15} & \dfrac{3}{28} & \dfrac{1}{15} & \dfrac{67}{420} 
\end{array}\right).$$
The above matrices with freely permutable columns are different, so
 $X,Y$ are distinguished by the Weighted Simplexwise Distributon $\WSD$ of order $h=1$.
\myskip 

Also, $\MSD[X;2]\neq\MSD[Y;2]$ because, for any base sequence $A=(p,q)\in X^2$, we have $\VSM[X;2](A;d_1,d_2)=(w(p),w(q))$ for $d_1,d_2<1$ since all other points have minimum distance $1$ from $p,q$, similarly for $Y$.
The unique points $p,q$ of weights $w(p)=\dfrac{23}{140}$ and $w(q)=\dfrac{67}{420}$ have different distances $d_X(p,q)=1$ and $d_Y(p,q)=2$.
Then $\MSD[X;2]\neq\MSD[Y;2]$ differ by the uniquely identifiable fields mapping $[0,1)^2$ to the constant vector $(w(p),w(q))$ with $\VID_X(A)=1\neq 2=\VID_Y(A)$.
\eexa
\end{exa}

We conjecture that any metric-measure spaces $X,Y$ are distinguished under isomorphism by Measured Simplexwise Distributions for a high enough $h$ depending on $X,Y$.

\bibliographystyle{plain}
\bibliography{Geometric-Data-Science-book}

%
%
%

\chapter{Complete and Lipschitz continuous invariants of unordered points in $\R^n$}
\label{chap:SCD} 

\abstract{
This chapter leverages the Euclidean structure of $\R^n$ to improve the Simplexwise Distance Distribution to a smaller Simplexwise Centred Distribution (SCD) for any unordered points.
The new invariant is complete under rigid motion and computable in polynomial time for a fixed dimension.
The key ingredient of Lipschitz continuity is the new strength of a simplex, which is a linear-growth analogue of the simplex volume.
}

\section{Geo-mapping problem under rigid motion in $\R^n$}
\label{sec:OSD}

This chapter follows paper \cite{widdowson2023recognizing} and its extension \cite{kurlin2023strength} to Euclidean spaces.
Problem~\ref{pro:Euclidean_unordered} adjusts Geo-Mapping Problem~\ref{pro:geocodes}
to finite clouds of unordered points in $\R^n$.
\myskip

The major difference with Problem~\ref{pro:metric_space}, which was stated under isometry in a metric space, is the full completeness under the stronger equivalence of rigid motion in $\R^n$.

\begin{pro}[partial case of Problem~\ref{pro:geocodes} for clouds under rigid motion in $\R^n$]
\label{pro:Euclidean_unordered}
Design a map $I$ on finite clouds of unordered points in $\R^n$ with values in a metric space satisfying the following conditions.
\myskip

\noindent
\tb{(a)} 
\emph{Completeness:} 
any clouds $A,B\subset\R^n$ are related by rigid motion ($A\cong B$)  in $\R^n$ if and only if $I(A)=I(B)$.
\myskip

\noindent
\tb{(b)} \emph{Metric:} 
there is a distance $d$ on the \emph{invariant space} $\{I(A) \vl A\subset\R^n\}$ satisfying all metric axioms in Definition~\ref{dfn:metrics}(a). 
\myskip

\nt
\tb{(c)} 
\emph{Continuity:} 
there is a constant $\la$ such that, for any $\ep>0$, if $B$ is obtained from a cloud $A\subset\R^n$ by perturbing every point up to Euclidean distance $\ep$, then $d(I(A),I(B))\leq \la\ep$.  
\myskip

\nt
\tb{(d)} 
\emph{Computability:} 
for a fixed dimension $n$, the invariant $I(A)$ and the metric $d(I(A),I(B))$ can be computed in times that depend polynomially on the maximum size $\max\{|A|,|B|\}$ of clouds $A,B\subset\R^n$.
\epro 
\end{pro} 

Problem~\ref{pro:Euclidean_unordered} will be fully solved by the Oriented Simplexwise  Distribution ($\OSD$), which we introduce in Definition~\ref{dfn:OSD} after a few auxiliary concepts below.

\begin{dfn}[matrices $D(A)$ and $M(C;A)$ for $A\subset C$]
\label{dfn:D(A)+M(C;A)}
Let $C$ be a cloud of $m$ unordered points in $\R^n$ with a fixed orientation.
Let $A=(p_1,\dots,p_n)$ be a \emph{base sequence} of $n$ distinct ordered points of $C$.
Let $D(A)$ be the $n\times n$ \emph{distance} matrix whose entry $D(A)_{i,j}$ is Euclidean distance $|p_i-p_j|$ for $1\leq i<j\leq n$, all other entries are zeros.
For any other point $q\in C-A$, write distances from $q$ to $p_1,\dots,p_n$ as a column.
Form the $n\times (m-n)$-matrix by these $m-n$ lexicographically ordered columns.
\sskip

At the bottom of the column of a $q\in C-A$, add the sign of the determinant consisting of the vectors $q-p_1,\dots,q-p_n$.  
The resulting $(n+1)\times(m-n)$-matrix with signs in the bottom $(n+1)$-st row is the \emph{oriented relative distance} matrix $M(C;A)$.
\edfn
\end{dfn}

Let $S_n$ denote the permutation group on indices $1,\dots,n$.
Any permutation $\xi\in S_n$ is a composition of some $t$ transpositions $i\lra j$ and has $\sign(\xi)=(-1)^t$. 

\begin{dfn}[oriented distributions $\ORD(C;A)$ and $\OSD(C)$ for a cloud $C\subset\R^n$]
\label{dfn:OSD}
Any permutation $\xi\in S_n$ acts on $D(A)$ by mapping $D(A)_{ij}$  to $D(A)_{kl}$, where $k\leq l$ is the pair $\xi(i),\xi(j)-1$ written in increasing order.
Then the permutation $\xi$ acts on $M(C;A)$ by mapping any $i$-th row to the $\xi(i)$-th row and by multiplying the $(n+1)$-st row by $\sign(\xi)$, after which all columns are written in the lexicographic order.
\smallskip

The \emph{Oriented Relative Distribution} $\ORD(C;A)$ is the equivalence class of the pair $[D(A);M(C;A)]$ under all permutations $\xi\in S_n$ acting on both $D(A)$ and $M(C;A)$. 
\sskip

The \emph{Oriented Simplexwise Distribution} $\OSD(C)$ is the unordered collection of $\ORD(C;A)$ for all $\binom{m}{n}$ unordered subsets $A\subset C$ of $n$ points. 
\edfn
\end{dfn}

Any mirror reflection in $\R^n$ reverses the sign of the $n\times n$ determinant consisting of vectors $v_1,\dots,v_n\in\R^n$, hence reverses all signs in the $(n+1)$-st rows of the matrices $M(C;A)$ in Oriented Relative Distributions.
$\mORD(C;A)$ and $\mOSD(C)$ denote the `mirror images' of $\ORD(C;A)$ and $\OSD(C)$, respectively, with all signs reversed.

\begin{figure}[h!]
\centering
\includegraphics[width=\linewidth]{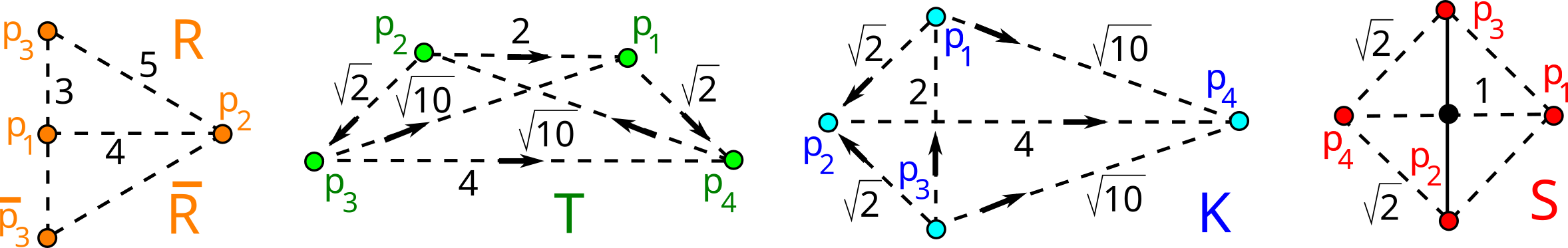}
\caption{\textbf{1st}: the right-angled cloud $R\subset\R^2$ consisting of points $p_1=(0,0)$, $p_2=(4,0)$, $p_3(0,3)$, and its mirror image $\bar R$ of $p_1,p_2$, and $\bar p_3=(0,-3)$ with respect to the $x$-axis. 
\textbf{2nd}: the trapezium cloud $T\subset\R^2$ consisting of points $p_1=(1,1)$, $p_2=(-1,1)$, $p_3=(-2,0)$, $p_4=(2,0)$. 
\textbf{3rd}: the kite cloud $K\subset\R^2$ consisting of points $p_1=(0,1)$, $p_2=(-1,0)$, $p_3=(0,-1)$, $p_4=(3,0)$. 
\textbf{4th}: the square cloud $S\subset\R^2$ consisting of points $p_1=(1,0)$, $p_2=(0,-1)$, $p_3=(-1,0)$, $p_4=(0,1)$. 
}
\label{fig:clouds_oriented}
\end{figure}
 
\begin{exa}[$\OSD$ for mirror images on right-angled clouds]
\label{exa:OSD+mirror}
In $\R^2$ with the counter-clockwise orientation, the right-angled cloud $R$ on the vertices $p_1=(0,0)$, $p_2=(4,0)$, $p_3=(0,3)$ of the triangle in Fig.~\ref{fig:clouds_oriented}~(1st) has the distribution $\OSD(R)$ consisting of 
$$\ORD(R;(p_1,p_2))=[4,\colthree{3}{5}{+}],$$
$$\ORD(R;(p_2,p_3))=[5,\colthree{4}{3}{+}],$$
$$\ORD(R;(p_3,p_1))=[3,\colthree{5}{4}{+}].$$
If we swap the points $p_1\lra p_3$, the last $\ORD$ above changes to the equivalent form $\ORD(R;(p_1,p_3))=[3,\colthree{4}{5}{-}]$, without affecting others.
If we reflect $R$ with respect to the $x$-axis, the mirror image $\bar R$ of $p_1,p_2,\bar p_3=(0,-3)$ has 
$\OSD(\bar R)=\mOSD(R)$ with
$$\ORD(\bar R;(p_1,p_2))=[4,\colthree{3}{5}{-}],$$
$$\ORD(\bar R;(p_2,\bar p_3))=[5,\colthree{4}{3}{-}],$$
$$\ORD(\bar R;(p'_3,p_1))=[3,\colthree{5}{4}{-}],$$
 which differs from $\OSD(R)$ even if we swap points in each pair.
\eexa
\end{exa}

\begin{exa}[$\OSD$ for $T,K$]
\label{exa:OSD+T+K}
Fix the \emph{counter-clockwise} orientation on $\R^2$ so that if a vector $\vec v$ is obtained from $\vec u$ by a counter-clockwise rotation, then $\det(u,v)>0$. 
Table~\ref{tab:OSD+TK} shows the Oriented Simplexwise Distributions for the clouds $T,K$ in Fig.~\ref{fig:clouds_oriented}.
Each row contains the most similar $\ORD$s whose differences are \hl{highlighted}.
\eexa
\end{exa}

\begin{table}[h!]
  \centering
  \begin{tabular}{@{}l|l@{}}
    $\ORD$s in $\OSD(T)$ & $\ORD$s in $\OSD(K)$ \\
    \hline
    \myskip
    
	$[\sqrt{2},\mathree{2}{\cbox{yellow}{\sqrt{10}}}{-}{\sqrt{10}}{4}{-}]$ &
	$[\sqrt{2},\mathree{2}{\cbox{yellow}{\sqrt{2}}}{-}{\sqrt{10}}{4}{-}]$ 
\myskip \\

	$[\sqrt{2},\mathree{2}{\cbox{yellow}{\sqrt{10}}}{+}{\sqrt{10}}{4}{+}]$ &
	$[\sqrt{2},\mathree{2}{\cbox{yellow}{\sqrt{2}}}{+}{\sqrt{10}}{4}{+}]$ 
\myskip \\
	
	    $[2,\mathree{\sqrt{2}}{\cbox{yellow}{\sqrt{10}}}{-}{\sqrt{10}}{\cbox{yellow}{\sqrt{2}}}{\cbox{yellow}{-}}]$ & 
    $[2,\mathree{\sqrt{2}}{\cbox{yellow}{\sqrt{2}}}{-}{\sqrt{10}}{\cbox{yellow}{\sqrt{10}}}{\cbox{yellow}{+}}]$
\myskip \\
    
	$[\sqrt{10},\mathree{\sqrt{2}}{\cbox{yellow}{2}}{\cbox{yellow}{+}}{\cbox{yellow}{4}}{\cbox{yellow}{\sqrt{2}}}{-}]$ &
	$[\sqrt{10},\mathree{\sqrt{2}}{\cbox{yellow}{4}}{\cbox{yellow}{-}}{\cbox{yellow}{2}}{\cbox{yellow}{\sqrt{10}}}{-}]$ \myskip \\
	
	$[\sqrt{10},\mathree{\sqrt{2}}{\cbox{yellow}{2}}{\cbox{yellow}{-}}{\cbox{yellow}{4}}{\cbox{yellow}{\sqrt{2}}}{+}]$ & 
	$[\sqrt{10},\mathree{\sqrt{2}}{\cbox{yellow}{4}}{\cbox{yellow}{+}}{\cbox{yellow}{2}}{\cbox{yellow}{\sqrt{10}}}{+}]$ \smallskip \\
	
	    $[4,\mathree{\sqrt{2}}{\sqrt{10}}{+}{\cbox{yellow}{\sqrt{10}}}{\cbox{yellow}{\sqrt{2}}}{\cbox{yellow}{+}}]$ &
	$[4,\mathree{\sqrt{2}}{\sqrt{10}}{+}{\cbox{yellow}{\sqrt{2}}}{\cbox{yellow}{\sqrt{10}}}{\cbox{yellow}{-}}]$ 
\end{tabular}
  \caption{The Oriented Simplexwise Distributions $\OSD$s from Definition~\ref{dfn:OSD} for the 4-point clouds $T,K\subset $ in Fig.~\ref{fig:clouds_oriented}. 
Forgetting all signs in the bottom rows of $\ORD$s gives $\SDD$s in Table~\ref{tab:SDD+TK}.
}
  \label{tab:OSD+TK}
\end{table}

Though Problem~\ref{pro:Euclidean_unordered} did not include the reconstruction condition as in \ref{pro:geocodes}(b), Lemma~\ref{lem:ORD_reconstruction} below proves this reconstruction in more detail than \cite[Lemma 3.6]{kurlin2023strength}. 

Recall that Definition~\ref{dfn:aff} introduces the affine dimension of a base sequence $A$ of $n$ ordered points $p_1,\dots,p_n$ as the maximum dimension of the vector space generated by all inter-point vectors $\vec p_i-\vec p_j$ for $i,j\in\{1,\dots,n\}$. 

\begin{lem}[reconstruction from $\ORD$]
\label{lem:ORD_reconstruction}
A cloud $C\subset\R^n$ of $m>n$ unordered points can be reconstructed, uniquely under rigid motion, from $\ORD(C;A)$ in Definition~\ref{dfn:OSD} for any base sequence $A\subseteq C$ with $\aff(A)=n-1$.  
\elem
\end{lem}
\begin{proof}
By Lemma~\ref{lem:seq_reconstruction}(b), any base sequence $A\subseteq C$ can be reconstructed, uniquely under rigid motion in $\R^n$, from the triangular distance matrix $D(A)$ in Definition~\ref{dfn:OSD}.
We may assume that the first $n$ points $p_1,\dots,p_n$ of $A\subseteq C$ span the subspace of the first $n-1$ coordinate axes of $\R^n$. 
We prove that any point $q\in C-A$ has a unique location in $\R^n$, determined by the $n$ distances $|q-p_1|,\dots,|q-p_n|$ written in a column of the matrix $\ORD(C;A)$.
Since the points of $A$ do not belong to any $(n-1)$-dimensional affine subspace of $\R^n$, the $n$ spheres $S(p_i;|q-p_i|)$ of radii $|q_i-p_i|$ and centres $p_i$, $i=1,\dots,n$, contain $q$ and their full intersections consists of one or two points. 
We can uniquely choose $q$ among these two options due to the sign of the determinant on the column vectors $\vec q-\vec p_1,\dots,\vec q-\vec p_n$ in the bottom row of $\ORD(C;A)$. 
\end{proof}

Lemma~\ref{lem:ORD_reconstruction} implies that $\ORD(C;A)$ can have identical columns only for degenerate subsets $A\subset C$ with $\aff(A)<n-1$.
For example, let $n=3$ and $A$ consist of three points $p_1,p_2,p_3$ in the same straight line $L\subset\R^3$.
The three distances $|q-p_i|$, $i=1,2,3$, to any other point $q\in C$ outside $L$ define three spheres $S(p_i;|q-p_i|)$ that share a common circle in $\R^3$, so the position of $q$ is not uniquely determined in this case.
\smallskip

Though one $\ORD(C;A)$ with $\aff(A)=n-1$ suffices to reconstruct $C\subset\R^n$ up to rigid motion, the dependence on a subset $A\subset C$ required us to consider the larger Oriented Simplexwise Distribution $\OSD(C)$ for all $n$-point subsets $A\subset C$ to get a complete invariant in Theorem~\ref{thm:OSD_complete}. 
Equality $\OSD(C)=\OSD(C')$ is interpreted as a bijection $\OSD(C)\to\OSD(C')$ matching all ORDs. 

\begin{thm}[completeness of $\OSD$, {\cite[Theorem~3.7]{kurlin2023strength}}]
\label{thm:OSD_complete}
The Oriented Simplexwise Distribution $\OSD(C)$ in Definition~\ref{dfn:OSD} is a complete isometry invariant and can be computed in time $O(m^{n+1}/(n-3)!)$. 
So any clouds $C,C'\subset\R^n$ of $m$ unlabelled points are related by rigid motion (isometry, respectively) \emph{if and only if} 
$\OSD(C)=\OSD(C')$ ($\OSD(C)=\OSD(C')$ or its mirror image $\mOSD(C')$, respectively).
\ethm
\end{thm}

\section{Simplexwise Centered Distributions of a cloud in $\R^n$}
\label{sec:SCD}

This section simplifies the $\OSD$ invariant to the Simplexwise Centred Distribution ($\SCD$) in Definition~\ref{dfn:SCD}.
The Euclidean structure of $\R^n$ allows us to translate the \emph{centre of mass} $\dfrac{1}{m}\sum\limits_{p\in C} p$ of a given $m$-point cloud $C\subset\R^n$ to the origin $0\in\R^n$.
Then Problem~\ref{pro:Euclidean_unordered} reduces to only rotations around $0$ from the orthogonal group $\Or(\R^n)$.
\myskip

Definition~\ref{dfn:OSD} introduced the Oriented Simplexwise Distribution (OSD) as an ordered collection of $\ORD(C;A)$ for all $\binom{m}{n}$ unordered subsets $A\subset C$ of $n$ points.
Including the centre of mass allows us to consider the smaller number of $\binom{m}{n-1}$ subsets $A\subset C$ of $n-1$ points instead of $n$.
\myskip

Though the centre of mass is uniquely determined for any cloud $C\subset\R^n$ of unordered points, real applications may offer one or several labelled points of $C$ that substantially speed up metrics on invariants.
For example, an atomic neighbourhood in a solid material is a cloud $C\subset\R^3$ of atoms around a central atom, which may not be the centre of mass of $C$, but can be an extra base point in all constructions below.
\smallskip

For any base sequence $A$ of $n-1$ ordered points $p_1,\dots,p_{n-1}\in C$, add the origin $0$ as the $n$-th point and consider the $n\times n$
distance matrix $D(A\cup\{0\})$ and the $(n+1)\times (m-n)$ matrix $M(C;A\cup\{0\})$ in Definition~\ref{dfn:D(A)+M(C;A)}.
Any $n$ vectors $\vec v_1,\dots,\vec v_n\in\R^n$ can be written as columns in the $n\times n$ matrix whose determinant has a \emph{sign} $\pm 1$ or $0$ (if the vectors $\vec v_1,\dots,\vec v_n$ are linearly dependent).
Any permutation $\xi\in S_{n-1}$ of $n-1$ points of $A$ acts on $D(A)$ by permuting the first $n-1$ rows of $M(C;A\cup\{0\})$ and by multiplying every sign in the $(n+1)$-st row by $\sign(\xi)$.

\begin{dfn}[Simplexwise Centred Distribution $\SCD$]
\label{dfn:SCD}
Let $C\subset\R^n$ be any cloud of $m$ unlabelled points.
For any base sequence $A$ of ordered $p_1,\dots,p_{n-1}$ in a cloud $C\subset\R^n$ with the center of mass at $0\in\R^n$,  the \emph{Oriented Centred Distribution} $\OCD(C;A)$ is the equivalence class of pairs $[D(A\cup\{0\}),M(C;A\cup\{0\})]$ considered up to permutations $\xi\in S_{n-1}$ of points of $A$.
The \emph{Simplexwise Centred Distribution} $\SCD(C)$ is the unordered set of distributions $\OCD(C;A)$ for all $\binom{m}{n-1}$ unordered $(n-1)$-point subsets $A\subset C$.
The mirror image $\mSCD(C)$ is obtained from $\SCD(C)$ by reversing all signs.
\edfn
\end{dfn}

Definition~\ref{dfn:SCD} needs no permutations for any $C\subset\R^2$ as $n-1=1$.
Columns of $M(C;A\cup\{0\})$ can be lexicographically ordered without affecting future metrics. 
\myskip

Some of the $\binom{m}{n-1}$ $\OCD$s in $\SCD(C)$ can be identical as in Example~\ref{exa:SCD}(b).
If we collapse any $l>1$ identical $\OCD$s into a single $\OCD$ with the \emph{weight} $l/\binom{m}{h}$, $\SCD$ can be considered as a weighted probability distribution 
of $\OCD$s.

\begin{exa}[Simplexwise Centered Distribution $\SCD$s for clouds in Fig.~\ref{fig:clouds_oriented}]
\label{exa:SCD}
\textbf{(a)}
Let $R\subset\R^2$ be the right-angled cloud of the points $p_1=(0,0)$, $p_2=(4,0)$, $p_3=(0,3)$  in Fig.~\ref{fig:clouds_oriented}~(1st).
Though $p_1=(0,0)$ is included in $R$ and is not its centre of mass,
$\SCD(R)$ still makes sense.
In
$\OCD(R;p_1)=\left[0,\left( \begin{array}{cc} 
4 & 3 \\
4 & 3 \\
0 & 0
\end{array}\right) \right]$,
the matrix $D(\{p_1,0\})$ is $|p_1-0|=0$,
the top row has $|p_2-p_1|=4$, $|p_3-p_1|=3$.
In $\OCD(R;p_2)=\left[4,\left( \begin{array}{cc} 
4 & 5 \\
0 & 3 \\
0 & -
\end{array}\right) \right]$,
the first row has $|p_1-p_2|=4$, $|p_3-p_2|=5$,
the second row has $|p_1-0|=0$, $|p_3-0|=3$,
$\det\mat{-4}{0}{3}{3}<0$.
In $\OCD(R;p_3)=\left[3,\left( \begin{array}{cc} 
3 & 5 \\
0 & 4 \\
0 & +
\end{array}\right) \right]$,
the first row has $|p_1-p_3|=3$, $|p_2-p_3|=5$,
the second row has $|p_1-0|=0$, $|p_2-0|=4$,
$\det\mat{4}{4}{-3}{0}>0$.
So $\SCD(R)$ consists of the three $\OCD$s above.
\myskip

If we reflect $R$ with respect to the $x$-axis, the new cloud $\bar R$ of the points $p_1,p_2,\bar p_3=(0,-3)$ has $\SCD(\bar R)=\mSCD(R)$ with
$$\OCD(\bar R;p_1)=\OCD(R),
\OCD(\bar R;p_2)=\left[4,\left( \begin{array}{cc} 
4 & 5 \\
0 & 3 \\
0 & +
\end{array}\right) \right],
\OCD(R;\bar p_3)=\left[3,\left( \begin{array}{cc} 
3 & 5 \\
0 & 4 \\
0 & -
\end{array}\right) \right]$$ whose signs changed under reflection, so $\SCD(R)\neq\SCD(\bar R)$.
\medskip

\noindent
\textbf{(b)}
Let $S\subset\R^2$ consist of $m=4$ points $(\pm 1,0),(0,\pm 1)$ that are vertices of the square in Fig.~\ref{fig:clouds_oriented}~(4th).
The centre of mass is $0\in\R^2$ and has a distance $1$ to each point of $S$.
\smallskip

For each 1-point subset $A=\{p\}\subset S$, the distance matrix $D(A\cup\{0\})$ on two points is the single number $1$.
The matrix $M(S;A\cup\{0\})$ has $m-n+1=3$ columns.
For $p_1=(1,0)$, we have   
$M\left(S;\vect{p_1}{0}\right)=\left(\begin{array}{ccc} 
\sqrt{2} & \sqrt{2} & 2 \\
1 & 1 & 1 \\
- & + & 0
\end{array}\right)$, where the columns are ordered according to
$p_2=(0,-1)$, $p_3=(0,1)$, $p_4=(-1,0)$ in Fig.~\ref{fig:clouds_oriented}~(4th).
The sign in the bottom right corner is 0 because the points $p_1,0,p_4$ are in a straight line. 
Due to the rotational symmetry, $M(S;\{p_i,0\})$ is independent of $i=1,2,3,4$.
So $\SCD(S)$ can be considered as one $\OCD=\left[1,M\left(S;\vect{p_1}{0}\right)\right]$ of weight 1. 
\bs 
\end{exa}

\begin{thm}[completeness of $\SCD$, {\cite[Theorem~3.10]{kurlin2023strength}}]
\label{thm:SCD_complete}
The Simplexwise Centred Distribution $\SCD(C)$ in Definition~\ref{dfn:SCD} is a complete isometry invariant of clouds $C\subset\R^n$ of $m$ unlabelled points with a centre of mass at the origin $0\in\R^n$, and can be computed in time $O(m^n/(n-4)!)$.
So any clouds $C,C'\subset\R^n$ are related by rigid motion (isometry, respectively) \emph{if and only if} 
$\SCD(C)=\SCD(C')$ ($\SCD(C)$ equals $\SCD(C')$ or its mirror image $\mSCD(C')$, respectively).
For any $m$-point clouds $C,C'\subset\R^n$, let 
$\SCD(C)$ and $\SCD(C')$ consist of $k=\binom{m}{n-1}$ $\OCD$s.
\ethm
\end{thm}

Corollary~\ref{cor:OCD_reconstruction}
follows from Lemma~\ref{lem:ORD_reconstruction} by adding the centre of mass of $C$ as an extra point to a base sequence $A\subset C$.

\begin{cor}[reconstruction from $\OCD$]
\label{cor:OCD_reconstruction}
A cloud $C\subset\R^n$ of $m>n$ unordered points with the centre of mass $O(A)$ can be reconstructed, uniquely under rigid motion, from $\OCD(C;A)$ in Definition~\ref{dfn:OSD} for any base sequence $A\subseteq C$ with $\aff(\{O(A)\}\cup A)=n-1$.  
\elem
\end{cor}

Example~\ref{exa:SCD}(b) illustrates the key discontinuity challenge:
if $p_4=(-1,0)$ is perturbed, the corresponding sign can discontinuously change to $+1$ or $-1$.
\myskip

To get a continuous metric on $\OCD$s, we will multiply each sign by 
a continuous \emph{strength} function, which vanishes for any zero sign, as defined in the next section.
 
\section{The Lipschitz continuous strength of a simplex in $\R^n$}
\label{sec:strength}

This section resolves the discontinuity of signs of determinants by introducing the multiplicative factor below.

\index{strength of a simplex}

\begin{dfn}[\emph{strength} $\si(A)$ of a simplex]
\label{dfn:strength_simplex}
For a set $A$ of $n+1$ points $q=p_0,p_1,\dots,p_n$ in $\R^n$,
let $p(A)=\frac{1}{2}\sum\limits_{i\neq j}^{n+1}|p_i-p_j|$ be half of the sum of all pairwise distances.
Let $V(A)$ denote the volume the $n$-dimensional simplex on the set $A$.
Define the \emph{strength} $\si(A)=V^2(A)/p^{2n-1}(A)$.
\edfn
\end{dfn}

\begin{exa}[strengths in dimensions $1,2$]
\label{exa:strength}
\tb{(a)}
For $n=1$ and a set $A={p_0,p_1}\subset\R$, the volume is $V(A)=|p_0-p_1|=2p(A)$, so $\si(A)=
2|p_0-p_1|$ is the double length.
\myskip

\nt
\tb{(b)}
For $n=2$ and a triangle $A\subset\R^2$ with sides $a,b,c$, Heron's formula gives $\si(A)=\dfrac{(p-a)(p-b)(p-c)}{p^2}$, $p=\dfrac{a+b+c}{2}=p(A)$ is the half-perimeter of $A$.
\eexa
\end{exa}
 
The strength $\si(A)$ depends only on the distance matrix $D(A)$ from Definition~\ref{dfn:D(A)+M(C;A)}, so the notation $\si(A)$ is used only for brevity.
In any $\R^n$, the squared volume $V^2(A)$ is expressed by the Cayley-Menger determinant \cite{sippl1986cayley} in pairwise distances between points of $A$.
The strength $\si(A)$ vanishes when the simplex on a set $A$ degenerates. 
\smallskip

Theorem~\ref{thm:OSD+SCD_continuous} will need the continuity of $s\si(A)$, when a sign $s\in\{\pm 1\}$ from a bottom row of ORD changes while passing through a degenerate set $A$.
In appendices, the proof of the continuity of $\si(A)$ in Theorem~\ref{thm:strength} gives an explicit upper bound for a Lipschitz constant $\la_n$ below. 

\index{strength of a simplex}
\index{Lipschitz continuity}

\begin{thm}[Lipschitz continuity of the strength $\si$, {\cite[Theorem~4.2]{kurlin2023strength}}]
\label{thm:strength}
Let a cloud $A'$ be obtained from another $(n+1)$-point cloud $A\subset\R^n$ by perturbing every point within its $\ep$-neighbourhood.
The strength $\si(A)$ from Definition~\ref{dfn:strength_simplex} is Lipschitz continuous so that $|\si(A')-\si(A)|\leq 2\ep \la_n$ for a Lipschitz constant $\la_n$.
\ethm
\end{thm}

\begin{exa}[approximates constants $\la_n$ of strength]
\label{exa:strength_bounds}
For $n\geq 2$, the proof of \cite[Theorem~4.2]{kurlin2023strength} implies the following approximate values for upper bounds of the Lipschitz constant of strength: $\la_2=2\sqrt{3}$, $\la_3\approx 0.43$, $\la_4\approx 0.01$, which quickly tend to 0 due to the `curse of dimensionality'. 
The plots in \cite[Fig.~4]{widdowson2023recognizing} illustrate that the strength $\si(A)$ behaves smoothly under perturbations and the  derivative $|\pd{\si}{x}|$ is much smaller than the proved bounds of $\la_n$ above. 
\eexa
\end{exa}

\section{Algorithms for continuous metrics on complete invariants}
\label{sec:algorithms_metrics}

This section introduces Lipschitz continuous metrics on the invariants $\OSD$ and $\SCD$ by using the strength of a simplex.
By Definition~\ref{dfn:OSD} an Oriented Relative Distribution $\ORD$
 is a pair $[D(A);M(C;A)]$ of matrices considered up to permutations $\xi\in S_n$ of $n$ ordered points of $A$.
Any column of $M(C;A)$ is a pair $(v,s)$, where $s\in\{\pm 1,0\}$ and $v\in\R^n$ is a vector of distances from $q\in C-A$ to $p_1,\dots,p_n\in A$.
\smallskip

For simplicity and similar to the case of a general metric space, we assume that a cloud $C\subset\R^n$ is given by a matrix of pairwise Euclidean distances.
If $C$ is given by Euclidean coordinates of points, then any distance requires $O(n)$ computations, and we should add the factor $n$ in all complexities below, keeping all times polynomial in $m$.
The $m-n$ permutable columns of the matrix $M(C;A)$ in $\ORD$ from Definition~\ref{dfn:OSD} can be interpreted as $m-n$ unordered points in $\R^n$.
Since any isometry is bijective, the simplest metric respecting bijections is the bottleneck distance $\BD$ from Example~\ref{dfn:metrics}(b). 

\begin{dfn}[max metric on $\ORD$s]
\label{dfn:ORD+metric}
Consider the bottleneck distance $\BD$ on the set of $m-n$ permutable columns of $M(C;A)$ as on a cloud of $m-n$ unordered points $(v,\frac{s}{\la_n}\si(A\cup\{q\}))\in\R^{n+1}$.
For another $\ORD'=[D(A');M(C';A')]$ and any permutation $\xi\in S_n$ of indices $1,\dots,n$ acting on $D(A)$ and rows of $M(C;A)$, set 
$$d_o(\xi)=\max\{L_{\infty}(\xi(D(A)),D(A')),\BD(\xi(M(C;A)),M(C';A')) \}.$$
Then the \emph{max metric} is defined as $M_\infty(\ORD,\ORD')=\min\limits_{\xi\in S_n} d_o(\xi)$.
\edfn
\end{dfn}

The coefficient $\dfrac{1}{\la_n}$ 
in front of the strength $\si(A\cup\{q\})$ in Definitions~\ref{dfn:ORD+metric} and~\ref{dfn:OCD+metric}
normalises the Lipschitz constant $\la_n$ of $\si$ to $1$ in line with changes of distances by at most $2\ep$ when points are perturbed within their $\ep$-neighbourhoods.

\index{max metric}
\index{strength of a simplex}

\begin{dfn}[max metric on $\OCD$s]
\label{dfn:OCD+metric}
Consider the bottleneck distance $\BD$ on the set of permutable $m-n+1$ columns of $M(C;A\cup\{0\})$ as on a cloud of $m-n+1$ unordered points $\left(v,\dfrac{s}{\la_n}\si(A\cup\{0,q\})\right)\in\R^{n+1}$.
For another $\OCD'=[D(A'\cup\{0\});M(C';A'\cup\{0\})]$ and any permutation $\xi\in S_{n-1}$ of indices $1,\dots,n-1$ acting on $D(A\cup\{0\})$ and the first $n-1$ rows of $M(C;A\cup\{0\})$, set $d_o(\xi)=\max\{L,W\}$, 
$$\text{ where }L=L_{\infty}\Big(\xi(D(A\cup\{0\})),D(A'\cup\{0\})\Big),$$
$$W=\BD\Big(\xi(M(C;A\cup\{0\})),M(C';A'\cup\{0\})\Big).$$
Then the max metric is defined as $M_\infty(\OCD,\OCD')=\min\limits_{\xi\in S_{n-1}} d_o(\xi)$.
\edfn
\end{dfn}

The max metric $M_{\infty}$ is used for intermediate costs to get metrics on unordered collections $\OSD$s and $\SCD$s by using the metrics LAC and $\EMD$ from Definitions~\ref{dfn:LAC_matrix}and~\ref{dfn:EMD}, respectively.  
Equality $\OSD(C)=\OSD(C')$ between unordered collections of ORDs is best verified by checking if LAC or EMD  between these $\OSD$s is 0.

\begin{thm}[times for metrics on $\OSD$s, {\cite[Theorem~5.6]{kurlin2023strength}}]
\label{thm:OSD_metrics}
\textbf{(a)}
For the $k\times k$ matrix of costs computed by the max metric $M_{\infty}$ between $\ORD$s from $\OSD(C)$ and $\OSD(C')$, 
$\LAC$ from Definition~\ref{dfn:LAC_matrix} satisfies all metric axioms on $\OSD$s and needs time $O(n!(n^2+m^{1.5}\log^{n+1} m)k^2+ k^3\log k)$.
\medskip

\noindent
\textbf{(b)}
Let $\OSD$s have a maximum size $l\leq k$ after collapsing identical $\ORD$s. The $\EMD$ from Definition~\ref{dfn:EMD} satisfies all metric axioms  on $\OSD$s and can be computed in time $O(n!(n^2 +m^{1.5}\log^{n+1} m) l^2 +l^3\log l)$.
\ethm
\end{thm}

Equality $\SCD(C)=\SCD(C')$ is interpreted as a bijection between unordered sets $\SCD(C)\to\SCD(C')$ matching all $\OCD$s, which is best verified by checking if the metrics between these $\SCD$s vanish in Theorem~\ref{thm:SCD_metrics}.

\begin{thm}[times for metrics on $\SCD$s, {\cite[Theorem~5.7]{kurlin2023strength}}]
\label{thm:SCD_metrics}
\textbf{(a)}
For the $k\times k$ matrix of costs computed by the max metric $M_{\infty}$ between $\OCD$s in $\SCD(C)$ and $\SCD(C')$, the metric $\LAC$ from Definition~\ref{dfn:LAC_matrix} satisfies all metric axioms on $\SCD$s and needs time $O((n-1)!(n^2+m^{1.5}\log^{n} m)k^2+ k^3\log k)$.
\medskip

\noindent
\textbf{(b)}
Let $\SCD$s have a maximum size $l\leq k$ after collapsing identical $\OCD$s. Then $\EMD$ from Definition~\ref{dfn:EMD} satisfies all metric axioms  on $\SCD$s and can be computed in time $O((n-1)!(n^2 +m^{1.5}\log^{n} m) l^2 +l^3\log l)$.
\ethm
\end{thm}

If we estimate $l\leq k=\binom{m}{n-1}=m(m-1)\dots(m-n+2)/n!$ as $O(m^{n-1}/n!)$, Theorem~\ref{thm:SCD_metrics} gives time 
$O(n(m^{n-1}/n!)^3\log m)$ for metrics on $\SCD$s, which is $O(m^3\log m)$ for $n=2$, and $O(m^6\log m)$ for $n=3$. 
Though the above estimates are very rough, 
the time $O(m^3\log m)$ in $\R^2$ is faster than the only past time $O(m^5\log m)$ for comparing $m$-point clouds by the Hausdorff distance minimized over isometries \cite{chew1997geometric}.

\index{strength of a simplex}
\index{Lipschitz continuity}

\begin{thm}[continuity of $\OSD$ and $\SCD$, {\cite[Corollary~6.1]{kurlin2023strength}}]
\label{thm:OSD+SCD_continuous}
For any cloud $C\subset\R^n$ of $m$ unordered points, perturbing any point within its $\ep$-neighbourhood changes $\OSD(C)$ and $\SCD(C)$ by at most $2\ep$ in the metrics $\LAC$ and $\EMD$.
\ethm
\end{thm}

Theorems~\ref{thm:SCD_complete}, \ref{thm:SCD_metrics}, and \ref{thm:OSD+SCD_continuous} imply that the Simplexwise Centred Distribution fully solves Problem~\ref{pro:Euclidean_unordered}.
In addition, Corollary~\ref{cor:OCD_reconstruction} proves that $\OCD(C;A)$ from $\SCD(C)$ suffices to reconstruct $C\subset\R^n$, uniquely under rigid motion, if $C$ has a base sequence $A$ of $n-2$ points with $\aff(\{O(A)\}\cup A)=n-1$.
\myskip

Forthcoming work will improve $\SCD$s to better invariants that allow a reconstruction in all degenerate cases and satisfy the much harder realisability condition as in \ref{pro:geocodes}(f).

\vspace*{-4mm}
 
\bibliographystyle{plain}
\bibliography{Geometric-Data-Science-book}

%
%
%

\begin{partbacktext}
\part{Geometric Data Science of periodic point sets}
\end{partbacktext}

%
%
%

\chapter{One-periodic sequences in high-dimensional Euclidean spaces}
\label{chap:1-periodic} 

\abstract{
This chapter studies high-dimensional data that is periodic in one direction.
These periodic sequences live in a high-dimensional space $\R\times\R^{n-1}$ for any dimension $n\geq 1$ and were indistinguishable by past invariants even in dimension $n=2$.
Experimental noise and atomic vibrations motivate a new continuous approach, because a minimal periodic pattern breaks down under almost any perturbation.
}

\section{One-periodic sequences under various isometries in $\R\times\R^{n-1}$}
\label{sec:1-periodic_isometries}

All sections in this chapter follow paper \cite{kurlin2025complete} with minor updates.

\index{1-periodic sequence}
\begin{dfn}[1-periodic sequences in $\R\times\R^{n-1}$]
\label{dfn:1-periodic}
Let $\vec e_1$ be the unit vector along the first axis in $\R\times\R^{n-1}$ for $n\geq 1$.
For a \emph{period} $l>0$, a \emph{motif} $M$ is a set of points $p_1,\dots,p_m$ in the \emph{slice} $[0,l)\times\R^{n-1}$ of the width $l>0$.
We assume that the \emph{time projections} $t(p_1),\dots,t(p_m)$ under $t:[0,l)\times\R^{n-1}\to[0,l)$ are distinct, while $v(p_1),\dots,v(p_m)$ under the \emph{value projection} $v:[0,l)\times\R^{n-1}\to\R^{n-1}$ are arbitrary.
\sskip

A \emph{1-periodic sequence} $S=M+l\vec e_1\Z$ is the infinite sequence of points $p(i+mj)=p_i+jl\vec e_1\in\R^n$, which are indexed by $i+mj$, where $j\in\Z$ and $i=1,\dots,m$.
\edfn
\end{dfn}

\begin{figure}[h]
\includegraphics[width=\textwidth]{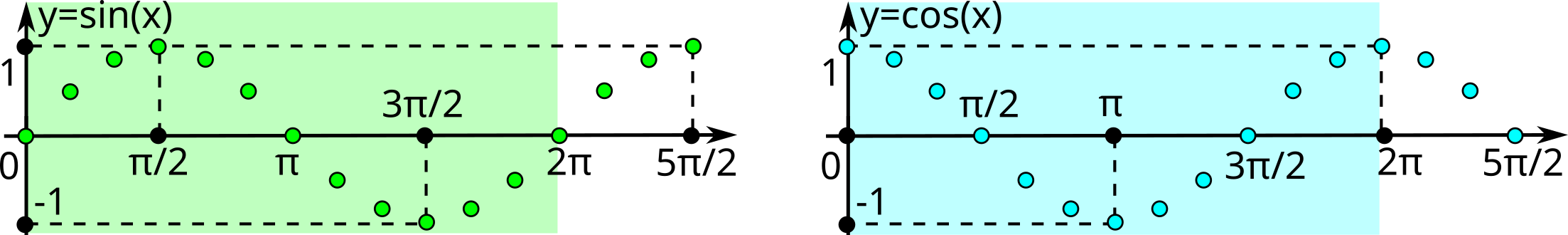}
\caption{The periodic sequences $C,S\subset\R\times\R$ are sampled from the sine and cosine graphs.
The motifs in the shaded slice $[0,2\pi)\times\R$ are non-isometric, but $S,C$ are related by translation.}
\label{fig:sampled-waves}
\end{figure}

The slice $[0,l)\times\R^{n-1}$ excludes all points with $t=l$, which are equivalent to points with $t=0$ by translation in the time factor $\R$.
Then all points $p_1,\dots,p_m\in[0,l)\times\R^{n-1}$ are counted once and ordered under the time projection $t:[0,l)\times\R^{n-1}\to[0,l)$.

\begin{exa}[1-periodic sequences in $\R\times\R$]
\label{exa:sampled_waves}
Fig.~\ref{fig:sampled-waves}~(left)  shows the 1-periodic sequence $S$ in $\R\times\R$ (from the sine graph) with the period $l=2\pi$ and motif $M_S$ of 

$$(0,0),\quad 
(\frac{\pi}{6},\frac{1}{2}),\quad 
(\frac{\pi}{3},\frac{\sqrt{3}}{2}),\quad (\frac{\pi}{2},1),\quad 
(\frac{2\pi}{3},\frac{\sqrt{3}}{2}),\quad (\frac{5\pi}{6},\frac{1}{2}),$$
$$(\pi,0),\quad
(\frac{7\pi}{6},-\frac{1}{2}),\quad
(\frac{4\pi}{3},-\frac{\sqrt{3}}{2}),\quad 
(\frac{3\pi}{2},-1),(\frac{5\pi}{3},\quad
\frac{\sqrt{3}}{2}),\quad
(\frac{11\pi}{6},-\frac{1}{2}).$$
Fig.~\ref{fig:sampled-waves}~(right) shows another sequence $C$ with the same period $l=2\pi$ and a different motif $M_C\neq M_S$.
However, $S,C$ are identical under translation: $\sin(x+\frac{\pi}{2})=\cos(x)$.
\eexa
\end{exa}

Example~\ref{exa:sampled_waves} illustrates the ambiguity of digital representations when many real objects look different in various coordinate systems despite being equivalent as rigid objects.
We adapt basic equivalences to sets in the product $\R\times\R^{n-1}$. 

\index{cyclic equivalence}
\index{dihedral equivalence}

\begin{dfn}[\emph{cyclic} vs \emph{dihedral} isometries and rigid motions]
\label{dfn:equivalences}
A \emph{cyclic} isometry of $\R\times\R^{n-1}$ is a composition of a translation in the time factor $\R$ and an isometry in the value factor $\R^{n-1}$.
If we allow compositions of a translation and symmetry $x\mapsto -x$ in the time factor $\R$, the resulting isometry of $\R\times\R^{n-1}$ is \emph{dihedral}.
\sskip

If we allow only isometries that preserve orientation in the value factor $\R^{n-1}$, the resulting equivalences are called cyclic and dihedral \emph{rigid motions}, 
respectively.
\edfn
\end{dfn}

The adjectives \emph{cyclic} and \emph{dihedral} are motivated by the names of the cyclic group $C_m$ and the dihedral group $D_m$ consisting of orientation-preserving isometries and all isometries in $\R^2$, respectively, that map the regular polygon on $m$ vertices to itself.
\medskip

The equivalences in Definition~\ref{dfn:equivalences} make sense for any finite sequence of points $T\subset\R\times\R^{n-1}$.
However, the periodicity substantially worsens the ambiguity of representations based a period $l$ and a motif $M$ as follows.
A translation in the time factor $\R$ allows us to fix any point $p$ of a motif $M$ at $t=0$, but this choice of $p$ is arbitrary, so a motif $M$ is defined only modulo cyclic permutations of its points.  
\myskip

The set of integers can be defined as $\Z$ with period $1$ or as $\{0,1\}+2\Z$ with period 2, or with any integer period $l>0$.
For any sequence $S=\{p_1,\dots,p_m\}+l\vec e_1\Z$, we can choose a \emph{minimal} period $l$ such that $S$ can not be represented with a smaller period.
\myskip

This classical approach in crystallography leads to an \emph{invariant} $I$ 
based on a minimum period (primitive cell) and defined as a set of numerical properties preserved under any rigid motion. 
Fixing a minimum period creates the following discontinuity.

\begin{exa}[discontinuity of a period]
\label{exa:discontinuous_period1}
For any small $\ep>0$ and integer $m\geq 1$, any point of $\Z$ is $\ep$-close to a unique point of the sequence $\{0,1+\ep,\dots, m+\ep\}+(m+1)\Z$, though their minimum periods $1$ and $m+1$ are arbitrarily different.
Hence comparing periodic sequences by their minimal motifs can miss infinitely many near-duplicates. 
\eexa
\end{exa}

We assume that the input for a 1-periodic sequence $S$ consists of a period $l>0$ and a motif of $m=|S|$ points in the high-dimensional slice $[0,l)\times\R^{n-1}$.

\begin{pro}[invariants of 1-periodic sequences in $\R\times\R^{n-1}$]
\label{pro:1-periodic}
Design an invariant $I$ of all 1-periodic sequences of points in $\R\times\R^{n-1}$ satisfying the following conditions.
\myskip

\noindent
\tb{(a)} 
\emph{Completeness:}
any 1-periodic sequences $S,Q\subset\R\times\R^{n-1}$ are related by cyclic isometry (denoted as $S\cong Q$) in Definition~\ref{dfn:equivalences} if and only if $I(S)=I(Q)$.
\myskip

\noindent
\tb{(b)} 
\emph{Reconstruction:}
any 1-periodic sequence $S\subset\R\times\R^{n-1}$ is reconstructable from its invariant value $I(S)$, uniquely under cyclic isometry. 
\myskip

\noindent
\tb{(c)} 
\emph{Metric:} 
there is a distance $d$ on the space $\{I(S) \vl \text{all 1-periodic sequences }S\subset\R\times\R^{n-1}\}$ satisfying all metric axioms in Definition~\ref{dfn:metrics}(a). 
\myskip

\noindent
\tb{(d)} 
\emph{Continuity:} 
there is a constant $\la>0$, such that, for all sufficiently small $\ep$, if a 1-periodic sequence $Q$ is obtained by perturbing every point of a 1-periodic sequence $S\subset\R\times\R^n$ up to Euclidean distance $\ep$, then $d(I(S),I(Q))\leq\la\ep$.
\myskip

\noindent
\tb{(e)} 
\emph{Computability:} the invariant $I$, a reconstruction of $S\subset\R\times\R^{n-1}$ from $I(S)$, and the metric $d(I(S),I(Q))$ can be computed in times that depend polynomially on the dimension $n$ and the maximum motif size of 1-periodic sequences $S,Q$.
\epro
\end{pro}
 
Further sections will develop invariants that solve Problem~\ref{pro:1-periodic} 
for all 1-periodic sequences under cyclic and dihedral isometries and rigid motions in $\R\times\R^{n-1}$. 
\myskip

Our constructions were inspired by the infinite family of 1-periodic sequences in \cite[Fig.~4]{pozdnyakov2022incompleteness} that were not distinguished by past invariants, see a review in \cite[section~2]{kurlin2025complete}. 

\section{Invariants and continuous metrics for finite sequences in $\R^n$}
\label{sec:Euclidean_finite}

This section studies complete invariants and metrics for isometry classes of finite sequences of ordered points in $\R^n$.
These invariants are easily extendable to the 1-periodic sequences in $\R\times\R^{n-1}$ and substantially differ from direction-based invariants in chapter~\ref{chap:directional} and the backbone invariants in chapter~\ref{chap:proteins}, which were defined only for non-degenerate sequences of triplets of points in $\R^3$.

\index{cyclic distance matrix}

\begin{dfn}[distance matrices $\DM$ and $\CDM$]
\label{dfn:CDM}
Let $T$ be an ordered sequence of $m$ points $p_1,\dots,p_m\in\R^{n}$.
\myskip

\nt
\tb{(a)}
In the \emph{distance matrix} $\DM(T)$ of the size $m\times m$, each element $\DM_{ij}(T)$ is the Euclidean distance $|p_j-p_{j}|$ for $i,j\in\{1,\dots,m\}$, so $d_{ii}=0$ for $i=1,\dots,m$.
\myskip

\nt
\tb{(b)}
To define the \emph{cyclic distance matrix} $\CDM(T)$ of the size $(m-1)\times m$, set element $\CDM_{ij}(T)$ to the Euclidean distance $|p_j-p_{i+j}|$ for $i\in\{1,\dots,m-1\}$ and $j\in\{1,\dots,m\}$, where all indices are considered modulo $m$, for example, $p_{m+1}=p_1$.
\edfn 
\end{dfn}

Any $m=3$ points in $\R^n$ with pairwise distances $d_{ij}$ have the distance matrix
$\DM=\left(\begin{array}{ccc}
0 & d_{12} & d_{13} \\
d_{12} & 0 & d_{23} \\
d_{13} & d_{23} & 0 
\end{array} \right)$ and 
the cyclic distance matrix
$\CDM=\left(\begin{array}{ccc}
d_{12} & d_{23} & d_{13} \\
d_{13} & d_{12} & d_{23} 
\end{array} \right)$.
$\CDM(T)$ is obtained from $\DM(T)$ by removing the zero diagonal and cyclically shifting each column so that the first row of $\CDM(T)$ has distances from $p_i$ to the next point $p_{i+1}$.

\begin{figure}[h]
\centering
\includegraphics[width=\textwidth]{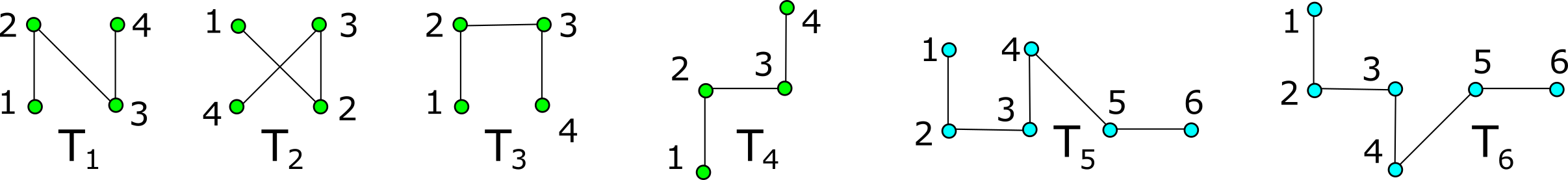}
\caption{These sequences are distinguished by their cyclic distance matrices in Example~\ref{exa:CDM}.}
\label{fig:finite_sequences}
\end{figure}

\begin{exa}[cyclic distance matrices]
\label{exa:CDM}
Fig.~\ref{fig:finite_sequences} shows the sequences $T_1,\dots,T_6\subset\R^2$ whose points are in the integer lattice $\Z^2$ so that the minimum inter-point distance is 1.
In each sequence, the points are connected by straight lines in the order $1\to 2\to\dots\to m$.
$\CDM(T_1)=\left(\begin{array}{cccc}
1 & \sqrt{2} & 1 & \sqrt{2} \\
1 & 1 & 1 & 1\\
\sqrt{2} & 1 & \sqrt{2} & 1
\end{array} \right)$,
$\CDM(T_2)=\left(\begin{array}{cccc}
\sqrt{2} & 1 & \sqrt{2} & 1 \\
1 & 1 & 1 & 1\\
1 & \sqrt{2} & 1 & \sqrt{2}
\end{array} \right)$ are different but related by a cyclic shift of columns. 
This shift of indices in $T_1$ gives a sequence isometric to $T_2$. 
Then
$\CDM(T_3)=\left(\begin{array}{cccc}
1 & 1 & 1 & 1 \\
\sqrt{2} & \sqrt{2} & \sqrt{2} & \sqrt{2}\\
1 & 1 & 1 & 1
\end{array} \right)$,
$\CDM(T_4)=\left(\begin{array}{cccc}
1 & 1 & 1 &  \sqrt{5}\\
\sqrt{2} & \sqrt{2} & \sqrt{2} & \sqrt{2}\\
\sqrt{5} & 1 & 1 & 1
\end{array} \right)$.
The CDMs of the sets $T_5,T_6$ differ only by distances $|p_1-p_4|=1$ in $T_5$ and $|p_1-p_4|=\sqrt{5}$ in the highlighted cells below.
If reduce the number $m-1$ of rows in $\CDM$ to the dimension $n=2$, the smaller matrices fail to distinguish the non-isometric sequences $T_5\not\cong T_6$.

$T_5:\; \left(\begin{array}{cccccc}
1 & 1 & 1 &  \sqrt{2} & 1 & \sqrt{10}\\
\sqrt{2} & \sqrt{2} & 1 & \sqrt{5} & \sqrt{5} & 3 \\
\mc{yellow}{1} & 2 & 2 & \mc{yellow}{1} & 2 & 2 \\
\sqrt{5} & 3 & \sqrt{2} & \sqrt{2} & 1 & \sqrt{5} \\
\sqrt{10} & 1 & 1 & 1 & \sqrt{2} & 1
\end{array} \right)$ and
$T_6:\; \left(\begin{array}{cccccc}
1 & 1 & 1 &  \sqrt{2} & 1 & \sqrt{10}\\
\sqrt{2} & \sqrt{2} & 1 & \sqrt{5} & \sqrt{5} & 3 \\
\mc{yellow}{\sqrt{5}} & 2 & 2 & \mc{yellow}{\sqrt{5}} & 2 & 2 \\
\sqrt{5} & 3 & \sqrt{2} & \sqrt{2} & 1 & \sqrt{5} \\
\sqrt{10} & 1 & 1 & 1 & \sqrt{2} & 1
\end{array} \right)$.
\eexa
\end{exa}

Recall that Definition~\ref{dfn:strength_simplex} introduced the strength $\si(A)=\dfrac{V^2(A)}{p^{2n-1}(A)}$ of a simplex $A$ on any set of $n+1$ points $q_0,q_1,\dots,q_n\in\R^n$, where $V(A)$ is the volume of $A$ and $p(A)=\dfrac{1}{2}\sum\limits_{0\leq i<j\leq n}|q_i-q_j|$ is the half-perimeter.

\begin{dfn}[\emph{cyclic distances with signs} $\CDS$]
\label{dfn:CDS}
For any sequence $T$ of $p_1,\dots,p_m\in\R^n$ and $i=1,\dots,m$, let $\si_i(T)$ be the strength of the simplex on the points $p_i,\dots,p_{i+n}$, where all indices are modulo $m$.
Let $\sign_i(T)$ be the sign ($\pm 1$ or $0$) of the $n\times n$ determinant with the columns $\vec p_{i+1}-\vec p_i,\vec p_{i+2}-\vec p_{i+1},\dots,\vec p_{i+n}-\vec p_{i+n-1}$.
The matrix $\CDS(T)$ of \emph{cyclic distances with signs} is obtained from $\CDM(T)$ in Definition~\ref{dfn:CDM} by attaching the extra $m$-th row $\sign(T)=(\sign_1(T),\dots,\sign_m(T))$.
\edfn 
\end{dfn}

\begin{exa}[strengths and signs]
\label{exa:signs+strengths}
For the first sequence $T_1$ in Fig.~\ref{fig:finite_sequences} with the points $p_1=(0,0)$, $p_2=(0,1)$, $p_3=(1,0)$, $p_4=(1,1)$, the first $2\times 2$ determinant with the columns $p_2-p_1=\vect{0}{1}$ and $p_3-p_2=(1,-1)$ is $\det\mat{0}{1}{1}{-1}$ has sign $-1$.
The further determinants for $i=2,3,4$ are $\det\mat{1}{1}{-1}{0}=+1$, $\det\mat{1}{-1}{0}{-1}=-1$, $\det\mat{-1}{1}{-1}{0}=+1$, so the row of signs is 
$\sign(T_1)=(-1,+1,-1,+1)$.
All triangles on 4 triples $p_i,p_{i+1},p_{i+2}$ for $i=1,2,3,4$ (with indices considered modulo 4) have the sides $1,1,\sqrt{2}$, half-perimeter $p=1+\dfrac{1}{\sqrt{2}}$, area $V=\dfrac{1}{2}$, and strength $\si=\dfrac{1}{\sqrt{2}(1+\sqrt{2})^3}$. 
\eexa
\end{exa}

Since the sign of a determinant discontinuously changes when a point set passes through a degenerate configuration, this sign will be multiplied by the Lipschitz continuous strength to get a metric satisfying condition~\ref{pro:1-periodic}(d), see  Theorem~\ref{thm:Euclidean_ordered}(d).
\smallskip

Section~\ref{sec:1-periodic_metrics} will adapt the matrices from Definitions~\ref{dfn:CDM} and~\ref{dfn:CDS} to 1-periodic sequences whose motifs of points should be considered under cyclic permutations.
The \emph{cyclic} group $C_m$ consists of $m$ permutations on $1,\dots,m$ generated by the \emph{shift permutation} $\ga_m:(1,2,\dots,m)\mapsto(2,\dots,m,1)$. 
The \emph{dihedral} group $D_m$ consists of $2m$ permutations generated by $\ga_m$ and the \emph{reverse permutation} $\iota_m:(1,2,\dots,m)\mapsto(m,\dots,2,1)$. 

\begin{lem}[actions on vectors and matrices]
\label{lem:actions_matrices}
The shift permutation $\ga_m\in C_m$ 
acts on the cyclic distance matrix $\CDM(T)$ by cyclically shifting its $m$ columns and keeping all rows.
The reverse permutation $\iota_m\in D_m$ reverses the order of columns and rows in $\CDM(T)$.
These permutations act on the row of signs in Definition~\ref{dfn:CDS} 
as $\ga_m(s_1,s_2\dots,s_m)=(s_2,\dots,s_m,s_1)$ and $\iota_m(s_1,s_2\dots,s_m)=(-1)^{[3n/2]}(s_m,\dots,s_2,s_1)$.
For any mirror image $\bar T$ of $T$, the matrix $\CDS(\bar T)$ is obtained from $\CDS(T)$ by reversing all signs in the last row.
Any element of the groups $C_m,D_m$ acts on any sequence of $m$ numbers as a composition of $\ga_m,\iota_m$.
\elem
\end{lem}

Any matrix $k\times m$ can be rewritten row-by-row as a vector $v\in\R^{km}$.
For any $q\in[1,+\infty]$, the Minkowski norm is $||v||_q=\left(\sum\limits_{i=1}^{km}|v_i|^q\right)^{1/q}$, where the limit case is $||v||_{\infty}=\max\limits_{i=1,\dots,km}|v_i|$.
Any power $a^{1/q}$ for $a>0$ is interpreted as $1$ in the case $q=+\infty$.

\begin{dfn}[metrics $\MCD_q$ and $\MCS_q$ for finite sequences]
\label{dfn:finite_metrics}
For any Minkowski norm with a parameter $q\in[1,+\infty]$ and ordered sequences $T,S\subset\R^{n-1}$ of $m$ points,
use the matrices from Definition~\ref{dfn:CDM}  
to define the \emph{metric based on cyclic distances} $\MCD_q(S,T)=\dfrac{||\CDM(S)-\CDM(T)||_q}{\big(m(m-1)\big)^{1/q}}$ 
and the \emph{metric based on cyclic distances with signs}
$\MCS_q(S,T)=\max\Big\{\MCD_q(S,T),\dfrac{2}{\la_n}\max\limits_{i=1,\dots,m}\big|\sign_i(S)\si_i(S)-\sign_i(T)\si_i(T) \big| \Big\}.$
\edfn
\end{dfn}

\begin{exa}[metric $\MCD_{q}$]
\label{exa:MCD}
For any $q\in[1,+\infty)$, we use cyclic distance matrices from Example~\ref{exa:CDM}
to compute
$$\begin{array}{l}
\MCD_q(T_1,T_3)=(\frac{2}{3})^{1/q}(\sqrt{2}-1),\\
\MCD_q(T_3,T_4)=(\frac{1}{6})^{1/q}(\sqrt{5}-1), \\
\MCD_q(T_1,T_4)=(\frac{1}{2}(\sqrt{2}-1)^q+\frac{1}{6}(\sqrt{5}-\sqrt{2})^q)^{1/q}.
\end{array}$$
The triangle inequality holds for $q\geq 1$ as follows:
$$\begin{array}{l}\big(\MCD_q(T_1,T_3)+ \MCD_q(T_3,T_4)\big)^q= \\
\big((\frac{2}{3})^{1/q}(\sqrt{2}-1)+(\frac{1}{6})^{1/q}(\sqrt{5}-1)\big)^q\geq \\
\big((\frac{1}{2})^{1/q}(\sqrt{2}-1)+(\frac{1}{6})^{1/q}(\sqrt{5}-\sqrt{2})\big)^q\geq \\
\frac{1}{2}(\sqrt{2}-1)^q+\frac{1}{6}(\sqrt{5}-\sqrt{2})^q=
\big(\MCD_q(T_1,T_4)\big)^q
\end{array}$$
 due to $(a+b)^q\geq a^q+b^q$ for $a,b>0$ and $q\geq 1$.
For $q=+\infty$, the inequality becomes $(\sqrt{2}-1)+(\sqrt{5}-1)\geq \sqrt{5}-\sqrt{2}$.
Then $T_5\not\cong T_6$ have
$\MCD_q(T_5,T_6)=2^{1/q}(\sqrt{5}-1)$.
\eexa
\end{exa}

We use the extra factors $\big(m(m-1)\big)^{1/q}$ and $\dfrac{2}{\la_n}$ in the definition above, where $\la_n$ is a Lipschitz constant $\si$ from Theorem~\ref{thm:strength}, to guarantee the Lipschitz constant $2$ for the new metrics. 
Indeed, perturbing any points up to $\ep$ changes the distance between them up to $2\ep$.
Instead of maxima in the formula for $\MCS_q(S,T)$, one can consider other metric transforms from \cite[section 4.1]{deza2009encyclopedia}, for example, sums of metrics.
\myskip

Theorem~\ref{thm:Euclidean_ordered} solves Problem~\ref{pro:Euclidean_ordered}   and will help solve Problem~\ref{pro:1-periodic}.

\begin{thm}[solution of Problem~\ref{pro:Euclidean_ordered} for ordered points in $\R^n$,
{\cite[Theorem~3.9]{kurlin2025complete}}]
\label{thm:Euclidean_ordered}
\textbf{(a)}
For any sequence $T\subset\R^n$ of $m$ ordered points, the matrices $\CDM(T)$ and $\CDS(T)$ are complete invariants of $T\subset\R^n$ under isometry and rigid motion in $\R^n$, which are computable in times $O(m^2n)$ and $O(m^2n+mn^3)$, respectively.
\medskip

\noindent
\textbf{(b)}
Any sequence $T\subset\R^n$ of $m$ points can be reconstructed from the invariants $\CDM(T)$ and $\CDS(T)$, uniquely under isometry and rigid motion, respectively, in time $O(m^3)$. 
\medskip

\noindent
\textbf{(c)}
For any sequences $S,T\subset\R^n$ of $m$ points, $\MCD_q(S,T),\MCS_q(S,T)$ satisfy all metric axioms and are computable in time {$O(m^2)$} and $O(m^2n+mn^3)$, respectively.
\medskip

\noindent
\textbf{(d)}
If $S$ is obtained from any sequence $T\subset\R^n$ by perturbing every point up to Euclidean distance $\ep$, then $\MCD_q(S,T)\leq 2\ep$ and $\MCS_q(S,T)\leq 2\ep$ for $q\in[1,+\infty]$.
\ethm
\end{thm}

\section{Discontinuity of a minimal period for 1-periodic sequences}
\label{sec:discontinuity}

The invariants and metrics from section~\ref{sec:Euclidean_finite} will be used for a motif of a 1-periodic sequence $S$ projected to the value factor $\R^{n-1}$.
To solve Problem~\ref{pro:1-periodic}, we first resolve the discontinuity of a period under perturbations by projecting $S$ to the time factor $\R$.

\begin{dfn}[time shift $\TS$] 
\label{dfn:time_shift}
Let $S\subset\R\times\R^{n-1}$ be a 1-periodic sequence with a period $l$ and a motif $M$ of points $p_1,\dots,p_m$, which have ordered time projection $t(p_1)<\cdots<t(p_m)$ in $[0,l)$ under $t:\R\times\R^{n-1}\to\R$, see Definition~\ref{dfn:1-periodic}.
Set $d_i=t(p_{i+1})-t(p_{i})$ for $i=1,\dots,m$, $t(p_{m+1})=t(p_1)+l$.
The \emph{time shift} of the pair (motif, period) of the 1-periodic sequence $S$ is $\TS(M;l)=(d_1,\dots,d_m)$.
\edfn
\end{dfn}

The sequences $S_2=\{0,1\}+3\Z$ and $3-S_2=\{0,2\}+3\Z$ are related by translation but have different time shifts $\TS(\{0,1\};3)=(1,2)$ and $\TS(\{0,2\};3)=(2,1)$.
To get isometry invariants, these shifts are considered modulo cyclic or dihedral permutations.

\begin{dfn}[cyclic and dihedral invariants under isometries]
\label{dfn:periodic_invariants}
\textbf{(a)}
For any 1-periodic sequence $S=M+l\vec e_1\Z\subset\R\times\R^{n-1}$ with a minimum motif $M$ of $m$ points, let $v(M)\subset\R^{n-1}$ be the image of $M$ under the value projection $v:\R\times\R^{n-1}\to\R^{n-1}$.
\medskip

\noindent
\textbf{(b)}
The \emph{cyclic} and \emph{dihedral isometry} invariants $\CI(S),\DI(S)$ are the classes of the pair $(\TS(M;l),\CDM(v(M)))$ under permutations $\ga$ from the groups $C_m,D_m$, respectively, acting simultaneously on the time shift $\TS(M;l)$ and the matrix $\CDM(v(M))$.
\medskip

\noindent
\textbf{(c)}
The \emph{cyclic} and \emph{dihedral rigid} invariants $\CR(S),\DR(S)$ are the classes of the pair $(\TS(M;l),\CDS(v(M)))$ under permutations $\ga$ from the groups $C_m,D_m$, respectively, acting simultaneously on the time shift $\TS(M;l)$ and the matrix $\CDS(v(M))$.
\edfn
\end{dfn}

The matrices $\CDM,\CDS$ are used for the projected motif $v(M)\subset\R^{n-1}$ and do not depend on a period $l$, because a shift along the time direction $\vec e_1$ keeps the value projection.
For $n=1$, when a periodic sequence $S=\{p_1,\dots,p_m\}+l\Z$ is in the line $\R$, Definition~\ref{dfn:periodic_invariants} simplifies to a single time shift obtained by lexicographic ordering. 
\medskip

Recall the \emph{lexicographic order} on vectors: $(d_1,\dots,d_m)<(d'_1,\dots,d'_m)$ if $d_1=d'_1,\dots,d_i=d'_i$ for some $0\leq i<m$, where $i=0$ means no identities, and $d_{i+1}<d'_{i+1}$.

\begin{dfn}[time invariants $\CT,\DT$] 
\label{dfn:time_invariants}
Let $S=\{p_1,\dots,p_m\}+l\Z$ be a 1-periodic sequence in $\R\times\R^{n-1}$ with a minimum period $l>0$.
Set $d_i=p_{i+1}-p_i$ for $i=1,\dots,m$, where $p_{m+1}=p_1+l$.
Apply all permutations of cyclic group $C_m$ to $(d_1,\dots,d_m)$, order all resulting lists lexicographically, and call the first (smallest) list the \emph{cyclic time} invariant $\CT(S)$.
Similarly, define the \emph{dihedral time} invariant $\DT(S)$ as the lexicographically smallest list obtained from $(d_1,\dots,d_m)$ by the action of $D_m$.
\edfn
\end{dfn}

\begin{exa}[time invariants $\CT,\DT$] 
\label{exa:time_invariants}
In $\R$, the periodic sequences $S=\{0,1,3\}+6\Z$ and $Q=6-S=\{0,3,5\}+6\Z$ are related by reflection $x\mapsto 6-x$ and not by translation.
Their time shifts are $\TS(\{0,1,3\};6)=(1,2,3)$ and $\TS(\{0,3,5\};6)=(3,2,1)$.
So the dihedral time invariants of both $S,Q$ are equal to $\DT=(1,2,3)$, but their cyclic time invariants differ: $\CT(S)=(1,2,3)\neq (1,3,2)=\CT(Q)$. 
\eexa
\end{exa}

\begin{figure}[h]
\centering
\includegraphics[height=20mm]{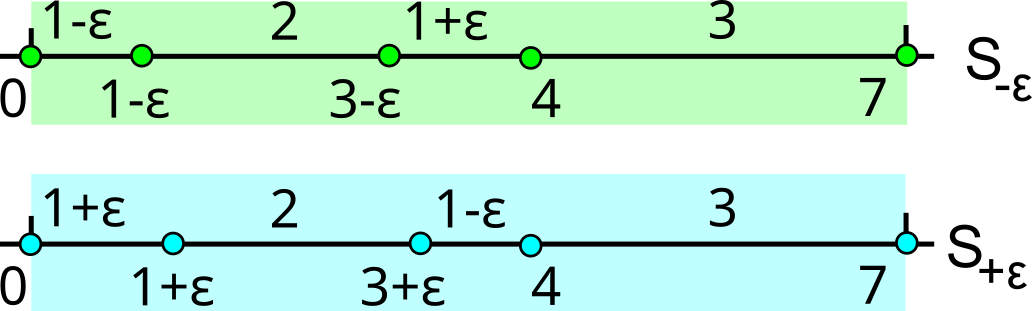}
\hspace*{2mm}
\includegraphics[height=20mm]{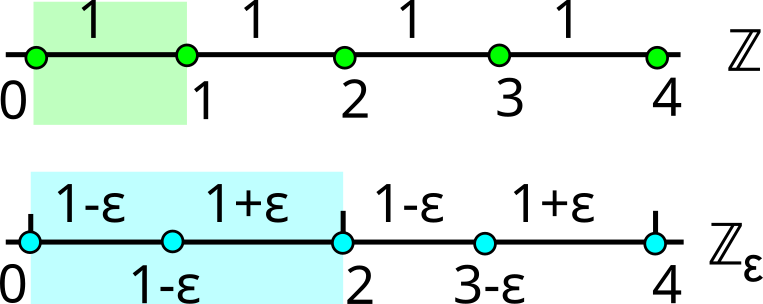}
\caption{\textbf{Left}: the near-duplicate periodic sequences $S_{\pm\ep}=\{0,1\pm\ep,3\pm\ep,4\}+7\Z$ have distant time invariants from Definition~\ref{dfn:time_shift}, see Example~\ref{exa:discontinuity}. 
\textbf{Right}: the periodic sequence $\Z$ and its $\ep$-perturbation $\Z_\ep$ have incomparable time shifts $\TS(\{0\};1)=(1)$ and $\TS(\{0,1-\ep\};2)=(1-\ep,1+\ep)$ of different lengths sizes.
This discontinuity motivates period-independent metrics in Definition~\ref{dfn:period_metrics}.
}
\label{fig:discontinuity}
\end{figure}

Though the time invariants from Definition~\ref{dfn:time_invariants} can be proved to be complete for sequences in $\R$,
Example~\ref{exa:discontinuity} and Fig.~\ref{fig:discontinuity} show their
discontinuity under noise. 

\begin{exa}[discontinuity of time shifts]
\label{exa:discontinuity} 
In $\R$, the periodic sequence $S_0=\{0,1,3,4\}+7\Z$ has two perturbations $S_{\pm\ep}=\{0,1\pm\ep,3\pm\ep,4\}+7\Z$ for any small $\ep>0$.
Rewriting the time shifts $\TS(\{0,1-\ep,3-\ep,4\};7)=(1-\ep,2,1+\ep,3)$ and $\TS(\{0,1+\ep,3+\ep,4\};7)=(1+\ep,2,1-\ep,3)$ in increasing order does not make them close, as the minimal distance $1-\ep$ is followed by the different distances $2<3$ in the nearly identical $S_{\pm\ep}$ for any $\ep>0$, see Fig.~\ref{fig:discontinuity}~(left).
This discontinuity will be resolved by minimising over cyclic permutations, but there is one more obstacle below.\eexa
\end{exa}

It seems natural to always use a minimum period $l>0$ of $S=\{p_1,\dots,p_m\}+l\vec e_1\Z\subset\R\times\R^{n-1}$.
However, the time shift $\TS=(d_1,\dots,d_m)$ of a fixed size $m$ cannot be directly used for comparing sequences that have different sizes of motifs, see Fig.~\ref{fig:discontinuity}~(right).

\section{Period-independent metrics for 1-periodic sequences}
\label{sec:1-periodic_metrics}

Definition~\ref{dfn:period_metrics} introduces continuous metrics after extending motifs to a common size.

\begin{dfn}[cyclic and dihedral metrics under isometry and rigid motion]
\label{dfn:period_metrics}
\tb{(a)}
For any 1-periodic sequences $S=M_S+l_S\vec e_1\Z$ and $Q=M_Q+l_Q\vec e_1\Z$ in $\R\times\R^{n-1}$, let $m=\lcm(|M_S|,|M_Q|)$ be the lowest common multiple of their motif sizes {(cardinalities)}.
For the integers $k_S=\dfrac{m}{|M_S|}$ and $k_Q=\dfrac{m}{|M_Q|}$, 
 the extended motifs defined as 
 $$k_S M_S=\bigcup\limits_{i=1,\dots,k_S}\big(M_S+il_S\vec e_1\big)\qquad\text{ and }\qquad
 k_Q M_Q=\bigcup\limits_{i=1,\dots,k_Q}\big(M_Q+il_Q\vec e_1\big)$$ have the same number $k_S |M_S|=m=k_Q |M_Q|$ of points.
\medskip
 
Any permutation $\ga$ from $C_m,D_m$ acts on the projected motif $v(k_Q M_Q)\subset\R^{n-1}$ as in Lemma~\ref{lem:actions_matrices}.
For any parameter $q\in[1,+\infty]$, the \emph{cyclic} and  \emph{dihedral isometry} metrics  are 
$\CIM_q(S,Q)=\min\limits_{\ga\in C_m}\max\{d_t,d_v\}$ and
$\DIM_q(S,Q)=\min\limits_{\ga\in D_m}\max\{d_t,d_v\}$, where
$$\begin{array}{l}
d_t=m^{-1/q}\big|\big|\TS(k_S M_S;k_S l_S)-\TS(\ga(k_Q M_Q);k_Q l_Q)\big|\big|_q, \\
d_v=\MCD_q\big(v(k_S M_S),\ga(v(k_Q M_Q))\big).
\end{array}$$

\nt
\tb{(b)}
The \emph{cyclic} and \emph{dihedral rigid} metrics $\CRM_q,\DRM_q$ are defined by the same formulae as $\CIM_q,\DIM_q$ after replacing $\MCD_q$ with $\MCS_q$ from Definition~\ref{dfn:finite_metrics}.\edfn
\end{dfn}

\begin{table}[h!]
\caption{Acronyms and references for the new invariants and metrics in sections~\ref{sec:Euclidean_finite}, \ref{sec:discontinuity},~\ref{sec:1-periodic_metrics}.}
\label{tab:acronyms_1-periodic}
\centering
\begin{tabular}{lll}
$\CDM(T)$ & Cyclic Distance Matrix of a finite sequence $T\subset\R^n$ & Definition~\ref{dfn:CDM} \\
$\CDS(T)$ & matrix of Cyclic Distances and Signs of a sequence $T\subset\R^n$ & Definition~\ref{dfn:CDS}  \\
$\MCD_q$ & Metric on Cyclic Distance matrices ($\CDM$) & Definition~\ref{dfn:finite_metrics} \\
$\MCS_q$ & Metric on matrices of Cyclic distances and Signs ($\CDS$) & Definition~\ref{dfn:finite_metrics} \\
$\TS(M;l)$ & Time Shift for a motif $M$ and period $l$ of a sequence & Definition~\ref{dfn:time_shift}  \\
$\CI(S)$ & Cyclic Isometry invariant of a sequence $S\subset\R\times\R^{n-1}$ & Definition~\ref{dfn:periodic_invariants}  \\
$\DI(S)$ & Dihedral Isometry invariant of a sequence $S\subset\R\times\R^{n-1}$ & Definition~\ref{dfn:periodic_invariants}  \\
$\CR(S)$ & Cyclic Rigid invariant of a sequence $S\subset\R\times\R^{n-1}$ & Definition~\ref{dfn:periodic_invariants} \\
$\DR(S)$ & Dihedral Rigid invariant of a sequence $S\subset\R\times\R^{n-1}$ & Definition~\ref{dfn:periodic_invariants} \\
$\CI(S)$ & Cyclic Isometry invariant of a sequence $S\subset\R\times\R^{n-1}$ & Definition~\ref{dfn:periodic_invariants}  \\
$\DI(S)$ & Dihedral Isometry invariant of a sequence $S\subset\R\times\R^{n-1}$ & Definition~\ref{dfn:periodic_invariants}  \\
$\CIM_q$ & Cyclic Isometry Metric on 1-periodic sequences in $\R\times\R^{n-1}$ & Definition~\ref{dfn:period_metrics} 
\\
$\DIM_q$ & Dihedral Isometry Metric on 1-periodic sequences in $\R\times\R^{n-1}$ & Definition~\ref{dfn:period_metrics} \\ 
$\CRM_q$ & Cyclic Rigid Metric on 1-periodic sequences in $\R\times\R^{n-1}$ & Definition~\ref{dfn:period_metrics} 
\\
$\DRM_q$ & Dihedral Rigid Metric on 1-periodic sequences in $\R\times\R^{n-1}$ & Definition~\ref{dfn:period_metrics}
\end{tabular}
\end{table}

In the limit case $q=+\infty$, any factor $a^{\pm 1/q}$ for $a>0$ is interpreted as $\lim\limits_{q\to+\infty}a^{\pm 1/q}=1$.
In Definition~\ref{dfn:period_metrics}, the extended periods $k_S l_S$ and $k_Q l_Q$ can be different.
For simplicity, the metrics $\MCD_q,\MCS_q$ were written via projected motifs as in Definition~\ref{dfn:finite_metrics} but will be computable via the complete invariants 
from Definition~\ref{dfn:periodic_invariants}.
\medskip

For $n=1$, the projected motifs are empty, so the cases of rigid motion and isometry in $\R^0$ trivially coincide.
In both cases, the metrics are obtained by minimizing only the differences $d_t$ between time shifts under cyclic and dihedral permutations.

\begin{exa}[invariant metrics]
\label{exa:period_metrics}
The periodic sequences $S=\{0,1\}+3\Z$ and $Q=\{0,1,3\}+6\Z$ have motifs $M_S=\{0,1\}$ and $M_Q=\{0,1,3\}$ of different sizes $m_S=2$ and $m_Q=3$ whose lowest common multiple is $m=6$.
In the notations of Definition~\ref{dfn:period_metrics}, we get $k_S=\dfrac{m}{|M_S|}=3$, $k_Q=\dfrac{m}{|M_Q|}=2$.
The extended motifs and periods are $3M_S=\{0,1,3,4,6,7\}$, $3l_S=9$, $2M_Q=\{0,1,3,6,7,9\}$, $2l_Q=12$.
Then $\TS(3M_S;9)=(1,2,1,2,1,2)$ and $\TS(2M_Q;12)=(1,2,3,1,2,3)$.
Any cyclic or dihedral permutation of the time shift $\TS(3M_S;9)$ relative to $\TS(2M_Q;12)$ gives the maximum component-wise distance $|1-3|=2$, so $\CIM_{+\infty}(S,Q)=2=\DIM_{+\infty}(S,Q)$.
\end{exa}

\index{cyclic rigid invariant}
\index{dihedral rigid invariant}
\index{polynomial-time complexity}

\begin{thm}[solution to Problem~\ref{pro:1-periodic} for 1-periodic sequences, {\cite[Theorem~4.8]{kurlin2025complete}}]
\label{thm:1-periodic}
\textbf{(a)}
For any 1-periodic sequence $S\subset\R\times\R^{n-1}$ with a motif of $m$ points, $\CI(S),\DI(S)$ from Definition~\ref{dfn:periodic_invariants} are complete invariants under cyclic and dihedral isometry in $\R\times\R^{n-1}$, respectively, and are computable in time $O(m^3n)$.
Then the invariants $\CR(S),\DR(S)$ are complete under cyclic and dihedral rigid motion in $\R\times\R^{n-1}$, respectively, and are computable in time $O(m^3 n+m^2 n^3)$.
\medskip

\noindent
\textbf{(b)}
Any 1-periodic sequence $S\subset\R\times\R^{n-1}$ with a motif of $m$ points can be reconstructed from its complete invariant {under} a relevant equivalence from part (a) in time $O(m^3n)$. 
\medskip

\noindent
\textbf{(c)}
The metrics in Definition~\ref{dfn:period_metrics} remain invariant if any 1-periodic sequence $S=M+l\vec e_1\Z$ is alternatively represented by an extended motif $kM$ and a period $kl$ for any integer $k>0$.
For any 1-periodic sequences $S,Q\subset\R\times\R^{n-1}$ with a lowest common multiple $m$ of motif sizes, the distances $\CIM_q,\DIM_q,\CRM_q,\DRM_q$ in Definition~\ref{dfn:period_metrics} satisfy all metric axioms and are computable in times $O(m^3n)$ and $O(m^3n+m^2 n^3)$ for the cases of isometry and rigid motion, respectively.
\medskip

\noindent
\textbf{(d)}
Let $Q$ be a 1-periodic sequence $S\subset\R\times\R^{n-1}$ after perturbing every point of $S$ up to some Euclidean distance $\ep$ that is smaller than a half-distance between any points of $t(S)$ and of $t(Q)$. 
Then $\CIM_q(S,T),\DIM_q(S,Q),\CRM_q(S,Q),\DRM_q(S,Q)\leq 2\ep$.
\ethm
\end{thm}

\begin{figure}[h!]
\centering
\includegraphics[width=\textwidth]{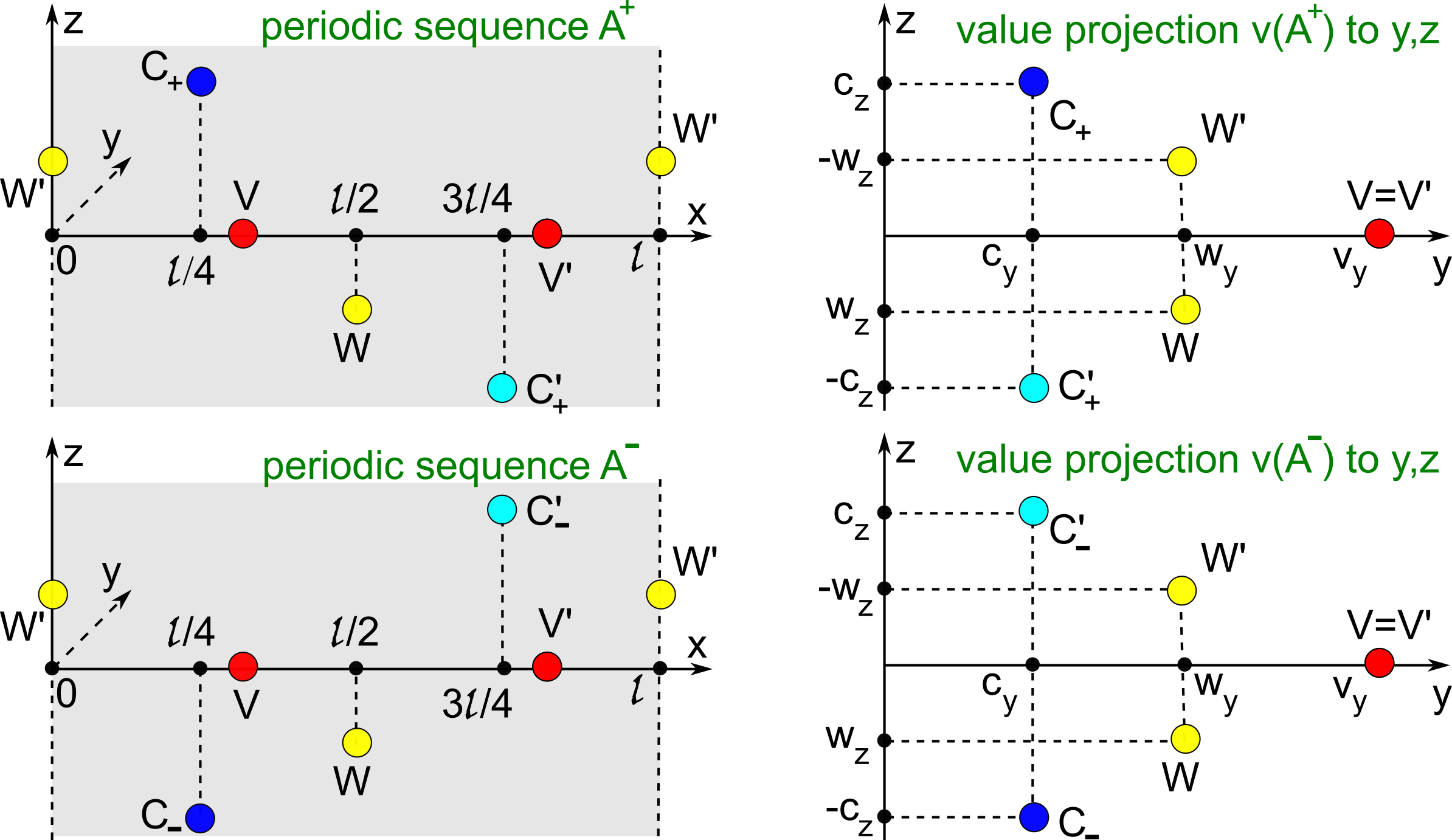}
\caption{These 1-periodic sequences $A^{\pm}\subset\R\times\R^2$ from \cite[Fig.~2]{pozdnyakov2022incompleteness} have identical past invariants.}
\label{fig:6-point_periodic_sequences}
\end{figure}

\begin{exa}[challenging 1-periodic sequences]
\label{exa:6-point_pair}
The infinite family of counter-examples in \cite[Fig.~4]{pozdnyakov2022incompleteness} to the completeness of past distance-based invariants includes the pairs of the 1-periodic sequences $A^{\pm}\subset\R\times\R^2$ with a period $l>0$ and 6-point motifs
$M^+=\{W',C_+,V,W,C'_+,V'\}$ and  $M^-=\{W',C_-,V,W,C'_-,V'\}$ with the points $V=(v_x,v_y,0)$,  $W=(\frac{l}{2},w_y,w_z)$, $C_{\pm}=(\frac{l}{4},c_y,\pm c_z)$, and free parameters $l,w_y,w_z,c_y,c_z>0$, and $v_x,v_y\in[0,\frac{l}{2}]$.
Any point denoted with a prime is obtained by $g(x,y,z)=(x+\frac{l}{2},y,-z)$.
The time projections are identical: $t(M^{\pm})=(0,\frac{l}{4},v_x,\frac{l}{2},\frac{3l}{4},\frac{l}{2}+v_x)$.
Assuming that $v_x\in(\frac{l}{4},\frac{l}{2})$ as in Fig.~\ref{fig:6-point_periodic_sequences}, the time shifts are 
$$\TS(M^{\pm};l)=\left(\frac{l}{4},v_x-\frac{l}{4},\frac{l}{2}-v_x,\frac{l}{4},v_x-\frac{l}{4},\frac{l}{2}-v_x\right).$$
{Order value projections along the $x$-axis from $\frac{l}{2}$ to the right}: 
$$v(M^\pm)=\{(w_y,-w_z),(c_y,\pm c_z),(v_y,0),(w_y,w_z),(c_y,\mp c_z),(v_y,0)\}.$$
The cyclic distance matrices of $M^+$ and $M^-$ are on the left and right, respectively:
\medskip

\noindent
$\left(\begin{array}{llllll}
\mc{yellow}{d_{11}} & d_{12} & d_{21} & \mc{yellow}{d_{11}} & d_{12} & d_{21} \\
d_{21} & \mc{yellow}{d_{22}} & d_{12} & d_{21} & \mc{yellow}{d_{22}} & d_{12}\\
2|w_z| & 2|c_z|  & 0 & 2|w_z| & 2|c_z|  & 0
\end{array}\right)\neq
\left(\begin{array}{llllll}
\mc{yellow}{d_{22}} & d_{12} & d_{21} & \mc{yellow}{d_{22}} & d_{12} & d_{21} \\
d_{21} & \mc{yellow}{d_{11}} & d_{12} & d_{21} & \mc{yellow}{d_{11}} & d_{12}\\
2|w_z| & 2|c_z|  & 0 & 2|w_z| & 2|c_z|  & 0
\end{array}\right)$.
\medskip

The differences in distances are highlighted in \hl{yellow}:
$$\begin{array}{ll}
d_{11}=\sqrt{(w_y-c_y)^2+(w_z\mc{yellow}{+c_z})^2}, &
d_{12}=\sqrt{(c_y-v_y)^2+c_z^2}, \\
d_{22}=\sqrt{(w_y-c_y)^2+(w_z\mc{yellow}{-c_z})^2}, &
d_{21}=\sqrt{(w_y-v_y)^2+w_z^2}.\end{array}$$
The matrix difference has the norm 
$||\CDM(M^+)-\CDM(M^-)||_{\infty}=|d_{11}-d_{22}|>0$ unless $c_z=0$ or $w_z=0$.
If $c_z=0$, $A^{\pm}$ are identical.
If $w_z=0$, then $A^{\pm}$ are isometric by $g(x,y,z)=(x+\frac{l}{2},y,-z)$.
If both $c_z,w_z\neq 0$, then $\CIM_{+\infty}(A^+,A^-)$ is obtained by minimizing over 6 cyclic permutations $\ga\in C_6$. 
The trivial permutation and the shift by 3 positions give $|d_{11}-d_{12}|$.
Any other permutation gives $d_t=\max\{v_x-\frac{l}{4},\frac{l}{2}-v_x\}$ from comparing $\TS(M^+;l)$ with $\ga(\TS(M^-;l))$ and $d_v=\max\{|a-b|\}$ maximized for all pairs of parameters $a,b\in\{d_{11},d_{12},d_{21},d_{22}\}$.
\medskip

In all cases, the metric is positive: $\CIM_{+\infty}(A^+,A^-)\geq|d_{11}-d_{22}|>0$.
Then the invariant $\CI$ from Definition~\ref{dfn:periodic_invariants} distinguishes these sequences $A^+\not\cong A^-$.
\eexa
\end{exa}

Geo-Mapping Problem~\ref{pro:geocodes} becomes much harder for point sets that are periodic in two directions.
The next chapter will solve the case of 2-dimensional lattices.

\bibliographystyle{plain}
\bibliography{Geometric-Data-Science-book}

%
%
%

\chapter{Moduli spaces of 2D lattices under isometry and rigid motion}
\label{chap:lattices2D} 

\abstract{
This chapter continuously parametrises moduli spaces of 2-dimensional lattices under Euclidean isometry, rigid motion, dilation, and homothety.
The new root invariants have easily computable metrics and settle the past discontinuity of reduced bases.
The moduli space of 2-dimensional lattices under rigid motion can be mapped to the sphere without one point.
Hence, any geographic location on Earth can be associated with a canonical lattice. 
We also define chiral distances that continuously measure deviations from higher-symmetry lattices.
}

\section{Representations of lattices by unit cells and reduced bases}
\label{sec:lattice_bases}

All sections in this chapter follow papers \cite{kurlin2024mathematics,bright2023geographic,bright2023continuous} with minor updates of notations.
\myskip

Recall that the most practical equivalence relations (rigid motion, isometry, dilation, homothety) on arbitrary subsets of $\R^n$ were introduced in Example~\ref{exa:isometry}.

\index{lattice}
\index{unit cell}

\begin{dfn}[a basis and a primitive unit cell $U(\vec v_1,\dots,\vec v_n)$ of a lattice $\La\subset\R^n$]
\label{dfn:lattice_cell}
\tb{(a)}
Let vectors $\vec v_1,\dots,\vec v_n$ form a linear {\em basis} in $\R^n$.
A {\em lattice} $\La\subset\R^n$ consists of all linear combinations 
$\sum\limits_{i=1}^n c_i \vec v_i$  with integer coefficients $c_i\in\Z$.
The parallelepiped $U(\vec v_1,\dots,\vec v_n)=\left\{ \sum\limits_{i=1}^n t_i \vec v_i \vl t_i\in[0,1) \right\}\subset\R^n$ is called a \emph{primitive unit cell} of $\La$.
\edfn
\end{dfn}

The inequalities $0\leq t_i<1$ in Definition~\ref{dfn:lattice_cell}(a) guarantee that the copies of primitive unit cells $U(\vec v_1,\dots,\vec v_n)$ translated by all $\vec v\in\La$ are disjoint and cover $\R^n$.
\myskip

\begin{figure}[h]
\includegraphics[height=28mm]{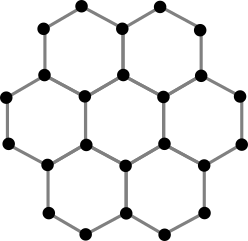}
\includegraphics[height=28mm]{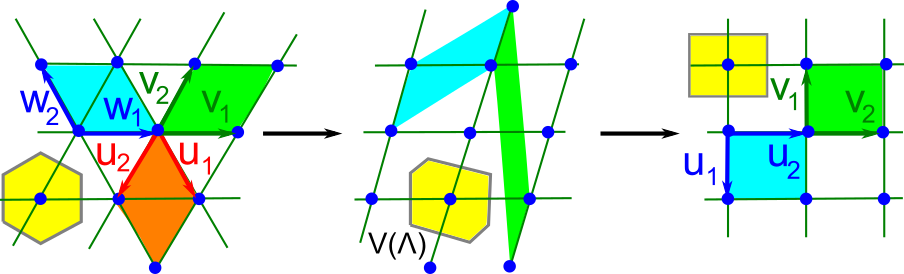}
\label{fig:hexagonal_graphene}
\caption{
\textbf{Left}: a 2-dimensional layer of graphene is formed by carbon atoms.
\textbf{Right}: one can generate a hexagonal lattice (as any other) by infinitely many bases and continuously deform into a rectangular lattice (far right) whose bases $\{\vec v_1,\vec v_2\}$ and $\{\vec u_1,\vec u_2\}$ are related by an orientation-reversing map.
The yellow Voronoi domain $\bar V(\La)$ of any point $\vec p$ in a lattice $\La$ consists of all points $q\in\R^2$ that are non-strictly closer to $\vec p$ than to other points $\La\setminus\{p\}$ in the lattice, see Definition~\ref{dfn:Voronoi_domain}(a).
}
\end{figure}

Recall that the \emph{special linear} group $\SL(\Z_n)$ consists of all $n\times n$ matrices with integer entries and determinant 1.
Let vectors $\vec v_1,\dots,\vec v_n\in\R^n$ form one basis of a lattice $\La\subset\R^n$.
Then, for any matrix $A\in\SL(\Z_n)$, the vectors $A\vec v_1,\dots,A \vec v_n$ form another basis of $\La$ with a very different primitive cell $U(A\vec v_1,\dots,A \vec v_n)$, see Fig.~\ref{fig:hexagonal_graphene}.
\myskip

The past approach to tackling the ambiguity of lattice representations was to consider a reduced basis, briefly as rectangular as possible.
In $\R^3$, \cite{gruber1989reduced} reviewed several reduced bases.
The most common is Niggli's cell \cite{niggli1928krystallographische}, whose 2-dimensional version is introduced below.
For any $\vec v_1=(a_1,a_2)$ and $\vec v_2=(b_1,b_2)$ in $\R^2$, let $\det(\vec v_1,\vec v_2)=a_1b_2-a_2b_1$ be the determinant of the matrix $\matfour{a_1}{b_1}{a_2}{b_2}$ with the columns $\vec v_1,\vec v_2$.

\begin{dfn}[reduced cell]
\label{dfn:reduced_cell}
\tb{(a)}
For a lattice $\La\subset\R^2$ under isometry, a basis and its unit cell $U(\vec v_1,\vec v_2)$ are \emph{reduced} (non-acute) if $|\vec v_1|\leq|\vec v_2|$ and $-\frac{1}{2}\vec v_1^2\leq \vec v_1\cdot \vec v_2\leq 0$.
\myskip

\nt
\tb{(b)}
Under rigid motion, the conditions are weaker: $|\vec v_1|\leq|\vec v_2|$ and $|\vec v_1\cdot \vec v_2|\leq\frac{1}{2}\vec v_1^2$, $\det(\vec v_1,\vec v_2)>0$, and the new \emph{special condition} : if $|\vec v_1|=|\vec v_2|$ then $\vec v_1\cdot \vec v_2\geq 0$.
\edfn
\end{dfn}

All bases in Fig.~\ref{fig:hexagonal_graphene} are reduced under rigid motion.
The condition $|\vec v_1\cdot \vec v_2|\leq\frac{1}{2}\vec v_1^2$ in Definition~\ref{dfn:reduced_cell} geometrically means that $\vec v_1,\vec v_2$ are close to being orthogonal: the projection of $\vec v_2$ to $\vec v_1$ is between $\pm\frac{1}{2}|\vec v_1|$.
The conditions $|\vec v_1|\leq|\vec v_2|$ and $-\frac{1}{2}\vec v_1^2\leq \vec v_1\cdot \vec v_2\leq 0$ in Definition~\ref{dfn:reduced_cell} coincide with the conventional definition from \cite[section 9.2.2]{aroyo2013international} for type II (non-acute) cells in $\R^3$ if we choose $\vec v_3$ to be very long and orthogonal to $\vec v_1,\vec v_2$.
Alternative type I cells with non-obtuse angles have $0\leq \vec v_1\cdot \vec v_2\leq\frac{1}{2}\vec v_1^2$.
\myskip

\cite[Proposition~3.10(a)]{kurlin2024mathematics} proves the uniqueness of a reduced basis under isometry. 
\medskip

Another well-known cell of a lattice $\La\subset\R^n$ is the \emph{Voronoi domain} \cite{voronoi1908nouvelles}, also called the \emph{Wigner-Seitz cell}, \emph{Brillouin zone} or \emph{Dirichlet cell}.
We use the word \emph{domain} not to confuse it with a unit cell in Definition~\ref{dfn:lattice_cell}. 
Though the Voronoi domain can be defined for any point of a lattice, it suffices to consider only the origin $0$.

\index{Voronoi domain}

\begin{dfn}[Voronoi domain $\bar V(\La)$]
\label{dfn:Voronoi_domain}
\tb{(a)}
The \emph{Voronoi domain} of a lattice $\La\subset\R^n$ is the neighbourhood $\bar V(\La)=\{\vec p\in\R^n: |\vec p|\leq|\vec p-\vec v| \text{ for any }\vec v\in\La\}$ of $0$ consisting of all $\vec p\in\R^n$ that are non-strictly closer to $0$ than to other points $\vec v\in\La$.
\myskip

\nt
\tb{(b)}
A vector $\vec v\in\La$ is called a \emph{Voronoi vector} if the bisector hyperspace $H(0,\vec v)=\{\vec p\in\R^n \vl \vec p\cdot \vec v=\frac{1}{2}\vec v^2\}$ between 0 and $\vec v$ intersects $\bar V(\La)$.
If $\bar V(\La)\cap H(0,\vec v)$ is an $(n-1)$-dimensional face of $\bar V(\La)$, then $\vec v$ is called a \emph{strict} Voronoi vector. 
\edfn
\end{dfn}

\begin{figure}[h]
\includegraphics[width=1.0\textwidth]{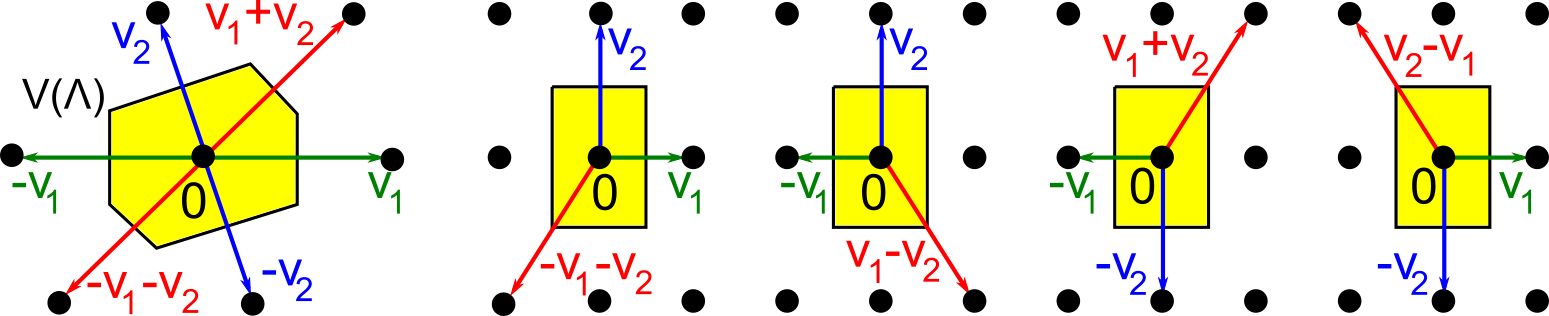}
\caption{\textbf{Left}: a generic lattice $\La\subset\R^2$ has a hexagonal Voronoi domain with an obtuse superbase $\vec v_1,\vec v_2,\vec v_0=-\vec v_1-\vec v_2$, which is unique under permutations and central symmetry.
\textbf{Other pictures}: two pairs of obtuse superbases (related by reflection) for a rectangular lattice.}
\label{fig:Voronoi2D}
\end{figure}

\index{Voronoi domain}

Fig.~\ref{fig:Voronoi2D} shows how the Voronoi domain $\bar V(\La)$ can be obtained as the intersection of the closed half-spaces $S(0,v)=\{\vec p\in\R^n \vl \vec p\cdot \vec v\leq\frac{1}{2}\vec v^2\}$ whose boundaries $H(0,v)$ are bisectors between $0$ and all strict Voronoi vectors $\vec v\in\La$.
A generic lattice $\La\subset\R^2$ has a hexagonal Voronoi domain $\bar V(\La)$ with six Voronoi vectors.
\medskip

Any lattice is determined by its Voronoi domain by \cite[Lemma~A.2]{kurlin2024mathematics}.
However, the combinatorial structure of $\bar V(\La)$ is discontinuous under perturbations.
Almost any perturbation of a rectangular basis in $\R^2$ gives a non-rectangular basis generating a lattice whose Voronoi domain $\bar V(\La)$ is hexagonal, not rectangular.
Hence, any integer-valued descriptors of $\bar V(\La)$, such as the numbers of vertices or edges, are always discontinuous and unsuitable for continuous quantifications.
\myskip

Lemma~\ref{lem:Voronoi_vectors} shows how to find all Voronoi vectors of any lattice $\La\subset\R^n$. 
The doubled lattice is $2\La=\{2\vec v \vl \vec v\in\La\}$.
Vectors $\vec u,\vec v\in\La$ are called \emph{$2\La$-equivalent} if $\vec u-\vec v\in 2\La$.
Then any vector $\vec v\in\La$ generates its $2\La$-class $\vec v+2\La=\{\vec v+2\vec u \vl \vec u\in\La\}$, which is $2\La$ translated by $\vec v$ and containing $-\vec v$.
All classes of $2\La$-equivalent vectors form the quotient space $\La/2\La$.
Any 1-dimensional lattice $\La$ generated by a vector $\vec v$ has the quotient $\La/2\La$ consisting of only two classes $\La$ and $\vec v+\La$.  
 
\begin{lem}[criterion for Voronoi vectors {\cite{minkowski1891ueber}, \cite[Theorem~2]{conway1992low}}]
\label{lem:Voronoi_vectors}
For any lattice $\La\subset\R^n$, a non-zero vector $\vec v\in\La$ is a Voronoi vector of $\La$ if and only if $\vec v$ is a shortest vector in its $2\La$-class $\vec v+2\La$.
Also, $\vec v$ is a strict Voronoi vector if and only if $\pm \vec v$ are the only shortest vectors in the $2\La$-class $\vec v+2\La$.
\elem
\end{lem}

Any lattice $\La\subset\R^2$ generated by $\vec v_1,\vec v_2$ has $\La/2\La=\{\vec v_1,\vec v_2,\vec v_1+\vec v_2\}+\La$.
Notice that the vectors $\vec v_1\pm \vec v_2$ belong to the same $2\La$-class.
Assume that $\vec v_1,\vec v_2$ are not longer than $\vec v_1+\vec v_2$, which holds if the angle $\angle(\vec v_1,\vec v_2)\in[60^{\circ},120^{\circ}]$.
If the sum $\vec v_1+\vec v_2$ is shorter than $\vec v_1-\vec v_2$ as in Fig.~\ref{fig:Voronoi2D}~(left), then $\La$ has three pairs of strict Voronoi vectors $\pm \vec v_1,\pm \vec v_2,\pm(\vec v_1+\vec v_2)$.
If $\vec v_1\pm \vec v_2$ have the same length, the unit cell spanned by $\vec v_1,\vec v_2$ degenerates to a rectangle, $\La$ has four non-strict Voronoi vectors $\pm \vec v_1\pm \vec v_2$. 
\medskip

The triple of vector pairs $\pm \vec v_1,\pm \vec v_2,\mp(\vec v_1+\vec v_2)$ in Fig.~\ref{fig:Voronoi2D} motivates the concept of a superbase with the extra vector $\vec v_0=-\vec v_1-\vec v_2$, which extends to any dimension $n$ by setting $\vec v_0 = -\sum\limits_{i=1}^n \vec v_n$.
For dimensions 2 and 3, \cite[Theorem~2.9]{kurlin2024mathematics} proved that any lattice has an obtuse superbase of vectors whose pairwise scalar products are non-positive and are called \emph{Selling parameters} \cite{selling1874ueber}.
For any superbase in $\R^n$, the negated parameters $p_{ij}=-\vec v_i\cdot \vec v_j$ can be interpreted as conorms of lattice characters, which are functions $\chi: \La\to\{\pm 1\}$ satisfying $\chi(\vec u+\vec v)=\chi(\vec u)\chi(\vec v)$), see \cite[Theorem~6]{conway1992low}.
So $p_{ij}$ will be defined as \emph{conorms} only for an obtuse superbase below. 

\index{superbase}
\index{obtuse superbase}
\index{conorm}

\begin{dfn}[obtuse superbase and conorms $p_{ij}$]
\label{dfn:conorms}
For any basis $\vec v_1,\dots,\vec v_n$ of $\R^n$, the \emph{superbase} $\vec v_0,\vec v_1,\dots,\vec v_n$ includes the vector $\vec v_0=-\sum\limits_{i=1}^n \vec v_i$.
The \emph{conorms} $p_{ij}=-\vec v_i\cdot \vec v_j$ are the negative scalar products of the vectors above. 
The superbase is \emph{obtuse} if all conorms $p_{ij}\geq 0$, so all angles between vectors $\vec v_i,\vec v_j$ are non-acute for distinct indices $i,j\in\{0,1,\dots,n\}$.
The superbase is called \emph{strict} if all $p_{ij}>0$.
\edfn
\end{dfn}

Formula (1) in \cite{conway1992low} has a typo initially defining $p_{ij}$ as exact Selling parameters, but later \cite[Theorems 3,~7,~8]{conway1992low} use the non-negative conorms $p_{ij}=-\vec v_i\cdot \vec v_j\geq 0$.
\medskip

The indices of a conorm $p_{ij}$ are distinct and unordered.
We set $p_{ij}=p_{ji}$ for all indices $i,j$.
For $n=1$, the 1-dimensional lattice generated by a vector $\vec v_1$ has the obtuse superbase consisting of the two vectors $\vec v_0=-\vec v_1$ and $\vec v_1$, so the only conorm $p_{01}=-\vec v_0\cdot \vec v_1=\vec v_1^2$ is the squared length of $\vec v_1$.
Any superbase of $\R^n$ has $\dfrac{n(n+1)}{2}$ conorms $p_{ij}$, for example, three conorms $p_{01},p_{02},p_{12}$ in dimension 2.

\index{vonorm}

\begin{dfn}[partial sums $\vec v_S$ and vonorms $\vec v_S^2$]
\label{dfn:vonorms}
Let a lattice $\La\subset\R^n$ have a superbase $B=\{\vec v_0,\vec v_1,\dots,\vec v_n\}$. 
For any proper subset $S\subset\{0,1,\dots,n\}$ of indices,
 consider its complement $\bar S=\{0,1,\dots,n\}\setminus S$ and the \emph{partial sum} $\vec v_S=\sum\limits_{i\in S} \vec v_i$ whose squared lengths $\vec v_S^2$ are called the \emph{vonorms} of $B$ and can be expressed as 
$$\vec v_S^2=\left(\sum\limits_{i\in S} \vec v_i\right)\left(-\sum\limits_{j\in\bar S}\vec v_j\right)=-\sum\limits_{i\in S,j\in\bar S}\vec v_{j}\cdot \vec v_j=\sum\limits_{i\in S,j\in\bar S}p_{ij}.$$
For $n=2$, we get the following simple formulae
$$
\vec v_0^2=p_{01}+p_{02},\qquad
\vec v_1^2=p_{01}+p_{12},\qquad 
\vec v_2^2=p_{02}+p_{12}.
\leqno{(\ref{dfn:vonorms}a)}$$
The above formulae allow us to express the conorms via vonorms as follows
$$
p_{12}=\dfrac{1}{2}(\vec v_1^2+\vec v_2^2-\vec v_0^2),\quad
p_{01}=\dfrac{1}{2}(\vec v_0^2+\vec v_1^2-\vec v_2^2),\quad
p_{02}=\dfrac{1}{2}(\vec v_0^2+\vec v_2^2-\vec v_1^2).
\leqno{(\ref{dfn:vonorms}b)}$$
So $p_{ij}=\dfrac{1}{2}(\vec v_i^2+\vec v_j^2-\vec v_k^2)$ for distinct $i,j\in\{0,1,2\}$ and $k=\{0,1,2\}-\{i,j\}$.
\edfn
\end{dfn}

Lemma~\ref{lem:partial_sums} will later help to prove that a lattice is uniquely determined under isometry by an obtuse superbase, hence by its vonorms or, equivalently, conorms.

\begin{lem}[Voronoi vectors $v_S$ {\cite[Theorem~3]{conway1992low}}]
\label{lem:partial_sums}
For any obtuse superbase $v_0,v_1,\dots,v_n$ of a lattice, all partial sums $v_S$ from Definition~\ref{dfn:vonorms} split into $2^n-1$ symmetric pairs $v_S=-v_{\bar S}$, which are Voronoi vectors representing distinct $2\La$-classes in $\La/2\La$.
All Voronoi vectors $v_S$ are strict if and only if all $p_{ij}>0$.
\elem
\end{lem}

By Conway and Sloane \cite[section~2]{conway1992low}, any lattice $\La\subset\R^n$ that has an obtuse superbase is called a \emph{lattice of Voronoi's first kind}.
Any lattice in dimensions 2 and 3 is of Voronoi's first kind due to \cite[p.~277]{voronoi1908nouvelles} for $n=2$ and \cite[Section~III.4.3]{delone1934mathematical} for $n=3$.

\begin{thm}[reduction to an obtuse superbase]
\label{thm:reduction_to_obtuse_superbase}
Any lattice $\La$ in dimensions $2$ and $3$ has an obtuse superbase $\{v_0,v_1,\dots,v_n\}$ so that $v_0=-\sum\limits_{i=1}^n v_i$ and all conorms $p_{ij}=-v_i\cdot v_j\geq 0$ for all distinct indices $i,j\in\{0,1,\dots,n\}$.
\ethm
\end{thm}

Conway and Sloane in \cite[section~7]{conway1992low} attempted to prove Theorem~\ref{thm:reduction_to_obtuse_superbase} for $n=3$ by example, which is corrected in \cite{kurlin2022complete}.
Theorem~\ref{thm:reduction_to_obtuse_superbase} for $n=2$ is proved in \cite[appendix~A]{kurlin2024mathematics}.
Proposition~\ref{prop:reduced_bases} will establish a 1-1 correspondence between obtuse superbases and reduced bases.
The latter bases are implemented by many fast algorithms in crystallography \cite{aroyo2011crystallography}.
So our lattice input will be any obtuse superbase.

\section{Geo-mapping problem for lattices under equivalences in $\R^2$}
\label{sec:geomapping_lattices2D}

\begin{figure}[h]
\includegraphics[width=1.0\textwidth]{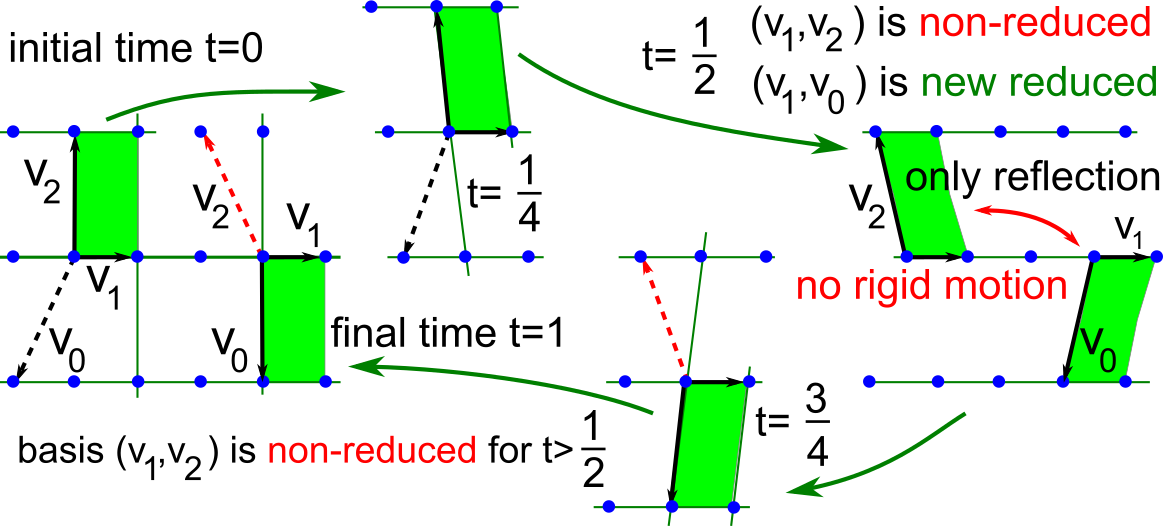}
\caption{
Discontinuity of obtuse superbases under rigid motion.
The obtuse superbase of the vectors $\vec v_1=(1,0)$, $\vec v_2(t)=(-t,2)$, $\vec v_0(t)=(t-1,-2)$ deforms for $t\in[0,1]$.
The initial and final superbases at $t=0$ and $t=1$ generate the same rectangular lattice but are not related by rigid motion.
}
\label{fig:rect_basis_discontinuity}
\end{figure}

This section shows how past approaches to lattice classifications remained discontinuous and states a suitable version of Geo-Mapping~\ref{pro:geocodes} for 2-dimensional lattices.
Fig.~\ref{fig:rect_basis_discontinuity} illustrates how both reduced basis and obtuse superbase discontinuously change under rigid motion.
As for many other objects, we need invariants, because \cite[Theorem~15]{widdowson2022average} proved that any reduced basis is discontinuous under coordinate-wise comparisons.

Any lattice $\La\subset\R^2$ with a basis $\vec v_1,\vec v_2$ defines the \emph{positive quadratic form} 
$$Q(x,y)=(x\vec v_1+y\vec v_2)^2=q_{11}x^2+2q_{12}xy+q_{22}y^2\geq 0 
\text{ for all }x,y\in\R,$$ 
where $q_{11}=\vec v_1^2$, $q_{22}=\vec v_2^2$, $q_{12}=\vec v_1\cdot \vec v_2$.
Changing the basis $\vec v_1,\vec v_2$ (possibly by reflection) is equivalent to replacing $x,y$ by the linear combinations of the coordinates of $x\vec v_1+y\vec v_2$ in a new basis.
Conversely, any positive quadratic form $Q(x,y)$ can be written as a sum $(a_1x+b_1y)^2+(a_2x+b_2y)^2$, see \cite[Theorem~2 on p.~116]{delone1975bravais}, and defines the lattice with the basis $\vec v_1=(a_1,a_2)$, $\vec v_2=(b_1,b_2)$.
\medskip

In 1773, Lagrange \cite{lagrange1773recherches} proved that any positive quadratic form can be rewritten so that $0<q_{11}\leq q_{22}$ and $-q_{11}\leq 2q_{12}\leq 0$. 
The resulting non-acute \emph{reduced} basis $\vec v_1,\vec v_2$ satisfies $0<\vec v_1^2\leq \vec v_2^2$ and $-\vec v_1^2\leq 2\vec v_1\cdot \vec v_2\leq 0$ without the new special conditions in Definition~\ref{dfn:reduced_cell}.
Alternatively, $0\leq 2q_{12}\leq q_{11}$ and $0\leq 2\vec v_1\cdot \vec v_2\leq \vec v_1^2$ define a non-obtuse reduced basis.
The mirror images $\La^{\pm}$ of $\La(\frac{1}{4})$ in Fig.~\ref{fig:rect_basis_discontinuity} (top) generated by the obtuse reduced bases $\vec v_1=(1,0)$, $\vec v_2^{\pm}=(-\frac{1}{4},\pm 2)$ have the same reduced form $Q(x,y)=x^2-\frac{1}{2}xy+4y^2$ not distinguishing $\La^{\pm}$ under rigid motion.
\medskip

If $\vec v_1,\vec v_2$ form a unique reduced basis, Lemma~\ref{lem:VF} shows that $\pm\{\vec v_1,\vec v_2,\vec v_1+\vec v_2\}$ are the three (pairs of) shortest Voronoi vectors.
Then the \emph{metric tensor} $(\vec v_1^2,\vec v_1\cdot \vec v_2,\vec v_2^2)$ is a complete isometry invariant but doesn't distinguish mirror images (enantiomorphs).
Instead of one scalar product and two squared lengths, Delone used the homogeneous parameters \cite[section~29]{delone1938geometry} equal to the conorms $p_{ij}$ from Definition~\ref{dfn:conorms}: \\
$$\begin{array}{l}
p_{01}=q_{11}+q_{12}=\vec v_1^2+\vec v_1\cdot \vec v_2=\vec v_1\cdot(\vec v_1+\vec v_2)=-\vec v_0\cdot \vec v_1,\\ 
p_{02}=q_{22}+q_{12}=\vec v_2^2+\vec v_1\cdot \vec v_2=\vec v_2\cdot(\vec v_1+\vec v_2)=-\vec v_0\cdot \vec v_2,\\
p_{12}=-q_{12}=-\vec v_1\cdot \vec v_2. \end{array}$$
The quadratic form becomes a sum of squares: $Q_\La=p_{01}x^2+p_{22}y^2+p_{12}(x-y)^2$.
The inequalities for $q_{ij}$ are equivalent to the simple ordering $0\leq p_{12}\leq p_{01}\leq p_{02}$, which Definition~\ref{dfn:RI} will use to introduce a more convenient root invariant.
\medskip

Bi-continuity conditions in general Problem~\ref{pro:geocodes}(d,e) become challenging for periodic sets already in dimension 1.
Any lattice under rigid motion in $\R$ is equivalent to a periodic sequence $l\Z$ for a period $l>0$.
However, for any small $\ep>0$, the $\ep$ perturbation of $l$ up to (say) $l+\ep$ makes the perturbed sequence $(l+\ep)\Z$ very different from $l\Z$. 
Indeed, extra $\ep$-shifts gradually move points $(l+\ep)n$ further and further away from $ln$ as $n\in Z$ increases.
Nonetheless, we will prove bi-continuity of bijections between invariant spaces of 2D lattices and the following metric spaces of obtuse superbases.

\label{obtuse superbase}

\begin{dfn}[spaces of obtuse superbases under equivalences]
\label{dfn:SBI}
Let $B=\{\vec v_i\}_{i=0}^n$ and $B'=\{\vec u_i\}_{i=0}^n$ be any obtuse superbases in $\R^n$.
\myskip

\nt
\textbf{(a)}
Let $\SBR(\R^n)$ be the space of equivalence classes obtuse of \emph{superbases under isometry} under the action of special orthogonal maps $f\in\SO(\R^n)$ with the \emph{superbase rigid metric} $\SRM_{\infty}(B,B')=\min\limits_{f\in\SO(\R^n)}\max\limits_{i=0,\dots,n}|f(\vec u_i)-\vec v_i|$.
\medskip

\noindent
\textbf{(b)}
Let $\SBI(\R^n)$ be the space of equivalence classes obtuse of \emph{superbases under isometry} under the action of orthogonal maps $f\in\Or(\R^n)$ with the \emph{superbase isometry metric} $\SIM_{\infty}(B,B')=\min\limits_{f\in\Or(\R^n)}\max\limits_{i=0,\dots,n}|f(\vec u_i)-\vec v_i|$.
\myskip

\noindent
\textbf{(c)}
Let $\SBD(\R^n)$ be the space of equivalence classes obtuse of \emph{superbases under dilation} under the action of dilation maps $f\in\SO(\R^n)\times\R_+$ with the \emph{superbase dilation metric} $\SDM_{\infty}(B,B')=\min\limits_{f\in\SO(\R^n)\times\R_+}\max\limits_{i=0,\dots,n}|f(\vec u_i)-\vec v_i|$.
\myskip

\noindent
\textbf{(d)}
Let $\SBH(\R^n)$ be the space of equivalence classes obtuse of \emph{superbases under homothety} under the action of homothety maps $f\in\Or(\R^n)\times\R_+$ with the \emph{superbase homothety metric} $\SHM_{\infty}(B,B')=\min\limits_{f\in\Or(\R^n)\times\R_+}\max\limits_{i=0,\dots,n}|f(\vec u_i)-\vec v_i|$.
\edfn
\end{dfn}

Since any continuous function over a compact domain achieves its minimum value and $\SO(\R^n),\Or(\R^n)$ are compact, the minima in Definition~\ref{dfn:SBI}(a,b) are achievable.
For fixed superbases $B,B'$, one can restrict uniform scaling by reasonable bounds to guarantee the existence of minima in Definition~\ref{dfn:SBI}(c,d).
\myskip

Problem~\ref{pro:lattices2D} is a case of Problem~\ref{pro:geocodes} for 2D lattices.
We use isometry as the main equivalence, but moduli spaces of lattices will also be parametrised under rigid motion and their compositions with uniform scaling (dilation and homothety). 

\index{lattice}
\index{geocode}
\index{geo-mapping}

\begin{pro}[geo-mapping for 2D lattices]
\label{pro:lattices2D}
Design a \emph{geocode} on the space of 2D lattices under isometry that is an invariant satisfying the following conditions.
\smallskip

\noindent
\tb{(a)} 
\emph{Completeness:} 
any lattices $\La\simeq\La'$ are isometric in $\R^2$ if and only if $I(\La)=I(\La')$.
\myskip

\noindent
\tb{(b)} 
\emph{Reconstruction:} 
any lattice $\La\subset\R^2$ can be reconstructed from its invariant value $I(\La)$, uniquely under isometry in $\R^2$.
\myskip

\noindent
\tb{(c)} 
\emph{Metric:} 
there is a metric $d$ on the invariant space $\{I(\La) \vl \text{lattices }\La\subset\R^2\}$, which satisfies all axioms in Definition~\ref{dfn:metrics}(a).
\myskip

\noindent
\tb{(d)} 
\emph{Continuity:} 
there is a constant $\la>0$ such that, for any $\ep>0$, if any lattices $\La,\La'\subset\R^2$ have obtuse have obtuse superbases $B,B'$, respectively with $\SIM_{\infty}(B,B')\leq\ep$, then $d(I(\La),I(\La'))\leq\la\sqrt{\ep}$.
\myskip

\noindent
\tb{(e)} 
\emph{Inverse continuity:} 
 for any $\ep>0$, there is $\de>0$ such that if lattices $\La,\La'\subset\R^2$ satisfy $d(I(\La),I(\La'))\leq\de$, then they have obtuse superbases $B,B'$, respectively, with $\SIM_{\infty}(B,B')\leq\ep$.
\myskip

\nt
\tb{(f)}
\emph{Realisability:} 
the invariant space $\{I(\La) \vl \text{lattices }\La\subset\R^2\}$ can be parametrised so that we can generate any value $I(\La)$ realisable by some lattice $\La\subset\R^2$.
\myskip

\nt
\tb{(g)}
\emph{Euclidean embedding:} 
the invariant space $\{I(\La) \vl \text{lattices }\La\subset\R^2\}$ with the metric $d$ allows a Lipschitz embedding into a suitable Euclidean space $\R^N$ for some $N$.
\myskip

\noindent
\tb{(h)}
\emph{Computability:} the invariant $I$ and the metric $d(I(\La),I(\La'))$ can be computed in a constant time $O(1)$ from reduced bases of $\La,\La'$.
\epro 
\end{pro}

A geocode that satisfies Problem~\ref{pro:lattices2D} continuously parametrises the space of 2D lattices under isometries, similarly defined for other equivalences below.

\index{lattice spaces}
\begin{dfn}[moduli spaces of lattices under four equivalences in $\R^n$]
\label{dfn:lattice_spaces}
We consider all lattices $\La\subset\R^n$ below.
\myskip

\noindent
\tb{(a)}
The \emph{Lattice Isometry Space} $\LIS(\R^n)$ is the space of lattices $\La$ under isometry.
\myskip

\noindent
\tb{(b)}
The \emph{Lattice Rigid Space} $\LRS(\R^n)$ is the space of lattices under rigid motion.
\myskip

\noindent
\tb{(c)}
The \emph{Lattice Homothety Space} $\LHS(\R^n)$ is the space of lattices under homothety.
\myskip

\noindent
\tb{(d)}
The \emph{Lattice Dilation Space} $\LDS(\R^n)$ is the space of lattices under dilation.
\edfn
\end{dfn}

The traditional approach to deciding if lattices are isometric is to compare their conventional or reduced cells.
Though this comparison theoretically gives a complete invariant, in practice, all real lattices in periodic crystals are non-isometric due to noise in measurements. 
Since all atoms vibrate, any real lattice basis is always perturbed.
The discontinuity of reduced bases under perturbations was experimentally known since 1965 \cite[p.~80]{lawton1965reduced} and was proved for all potential reductions in \cite[Theorem~15]{widdowson2022average}.
\medskip

A more practical goal is to design a complete invariant that is continuous under any perturbations of (bases of) lattices.
Such a \emph{geocode}, which is more generally defined in Problem~\ref{pro:geocodes}, will unambiguously parametrise the Lattice Rigid Space $\LRS(\R^n)$ consisting of infinitely many equivalence classes of lattices under rigid motion in $\R^n$.
For example, our Earth is continuously parametrised by the latitude and longitude, very similar to the Lattice Dilation Space $\LDS(\R^2)$ as will become clear soon. 
\medskip

The Lattice Rigid Space $\LRS(\R^n)$ is continuous and connected because any two lattices can be joined by a continuous deformation of their bases as in Fig.~\ref{fig:hexagonal_graphene}.
Such deformation can always be visualised as a continuous path in the space $\LRS(\R^n)$.
\myskip

The Euclidean embeddability in \ref{pro:lattices2D}(g) raises Problem~\ref{pro:lattices2D} above metric geometry to define a simpler Euclidean structure on $\LIS(\R^n)$.
It is easy to multiply any lattice by a fixed scalar, but a sum of any two lattices is harder to define in a meaningful way independent of lattice bases.
We will overcome this obstacle by a natural embedding of invariant spaces of lattices into $\R^3$ and $\R^4$ to fully solve Problem~\ref{pro:lattices2D}.
\myskip

The isometry classification condition in \ref{pro:lattices2D}(a) can be interpreted via group actions as follows, see \cite{engel2004lattice} and \cite{zhilinskii2016introduction}.
Let $\mathcal{B}_n$ be the space of all linear bases in $\R^n$.
\myskip

Under a change of basis, all lattices in $\R^n$ form the $n^2$-dimensional orbit space $\mathcal{L}_n=\mathcal{B}_n/\GL(\Z^n)$, see \cite[formula (1.38), p.~34]{engel2004lattice}.
Under homothety, the orbit space $\mathcal{L}_n^d=\mathcal{L}_n/\R_+^\times$ becomes $(n^2-1)$-dimensional.
Under orthogonal maps from the group $\Or(\R^n)$, the orbit space of lattices can be identified with the cone $\mathcal{C}_+(\mathcal{Q}_n)=\mathcal{B}_n/\Or(\R^n)$ of positive quadratic forms, where $\mathcal{Q}_n$ denotes the space of real symmetric $n\times n$ matrices, see \cite[formula (1.67), p.~41]{engel2004lattice}.
The Lattice Isometry Space $\LIS(\R^n)$ was called the space of \emph{intrinsic} lattices $\mathcal{L}_n^o=\mathcal{C}_+(\mathcal{Q}_n)/\GL(\Z^n)$ in \cite[formula (1.70), p.~42]{engel2004lattice}.
\myskip

Another approach to identify an intrinsic lattice (isometry class), say for $n=2$,  was to choose a fundamental domain of the action of $\GL(\Z^2)$ on the cone $\mathcal{C}_+(\mathcal{Q}_2)$.
This choice is equivalent to a choice of a reduced basis, which can be discontinuous.
\myskip

\begin{figure}[h]
\includegraphics[width=1.0\textwidth]{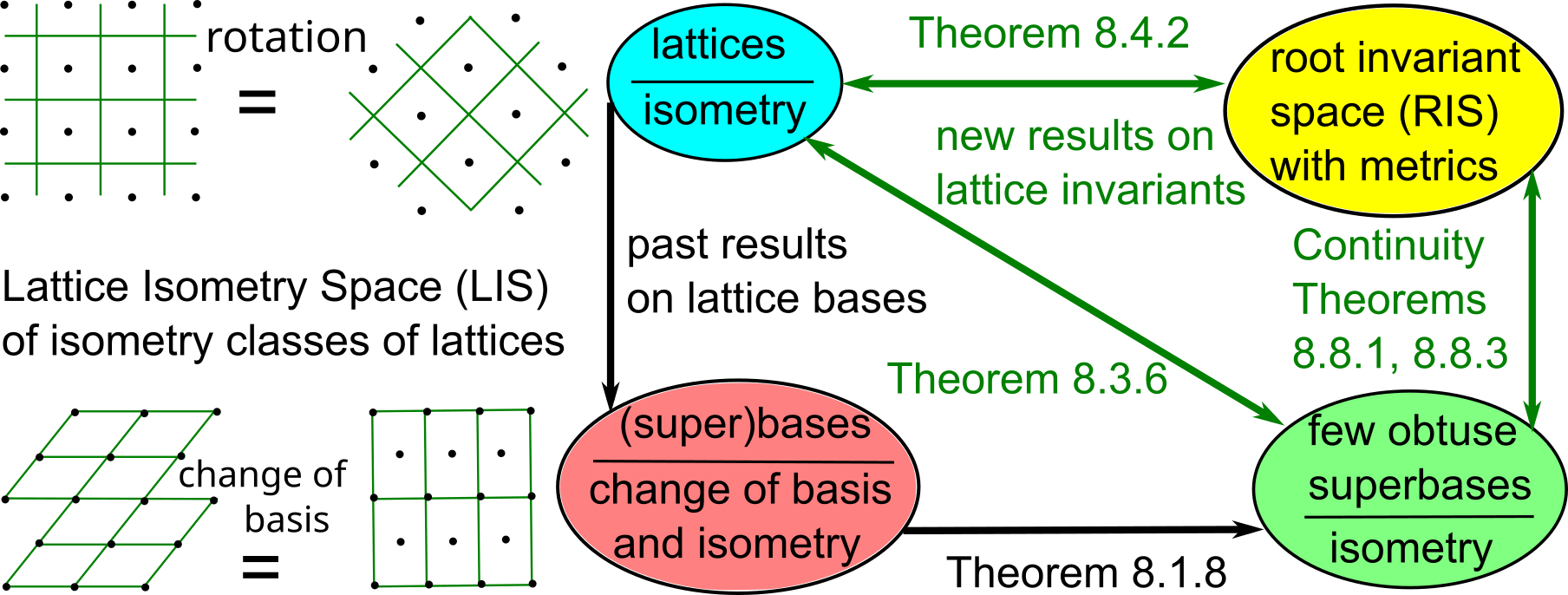}
\caption{
$\LIS(\R^2)$ is bijectively and continuously mapped to the space of root invariants, which are ordered triples of square roots of scalar products of vectors of an obtuse superbase of a lattice $\La\subset\R^2$.}
\label{fig:lattice_classification}
\end{figure}

Mirror reflections of any lattice $\La$ correspond to quadratic forms $q_{11}x^2\pm 2q_{12}xy+q_{22}y^2$ that differ by a sign of $q_{12}$.
To distinguish mirror images of lattices, Definition~\ref{dfn:sign} will introduce $\sign(\La)$.
Then continuous deformations of lattices become continuous paths in a space of invariants, see Remark~\ref{rem:group_action}.
\myskip
 
Fig.~\ref{fig:lattice_classification} summarises the past obstacles and a full solution to Problem~\ref{pro:lattices2D}.
The Root Invariant Space $\RIS(\R^2)$ consists of ordered triples of square roots of conorms from Definition~\ref{dfn:conorms}.
Related invariants will continuously parametrise
spaces of lattices under rigid motion, dilation, and homothety, as defined below.

\section{Invariants of an obtuse superbase of a 2-dimensional lattice}
\label{sec:lattices2D_forms}

Definition~\ref{dfn:RI} introduces voforms $\VF$ and coforms $\CF$, which are triangular cycles whose three nodes are marked by vonorms and conorms, respectively.
We start from any obtuse superbase $B$ of a lattice $\La\subset\R^2$ to define $\VF$, $\CF$, and a root invariant $\RI$. Lemma~\ref{lem:lattice_invariants}(a) will justify that $\RI$ depends only on $\La$, not on $B$.

\begin{figure}[h]
\includegraphics[height=35mm]{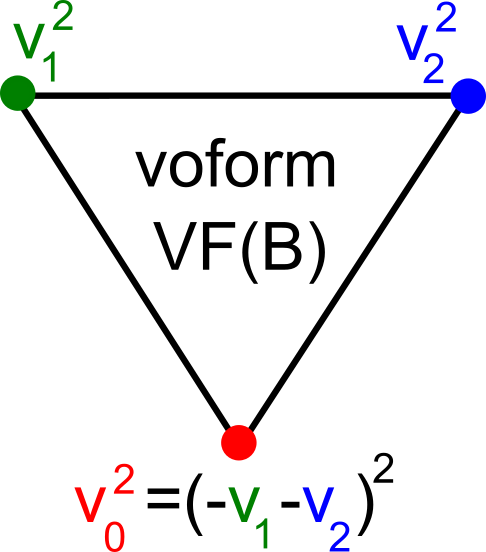}
\includegraphics[height=35mm]{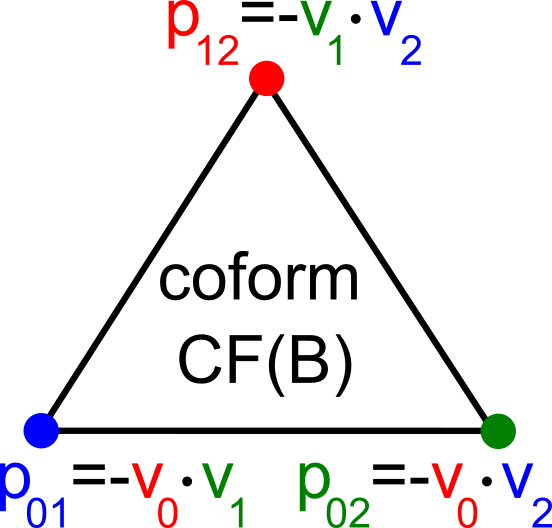}
\hspace*{2mm}
\includegraphics[height=35mm]{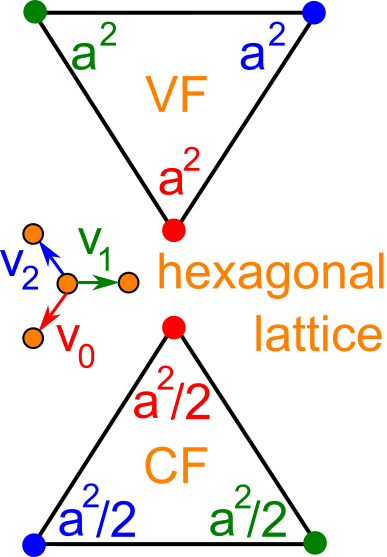}
\hspace*{2mm}
\includegraphics[height=35mm]{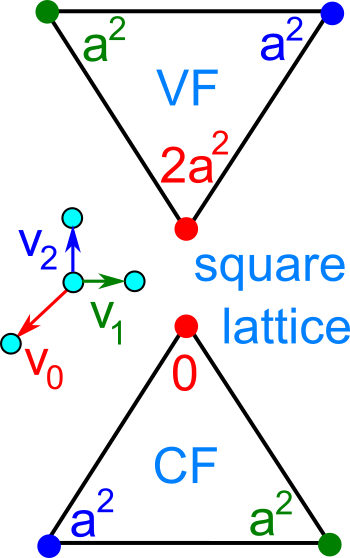}
\caption{\textbf{1st picture}: a voform $\VF(B)$ of a 2D lattice with an obtuse superbase $B=\{\vec v_0,\vec v_1,\vec v_2\}$.
\textbf{2nd picture}: nodes of a coform $\CF(B)$ are marked by conorms $p_{ij}$.
\textbf{3rd and 4th pictures}: $\VF$ and $\CF$ of the hexagonal and square lattice with a minimum inter-point distance $a$.}
\label{fig:forms2d}
\end{figure}

\index{root invariant} 

\begin{dfn}[voform $\VF$, coform $\CF$, ordered \emph{root invariant} $\RI$]
\label{dfn:RI}
For any ordered obtuse superbase $B$ in $\R^2$, the \emph{voform} $\VF(B)$ is the cycle on three nodes marked by the vonorms $\vec v_0^2,\vec v_1^2,\vec v_2^2$, see Fig.~\ref{fig:forms2d}.
The \emph{coform} $\CF(B)$ is the cycle on three nodes marked by the conorms $p_{12},p_{02},p_{01}$.
Since all conorms $p_{ij}\geq 0$, we can define the \emph{root products} $r_{ij}=\sqrt{p_{ij}}$.
The \emph{root invariant} $\RI(B)$ is obtained by writing the three root products $r_{12},r_{01},r_{02}$ in the increasing order.
\edfn
\end{dfn}

By Theorem~\ref{thm:reduction_to_obtuse_superbase} any lattice $\La\subset\R^2$ has an obtuse superbase with all $p_{ij}\geq 0$.
At least two root products $r_{ij}$ should be positive, otherwise one vonorm vanishes, but there are no other restrictions on $r_{ij}\geq 0$.
The vonorms $\vec v_0^2,\vec v_1^2,\vec v_2^2>0$ should satisfy three triangle inequalities such as $\vec v_0^2\leq \vec v_1^2+\vec v_2^2$, only one of them can be an equality.
The ordering $r_{12}\leq r_{01}\leq r_{02}$ is equivalent to $\vec v_1^2\leq \vec v_2^2\leq \vec v_0^2$ by formulae~(\ref{dfn:vonorms}a).
Root products have the same units as original coordinates of basis vectors, for example, Angstroms: $1\AA=10^{-10}$m.
The ordered root invariant $\RI(B)$ is more convenient than $\VF(B)$ and $\CF(B)$, which depend on an order of vectors of $B$.

\begin{figure}[h]
\includegraphics[width=1.0\textwidth]{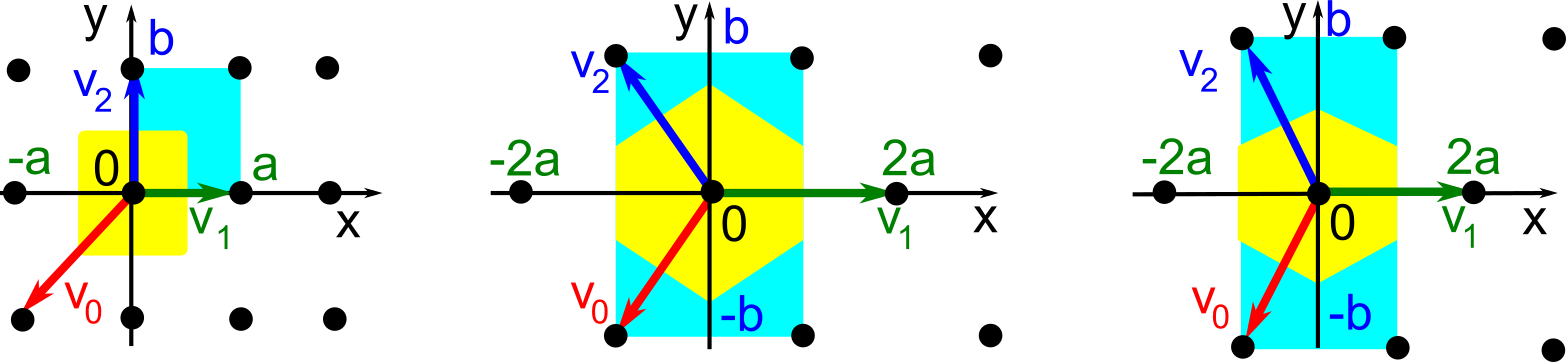}
\caption{
\textbf{Left}: $\La$ has a rectangular cell and obtuse superbase $B$ with $\vec v_1=(a,0)$, $\vec v_2=(0,b)$, $\vec v_0=(-a,-b)$, see Example~\ref{exa:achiral_lattices} and Lemma~\ref{lem:achiral_lattices}.
Other lattices $\La$ have a rectangular cell $2a\times 2b$ and an obtuse superbase $B$ with $\vec v_1=(2a,0)$, $\vec v_2=(-a,b)$, $\vec v_0=(-a,-b)$.
\textbf{Middle}: $\RI(B)=(\sqrt{b^2-a^2}, a\sqrt{2},a\sqrt{2})$, $a\leq b\leq a\sqrt{3}$.
\textbf{Right}: $\RI(B)=(a\sqrt{2},a\sqrt{2},\sqrt{b^2-a^2})$, $a\sqrt{3}\leq b$.
}
\label{fig:achiral_lattices}
\end{figure}

\begin{exa}
\label{exa:achiral_lattices}
\textbf{(a)} 
A lattice $\La$ with a rectangular cell of sides $a\leq b$ has an obtuse superbase $B$ with $\vec v_1=(a,0)$, $\vec v_2=(0,b)$, $\vec v_0=(-a,-b)$, and $\RI(B)=(0,a,b)$.
\medskip

\noindent
\textbf{(b)} 
For any lattice $\La\subset\R^2$ whose Voronoi domain $\bar V(\La)$ is a mirror-symmetric hexagon, assume that the $x$-axis is its line of symmetry.
Since $\bar V(\La)$ is centrally symmetric with respect to the origin $0$, the $y$-axis is also its line of symmetry, see Fig.~\ref{fig:achiral_lattices}.
\sskip

Then $\La$ has the centred rectangular (non-primitive) cell with sides $2a\leq 2b$.
The obtuse superbase $B$ with $\vec v_1=(2a,0)$, $\vec v_2=(-a,b)$, $\vec v_0=(-a,-b)$ has $\RI(B)=(a\sqrt{2},a\sqrt{2},\sqrt{b^2-a^2})$ for $b\geq a\sqrt{3}$.
For $a\leq b<a\sqrt{3}$, we should swap $r_{02}=\sqrt{b^2-a^2}$ with $r_{12}=a\sqrt{2}$ to get an ordered root invariant $\RI(B)$.
\eexa
\end{exa}

A lattice $\La\subset\R^n$ that can be mapped to itself by a mirror reflection with respect to a $(n-1)$-dimensional hyperspace can be called \emph{mirror-symmetric} or \emph{achiral}.
Since a mirror reflection of any lattice $\La\subset\R^2$ with respect to a line $L\subset\R^2$ can be realised by a rotation in $\R^3$ around $L$ through $180^{\circ}$, the term \emph{achiral} sometimes applies to all 2D lattices and becomes non-trivial only for 3D lattices.
This paper for 2D lattices uses the clearer adjective \emph{mirror-symmetric}.

\index{root invariant} 
\index{Voronoi domain}
 
\begin{lem}[root invariants of mirror-symmetric lattices $\La\subset\R^2$,
{\cite[Lemma 3.3]{kurlin2024mathematics}}]
\label{lem:achiral_lattices}
An obtuse superbase $B$ generates a \emph{mirror-symmetric} lattice $\La(B)$ \emph{if and only if} 
\smallskip

\noindent
\tb{(a)} 
the root invariant $\RI(B)$ contains a zero value and $\La(B)$ is rectangular, or
\smallskip

\noindent
\tb{(b)}
$\RI(B)$ has equal root products and the Voronoi domain of $\La(B)$ is a square or a hexagon whose symmetry group has two orthogonal axes of symmetry.
\elem
\end{lem}

\begin{figure}[h]
\label{fig:root_forms2d_reflection}
\caption{The lattices $\La,\La'$ are mirror reflections of each other and have oriented root invariants $\RI^o=(\sqrt{3},\sqrt{6},\sqrt{7})_{\pm}$ with opposite signs introduced in Definition~\ref{dfn:sign}, see Example~\ref{exa:signs}.} 
\includegraphics[width=\textwidth]{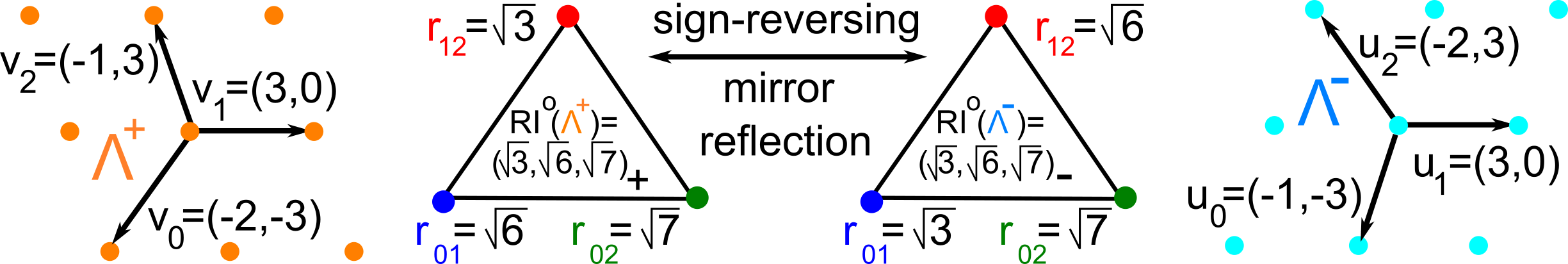}
\end{figure}

\begin{dfn}[$\sign(B)$, the oriented root invariant $\RI^o(B)$]
\label{dfn:sign}
If an obtuse superbase $B$ generates a mirror-symmetric lattice, set $\sign(B)=0$.
Else all vectors of $B$ have different lengths and angles not equal to $90^{\circ}$ by Lemma~\ref{lem:achiral_lattices}.
Let $\vec v_1,\vec v_2$ be the shortest vectors of $B$ so that $|\vec v_1|<|\vec v_2|$.
Then $\sign(B)=\pm 1$ is the sign of the determinant $\det(\vec v_1,\vec v_2)$ of the matrix with the columns $\vec v_1,\vec v_2$.
The \emph{oriented root invariant} $\RI^o(B)$ is obtained by adding $\sign(B)$ as a superscript to $\RI(B)$, see Fig.~\ref{fig:root_forms2d_reflection}. 
\edfn
\end{dfn}

If $\sign(B)=0$, this zero superscript in $\RI^o(B)$ can be skipped for simplicity, so $\RI^o(B)=\RI(B)$ in this case.
Theorem~\ref{thm:isometric_superbases} will show that $\sign(B)$ can be considered as an invariant of a lattice $\La$ under dilation.
\medskip

In Definition~\ref{dfn:sign} the determinant $\det(\vec v_1,\vec v_2)$ is the signed area of the unit cell $U(\vec v_1,\vec v_2)$ equal to $|\vec v_1|\cdot|\vec v_2|\sin\angle(\vec v_1,\vec v_2)$, where the angle is measured from $\vec v_1$ to $\vec v_2$ in the anticlockwise direction around the origin $0\in\R^2$.
For a strict obtuse superbase $B$, all angles between its basis vectors are strictly obtuse.
Then $\sign(B)=+1$ if $\angle(\vec v_1,\vec v_2)$ is in the positive range $(90^\circ,180^\circ)$, else $\sign(B)=-1$.

\begin{exa}[signs of lattices in Fig.~\ref{fig:root_forms2d_reflection}]
\label{exa:signs}
\tb{(a)}
The lattice $\La^+$ in the first picture of Fig.~\ref{fig:root_forms2d_reflection} has the obtuse superbase $B$ with $\vec v_1=(3,0)$, $\vec v_2=(-1,3)$, $\vec v_0=(-2,-3)$ of lengths $3,\sqrt{10},\sqrt{13}$, respectively, so $\La^+$ is not mirror-symmetric.
Since $\vec v_1,\vec v_2$ are the two shortest vectors of $B^+$ and $\det(\vec v_1,\vec v_2)=\det\matfour{3}{-1}{0}{3}>0$, we get $\sign(B^+)=+1$.
The anticlockwise angle is $\angle(\vec v_1,\vec v_2)=180^\circ-\arcsin\frac{3}{\sqrt{10}}\approx 108^\circ$.
\medskip

\nt
\tb{(b)}
The lattice $\La^-$ in the last picture of Fig.~\ref{fig:root_forms2d_reflection} is obtained from $\La^+$ by a mirror reflection and has the obtuse superbase $B^-$ with $\vec u_1=\vec v_1$, $\vec u_2=(-2,3)$, $\vec u_0=(-1,-3)$ of lengths
$3,\sqrt{13},\sqrt{10}$, respectively, so $\La^-$ is not mirror-symmetric. 
Since $\vec u_1,\vec u_0$ are the shortest vectors, $\det(\vec u_1,\vec u_0)=\det\matfour{3}{-1}{0}{-3}<0$, we get $\sign(B^-)=-1$.
The anticlockwise angle is $\angle(\vec u_1,\vec u_0)=\arcsin\frac{3}{\sqrt{10}}-180^\circ\approx -108^\circ$.
\eexa
\end{exa}

Theorem~\ref{thm:isometric_superbases}  below is crucial for a complete classification of 2D lattices in Theorem~\ref{thm:lattices2D_classification} and Corollary~\ref{cor:lattices2D_classification}. 
Theorem~\ref{thm:isometric_superbases} highlights that mirror-symmetric lattices have more options for obtuse superbases under rigid motion.
The same rectangular lattice can have two obtuse bases with $\vec v_1=(1,0)$, $\vec v_2=(0,\pm 2)$, which are related by reflection in the $x$-axis, not by rigid motion.
This symmetry-related ambiguity is much harder to resolve for 3D lattices even under isometry, see \cite{kurlin2022complete}.

\index{obtuse superbase}

\begin{thm}[isometric obtuse superbases, {\cite[Theorem 3.7]{kurlin2024mathematics}}]
\label{thm:isometric_superbases}  
Any lattices $\La,\La'\subset\R^2$ are isometric if and only if any obtuse superbases of $\La,\La'$ are isometric.
If $\La,\La'$ are not rectangular, the same conclusion holds for rigid motion instead of isometry.
Any rectangular (non-square) lattice has two obtuse superbases related by reflection.
\ethm
\end{thm}

\index{obtuse superbase}

\begin{lem}[lattice invariants, {\cite[Lemma 3.8]{kurlin2024mathematics}}]
\label{lem:lattice_invariants}
\textbf{(a)}
For any obtuse superbase $B$ of a lattice $\La\subset\R^2$,
 the root  invariant $\RI(B)$ is an isometry invariant of $\La$ and can be denoted by $\RI(\La)$.
Similarly, $\RI^o(\La)$ and $\sign(\La)$ are invariants of a lattice $\La$ under rigid motion and dilation, respectively.
\medskip

\noindent
\textbf{(b)}
A lattice $\La\subset\R^2$ is mirror-symmetric if and only if $\sign(\La)=0$.
\elem
\end{lem}

\begin{figure}[h!]
\includegraphics[width=\textwidth]{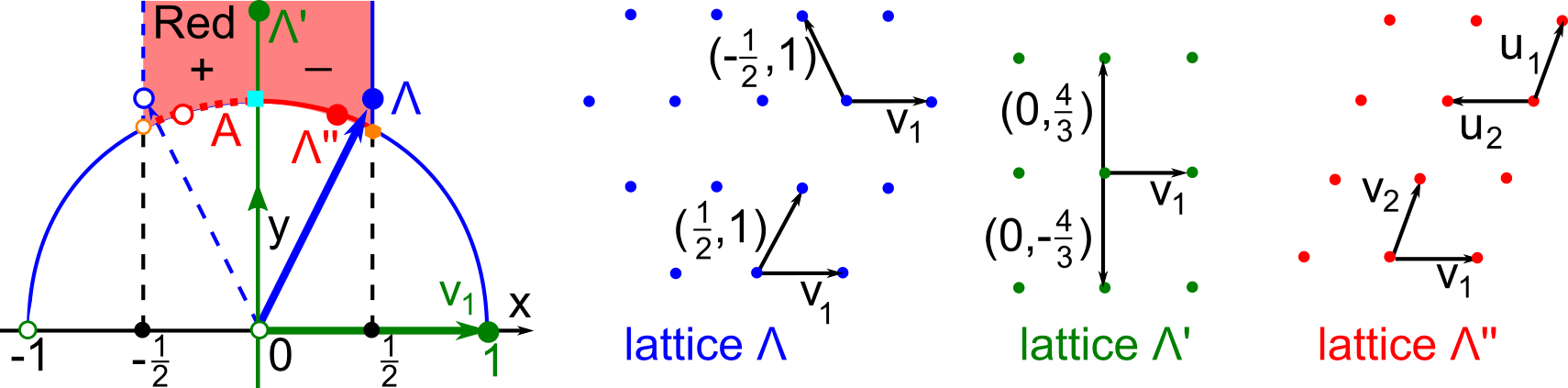}
\caption{\textbf{Left}: any reduced basis in Definition~\ref{dfn:reduced_cell} can be mapped under dilation to the basis of $\vec v_1=(1,0)$ and $\vec v_2=(x,y)\in\RB$ from Proposition~\ref{prop:reduced_bases}.
\textbf{Right}: for each of the lattices $\La,\La',\La''$ represented by small blue, green, red circles/disks on the right, the conditions of Definition~\ref{dfn:reduced_cell} choose one reduced basis among two bases that differ under rigid motion.
}
\label{fig:reduced_bases}
\end{figure}

\index{reduced basis}

\begin{prop}[reduced bases, {\cite[Proposition 3.10]{kurlin2024mathematics}}]
\label{prop:reduced_bases}
\textbf{(a)}
Under isometry in $\R^2$, all reduced bases $\vec v_1,\vec v_2$ from Definition~\ref{dfn:reduced_cell} are in a 1-1 correspondence with all obtuse superbases $B=\{\vec v_0,\vec v_1,\vec v_2\}$ such that $|\vec v_1|\leq|\vec v_2|\leq|\vec v_0|$.
Under isometry, any lattice $\La\subset\R^2$ has a unique reduced basis in the conditions of Definition~\ref{dfn:reduced_cell}, see Fig.~\ref{fig:reduced_bases}.
\medskip

\noindent
\textbf{(b)}
Under rigid motion, any lattice has a unique reduced basis in Definition~\ref{dfn:reduced_cell}.
\ethm
\end{prop}

\section{Complete invariants of 2D lattices under four equivalences}
\label{sec:lattices2D_invariants}

Lemma~\ref{lem:lattice_invariants} showed that $\RI(\La),\RI^o(\La)$ are invariants of lattices under isometry and rigid motion, respectively. 
To prove completeness of the invariants in Theorem~\ref{thm:lattices2D_classification}, Lemma~\ref{lem:superbase_reconstruction} reconstructs an obtuse superbase of $\La$.
Corollary~\ref{cor:lattices2D_classification} will classify lattices under homothety by projected invariants introduced in Definition~\ref{dfn:PI}.

\begin{lem}[superbase reconstruction, {\cite[Lemma 4.1]{kurlin2024mathematics}}]
\label{lem:superbase_reconstruction}
An obtuse superbase $B=\{\vec v_0,\vec v_1,\vec v_2\}$ of a lattice $\La\subset\R^2$ can be uniquely reconstructed under isometry and under rigid motion from its root invariant $\RI(\La)$ and its oriented root invariant $\RI^o(\La)$, respectively.
If $\RI(\La)=(r_{12},r_{01},r_{02})$, the basis vectors $\vec v_1,\vec v_2$ are determined by 
$$|\vec v_1|=\sqrt{r_{12}^2+r_{01}^2},\quad 
|\vec v_2|=\sqrt{r_{12}^2+r_{02}^2},\quad
\cos\angle(\vec v_1,\vec v_2)=\dfrac{-r_{12}^2}{\sqrt{r_{12}^2+r_{01}^2}\sqrt{r_{12}^2+r_{02}^2}},
$$ and span a primitive unit cell of the area $A(\La)=\sqrt{r_{12}^2 r_{01}^2 + r_{12}^2 r_{02}^2 + r_{01}^2 r_{02}^2}$.
\elem
\end{lem}

\index{rigid motion}
\index{isometry}

\begin{thm}[completeness of root invariants, {\cite[Theorem~4.2]{kurlin2024mathematics}}]
\label{thm:lattices2D_classification}
\tb{(a)}
Any lattices $\La,\La'\subset\R^2$ are isometric if and only if their root invariants coincide: $\RI(\La)=\RI(\La')$.
\myskip

\nt
\tb{(b)}
Any lattices $\La,\La'$ are related by rigid motion
if and only if $\RI^o(\La)=\RI^o(\La')$.
\ethm 
\end{thm}

The above classification helps prove that some other isometry invariants of lattices are also complete and continuous.
By (\ref{dfn:vonorms}ab) the voform $\VF=(\vec v_0^2,\vec v_1^2,\vec v_2^2)$ and coform $\CF=(p_{12},p_{01},p_{02})$ are both complete if considered under $3!$ permutations.
The root invariant $\RI$ is a uniquely ordered version of $\CF$ and deserves its own name.
The square roots $r_{ij}=\sqrt{p_{ij}}$ have original units of the vector coordinates.
\medskip

\begin{figure}[h]
\includegraphics[height=34mm]{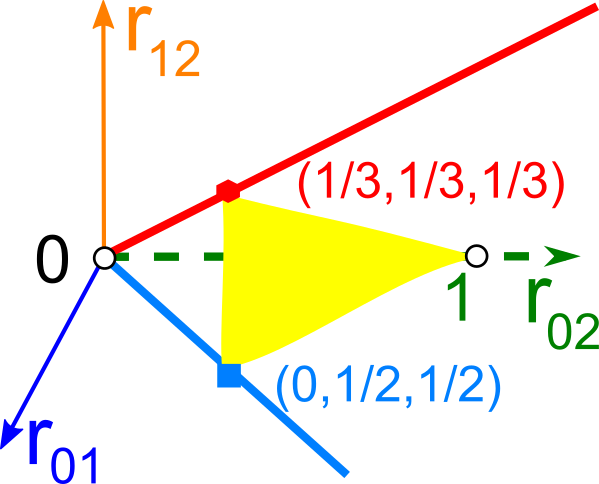}\hspace*{2mm}
\includegraphics[height=34mm]{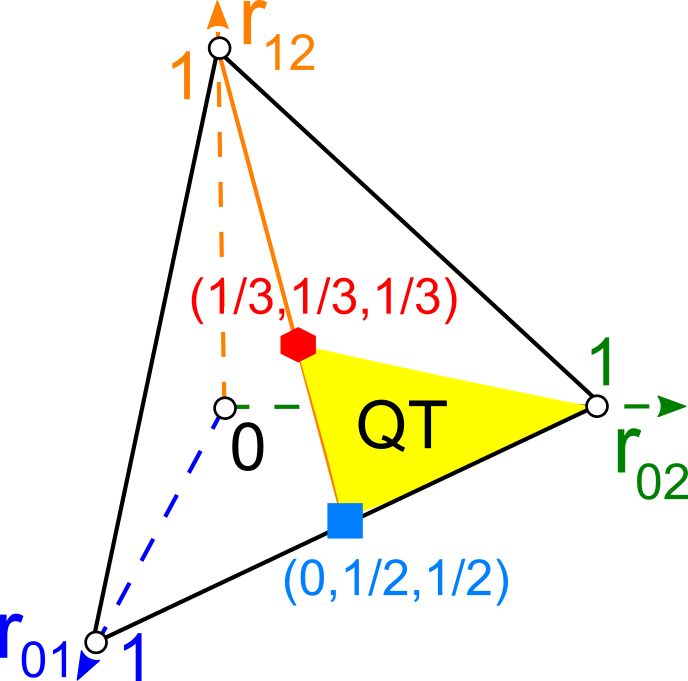}\hspace*{2mm}
\includegraphics[height=34mm]{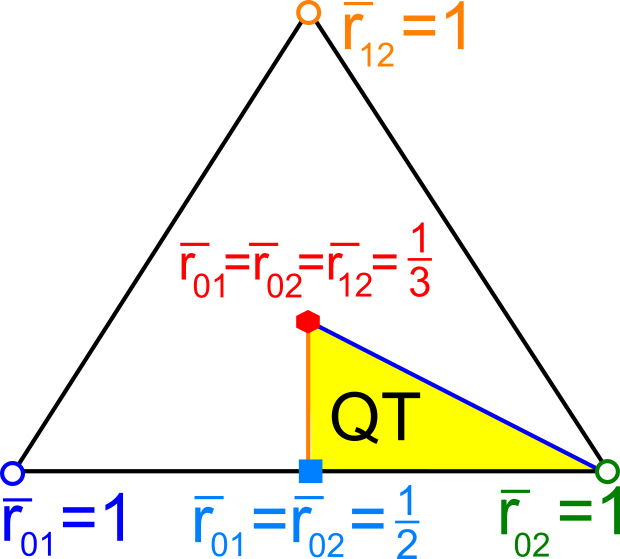}
\caption{\textbf{Left}: the triangular cone $\TC=\{(r_{12},r_{01},r_{02})\in\R^3 \mid 0\leq r_{12}\leq r_{01}\leq r_{02}\neq 0\}$ represents the space $\RIS$ of all root invariants of 2D lattices, see Definition~\ref{dfn:TC}. 
\textbf{Middle}: 
$\TC$ projects to the quotient triangle $\QT=\TC\cap\{r_{12}+r_{01}+r_{02}=1\}$ representing the space $\LHS$ of 2D lattices under homothety, see Corollary~\ref{cor:lattices2D_classification}.
\textbf{Right}: the quotient triangle $\QT$ can be parametrised by $x=\bar r_{02}-\bar r_{01}\in[0,1)$ and $y=3\bar r_{12}\in[0,1]$, see $\QT$ also in Fig.~\ref{fig:QT+QS}.
}
\label{fig:TC}
\end{figure}

\begin{figure}[h]
\includegraphics[height=55mm]{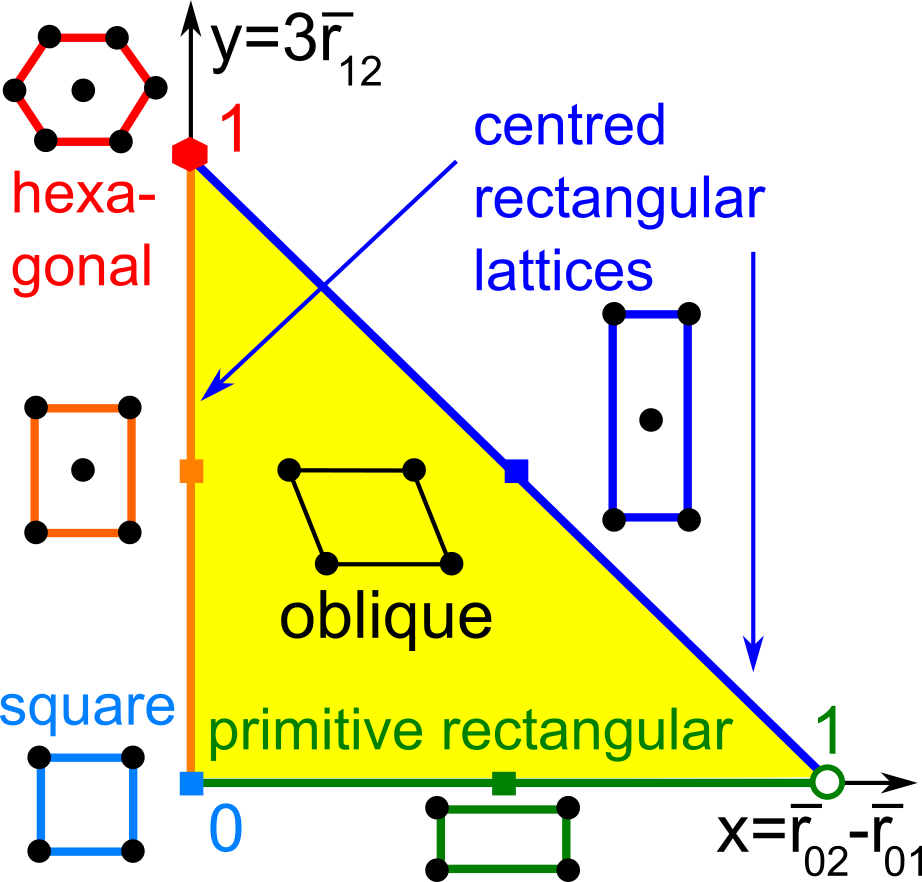}
\hspace*{3mm}
\includegraphics[height=55mm]{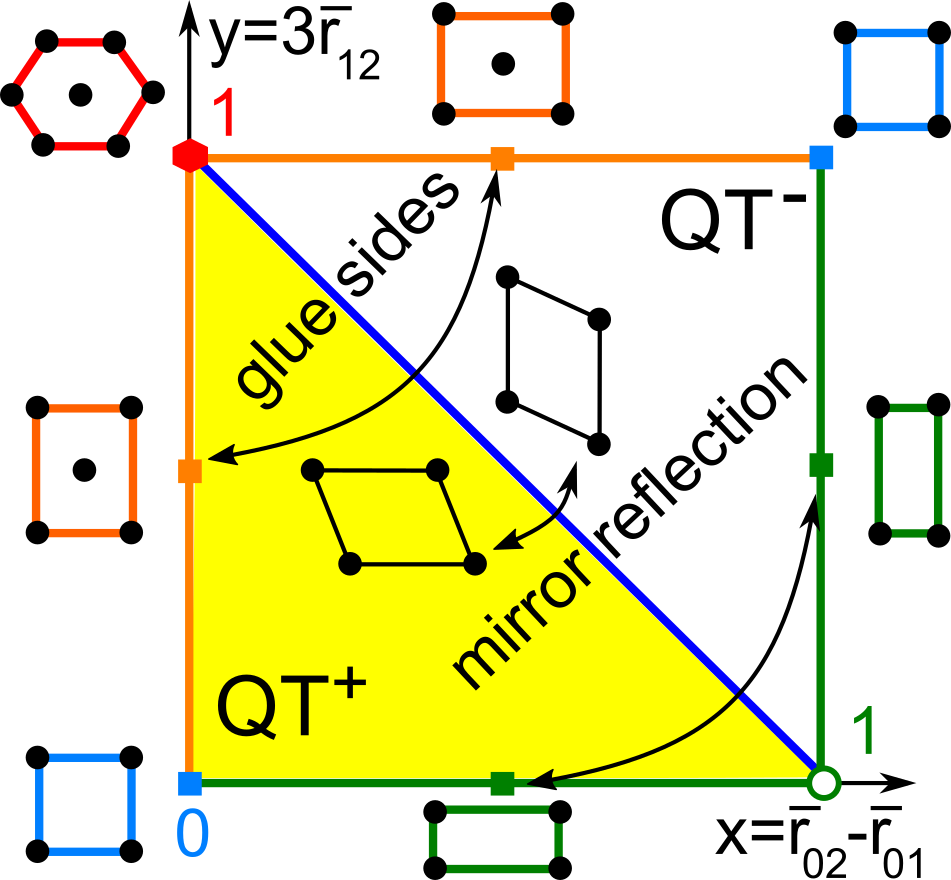}
\caption{\textbf{Left}:
all projected invariants $\PI(\La)$ of lattices $\La\subset\R^2$ live in the quotient triangle $\QT$ from Fig.~\ref{fig:TC}, which is parameterised by $x=\bar r_{02}-\bar r_{01}\in[0,1)$ and $y=3\bar r_{12}\in[0,1]$.
\textbf{Right}: mirror reflections $\La^{\pm}$ of any non-mirror-symmetric lattice can be represented by a pair of points in the quotient square $\QS=\QT^+\cup\QT^-$ symmetric in the diagonal $x+y=1$.}
\label{fig:QT+QS}
\end{figure}

Theorem~\ref{thm:lattices2D_classification} and Lemma~\ref{lem:VF} imply that, after taking square roots of vonorms, the ordered lengths, say $|\vec v_1|\leq|\vec v_2|\leq|\vec v_0|$, form a complete invariant that should satisfy the triangle inequality $|\vec v_1|+|\vec v_2|\geq|\vec v_0|$.
This inequality is the only disadvantage of the complete invariant $|\vec v_1|\leq|\vec v_2|\leq|\vec v_0|$ in comparison with ordered root products $r_{12}\leq r_{01}\leq r_{02}$, which are easier to visualise in Fig.~\ref{fig:TC},~\ref{fig:QT+QS}.
\medskip

Classification Theorem~\ref{thm:lattices2D_classification} says that all isometry classes of lattices $\La\subset\R^2$ are in a 1-1 correspondence with all ordered triples $0\leq r_{12}\leq r_{01}\leq r_{02}$ of root products in $\RI(\La)$.
Only the smallest root product $r_{12}$ can be zero, two others $r_{01}\leq r_{02}$ should be positive, otherwise $\vec v_1^2=r_{12}^2+r_{01}^2=0$ by formulae (\ref{dfn:vonorms}a).
\medskip

We explicitly describe the set of all possible root invariants, which will be later converted into metric spaces with continuous metrics in Definitions~\ref{dfn:RM} and~\ref{dfn:RMo}.

\index{triangular cone } 
\index{quotient triangle}
\index{quotient square}

\begin{dfn}[triangular cone $\TC$]
\label{dfn:TC}
All root invariants $\RI(\La)=(r_{12},r_{01},r_{02})$ of lattices $\La\subset\R^2$ live in the \emph{triangular cone} $\TC=\{0\leq r_{12}\leq r_{01}\leq r_{02}\}$ within the octant $\Oct=[0,+\infty)^3$ excluding the axes in the coordinates $r_{12},r_{01},r_{02}$, see Fig.~\ref{fig:TC}~(left).
\sskip

The boundary $\bd(\TC)$ of the cone $\TC$ consists of root invariants of all mirror-symmetric lattices from Lemma~\ref{lem:achiral_lattices}: 
the bisector planes $\{r_{01}=r_{02}\}$ and $\{r_{12}=r_{01}\}$ within $\TC$. 
The orange line $\{0<r_{12}=r_{01}=r_{02}\}\subset\bd(\TC)$ in Fig.~\ref{fig:TC}~(left) consists of root invariants of hexagonal lattices with a minimum inter-point distance $r_{12}\sqrt{2}$.
The blue line $\{r_{12}=0<r_{01}=r_{02}\}\subset\bd(\TC)$ consists of root invariants of square lattices with a minimum inter-point distance $r_{01}$.
\edfn
\end{dfn}

To classify lattices under homothety, it is convenient to scale them by the \emph{size} $|\La|=r_{12}+r_{01}+r_{02}$.
This sum is a simpler uniform measure of size than (say) the unit cell area $A(\La)$ from Lemma~\ref{lem:superbase_reconstruction}, which can be small even for long cells. 

\index{projected invariant} 

\begin{dfn}[projected invariants $\PI(\La)$ and $\PI^o(\La)$]
\label{dfn:PI}
The \emph{triangular projection} $\TP:\TC\to\{r_{12}+r_{01}+r_{02}=1\}$ divides each coordinate by the \emph{size} $|\La|=r_{12}+r_{01}+r_{02}$ and gives $\PRF(\La)=(\bar r_{12},\bar r_{01},\bar r_{02})=\dfrac{(r_{12},r_{01},r_{02})}{r_{12}+r_{01}+r_{02}}$ in $\TC\cap\{r_{12}+r_{01}+r_{02}=1\}$.
Then we map $(\bar r_{12},\bar r_{01},\bar r_{02})$ to the \emph{projected invariant} $\PI(\La)=(x,y)$ with $x=\bar r_{02}-\bar r_{01}\in[0,1)$ and $y=3\bar r_{12}\in[0,1]$ in the \emph{quotient triangle} 
$$\QT=\{(x,y)\in\R^2\mid 0\leq x<1,\, 0\leq y\leq 1,\, x+y\leq 1\},$$ see Fig.~\ref{fig:QT+QS}.
All oriented root invariants $\RI^o(\La)$ live in the \emph{doubled cone} $\DC$ that is the union of two triangular cones $\TC^{\pm}$, where we identify any two boundary points representing the same root invariant $\RI(\La)$ with $\sign(\La)=0$.
The \emph{oriented projected invariant} $\PI^o(\La)=(x,y)^{\pm}$ is $\PI(\La)$ with the superscript from $\sign(\La)$.
\edfn
\end{dfn}

The inequality $1\geq x+y=(\bar r_{02}-\bar r_{01})+3\bar r_{12}$ follows after multiplying both sides by the size $|\La|$, because $r_{12}+r_{01}+r_{02}\geq (r_{02}-r_{01})+3r_{12}$ becomes $r_{01}\geq r_{12}$.
\medskip

The set of oriented projected invariants $\PI^o$ is visualised in Fig.~\ref{fig:QT+QS} (right) as the \emph{quotient square} $\QS$ obtained by gluing the quotient triangle $\QT^+$ with its mirror image $\QT^-$.
The boundaries of both triangles excluding the vertex $(x,y)=(1,0)$ are glued by the diagonal reflection $(x,y)\lra(1-y,1-x)$.
Any pair of points $(x,y)\in\QT^+$ and $(1-y,1-x)\in\QT^-$ in Fig.~\ref{fig:QT+QS}~(right) represent mirror images of a lattice under homothety, see Corollary~\ref{cor:lattices2D_classification}.
So $\QS$ is a topological sphere without a single point and will be parameterised by geographic-style coordinates in \cite{bright2023geographic}.
\medskip

Following Fig.~\ref{fig:forms2d}, any square lattice 
has a root invariant $\RI=(0,a,a)$, so its projected invariant $\PI=(0,0)$ is at the bottom left vertex of $\QT$ in Fig.~\ref{fig:QT+QS}~(left), identified with top right vertex of $\QS$ in Fig.~\ref{fig:QT+QS}~(right).
Any hexagonal lattice has a root invariant $\RI=(a,a,a)$, so its projected invariant $\PI=(0,1)$ is at the top left vertex of $\QT$ in Fig.~\ref{fig:QT+QS}~(left), identified with bottom right vertex of $\QS$.
\medskip

By Example~\ref{exa:achiral_lattices}(a) any rectangular lattice has $\RI=(0,a,b)$ for $a<b$, hence its projected invariant $\PI=(\frac{b-a}{a+b},0)$ belongs to the bottom edge of $\QT$ identified with the top edge of $\QS$.
By Example~\ref{exa:achiral_lattices}(b) any lattice with a mirror-symmetric Voronoi domain has $\RI$ with 0 or two equal root products.
Such lattices have a rhombic unit cell and form the centred rectangular Bravais class.
Their projected invariants belong to the vertical edges and diagonal of $\QS$ in Fig.~\ref{fig:QT+QS}~(right).
The companion paper \cite{bright2023geographic} discusses Bravais classes of 2-dimensional lattices in detail.

\index{projected invariant}
\index{dilation}
\index{homothety}

\begin{cor}[completeness of $\PI$, {\cite[Corollary 4.6]{kurlin2024mathematics}}]
\label{cor:lattices2D_classification}
Any lattices $\La,\La'\subset\R^2$ are related by homothety if and only if their projected invariants are equal: $\PI(\La)=\PI(\La')$.
Any lattices $\La,\La'$ are related by dilation if and only if $\PI^o(\La)=\PI^o(\La')$.
\ethm
\end{cor}

\begin{lem}[criteria of mirror-symmetric lattices in $\R^2$, {\cite[Lemma 4.7]{kurlin2024mathematics}}]
\label{lem:sign0}
A lattice $\La$ in $\R^2$ is mirror-symmetric if and only if one of the following equivalent conditions holds:
$\sign(\La)=0$ or $\RI(\La)\in\bd\TC$ or $\PI(\La)\in\bd\QT$.
So the boundaries of the triangular cone $\TC$ and the quotient triangle $\QT$ consist of root invariants and projected invariants, respectively, of all mirror-symmetric lattices $\La\subset\R^2$.
\elem
\end{lem}

\section{Inverse design and a spherical map of 2D lattices}
\label{sec:lattices2D_spaces}

This section discusses lattices in terms of group actions, inversely designs lattices from invariants, and embeds the Lattice Rigid Space $\LRS(\R^2)$ in a 2-dimensional sphere.
\myskip    

In the theory of complex functions, any lattice $\La\subset\R^2$ can be considered as a subgroup of the complex plane $\C$ whose quotient $\C/\La$ is a torus.
By the Riemann mapping theorem any compact Riemann surface of genus 1 is conformally equivalent (holomorphically homeomorphic) to the quotient $\C/\La$ for some lattice $\La$, see \cite[Section 5.3]{jost2013compact}.
These tori $\C/\La$ and $\C/\La'$ are  \emph{conformally equivalent} if and only if $\La,\La'$ are related by homothety, see \cite[Theorem~6.1.4]{jones1987complex}.
The spaces $\LHS(\R^2)$ and $\LDS(\R^2)$ of all lattices $\La\subset\C=\R^2$ under homothety and dilation are the quotient triangle $\QT$ and the quotient square $\QS$, respectively, see Fig.~\ref{fig:QT+QS}.

\begin{rem}[lattices via group actions]
\label{rem:group_action}
Another parametrisation of the Lattice Homothety Space $\LHS(\R^2)$ can be obtained from a fundamental domain of the action of $\GL(\Z^2)\times\R_+^\times$ on the cone $\mathcal{C}_+(\mathcal{Q}_2)$ of positive quadratic forms.
Recall that any lattice $\La\subset\R^2$ with a basis $\vec v_1,\vec v_2$ defines the positive quadratic form
$$Q_\La(x,y)=(x\vec v_1+y\vec v_2)^2=\vec v_1^2x^2+2\vec v_1\vec v_2xy+\vec v_2^2y^2=q_{11}x^2+2q_{12}xy+q_{22}y^2\geq 0$$ 
whose positivity for all $(x,y)\in\R^2-0$ means that $q_{12}^2<q_{11}q_{22}$.
The cone $\mathcal{C}_+(\mathcal{Q}_2)$ of all positive quadratic forms projects to the unit disk $\xi^2+\eta^2<1$ parameterised by $\xi=\dfrac{q_{22}-q_{11}}{q_{11}+q_{22}}$ and $\eta=\dfrac{-2q_{12}}{q_{11}+q_{22}}$.
Indeed, the positivity 
condition $q_{12}^2<q_{11}q_{22}$ for the form $Q_\La(x,y)$ is equivalent to $\xi^2+\eta^2<1$ in the coordinates above.
\medskip

The quadratic form $Q_\La$ is in a \emph{reduced} (non-acute) form if $0\leq -2q_{12}\leq q_{11}\leq q_{22}$ and $q_{11}>0$, see \cite[formula~(1.130) on p.~75]{engel2004lattice}.
The above conditions define the fundamental domain $T=\{0\leq\xi<1,\; 0\leq\eta\leq\frac{1}{2},\; \xi+2\eta\leq 1\}$, see  \cite[Fig.~8.1]{zhilinskii2016introduction}. 
This non-isosceles triangle is one of the infinitely many triangular domains within the disk $\xi^2+\eta^2<1$ in \cite[Fig.~1.2 on p.~82]{engel2004lattice} or \cite[Fig.~6.2]{zhilinskii2016introduction}.
Choosing one triangular domain is equivalent to choosing a reduced basis under isometry, not under rigid motion.
\myskip

For instance, the mirror-symmetric bases $\vec v_1=(1,0)$, $\vec v_2^{\pm}=(-\frac{1}{2},\pm 1)$ have the same reduced non-acute form $x^2-xy+\frac{5}{4}y^2$ represented only by $(\xi,\eta)=(\frac{1}{9},\frac{4}{9})$.
The above ambiguity under rigid motion is resolved by $\sign(\La)$ in the twice larger space $\LDS(\R^2)$ visualised as the quotient square $\QS$, see Example~\ref{exa:rect_deformation}.
\medskip

More importantly, the inverse map from $\RI^o(\La)$ to a reduced basis  is discontinuous at any rectangular lattice $\La$ with a unit cell $a\times b$.
Indeed, slight perturbations of $\La$ have unique reduced bases that are not close to each other, being close to the distant bases $(a,0),(0,\pm b)$, which are not equivalent under rigid motion for $a<b$.
This discontinuity of lattice bases will emerge in $\R^3$ even under isometry \cite{kurlin2022complete}.
In $\R^2$, Corollary~\ref{cor:rect_discontinuity} will completely settle the basis discontinuity under rigid motion.
\medskip

Another complete invariant is the ordered voform $\vec v_1^2\leq \vec v_2^2\leq \vec v_0^2$ or the lengths $|\vec v_1|\leq|\vec v_2|\leq|\vec v_0|$ of the three shortest Voronoi vectors from Lemma~\ref{lem:VF} below.
However, this invariant doesn't extend even to dimension $n=3$ due to a 6-parameter family of pairs of non-isometric lattices $\La_1\not\cong\La_2$ that have the same lengths of seven shortest Voronoi vectors in $\R^3$, see \cite{kurlin2022complete}.
The above reasons justify the choice of homogeneous coordinates $r_{ij}$, which easily extend to higher dimensions. 
\erem
\end{rem}

\begin{lem}[{\cite[Theorem~7]{conway1992low}}]
\label{lem:VF}
For any obtuse superbase $(v_0,v_1,v_2)$ of a lattice $\La\subset\R^2$, the vonorms $v_0^2,v_1^2,v_2^2$ are squared lengths of shortest Voronoi vectors. 
\elem
\end{lem}

The projected invariant $\PI=(x,y)$ obtained from $\RI$ is preferable to the coordinates $(\xi,\eta)$, which define a non-isosceles triangle, while the isosceles quotient triangle $\QT$ will lead to easier formulae for metrics in the next section.
Since the metric tensor $(\vec v_1^2,\vec v_1\cdot \vec v_2,\vec v_2^2)=(q_{11},q_{12},q_{22})$ and its 3-dimensional analogue are more familiar to crystallographers, we will rephrase key results from sections~\ref{sec:lattices2D_metrics}-\ref{sec:chiral_distances} by using these non-homogeneous coordinates in the companion paper \cite{bright2023geographic}.

\begin{prop}[inverse design of 2D lattices, {\cite[Proposition 4.9]{kurlin2024mathematics}}]
\label{prop:inverse_design}
For $s>0$ and any point $(x,y)$ in the quotient triangle $\QT$, there is a unique (under isometry) lattice $\La$ with the projected invariant $\PI(\La)=(x,y)$ and size $|\La|=s=r_{12}+r_{01}+r_{02}$. 
Then
$$\RI(\La)=(r_{12},r_{01},r_{02})=
\left(\frac{s}{3} y,\; \frac{s}{6}(3-3x-y),\; \frac{s}{6}(3+3x-y)\right).
\leqno{(\ref{prop:inverse_design}a)}$$ 
If $(x,y)$ is in the interior of $\QT$, the invariant $\RI$ defines a pair of lattices $\La^{\pm}$ that have opposite signs and 
unique (under isometry) reduced basis vectors $\vec v_1,\vec v_2$ with the lengths $|\vec v_1|=\sqrt{r_{12}^2+r_{01}^2}$, $|\vec v_2|=\sqrt{r_{12}^2+r_{02}^2}$ and the anticlockwise angle
$(\ref{prop:inverse_design}b)\qquad \angle(\vec v_1,\vec v_2)=\arccos\dfrac{-4y^2}{\sqrt{(9x^2+5y^2-6y+9)^2-36x^2(3-y)^2}}$.
\ethm
\end{prop}

Example~\ref{exa:inverse_design} shows the power of Proposition~\ref{prop:inverse_design} based on Theorem~\ref{thm:lattices2D_classification} and Corollary~\ref{cor:lattices2D_classification} for inverse design by sampling the square $\QS$ at interesting places.
\medskip

Fig.~\ref{fig:QS+DC}~(right) visualises the doubled cone $\DC$ of oriented root invariants $\RI^o$ from Definition~\ref{dfn:sign} 
 by uniting the triangular cone $\TC=\{0\leq r_{12}\leq r_{01}\leq r_{02}\}$ with its mirror reflection in the vertical plane $\{r_{01}=r_{02}\}$ including the $r_{12}$-axis.
\medskip

The lattice $\La_0$ with $\RI=(1,1,4)$ is represented by two boundary points of $\DC$ identified by $(r_{01},r_{02})\lra(r_{02},r_{01})$.
The lattices $\La_{\infty}^{\pm}$ with the root invariant $\RI=(r_{12},r_{01},r_{02})=(1,4,7)$ are represented by $(1,4,7)$ and its mirror image $(1,7,4)$ in $\DC$ related by the reflection in the vertical bisector plane $r_{01}=r_{02}$ containing the root invariants of $\La_4,\La_6$.
The superscript shows $\sign(\La_{\infty}^{\pm})=\pm 1$. 
  

\begin{exa}[inverse design of 2D lattices]
\label{exa:inverse_design}
We will inversely design the lattices $\La_4,\La_6,\La_0,\La_2^{\pm},\La_{\infty}^\pm$, see their visualised invariants in Fig.~\ref{fig:QS+DC}~(right).
\medskip

\noindent
\emph{($\pmb{\La_4}$)}
We design the square lattice $\La_4$ starting from its projected invariant  at the origin $\PI(\La_4)=(0,0)\in\QT$, which is identified with the top right vertex $(1,1)\in\QS$ in Fig.~\ref{fig:QS+DC}~(left).
Formula~(\ref{prop:inverse_design}a) for the size $|\La_4|=2$ (only to get simplest integers) gives $\RI(\La_4)=(0,1,1)$.
An obtuse superbase $\{\vec v_0,\vec v_1,\vec v_2\}$ can be reconstructed by Lemma~\ref{lem:superbase_reconstruction}.
The vonorms are $\vec v_1^2=\vec v_2^2=0^2+1^2=1$, $\vec v_0^2=1^2+1^2=2$.
We can choose the standard obtuse superbase $\vec v_1=(1,0)$, $\vec v_2=(0,1)$, $\vec v_0=(-1,-1)$.
\medskip

\noindent
\emph{($\pmb{\La_6}$)}
We design the hexagonal lattice $\La_6$ starting from the projected invariant at the top left vertex $\PI(\La_6)=(0,1)\in\QT$, which is identified with the bottom right vertex $(1,0)\in\QS$ in Fig.~\ref{fig:QS+DC}~(left).
Formula~(\ref{prop:inverse_design}a) for the size $|\La_6|=3$ (only to get simplest integers) gives $\RI(\La_6)=(1,1,1)$.
To reconstruct an obtuse superbase $\{\vec v_0,\vec v_1,\vec v_2\}$ by Lemma~\ref{lem:superbase_reconstruction}, find the vonorms $\vec v_1^2=\vec v_2^2=\vec v_0^2=1^2+1^2=2$.
Formula~(\ref{prop:inverse_design}b) gives the angle $\angle(\vec v_1,\vec v_2)=\arccos
\frac{-4}{\sqrt{(5-6+9)^2}}=
\arccos\left(-\frac{1}{2}\right)=120^\circ$.
We can choose the superbase $\vec v_1=(\sqrt{2},0)$, $\vec v_2=(-\frac{1}{\sqrt{2}},\frac{\sqrt{3}}{\sqrt{2}})$,
$\vec v_0=(-\frac{1}{\sqrt{2}},-\frac{\sqrt{3}}{\sqrt{2}})$.
\medskip

\begin{figure}[h]
\includegraphics[height=50mm]{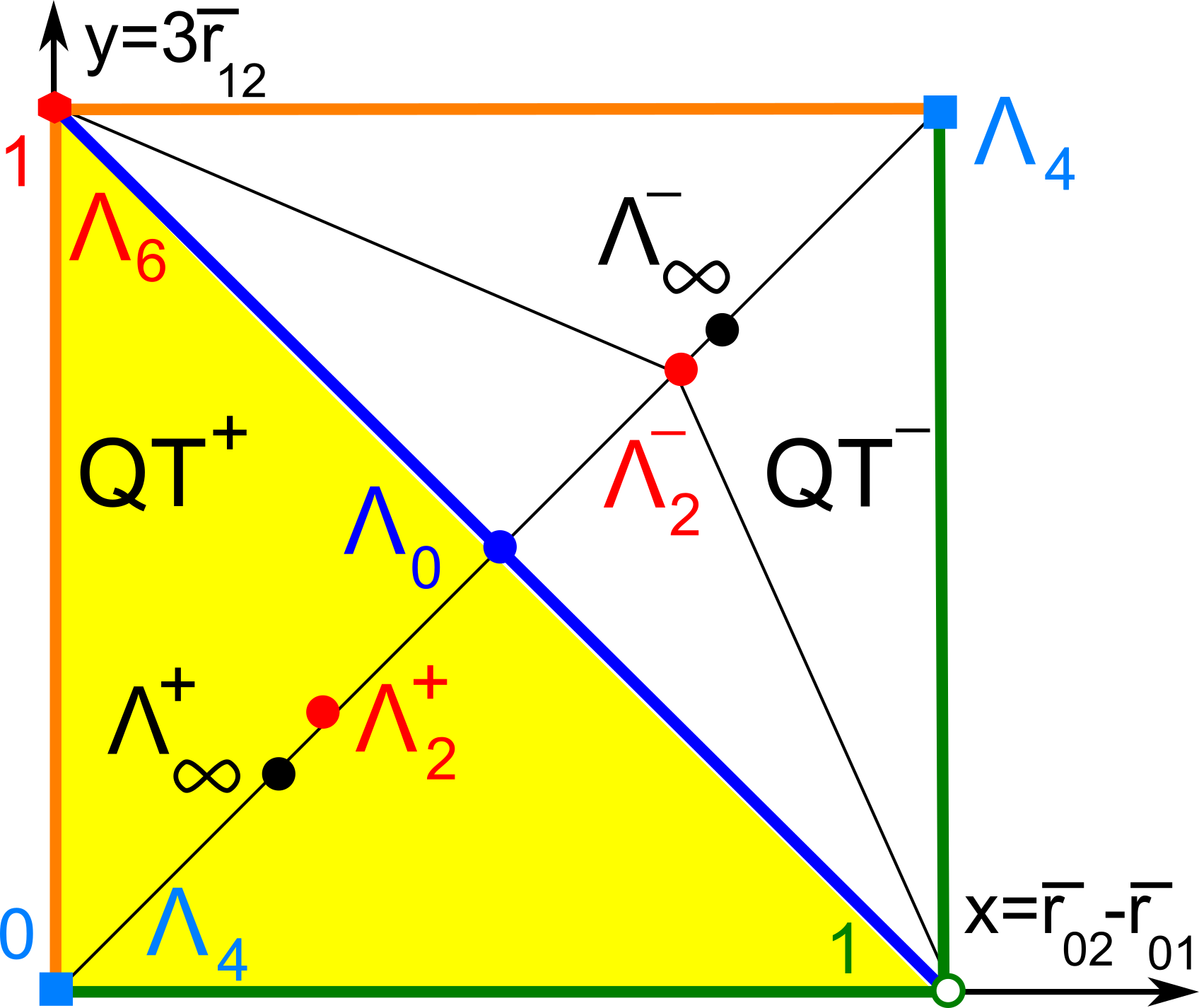}
\hspace*{2mm}
\includegraphics[height=50mm]{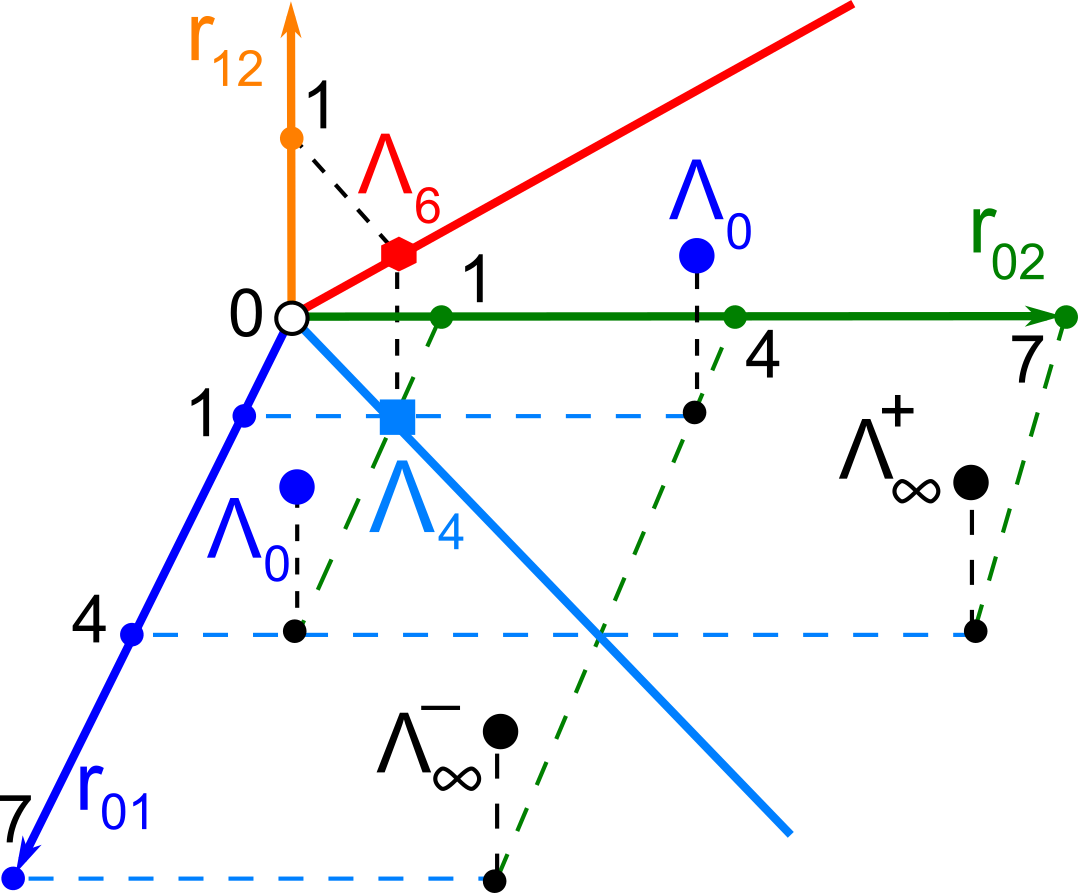}
\caption{\textbf{Left}: $\QS=\QT^+\cup\QT^-$ includes mirror-symmetric lattices $\La_4,\La_6,\La_0$ and non-mirror-symmetric lattices $\La_{\infty}^{\pm}$, see Example~\ref{exa:RM} and Table~\ref{tab:forms460inf} later.
\textbf{Right}: the doubled cone $\DC$ is visualised as $\{0\leq r_{12}\leq\min\{r_{01},r_{02}\}>0\}$ bounded by the planes $\{r_{12}=0\}$, $\{r_{12}=r_{01}\}$, $\{r_{12}=r_{02}\}$ with the identifications $(r_{12},r_{01},r_{02})\lra (r_{12},r_{02},r_{01})$ on the boundary $\bd\DC$.
}
\label{fig:QS+DC}
\end{figure}

\noindent
\emph{($\pmb{\La_0}$)}
We inversely design the lattice $\La_0$ in Fig.~\ref{fig:QS+DC} starting from $\PI(\La_0)=(x,y)$ at the centre $(\frac{1}{2},\frac{1}{2})\in\QS$.
Formula~(\ref{prop:inverse_design}a) for the size $|\La_0|=6$ (only to get simplest integers) gives $\RI(\La_0)=(1,1,4)$.
To reconstruct an obtuse superbase $\{\vec v_0,\vec v_1,\vec v_2\}$ by Lemma~\ref{lem:superbase_reconstruction}, find the vonorms $\vec v_1^2=1^2+1^2=2$, $\vec v_0^2=\vec v_2^2=1^2+4^2=17$.
Formula~(\ref{prop:inverse_design}b) gives the angle
$\angle(\vec v_1,\vec v_2)=\arccos
\dfrac{-4}{\sqrt{(\frac{9}{4}+\frac{5}{4}-3+9)^2-9(\frac{5}{4})^2}}
=\arccos(-\frac{1}{\sqrt{34}})\approx 99.9^\circ$.
We can choose the following superbase, see Fig.~\ref{fig:DT}:  
$\vec v_1=(\sqrt{2},0)$, 
$\vec v_2=|\vec v_2|(\cos\angle(\vec v_1,\vec v_2),\sin\angle(\vec v_1,\vec v_2))
=(-\frac{1}{\sqrt{2}},\frac{\sqrt{33}}{\sqrt{2}})$,
$\vec v_0=(-\frac{1}{\sqrt{2}},-\frac{\sqrt{33}}{\sqrt{2}})$.
\medskip

\noindent
\emph{($\pmb{\La_2}$)}
We inversely design the lattice $\La_2$ in Fig.~\ref{fig:QS+DC} 
starting from their projected invariants $\PI(\La_2)=(\frac{1}{2+\sqrt{2}},\frac{1}{2+\sqrt{2}})$, which will maximise the chiral distance $\PC[D_2]$ in Theorem~\ref{prop:PC}(a).
Formula~(\ref{prop:inverse_design}a) for the size $|\La_2|=6$ (only to simplify the root invariant) gives $\RI(\La_2)=(2-\sqrt{2},2\sqrt{2}-1,5-\sqrt{2})$.
Since all root products are non-zero and distinct, by Lemma~\ref{lem:achiral_lattices} there is a pair of lattices $\La_2^\pm$ with $\sign(\La_2^\pm)=\pm 1$.
The lattices $\La_2^\pm$ are related by reflection, not by rigid motion.
\medskip

To reconstruct an obtuse superbase $\{\vec v_0,\vec v_1,\vec v_2\}$ of $\La_2^\pm$ by Lemma~\ref{lem:superbase_reconstruction}, find \\
$\vec v_0^2=(2\sqrt{2}-1)^2+(5-\sqrt{2})^2=(9-4\sqrt{2})+(27-10\sqrt{2})=36-14\sqrt{2}\approx 16.2,$ \\
$\vec v_1^2=(2-\sqrt{2})^2+(2\sqrt{2}-1)^2=(6-4\sqrt{2})+(9-4\sqrt{2})=15-8\sqrt{2}\approx 3.7,$ \\
$\vec v_2^2=(2-\sqrt{2})^2+(5-\sqrt{2})^2=(6-4\sqrt{2})+(27-10\sqrt{2})=33-14\sqrt{2}\approx 13.2,$  \\
and the anticlockwise angle
$\angle(\vec v_1,\vec v_2)=\arccos
\dfrac{-r_{12}^2}{|\vec v_1|\cdot|\vec v_2|}\approx 92.8^\circ$.
Then $\La_2^\pm$ have the following obtuse superbases in
Fig.~\ref{fig:DT}:  
$\vec v_1=(\sqrt{15-8\sqrt{2}},0)\approx(1.9,0)$, 
$\vec v_2=|\vec v_2|(\cos\angle(\vec v_1,\vec v_2),\sin\angle(\vec v_1,\vec v_2))
\approx(-0.18,3.63)$,
$\vec v_0
\approx(-1.72,-3.63)$.
\medskip

\begin{figure}[h]
\includegraphics[width=\textwidth]{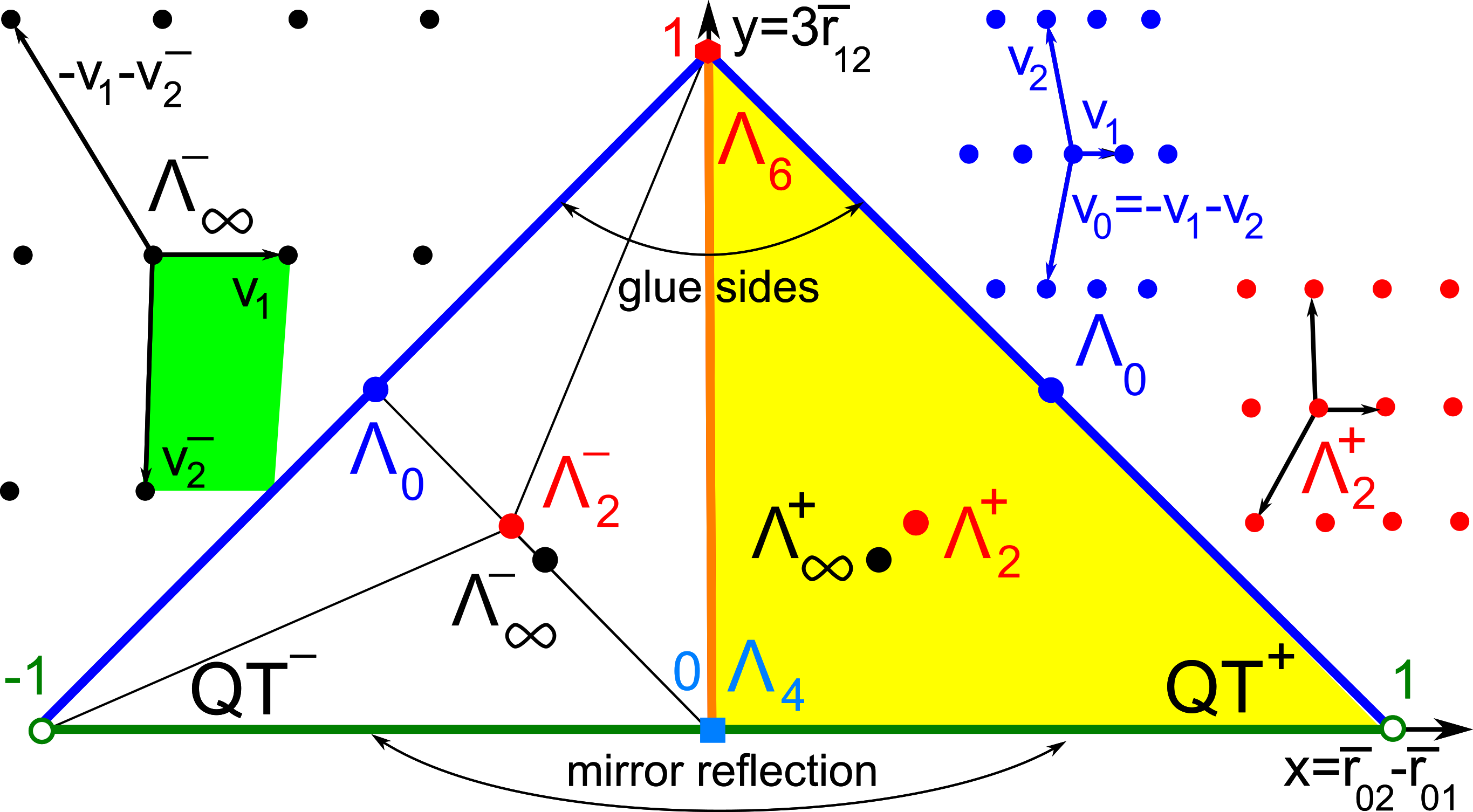}
\caption{
The doubled cone $\DC$ in Fig.~\ref{fig:QS+DC} (right) projects to the doubled triangle $\DT$ parameterised by $x\in(-1,1)$, $y\in[0,1]$ and obtained by gluing two copies $\QT^{\pm}$ of the quotient triangle along vertical sides instead of hypotenuses as in $\QS$, see Example~\ref{exa:inverse_design} and Table~\ref{tab:forms460inf}.
}
\label{fig:DT}
\end{figure}

\noindent
\emph{($\pmb{\La_{\infty}}$)}
We inversely design the lattice $\La_{\infty}$ in Fig.~\ref{fig:DT} starting from $\PI(\La_{\infty})=(x,y)$ at the mid-point $(\frac{1}{4},\frac{1}{4})$ of the segment between $\PI(\La_4),\PI(\La_0)\in\QT$.
Formula~(\ref{prop:inverse_design}a) for the size $|\La_{\infty}|=12$ (only to simplify the root invariant) gives 
 $\RI(\La_{\infty})=(1,4,7)$.
Since all root products are non-zero and distinct, by Lemma~\ref{lem:achiral_lattices} there is a pair of lattices $\La_{\infty}^\pm$ of opposite signs $\sign(\La_{\infty}^\pm)=\pm 1$.
\medskip

To reconstruct an obtuse superbase $\{\vec v_0,\vec v_1,\vec v_2\}$ of $\La_{\infty}^\pm$ by Lemma~\ref{lem:superbase_reconstruction}, find the vonorms $\vec v_0^2=4^2+7^2=65$, $\vec v_1^2=1^2+4^2=17$, $\vec v_2^2=1^2+7^2=50$, and the anticlockwise angle
$\angle(\vec v_1,\vec v_2)=\arccos
\dfrac{-r_{12}^2}{|\vec v_1|\cdot|\vec v_2|}
=\arccos(-\frac{1}{\sqrt{850}})\approx 92^\circ$.
Then $\La_{\infty}^\pm$ have the following obtuse superbases in
Fig.~\ref{fig:DT}:  
$\vec v_1=(\sqrt{17},0)\approx(4.12,0)$, 
$$\vec v_2^{\pm}=|\vec v_2|(\cos\angle(\vec v_1,\vec v_2),\sin\angle(\vec v_1,\vec v_2))
=\left(-\frac{1}{\sqrt{17}},\pm\frac{\sqrt{849}}{\sqrt{17}}\right)\approx(-0.24,\pm7.1),$$
$\vec v_0^\pm=-\vec v_1-\vec v_2^{\pm}=(-\frac{16}{\sqrt{17}},\mp\frac{\sqrt{849}}{\sqrt{17}})\approx(-3.88,\mp7.1)$, see all forms in
Table~\ref{tab:forms460inf}.
\eexa
\end{exa}

\begin{table}[h]
\caption{Various invariants of the lattices computed in Example~\ref{exa:inverse_design}, see Fig.~\ref{fig:QS+DC} and~\ref{fig:DT}. }
\medskip

\hspace*{-2mm}
\begin{tabular}{l|ccccc}      
$\La$ & 
$\La_4$ & 
$\La_6$ & 
$\La_0$ & 
$\La_2^{\pm}$ & 
$\La_{\infty}^{\pm}$ \\
\hline

$|\La)$ & 2 & 3 & 6 & 6 & 12 \\

$\PI(\La)$ 
& $(0,0)$
& $(0,1)$
& $\left(\dfrac{1}{2},\dfrac{1}{2}\right)$
& $\left(\dfrac{1}{2+\sqrt{2}},\dfrac{1}{2+\sqrt{2}}\right)$
& $\left(\dfrac{1}{4},\dfrac{1}{4}\right)$ \\


$\RI^o(\La)$  & 
(0,1,1) & 
(1,1,1) & 
(1,1,4) & 
$(2-\sqrt{2},2\sqrt{2}-1,5-\sqrt{2})^{\pm}$ & 
$(1,4,7)^{\pm}$ \\

$\VF(\La)$  & 
(2,1,1) & 
(2,2,2) & (17,2,17) & 
$(15-8\sqrt{2},33-14\sqrt{2},36-14\sqrt{2})$ & 
(65,17,50)
\end{tabular}
\label{tab:forms460inf}
\end{table}

Since the quotient square $\QS=\QT^+\cup\QT^-$ with identified sides is a {punctured sphere}, it is natural to visualise $\QS$ as the round surface of Earth with $\QT^\pm$ as the north/south hemispheres separated by the equator along their common boundary of $\QT$ represented by projected invariants $\PI(\La)$ of all mirror-symmetric lattices $\La$.
\smallskip

We can choose any internal point of the quotient triangle $\QT$ as the north pole.
The most natural choice is the incentre $P^+$ (pole), the centre of the circle inscribed into $\QT^+$ because the rays from $P^+$ to the vertices of $\QT^+$ equally bisect the angles $90^\circ, 45^\circ,45^\circ$. 
The incentre of $\QT^+$ has the coordinates $(x,x)$, where $x=1-\frac{1}{\sqrt{2}}=\frac{1}{2+\sqrt{2}}$.
The lattice $\La_2^+$ with the projected invariant $\PI(\La_2^+)=(x,x)$ has the basis $v_1\approx(1.9,0)$, $v_2\approx(-0.18,3.63)$ inversely designed in Example~\ref{exa:inverse_design}~($\La_2$).

\begin{figure}
\label{fig:QT+HS}
\caption{
\textbf{Left}: in $\QT^+$, the Greenwich line goes from the `empty' point (1,0) through incentre $P^+$ to the point $G=(0,\sqrt{2}-1)$.
\textbf{Middle}: the hemisphere $\HS^+$ has the north pole at $P^+$, the equator $\bd\QT^+$ of mirror-symmetric lattices.
\textbf{Right}: the longitude $\mu\in(-180^\circ,+180^\circ]$ anticlockwise measures angles from the Greenwich line, the latitude $\varphi\in[-90^\circ,+90^\circ]$ measures angles from the equator to the north.}
\includegraphics[height=36mm]{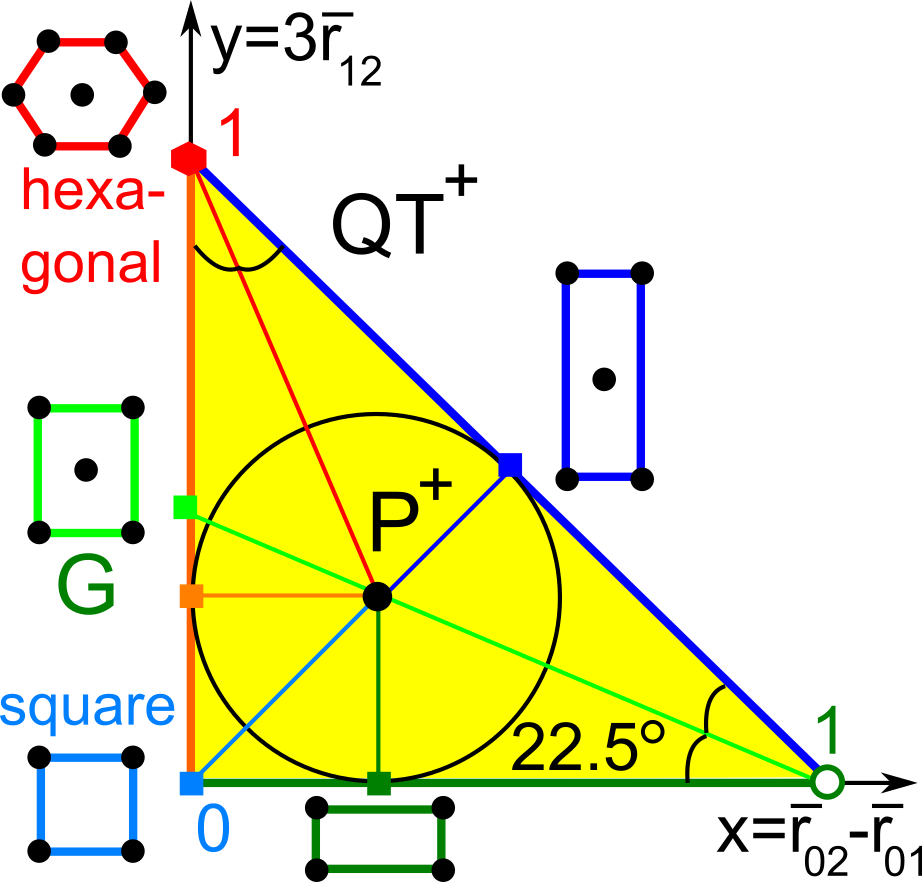}
\hspace*{1mm}
\includegraphics[height=36mm]{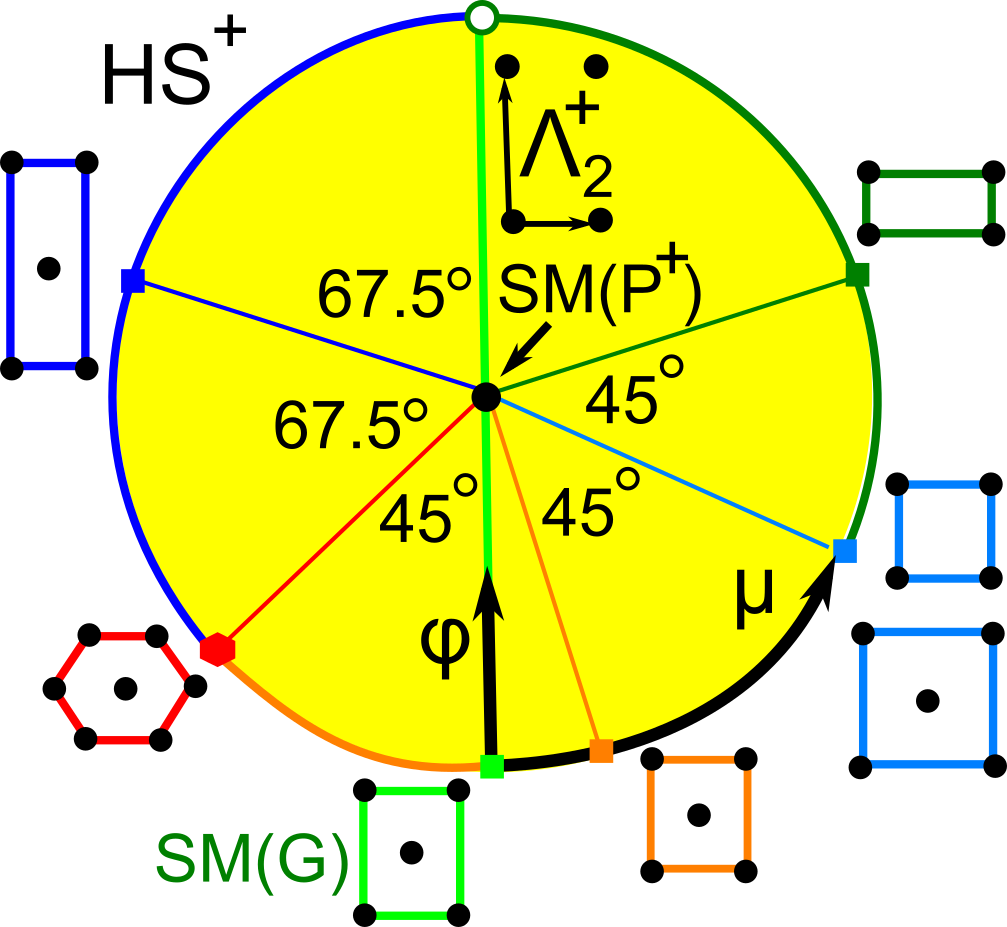}
\hspace*{1mm}
\includegraphics[height=36mm]{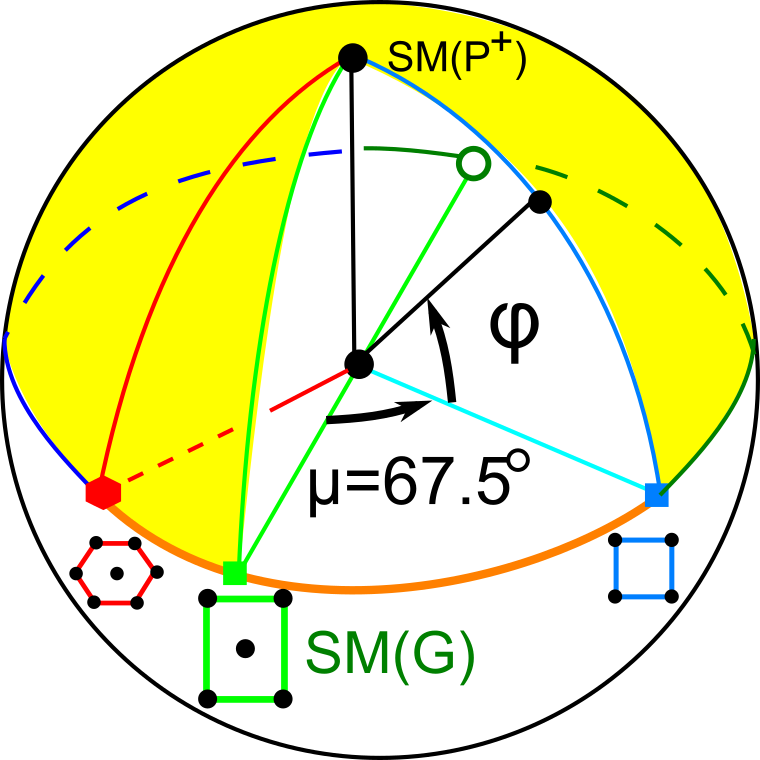}
\end{figure}

\begin{dfn}[spherical lattice map $\SLM:\QS\to S^2$] 
\label{dfn:SLM}
\textbf{(a)}
The \emph{spherical map} $\SLM$ sends the incentre $P^+$ of $\QT$ to the north pole of the hemisphere $\HS^+$ and the boundary $\bd\QT$ to the equator of $\HS^+$, see Fig.~\ref{fig:QT+HS}~(middle).
Linearly map the line segment between $P^+$ and
any point $(x,y)$ in the boundary $\bd\QT$ to the shortest arc connecting the north pole $\SLM(P^+)$ to $\SLM(x,y)$ in the equator of $\HS^+$.
Extend the \emph{spherical map} to $\SLM:\QS\to S^2$ by sending any pair of invariants $\PI^o(\La^{\pm})$ with $\sign(\La^{\pm})=\pm 1$ to the northern/southern hemispheres of the 2-dimensional sphere $S^2$, respectively.
\medskip

\noindent
\textbf{(b)}
For any lattice $\La\subset\R^2$, the \emph{latitude} $\ph(\La)\in[-90^\circ,+90^\circ]$ is the angle from the equatorial plane $\EP$ of $S^2$ to the radius-vector to the point $\SLM(\PI^o(\La))\in S^2$ in the upwards direction.
Let $v(\La)$ be the orthogonal projection of this radius-vector to $\EP$.
{Define the \emph{Greenwich} point as $G=(0,\sqrt{2}-1)\in\bd\QT$ in the line through $P^+$ and $(1,0)$.
This $G$ represents all centred rectangular lattices with a conventional unit cell $2a\times 2b$ whose ratio $r=\frac{b}{a}$ can be found from} Example~\ref{exa:achiral_lattices}(b): {$\sqrt{2}-1=\frac{3\sqrt{b^2-a^2}}{2a\sqrt{2}+\sqrt{b^2-a^2}}$.
Setting $s=\sqrt{r^2-1}$, we get $\sqrt{2}-1=\frac{3s}{2\sqrt{2}+s}$, 
$s=\frac{4-2\sqrt{2}}{4-\sqrt{2}}$, $r=\sqrt{s^2+1}\approx 1.1$. 
The \emph{Greenwich meridian} is the great circle on $S^2$ passing through the point $\SLM(G)$ in the equator $E$.} 
The longitude $\mu(\La)\in(-180^\circ,180^\circ]$ is the anticlockwise angle from the \emph{Greenwich plane} through the Greenwich meridian to the vector $v(\La)$ above.
\edfn
\end{dfn}

For lattices with $\PI(\La)$ in the straight-line segment between the excluded vertex $(1,0)$ and the incentre $P^+$, we choose the longitude $\mu=+180^\circ$ rather than $-180^\circ$.
Proposition~\ref{prop:SLM} computes 
the longitude  and latitude coordinates 
$\mu(\La),\ph(\La)$ via $\PI(\La)=(x,y)$
in terms of the projected invariant $\PI(\La)=(x,y)$ .

\begin{prop}[formulae for spherical lattice map $\SLM$, {\cite[Proposition 5.2]{bright2023geographic}}]
\label{prop:SLM}
For any lattice $\La\subset\R^2$ with 
$\PI(\La)=(x,y)\in\QT$,
 if $x\neq t=1-\dfrac{1}{\sqrt{2}}$, then set $\psi=\arctan\dfrac{y-t}{x-t}$, otherwise $\psi=\sign(y-t)90^\circ$.
\myskip

\noindent
\tb{(a)}
The longitude of the lattice $\La$ is 
$\mu(\La)=\left\{\begin{array}{l} 
\psi+22.5^\circ \text{ if } x<t, \\
\psi-157.5^\circ \text{ if } x\geq t, \psi\geq-22.5^\circ, \\
\psi+202.5^\circ \text{ if } x\geq t, \psi\leq-22.5^\circ.
\end{array} \right.$

\noindent
\tb{(b)}
The latitude is 
$\ph(\La)=\sign(\La)\cdot\left\{\begin{array}{l} 
\frac{x\sqrt{2}}{\sqrt{2}-1}90^\circ \text{ if } \mu(\La)\in[-45^\circ,+67.5^\circ], \\
\frac{y\sqrt{2}}{\sqrt{2}-1}90^\circ \text{ if }\mu(\La)\in[+67.5^\circ,+180^\circ], \\
\frac{1-x-y}{\sqrt{2}-1}90^\circ \text{ if } \mu(\La)\in[-180^\circ,-45^\circ].
\end{array} \right.$

\noindent
The incentres $P^{\pm}\in\QT^{\pm}$ have $\psi=0$ and {$\ph=\pm 90^{\circ}$, respectively, $\mu$ is undefined}.
\ethm
\end{prop}

\begin{exa}[prominent lattices]
\label{exa:SLM}
Any mirror-symmetric lattice $\La\subset\R^2$ has $\sign(\La)=0$, hence belongs to the equator $E$ of $S^2$ and has $\ph(\La)=0$ by (\ref{prop:SLM}b).
Any square lattice $\La_4$ with $\PI(\La_4)=(0,0)$ has $\mu(\La_4)=\arctan 1+22.5^\circ=67.5^\circ$ by (\ref{prop:SLM}a).
Any hexagonal lattice $\La_6$ with $\PI(\La_4)=(0,1)$ has $\mu(\La_4)=\arctan\frac{1}{1-\sqrt{2}}+22.5^\circ=-45^\circ$. 
Any rectangular lattice $\La$ with $\PI(\La)=(1-\frac{1}{\sqrt{2}},0)$ has $\mu(\La)=-90^\circ+202.5^\circ =112.5^\circ$.
Any centered rectangular lattice $\La$ with $\PI(\La)=(\frac{1}{2},\frac{1}{2})$ at the mid-point of the diagonal of $\QT$ has $\mu(\La)=\arctan 1-157.5^\circ =-112.5^\circ$.
Any \emph{Greenwich} lattice $\La_G$ with $\PI(\La_G)=G=(0,\sqrt{2}-1)$ has $\mu(\La_G)=\arctan(1-\sqrt{2})+22.5^\circ=0$.
\eexa
\end{exa}

In addition to the original paper \cite{kurlin2024mathematics}, we add new Corollary~\ref{cor:lattices2D_spaces} to fulfill the realisability and Euclidean embeddability conditions in Problem~\ref{pro:lattices2D}(f,g).

\index{embedding}

\begin{cor}[Euclidean embeddings of lattice spaces]
\label{cor:lattices2D_spaces}
\tb{(a)}
For all lattices $\La\subset\R^2$ under isometry, the Root Invariant Space $\RIS(\R^2)=\{\RI(\La) \vl \text{ lattices }\La\subset\R^2\}$ consists of all ordered triples $0\leq r_{12}\leq r_{01}\leq r_{02}$, where the smallest root product $r_{12}$ can be zero.
Then $\RIS(\R^2)$ is embedded into $\R^3$ as the triangular cone $\TC$ in Definition~\ref{dfn:TC}.
\myskip

\nt
\tb{(b)}
For all lattices $\La\subset\R^2$ under homothety, the invariant space $\{\PI(\La) \vl \text{ lattices }\La\subset\R^2\}$ consists of all points $(x,y)$ in the quotient triangle $\QT\subset\R^2$ from Definition~\ref{dfn:PI}.
\myskip

\nt
\tb{(c)}
For all lattices $\La\subset\R^2$ under dilation, the invariant space $\{\PI^o(\La) \vl \text{ lattices }\La\subset\R^2\}$ can be embedded onto $S^2\setminus\{\text{one point}\}\subset\R^3$ by the map $\SLM$ from Definition~\ref{dfn:SLM}.
\myskip

\nt
\tb{(d)}
For all lattices $\La\subset\R^2$ under rigid motion, the space $\{\PI^o(\La) \vl \text{ lattices }\La\subset\R^2\}$ can be embedded onto $(\S^2\setminus\{\text{one point}\})\times(0,+\infty)\subset\R^4$.
\ethm
\end{cor}
\begin{proof}
Parts (a,b,c) follow Definitions~\ref{dfn:TC}, \ref{dfn:PI}, and~\ref{dfn:SLM}.
Part (d) is obtained by extending the embedding $S^2\subset\R^3$ from part~(c) by a scaling factor $s\in(0,+\infty)$
All embeddings are bi-Lipschitz because all invariant spaces can be compactified and all underlying maps are expressed via elementary functions as in \ref{prop:SLM}.
\end{proof}

\section{Metrics on spaces of 2D lattices under all four equivalences}
\label{sec:lattices2D_metrics}

All lattices $\La\subset\R^2$ are uniquely represented under isometry and homothety by their invariants $\RI\in\TC$ and $\PI\in\QT$, respectively.
Then any metric $d$ on the triangular cone $\TC\subset\R^3$ or the quotient triangle $\QT\subset\R^2$ gives rise to a metric in Definition~\ref{dfn:RM} on the spaces $\LIS$ and $\LHS$, respectively. 
The oriented case in Definition~\ref{dfn:RMo} will be harder because of identifications on the boundary $\bd\TC$.

\begin{dfn}[root metrics $\RM$, projected metrics $\PM$]
\label{dfn:RM}
Any metric $d$ on $\R^3$ defines the \emph{root metric} 
$\RM(\La_1,\La_2)=d(\RI(\La_1),\RI(\La_2))$ on lattices $\La_1,\La_2\subset\R^2$  under isometry.
The \emph{Root Invariant Space} $\RIS=(\TC,d)$ is the triangular cone with a fixed metric $d$.
If we use the Minkowski norm $M_q(v)=||v||_q=(\sum\limits_{i=1}^n |x_i|^q)^{1/q}$ of a vector $\vec v=(x_1,\dots,x_n)\in\R^n$  for any real $q\in[1,+\infty]$, the root metric is denoted by $\RM_q(\La_1,\La_2)=||\RI(\La_1)-\RI(\La_2)||_q$.
The limit case $q=+\infty$ uses 
$||v||_{\infty}=\max\limits_{i=1,\dots,n}|x_i|$.
The \emph{projected metric} $\PM(\La_1,\La_2)=d(\PI(\La_1),\PI(\La_2))$ is on lattices under homothety for any metric $d$ on $\R^2$.
The space of \emph{projected invariants} $\PIN=(\QT,d)$ is the quotient triangle with a metric $d$.
The notation $\PM_q(\La_1,\La_2)=||\PI(\La_1)-\PI(\La_2)||_q$ includes a parameter $q\in[1,+\infty]$ of $M_q$.
\edfn
\end{dfn}

The Minkowski distance $M_q$ for $q=2$ is Euclidean.
The root metric $\RM_q$ can take any large values in original units of vector coordinates such as Angstroms.
The projected metric $\PM_q$ is unitless and the space $\PIN=(\QT,d)$ is bounded.

\begin{table}[h]
\caption{Metrics $\RM_q$ and $\PM_q$ for the lattices from Example~\ref{exa:RM} and shown Fig.~\ref{fig:QS+DC} and~\ref{fig:DT}. }
\medskip

\begin{tabular}{|l|cccc|}  
\hline    
$\RM_{\infty}$ & $\La_4$ & $\La_6$ & $\La_0$ & $\La_{\infty}^{\pm}$\\
\hline

$\RI(\La_4)=(0,1,1)$ & $0$ & $1$ & $3$ & $6$ \\

$\RI(\La_6)=(1,1,1)$ & $1$ & $0$ & $3$ & $6$ \\

$\RI(\La_0)=(1,1,4)$ & $3$ & $3$ & $0$ & $3$  \\

$\RI(\La_{\infty}^{\pm})=(1,4,7)$ & $6$ & $6$ & $3$ & $0$  \\
\hline
\end{tabular}
\hspace*{4mm}
\begin{tabular}{|l|cccc|}
\hline      
$\PM_{\infty}$ & $\La_4$ & $\La_6$ & $\La_0$ & $\La_{\infty}^{\pm}$\\
\hline

$\PI(\La_4)=(0,0)$ & $0$ & $1$ & $\frac{1}{2}$ & $\frac{1}{4}$ \\

$\PI(\La_6)=(0,1)$ & $1$ & $0$ & $\frac{1}{2}$ & $\frac{3}{4}$ \\

$\PI(\La_0)=(\frac{1}{2},\frac{1}{2})$ & $\frac{1}{2}$ & $\frac{1}{2}$ & $0$ & $\frac{1}{4}$ \\

$\PI(\La_{\infty}^{\pm})=(\frac{1}{4},\frac{1}{4})$ & $\frac{1}{4}$ & $\frac{3}{4}$ & $\frac{1}{4}$ & $0$  \\
\hline
\end{tabular}
\medskip

\begin{tabular}{|l|cccc|}  
\hline    
$\RM_q$ for $q\in[1,+\infty)$ & $\La_4$ & $\La_6$ & $\La_0$ & $\La_{\infty}^{\pm}$\\
\hline

$\RI(\La_4)=(0,1,1)$ & $0$ & $1$ & $(1+3^q)^{1/q}$ & $(1+3^q+6^q)^{1/q}$ \\

$\RI(\La_6)=(1,1,1)$ & $1$ & $0$ & $3$ & $(3^q+6^q)^{1/q}$ \\

$\RI(\La_0)=(1,1,4)$ & $(1+3^q)^{1/q}$ & $3$ & $0$ & $3\cdot 2^{1/q}$  \\

$\RI(\La_{\infty}^{\pm})=(1,4,7)$ & $(1+3^q+6^q)^{1/q}$ & $(3^q+6^q)^{1/q}$ & $3\cdot 2^{1/q}$ & $0$  \\
\hline
\end{tabular}
\medskip

\begin{tabular}{|l|c|c|c|c|}
\hline      
$\PM_q$ for $q\in[1,+\infty)$ & $\La_4$ & $\La_6$ & $\La_0$ & $\La_{\infty}^{\pm}$\\
\hline

$\PI(\La_4)=(0,0)$ & $0$ & $1$ & $2^{(1/q)-1}$ & $2^{(1/q)-2}$ \\

$\PI(\La_6)=(0,1)$ & $1$ & $0$ & $2^{(1/q)-1}$ & $\frac{1}{4}(1+3^q)^{1/q}$ \\

$\PI(\La_0)=(\frac{1}{2},\frac{1}{2})$ & $2^{(1/q)-1}$ & $2^{(1/q)-1}$ & $0$ & $2^{(1/q)-2}$  \\

$\PI(\La_{\infty}^{\pm})=(\frac{1}{4},\frac{1}{4})$ & $2^{(1/q)-2}$ & $\frac{1}{4}(1+3^q)^{1/q}$ & $2^{(1/q)-2}$ & $0$  \\
\hline
\end{tabular}
\label{tab:RM}
\end{table}

\begin{exa}[metrics $\RM_q,\PM_q$]
\label{exa:RM}
Table~\ref{tab:RM} summarises metric computations for the lattices $\La_4,\La_6,\La_0,\La_{\infty}^{\pm}$, which  
were inversely designed in Example~\ref{exa:inverse_design}.
\eexa
\end{exa}

Since the mirror images $\La_{\infty}^{\pm}$ have the same root invariant $\RI(\La_{\infty}^{\pm})=(1,4,7)$, for any lattice $\La$, the distances $\RM(\La,\La_{\infty}^{\pm})$ and $\PM(\La,\La_{\infty}^{\pm})$ are independent of $\sign(\La_{\infty}^{\pm})=\pm 1$.
Any mirror images $\La^{\pm}$ have $\RM(\La^+,\La^-)=0=\PM(\La^+,\La^-)$ because $\La^{\pm}$ are isometric to each other.
The metric $\RM$ from Definition~\ref{dfn:RM} is well-defined only for lattices under any isometry including reflections.
\medskip

\begin{figure}[h]
\includegraphics[width=1.0\textwidth]{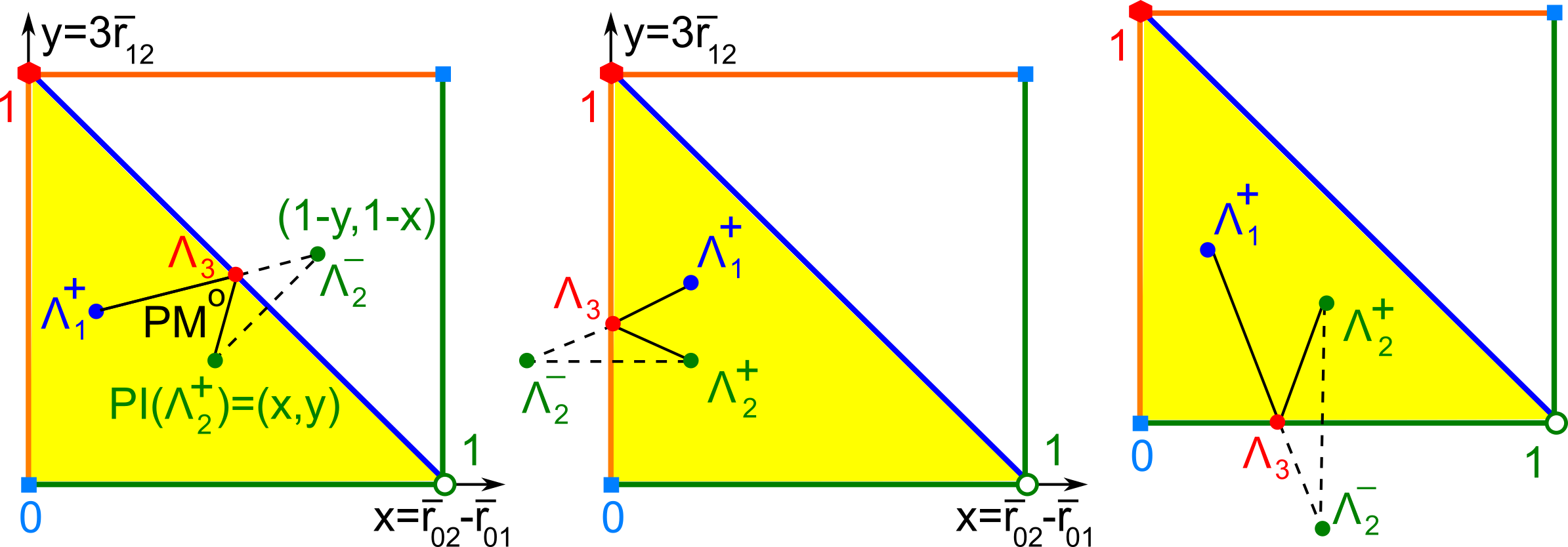}
\caption{By Definition~\ref{dfn:RMo}, the projected metric $\PM_2(\La_1^+,\La_2^-)$ is the minimum sum $\PM_2(\La_1^+,\La_3)+\PM_2(\La_3\La_2^-)$ achieved in the left image, see computations in Proposition~\ref{prop:RMq}.}
\label{fig:QS+RMo}
\end{figure}

Definition~\ref{dfn:RMo} introduces the metric $\RM^o$ on lattices under rigid motion so that $\RM^o(\La^+,\La^-)>0$ on mirror images of a non-mirror-symmetric lattice, see Fig.~\ref{fig:QS+RMo}.

\begin{dfn}[orientation-aware metrics $\RM^o,\PM^o$]
\label{dfn:RMo}
For lattices $\La_1,\La_2\subset\R^2$ with $\sign(\La_1)\sign(\La_2)\geq 0$, the \emph{orientation-aware} root metric is $\RM^o(\La_1,\La_2)=\RM(\La_1,\La_2)$ as in Definition~\ref{dfn:RM}.
If any lattices $\La_1,\La_2$ have opposite signs, set $\RM^o(\La_1,\La_2)=\inf\limits_{\sign(\La_3)=0}(\RM(\La_1,\La_3)+\RM(\La_2,\La_3))$.
The \emph{orientation-based} metric $\PM^o(\La_1,\La_2)$ is defined by the same formula, where we replace $\RM$ by $\PM$.
\edfn
\end{dfn}

The infimum in $\RM^o(\La_1,\La_2)$ is the greatest lower bound 
defining a metric on a union of metric spaces glued by isometries.
Theoretically, this bound may not be achieved over a non-compact domain.
When using a Minkowski base metric $M_q$, Propositions~\ref{prop:RMq}-\ref{prop:PMq} explicitly compute $\RM_q^o,\PM_q^o$ for $q=2,+\infty$, so the infimum in Definition~\ref{dfn:RMo} can be replaced by a minimum in practice.
\medskip

The \emph{oriented} root invariant space $\RIS^o$ and the space of \emph{oriented} projected invariants $\PIN^o$ can be defined similarly to $\RIS$ and $\PIN$ in Definition~\ref{dfn:RM} as the doubled cone $\DC$ and quotient square $\QS$ with any metrics from Definition~\ref{dfn:RMo}.
\cite[Lemmas 5.3 and 5.5]{kurlin2024mathematics} prove the metric axioms for $\RM,\PM$ and $\RM^o,\PM^o$, respectively. 
Lemma~\ref{lem:reversed_signs} speeds up computations in the oriented case, see Example~\ref{exa:RMo}.

\begin{lem}[reversed signs, {\cite[Lemma 5.6]{kurlin2024mathematics}}]
\label{lem:reversed_signs}
If lattices $\La_1^{\pm},\La_2^{\pm}\subset\R^2$ have specified signs, then 
$\RM^o(\La_1^+,\La_2^-)=\RM^o(\La_1^-,\La_2^+)$ and
$\PM^o(\La_1^+,\La_2^-)=\PM^o(\La_1^-,\La_2^+)$.
\elem
\end{lem}

If lattices $\La_1,\La_2$ have non-opposite signs, so $\sign(\La_1)\sign(\La_2)\geq 0$, then the metrics $\RM_q^o$ and $\PM_q^o$ from Definition~\ref{dfn:RMo} coincide with the easily computable unoriented metrics $\RM_q,\PM_q$ from Definition~\ref{dfn:RM}.
Hence Propositions~\ref{prop:RMq} and~\ref{prop:PMq} compute $\RM_q^o(\La_1,\La_2)$ and $\PM_q^o(\La_1,\La_2)$ only for lattices of opposite signs.
 
\begin{prop}[root metrics for $q=2,+\infty$, {\cite[Proposition 5.8]{kurlin2024mathematics}}]
\label{prop:RMq}
Let $\La_1,\La_2\subset\R^2$ be lattices of opposite signs with $\RI(\La_1)=(r_{12},r_{01},r_{02})$, 
$\RI(\La_2)=(s_{12},s_{01},s_{02})$.
\myskip

\noindent
\textbf{(a)}
$\RM_2^o(\La_1,\La_2)$ is the minimum of the Euclidean distances from the point $\RI(\La_1)$ to the three points $(-s_{12},s_{01},s_{02})$, $(s_{01},s_{12},s_{02})$, and $(s_{12},s_{02},s_{01})$ in $\R^3$.
\myskip

\noindent
\textbf{(b)}
$\RM_{\infty}^o(\La_1,\La_2)=\min\{ d_0, d_1, d_2 \}$, where \\  
$d_0=\max\{r_{12}+s_{12}, |r_{01}-s_{01}|, |r_{02}-s_{02}|\},$\\
$d_1=\max\{\MS(r_{12},r_{01},s_{12},s_{01}),|r_{02}-s_{02}|\}$,\\
$d_2=\max\{|r_{12}-s_{12}|,\MS(r_{01},r_{02},s_{01},s_{02})\}$, \\
where $\MS(a,b,c,d)=\max\{|a-b|,|c-d|,\frac{1}{2}|a+b-c-d|\}$.
\ethm
\end{prop}

\begin{prop}[projected metrics for $q=2,+\infty$, {\cite[Proposition 5.9]{kurlin2024mathematics}}]
\label{prop:PMq}
Let $\La_1,\La_2$ be lattices with opposite signs and invariants $\PI(\La_1)=(x_1,y_1)$, $\PI(\La_2)=(x_2,y_2)$.
\myskip

\noindent
\textbf{(a)}
$\PM_2^o(\La_1,\La_2)$ is the minimum of the Euclidean distances from $\PI(\La_1)=(x_1,y_1)$ to the three points $(-x_2,y_2),(x_2,-y_2),(1-y_2,1-x_2)$ in $\R^2$.  
\myskip

\noindent
\textbf{(b)}
For $x_1\leq x_2$, $\PM_{\infty}^o(\La_1,\La_2)=\min\{ d_x, d_y, d_{xy} \}$ for $d_x=\max\{x_2-x_1,y_2+y_1\}$, $d_y=\max\{x_2+x_1,|y_2-y_1|\}$, $d_{xy}=\max\{x_2-x_1,1-x_2-y_2+|1-y_1-x_2|\}$.
\bt
\end{prop}

\section{Real-valued chiral distances measure asymmetry of lattices}
\label{sec:chiral_distances}

The classical concept of chirality is a binary property distinguishing mirror images of the same object such as a molecule or a periodic crystal.
Continuous classifications in Theorem~\ref{thm:lattices2D_classification} and Corollary~\ref{cor:lattices2D_classification} imply that the binary chirality is discontinuous under almost any perturbations similar to other discrete invariants such as symmetry groups.
To avoid arbitrary thresholds, it makes more sense to continuously quantify a deviation of a lattice from a higher-symmetry neighbour.  
\medskip

The term \emph{chirality} often refers to 3-dimensional molecules or crystal lattices.
One reason is the fact that in $\R^2$ a reflection with respect to a line $L$ is realised by the rotation in $\R^3$ around $L$ through $180^\circ$. 
However, if our ambient space is only $\R^2$, the concepts of isometry and rigid motion differ.
For example, Lemma~\ref{lem:achiral_lattices} described root invariants of all lattices that are related to their mirror images by rigid motion.
Such lattices can be called \emph{achiral}. 
We call them mirror-symmetric to avoid a potential confusion with the literature in crystallography.
Definition~\ref{dfn:RC} introduces the real-valued $G$-chiral distances of a lattice $\La\subset\R^2$.
Proposition~\ref{prop:cont_chiral_distances} proves the continuity of these functions $\RC[G]:\LIS(\R^2)\to\R$ and $\PC[G]:\LHS(\R^2)\to\R$.
\medskip

Recall that the \emph{crystallographic point group} $G$ of a lattice $\La\subset\R^2$ containing the origin $0$ consists of all symmetry operations that keep $0$ and map $\La$ to itself.
For example, any such group $G$ includes the central symmetry with respect to $0\in\La\subset\R^2$.
If $G$ has no other non-trivial symmetries, we get $G=C_2$ in Schonflies notations.
All 2D lattices split into four crystal families by their point groups: oblique ($C_2$), orthorhombic ($D_2$), tetragonal or square ($D_4$) and hexagonal ($D_6$).
Orthorhombic lattices split into rectangular and centred rectangular, see Fig.~\ref{fig:QT+QS}.

\begin{figure}[h]
\includegraphics[width=1.0\textwidth]{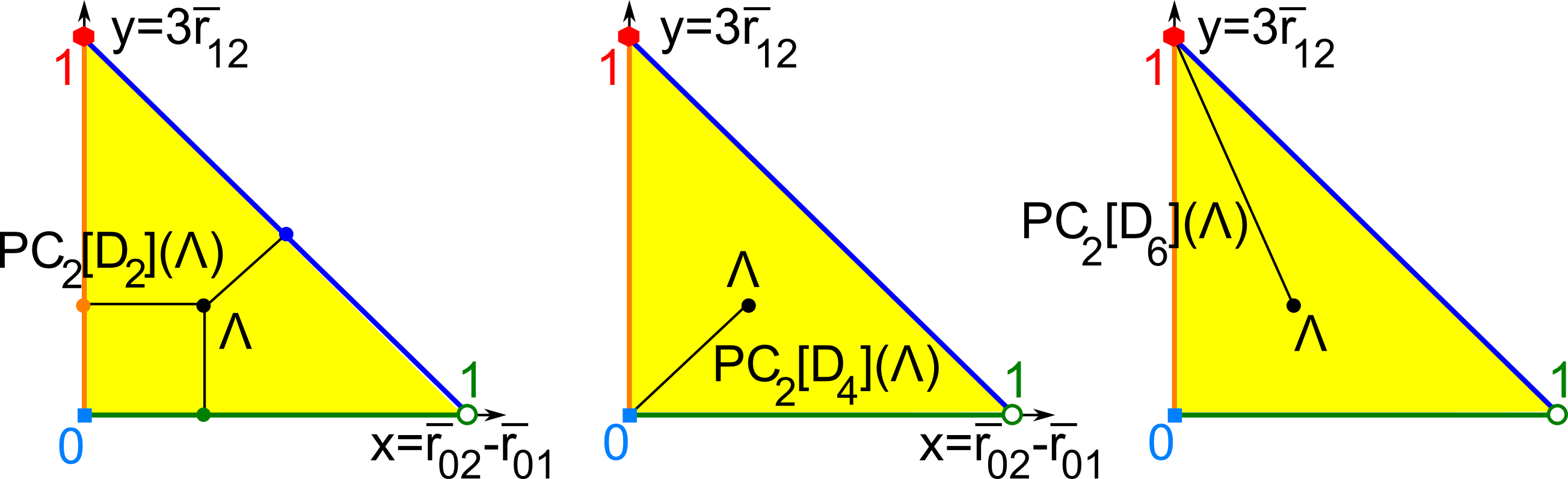}
\caption{\textbf{Left}: by Definition~\ref{dfn:RC}, the projected $D_2$ chiral distance $\PC_2[D_2](\La)$ is the minimum Euclidean distance from $\PI(\La)\in\QT$ to the boundary $\bd\QT$.
\textbf{Middle}: $\PC_2[D_4](\La)$ is the distance from $\PI(\La)$ to $(0,0)$.
\textbf{Right}: $\PC_2[D_6](\La)$ is the distance from $\PI(\La)$ to $(0,1)$.
}
\label{fig:QT+PC}
\end{figure}
 
\begin{dfn}[$G$-chiral distances {$\RC[G]$ and $\PC[G]$}]
\label{dfn:RC}
For any crystallographic point group $G$ in $\R^2$, let $\LIS[G]\subset\LIS(\R^2)$ be the closure of the subspace of all (isometry classes of) lattices that have the crystallographic  point group $G$.
For $G=D_2$ or $G=D_4$ or $G=D_6$, the root and projected \emph{$G$-chiral} distances are 
$$\RC[G](\La)=\min\limits_{\La'\in\LIS[G]}\RM(\La,\La')\geq 0 \text{ and }
\PC[G](\La)=\min\limits_{\La'\in\LIS[G]}\PM(\La,\La')\geq 0,$$
where $\RM$, $\PM$ are any metrics from Definition~\ref{dfn:RM} with a base metric $d$.
If $d=M_q$ for $q\in[1,+\infty]$, denote the $G$-chiral distances by $\RC_q[G]$ and $\PC_q[G]$.
\edfn
\end{dfn}

Since any lattice $\La$ is symmetric with respect to the origin $0\in\La$,
the closed subspace $\LIS[C_2]$ coincides with the 3-dimensional Lattice Isometry Space $\LIS(\R^2)$.
The 2-dimensional subspace $\LIS[D_2]$ consists of all mirror-symmetric lattices (rectangular and centred-rectangular) represented by root invariants $\RI$ on the boundary $\bd\TC$ of the triangular cone in Definition~\ref{dfn:TC}, see  Fig.~\ref{fig:TC}.
The 1-dimensional subspaces $\LIS[D_4],\LIS[D_6]\subset\LIS[D_2]$ can be viewed as the blue and orange rays $\{r_{12}=0<r_{01}=r_{02}\}$ and  $\{0<r_{12}=r_{01}=r_{02}\}$, respectively.   
\medskip

The $G$-chiral distance $\RC[G]$ in Definition~\ref{dfn:RC} measures a distance from $\RI(\La)$ to the root invariant of a closest neighbour in the subspace $\LIS[G]$.
Any $\RC[G](\La)$ is invariant under isometry and measures a distance from $\La$ to its nearest neighbour $\La'\in\LIS[G]$.
The \emph{signed chiral distances} $\sign(\La)\RC(\La)$ and $\sign(\La)\PC(\La)$ are invariant under rigid motion.
Since $\LIS[G]$ is a closed subspace within $\LIS(\R^2)$, the continuous distances $\RM,\PM$ achieve their minima if their base distances $d$ are continuous.
If $\LIS[D_2]$ was defined as an open subspace of only lattices that have the point group $D_2$ (not $D_4$ or $D_6$), then $\RC[G],\PC[G]$ should be defined via infima instead of simpler minima.
Indeed, any square or hexagonal lattice $\La$ can be approximated by infinitely many closer and closer orthorhombic lattices $\La'$, but the expected distance $\RM(\La,\La')=0$ will not be achieved on an open set.
\medskip

For 
$q=2,+\infty$, the distances $\RC_q,\PC_q$ are computed in Propositions~\ref{prop:RC},~\ref{prop:PC}.

\begin{lem}[properties of chiral distances, {\cite[Lemma 6.2]{kurlin2024mathematics}}]
\label{lem:RC}
\textbf{(a)}
A lattice $\La\subset\R^2$ is mirror-symmetric if and only if $\RC[D_2](\La)=0$ or, equivalently, $\PC[D_2](\La)=0$.
\medskip

\noindent
\textbf{(b)}
For any crystallographic point group $G$ in $\R^2$, mirror reflections $\La^{\pm}\subset\R^2$ have equal $G$-chiral distances: $\RC[G](\La^+)=\RC[G](\La^-)$,
$\PC[G](\La^+)=\PC[G](\La^-)$.
\elem
\end{lem}

\begin{lem}[lower bounds, {\cite[Lemma 6.3]{kurlin2024mathematics}}]
\label{lem:RCbounds}
\textbf{(a)}
If lattices $\La_1,\La_2$ have opposite signs, then
$\RM^o(\La_1,\La_2)\geq
\RC[D_2](\La_1)+\RC[D_2](\La_2)$ and \\
$\PM^o(\La_1,\La_2)\geq
\PC[D_2](\La_1)+\PC[D_2](\La_2)$.
\medskip

\noindent
\textbf{(b)}
For the mirror images $\La^{\pm}$ of any lattice $\La$, the lower bounds in part (a) become equalities:
$\RM^o(\La^+,\La^-)=2\RC[D_2](\La)$ and
$\PM^o(\La^+,\La^-)=2\PC[D_2](\La)$.
\elem
\end{lem}

\begin{prop}[chiral distances {$\RC_q[G]$} for $q=2,+\infty$, {\cite[Proposition 6.5]{kurlin2024mathematics}}]
\label{prop:RC}
Let a lattice $\La\subset\R^2$ have a root invariant $\RI(\La)=(r_{12},r_{01},r_{02})$ with $0\leq r_{12}\leq r_{01}\leq r_{02}$.
\medskip

\hspace*{-7mm}
$\begin{array}{llr}
\tb{(a)} &
\RC_2[D_2](\La)=\min\left\{ r_{12}, \dfrac{r_{01}-r_{12}}{\sqrt{2}}, \dfrac{r_{02}-r_{01}}{\sqrt{2}}\right\}; & \\
\smallskip

&
\RC_2[D_4](\La)=\sqrt{r_{12}^2+\frac{1}{4}(r_{02}-r_{01})^2}; & \\
\smallskip

&
\RC_2[D_6](\La)=\sqrt{\frac{2}{3}(r_{12}^2+r_{01}^2+r_{02}^2-r_{12}r_{01}-r_{12}r_{02}-r_{01}r_{02})}; & \\

\tb{(b)} &
\RC_{\infty}[D_2](\La)=\min\left\{ r_{12}, \dfrac{r_{01}-r_{12}}{2}, \dfrac{r_{02}-r_{01}}{2}\right\}. & \\

&
\RC_{\infty}[D_4](\La)=\min\{r_{12},\dfrac{r_{02}-r_{01}}{2}\}; & \\
\smallskip

&
\RC_{\infty}[D_6](\La)=\dfrac{r_{02}-r_{12}}{2}. & \blacktriangle
\end{array}$
\ethm
\end{prop}

When considering lattices under homothety, the subspace $\LHS[D_4]$ consists of a single class of all square lattices, which are all equivalent under isometry and uniform scaling.
The subspace $\LHS[D_6]$ is also a single point representing all hexagonal lattices.
Then $\PC[D_4]$ and $\PC[D_6]$ are distances to these single points. 

\begin{prop}[chiral distances $\PC_q$ for $q=2,+\infty$, {\cite[Proposition 6.6]{kurlin2024mathematics}}]
\label{prop:PC}
Let a lattice $\La$ have $\PI(\La)=(x,y)\in\QT$ so that $x\in[0,1)$, $y\in[0,1]$, $x+y\leq 1$.
\medskip

\hspace*{-7mm}
$\begin{array}{llr}
\tb{(a)} &
\PC_2[D_2](\La)=\min\left\{x,y,\dfrac{1-x-y}{\sqrt{2}}\right\}, & \\
\smallskip

&
\PC_q[D_4](\La)=(x^q+y^q)^{1/q} \text{ for any } q\in[1,+\infty), & \\
\smallskip

&
\PC_q[D_6](\La)=(x^q+(1-y)^q)^{1/q} \text{ for any } q\in[1,+\infty); & \\
\smallskip

\tb{(b)} &
\PC_{\infty}[D_2](\La)=\min\left\{x,y,\dfrac{1-x-y}{2}\right\}, & \\
\smallskip

&
\PC_{\infty}[D_4](\La)=x, & \\
\smallskip

&
\PC_{\infty}[D_6](\La)=1-y. & 
\end{array}$
\medskip

\noindent
\tb{(c)}
The upper bounds 
$\PC_{2}[D_2](\La)\leq\frac{1}{2+\sqrt{2}}$,
$\PC_{\infty}[D_2](\La)\leq\frac{1}{4}$ hold for any $\La$, achieved for lattices with $\PI(\La_2)=(\frac{1}{2+\sqrt{2}},\frac{1}{2+\sqrt{2}})$, $\PI(\La_{\infty})=(\frac{1}{4},\frac{1}{4})$, respectively.
For $q\in[1,+\infty]$, the bound $\PC_{q}[D_4](\La)\leq 1$ holds for any $\La$ and is achieved for any hexagonal lattice.
For $q\in[1,+\infty)$, the upper bound $\PC_{q}[D_6](\La)<2^{1/q}$ holds for any $\La$ and is approached but not achieved as $x\to 1$.
The bound $\PC_{\infty}[D_6](\La)\leq 1$ holds for any $\La$ and is achieved for any square and rectangular lattice. 
\ethm
\end{prop}

\begin{exa}[distances $\RC_q,\PC_q$]
\label{exa:RC}
Table~\ref{tab:RC} shows the chiral distances computed by Propositions~\ref{prop:RC},~\ref{prop:PC} for the prominent lattices $\La_2^{\pm}$, $\La_{\infty}^{\pm}$ in
Example~\ref{exa:inverse_design}.
\eexa
\end{exa}

\begin{table}[h]
\caption{Chiral distances $\PC_q,\RC_q$ for the lattices $\La_2^{\pm},\La_{\infty}^{\pm}$ in Fig.~\ref{fig:QS+DC} and~\ref{fig:DT}, see Example~\ref{exa:RC}. }
\medskip

$\begin{array}{|l|cc|}  
\hline    
\La & \La_{\infty} & \La_2 \\

\PI(\La) & 
\left(\dfrac{1}{4},\dfrac{1}{4}\right) & 
\left(\dfrac{1}{2+\sqrt{2}},\dfrac{1}{2+\sqrt{2}}\right) \\
\hline
 
\PC_2[D_2] & 
\dfrac{1}{4} & 
\dfrac{1}{2+\sqrt{2}} \\

\PC_2[D_4]  & 
\dfrac{\sqrt{2}}{4} & 
\sqrt{2}-1 \\

\PC_2[D_6] & 
\dfrac{\sqrt{10}}{4} & 
\sqrt{2-\sqrt{2}} \\
\hline

\PC_{\infty}[D_2]  & 
\dfrac{1}{4} &
\dfrac{1}{2+\sqrt{2}} \\

\PC_{\infty}[D_4] & 
\dfrac{1}{4} &
\dfrac{1}{2+\sqrt{2}} \\

\PC_{\infty}[D_6]  & 
\dfrac{3}{4} &
\dfrac{1}{\sqrt{2}} \\
\hline
\end{array}$
\hspace*{4mm}
\begin{tabular}{|l|cc|}  
\hline    
$\La$ & $\La_{\infty}$ & $\La_2$ \\

$\RI(\La)$ & 
$(1,4,7)$ & 
$(2-\sqrt{2},2\sqrt{2}-1,5-\sqrt{2})$ \\
\hline
 
$\RC_2[D_2]$ & 
$1$ & 
$2-\sqrt{2}$ \\

$\RC_2[D_4]$  & 
$\dfrac{\sqrt{13}}{2}$ & 
$(2-\sqrt{2})\dfrac{\sqrt{13}}{2}$ \\

$\RC_2[D_6]$ & 
$3\sqrt{2}$ & 
$\sqrt{2(13-3\sqrt{2})}$ \\
\hline

$\RC_{\infty}[D_2]$  & 
$1$ & 
$2-\sqrt{2}$ \\

$\RC_{\infty}[D_4]$ & 
$1$ & 
$2-\sqrt{2}$ \\

$\RC_{\infty}[D_6]$  & 
$3$ &
$\dfrac{3}{2}$ \\
\hline
\end{tabular}
\label{tab:RC}
\end{table}

\begin{exa}[metrics $\RM_q^o,\PM_q^o$]
\label{exa:RMo}
Table~\ref{tab:RMo} has $\RM_q^o,\RM_q^o$ for $q=2,+\infty$ and
 the prominent lattices $\La_2^{\pm}$, $\La_{\infty}^{\pm}$, which were inversely designed in 
Example~\ref{exa:inverse_design}.
\medskip

If lattices have the same sign, then $\RM^o,\PM^o$ coincide with their unoriented versions by Definition~\ref{dfn:RMo}.
For example, $\PM_q^o(\La_2^+,\La_{\infty}^+)$ is the distance $M_q$ between the invariants $\PI(\La_{\infty})=(\frac{1}{4},\frac{1}{4})$ and $\PI(\La_2)=(\frac{1}{2+\sqrt{2}},\frac{1}{2+\sqrt{2}})=(1-\frac{1}{\sqrt{2}},1-\frac{1}{\sqrt{2}})$, so
$\PM_{\infty}^o(\La_2^+,\La_{\infty}^+)=\frac{3}{4}-\frac{1}{\sqrt{2}}\approx 0.04$ and $\PM_{2}^o(\La_2^+,\La_{\infty}^+)=\frac{3}{4}\sqrt{2}-1\approx 0.06$.
\medskip

Similarly, 
$\RM_q^o(\La_2^+,\La_{\infty}^+)$ is the $M_q$ distance between the root invariants $\PI(\La_{\infty})=(1,4,7)$ and $\RI(\La_2)=(2-\sqrt{2},2\sqrt{2}-1,5-\sqrt{2})$, so
$\RM_{\infty}^o(\La_2^+,\La_{\infty}^+)=\max\{\sqrt{2}-1,5-2\sqrt{2},2+\sqrt{2}\}=2+\sqrt{2}\approx 3.41$ and $\RM_{2}^o(\La_2^+,\La_{\infty}^+)=\sqrt{6(7-3\sqrt{2})}\approx 4.1$.
\medskip

By Lemma~\ref{lem:RCbounds}(b) the distance between mirror images of the same lattice equals the doubled $D_2$-chiral distance.
For example, $\PM_q^o(\La_{\infty}^+,\La_{\infty}^-)=2\PC_q[D_2](\La_{\infty})=\frac{1}{2}$ and 
$\PM_q^o(\La_2^+,\La_2^-)=2\PC_q[D_2](\La_2)=\frac{2}{2+\sqrt{2}}=2-\sqrt{2}\approx 0.59$ for $q=2,+\infty$.
\medskip

Lemma~\ref{lem:RCbounds}(b) and Table~\ref{tab:RC} also give $\RM_q^o(\La_{\infty}^+,\La_{\infty}^-)=2\RC_q[D_2](\La_{\infty})=2$ and 
$\RM_q^o(\La_2^+,\La_2^-)=2\RC_q[D_2](\La_2)=2(2-\sqrt{2})\approx 1.17$ for $q=2,+\infty$.
\medskip

Lemma~\ref{lem:reversed_signs} says that $\RM^o(\La_2^+,\La_{\infty}^-)=\RM^o(\La_2^-,\La_{\infty}^+)$ and $\PM^o(\La_2^+,\La_{\infty}^-)=\PM^o(\La_2^-,\La_{\infty}^+)$.
Using the above properties, it remains to find four distances. 
\medskip

\begin{table}[h]
\caption{Metrics $\PM_q^o$ and $\RM_q^o$ for the lattices given by their invariants in Table~\ref{tab:RC}, see Fig.~\ref{fig:QS+DC}. }
\medskip

\begin{tabular}{|l|cccc|}
\hline      
$\PM_2^o$ & $\La_{\infty}^+$ & $\La_{\infty}^-$ & $\La_2^+$ & $\La_2^-$\\
\hline

$\La_{\infty}^+$ & 0 & $\frac{1}{2}$ & $\frac{3}{4}\sqrt{2}-1\approx 0.06$ & $\frac{\sqrt{25-16\sqrt{2}}}{2\sqrt{2}}\approx 0.54$ \\

$\La_{\infty}^-$  & $\frac{1}{2}$ & 0 & $\frac{\sqrt{25-16\sqrt{2}}}{2\sqrt{2}}\approx 0.54$ & $\frac{3}{4}\sqrt{2}-1\approx 0.06$  \\

$\La_2^+$ & $\frac{3}{4}\sqrt{2}-1\approx 0.06$
& $\frac{\sqrt{25-16\sqrt{2}}}{2\sqrt{2}}\approx 0.54$ & 0 & $2-\sqrt{2}\approx 0.59$ \\

$\La_2^-$ & $\frac{\sqrt{25-16\sqrt{2}}}{2\sqrt{2}}\approx 0.54$
 & $\frac{3}{4}\sqrt{2}-1\approx 0.06$ & $2-\sqrt{2}\approx 0.59$ & 0 \\
\hline
\end{tabular}
\medskip

\begin{tabular}{|l|cccc|}
\hline      
$\PM_{\infty}^o$ & $\La_{\infty}^+$ & $\La_{\infty}^-$ & $\La_2^+$ & $\La_2^-$\\
\hline

$\La_{\infty}^+$ & 0 & $\frac{1}{2}$ & $\frac{3}{4}-\frac{1}{\sqrt{2}}\approx 0.04$ & $\frac{5}{4}-\frac{1}{\sqrt{2}}\approx 0.54$ \\

$\La_{\infty}^-$  & 
$\frac{1}{2}$ & 
0 & 
$\frac{5}{4}-\frac{1}{\sqrt{2}}\approx 0.54$ & 
$\frac{3}{4}-\frac{1}{\sqrt{2}}\approx 0.04$  \\

$\La_2^+$ & 
$\frac{3}{4}-\frac{1}{\sqrt{2}}\approx 0.04$ & 
$\frac{5}{4}-\frac{1}{\sqrt{2}}\approx 0.54$ & 
0 & 
$2-\sqrt{2}\approx 0.59$ \\

$\La_2^-$ & 
$\frac{5}{4}-\frac{1}{\sqrt{2}}\approx 0.54$ & 
$\frac{3}{4}-\frac{1}{\sqrt{2}}\approx 0.04$ & 
$2-\sqrt{2}\approx 0.59$ & 0 \\
\hline
\end{tabular}
\medskip

\begin{tabular}{|l|cccc|}
\hline      
$\RM_{2}^o$ & $\La_{\infty}^+$ & $\La_{\infty}^-$ & $\La_2^+$ & $\La_2^-$\\
\hline

$\La_{\infty}^+$ & 0 & 
2 & 
$\sqrt{6(7-3\sqrt{2})}\approx 4.1$ & 
$\sqrt{50-22\sqrt{2}}$ \\

$\La_{\infty}^-$  & 2 & 0 & $\sqrt{50-22\sqrt{2}}\approx 4.3$ & $\sqrt{6(7-3\sqrt{2})}$  \\

$\La_2^+$ & $\sqrt{6(7-3\sqrt{2})}$
& $\sqrt{50-22\sqrt{2}}\approx 4.3$ & 
0 & 
$2(2-\sqrt{2})$ \\

$\La_2^-$ & 
$\sqrt{50-22\sqrt{2}}$ & 
$\sqrt{6(7-3\sqrt{2})}\approx 4.1$ & 
$2(2-\sqrt{2})\approx 1.17$ & 0 \\
\hline
\end{tabular}
\medskip

\begin{tabular}{|l|cccc|}
\hline      
$\RM_{\infty}^o$ & $\La_{\infty}^+$ & $\La_{\infty}^-$ & $\La_2^+$ & $\La_2^-$\\
\hline

$\La_{\infty}^+$ & 0 & 2 & $2+\sqrt{2}\approx 3.41$ & 3 \\

$\La_{\infty}^-$  & 2 & 0 & 3 & $2+\sqrt{2}\approx 3.41$  \\

$\La_2^+$ & 
$2+\sqrt{2}\approx 3.41$ & 
3 & 
0 & 
$2(2-\sqrt{2})\approx 1.17$ \\

$\La_2^-$ & 
3 & 
$2+\sqrt{2}\approx 3.41$ & 
$2(2-\sqrt{2})\approx 1.17$ & 0 \\
\hline
\end{tabular}
\label{tab:RMo}
\end{table}

Proposition~\ref{prop:PMq}(a) finds
$\PM_2^o(\La_2^+,\La_{\infty}^-)$ as the minimum of the Euclidean distances from $\PI(\La_2)=(\frac{1}{2+\sqrt{2}},\frac{1}{2+\sqrt{2}})=(1-\frac{1}{\sqrt{2}},1-\frac{1}{\sqrt{2}})$ to the three points $(-\frac{1}{4},\frac{1}{4})$, $(-\frac{1}{4},\frac{1}{4})$, $(\frac{3}{4},\frac{3}{4})$ obtained from $\PI(\La_{\infty})=(\frac{1}{4},\frac{1}{4})$ by reflections in the edges of $\QT$.
The first two distances equal to $\frac{\sqrt{25-16\sqrt{2}}}{2\sqrt{2}}\approx 0.54$ are larger than the third.
\medskip

Given $\PI(\La_2)=(x_1,y_1)=(1-\frac{1}{\sqrt{2}},1-\frac{1}{\sqrt{2}})$ and $\PI(\La_{\infty})=(x_2,y_2)=(\frac{1}{4},\frac{1}{4})$, Proposition~\ref{prop:PMq}(b) computes
$\PM_{\infty}^o(\La_2^+,\La_{\infty}^-)$ for as the minimum of
$d_x=\max\{x_2-x_1,y_2+y_1\}=\frac{5}{4}-\frac{1}{\sqrt{2}}$, 
$d_y=\max\{x_2+x_1,|y_2-y_1|\}=\frac{5}{4}-\frac{1}{\sqrt{2}}$, 
$d_{xy}=\max\{x_2-x_1,1-x_2-y_2+|1-y_1-x_2|\}=\frac{1}{4}+\frac{1}{\sqrt{2}}$, so $\PM_{\infty}^o(\La_2^+,\La_{\infty}^-)=\frac{5}{4}-\frac{1}{\sqrt{2}}\approx 0.54$
\medskip

Proposition~\ref{prop:RMq}(a) computes
$\RM_2^o(\La_2^+,\La_{\infty}^-)$ as the minimum of the Euclidean distances from $\RI(\La_2)=(2-\sqrt{2},2\sqrt{2}-1,5-\sqrt{2})$ to the three points $(-1,4,7)$, $(4,1,7)$, $(1,7,4)$ obtained from $\RI(\La_{\infty})=(1,4,7)$ by reflections in the boundaries of $\TC$.
The first distance is the smallest, so $\RM_2^o(\La_2^+,\La_{\infty}^-)=\sqrt{50-22\sqrt{2}}\approx 4.3$. 
\medskip

Given $\RI(\La_2)=(r_{12},r_{01},r_{02})=(2-\sqrt{2},2\sqrt{2}-1,5-\sqrt{2})$ and
$\RI(\La_{\infty})=(s_{12},s_{01},s_{02})=(1,4,7)$,
by Proposition~\ref{prop:RMq}(b) 
$\RM_{\infty}^o(\La_2^+,\La_{\infty}^-)=\min\{d_0,d_1,d_2\}$.
Using $\MS(a,b,c,d)=\max\{|a-b|,|c-d|,\frac{1}{2}|a+b-c-d|\}$, we compute 
\medskip

\noindent
$d_0=\max\{r_{12}+s_{12}, |r_{01}-s_{01}|, |r_{02}-s_{02}|\}$ \\
$=\max\{3-\sqrt{2},5-2\sqrt{2},2+\sqrt{2}\}
=2+\sqrt{2}\approx 3.4,$ 
\medskip

\noindent
$d_1=
\max\{\MS(r_{12},r_{01},s_{12},s_{01}),|r_{02}-s_{02}|\}=$\\
$=\max\{\MS(r_{12},r_{01},s_{12},s_{01}),2+\sqrt{2}\}=$ \\
$=\max\{\MS(2-\sqrt{2},2\sqrt{2}-1,1,4),2+\sqrt{2}\}=$ \\
$=\max\{\max\{3(\sqrt{2}-1),3,2-\frac{1}{\sqrt{2}}\}, 2+\sqrt{2}\}=\max\{3,2+\sqrt{2}\}=2+\sqrt{2},$
\medskip

\noindent
$d_2=\max\{|r_{12}-s_{12}|,\MS(r_{01},r_{02},s_{01},s_{02})\}=$ \\
$=\max\{\sqrt{2}-1, \MS(2\sqrt{2}-1,5-\sqrt{2},4,7)\}=$ \\
$=\max\{\sqrt{2}-1, \max\{6-3\sqrt{2},3,\frac{7-\sqrt{2}}{2})\}=3$, hence $\RM_{\infty}^o(\La_2^+,\La_{\infty}^-)=3$.
\eexa
\end{exa}

\section{Bi-continuity of the map from obtuse superbases to invariants}
\label{sec:continuity}

\index{obtuse superbase}

This section studies continuity of the bijection $B\mapsto\La(B)$, where an obtuse superbase $B$ and its lattice $\La(B)$ are considered under isometry, rigid motion, dilation, and homothety.
Theorems~\ref{thm:isometric_superbases} and~\ref{thm:lattices2D_classification} established the bijections $\LIS\to\SBI\to\RIS$, $\La\mapsto B\to\RI(B)=\RI(\La)$, mapping
any lattice $\La\subset\R^2$ to its obtuse superbase $B$ (unique under isometry) and then to the complete invariant $\RI(\La)$.
Hence, the Lattice Isometry Space $\LIS(\R^2)$ with a root metric $\RM$ can be identified with the Root Invariant Space $\RIS=(\TC,d)$ with a metric $d$ on the triangular cone $\TC\subset\R^3$.

\index{obtuse superbase}

\begin{thm}[continuity of $\SBI\to\LIS=\RIS$, {\cite[Theorem 7.5]{kurlin2024mathematics}}]
\label{thm:SBI->RIS}
\textbf{(a)}
Let $q\in[1,+\infty]$ and lattices $\La,\La'\subset\R^2$ have obtuse superbases $B$ and $B'$ whose vectors have a maximum length $l$.
If $\SIM_{\infty}(B,B')=\de\geq 0$,
then $\RM_q(\La,\La')\leq 3^{1/q}\sqrt{2l\de}$.
Hence the bijection $\SBI(\R^2)\to\LIS(\R^2)$ is continuous in the metrics $\SIM_{\infty}$ and $\RM_q$.
\medskip

\noindent
\textbf{(b)}
In the conditions above, the projected metric satisfies
 $\PM_q(\La,\La')\leq 2^{1/q}3\sqrt{2\de/l}$, so the bijection $\SBH(\R^2)\to\LHS(\R^2)$ is continuous in the metrics $\SHM_{\infty},\PM_q$.
\medskip

\noindent
\textbf{(c)}
In the oriented case, if $\de\to 0$, then $\RM_{q}^o(\La,\La')\to 0$ and $\PM_{q}^o(\La,\La')\to 0$, so both maps $\SBR(\R^2)\to\LRS(\R^2)$ and $\SBD(\R^2)\to\LDS(\R^2)$ are continuous.
\ethm
\end{thm}

Theorem~\ref{thm:SBI->RIS} is proved for the metrics $\RM_q,\PM_q$ only to give explicit upper bounds.
A similar argument proves continuity for any metrics $\RM,\PM$ in Definition~\ref{dfn:RM} based on a metric $d$ satisfying $d(u,v)\to 0$ when $\vec u\to v$ coordinate-wise. 
All Minkowski norms in $\R^n$ are topologically equivalent \cite{norms} due to the bounds  $||v||_q\leq ||v||_{r}\leq n^{\frac{1}{q}-\frac{1}{r}}||v||_q$ for any $1\leq q\leq r$, hence continuity for one value of $q$ is enough.
Theorem~\ref{thm:SBI->RIS} implies continuity of $\SBR\to\RIS^o$, because closeness of superbases under rigid motion is a stronger condition than under isometry.
\medskip

Example~\ref{exa:rect_deformation} illustrates Theorem~\ref{thm:SBI->RIS} and shows that the root invariant changes continuously for a deformation when a reduced basis changes discontinuously.

\begin{exa}[continuity of root invariants under deformaiton]
\label{exa:rect_deformation}
The obuse superbase $\vec v_1=(1,0)$, $\vec v_2(t)=(-t,2)$, $\vec v_0(t)=(t-1,-2)$ continuously deforms for $t\in[0,1]$ in Fig.~\ref{fig:rect_basis_discontinuity}.
The basis of $\vec v_1,\vec v_2(t)$ is reduced (non-acute) for $t\in[0,\frac{1}{2}]$ and at the critical moment $t=\frac{1}{2}$ changes to its mirror image $\vec v_1,\vec v_0(t)$, which remains reduced for $t\in[\frac{1}{2},1]$.
The obtuse superbase of unordered vectors $\{\vec v_1,\vec v_2(t),\vec v_0(t)\}$ keeps changing continuously because $\vec v_2(t),\vec v_0(t)$ only swap their places at $t=\frac{1}{2}$.
\myskip

The discontinuity of the obtuse superbases in the above deformation emerges at $t=1$ when the final superbase of $(1,0),(-1,2),(0,-2)$ becomes a mirror image of the initial superbase of $(1,0),(0,2),(-1,-2)$, not related by rigid motion, though both (unordered) superbases at $t=0,1$ generate the same lattice with the rectangular cell $1\times 2$.  
The root invariants are $r_{12}=\sqrt{t}$, $r_{01}=\sqrt{1-t}$, $r_{02}=\sqrt{4-t+t^2}$.
Since $4-t+t^2\geq\frac{15}{4}\geq\max\{t,1-t\}$ for $t\in[0,1]$, the root invariant can be written as 
$$\RI(\La(t))=\left\{\begin{array}{ll}
(\sqrt{t},\sqrt{1-t},\sqrt{4-t+t^2}) \text{ for } t\in[0,\frac{1}{2}],\\
(\sqrt{1-t},\sqrt{t},\sqrt{4-t+t^2}) \text{ for } t\in[\frac{1}{2},1].
\end{array} \right.$$

\begin{figure}[h]
\includegraphics[width=\textwidth]{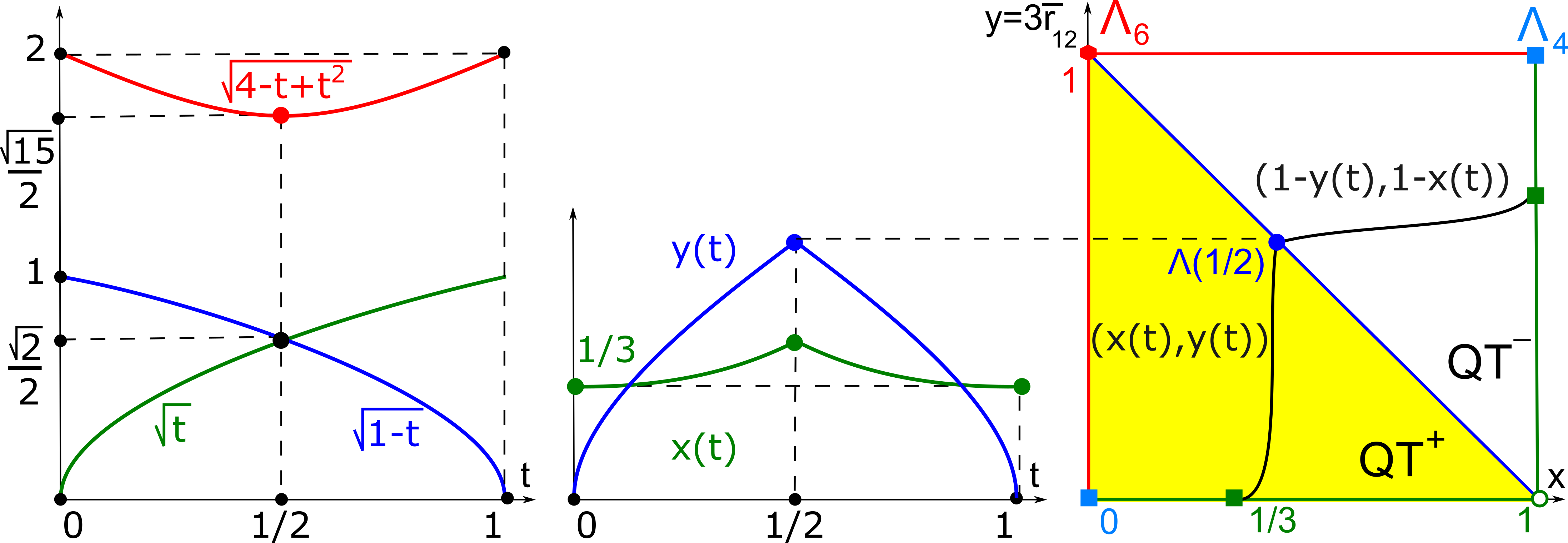}
\caption{
\textbf{Left}: graphs of root products in $\RI(\La(t))$, see Example~\ref{exa:rect_deformation}.
\textbf{Middle}: graphs of the components in $\PI(\La(t))$.
\textbf{Right}: the continuous path of $\PI(\La(t))$ in the quotient square $QS$.}
\label{fig:rect_deformation}
\end{figure}

By Definition~\ref{dfn:PI} the size is $|\La(t)|=r_{12}+r_{01}+r_{02}=\sqrt{t}+\sqrt{1-t}+\sqrt{4-t+t^2}$.
The projected invariant is $\PI(\La(t))=(x(t),y(t))$, see Fig.~\ref{fig:rect_deformation}, where
$$x(t)=\frac{\sqrt{4-t+t^2}-\max\{\sqrt{t},\sqrt{1-t}\}}{\sqrt{t}+\sqrt{1-t}+\sqrt{4-t+t^2}},\quad
y(t)=\frac{3\min\{\sqrt{t},\sqrt{1-t}\}}{\sqrt{t}+\sqrt{1-t}+\sqrt{4-t+t^2}}.$$
If $t=0$ or $t=1$, then $\RI(\La(t))=(0,1,2)$, $|\La(t))=3$, $\PI(\La(t))=(\frac{1}{3},0)$.
If $t=\frac{1}{2}$, then $\sqrt{t}=\sqrt{1-t}=\frac{\sqrt{2}}{2}$, $\sqrt{4-t+t^2}=\frac{\sqrt{15}}{2}$, $|\La(\frac{1}{2}))=\sqrt{2}+\frac{\sqrt{15}}{2}$.
So
$$\RI\left(\La\Big(\frac{1}{2}\Big)\right)=\left(\frac{\sqrt{2}}{2},\frac{\sqrt{2}}{2},\frac{\sqrt{15}}{2}\right),\; 
\PI\left(\La\Big(\frac{1}{2}\Big)\right)=\left(\frac{\sqrt{15}-\sqrt{2}}{\sqrt{15}+2\sqrt{2}},\frac{3\sqrt{2}}{\sqrt{15}+2\sqrt{2}}\right).$$
The last point is approximately $(0.37,0.63)$ in the diagonal $x+y=1$ of $\QS$.
Under the symmetry $t\lra 1-t$, all the functions above remain invariant and $\La(t)$ changes its sign.
Both paths $\RI(\La(t))$ and $\PI^o(\La(t))\in\QS$ are continuous everywhere, while the obtuse superbasis is discontinuous (under rigid motion) at $t=0,1$.
\eexa
\end{exa}

Theorem~\ref{thm:RIS->SBI} below proves the inverse continuity 
of $\RIS\to\SBI$ and a weaker claim in the oriented case, saying that we can choose an obtuse superbase $B'$ of a perturbed lattice arbitrarily close to a given superbase $B$ of an original lattice.

\begin{thm}[continuity of $\LIS\to\SBI$, {\cite[Theorem 7.7]{kurlin2024mathematics}}]
\label{thm:RIS->SBI}
\textbf{(a)}
For $q\in[1,+\infty]$, let lattices $\La,\La'$ in $\R^2$ satisfy $\RM_q(\La,\La')\leq\de$.
For any obtuse superbase $B$ of $\La$, there is an obtuse superbase $B'$ of $\La'$ such that $\SIM_{\infty}(B,B')\leq\SRM_{\infty}(B,B')\to 0$ as $\de\to 0$.
\medskip

\noindent
\textbf{(b)}
The bijection $\LIS(\R^2)\to\SBI(\R^2)$ is continuous in the metrics $\RM_q,\SIM_{\infty}$.
$\LRS(\R^2)\to\SBR(\R^2)$ is continuous in $\RM_q^o,\SIM_{\infty}^o$ at non-rectangular lattices.
\medskip

\noindent
\textbf{(c)}
The above conclusions hold for lattices under dilation and homothety.
\ethm
\end{thm}

Corollary~\ref{cor:rect_discontinuity} shows that Theorem~\ref{thm:RIS->SBI}(b) is the strongest possible continuity in the oriented case.
In $\R^3$, a similar discontinuity around high-symmetry lattices will be much harder to resolve for continuous invariants even under isometry \cite{kurlin2022complete}.

\begin{cor}[partial discontinuity of $\LRS\to\SBR$, {\cite[Cor.~7.9]{kurlin2024mathematics}}]
\label{cor:rect_discontinuity}
The bijection $\LRS\to\SBR$ is discontinuous in the metrics $\RM_{\infty},\SIM_{\infty}^o$ at all rectangular lattices.
\ethm
\end{cor}

Corollary~\ref{cor:rect_discontinuity} should be positively interpreted in the sense that we need to study lattices under rigid motion by their complete oriented root invariants in the continuous space $\LRS(\R^2)$ rather than in terms of reduced bases (or, equivalently, obtuse superbases due to Proposition~\ref{prop:reduced_bases}b), which are inevitably discontinuous.
\medskip

Proposition~\ref{prop:cont_chiral_distances} shows that all $G$-chiral distances $\RC[G]:\LIS(\R^2)\to\R$ and $\PC[G]:\LHS(\R^2)\to\R$ are continuous in any metrics $\RM,\PM$ from Definition~\ref{dfn:RM}.

\begin{prop}[continuous chiral distances, {\cite[Proposition 7.10]{kurlin2024mathematics}}]
\label{prop:cont_chiral_distances}
For a crystallographic point group $G$ and lattices $\La_1,\La_2$ in $\R^2$, we have 
$$|\RC[G](\La_1)-\RC[G](\La_2)|\leq \RM(\La_1,\La_2),$$
$$|\PC[G](\La_1)-\PC[G](\La_2)|\leq \PM(\La_1,\La_2)$$
for any metrics $\RM$ and $\PM$. 
\ethm
\end{prop}

Now Remark~\ref{rem:structures} summarises a wide range of rich mathematical structures that can be considered on the lattice spaces in addition to continuous metrics.

\begin{rem}[linear structure, scalar product on lattices]
\label{rem:structures}
Since the triangular cone $\TC$ in Fig.~\ref{fig:TC} is convex, we can consider any convex linear combination of root invariants $t\RI(\La_1)+(1-t)\RI(\La_2)\in\TC$, $t\in[0,1]$.
The resulting root invariant determines (an isometry class of) the new lattice that can be denoted by $t\La_1+(1-t)\La_2$.
The average of the square and hexagonal lattices with $\RI(\La_4)=(0,1,1)$, $\RI(\La_6)=(1,1,1)$ has $\RI=(\frac{1}{2},1,1)$.
The new lattice $\frac{1}{2}(\La_4+\La_6)$ is centred rectangular and has the basis $\vec v_1=(\sqrt{\frac{3}{2}},0)$ and $\vec v_2=(-\frac{1}{9}\sqrt{\frac{3}{2}},\frac{4}{9}\sqrt{\frac{15}{2}})$.
We can define similar sums in $\LHS(\R^2)$ due to the convexity of the triangle $\QT$. 
The usual scalar product of vectors in $\R^3$ defines the positive product of root invariants: $\RI(\La_4)\cdot\RI(\La_6)=(0,1,1)\cdot(1,1,1)=2$.
\erem
\end{rem}

In conclusion, Problem~\ref{pro:lattices2D} was resolved by the new invariants $\RI,\RI^o,\PI,\PI^o$ classifying all 2D lattices under four equivalences, see a summary in Table~\ref{tab:lattices2D_geocodes}.
\myskip

\noindent
\ref{pro:lattices2D}(a) 
Completeness of invariants: 
Theorem~\ref{thm:lattices2D_classification} and Corollary~\ref{cor:lattices2D_classification}.
\myskip

\noindent
\ref{pro:lattices2D}(b) Reconstruction: 
Lemma~\ref{lem:superbase_reconstruction} and
Proposition~\ref{prop:inverse_design} with Example~\ref{exa:inverse_design}.
\medskip

\noindent
\ref{pro:lattices2D}(c-e) Continuous metrics: 
Definitions~\ref{dfn:RM} and~\ref{dfn:RMo}, Theorems~\ref{thm:SBI->RIS} and~\ref{thm:RIS->SBI}.
\myskip

\noindent
\ref{pro:lattices2D}(f,g) Realisability and Euclidean embeddability: Corollary~\ref{cor:lattices2D_spaces}.
\myskip

\noindent
\ref{pro:lattices2D}(h) Computability of metrics: 
Propositions~\ref{prop:RMq},\ref{prop:PMq} and Examples~\ref{exa:RM},~\ref{exa:RC}. 
\myskip

\begin{table}[h]
\caption{A summary of classifications of all lattices $\La\subset\R^2$ under four equivalence relations.}

\begin{tabular}{p{20mm}|p{25mm}|c|c|c}      
equivalence  & complete invariant & configuration space & continuous metric & visual results \\
\hline

isometry  & root invariant  & $\LIS(\R^2)\lra\TC$ & root metric &  
Theorem~\ref{thm:lattices2D_classification} \\

 & $\RI(\La)$ & triangular cone & $\RM$ & Fig.~\ref{fig:TC} (left) \\
\hline

rigid & oriented  & $\LRS(\R^2)\lra\DC$ & oriented & Theorem~\ref{thm:lattices2D_classification} \\

motion & invariant $\RI^o(\La)$ & doubled cone & metric $\RM^o$ & Fig.~\ref{fig:QS+DC} (right) \\
\hline

homothety & projected  & $\LHS(\R^2)\lra\QT$ & projected & Corollary~\ref{cor:lattices2D_classification} 
\\

 & invariant $\PI(\La)$ & quotient triangle & metric $\PM$ & Fig.~\ref{fig:QT+QS} (left) \\
\hline

 & oriented & $\LDS(\R^2)\lra\QS$ & oriented & Corollary~\ref{cor:lattices2D_classification} \\

dilation & projected & quotient & projected & Fig.~\ref{fig:QT+QS} (right) \\

 & invariant $\PI^o(\La)$ & square & metric $\PM^o$ & Fig.~\ref{fig:DT}
\end{tabular}
\label{tab:lattices2D_geocodes}
\end{table}
\medskip

The key contributions of this chapter are the easily computable metrics in Definitions~\ref{dfn:RM},\ref{dfn:RMo}, which led to continuous real-valued deviations of lattices from their higher symmetry neighbours.
The chiral distances in Definition~\ref{dfn:RC} continuously extend the binary chirality by explicit formulae in Propositions~\ref{prop:RC} and \ref{prop:PC}.
\medskip

The discontinuity of basis reductions in \cite[Theorem~15]{widdowson2022average} was proved with a simple metric on lattice bases without isometry. 
When we consider obtuse superbases under isometry, the continuity holds in Theorem~\ref{thm:RIS->SBI} under isometry for all lattices and under rigid motion for non-rectangular lattices.
For rigid motion, when orientation is preserved, Corollary~\ref{cor:rect_discontinuity} proves discontinuity at any rectangular lattice in $\R^2$.
\medskip

It was important to clarify the above discontinuity of bases in Corollary~\ref{cor:rect_discontinuity} since the 3-dimensional case is much harder to resolve even under isometry \cite{kurlin2022complete}. 
\medskip

The structures in Remark~\ref{rem:structures} help treat lattices as vectors in a meaningful way (independent of a basis), for example, as inputs or outputs in machine learning algorithms.
Paper \cite{bright2023geographic} visualises for the first time millions of 2D lattices extracted from real crystals in the Cambridge Structural Database (CSD), see the Python code at https://github.com/MattB-242/Lattice\_Invariance.
\medskip
 
Lattice invariants can be used as a first ultra-fast step to find (near-)duplicates of a potentially new material in all existing experimental datasets.
The forthcoming work \cite{kurlin2022complete,bright2021complete} extends the isometry classification of Theorem~\ref{thm:lattices2D_classification} to $\R^3$.
The next chapter will introduce more advanced distance-based invariants of general periodic point sets.

\bibliographystyle{plain}
\bibliography{Geometric-Data-Science-book}

%
%
%



\chapter{Density functions of periodic sets of points in $\R^n$ and intervals in $\R$}
\label{chap:densities} 

\abstract{
This chapter adapts the general geo-mapping problem to periodic sets of points under isometry in high dimensions, motivated by periodic crystals in dimension 3.
We introduce density functions, which extend the scalar point density to Lipschitz continuous isometry invariants depending on a variable radius of balls centred at given points.
These functions can be efficiently computed at discrete radii in low dimensions and are generically complete for periodic point sets under isometry in $\R^3$.
In dimension $n=1$, the density functions are analytically computable for periodic sets of intervals.  
 }

\section{Periodic point sets in $\R^n$ and their invariant density functions}
\label{sec:periodic_sets}

This section follows \cite{widdowson2022resolving} to adapt Geo-Mapping Problem~\ref{pro:geocodes} to periodic point sets in $\R^n$. These sets extend lattices from Definition~\ref{dfn:lattice_cell}, as introduced below.

\index{motif}
\index{periodic point set}

\begin{dfn}[motif, periodic point set in $\R^n$] 
\label{dfn:periodic}
Let vectors $\vec v_1,\dots,\vec v_n\in\R^n$ form a basis of $\R^n$,  define the lattice $\La=\{\sum\limits_{i=1}^l c_i\vec  v_i \mid c_1,\dots,c_l\in\Z\}$ and the unit cell $U=\{\sum\limits_{i=1}^n x_i\vec  v_i \mid x_1,\dots,x_n\in[0,1)\}\subset\R^n$. 
For any finite set of points (called a \emph{motif}) $M\subset U$, the sum $S=M+\La=\{\vec p+\vec v \mid \vec p\in M, \vec v\in\La\}$ is an \emph{periodic point set}.
\edfn
\end{dfn}

\begin{figure}[h!]
\centering
\includegraphics[height=21mm]{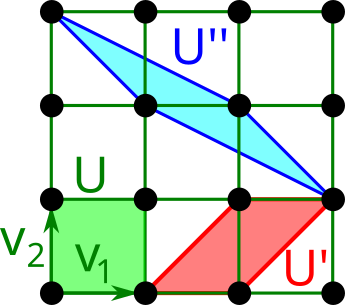}
\hspace*{0.5mm}
\includegraphics[height=21mm]{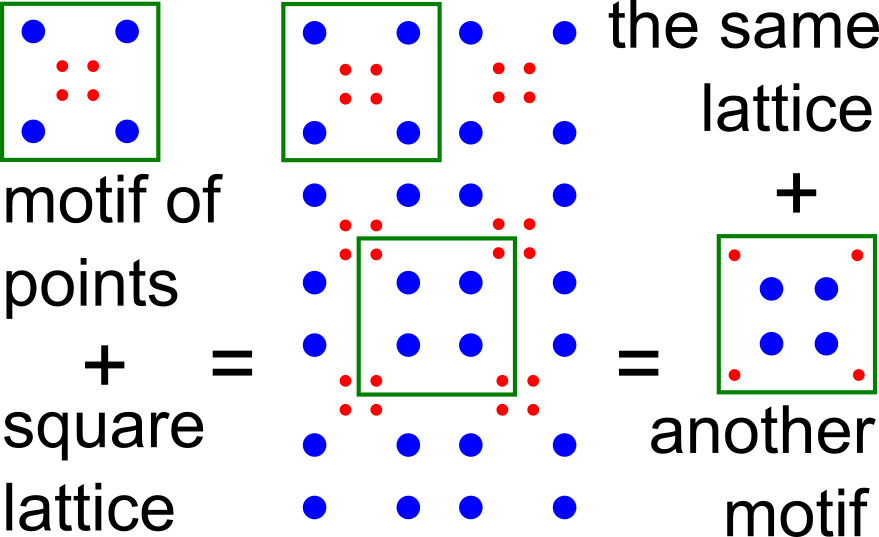}
\hspace*{0.5mm}
\includegraphics[height=21mm]{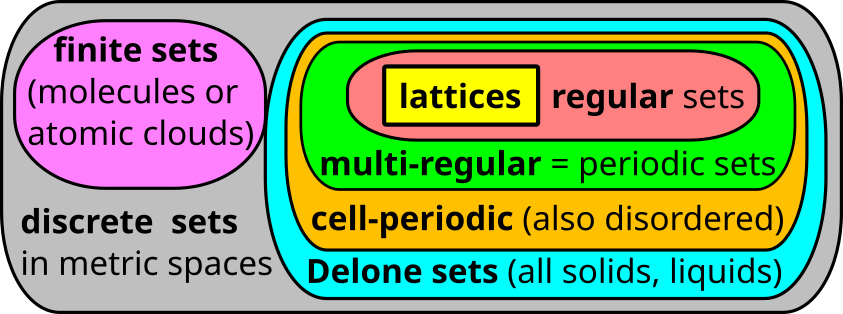}
\caption{\textbf{Left}: a lattice can be defined by many primitive bases.
\textbf{Middle}: a periodic set can be defined by different pairs (basis, motif).
\textbf{Right}: a hierarchy of discrete sets, which model all crystalline materials (periodic crystals) and amorphous solids with points at atomic centers, see Definition~\ref{dfn:periodic}.  }
\label{fig:lattice_periodic_set_hierarchy}
\end{figure}

Any unit cell $U$ includes only a partial boundary: we exclude the points with any coefficient $t_i=1$, $i=1,\dots,l$, for convenience.
Then $\R^n$ is tiled by the shifted cells $\{U+\vec v \mid \vec v\in\La\}$ without overlaps.
Any lattice is an example of a periodic set with one point in a motif.
Any periodic point set $S=M+\La$ can be considered a finite union $\bigcup_{p\in M}(\vec p+\La)$ of lattices whose origins are shifted to all $p\in M=S\cap U$.
\myskip

If we double a unit cell in one direction, e.g. by taking the basis $2\vec v_1,\vec v_2,\dots,\vec v_n$, the doubled motif $M\cup(M+\vec v_1)$ with the sublattice on the new basis defines the original periodic point set $S=M+\La$.
A basis and its cell $U$ of $S$ are called \emph{primitive} if $S\cap U$ has the smallest size among all unit cells $U$ of $S$.
Fig.~\ref{fig:lattice_periodic_set_hierarchy}~(left) shows a square lattice in $\R^2$, which (as any lattice) can be generated by infinitely many primitive bases.
Even if we fix a basis, Fig.~\ref{fig:lattice_periodic_set_hierarchy}~(middle) shows that different motifs in the same primitive cell $U$ define equivalent periodic sets, which differ only by translation. 
\smallskip

Finite and periodic point sets represent molecules and periodic crystals at the atomic scale by considering zero-sized points at all atomic centers.
Chemical bonds can be modelled by straight-line edges between atomic centers.
However, even the strongest covalent bonds within a molecule depend on various thresholds for distances and angles.
So these bonds are not real sticks and only abstractly represent inter-atomic interactions, while atomic nuclei are real.
We model all materials at the fundamental level of atoms. 
\myskip

Definition~\ref{dfn:crystal_spaces} extends moduli spaces of lattices in Definition~\ref{dfn:lattice_spaces} to periodic point sets with up to $m$ points in their motifs.
One physically justified metric of all these spaces is the bottleneck distance $\BD$ from Example~\ref{exa:metrics}(b), which quantifies atomic vibrations as a maximum deviation of all points from their original positions.
However, $\BD$ involves a minimisation over bijections between infinite periodic sets and also over infinitely many equivalences, such as isometries.
Hence, periodic sets need a simpler distance metric that should be efficiently computable.
Nonetheless, the bottleneck distance allows us to define the concept of a generic set in crystal spaces below. 
We use the word \emph{crystal} instead of the \emph{periodic point set} to keep all names short. 

\index{crystal spaces}
\index{Crystal Rigid Space}
\index{Crystal Isometry Space}
\index{Crystal Dilation Space}
\index{Crystal Homothety Space}

\begin{dfn}[moduli spaces of periodic point sets and generic subspaces]
\label{dfn:crystal_spaces}
In all cases below, we consider all periodic point sets $S\subset\R^n$ with motifs of up to $m$ points.
\myskip

\nt
\tb{(a)}
The \emph{Crystal Rigid Space} $\CIMS(\R^n;m)$: periodic point sets under rigid motion.
\myskip

\nt
\tb{(b)}
The \emph{Crystal Isometry Space} $\CIMS(\R^n;m)$: periodic point sets under isometry.
\myskip

\nt
\tb{(c)}
The \emph{Crystal Dilation Space} $\CRDS(\R^n;m)$: periodic point sets under dilation.
\myskip

\nt
\tb{(d)}
The \emph{Crystal Homothety Space} $\CRHS(\R^n;m)$: periodic point sets under homothety.
\myskip

\nt
\tb{(e)}
For any space $X$ above, a subspace $Y\subset X$ is \emph{dense} (or \emph{generic}) if, for any $\ep>0$, any periodic point set $Q$ representing a class in $X$ can be obtained from some $S\subset\R^n$ representing a class in $Y$ by perturbing any point of $S$ up to Euclidean distance $\ep$.
\edfn
\end{dfn}

\begin{figure}[h!]
\centering
\includegraphics[width=\textwidth]{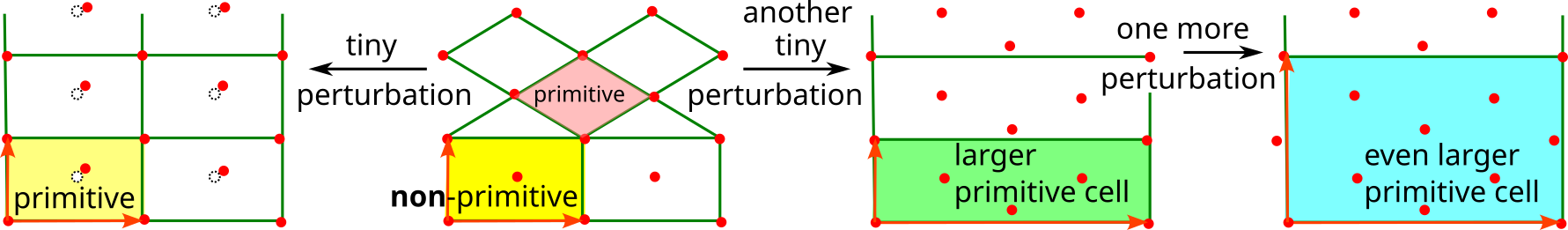}
\caption{Almost any noise arbitrarily scales up a primitive unit cell of any periodic point set.  }
\label{fig:hexagonal-lattice-cell-scales-up}
\end{figure}

Fig.~\ref{fig:hexagonal-lattice-cell-scales-up} illustrates for an initial hexagonal lattice of red points that the discontinuity of cell-based representations substantially worsens for periodic point sets in comparison with lattices, which have only point in a motif.
Indeed, if we arbitrarily extend a given primitive unit cell $U$ of $S$ to a larger cell $U'$, then almost any perturbation of one atom in $U'$ (and all its periodic copies obtained by translations along the edges of $U'$) makes $U'$ primitive.
Problem~\ref{pro:periodic_density} extends Problem~\ref{pro:1-periodic} to general periodic point sets under isometry.
The stronger equivalence of rigid motion will be considered in Chapter~\ref{chap:isosets}.

\begin{pro}[isometry invariants of periodic point sets in $\R^{n}$]
\label{pro:periodic_density}
Design an invariant $I$ on the Crystal Isometry Space $\CIMS(\R^n;m)$
satisfying the following conditions.
\smallskip

\noindent
\tb{(a)} 
\emph{Generic completeness:}
let $S,Q$ be any \emph{generic} sets whose isometry classes are in a dense subspace $Y\subset\CIMS(\R^n;m)$, then $S,Q$ are isometric if and only if $I(S)=I(Q)$.
\myskip

\noindent
\tb{(b)} 
\emph{Metric:} 
there is a distance $d$ on the Crystal Isometry Space $\CIMS(\R^n;m)$ satisfying all metric axioms in Definition~\ref{dfn:metrics}(a). 
\myskip

\noindent
\tb{(c)} 
\emph{Continuity:} 
there is a constant $\la>0$, such that, for all sufficiently small $\ep>0$, if a periodic point set $Q$ is obtained by perturbing every point of a periodic point set $S\subset\R^n$ up to Euclidean distance $\ep$, then $d(I(S),I(Q))\leq\la\ep$.
\myskip

\noindent
\tb{(d)} 
\emph{Computability:} the invariant $I$, a reconstruction of $S\subset\R^{n}$ from $I(S)$, and the metric $d(I(S),I(Q))$ can be computed in times that depend polynomially on the dimension $n$ and the maximum motif size of periodic point sets $S,Q$.
\epro
\end{pro}
 
Condition~\ref{pro:periodic_density}(a) includes only generic completeness, which will be tackled in this and the next chapter, while the full completeness will be resolved in Chapter~\ref{chap:isosets}.
\myskip

Now we introduce density functions that satisfy all conditions of Problem~\ref{pro:periodic_density} in $\R^3$, though polynomial-time algorithms in \ref{pro:periodic_density}(d) will be approximate.
Chapter~\ref{chap:PDD-periodic} will introduce newer isometry invariants and metrics with exact polynomial-time algorithms. 
For any $t\geq 0$ and $p\in\R^n$, let $\bar B(p;t)\subset\R^n$ be the closed ball of radius $t$ centred at $p$.

\index{density function}

\begin{dfn}[density functions of a periodic point set in $\R^n$]
\label{dfn:densities}
Let a periodic set $S=\La+M\subset\R^n$ have a unit cell $U$ of volume $\vol[U]$.
For any integer $k\geq 0$, let $U_k(t)$ be the region within the cell $U$ covered by exactly $k$ closed balls $\bar B(p;t)$ with a radius $t\geq 0$ and centres at all points of $S$.
The $k$-th \emph{density function} is $\psi_k[S](t)=\dfrac{\vol[U_k(t)]}{\vol[U]}$.
The \emph{density fingerprint} is the infinite sequence $\Psi[S]=\{\psi_k[S](t)\}_{k=0}^{+\infty}$.
\edfn
\end{dfn}

Notice that closed balls $\bar B(p;t)$ are considered for all points $p\in S$, not restricted to the motif $M=S\cap U$. 
The $0$-th density $\psi_0[S](t)$ measures the subset of $U$ that is not covered by any balls $\bar B(p;t)$ for $p\in S$, i.e. the fractional volume subset of all points $q\in U$ that are more than $t$ away from all points of $S$.
For $k\geq 1$, 
$\psi_k[S](t)$ measures the fractional volume of all $k$-fold intersections $\bigcap\limits_{p_1,\dots,p_k\in S}\bar B(p_i;t)$ within 
$U$.
\myskip

Since each density function $\psi_k[S](t)$ depends on $t\in[0,+\infty)$, our computations in  dimensions $n=2,3$ use uniformly sampled radii $t$.
For any fixed radius $t$, the density function $\psi_k[S](t)$ will be efficiently computed in section~\ref{sec:Brillouin}.
In dimension $n=1$, all density functions will be analytically computed, also for more general periodic sequences of intervals.
Fig.~\ref{fig:densities_sq} and~\ref{fig:densities_hex} illustrate the \emph{densigrams} that combine several density functions in one diagram for the square and hexagonal lattices in $\R^2$. 
 
\newcommand{\cheight}{23mm}
\newcommand{\bwidth}{49mm}
\newcommand{\dwidth}{70mm}
\begin{figure}[h!]
\parbox{\bwidth}{
  \includegraphics[height=\cheight]{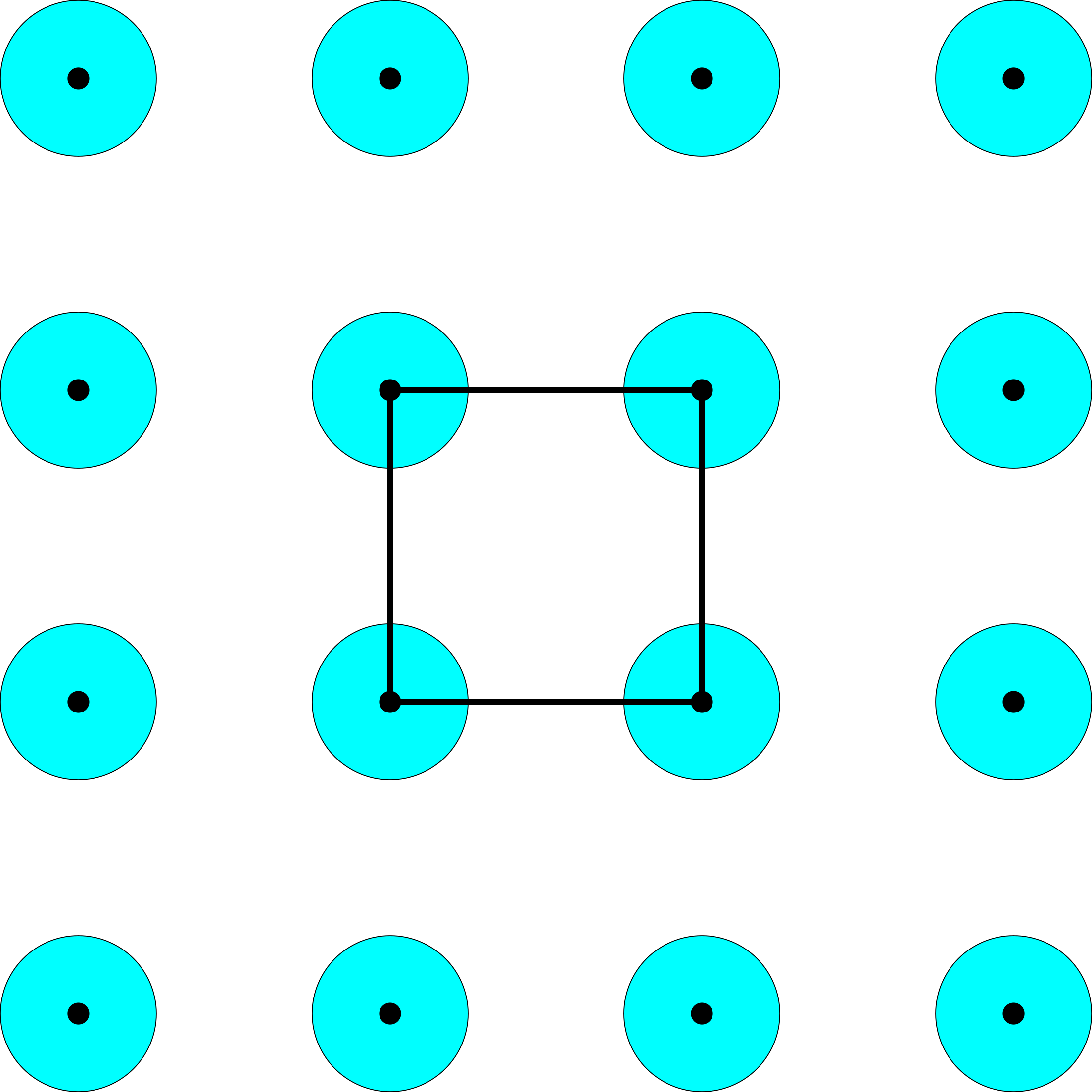}
  \hspace*{1mm}
  \includegraphics[height=\cheight]{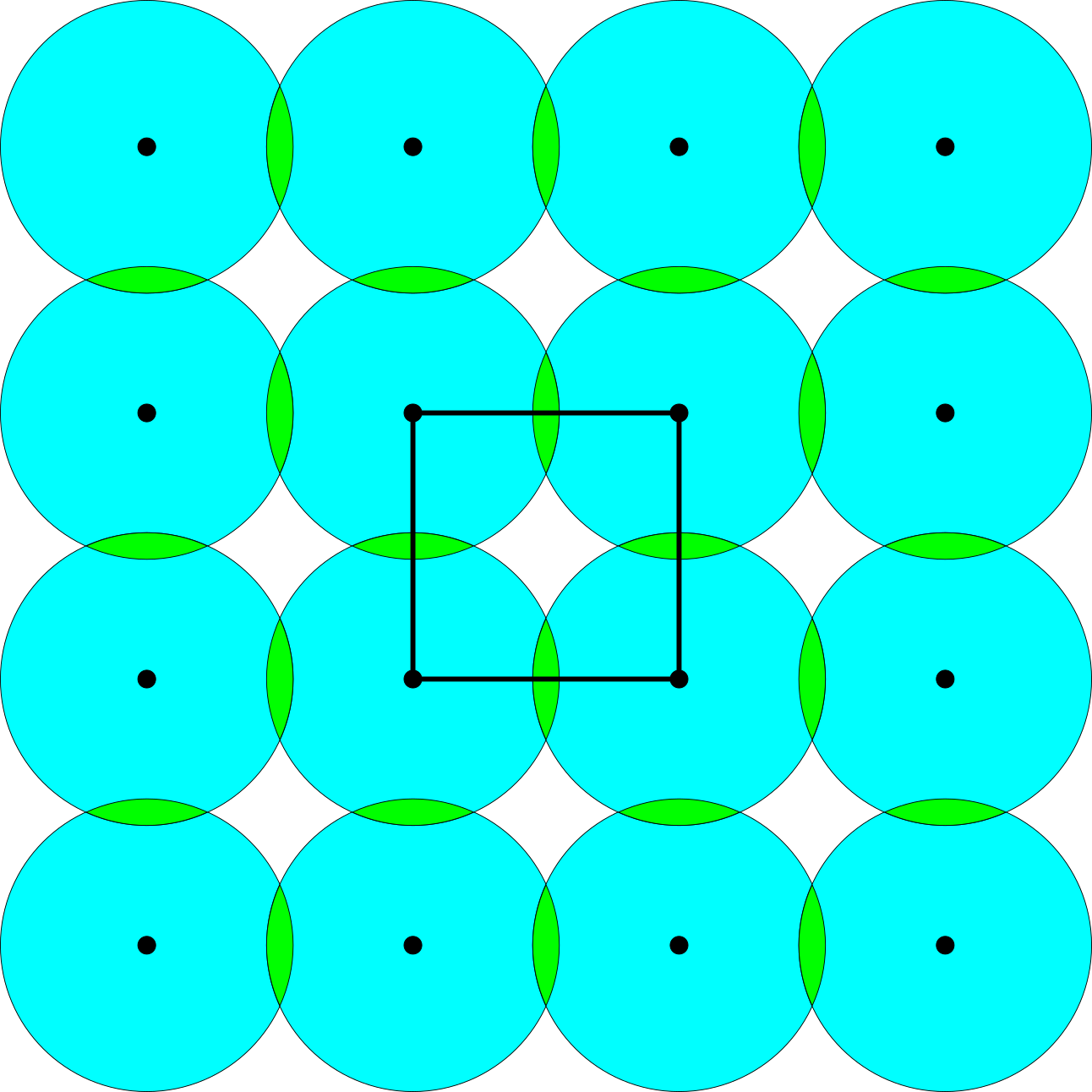}
  \medskip
  
  \includegraphics[height=\cheight]{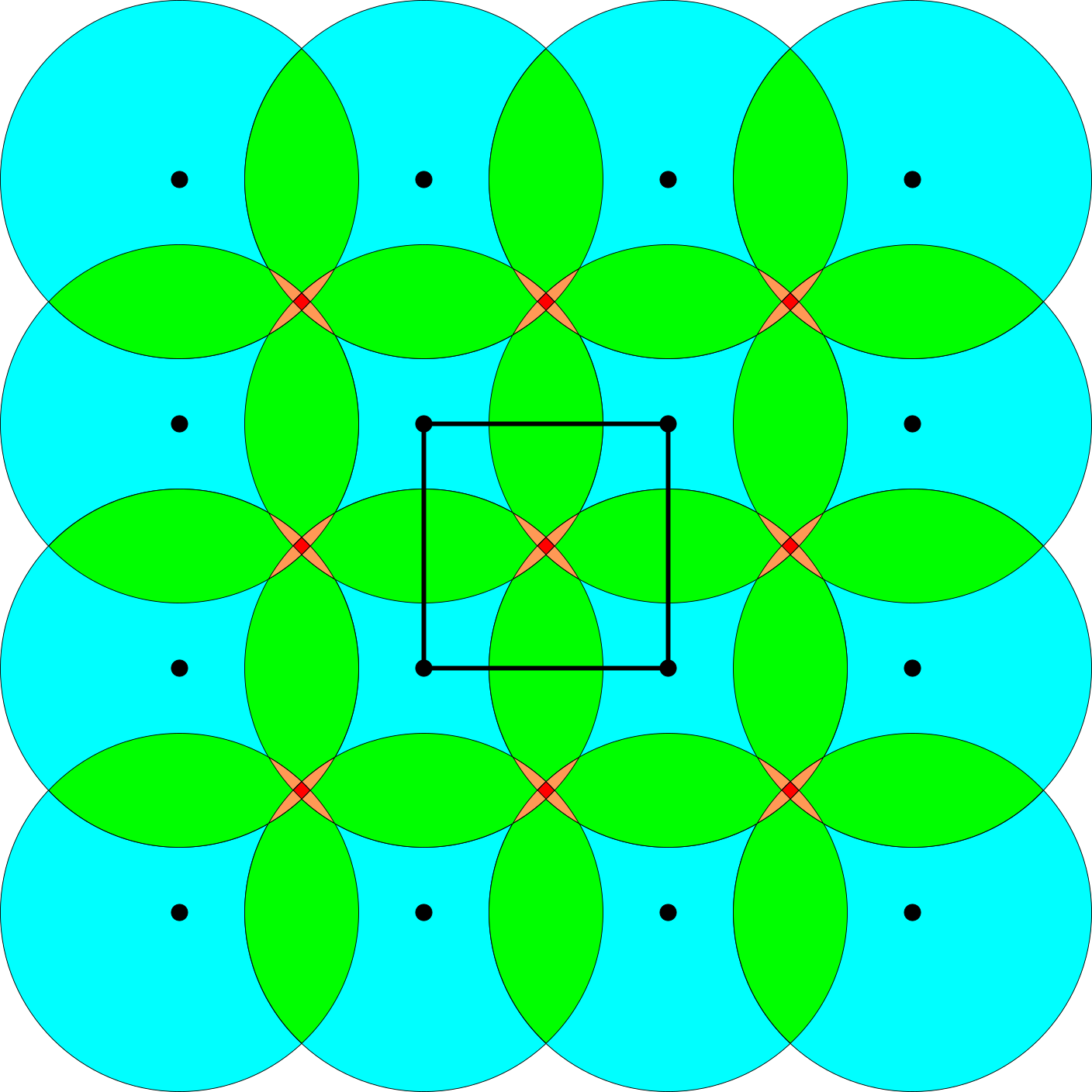}
  \hspace*{1mm}
  \includegraphics[height=\cheight]{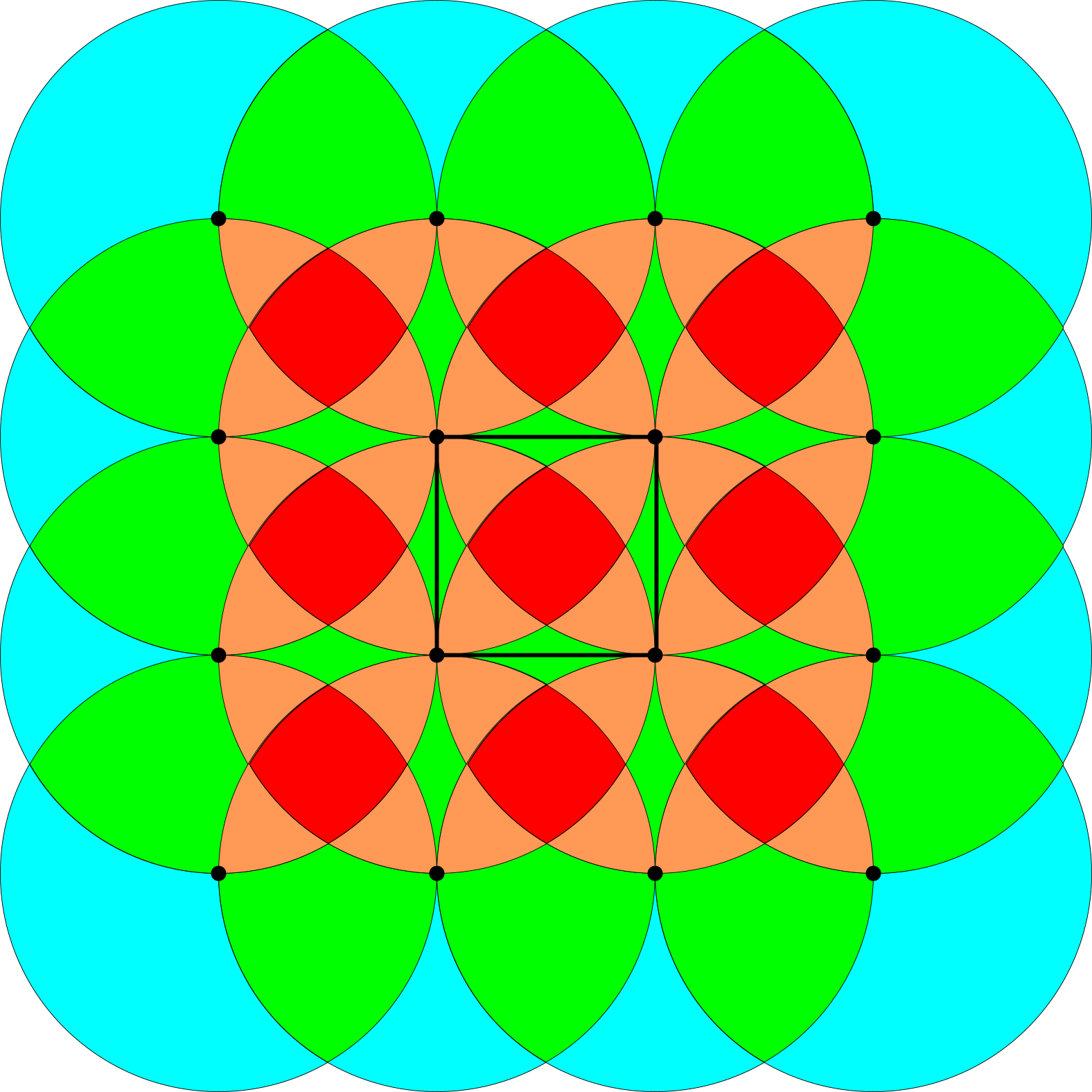}}
  \parbox{80mm}{
  \includegraphics[width=\dwidth]{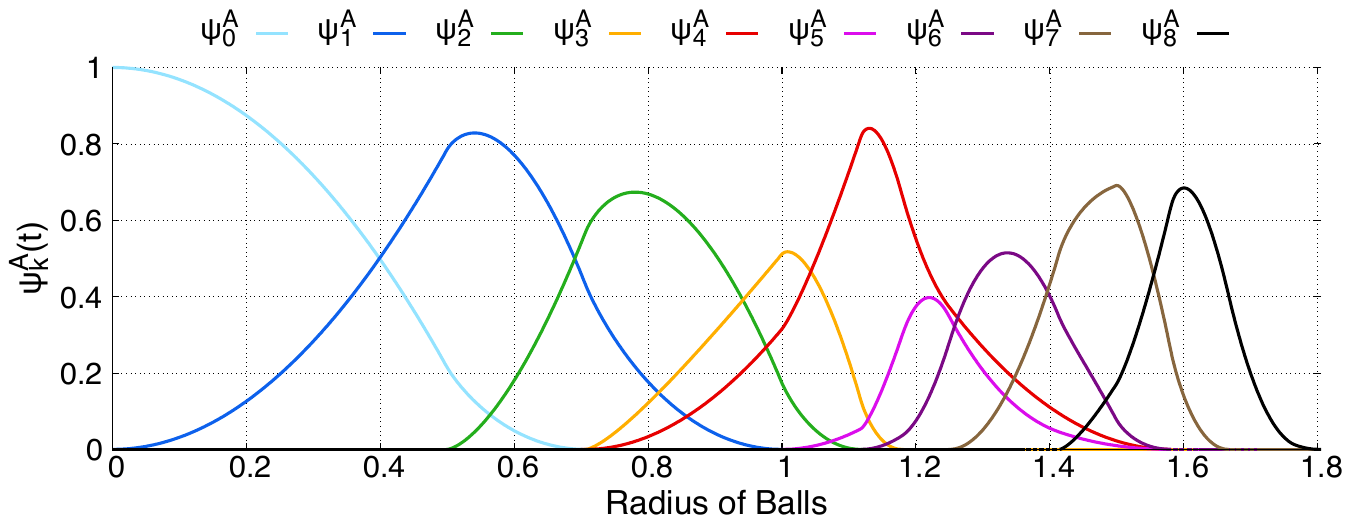}
  \includegraphics[width=\dwidth]{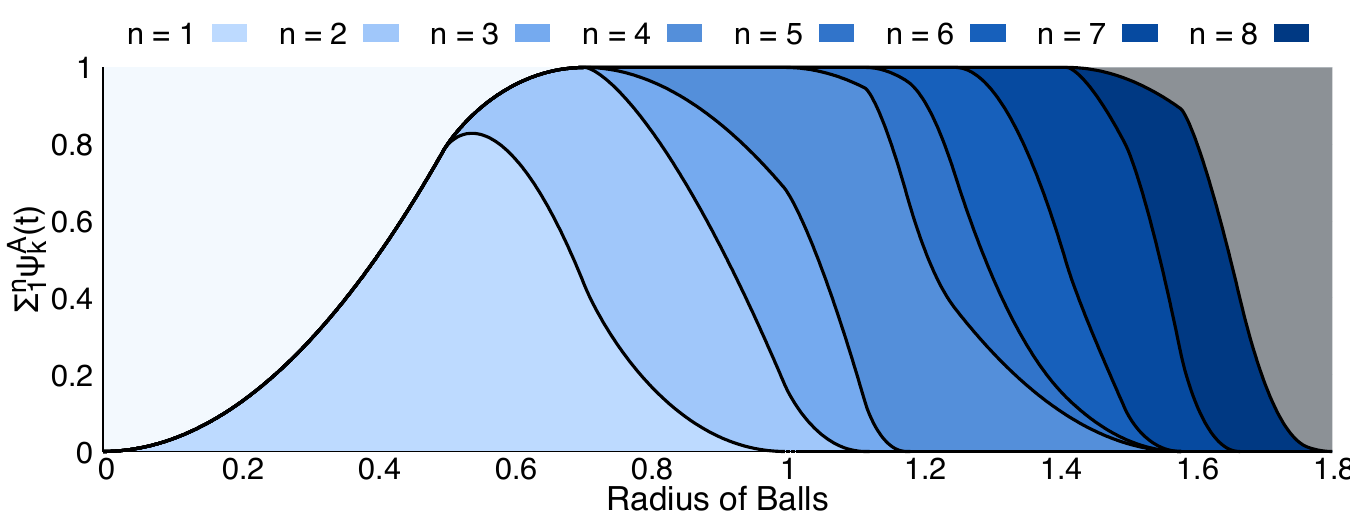}}
  \caption{Illustration of Definition~\ref{dfn:densities} for the square lattice $\La_4$.
  \textbf{Left}: subregions $U_k(t)$ are covered by $k$ disks for the radii $t=0.25, 0.55, 0.75, 1$.
  \textbf{Right}: the nine density functions are above the corresponding \emph{densigram} of 
  accumulated functions 
  $\sum\limits_{i=1}^k\psi_i(\La_4;t)$ \cite[Fig.~2]{edelsbrunner2021density}.}
  \label{fig:densities_sq}
\end{figure}

\begin{figure}[h!]
  \parbox{\bwidth}{
  \includegraphics[height=\cheight]{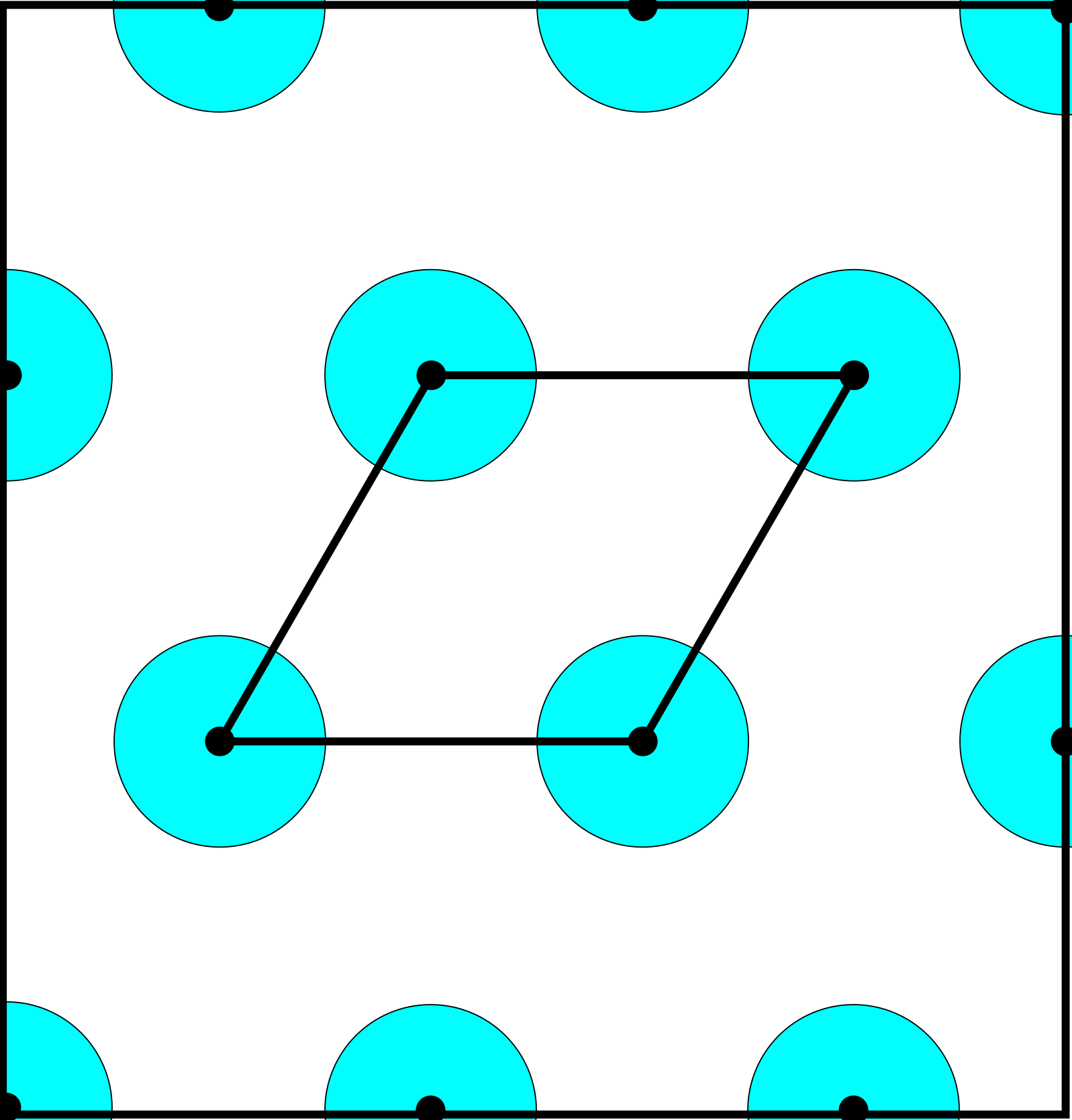}
  \hspace*{1mm}
  \includegraphics[height=\cheight]{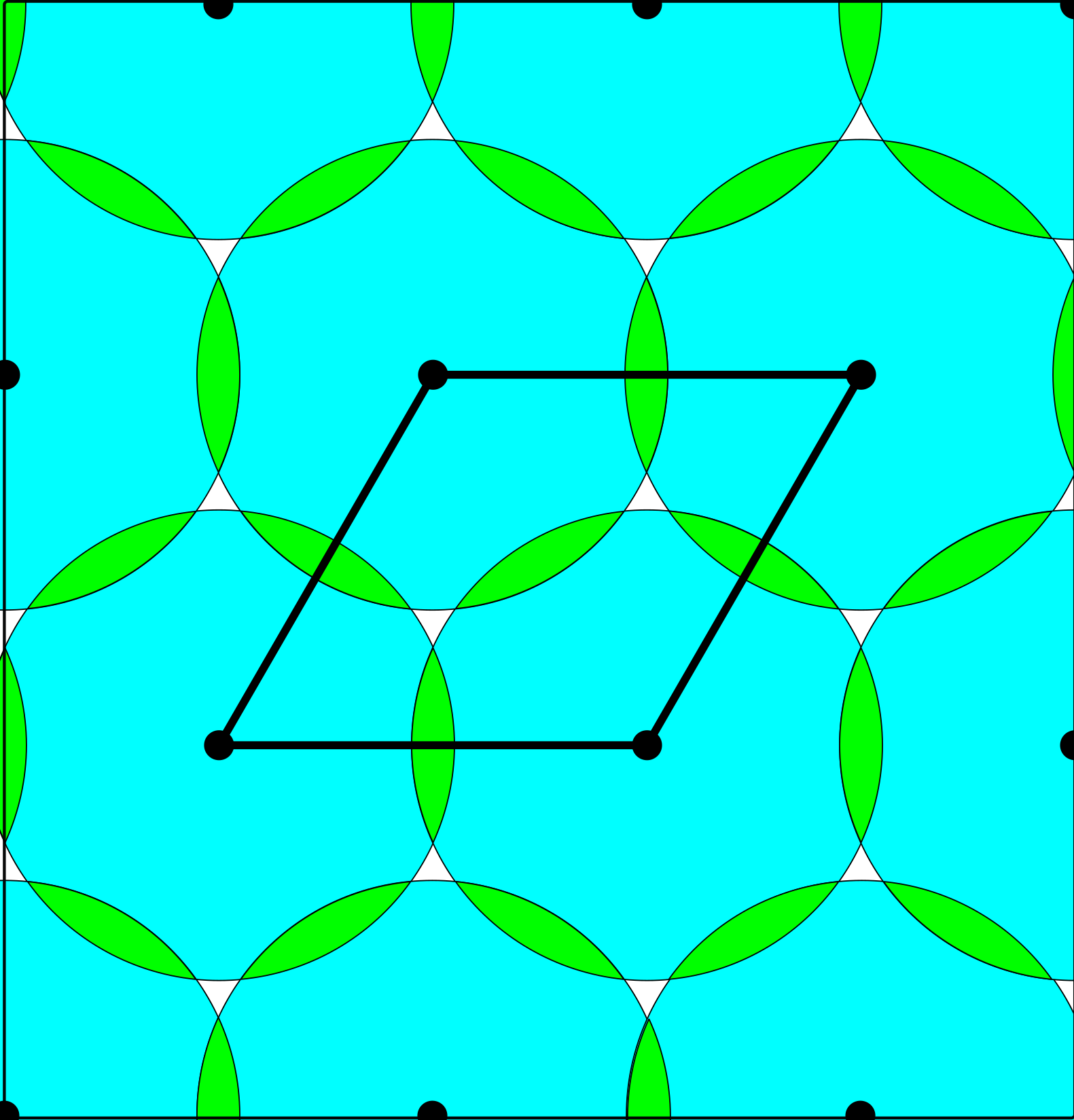}
  \medskip
  
  \includegraphics[height=\cheight]{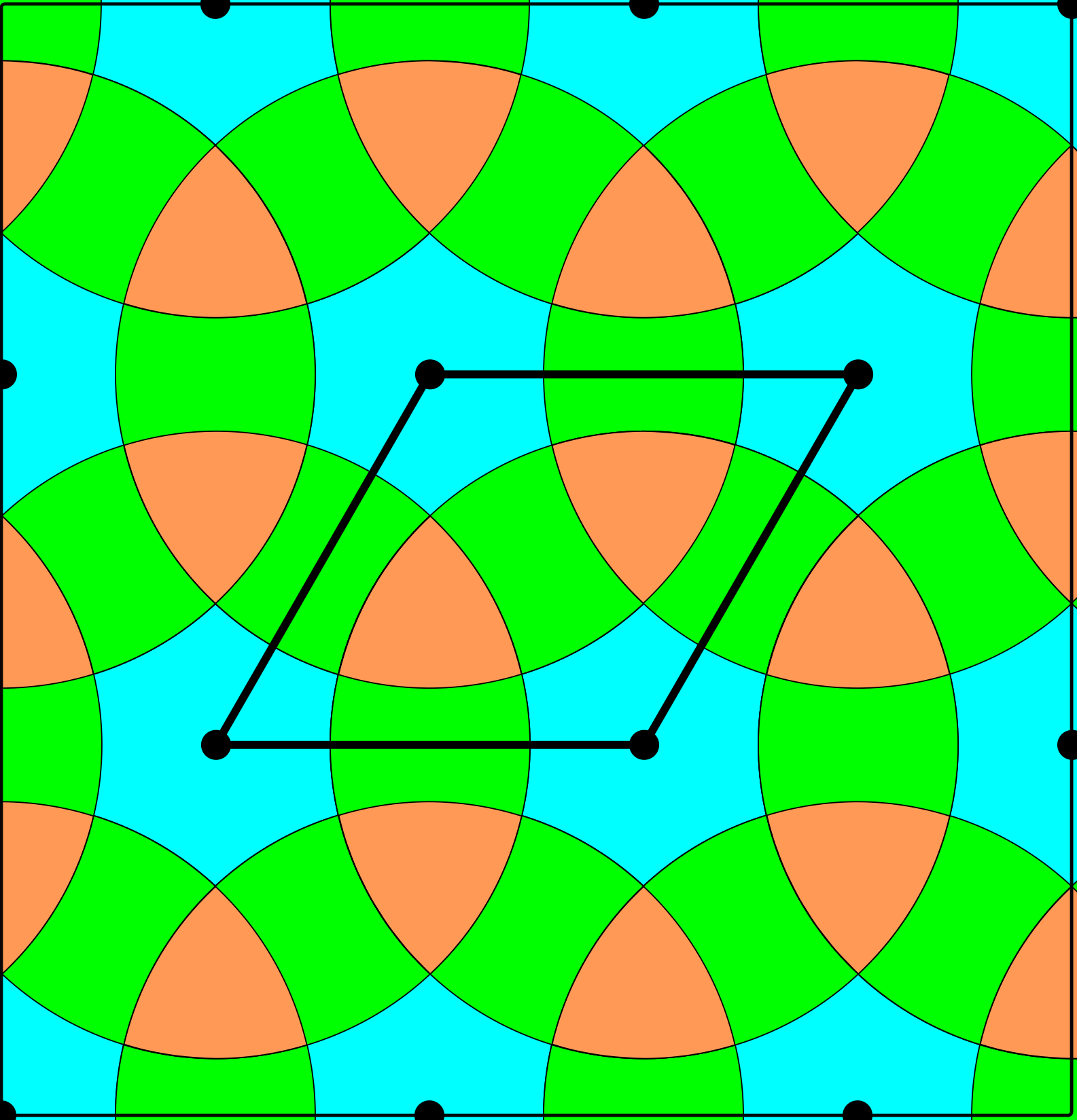}
  \hspace*{1mm}
  \includegraphics[height=\cheight]{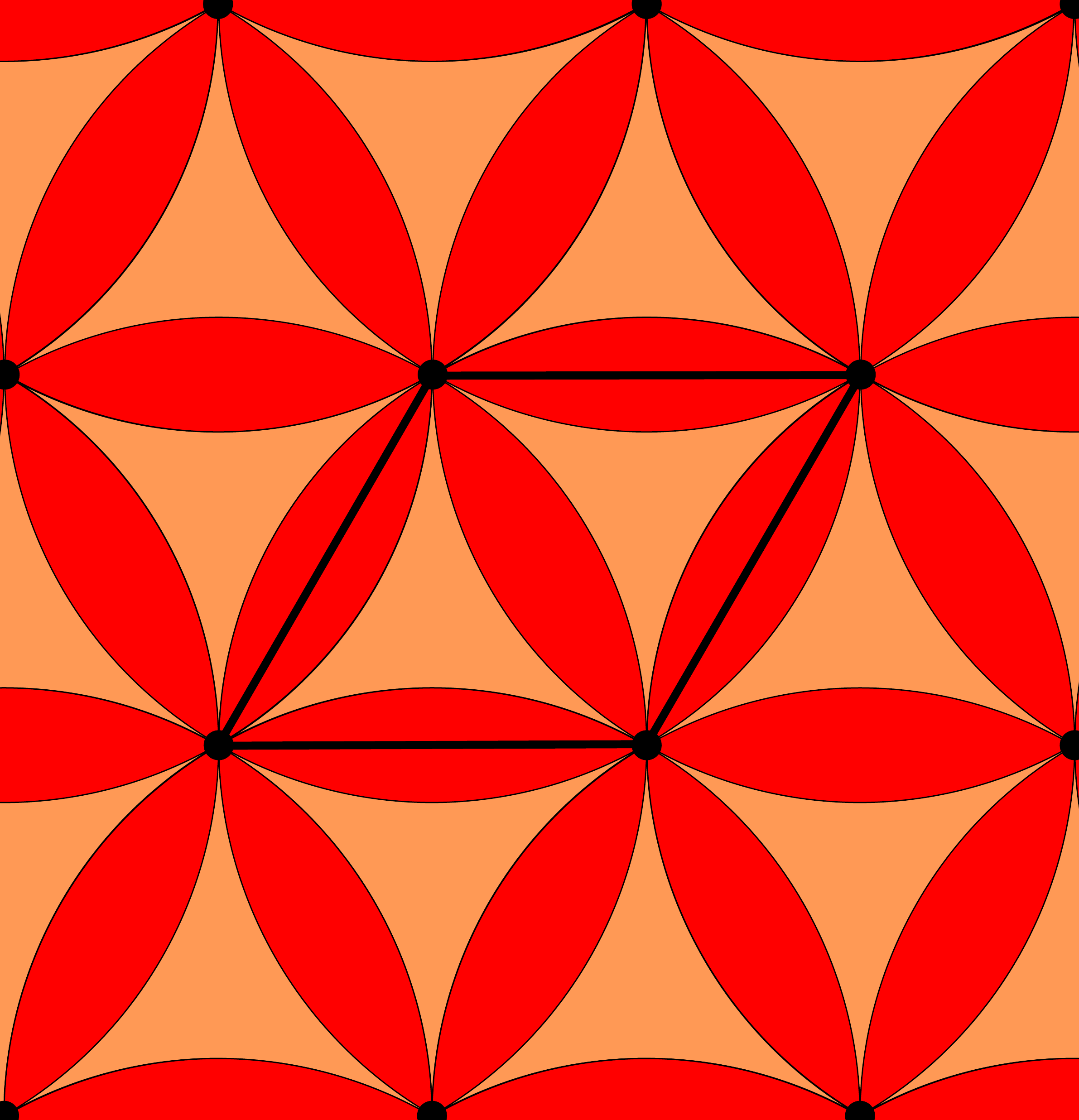}}
  \parbox{80mm}{
  \includegraphics[width=\dwidth]{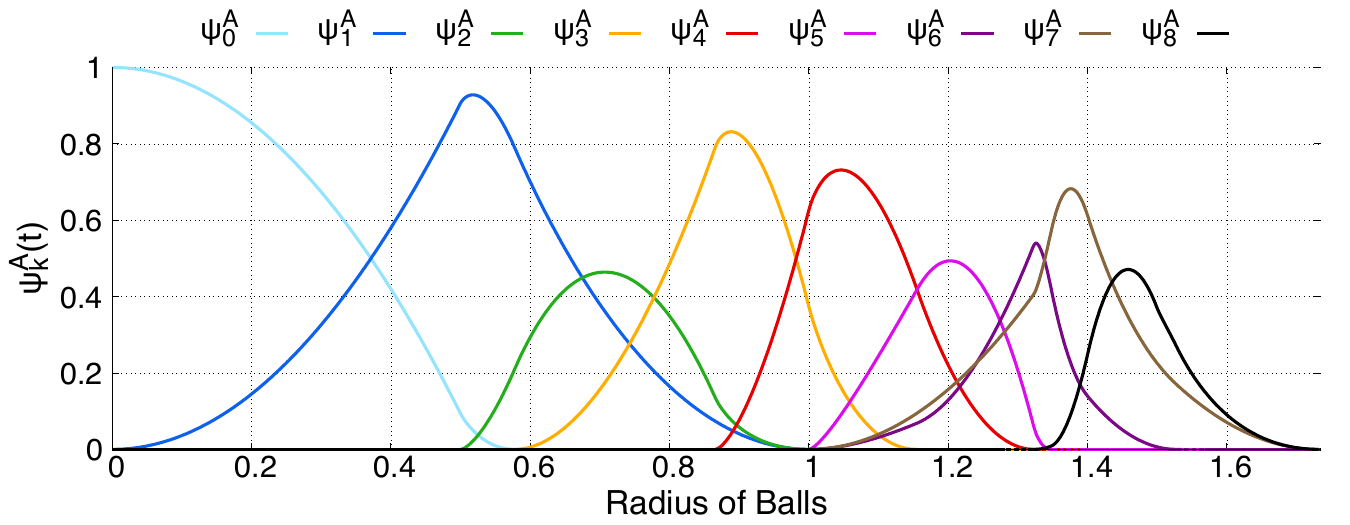}
  \includegraphics[width=\dwidth]{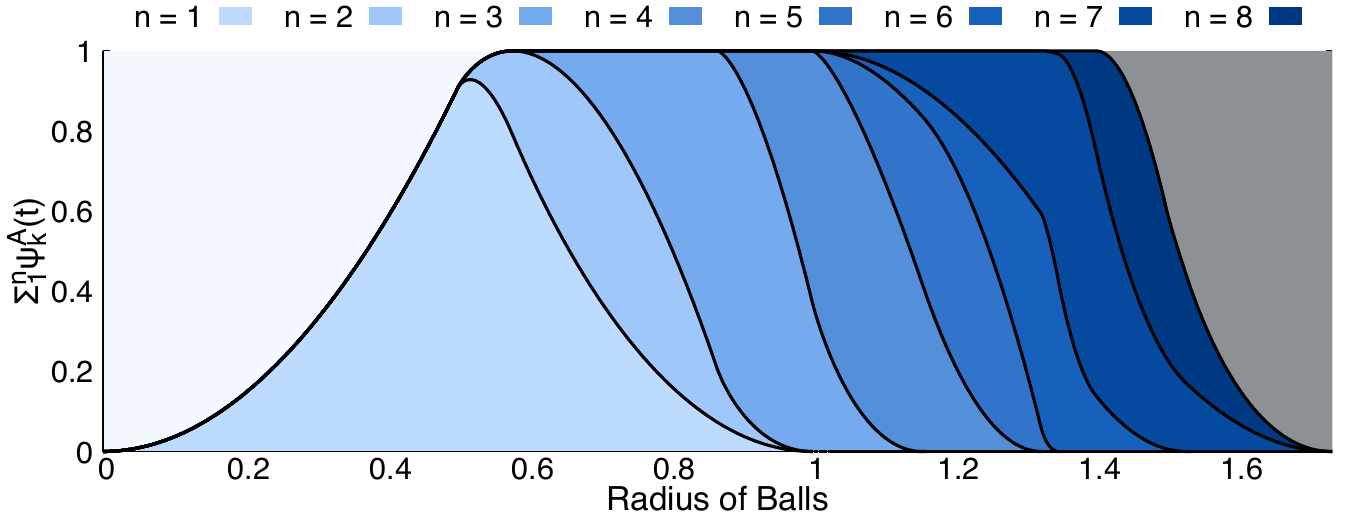}}
  \caption{illustration of Definition~\ref{dfn:densities} for the hexagonal lattice $\La_6$.
  \textbf{Left}: subregions $U_k(t)$ are covered by $k$ disks for the radii $t=0.25, 0.55, 0.75, 1$.
  \textbf{Right}: the nine density functions are above the corresponding \emph{densigram} of 
  accumulated functions 
  $\sum\limits_{i=1}^k\psi_i(\La_6;t)$ \cite[Fig.~2]{edelsbrunner2021density}.}
  \label{fig:densities_hex}
\end{figure}

Any density function $\psi_k[S](t)$ can also be interpreted as the probability that a random point $p\in\R^n$ is at a maximum distance $t$ to exactly $k$ points of $S$.
However, $\psi_k[S](t)$ is not a probability density function, so a potentially better name for $\psi_k[S](t)$ can be called the \emph{$k$-fold overlapping function}, but it is now a bit late to change this. 
\myskip

Since any isometry preserves distances and hence volumes of areas in $\R^n$, any density function is invariant under isometry and independent of a unit cell choice, see \cite[Lemma 1]{edelsbrunner2021density}.
Their further properties will be discussed in the next section.

\section{Continuity and generic completeness of density functions in $\R^3$}
\label{sec:densities_gen_comp}

This section presents the key results from paper \cite{edelsbrunner2021density}.
Definition~\ref{dfn:densities_metrics} introduces metrics on density functions to prove their Lipschitz continuity in Theorem~\ref{thm:densities_stable}.

\index{density function}

\begin{dfn}[metrics on density functions]
\label{dfn:densities_metrics}
For any $k\geq 1$ and periodic point sets $S,Q\subset\R^n$, 
define the \emph{max metrics} between their density functions as
$$|\psi_k[S]-\psi_k[Q]|_\infty=\sup\limits_{t\geq 0}|\psi_k[S](t)-\psi_k[Q](t)|$$ 
and between their fingerprints
$d_\infty(\Psi[S],\Psi[Q])=
\sup\limits_{k\geq 0} \dfrac{|\psi_k[S]-\psi_k[Q]|_\infty}{(\sqrt[3]{k+1})^2}$.
\edfn
\end{dfn}

In practice, the max distance $|\psi_k[S]-\psi_k[Q]|_\infty$ is approximated for uniformly sampled radii $t$, while $d_\infty(\Psi[S],\Psi[Q])$ is approximated by considering $k=1,\dots,8$.
\myskip

For any fixed index $k\geq 0$, the density function $\psi_k[S](t)$ eventually vanishes for large $t$ because the whole unit cell $U$ of $S$ includes only higher $l$-fold intersections of ball $B(p;t)$ for $l>k$.
The extra factor $\dfrac{1}{(\sqrt[3]{k+1})^2}$ reduces the apparently growing sensitivity of the density functions $\psi_k[S](t)$ to perturbations of points.

\index{packing radius}
\index{covering radius}

\begin{dfn}[radii $r(S)$ and $R(S)$]
\label{dfn:packing+covering}
Let $S\subset\R^n$ be any periodic point set.
\myskip

\nt
\tb{(a)}
The \emph{packing} radius $r(S)$ is the maximum $r$ such that the open balls $B(p;r)$ are disjoint for all $p\in S$, or $r(S)$ is the minimum half-distance between any points of $S$.
\myskip

\nt
\tb{(b)}
The \emph{covering} radius $R(S)$ is the minimum radius $R$ such that the union of closed balls $\bar B(p;R)$ for all $p\in S$ covers $\R^n$.
Alternatively, $R(S)$ is the maximum distance from any point $q\in\R^n$ to its nearest neighbour in $S$. 
\edfn
\end{dfn}

In terms of density functions, $r(S)$ is the maximum radius $t\geq 0$ such that $\psi_1[S](t)=0$, while $R(S)$ is the maximum radius $R$ such that $\psi_1[S](t)=0$ for all $t\in[0,R]$.

\begin{thm}[continuity of density functions, {\cite[Theorem~1]{edelsbrunner2021density}}]
\label{thm:densities_stable}
Let $Q\subset\R^3$ be a periodic point set obtained from another periodic point set $S\subset\R^3$ by perturbing any point of $S$ up to Euclidean distance $\ep$ such that $0\leq\ep<r=\min\{r(S),r(Q)\}$.
Set $R=\max\{R(S),R(Q)\}$ and $\la=13\dfrac{R^2}{r^3}$.
Then 
$d_\infty(\Psi[S],\Psi[Q])\leq\la\ep$.
\ethm
\end{thm}

\cite[section~5.1]{edelsbrunner2021density} describes technical conditions defining \emph{generic} periodic point sets $S\subset\R^3$, which satisfy Theorem~\ref{densities_gen_comp}.

\begin{thm}[generic completeness of density functions, {\cite[Theorem~1]{edelsbrunner2021density}}]
\label{densities_gen_comp}
If any generic periodic point sets $S,Q\subset\R^3$ are not isometric ($S\not\simeq Q$), then $\Psi[S]\neq\Psi[Q]$.
\ethm
\end{thm}

Now we describe the first practical impact of density functions.
Crystal Structure Prediction (CSP) aims to predict whether a selected molecule can be crystallised into a functional material, i.e. 
a crystal with useful functions or properties. In theory, CSP seeks to answer the question of whether copies of a molecule can be arranged in such a way that the resulting crystal is thermodynamically stable as well as useful.
Crucially, CSP aims to answer this question purely computationally to streamline the trial-and-error in molecular synthesise.

\begin{figure}[ht]
\vspace{0.1in}
\includegraphics[height=32mm]{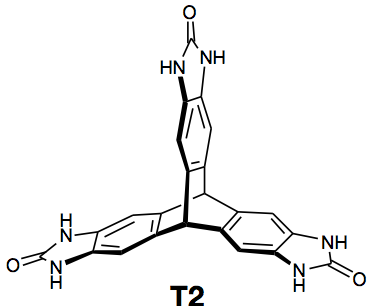}
\hspace{1mm}
\includegraphics[height=32mm]{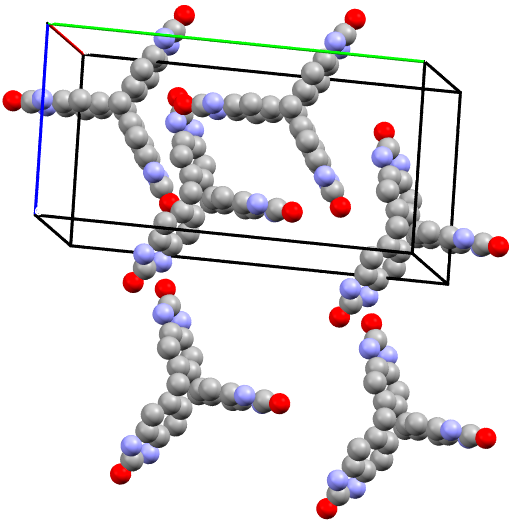} 
\hspace{1mm}
\includegraphics[height=32mm]{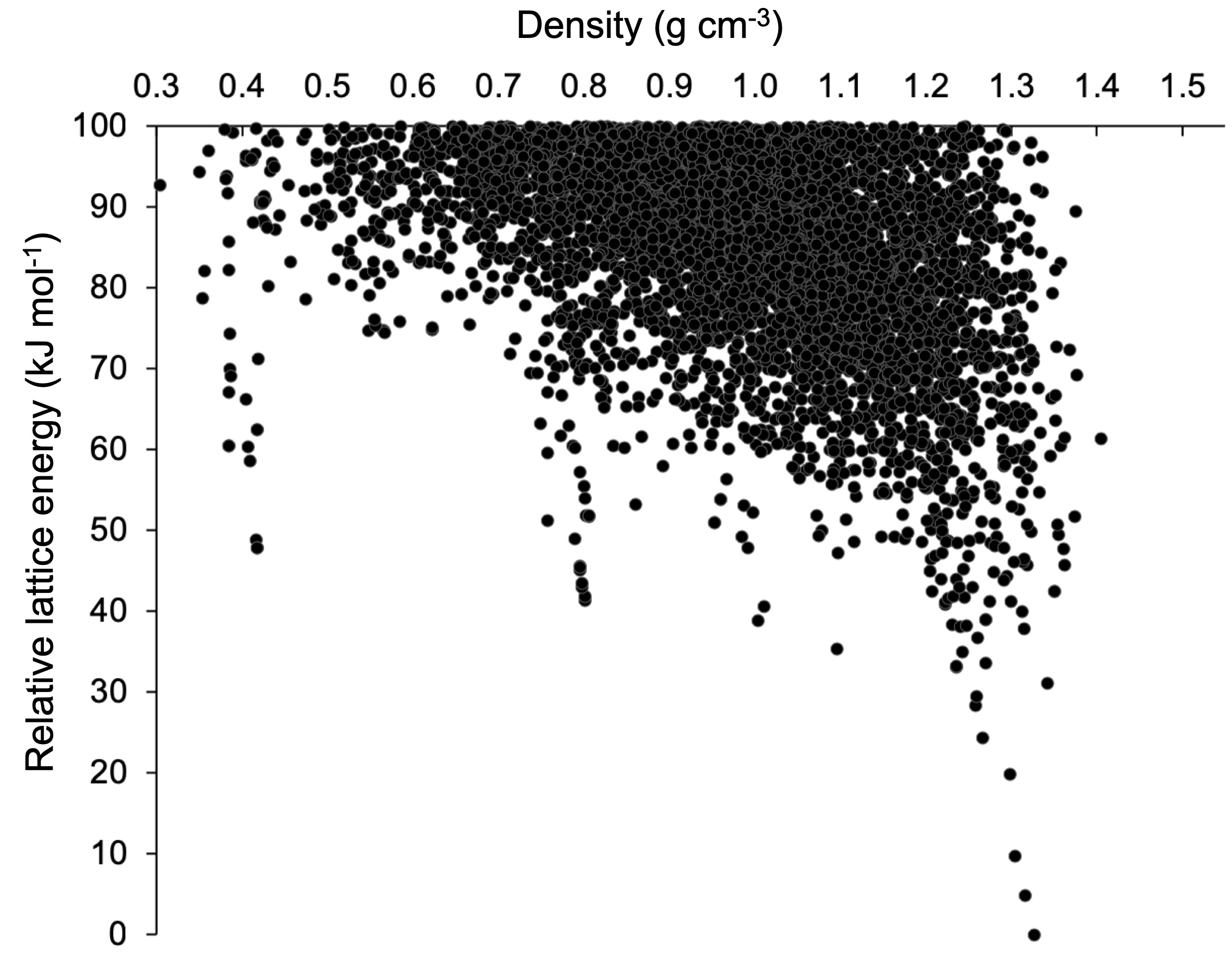}
\caption{\textbf{Left}: a T2 molecule. 
  \textbf{Middle}: the T2-$\delta$ crystal with highlighted unit cell. 
  \textbf{Right}: the output of CSP for the T2 molecule.
  It is a plot of 5679 simulated T2 crystal structures \cite[Fig.~2d]{pulido2017functional}, each represented by two coordinates: the physical density (atomic mass within a unit cell divided by the unit cell volume) and energy (determining the crystal's thermodynamic stability). 
  Structures at the bottom of the `downward spikes' are likely to be stable.}
\label{fig:T2molecule}
\end{figure}

Our colleagues at Liverpool's Materials Innovation Factory \cite{pulido2017functional} computationally predicted that the T2 molecule in Fig.~\ref{fig:T2molecule} can be crystallised into a new nanoporous material for gas storage. 
As part of this process, they also identified four other structures of interest.
Following the CSP predictions, they synthesised $5$ families of T2-crystals in the laboratory by varying parameters like temperature and pressure, calling them T2-$\alpha$, T2-$\beta$, $\dots$, T2-$\epsilon$. One of them, T2-$\gamma$, indeed had the desired property of having only half the physical density of the previously known structure T2-$\alpha$.
They scanned the synthesised crystals using X-ray powder diffraction yielding Crystallographic Information Files, 
each containing the unit cell and the motif points representing the atoms. 
These files were then compared with the results of the simulations,  
either by using their physical density alongside the \textsc{Compack} algorithm---which compares only a finite portion of the structure---or by looking at visualisations of the crystal structures.
This comparison showed that the synthesised crystals matched the prediction well.
Our colleagues deposited these structures into the Cambridge Structural Database (CSD).
\medskip

At a later time, we used our density functions to verify our collaborators' matchings between the synthesised crystals T2-$\alpha$ to T2-$\epsilon$ and the simulated crystals entry 99, 28, 62, 09, 01.
We did so by computing, for each of the five matches, the distance between the density functions of the synthesised and the simulated crystal.
As one is the prediction of the other, we expected to see small distances.
And for four of the five structures this was true: T2-$\gamma$, for example, always has an $L_{\infty}$-distance of less than $0.04$ over the first eight pairs of corresponding density functions; see Table~\ref{tab:T2distances}.
However, when we came to check the distances between density functions of T2-$\delta$ with its predicted structure, we were surprised to see large distances (the final row of Table~\ref{tab:T2distances}).
It turned out that a mix-up of files had happened, and what was uploaded to the Cambridge Structural Database as T2-$\delta$ was in fact T2-$\beta'$ (a crystal from the T2-$\beta$ family). 
The density fingerprint revealed this error, which was verified by chemists upon a visual inspection, and it is 
because of this that T2-$\delta$ was subsequently correctly deposited. 

\begin{table}[ht]
  \begin{center}
  \begin{tabular}{c||cccccccc}
    $|\psi_{k}[S]-\psi_{k}[Q]|_\infty$ & $k=0$ & 1 & 2 & 3 & 4 & 5 & 6 & 7 \\
    \hline \hline
    T2-$\alpha$ vs entry 99   & 0.0042 & 0.0092 & 0.0125 & 0.0056 & 0.0099 & 0.0088 & 0.0127  & 0.0099 \\
    T2-$\beta$ vs entry 28    & 0.0157 & 0.0156 & 0.0159 & 0.0224 & 0.0334 & 0.0396 & 0.0357  & 0.0454 \\
    T2-$\gamma$ vs entry 62   & 0.0020 & 0.0080 & 0.0128 & 0.0155 & 0.0153 & 0.0250 & 0.0296  & 0.0391 \\
    T2-$\delta$ vs entry 09   & 0.0610 & 0.0884 & 0.1267 & 0.0676 & 0.0915 & 0.0801 & 0.0733  & 0.0388 \\
    T2-$\epsilon$ vs entry 01 & 0.0132 & 0.0152 & 0.0207 & 0.0571 & 0.0514 & 0.0431 & 0.0468  & 0.0550 \\ \hline
    T2-$\beta'$ vs entry 09   & 0.2981 & 0.2631 & 0.3718 & 0.3747 & 0.2563 & 0.2360 & 0.3161  & 0.3232
 \end{tabular}
 \end{center}
  \caption{\textit{First five rows:} the $L_\infty$-distances between the first eight pairs of corresponding density functions of physically synthesised T2 crystals (T2-$\alpha$, T2-$\beta$, etc.) and the simulated structures that had predicted them from the CSP output dataset (entry XX).
  \textit{Last row:} the suspiciously larger numbers revealed the mix-up of the files T2-$\delta$ and T2-$\beta'$ and thus led to depositing the initially omitted Crystallographic Information File of the T2-$\delta$ crystal into the Cambridge Structural Database. 
  }
  \label{tab:T2distances}
\end{table}

Plots of the density functions of correctly matched synthesised and simulated structures can be seen in Figure~\ref{fig:T2densities}. 
As another application, we expect that the fingerprint will 
be used 
to simplify the large output data sets produced by CSP by comparing simulated structures with each other, thus speeding up what is currently a slow process.

\newcommand{\cdheight}{30mm}
\begin{figure}[!ht]
  \centering
\includegraphics[height=\cdheight]{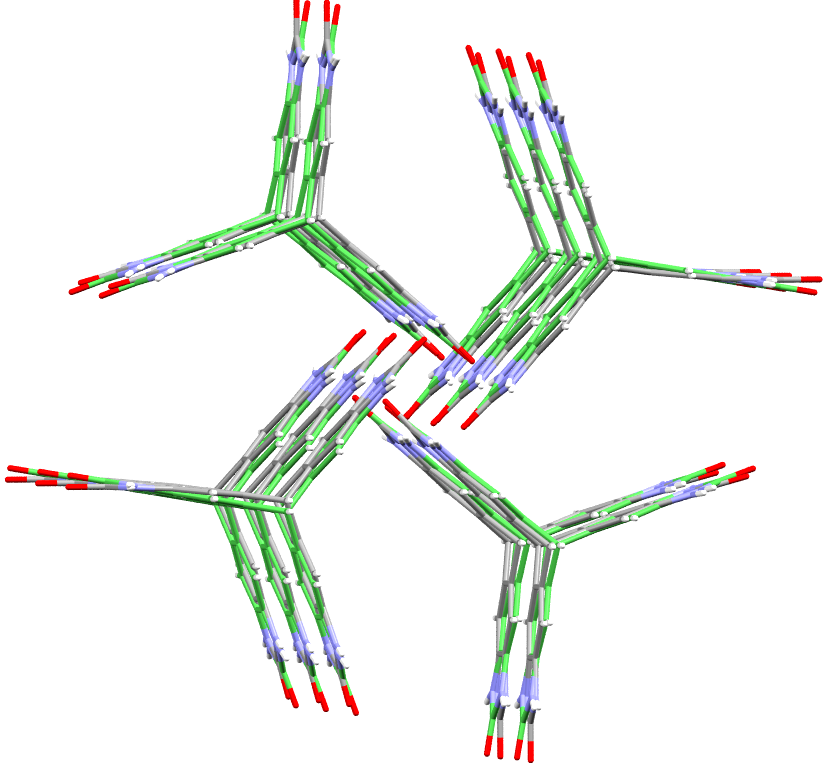} \hfill
  \includegraphics[height=\cdheight]{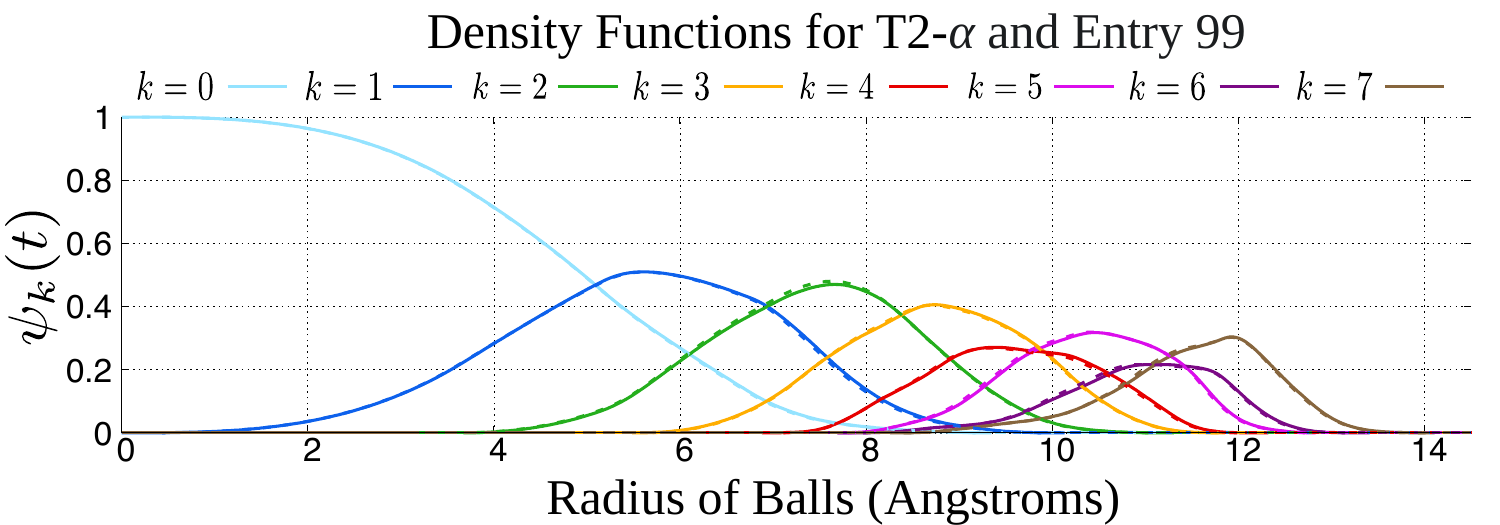} 
\smallskip
  
\includegraphics[height=\cdheight]{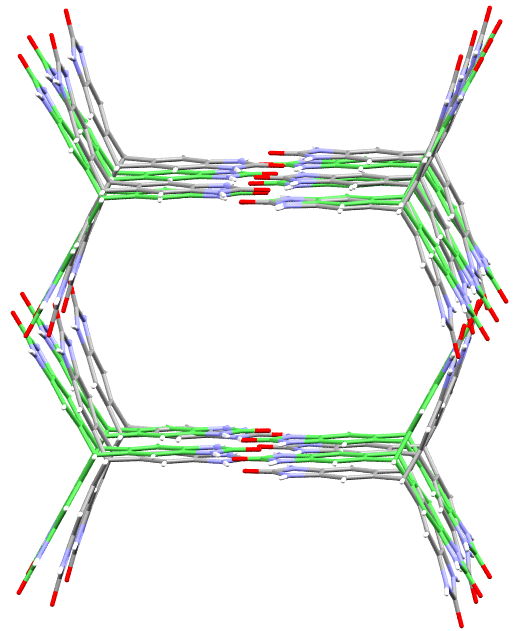} \hfill
  \includegraphics[height=\cdheight]{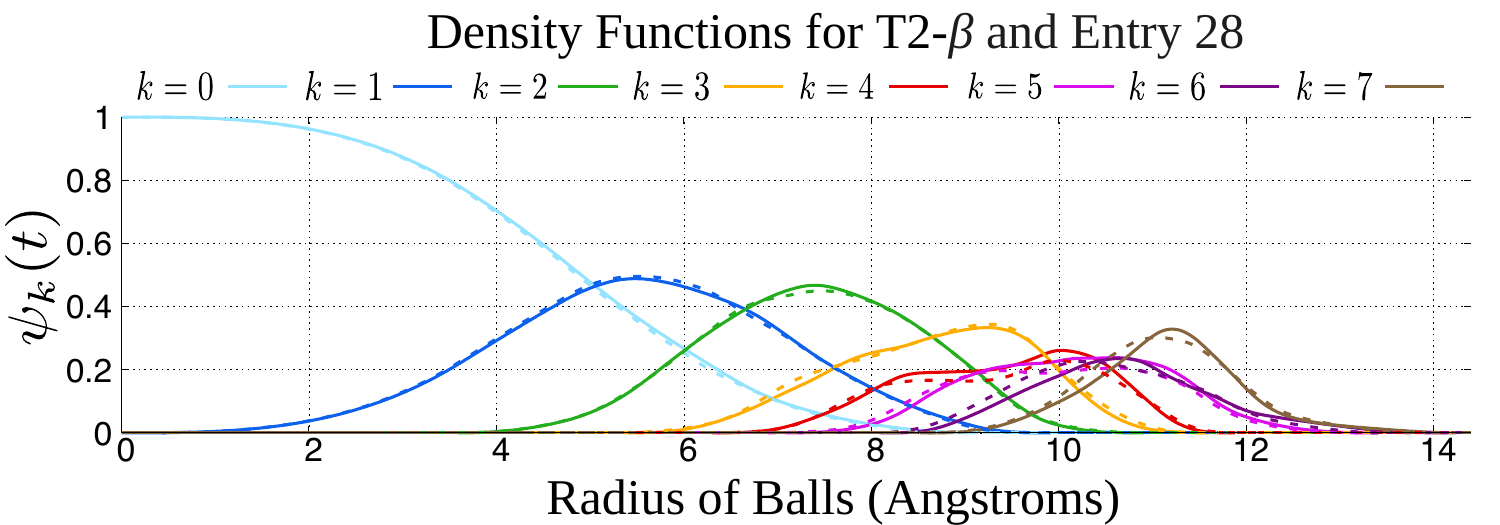} 
\medskip

\includegraphics[height=\cdheight]{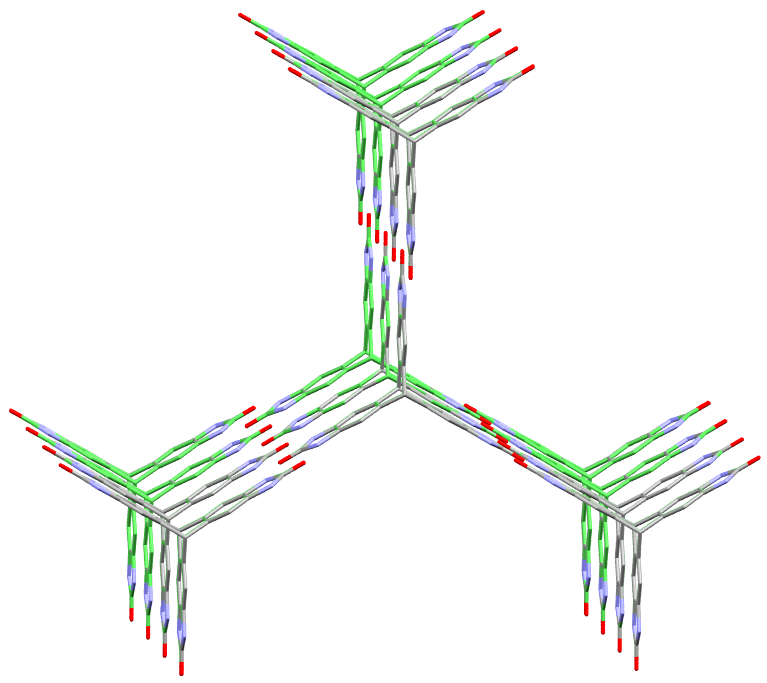} \hfill
  \includegraphics[height=\cdheight]{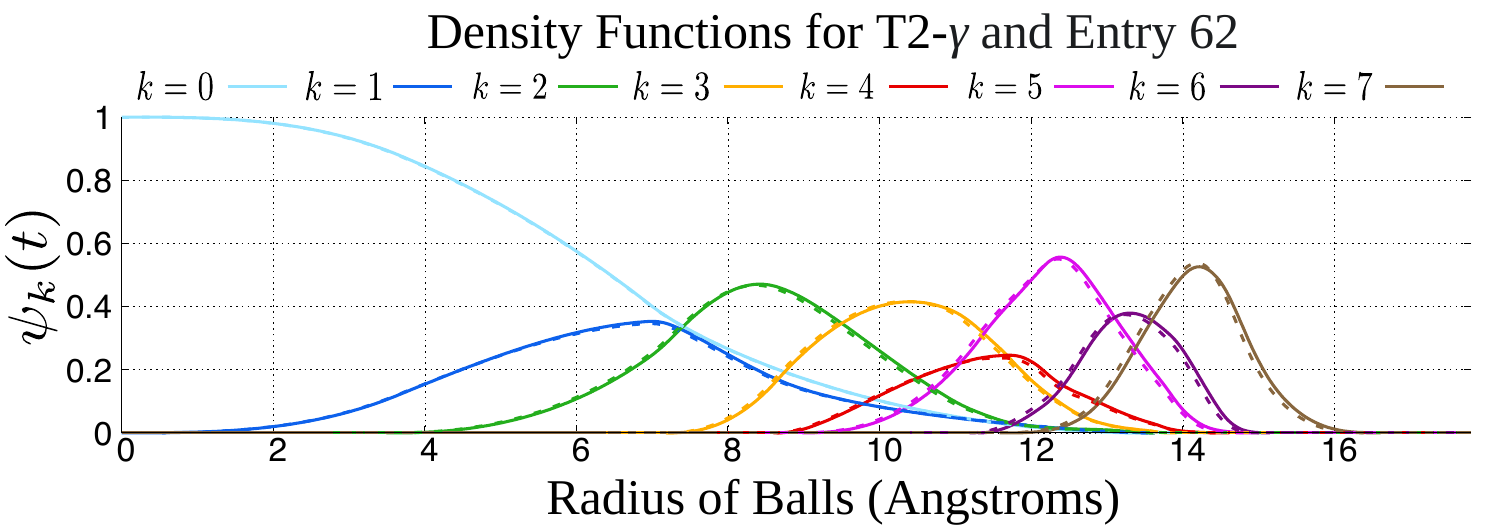} 
\medskip
  
\includegraphics[height=\cdheight]{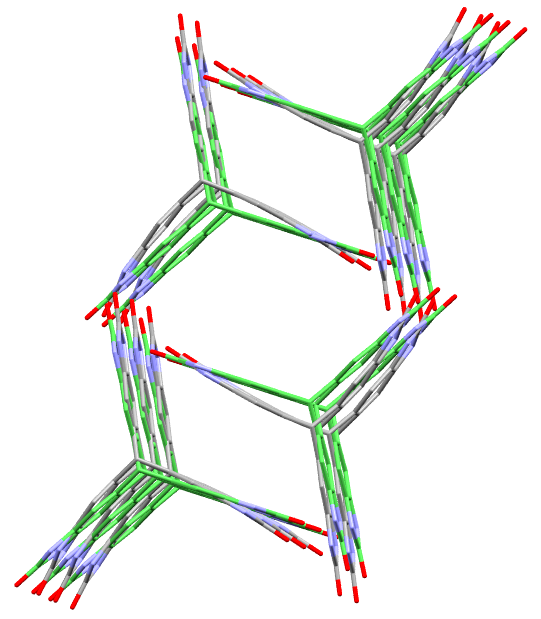} \hfill
  \includegraphics[height=\cdheight]{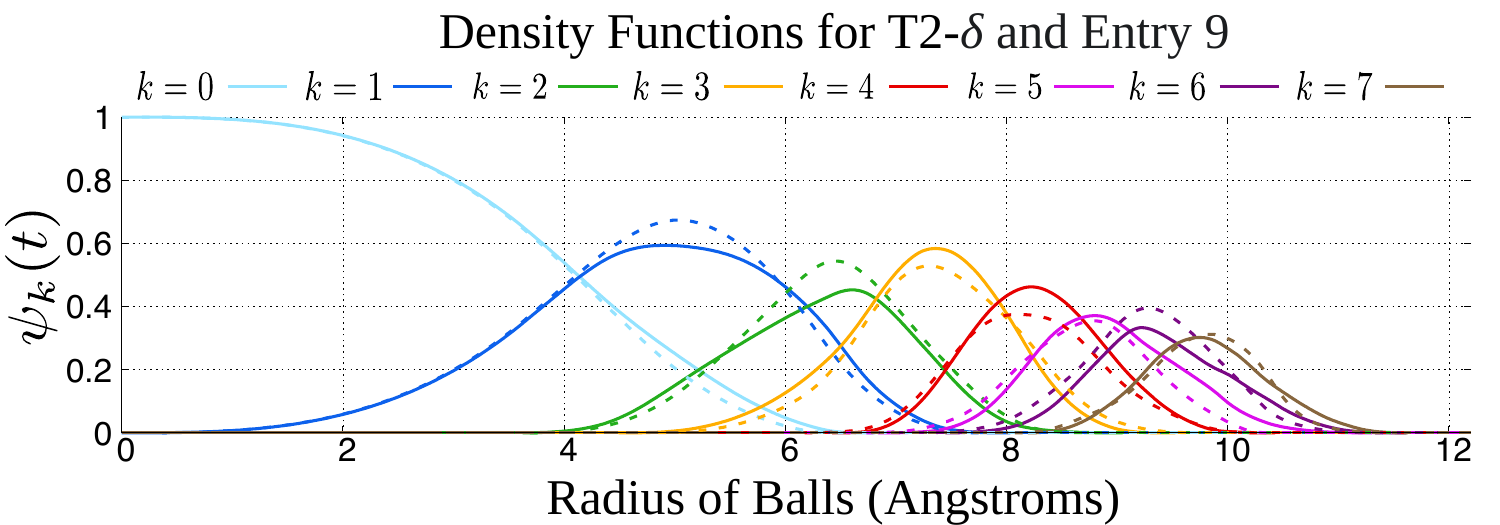}
\medskip
  
\includegraphics[height=\cdheight]{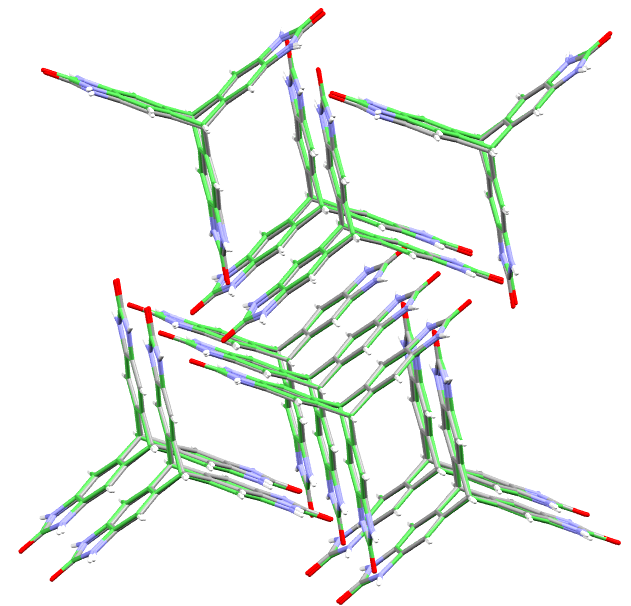} \hfill
  \includegraphics[height=\cdheight]{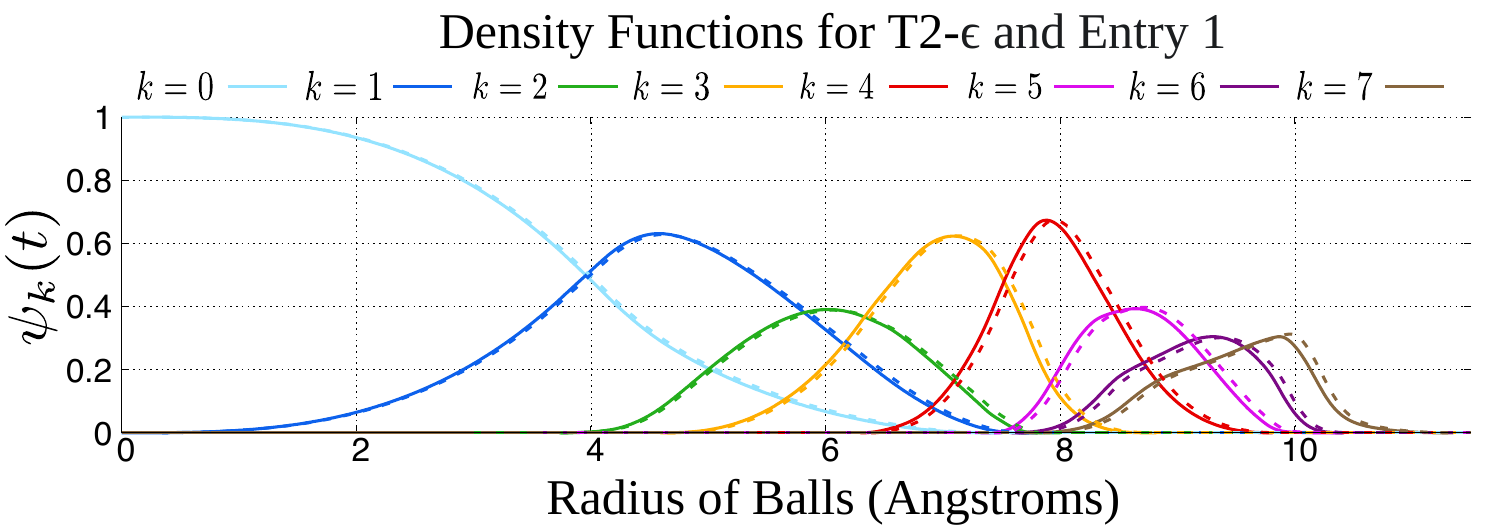}
\caption{
  \emph{Left}: experimental T2 crystals (curved gray molecules) and their simulated versions (straight green molecules) overlaid. \emph{Right}: the density functions of the periodic sets of molecular centres of the experimental T2 crystals (solid curves) vs.\ simulated crystals (dashed curves).}
\label{fig:T2densities}
\end{figure}

\section{Density functions via degree-$k$ Voronoi domains}
\label{sec:Brillouin}

This section follows papers \cite{edelsbrunner2021density,smith2022practical} with minor updates.
To compute any density function $\psi_k[S](t)$, we introduce degree-$k$ Voronoi domains, also called the $k$-th Brillouin zone.
These domains generalise the Voronoi domain $\bar V(\La)$ of a lattice from Definition~\ref{dfn:Voronoi_domain}(a) and differ from the \emph{order-$k$ Voronoi} domain \cite{edelsbrunner2023simple} defined for a $k$-point subset $A\subset S\subset\R^n$ and as the set of all points for whom the points in $A$ are the closest $k$ points in $S$.

\begin{dfn}[index-$k$ Voronoi domain $V_k(S;p)$ and degree-$k$ Voronoi domain $Z_k(S; p)$]
\label{def:Voronoi_domains_ind+deg}
\tb{(a)}
For a finite or a periodic point set $S\subset\R^n$ and a point $p\in S$, the \emph{$k$-th Voronoi domain} $V_k(C; p)$ is the closure of the set of all points $q\in\R^n$ such that $p$ is among the $k$ nearest points of $S$ to $q$. 
\myskip

\nt
\tb{(b)}
For any periodic point set $S\subset\R^n$ and $p\in S$, the \emph{degree-$k$ Voronoi domain} is the difference between successive closed index-$k$ Voronoi domains, i.e. 
$Z_k(C; p)=V_k(C; p) - V_{k - 1}(C; p)$ for $k\geq 1$, where we set $V_0(C; p) = \emptyset$.
\edfn
\end{dfn}

If $k=1$, $V_1(S;p)=Z_1(S;p)$ is the classical Voronoi domain for a point $p\in S$.
\myskip

The index-$k$ Voronoi domain $V_k(C;p)\subset\R^n$ is defined as a closed set above to cover all cases where $p$ has equal distances to several neighbours, so a $k$-th neighbour of $p$ may not be unique. 
Unlike order-$k$ Voronoi domains, which tile $\mathbb{R}^n$ \cite{edelsbrunner1986voronoi}, index-$k$ Voronoi domains form a nested sequence.
Any $V_k(C; p)$ is \emph{star-convex}, which means it contains all line segments connecting $\bd V_k(C; p)$ to $p$.
Indeed, if $p\in C$ is among the $k$ nearest to $q\in \bd V_k(C; p)$, then any intermediate point in the line segment $[p,q]$ has $p$ among its $k$ nearest neighbours of $C$.
\myskip

\begin{figure}[h!]
\centering
\includegraphics[width=\textwidth]{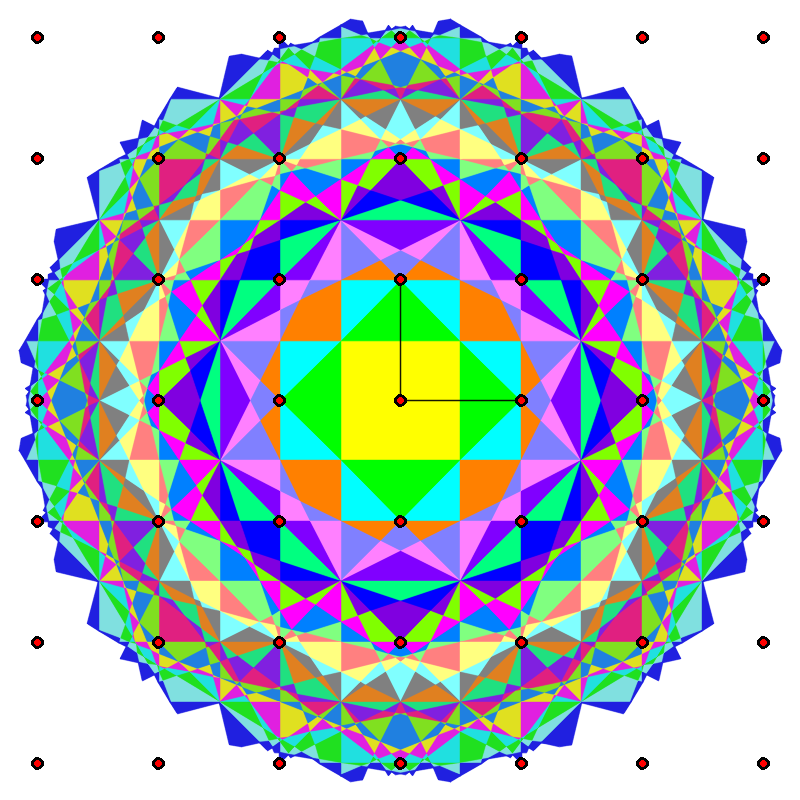}
\caption{
Degree-$k$ Voronoi domains $Z_k(\La_4;0)$ of the square lattice $\La_4$ of red points for $k=1,\dots,30$.}
\label{fig:Brillouin_30_Square}
\end{figure}

\begin{figure}[h!]
\centering
\includegraphics[width=\textwidth]{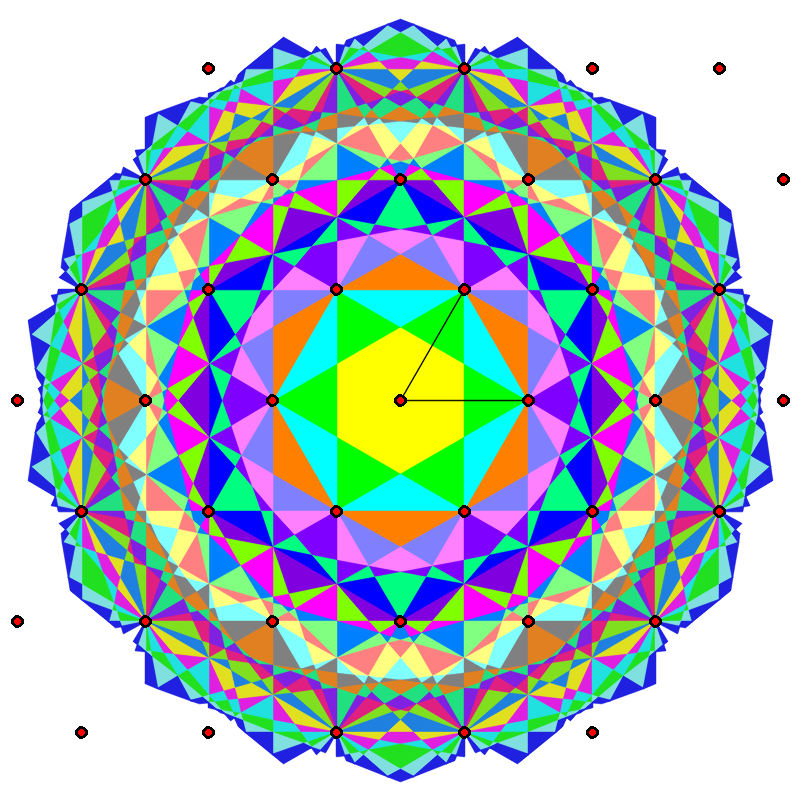}
\caption{
Degree-$k$ Voronoi domains $Z_k(\La_6;0)$ of the hexagonal lattice $\La_6$ for $k=1,\dots,30$.}
\label{fig:Brillouin_30_Hexagonal}
\end{figure}

Fig.~\ref{fig:Brillouin_30_Square} and~\ref{fig:Brillouin_30_Hexagonal} show degree-$k$ Voronoi domains for the square and hexagonal lattices, where any degree-$k$ Voronoi domain is the union of polygons of the same colour, and has the origin as its $k$-th nearest neighbour among all lattice points.
\myskip

Theorem~\ref{thm:Voronoi_volume} implies that Fig.~\ref{fig:Brillouin_30_Square} and~\ref{fig:Brillouin_30_Hexagonal} have the same total area of each colour.

\begin{thm}[volumes of degree-$k$ Voronoi domains, {\cite[Theorem 7]{smith2022practical}}]
\label{thm:Voronoi_volume}
For any periodic point set $S\subset\R^n$ with a motif $M$, the sum of the volumes of the degree-$k$ Voronoi domains $Z_k(S;p)$ over all motif points $p\in M$ is independent of $k$.
\ethm
\end{thm}


Theorem~\ref{thm:density_from_Voronoi} computes any density function $\psi_k[S]$ via  degree-$k$ Voronoi domains of $S$.
Let $\vol[C]$ denote the volume of any compact region $C\subset\R^n$.

\begin{thm}[formula for density functions, {\cite[Theorem~2]{edelsbrunner2021density}}]
\label{thm:density_from_Voronoi}
Let $S=\La+M\subset\R^n$ be a periodic point set with a motif $M\subseteq U$ in a unit cell $U$ of a lattice $\La\subset\R^n$.
Fix any integer $k\geq 1$.
Then $\psi_k[S](t)=\dfrac{1}{\vol[U]} \sum_{p\in M} \vol[V_k(S;p)\cap \bar B(p;t)]$.
\ethm
\end{thm}

Lemma~\ref{lem:Minkowski-reduced} can be considered a definition of a Minkowski-reduced basis, which is the last auxiliary concept needed to estimate the time for computing degree-$k$ Voronoi domains in Theorem~\ref{thm:deg-k_Voronoi_complexity}.

\begin{lem}[\emph{Minkowski-reduced} basis,  {\cite[Lemma 2.2.1]{nguyen2009low}}]
\label{lem:Minkowski-reduced}
A basis $\vec v_1,\dots,\vec v_n$ of a lattice $\La\subset\R^n$ is  \emph{Minkowski-reduced} if and only if, for any $i=1,\dots,n$ and integers $c_1,\dots,c_n\in\Z$ such that $c_i,\dots,c_n$ have no common integer factor $c>1$, the inequality $|\sum\limits_{i=j}^n c_j\vec v_j|\geq|\vec v_j|$ holds.  
\elem
\end{lem}

\begin{thm}[time of degree-$k$ Voronoi domains, {\cite[Theorem 7]{smith2022practical}}]
\label{thm:deg-k_Voronoi_complexity}
For $n=2,3$, let a periodic point set $S=\La+M\subset \R^{n}$ have a motif $M$ of $m$ points and a lattice $\La$ with a Minkowski-reduced basis.
For any point $p\in S$, the time to compute all degree-$i$ Voronoi domains $Z_i(S;p)$ for $i=1,\dots,k$ is $O(m^n(4k)^{n^2}(n\log(4k)+\log m))$.
\ethm
\end{thm}

\section{A description of density functions of periodic sequences}
\label{sec:densities1D_points}

All quoted results in this section have detailed proofs in \cite{anosova2022density}.
\myskip

\begin{figure}[ht]
\includegraphics[width=\textwidth]{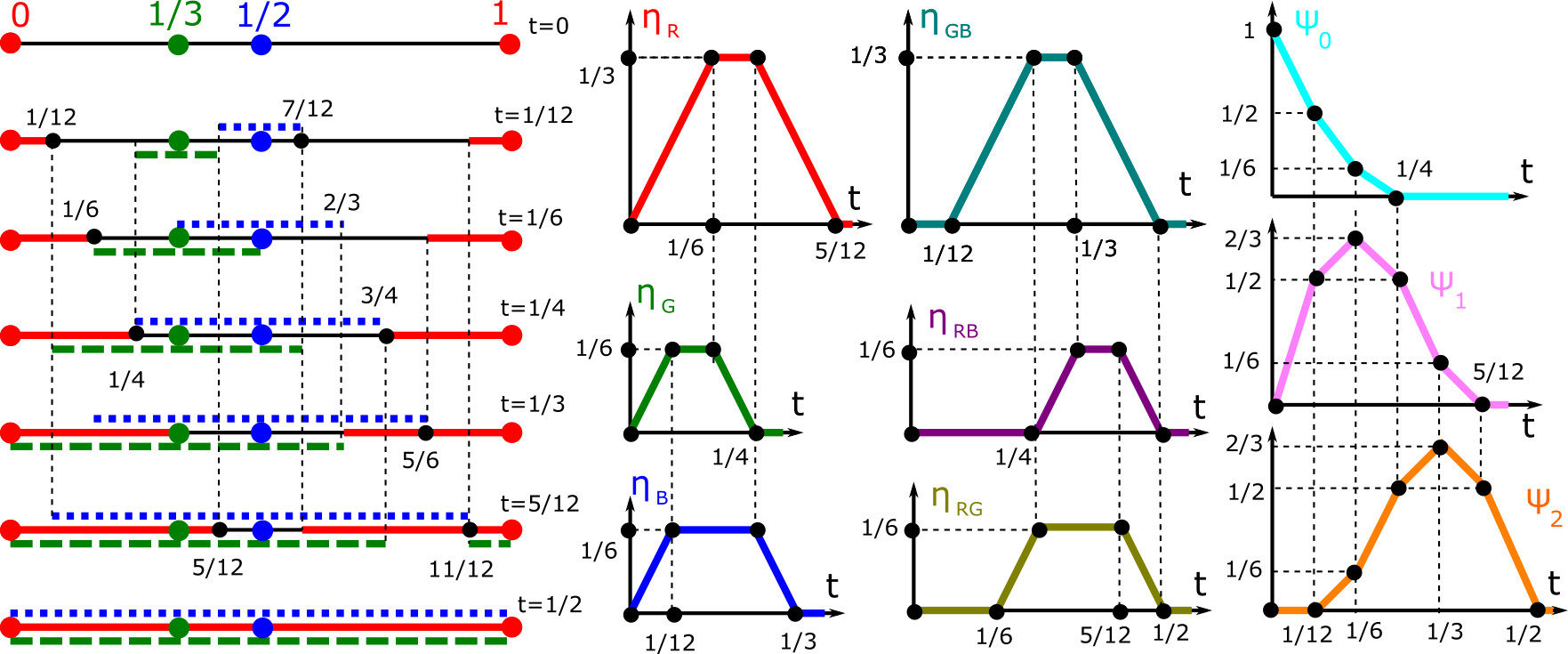}
\caption{\textbf{Left}: the periodic sequence $S=\{0,\frac{1}{3},\frac{1}{2}\}+\Z$ with points of three colors. 
The growing intervals around the red point $0\equiv 1\pmod{1}$, green point $\frac{1}{3}$, blue point $\frac{1}{2}$ have the same color for various radii $t$.
\textbf{Right}: the trapezium functions $\eta$ from Example~\ref{exa:densities}.
}
\label{fig:densities1D}      
\end{figure}

For convenience, scale any periodic sequence to period 1 so that $S=\{p_1,\dots,p_m\}+\Z$.
Since the expanding balls in $\R$ are growing intervals, volumes of their intersections linearly change in the variable radius $t$.
Hence, any density function $\psi_k(t)$ is piecewise linear and uniquely determined by \emph{corner} points $(a_j,b_j)$ where the gradient changes.
Examples~\ref{exa:0th_density} and~\ref{exa:densities} explain
how the density functions $\psi_k(t)$ are computed for the periodic sequence $S=\{0,\frac{1}{3},\frac{1}{2}\}+\Z$, see all graphs in Fig.~\ref{fig:densities1D}.

\begin{exa}[$0$-th density $\psi_0(t)$ for $S=\{0,\frac{1}{3},\frac{1}{2}\}+\Z$]
\label{exa:0th_density}
By Definition~\ref{dfn:densities} $\psi_0(t)$ is the fractional length within the period interval $[0,1]$ not covered by the intervals of radius $t$ (length $2t$), which are the red intervals $[0,t]\cup[1-t,1]$, green dashed interval $[\frac{1}{3}-t,\frac{1}{3}+t]$ and blue dotted interval $[\frac{1}{2}-t,\frac{1}{2}+t]$.
The graph of $\psi_0(t)$ starts from the point $(0,1)$ at $t=0$.
Then $\psi_0(t)$ linearly drops to the point $(\frac{1}{12},\frac{1}{2})$ at $t=\frac{1}{12}$ when a half of the interval $[0,1]$ remains uncovered.
\medskip

The next linear piece of $\psi_0(t)$ continues to the point $(\frac{1}{6},\frac{1}{6})$ at $t=\frac{1}{6}$ when only $[\frac{2}{3},\frac{5}{6}]$ is uncovered.
The graph of $\psi_0(t)$ finally returns to the $t$-axis at the point $(\frac{1}{4},0)$ and remains there for $t\geq \frac{1}{4}$.
The piecewise linear behaviour of $\psi_0(t)$ can be briefly described via the \emph{corner} points $(0,1)$, $(\frac{1}{12},\frac{1}{3})$, $(\frac{1}{6},\frac{1}{6})$, $(\frac{1}{4},0)$.
\eexa
\end{exa}

Theorem~\ref{thm:0th_density} extends Example~\ref{exa:0th_density}
 to any periodic sequence $S$ and implies that $\psi_0(t)$ is uniquely determined by the ordered distances within a unit cell of $S$.

\begin{thm}[description of $\psi_0$, {\cite[Theorem 5]{anosova2022density}}]
\label{thm:0th_density}
For any periodic sequence $S=\{p_1,\dots,p_m\}+\Z$ with motif points $0\leq p_1<\dots<p_m<1$, set $d_i=p_{i+1}-p_i\in(0,1)$, where $i=1,\dots,m$ and $p_{m+1}=p_1+1$. 
Put the distances in the increasing order $d_{[1]}\leq d_{[2]}\leq\dots\leq d_{[m]}$.
Then the 0-th density function $\psi_0$ is piecewise linear with the following (unordered) corners: $(0,1)$ and $(\frac{1}{2}d_{[i]}, 1-\sum\limits_{j=1}^{i-1} d_{[j]}-(m-i+1)d_{[i]})$ for $i=1,\dots,m$, so the last corner is $(\frac{1}{2}d_{[m]},0)$.
If any corner points are repeated, e.g. if $d_{[i-1]}=d_{[i]}$, these corners are collapsed into one. 
\ethm
\end{thm}

Theorem~\ref{thm:0th_density} for the sequence $S=\{0,\frac{1}{3},\frac{1}{2}\}+\Z$ gives the ordered distances $d_{[1]}=\frac{1}{6}<d_{[2]}=\frac{1}{3}<d_{[3]}=\frac{1}{2}$, which determine the corner points $(0,1)$, $(\frac{1}{12},\frac{1}{2})$, $(\frac{1}{6},\frac{1}{6})$, $(\frac{1}{4},0)$ of the density function $\psi_0(t)$ in Fig.~\ref{fig:densities1D}, see Example~\ref{exa:0th_density}.
\medskip

For any periodic sequence with $m$ points in a unit cell, by Theorem~\ref{thm:0th_density}, any 0th density function $\psi_0(t)$ is uniquely determined by the (unordered) set of lengths of intervals between successive points.
Hence, we can reorder these intervals without changing $\psi_0(t)$.
For instance, the periodic sequence $Q=\{0,\frac{1}{2},\frac{2}{3}\}+\Z$ has the same set of interval lengths $d_{[1]}=\frac{1}{6}$, $d_{[2]}=\frac{1}{3}$, $d_{[3]}=\frac{1}{2}$ as the periodic sequence $S=\{0,\frac{1}{3},\frac{1}{2}\}+\Z$ in Example~\ref{exa:0th_density}.
\medskip

The above sequences $S,Q$ are related by the mirror reflection $t\mapsto 1-t$.
One can easily construct many non-isometric sequences with $\psi_0[S](t)=\psi_0[Q](t)$.
For any $1\leq i\leq m-3$, the sequences $S_{m,i}=\{0,2,3,\dots,i+2,i+4,i+5,\dots,m+2\}+(m+2)\Z$ have the same interval lengths $d_{[1]}=\dots=d_{[m-2]}=1$, $d_{[m-1]}=d_{[m]}=2$ but are not related by isometry (translations and reflections in $\R$) because the intervals of length 2 are separated by $i-1$ intervals of length 1 in $S_{m,i}$.  
\medskip

Corollary~\ref{cor:gen_complete} will prove that the 1st density function $\psi_1[S](t)$  uniquely determines a periodic sequence $S\subset\R$ in general position up to isometry of $\R$.

\begin{exa}[functions $\psi_k(t)$ for $S=\{0,\frac{1}{3},\frac{1}{2}\}+\Z$]
\label{exa:densities}
The 1st density function $\psi_1(t)$ can be obtained as a sum of the three \emph{trapezium} functions $\eta_R$, $\eta_G$, $\eta_B$, each measuring the length of a region covered by a single interval (of one color).
The red intervals $[0,t]\cup[1-t,1]$ grow until $t=\frac{1}{6}$ when they touch the green interval $[\frac{1}{6},\frac{1}{2}]$.
So the length $\eta_R(t)$ of this interval linearly grows from the origin $(0,0)$ to the corner point $(\frac{1}{6},\frac{1}{3})$.
For $t\in[\frac{1}{6},\frac{1}{4}]$, the left red interval is shrinking at the same rate due to the overlapping green interval, while the right red interval continues to grow until $t=\frac{1}{4}$, when it touches the blue interval $[\frac{1}{4},\frac{3}{4}]$.  
Hence the graph of $\eta_R(t)$ remains constant up to the corner point $(\frac{1}{4},\frac{1}{3})$.
After that $\eta_R(t)$ linearly returns to the $t$-axis at $t=\frac{5}{12}$.
Hence the trapezium function $\eta_R$ has the piecewise linear graph through the corner points $(0,0)$, $(\frac{1}{6},\frac{1}{3})$, $(\frac{1}{4},\frac{1}{3})$, $(\frac{5}{12},0)$. 
\medskip

The 2nd function $\psi_2(t)$ is the sum of the \emph{trapezium} functions $\eta_{GB},\eta_{RG},\eta_{RB}$, each measuring the length of a double intersection.
For the green interval $[\frac{1}{3}-t,\frac{1}{3}+t]$ and the blue interval $[\frac{1}{2}-t,\frac{1}{2}+t]$, the graph of the trapezium function $\eta_{GB}(t)$ is piecewise linear and starts at the point $(\frac{1}{12},0)$, where the intervals touch.
The green-blue intersection interval $[\frac{1}{2}-t,\frac{1}{3}+t]$ grows until $t=\frac{1}{4}$, when $[\frac{1}{4},\frac{7}{12}]$ touches the red interval on the left.
At the same time $\eta_{GB}(t)$ is linearly growing to the point $(\frac{1}{4},\frac{1}{3})$.
For $t\in[\frac{1}{4},\frac{1}{3}]$, the green-blue intersection interval becomes shorter on the left, but grows at the same rate on the right until $[\frac{1}{3},\frac{2}{3}]$ touches the red interval $[\frac{2}{3},1]$.
Then $\eta_{GB}(t)$ remains constant up to the point $(\frac{1}{3},\frac{1}{3})$.
For $t\in[\frac{1}{3},\frac{1}{2}]$ the green-blue intersection interval is shortening from both sides.
The graph of $\eta_{GB}(t)$ returns to the $t$-axis at $(\frac{1}{2},0)$, see Fig.~\ref{fig:densities1D}.
\eexa
\end{exa}

Theorem~\ref{thm:densities} extends Example~\ref{exa:densities} 
and proves that any $\psi_k(t)$ is a sum of trapezium functions whose corners are explicitly described. 
We consider any index $i=1,\dots,m$ (of a point $p_i$ or a distance $d_i$) modulo $m$ so that $m+1\equiv 1\pmod{m}$.

\begin{thm}[description of $\psi_k$ for $k>0$, {\cite[Theorem~7]{anosova2022density}}]
\label{thm:densities}
For any periodic sequence $S=\{p_1,\dots,p_m\}+\Z$ with points $0\leq p_1<\dots<p_m<1$ in a motif, set $d_i=p_{i+1}-p_i\in(0,1)$, where $i=1,\dots,m$ and $p_{m+1}=p_1+1$. 
Any interval $[p_i-t,p_i+t]$ is projected to $[0,1]$ modulo $\Z$.
For $1\leq k\leq m$, the density function $\psi_k(t)$ is the sum of $m$ \emph{trapezium} functions $\eta_{k,i}$ with the corner points 
$(\frac{s}{2},0)$, 
$(\frac{d_{i-1}+s}{2},d)$, 
$(\frac{s+d_{i+k-1}}{2},d)$, 
$(\frac{d_{i-1}+s+d_{i+k-1}}{2},0)$, where 
$d=\min\{d_{i-1},d_{i+k-1}\}$, $s=\sum\limits_{j=i}^{i+k-2}d_j$,
$i=2,\dots,m+1$.
If $k=1$, then $s=0$ is the empty sum.
So $\psi_k(t)$ is determined by the unordered set of triples $(d_{i-1},s,d_{i+k-1})$ whose first and last entries are swappable.
\ethm
\end{thm}

In Example~\ref{exa:densities} for $S=\{0,\frac{1}{3},\frac{1}{2}\}+\Z$, we have $d_{1}=\frac{1}{3}$, $d_{2}=\frac{1}{6}$, $d_{3}=\frac{1}{2}=d_0$.
For $k=2$, $i=2$, we get $d_{i-1}=d_1=\frac{1}{3}$, $d_{i+k-1}=d_3=\frac{1}{2}$, i.e. $d=\min\{d_1,d_3\}=\frac{1}{3}$, $s=d_{2}=\frac{1}{6}$.
Then $\eta_{22}=\eta_{GB}$ has the corner points $(\frac{1}{12},0)$, $(\frac{1}{4},\frac{1}{3})$, $(\frac{1}{3},\frac{1}{3})$, $(\frac{1}{2},0)$. 

\begin{thm}[symmetries of $\psi_k(t)$, {\cite[Theorem 8]{anosova2022density}}]
\label{thm:symmetries}
For any periodic sequence $S\subset\R$ with a unit cell $[0,1]$, we have the \emph{periodicity} 
$\psi_{k+m}(t+\frac{1}{2})=\psi_{k}(t)$ for any $k\geq 0$, $t\geq 0$, and the \emph{symmetry} $\psi_{m-k}(\frac{1}{2}-t)=\psi_k(t)$ for $k=0,\dots,[\frac{m}{2}]$, and $t\in[0,\frac{1}{2}]$.
\ethm
\end{thm}

\begin{cor}[time of $\psi_k(t)$ for periodic sequences of points, {\cite[Corollary~9]{anosova2022density}}]
\label{cor:densities1D_points_time}
Let $S,Q\subset\R$ be periodic sequences with at most $m$ points in motifs.
For $k\geq 1$, one can draw the graph of the $k$-th density function $\psi_k[S]$ in time $O(m^2)$.
One can check in time $O(m^3)$ if 
the full density fingerprints coincide: 
$\Psi[S]=\Psi[Q]$.
\ethm
\end{cor}

To illustrate Corollary~\ref{cor:densities1D_points_time}, Example~\ref{exa:SQ15} will justify that the periodic sequences $S_{15}$ and $Q_{15}$ in Fig.~\ref{fig:SQ15} have identical density fingerprints $\Psi[S_{15}]=\Psi[Q_{15}]$.

\begin{exa}[periodic sequences $S_{15},Q_{15}\subset\R$]
\label{exa:SQ15}
\cite[Appendix~B]{widdowson2022average} discusses homometric periodic sets that can be distinguished by the invariant AMD (Average Minimum Distances) and not by inter-point distance distributions.
The periodic sequences
$$S_{15} = \{0,1,3,4,5,7,9,10,12\}+15\Z,\quad
Q_{15} = \{0,1,3,4,6,8,9,12,14\}+15\Z$$ have period 15 and unit cell $[0,15]$ shown as a circle in Fig.~\ref{fig:SQ15}.
\myskip

\begin{figure}[h!]
\includegraphics[width=\textwidth]{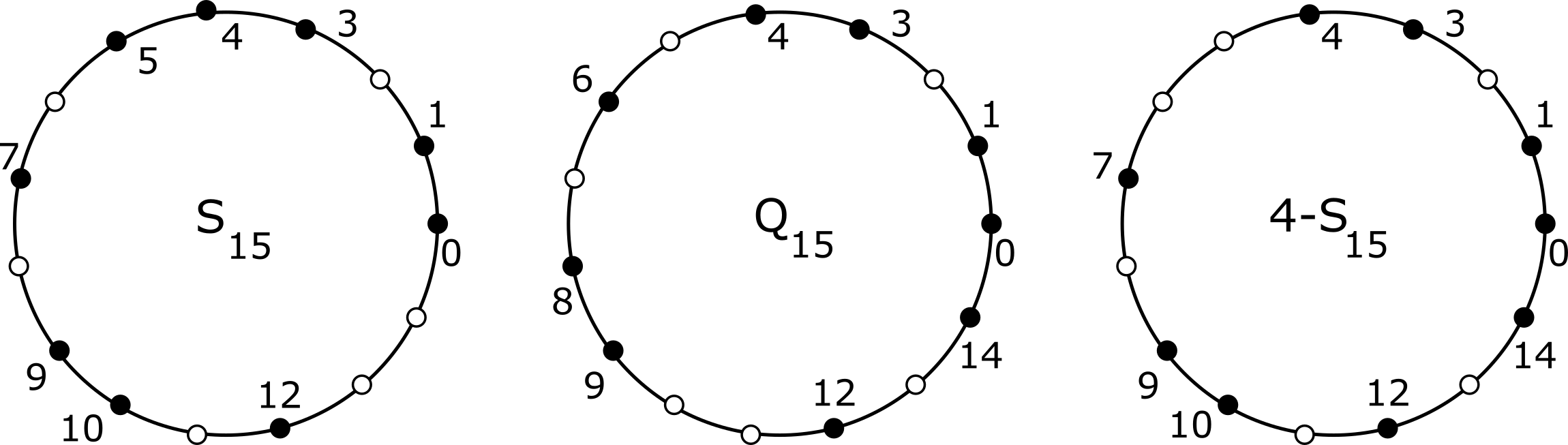}
\caption{Circular versions of the periodic sets $S_{15},Q_{15}$.
Distances are along round arcs.}
\label{fig:SQ15}      
\end{figure}

These periodic sequences \cite{grunbaum1995use} are obtained as Minkowski sums $S_{15}=U+V+15\Z$ and $Q_{15}=U-V+15\Z$ for 
$U = \{0, 4, 9\}$ and $V = \{0, 1, 3\}$.
The last picture in Fig.~\ref{fig:SQ15} shows the periodic set $4-S_{15}$ isometric to $S_{15}$.
Now the difference between $Q_{15}$ and $4-S_{15}$ is better visible: points $0,1,3,4,5,12,14$ are common, but points $6,8,9\in Q_{15}$ are shifted to $7,9,10$ in the circular set $4-S_{15}$.
\myskip

To avoid fractions, we keep the unit cell $[0,15]$ of the sequences $S_{15},Q_{15}$ without scaling it down to $[0,1]$ 
because all quantities in Theorem~\ref{thm:densities} can be scaled up by factor 15. 
To conclude that $\psi_0[S_{15}]=\psi_0[Q_{15}]$, by 
Theorem~\ref{thm:0th_density} we check that $S_{15},Q_{15}$ have the same set of the ordered distances $d_{[i]}$ between successive points, which is shown in identical rows 3 of Tables~\ref{tab:S15d} and~\ref{tab:Q15d}.

\begin{table}[h!]
\caption{\textbf{Row 1}: points $p_i$ from $S_{15}$ in Fig.~\ref{fig:SQ15}.
\textbf{Row 2}: the distances $d_i$ between successive points of $S_{15}$. 
\textbf{Row 3}: the distances $d_{[i]}$ are in the increasing order.
\textbf{Row 4}: the unordered set of these pairs determines the density function $\psi_1$ by Theorem~\ref{thm:densities}. 
\textbf{Row 5}: the pairs are lexicographically ordered for comparison with row 5 in Table~\ref{tab:Q15d}. 
\textbf{Rows 6,8,10}: the unordered sets of these triples determine the density functions $\psi_2,\psi_3,\psi_4$ by Theorem~\ref{thm:densities} for $k=2,3,4$.
\textbf{Rows 7,9,11}: the triples from rows 6,8,10 are ordered for easier comparison with corresponding rows 7,9,11 in Table~\ref{tab:Q15d}, see  Example~\ref{exa:SQ15}.
}
\label{tab:S15d}
\begin{tabular}{C{30mm}|C{9mm}C{9mm}C{9mm}C{9mm}C{9mm}C{9mm}C{9mm}C{9mm}C{9mm}}
$p_i$                      & 0 & 1 & 3 & 4 & 5 & 7 & 9 & 10 & 12 \\
\hline
\hline
$d_i=p_{i+1}-p_i$ & 1 & 2 & 1 & 1 & 2 & 2 & 1 &  2 & 3 \\
ordered $d_{[i]}$ & 1 & 1 & 1 & 1 & 2 & 2 & 2 &  2 & 3 \\
\hline
\hline
$(d_{i-1},d_i)$  & (3,1) & (1,2) & (2,1) & (1,1) & (1,2) & (2,2) & (2,1) & (1,2) & (2,3) \\
order $(d_{i-1},d_i)$  & (1,1) & (1,2) & (1,2) & (1,2) & (1,2) & (1,2) & (1,3) & (2,2) & (2,3) \\
\hline
\hline
$(d_{i-1},\mathbf{d_i},d_{i+1})$  & (3,{\bf 1},2) & (1,{\bf 2},1) & (2,{\bf 1},1) & (1,{\bf 1},2) & (1,{\bf 2},2) & (2,{\bf 2},1) & (2,{\bf 1},2) & (1,{\bf 2},3) & (2,{\bf 3},1) \\
order $(d_{i-1},\mathbf{d_i},d_{i+1})$ & (1,{\bf 1},2)  & (1,{\bf 1},2) & (2,{\bf 1},2) & (2,{\bf 1},3)  & (1,{\bf 2},1) & (1,{\bf 2},2) & (1,{\bf 2},2) & (1,{\bf 2},3) & (1,{\bf 3},2) \\
\hline
\hline
$(d_{i-1},\mathbf{s},d_{i+2})$  & (3,{\bf 3},1) & (1,{\bf 3},1) & (2,{\bf 2},2) & (1,{\bf 3},2) & (1,{\bf 4},1) & (2,{\bf 3},2) & (2,{\bf 3},3) & (1,{\bf 5},1) & (2,{\bf 4},2) \\
order $(d_{i-1},\mathbf{s},d_{i+2})$ & (2,{\bf 2},2) & (1,{\bf 3},1) & (1,{\bf 3},2) & (1,{\bf 3},3)  & (2,{\bf 3},2)  & (2,{\bf 3},3) & (1,{\bf 4},1) & (2,{\bf 4},2) & (1,{\bf 5},1) \\
\hline
\hline
$(d_{i-1},\mathbf{s},d_{i+3})$  & (3,{\bf 4},1) & (1,{\bf 4},2) & (2,{\bf 4},2) & (1,{\bf 5},1) & (1,{\bf 5},2) & (2,{\bf 5},3) & (2,{\bf 6},1) & (1,{\bf 6},2) & (2,{\bf 6},1) \\
order $(d_{i-1},\mathbf{s},d_{i+3})$  & (1,{\bf 4},2) & (1,{\bf 4},3) & (2,{\bf 4},2) & (1,{\bf 5},1) & (1,{\bf 5},2) & (2,{\bf 5},3) & (1,{\bf 6},2) & (1,{\bf 6},2) & (1,{\bf 6},2) 
\end{tabular}
\end{table}

\begin{table}[h!]
\caption{\textbf{Row 1}: points $p_i$ from $Q_{15}$ in Fig.~\ref{fig:SQ15}.
\textbf{Row 2}: the distances $d_i$ between successive points of $Q_{15}$. 
\textbf{Row 3}: the distances $d_{[i]}$ are in the increasing order.
\textbf{Row 4}: the unordered set of these pairs determines the density function $\psi_1$ by Theorem~\ref{thm:densities}b. 
\textbf{Row 5}: the pairs are lexicographically ordered for comparison with row 5 in Table~\ref{tab:S15d}. 
\textbf{Rows 6,8,10}: the unordered sets of these triples determine the density functions $\psi_2,\psi_3,\psi_4$ by Theorem~\ref{thm:densities} for $k=2,3,4$.
\textbf{Rows 7,9,11}: the triples from rows 6,8,10 are ordered for comparison with corresponding rows 7,9,11 in Table~\ref{tab:S15d}, see  Example~\ref{exa:SQ15}.
}
\label{tab:Q15d}
\begin{tabular}{C{30mm}|C{9mm}C{9mm}C{9mm}C{9mm}C{9mm}C{9mm}C{9mm}C{9mm}C{9mm}}
\hline
$p_i$                      & 0 & 1 & 3 & 4 & 6 & 8 & 9 & 12 & 14 \\
\hline
\hline
$d_i=p_{i+1}-p_i$ & 1 & 2 & 1 & 2 & 2 & 1 & 3 &  2 & 1 \\
ordered $d_{[i]}$ & 1 & 1 & 1 & 1 & 2 & 2 & 2 &  2 & 3 \\
\hline
\hline
$(d_{i-1},d_i)$  & (1,1) & (1,2) & (2,1) & (1,2) & (2,2) & (2,1) & (1,3) & (3,2) & (2,1) \\
ordered $(d_{i-1},d_i)$  & (1,1) & (1,2) & (1,2) & (1,2) & (1,2) & (1,2) & (1,3) & (2,2) & (2,3) \\
\hline
\hline
$(d_{i-1},\mathbf{d_i},d_{i+1})$  & (1,{\bf 1},2) & (1,{\bf 2},1) & (2,{\bf 1},2) & (1,{\bf 2},2) & (2,{\bf 2},1) & (2,{\bf 1},3) & (1,{\bf 3},2) & (3,{\bf 2},1) & (2,{\bf 1},1) \\
order $(d_{i-1},\mathbf{d_i},d_{i+1})$ & (1,{\bf 1},2)  & (1,{\bf 1},2) & (2,{\bf 1},2) & (2,{\bf 1},3)  & (1,{\bf 2},1) & (1,{\bf 2},2) & (1,{\bf 2},2) & (1,{\bf 2},3) & (1,{\bf 3},2) \\
\hline
\hline
$(d_{i-1},\mathbf{s},d_{i+2})$  & (1,{\bf 3},1) & (1,{\bf 3},2) & (2,{\bf 3},2) & (1,{\bf 4},1) & (2,{\bf 3},3) & (2,{\bf 4},2) & (1,{\bf 5},1) & (3,{\bf 3},1) & (2,{\bf 2},2) \\
order $(d_{i-1},\mathbf{s},d_{i+2})$ & (2,{\bf 2},2) & (1,{\bf 3},1) & (1,{\bf 3},2) & (1,{\bf 3},3)  & (2,{\bf 3},2)  & (2,{\bf 3},3) & (1,{\bf 4},1) & (2,{\bf 4},2) & (1,{\bf 5},1) \\
\hline
\hline
$(d_{i-1},\mathbf{s},d_{i+3})$  & (1,{\bf 4},2) & (1,{\bf 5},2) & (2,{\bf 5},1) & (1,{\bf 5},3) & (2,{\bf 6},2) & (2,{\bf 6},1) & (1,{\bf 6},1) & (3,{\bf 4},2) & (2,{\bf 4},1) \\
order $(d_{i-1},\mathbf{s},d_{i+3})$  & (1,{\bf 4},2) & (1,{\bf 4},2) & (2,{\bf 4},3) & (1,{\bf 5},2) & (1,{\bf 5},2) & (1,{\bf 5},3)  & (1,{\bf 6},1) & (1,{\bf 6},2) & (2,{\bf 6},2)
\end{tabular}
\end{table} 

To conclude that $\psi_1[S_{15}]=\psi_1[Q_{15}]$ by Theorem~\ref{thm:densities}, we check that $S_{15},Q_{15}$ have the same set of unordered pairs $(d_{i-1},d_i)$ of distances between successive points.
Indeed, Tables~\ref{tab:S15d} and~\ref{tab:Q15d} have identical rows 5, where pairs are \emph{lexicograpically} ordered for comparison: $(a,b)<(c,d)$ if $a<b$ or $a=b$ and $c<d$.
\medskip

To conclude that $\psi_k[S_{15}]=\psi_k[Q_{15}]$ for $k=2,3,4$, 
 we compare the triples $(d_{i-1},\mathbf{s},d_{i+k-1})$ from Theorem~\ref{thm:densities} for $S_{15},Q_{15}$.
For $k=2$ and $k=3$, Tables~\ref{tab:S15d} and~\ref{tab:Q15d} have identical rows 7 and 9, where the triples are ordered for easier comparison as follows.
If needed, we swap $d_{i-1},d_{i+k-1}$ to make sure that the first entry is not larger than the last.
Then we order by the middle bold number $\mathbf s$.
Finally, we lexicographically order the triples with the same middle value $s$.
\medskip

Final rows 11 of Tables~\ref{tab:S15d} and~\ref{tab:Q15d} look different for $k=4$.
More exactly, the rows share three triples (1,{\bf 4},2), (1,{\bf 5},2), (1,{\bf 6},4), but the remaining six triples differ.
However, the density function $\psi_4$ is the \emph{sum} of nine trapezium functions.
Fig.~\ref{fig:SQ15density4} shows that these sums are equal for $S_{15},Q_{15}$.
Then the sequences $S_{15},Q_{15}$ have identical density functions $\psi_k$ for $k=0,1,2,3,4$, hence for all $k$ by the symmetry and periodicity from Theorem~\ref{thm:symmetries}.
Fig.~\ref{fig:SQ15densities} shows the density functions $\psi_k$ for $k=0,1,\dots,9$.
\eexa
\end{exa}


\begin{figure}[ht]
\includegraphics[width=\textwidth]{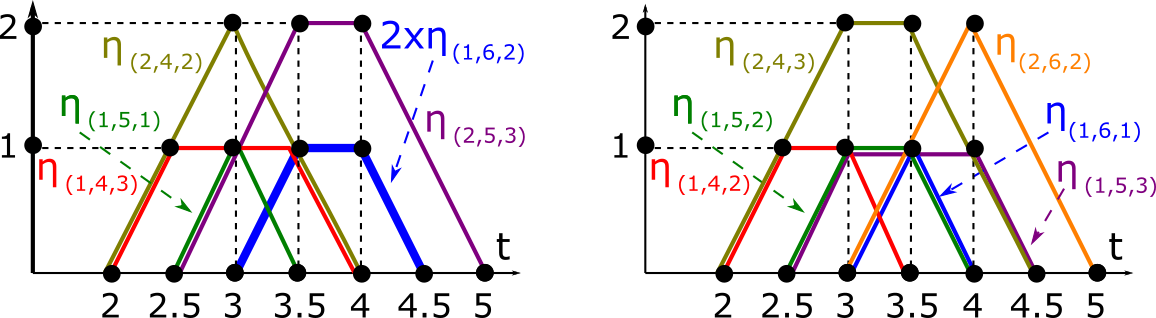}
\caption{The 4th-density function $\psi_4[S_{15}]$ includes the six trapezium functions on the left, which are replaced by other six trapezium functions in $\psi_4[Q_{15}]$ on the right, compare the last rows of Tables~\ref{tab:S15d} and~\ref{tab:Q15d}.
However, the sums of these six functions are equal, which can be checked at critical radii: both sums of six functions have $\eta(2.5)=2$, $\eta(3)=5$, $\eta(3.5)=6$, $\eta(4)=4$, $\eta(4.5)=1$.
Hence, the sequences $S_{15},Q_{15}$ in Fig.~\ref{fig:SQ15} have identical density functions $\psi_k$ for all $k\geq 0$, see Example~\ref{exa:SQ15}.} 
\label{fig:SQ15density4}      
\end{figure}

\begin{figure}[h!]
\includegraphics[width=\textwidth]{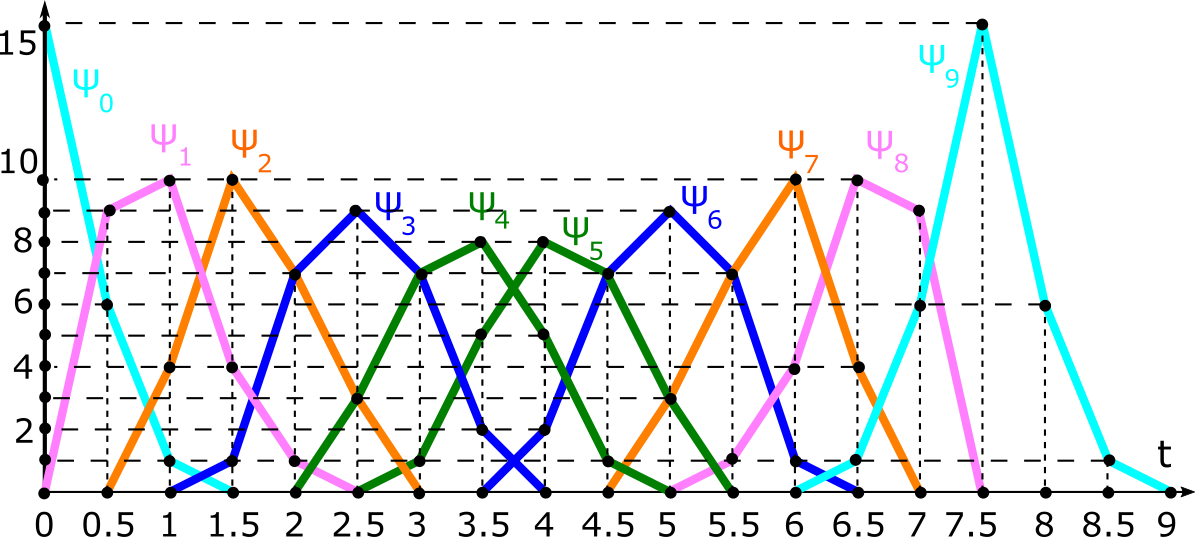}
\caption{The periodic sequences $S_{15},Q_{15}$ in Fig.~\ref{fig:SQ15} have identical density functions $\psi_k(t)$ for all $k\geq 0$.
Both axes are scaled by factor 15.
Theorem~\ref{thm:symmetries} implies
the symmetry $\psi_k(\frac{15}{2}-t)=\psi_{9-k}(t)$, $t\in[0,\frac{15}{2}]$, and periodicity $\psi_9(t+\frac{15}{2})=\psi_0(t)$, $t\geq 0$.}
\label{fig:SQ15densities}      
\end{figure}

Recall that all indices $i$ of distances $d_i$ are considered modulo $m$.  

\begin{cor}[$k$-th density $\rho_k$, {\cite[Corollary~11]{anosova2022density}}]
\label{cor:rho-densities}
For any periodic sequence $S=\{p_1,\dots,p_m\}+\Z$ with inter-point distances $d_i=p_{i+1}-p_i$, where $i=1,\dots,m$ and $p_{m+1}=p_1+1$, the $k$-th \emph{density} $\rho_k[S]=\int\limits_{-\infty}^{+\infty} \psi_k(t)dt$ defined as the area under the graph of $\psi_k(t)$ over $\R$ equals
$\rho_k[S]=\dfrac{1}{2}\sum\limits_{i=1}^m d_{i-1} d_{i+k-1}$ for any $k>0$ and $\rho_0[S]=\dfrac{1}{4}\sum\limits_{i=1}^md_i^2$.
\ethm
\end{cor}

For $S=\{0,\frac{1}{3},\frac{1}{2}\}+\Z$,
Corollary~\ref{cor:rho-densities} gives $\rho_0=\frac{7}{72}$, 
$\rho_1=\rho_2=\frac{11}{12^2}$ as in Fig.~\ref{fig:densities1D}.  

\begin{cor}[generic completeness of $\psi_1$, {\cite[Corollary~12]{anosova2022density}}]
\label{cor:gen_complete}
Let $S\subset\R$ be a sequence with period 1 and $m$ points $0\leq p_1<\dots<p_m<1$.
The sequence $S$ is called \emph{generic} if $d_i=p_{i+1}-p_i$ are distinct, where $i=1,\dots,m$ and $p_{m+1}=p_1+1$.
Then any generic $S$ can be reconstructed from the 1st density function $\psi_1[S](t)$ up to isometry in $\R$.
Hence $\psi_1(t)$ is a complete isometry invariant for all generic $S$. 
\ethm
\end{cor}

\section{Density functions of periodic sequences of intervals in $\R$}
\label{sec:densities1D_intervals}

This section follows \cite[sections 2-5]{anosova2023density} by extending density functions to periodic sets of points with radii, including periodic sequences of disjoint intervals in $\R$.

\index{density function}

\begin{dfn}[density functions for periodic sets of points with radii]
\label{dfn:densities_radii}
Let a periodic set $S=\La+M\subset\R^n$ have a unit cell $U$.
For every point $p\in M$, fix a radius $r(p)\geq 0$.
For any integer $k\geq 0$, let $U_k(t)$ be the region within the cell $U$ covered by exactly $k$ closed balls $\bar B(p;r(p)+t)$ for $t\geq 0$ and all  points $p\in M$ and their translations by $\La$.
The $k$-th \emph{density} function $\psi_k[S](t)=\vol[U_k(t)]/\vol[U]$ is the fractional volume of the $k$-fold intersections of these balls within $U$.
\ethm
\end{dfn}

In Definition~\ref{dfn:densities_radii}, the balls are growing at all points of $S$, because centers $p\in M$ are translated by all lattice vectors $v\in\La$.
The initially different radii $r_i$ are motivated by real lengths of continuous events in periodic time series for $n=1$ and also by atomic radii of different chemical elements for $n=3$.  
Another (possibly, non-linear) growth of radii lead to more complicated density functions.
\medskip

The density $\psi_k[S](t)$ can be interpreted as the probability that a random (uniformly chosen in $U$) point $q$ is at a maximum distance $t$ to exactly $k$ balls with initial radii $r(p)$ and all centers $p\in S$.
For $k=0$, the $0$-th density $\psi_0[S](t)$ measures the fractional volume of the empty space not covered by any expanding balls $\bar B(p;r(p)+t)$
\medskip

For $k=1$ and small $t>0$ while all equal-sized balls $\bar B(p;t)$ remain disjoint, the 1st density $\psi_1[S](t)$ increases proportionally to $t^n$ but later reaches a maximum and eventually drops back to $0$ when all points of $\R^n$ are covered of by at least two balls.
\myskip

The original densities helped find a missing crystal in the Cambridge Structural Database, which was accidentally confused with a slight perturbation (measured at a different temperature) of another crystal (polymorph) with the same chemical composition, see \cite[section~7]{edelsbrunner2021density}. 
\medskip

The new weighted case with radii $r(p)\geq 0$ in Definition~\ref{dfn:densities_radii} is even more practically important due to different Van der Waals radii, which are individually defined for all chemical elements.
\medskip

The key advantage of density functions over other isometry invariants of periodic crystals (such as symmetries or conventional representations based on a geometry of a minimal cell) is their continuity under perturbations.
The only limitation is the infinite size of densities $\psi_k(t)$ due to the unbounded parameters: integer index $k\geq 0$ and continuous radius $t\geq 0$. 
\myskip

Theorem~\ref{thm:0-th_density} will explicitly describing the 0-th density function $\psi_0[S](t)$ for any periodic sequence $S\subset\R$ of intervals.
All intervals are considered closed and called \emph{disjoint} if their open interiors (not endpoints) have no common points.
\medskip

For convenience, scale any periodic sequence $S$ to period 1 so that $S$ is given by points $0\leq p_1<\cdots<p_m<1$ with radii $r_1,\dots,r_m$, respectively.
Since the expanding balls in $\R$ are growing intervals, volumes of their intersections linearly change with respect to the variable radius $t$.
Hence any density function $\psi_k(t)$ is piecewise linear and uniquely determined by \emph{corner} points $(a_j,b_j)$ where the gradient of $\psi_k(t)$ changes.
To illustrate Theorem~\ref{thm:0-th_density}, we consider 
Example~\ref{exa:0-th_density} for the simple sequence $S$.

\begin{exa}[$0$-th density function $\psi_0$]
\label{exa:0-th_density}
Let the periodic sequence $S=\{0,\frac{1}{3},\frac{1}{2}\}+\Z$ have three points $p_1=0$, $p_2=\frac{1}{3}$, $p_3=\frac{1}{2}$ of radii $r_1=\frac{1}{12}$, $r_2=0$, $r_3=\frac{1}{12}$, respectively.
Fig.~\ref{fig:growing_intervals} shows each point $p_i$ and its growing interval 
$$L_i(t)=[(p_i-r_i)-t,(p_i+r_i)+t] \text{ of the length }2r_i+2t$$ 
for $i=1,2,3$ in its own color: red, green, blue.
\medskip

\begin{figure}[h!]
\includegraphics[width=\linewidth]{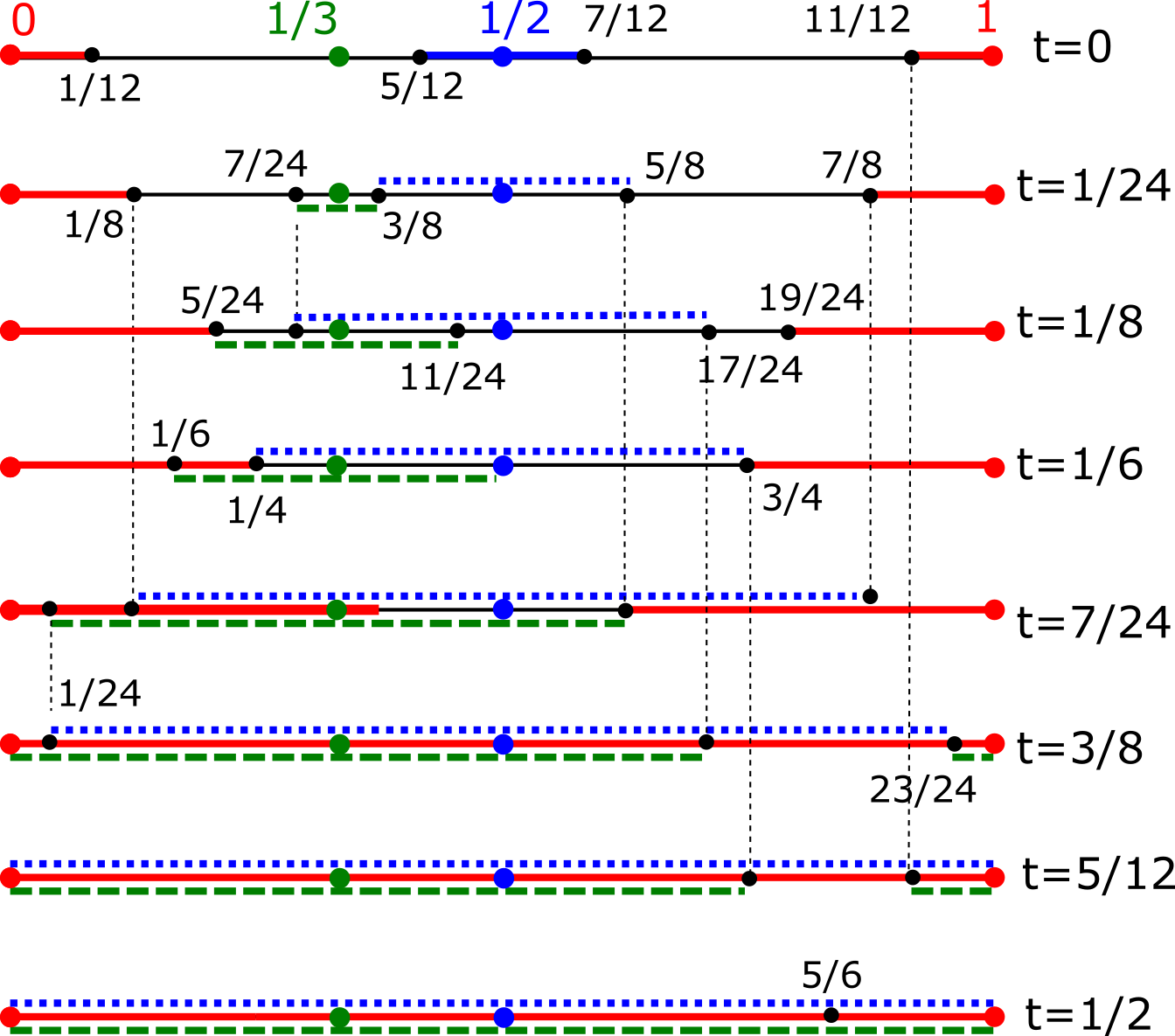}
\caption{The sequence $S=\{0,\frac{1}{3},\frac{1}{2}\}+\Z$ has the points of weights $\frac{1}{12},0,\frac{1}{12}$, respectively. 
The intervals around the red point $0\equiv 1\pmod{1}$, green point $\frac{1}{3}$, blue point $\frac{1}{2}$ have the same color for various radii $t$, see Examples~\ref{exa:0-th_density},~\ref{exa:1st_density},~\ref{exa:2nd_density}.
}
\label{fig:growing_intervals}      
\end{figure}

By Definition~\ref{dfn:densities_radii}, each density function $\psi_k[S](t)$ measures a fractional length covered by exactly $k$ intervals within the unit cell $[0,1]$.
It is convenient to periodically map the endpoints of each growing interval to the unit cell $[0,1]$.
\medskip

For instance, the interval $[-\frac{1}{12}-t,\frac{1}{12}+t]$ of the point $p_1=0\equiv 1\pmod{1}$ maps to the red intervals $[0,\frac{1}{12}+t]\cup[\frac{11}{12}-t,1]$ shown by solid red lines in Fig.~\ref{fig:growing_intervals}.
The same image shows the green interval $[\frac{1}{3}-t,\frac{1}{3}+t]$ by dashed lines and the blue interval $[\frac{5}{12}-t,\frac{7}{12}+t]$ by dotted lines.
\medskip

At the moment $t=0$, since the starting intervals are disjoint, they cover the length $l=2(\frac{1}{12} + 0 +\frac{1}{12})=\frac{1}{3}$.
The non-covered part of $[0,1]$ has length $1-\frac{1}{3}=\frac{2}{3}$.
So the graph of $\psi_0(t)$ at $t=0$ starts from the point $(0,\frac{2}{3})$, see Fig.~\ref{fig:0-th_density}~(right).
\medskip

\begin{figure}[h!]
\includegraphics[height=50mm]{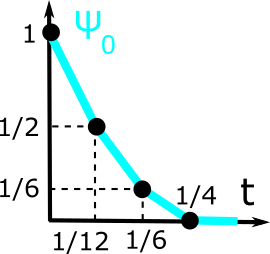}
\hspace*{1pt}
\includegraphics[height=50mm]{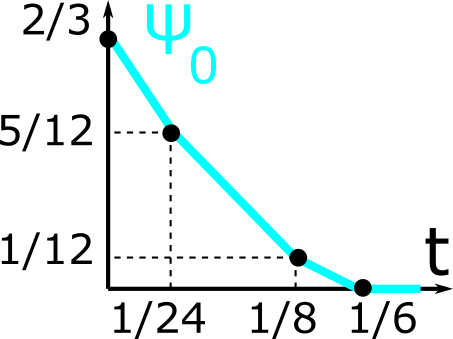}
\caption{
\textbf{Left}: the 0-th density function $\psi_0(t)$ for the 1-period sequence $S=\{0,\frac{1}{3},\frac{1}{2}\}+\Z$ with radii 0.
\textbf{Right}: the 0-th density $\psi_0(t)$ for the 1-period sequence $S$ whose points $0,\frac{1}{3},\frac{1}{2}$ have radii $\frac{1}{12},0,\frac{1}{12}$, respectively, see Example~\ref{exa:0-th_density}. }
\label{fig:0-th_density}      
\end{figure}

At the first critical moment $t=\frac{1}{24}$ when the green and blue intervals collide at $p=\frac{3}{8}$, only the intervals $[\frac{1}{8},\frac{7}{24}]\cup[\frac{5}{8},\frac{7}{8}]$ of total length $\frac{5}{12}$ remain uncovered.
Hence $\psi_0(t)$ linearly drops to the point $(\frac{1}{12},\frac{5}{12})$. 
At the next critical moment $t=\frac{1}{8}$ when the red and green intervals collide at $p=\frac{5}{24}$, only the interval $[\frac{17}{24},\frac{19}{24}]$ of length $\frac{1}{12}$ remain uncovered, so $\psi_0(t)$ continues to $(\frac{1}{8},\frac{1}{12})$.
\medskip

The graph of $\psi_0(t)$ finally returns to the $t$-axis at the point $(\frac{1}{6},0)$ and remains there for $t\geq \frac{1}{6}$.
The piecewise linear behaviour of $\psi_0(t)$ can be described by specifying the \emph{corner} points in Fig.~\ref{fig:0-th_density}: 
$(0,\frac{2}{3})$, 
$(\frac{1}{24},\frac{5}{12})$, 
$(\frac{1}{8},\frac{1}{12})$, 
$(\frac{1}{6},0)$.
\eexa
\end{exa}

Theorem~\ref{thm:0-th_density} extends Example~\ref{exa:0-th_density}
 to any periodic sequence $S$ and implies that the 0-th density function $\psi_0(t)$ is uniquely determined by the ordered gap lengths between successive intervals.

\begin{thm}[description of $\psi_0$, {\cite[Theorem~3.2]{anosova2023density}}]
\label{thm:0-th_density}
Let a periodic sequence $S=\{p_1,\dots,p_m\}+\Z$ consist of disjoint intervals with centers $0\leq p_1<\dots<p_m<1$ and radii $r_1,\dots,r_m\geq 0$. 
Consider the \emph{total length} $l=2\sum\limits_{i=1}^m r_i$ and \emph{gaps} between successive intervals $g_i=(p_{i}-r_{i})-(p_{i-1}+r_{i-1})$, where $i=1,\dots,m$ and $p_{0}=p_m-1$, $r_0=r_m$. 
Put the gaps in increasing order: $g_{[1]}\leq g_{[2]}\leq\dots\leq g_{[m]}$.
Then the 0-th density $\psi_0[S](t)$ is piecewise linear with the following (unordered) corner points: $(0,1-l)$ and $(\frac{g_{[i]}}{2},\; 1-l-\sum\limits_{j=1}^{i-1} g_{[j]}-(m-i+1)g_{[i]})$ for $i=1,\dots,m$, so the last corner is $(\frac{g_{[m]}}{2},0)$.
If any corners are repeated, e.g. if $g_{[i-1]}=g_{[i]}$, these corners are collapsed into one. 
\ethm
\end{thm}

Example~\ref{exa:revisit_0-th_density} applies Theorem~\ref{thm:0-th_density} to get $\psi_0$ for the sequence $S$ in Example~\ref{exa:0-th_density}.

\begin{exa}[using Theorem~\ref{thm:0-th_density}]
\label{exa:revisit_0-th_density}
The sequence $S=\{0,\frac{1}{3},\frac{1}{2}\}+\Z$ in Example~\ref{exa:0-th_density} with points $p_1=0$, $p_2=\frac{1}{3}$, $p_3=\frac{1}{2}$ of radii $r_1=\frac{1}{12}$, $r_2=0$, $r_3=\frac{1}{12}$, respectively, has 
$l=2(r_1+r_2+r_3)=\frac{1}{3}$ and the initial gaps 
between successive intervals 
$$\begin{array}{l}
g_1=p_{1}-r_{1}-p_{3}-r_{3}=(1-\frac{1}{12})-(\frac{1}{2}+\frac{1}{12})=\frac{1}{3}, \\
g_2=p_{2}-r_{2}-p_{1}-r_{1}=(\frac{1}{3}-0)-(0+\frac{1}{12})=\frac{1}{4}, \\
g_3=p_{3}-r_{3}-p_{2}-r_{2}=(\frac{1}{2}-\frac{1}{12})-(\frac{1}{3}+0)=\frac{1}{12}.
\end{array}$$
Order the gaps: 
$g_{[1]}=\frac{1}{12}<g_{[2]}=\frac{1}{4}<g_{[3]}=\frac{1}{3}$.
Then
$$\begin{array}{l}
1-l=1-\frac{1}{3}=\frac{2}{3}, \\
1-l-3g_{[1]}=\frac{2}{3}-\frac{3}{12}=\frac{5}{12}, \\
1-l-g_{[1]}-2g_{[2]}=\frac{2}{3}-\frac{1}{12}-\frac{2}{4}=\frac{1}{12}, \\
1-l-g_{[1]}-g_{[2]}-g_{[3]}=\frac{2}{3}-\frac{1}{12}-\frac{1}{4}-\frac{1}{3}=0. \end{array}$$ 
By Theorem~\ref{thm:0-th_density} $\psi_0(t)$ has the corner points
$$\begin{array}{l}
(0,1-l)=(0,\frac{2}{3}), \\
(\frac{1}{2}g_{[1]},1-l-3g_{[1]})=(\frac{1}{24},\frac{5}{12}), \\
(\frac{1}{2}g_{[2]},1-l-g_{[1]}-2g_{[2]})=(\frac{1}{8},\frac{1}{12}), \\
(\frac{1}{2}g_{[3]},1-l-g_{[1]}-g_{[2]}-g_{[3]})=(\frac{1}{6},0).
\end{array}$$ 

See the graph of the 0-th density $\psi_0(t)$ in Fig.~\ref{fig:0-th_density}.
\eexa
\end{exa}

By Theorem~\ref{thm:0-th_density} any 0-th density function $\psi_0(t)$ is uniquely determined by the (unordered) set of gap lengths between successive intervals.
Hence we can re-order these intervals without changing $\psi_0(t)$.
For instance, the periodic sequence $Q=\{0,\frac{1}{2},\frac{2}{3}\}+\Z$ with points $0,\frac{1}{2},\frac{2}{3}$ of weights $\frac{1}{12},\frac{1}{12},0$ has the same set ordered gaps $g_{[1]}=\frac{1}{12}$, $d_{[2]}=\frac{1}{3}$, $d_{[3]}=\frac{1}{2}$ as the periodic sequence $S=\{0,\frac{1}{3},\frac{1}{2}\}+\Z$ in Example~\ref{exa:0-th_density}.
\medskip

The above sequences $S,Q$ are related by the mirror reflection $t\mapsto 1-t$.
One can easily construct many non-isometric sequences with $\psi_0[S](t)=\psi_0[Q](t)$.
For any $1\leq i\leq m-3$, the sequences $S_{m,i}=\{0,2,3,\dots,i+2,i+4,i+5,\dots,m+2\}+(m+2)\Z$ have the same interval lengths $d_{[1]}=\dots=d_{[m-2]}=1$, $d_{[m-1]}=d_{[m]}=2$ but are not related by isometry (translations and reflections in $\R$) because the intervals of length 2 are separated by $i-1$ intervals of length 1 in $S_{m,i}$.  
\bskip

Theorem~\ref{thm:1st_density} will explicitly describe the 1st density function $\psi_1[S](t)$ for any periodic sequence $S$ of disjoint intervals.
To illustrate Theorem~\ref{thm:1st_density}, 
Example~\ref{exa:1st_density} finds $\psi_1[S]$ for the sequence $S$ from Example~\ref{exa:0-th_density}.

\begin{figure}[h!]
\includegraphics[height=130mm]{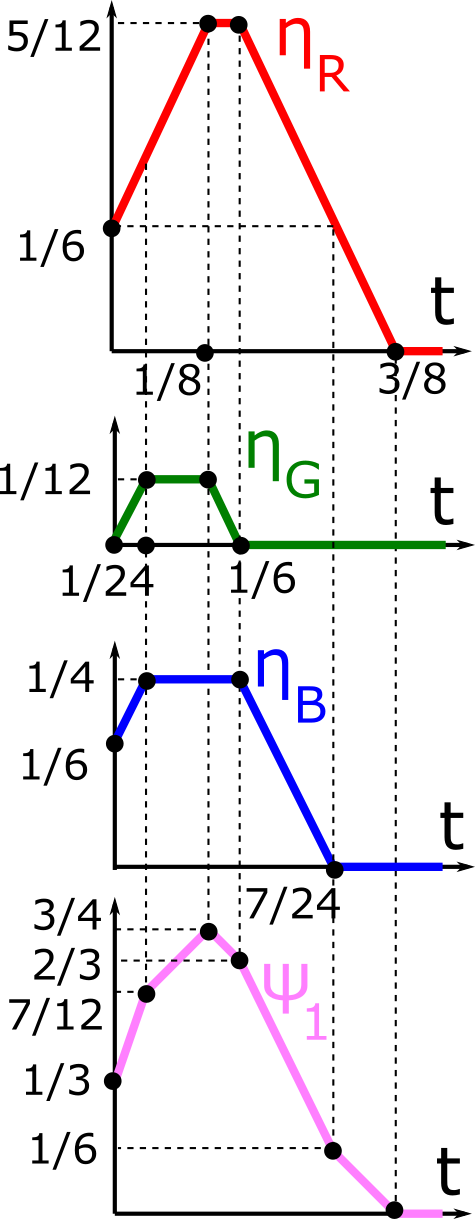}
\hspace*{4mm}
\includegraphics[height=130mm]{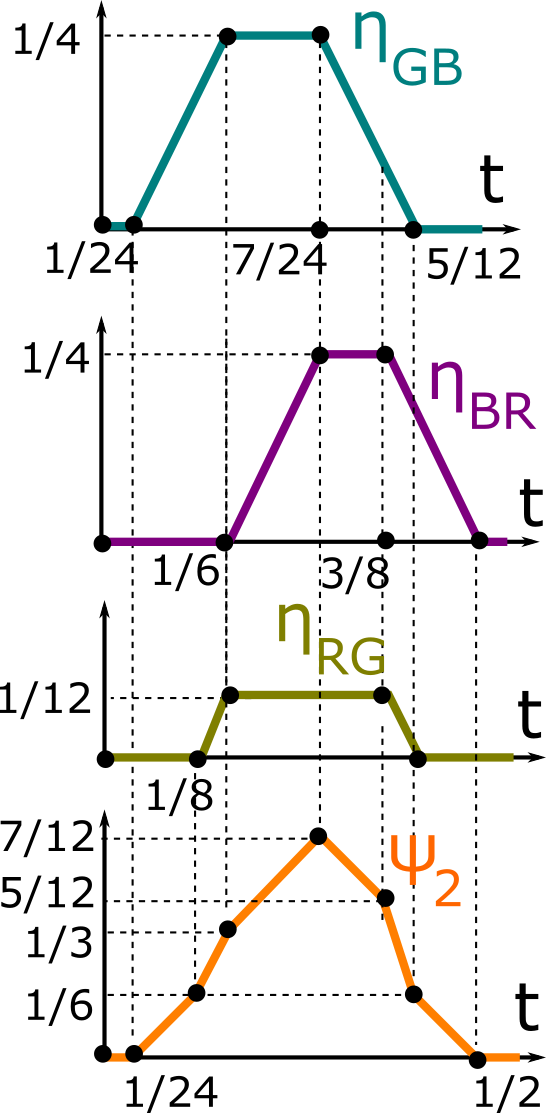}
\caption{\textbf{Left}: the trapezium functions $\eta_R,\eta_G,\eta_B$ and the 1st density function $\psi_1(t)$ for the 1-period sequence $S$  whose points $0,\frac{1}{3},\frac{1}{2}$ have radii $\frac{1}{12},0,\frac{1}{12}$, see Example~\ref{exa:1st_density}. 
\textbf{Right}: The trapezium functions $\eta_{GB},\eta_{BR},\eta_{RG}$ and the 2nd density function $\psi_2(t)$ for the 1-period sequence $S$ whose points $0,\frac{1}{3},\frac{1}{2}$ have radii $\frac{1}{12},0,\frac{1}{12}$, see Example~\ref{exa:2nd_density}.
}
\label{fig:1st_density}      
\end{figure}

\begin{exa}[$\psi_1$ for the sequence $S=\{0,\frac{1}{3},\frac{1}{2}\}+\Z$]
\label{exa:1st_density}
The 1st density function $\psi_1(t)$ can be obtained as a sum of the three \emph{trapezium} functions $\eta_R$, $\eta_G$, $\eta_B$, each measuring the length of a region covered by a single interval of one color, see Fig.~\ref{fig:growing_intervals}.
\medskip

At the initial moment $t=0$, the red intervals $[0,\frac{1}{12}]\cup[\frac{11}{12},1]$ have the total length $\eta_R(0)=\frac{1}{6}$.
These red intervals $[0,\frac{1}{12}+t]\cup[\frac{11}{12}-t,1]$ for $t\in[0,\frac{1}{8}]$ grow until they touch the green interval $[\frac{7}{24},\frac{3}{8}]$ and have the total length $\eta_R(\frac{1}{8})=\frac{1}{6}+\frac{2}{8}=\frac{5}{12}$ in the second picture of Fig.~\ref{fig:growing_intervals}.
So the graph of the red length $\eta_R(t)$ linearly grows with gradient 2 from the point $(0,\frac{1}{6})$ to the corner point $(\frac{1}{8},\frac{5}{12})$.
\medskip

For $t\in[\frac{1}{8},\frac{1}{6}]$, the left red interval is shrinking at the same rate (due to the overlapping green interval) as the right red interval continues to grow until $t=\frac{1}{6}$, when it touches the blue interval $[\frac{1}{4},\frac{3}{4}]$.  
Hence the graph of $\eta_R(t)$ remains constant for $t\in[\frac{1}{8},\frac{1}{6}]$ up to the corner point $(\frac{1}{6},\frac{5}{12})$.
After that, the graph of $\eta_R(t)$ linearly decreases (with gradient $-2$) until all red intervals are fully covered by the green and blue intervals at moment $t=\frac{3}{8}$, see the 6th picture in Fig.~\ref{fig:growing_intervals}.
\medskip

Hence, the trapezium function $\eta_R$ has the piecewise linear graph through the corner points $(0,\frac{1}{6})$, $(\frac{1}{8},\frac{5}{12})$, $(\frac{1}{6},\frac{5}{12})$, $(\frac{3}{8},0)$.
After that, $\eta_R(t)=0$ remains constant for $t\geq \frac{3}{8}$.
Fig.~\ref{fig:1st_density} shows the graphs of $\eta_R,\eta_G,\eta_B$ and $\psi_1=\eta_R+\eta_G+\eta_B$. 
\eexa
\end{exa}

Theorem~\ref{thm:1st_density} extends Example~\ref{exa:1st_density} 
and proves that any $\psi_1(t)$ is a sum of trapezium functions whose corners are explicitly described. 
We consider any index $i=1,\dots,m$ (of a point $p_i$ or a gap $g_i$) modulo $m$ so that $m+1\equiv 1\pmod{m}$.

\begin{thm}[description of $\psi_1$, {\cite[Theorem~4.2]{anosova2023density}}]
\label{thm:1st_density}
Let a periodic sequence $S=\{p_1,\dots,p_m\}+\Z$ consist of disjoint intervals with centers $0\leq p_1<\dots<p_m<1$ and radii $r_1,\dots,r_m\geq 0$, respectively. 
Consider the \emph{gaps} $g_i=(p_{i}-r_{i})-(p_{i-1}+r_{i-1})$, between successive intervals, 
where $i=1,\dots,m$ and $p_{0}=p_m-1$, $r_0=r_m$. 
Then the 1st density 
$\psi_1(t)$ is the sum of $m$ \emph{trapezium} functions $\eta_{i}$, $i=1,\dots,m$, with the corners 
$(0,2r_{i})$,
$(\frac{g_{i}}{2}, g+2r_i)$,
$(\frac{g_{i+1}}{2}, g+2r_i)$,  
$(\frac{g_{i}+g_{i+1}}{2}+r_i,0),$ 
where 
$g=\min\{g_{i},g_{i+1}\}$.
Hence, $\psi_1(t)$ is determined by the unordered set of unordered pairs $(g_{i},g_{i+1})$, $i=1,\dots,m$.
\ethm 
\end{thm}

\begin{figure}[h!]
\includegraphics[width=\linewidth]{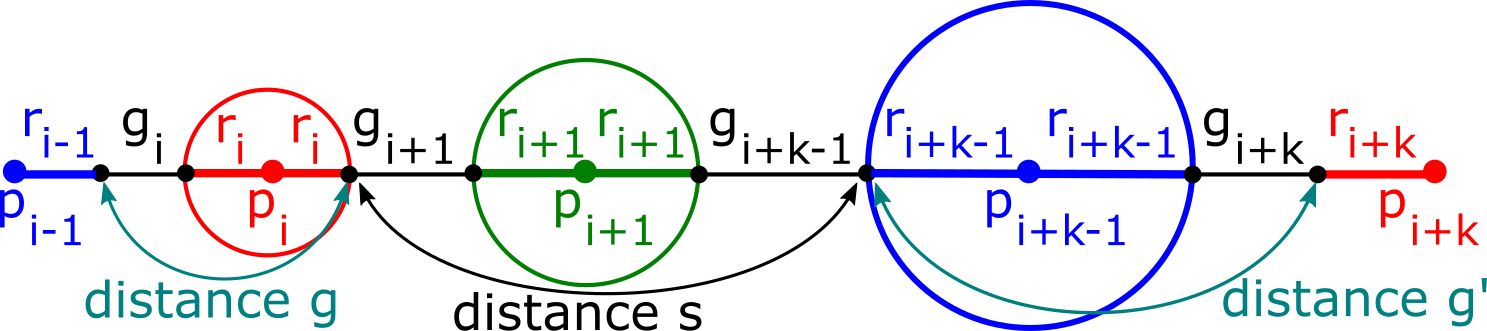}
\caption{The distances $g,s,g'$ between line intervals in Theorems~\ref{thm:1st_density} and~\ref{thm:k-th_density}, shown for $k=3$. }
\label{fig:intervals_gaps}      
\end{figure}

Example~\ref{exa:revisit_1st_density} applies Theorem~\ref{thm:1st_density} to get $\psi_1$ for the sequence $S$ in Example~\ref{exa:1st_density}.

\begin{exa}[using Theorem~\ref{thm:1st_density} for $\psi_1$]
\label{exa:revisit_1st_density}
The sequence $S=\{0,\frac{1}{3},\frac{1}{2}\}+\Z$ in Example~\ref{exa:1st_density} with points $p_1=0$, $p_2=\frac{1}{3}$, $p_3=\frac{1}{2}$ of radii $r_1=\frac{1}{12}$, $r_2=0$, $r_3=\frac{1}{12}$, respectively, has 
 the initial gaps between successive intervals 
$g_1=\frac{1}{3}$, 
$g_2=\frac{1}{4}$,
$g_3=\frac{1}{12}$, see all the computations in Example~\ref{exa:revisit_0-th_density}.
\medskip

\noindent
\textbf{Case (R)}.
In Theorem~\ref{thm:1st_density} for the trapezium function $\eta_R=\eta_1$ measuring the fractional length covered only by the red interval, we set $i=1$.
Then $r_i=\frac{1}{12}$, $g_i=\frac{1}{3}$ and $g_{i+1}=\frac{1}{4}$, 
$$\begin{array}{l}
\frac{g_i+g_{i+1}}{2}+r_i=\frac{1}{2}(\frac{1}{3}+\frac{1}{4})+\frac{1}{12}=\frac{3}{8}, \\
g=\min\{g_i,g_{i+1}\}=\frac{1}{4},\quad g+2r_i=\frac{1}{4}+\frac{2}{12}=\frac{5}{12}.
\end{array}$$ 

Then $\eta_{R}=\eta_{1}$ has the following corner points:
$$\begin{array}{l}
(0,2r_i)=(0,\frac{1}{6}), \quad
(\frac{g_{i}}{2},g+2r_i)=(\frac{1}{6},\frac{5}{12}), \\
(\frac{g_{i+1}}{2},g+2r_i)=(\frac{1}{8},\frac{5}{12}), \\ 
(\frac{g_{i}+g_{i+1}}{2}+r_i,0)=(\frac{3}{8},0),
\end{array}$$
where the two middle corners are accidentally swapped due to $g_i>g_{i+1}$ but they define the same trapezium function as in the first picture of Fig.~\ref{fig:1st_density}.
\medskip

\noindent
\textbf{Case (G)}.
In Theorem~\ref{thm:1st_density} for the trapezium function $\eta_G=\eta_2$ measuring the fractional length covered only by the green interval, we set $i=2$.
Then $r_i=0$, $g_i=\frac{1}{4}$ and $g_{i+1}=\frac{1}{12}$, 
$\begin{array}{l}
\frac{g_i+g_{i+1}}{2}+r_i=\frac{1}{2}(\frac{1}{4}+\frac{1}{12})+0=\frac{1}{6},\;
g=\min\{g_i,g_{i+1}\}=\frac{1}{12},\; g+2r_i=\frac{1}{12}+0=\frac{1}{12}.
\end{array}$
\myskip

Then $\eta_{G}=\eta_{2}$ has the following corner points
exactly as shown in the second picture of Fig.~\ref{fig:1st_density}~(left):
$$\begin{array}{ll}
(0,2r_i)=(0,0), &
(\frac{g_{i}}{2},g+2r_i)=(\frac{1}{8},\frac{1}{12}), \\
(\frac{g_{i+1}}{2},g+2r_i)=(\frac{1}{24},\frac{5}{12}), &
(\frac{g_{i}+g_{i+1}}{2}+r_i,0)=(\frac{1}{6},0).
\end{array}$$

\noindent
\textbf{Case (B)}.
In Theorem~\ref{thm:1st_density} for the trapezium function $\eta_B=\eta_3$ measuring the fractional length covered only by the blue interval, we set $i=3$.
Then $r_i=\frac{1}{12}$, $g_i=\frac{1}{12}$ and $g_{i+1}=\frac{1}{3}$, 
$$\begin{array}{l}
\frac{g_i+g_{i+1}}{2}+r_i=\frac{1}{2}(\frac{1}{12}+\frac{1}{3})+\frac{1}{12}=\frac{7}{24}, \;
g=\min\{g_i,g_{i+1}\}=\frac{1}{12},\;
g+2r_i=\frac{1}{12}+\frac{2}{12}=\frac{1}{4}.
\end{array}$$ 

Then $\eta_{B}=\eta_{3}$ has the following corner points:
$$\begin{array}{ll}
(0,2r_i)=(0,\frac{1}{6}), &
(\frac{g_{i}}{2},g+2r_i)=(\frac{1}{24},\frac{1}{4}), \\
(\frac{g_{i+1}}{2},g+2r_i)=(\frac{1}{6},\frac{1}{4}), &
(\frac{g_{i}+g_{i+1}}{2}+r_i,0)=(\frac{7}{24},0)
\end{array}$$
exactly as shown in the third picture of Fig.~\ref{fig:1st_density}.
\eexa
\end{exa}

Theorem~\ref{thm:k-th_density} describing the $k$-th density function $\psi_k[S](t)$ for any $k\geq 2$ and a periodic sequence $S$ of disjoint intervals.
To illustrate Theorem~\ref{thm:k-th_density}, Example~\ref{exa:2nd_density} computes $\psi_2[S]$ for the periodic sequence $S$ from Example~\ref{exa:0-th_density}.

\begin{exa}[$\psi_2$ for the sequence $S=\{0,\frac{1}{3},\frac{1}{2}\}+\Z$]
\label{exa:2nd_density}
The density $\psi_2(t)$ can be found as the sum of the \emph{trapezium} functions $\eta_{GB},\eta_{BR},\eta_{RG}$, each measuring the length of a double intersection, see Fig.~\ref{fig:growing_intervals}.
For the green interval $[\frac{1}{3}-t,\frac{1}{3}+t]$ and the blue interval $[\frac{5}{12}-t,\frac{7}{12}+t]$, the graph of the function $\eta_{GB}(t)$ is piecewise linear and starts at the point $(\frac{1}{24},0)$ because these intervals touch at $t=\frac{1}{24}$.
\medskip

The green-blue intersection $[\frac{5}{12}-t,\frac{1}{3}+t]$ grows until $t=\frac{1}{6}$, when the resulting interval $[\frac{1}{4},\frac{1}{2}]$ touches the red interval on the left.
At the same time, the graph of $\eta_{GB}(t)$ is linearly growing (with gradient 2) to the corner $(\frac{1}{6},\frac{1}{4})$, see Fig,~\ref{fig:1st_density}.
\medskip

For $t\in[\frac{1}{6},\frac{7}{24}]$, the green-blue intersection interval becomes shorter on the left, but grows at the same rate on the right until $t=\frac{7}{24}$ when $[\frac{1}{8},\frac{5}{8}]$ touches the red interval $[\frac{5}{8},1]$ on the right, see the 5th picture in Fig.~\ref{fig:growing_intervals}.
So the graph of $\eta_{GB}(t)$ remains constant up to the point $(\frac{7}{24},\frac{1}{4})$.
\medskip

For $t\in[\frac{7}{24},\frac{5}{12}]$ the green-blue intersection interval is shortening from both sides.
So the graph of $\eta_{GB}(t)$ linearly decreases (with gradient $-2$) and returns to the $t$-axis at the corner $(\frac{5}{12},0)$, then remains constant $\eta_{GB}(t)=0$ for $t\geq \frac{5}{12}$.
Fig.~\ref{fig:1st_density} shows all trapezium functions for double intersections and $\psi_2=\eta_{GB}+\eta_{BR}+\eta_{RG}$.
\eexa
\end{exa}

\begin{thm}[description of $\psi_k$ for $k\geq 2$, {\cite[Theorem~5.2]{anosova2023density}}]
\label{thm:k-th_density}
Let a periodic sequence $S=\{p_1,\dots,p_m\}+\Z$ consist of disjoint intervals with centers $0\leq p_1<\dots<p_m<1$ and radii $r_1,\dots,r_m\geq 0$, respectively. 
Consider the \emph{gaps} $g_i=(p_{i}-r_{i})-(p_{i-1}+r_{i-1})$ between the successive intervals of $S$, where $i=1,\dots,m$ and $p_{0}=p_m-1$, $r_0=r_m$. 
\medskip

For $k\geq 2$, the density function $\psi_k(t)$ equals the sum of $m$ \emph{trapezium} functions $\eta_{k,i}(t)$, $i=1,\dots,m$, each having the corner points 
$(\frac{s}{2},0),  
(\frac{g+s}{2},g),
(\frac{s+g'}{2},g),
(\frac{g+s+g'}{2},0),$
where 
$g,g'$ are the minimum and maximum values in the pair $\{g_{i}+2r_i,g_{i+k}+2r_{i+k-1}\}$, and $s=\sum\limits_{j=i+1}^{i+k-1}g_j+2\sum\limits_{j=i+1}^{i+k-2} r_j$.
For $k=2$, we have $s=g_{i+1}$.
Hence, $\psi_k(t)$ is determined by the unordered set of the ordered tuples $(g,s,g')$, $i=1,\dots,m$.
\ethm
\end{thm}

Example~\ref{exa:revisit_2nd_density} applies Theorem~\ref{thm:k-th_density} to get $\psi_2$ for the sequence $S$ in Example~\ref{exa:0-th_density}.

\begin{exa}[using Theorem~\ref{thm:k-th_density} for $\psi_2$]
\label{exa:revisit_2nd_density}
The sequence $S=\{0,\frac{1}{3},\frac{1}{2}\}+\Z$ in Example~\ref{exa:1st_density} with points $p_1=0$, $p_2=\frac{1}{3}$, $p_3=\frac{1}{2}$ of radii $r_1=\frac{1}{12}$, $r_2=0$, $r_3=\frac{1}{12}$, respectively, has 
 the initial gaps 
$g_1=\frac{1}{3}$, 
$g_2=\frac{1}{4}$,
$g_3=\frac{1}{12}$, see Example~\ref{exa:revisit_0-th_density}.
\medskip

In Theorem~\ref{thm:k-th_density}, the 2nd density function $\psi_2[S](t)$ is expressed as a sum of the trapezium functions computed via their corners below.
\medskip

\noindent
\textbf{Case (GB)}.
For the function $\eta_{GB}$ measuring the double intersections of the green and blue intervals centered at $p_2=p_i$ and $p_3=p_{i+k-1}$, we set $k=2$ and $i=2$.
Then we have the radii $r_i=0$ and $r_{i+1}=\frac{1}{12}$, the gaps $g_i=\frac{1}{4}$, $g_{i+1}=\frac{1}{12}$, $g_{i+2}=\frac{1}{3}$, and the sum $s=g_{i+1}=\frac{1}{12}$.
The pair 
$
\{g_i+2r_i,g_{i+2}+2r_{i+1}\}=\{\frac{1}{4}+0,\frac{1}{3}+\frac{2}{12}\}
$ 
has the minimum value $g=\frac{1}{4}$ and maximum value $g'=\frac{1}{2}$.
Then $\eta_{2,2}[S](t)=\eta_{GB}$ has the following corners as in the top picture of Fig.~\ref{fig:1st_density}~(right): 

$$\begin{array}{l}
(\frac{s}{2},0)=(\frac{1}{24},0),\\
(\frac{g+s}{2},g)
=(\frac{1}{2}(\frac{1}{4}+\frac{1}{12}),\frac{1}{4})=(\frac{1}{6},\frac{1}{4}), \\
(\frac{s+g'}{2},g)
=(\frac{1}{2}(\frac{1}{12}+\frac{1}{2}),\frac{1}{4})=(\frac{7}{24},\frac{1}{4}), \\
(\frac{g+s+g'}{2},0)
=(\frac{1}{2}(\frac{1}{4}+\frac{1}{12}+\frac{1}{2}),0)=(\frac{5}{12},0).
\end{array}$$

\noindent
\textbf{Case (BR)}.
For the trapezium function $\eta_{BR}$ measuring the double intersections of the blue and red intervals centered at $p_3=p_i$ and $p_1=p_{i+k-1}$, we set $k=2$ and $i=3$.
Then we have the radii $r_i=\frac{1}{12}=r_{i+1}$, the gaps $g_i=\frac{1}{12}$, $g_{i+1}=\frac{1}{3}$, $g_{i+2}=\frac{1}{4}$, and $s=g_{i+1}=\frac{1}{3}$.
The pair 
$
\{g_i+2r_i,g_{i+2}+2r_{i+1}\}=\{\frac{1}{12}+\frac{2}{12},\frac{1}{4}+\frac{2}{12}\}
$ 
has the minimum $g=\frac{1}{4}$ and maximum $g'=\frac{5}{12}$.
Then $\eta_{2,3}[S](t)=\eta_{BR}$ has the following corners as expected in the second picture of Fig.~\ref{fig:1st_density}~(right): 
$$\begin{array}{l}
(\frac{s}{2},0)=(\frac{1}{6},0),\\
(\frac{g+s}{2},g)
=(\frac{1}{2}(\frac{1}{4}+\frac{1}{3}),\frac{1}{4})=(\frac{7}{24},\frac{1}{4}), \\
(\frac{s+g'}{2},g)
=(\frac{1}{2}(\frac{1}{3}+\frac{5}{12}),\frac{1}{4})=(\frac{3}{8},\frac{1}{4}), \\
(\frac{g+s+g'}{2},0)
=(\frac{1}{2}(\frac{1}{4}+\frac{1}{3}+\frac{5}{12}),0)=(\frac{1}{2},0).
\end{array}$$

\noindent
\textbf{Case (RG)}.
For the trapezium function $\eta_{RG}$ measuring the double intersections of the red and green intervals centered at $p_1=p_i$ and $p_2=p_{i+k-1}$, we set $k=2$ and $i=1$.
Then we have the radii $r_i=\frac{1}{12}$ and $r_{i+1}=0$, the gaps $g_i=\frac{1}{3}$, $g_{i+1}=\frac{1}{4}$, $g_{i+2}=\frac{1}{12}$, and $s=g_{i+1}=\frac{1}{4}$.
The pair 
$
\{g_i+2r_i,g_{i+2}+2r_{i+1}\}=\{\frac{1}{3}+\frac{2}{12},\frac{1}{12}+0\}
$ 
has the minimum $g=\frac{1}{12}$ and maximum $g'=\frac{1}{2}$.
Then $\eta_{2,1}[S](t)=\eta_{RG}$ has the following corners: 
$$\begin{array}{l}
(\frac{s}{2},0)=(\frac{1}{8},0),\\
(\frac{g+s}{2},g)
=(\frac{1}{2}(\frac{1}{12}+\frac{1}{4}),\frac{1}{12})=(\frac{1}{6},\frac{1}{12}), \\
(\frac{s+g'}{2},g)
=(\frac{1}{2}(\frac{1}{4}+\frac{1}{2}),\frac{1}{12})=(\frac{3}{8},\frac{1}{12}), \\
(\frac{g+s+g'}{2},0)
=(\frac{1}{2}(\frac{1}{12}+\frac{1}{4}+\frac{1}{2}),0)=(\frac{5}{12},0).
\end{array}$$
 as expected in the third picture of Fig.~\ref{fig:1st_density}~(right).
\eexa
\end{exa}

\begin{figure}[h!]
\includegraphics[width=\linewidth]{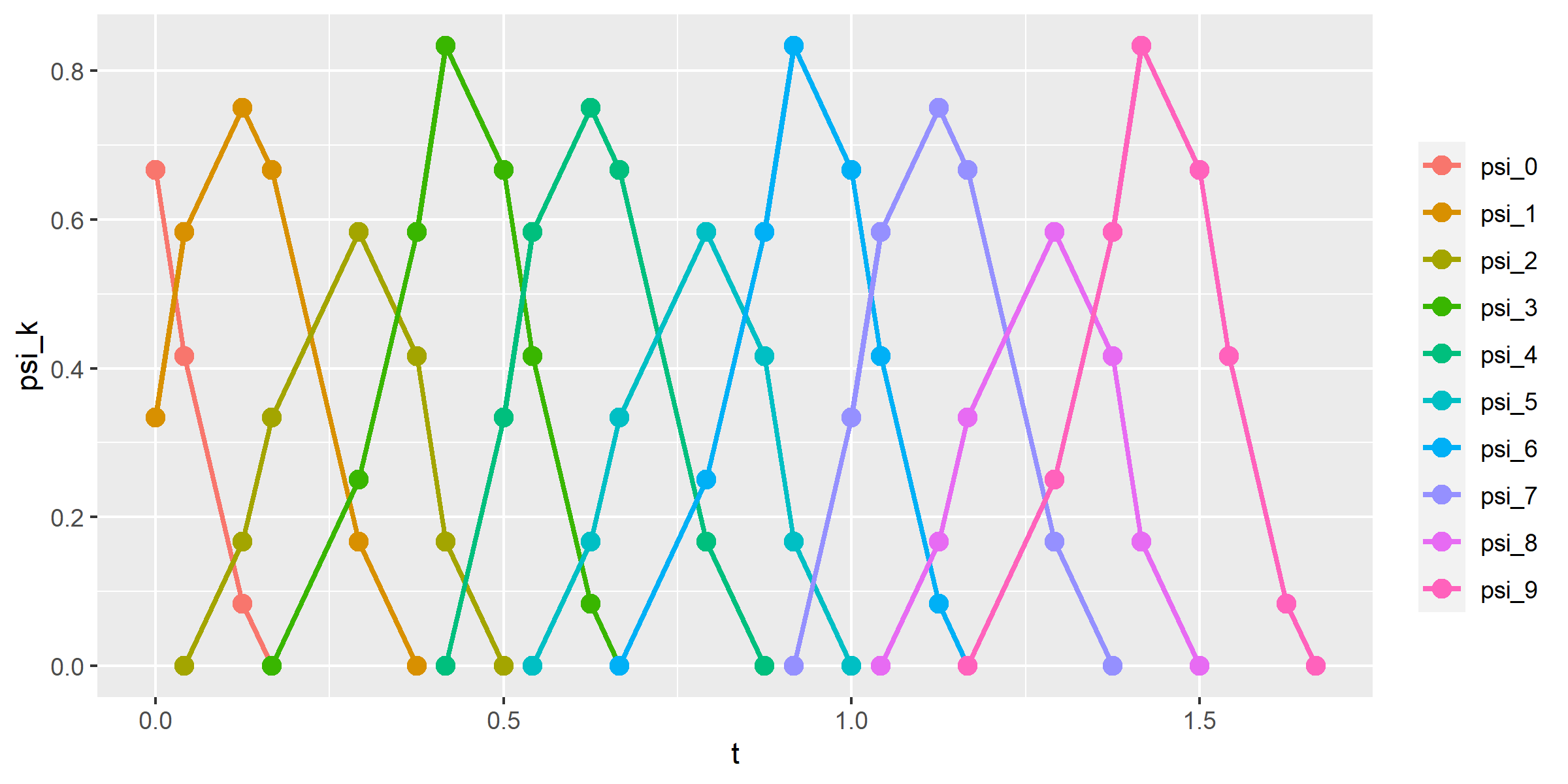}
\caption{The densities $\psi_k$, $k=0,\dots,9$ for the 1-period sequence $S$ whose points $0,\frac{1}{3},\frac{1}{2}$ have radii $\frac{1}{12},0,\frac{1}{12}$, respectively. 
The densities $\psi_0,\psi_1,\psi_2$ are described in Examples~\ref{exa:0-th_density},~\ref{exa:1st_density},~\ref{exa:2nd_density} and determine all other densities by periodicity in Theorem~\ref{thm:periodicity}. }
\label{fig:3-point_set_densities9}      
\end{figure}

\section{Properties of density functions of periodic sequences of intervals}
\label{sec:densities1D_properties}

All results in this section have detailed proofs in \cite[sections 6]{anosova2023density}.
Now we study the periodicity of the sequence $\{\psi_k\}$ with respect to the index $k\geq 0$ in Theorem~\ref{thm:periodicity}, which was a bit unexpected from Definition~\ref{dfn:densities_radii}.
We start with the simpler example for the familiar 3-point sequence in Fig.~\ref{fig:growing_intervals}.

\begin{figure}[H]
\includegraphics[width=\linewidth]{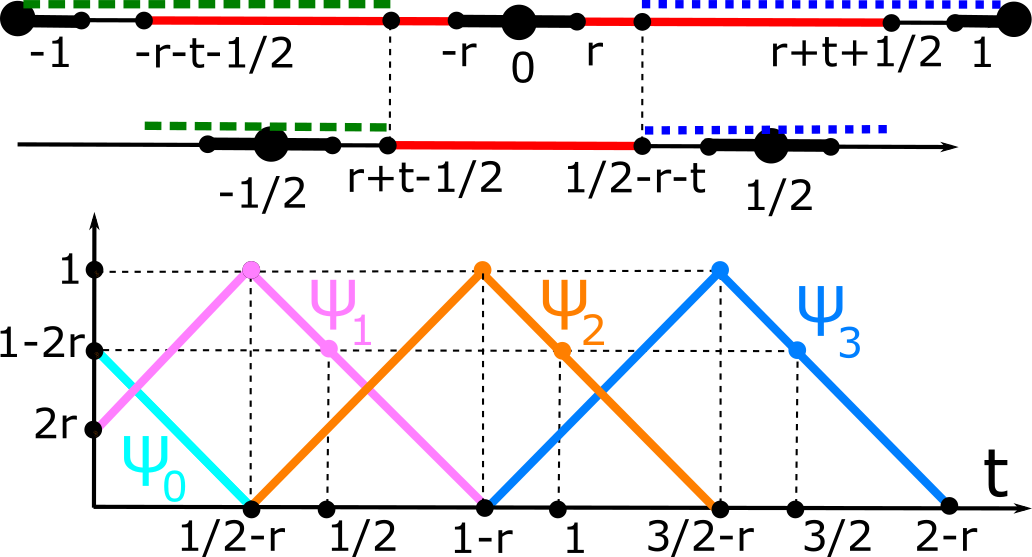}
\caption{
\textbf{Top}: Example~\ref{exa:periodicity_m=1} illustrates the proof of Theorem~\ref{thm:periodicity} for $m=1$.
\textbf{Bottom}: the density functions $\psi_k$ of $S=\Z$ whose points have a radius $0<r<\frac{1}{4}$ satisfy the periodicity $\psi_{k+1}(t+\frac{1}{2})=\psi_{k}(t)$ for any $k\geq 0$ and $t\geq 0$.} 
\label{fig:densities_1point_radius}      
\end{figure}

\begin{exa}[periodicity of $\psi_k$ in the index $k$]
\label{exa:periodicity}
Let the periodic sequence $S=\{0,\frac{1}{3},\frac{1}{2}\}+\Z$ have three points $p_1=0$, $p_2=\frac{1}{3}$, $p_3=\frac{1}{2}$ of radii $r_1=\frac{1}{12}$, $r_2=0$, $r_3=\frac{1}{12}$, respectively.
The initial intervals $L_1(0)=[-\frac{1}{12},\frac{1}{12}]$, $L_2(0)=[\frac{1}{3},\frac{1}{3}]$, $L_3(0)=[\frac{5}{12},\frac{7}{12}]$ have the 0-fold intersection measured by $\psi_0(0)=\frac{2}{3}$ and the 1-fold intersection measured by $\psi_1(0)=\frac{1}{3}$, see Fig.~\ref{fig:0-th_density} and~\ref{fig:1st_density}.
\medskip

By the time $t=\frac{1}{2}$ the initial intervals will grow to 
$L_1(\frac{1}{2})=[-\frac{7}{12},\frac{7}{12}]$, 
$L_2(\frac{1}{2})=[-\frac{1}{6},\frac{5}{6}]$, 
$L_3(\frac{1}{2})=[-\frac{1}{12},\frac{13}{12}]$.
The grown intervals at the radius $t=\frac{1}{2}$ have the 3-fold intersection $[-\frac{1}{12},\frac{7}{12}]$ of the length $\psi_3(\frac{1}{2})=\frac{2}{3}$, which coincides with $\psi_0(0)=\frac{2}{3}$.
\medskip

With the extra interval $L_4(\frac{1}{2})=[\frac{5}{12},\frac{19}{12}]$ centered at $p_4=1$, the 4-fold intersection is $L_1\cap L_2\cap L_3\cap L_4=[\frac{5}{12},\frac{7}{12}]$.
With the extra interval $L_{5}(\frac{1}{2})=[\frac{5}{6},\frac{11}{6}]$ centered at $p_{5}=\frac{4}{3}$, the 4-fold intersection $L_2\cap L_3\cap L_4\cap L_5$ is the single point $\frac{5}{6}$.
With the extra interval $L_{6}(\frac{1}{2})=[\frac{11}{12},\frac{13}{12}]$ centered at $p_6=\frac{3}{2}$, the 4-fold intersection is 
$L_3\cap L_4\cap L_5\cap L_6=[\frac{11}{12},\frac{13}{12}]$.
Hence the total length of the 4-fold intersection at $t=\frac{1}{2}$ is $\psi_4(\frac{1}{2})=\frac{1}{3}$, which coincides with $\psi_1(0)=\frac{1}{3}$.
\medskip

For the larger $t=1$, the six grown intervals
$$\begin{array}{ll}
L_1(1)=\left[-\frac{13}{12},\frac{13}{12}\right], & 
L_2(1)=\left[-\frac{2}{3},\frac{4}{3}\right], \\
L_3(1)=\left[-\frac{7}{12},\frac{19}{12}\right], & 
L_4(1)=\left[-\frac{1}{12},\frac{25}{12}\right], \\
L_5(1)=\left[\frac{1}{3},\frac{7}{3}\right], &
L_6(1)=\left[\frac{5}{12},\frac{31}{12}\right]
\end{array}$$

have the 6-fold intersection 
$\left[\frac{5}{12},\frac{13}{12}\right]$ of length $\psi_6(1)=\frac{2}{3}$ coinciding with $\psi_0(0)=\psi_3(\frac{1}{2})=\frac{2}{3}$.
\eexa 
\end{exa}

Corollary~\ref{thm:periodicity} says that the coincidences in Example~\ref{exa:periodicity} are not accidental. 
The periodicity of $\psi_k$ with respect to $k$ is illustrated by Fig.~\ref{fig:3-point_set_densities9}.

\begin{thm}[periodicity of $\psi_k$ in the index $k$, {\cite[Theorem~6.2]{anosova2023density}}]
\label{thm:periodicity}
The density functions $\psi_k[S]$ of a periodic sequence $S=\{p_1,\dots,p_m\}+\Z$ consisting of disjoint intervals with centers $0\leq p_1<\dots<p_m<1$ and radii $r_1,\dots,r_m\geq 0$, respectively, satisfy the \emph{periodicity} 
$\psi_{k+m}(t+\frac{1}{2})=\psi_{k}(t)$ for any $k\geq 0$ and $t\geq 0$.
\ethm
\end{thm}

\begin{figure}[H]
\includegraphics[width=\linewidth]{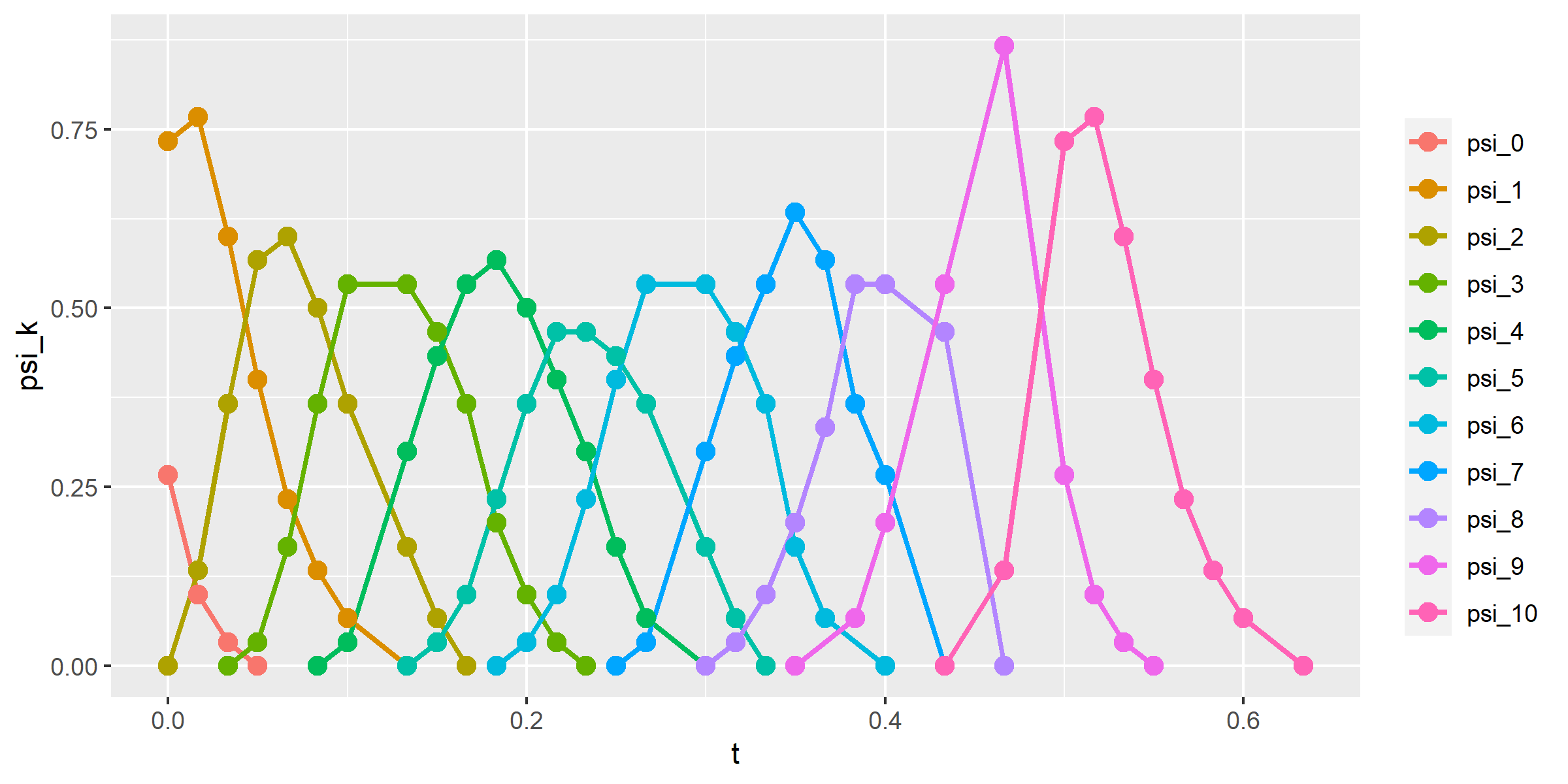}
\includegraphics[width=\linewidth]{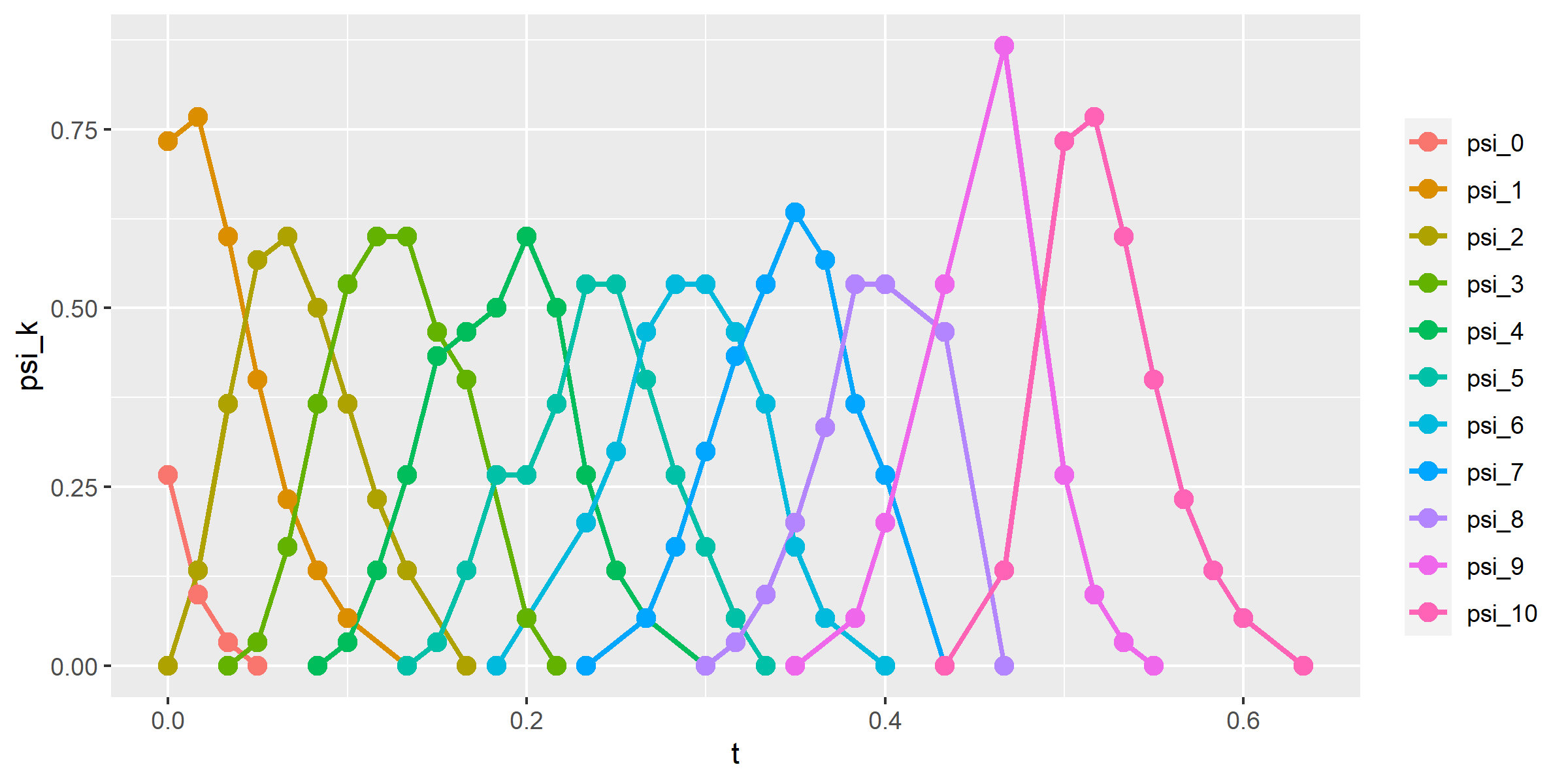}
\caption{The densities $\psi_k$, $k=0,\dots,10$, distinguish (already for $k\geq 2$) the sequences (scaled down by period $15$) $S_{15} = \{0,1,3,4,5,7,9,10,12\}+15\Z$ (\textbf{top}) and $Q_{15} = \{0,1,3,4,6,8,9,12,14\}+15\Z$ (\textbf{bottom}), where the radius $r_i$ of any point is the half-distance to its closest neighbour, see Example~\ref{exa:SQ15}.}
\label{fig:S5+Q15}      
\end{figure}

\begin{exa}[Theorem~\ref{thm:periodicity} for $m=1$ in Fig.~\ref{fig:densities_1point_radius}]
\label{exa:periodicity_m=1}
Let a 1-period sequence $S$ have one point $p_1=0$ of a radius $0<r<\frac{1}{2}$.
The grown interval $[-r-t-\frac{1}{2},r+t+\frac{1}{2}]$ around $0$ has the 1-fold intersection $I=[r+t-\frac{1}{2},\frac{1}{2}-r-t]$ centered at $p=0$ and not covered by the adjacent intervals centered at $\pm 1$, so $\psi_1(t+\frac{1}{2})=1-2(t+r)$.
\medskip

After collapsing $[-\frac{1}{2},\frac{1}{2}]$ to $0$, which is excluded from $S$, the periodic sequence has new points $\pm\frac{1}{2}$ of the smaller radius $r+t$.
The new shorter intervals have the same endpoints $-\frac{1}{2}+(r+t)$ and $\frac{1}{2}-(r+t)$ around $p=0$.
Now $I=[r+t-\frac{1}{2},\frac{1}{2}-r-t]$ is not covered by any
shorter intervals, so the get the same length of the 0-fold intersection: $\psi_0(t)=1-2(t+r)$.
\eexa
\end{exa}

The symmetry $\psi_{m-k}(\frac{1}{2}-t)=\psi_k(t)$ for $k=0,\dots,[\frac{m}{2}]$, and $t\in[0,\frac{1}{2}]$ from \cite[Theorem~8]{anosova2022density} no longer holds for points with different radii.
For example, $\psi_1(t)\neq \psi_2(\frac{1}{2}-t)$ for the periodic sequence $S=\{0,\frac{1}{3},\frac{1}{2}\}+\Z$, see Fig.~\ref{fig:1st_density}.
If all points have the same radius $r$, \cite[Theorem~8]{anosova2022density} implies the symmetry after replacing $t$ by $t+2r$.
\medskip

Example~\ref{exa:S15+Q15} justified that all density functions cannot distinguish the non-isometric sequences 
$S_{15} = \{0,1,3,4,5,7,9,10,12\}+15\Z$ and $Q_{15} = \{0,1,3,4,6,8,9,12,14\}+15\Z$ of points with zero radii.
Example~\ref{exa:S15+Q15} shows that the densities for sequences with non-zero radii are strictly stronger and distinguish the sequences $S_{15}\not\cong Q_{15}$.

\begin{exa}[$\psi_k$ for $S_{15},Q_{15}$ with radii]
\label{exa:S15+Q15}
For any point $p$ in a periodic sequence $S\subset\R$, define its \emph{neighbour} radius as the half-distance to a closest neighbour of $p$ within the sequence $S$.
This choice of radii respects the isometry in the sense that periodic sequences $S,Q$ with zero-sized radii are isometric if and only if $S,Q$ with neighbour radii are isometric.
Fig.~\ref{fig:S5+Q15} shows that the densities $\psi_k$ for $k\geq 2$ distinguish the non-isometric sequences $S_{15}$ and $Q_{15}$ scaled down by factor 15 to the unit cell $[0,1]$, see 
Example~\ref{exa:SQ15}.
\eexa
\end{exa}

\begin{figure}[H]
\includegraphics[width=\linewidth]{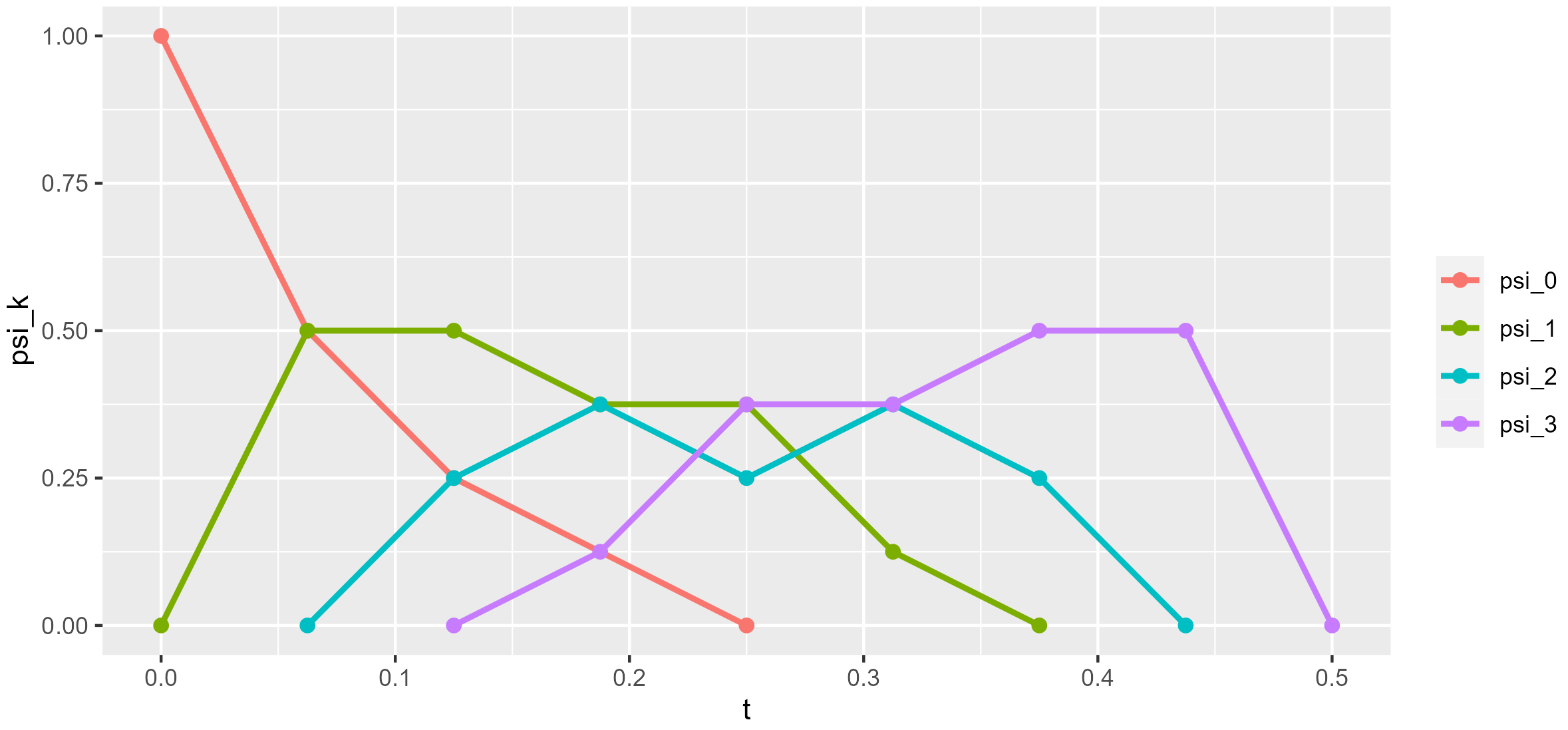}
\caption{For the periodic sequence $S=\{0,\frac{1}{8},\frac{1}{4},\frac{3}{4}\}+\Z$ whose all points have radii 0, the 2nd density $\psi_2[S](t)$ has the local minimum at $t=\frac{1}{4}$ between two local maxima.
}
\label{fig:set0_0125_025_075_densities3}      
\end{figure}

\begin{cor}[time of $\psi_k(t)$, {\cite[Corollary~6.5]{anosova2023density}}]
\label{cor:densities1D_intervals_computation}
Let $S,Q\subset\R$ be periodic sequences with at most $m$ motif points.
For $k\geq 1$, one can draw the graph of the $k$-th density function $\psi_k[S]$ in time $O(m^2)$.
One can check in time $O(m^3)$ if 
$\Psi[S]=\Psi[Q]$.
\ethm
\end{cor}

All previous examples show densities with a single local maximum.
However, the new R code \cite{anosova2023R} helped us discover the opposite examples.

\begin{figure}[H]
\includegraphics[width=\linewidth]{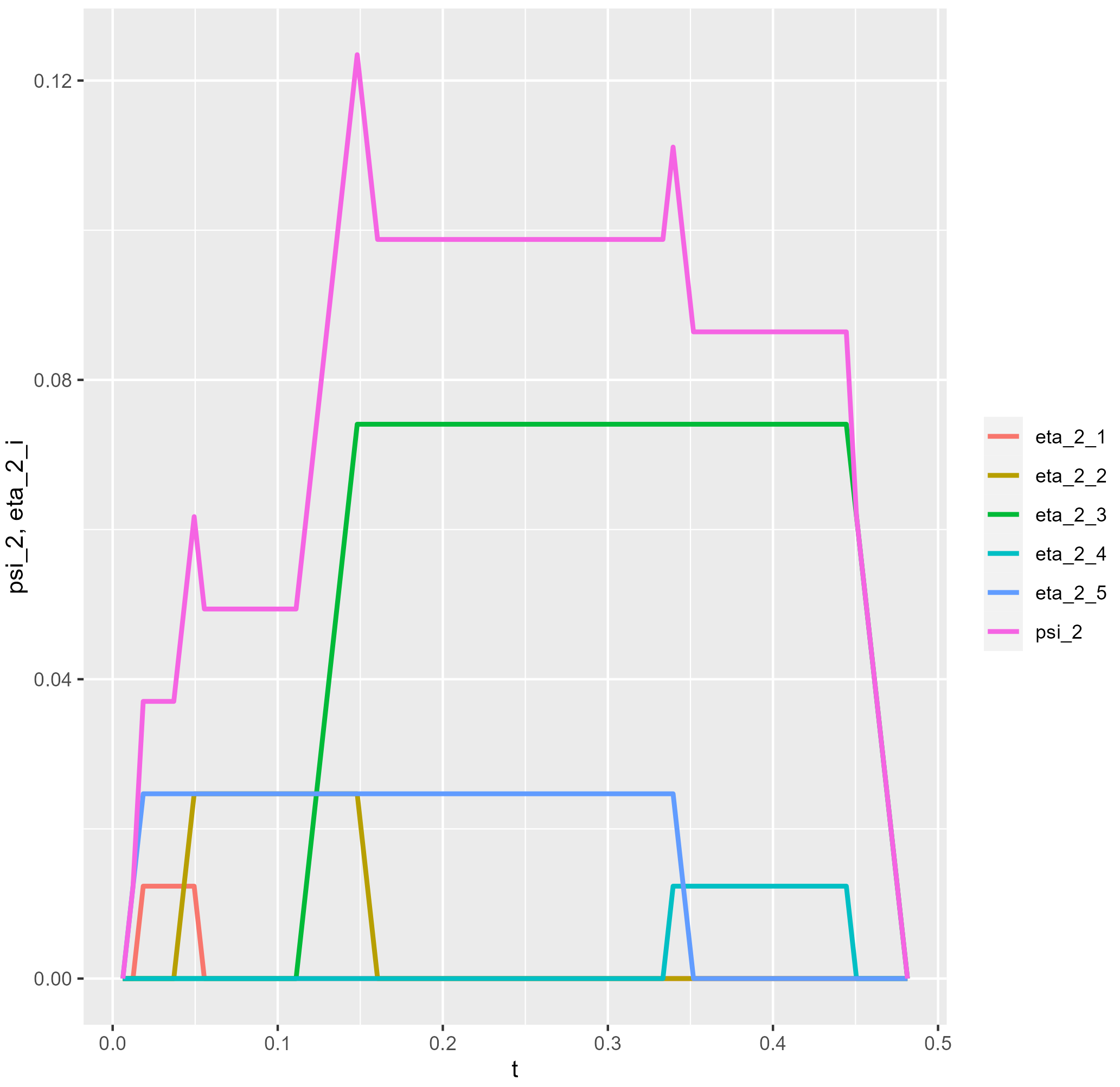}
\caption{For the sequence $S=\Big\{0,\frac{1}{81},\frac{1}{27},\frac{1}{9},\frac{1}{3}\Big\}+\Z$ whose all points have radii 0, $\psi_2[S]$ 
equal to the sum of the shown five trapezium functions has three maxima.
}
\label{fig:set0_powers3_psi_2_eta}      
\end{figure}

\begin{figure}[H]
\includegraphics[width=\linewidth]{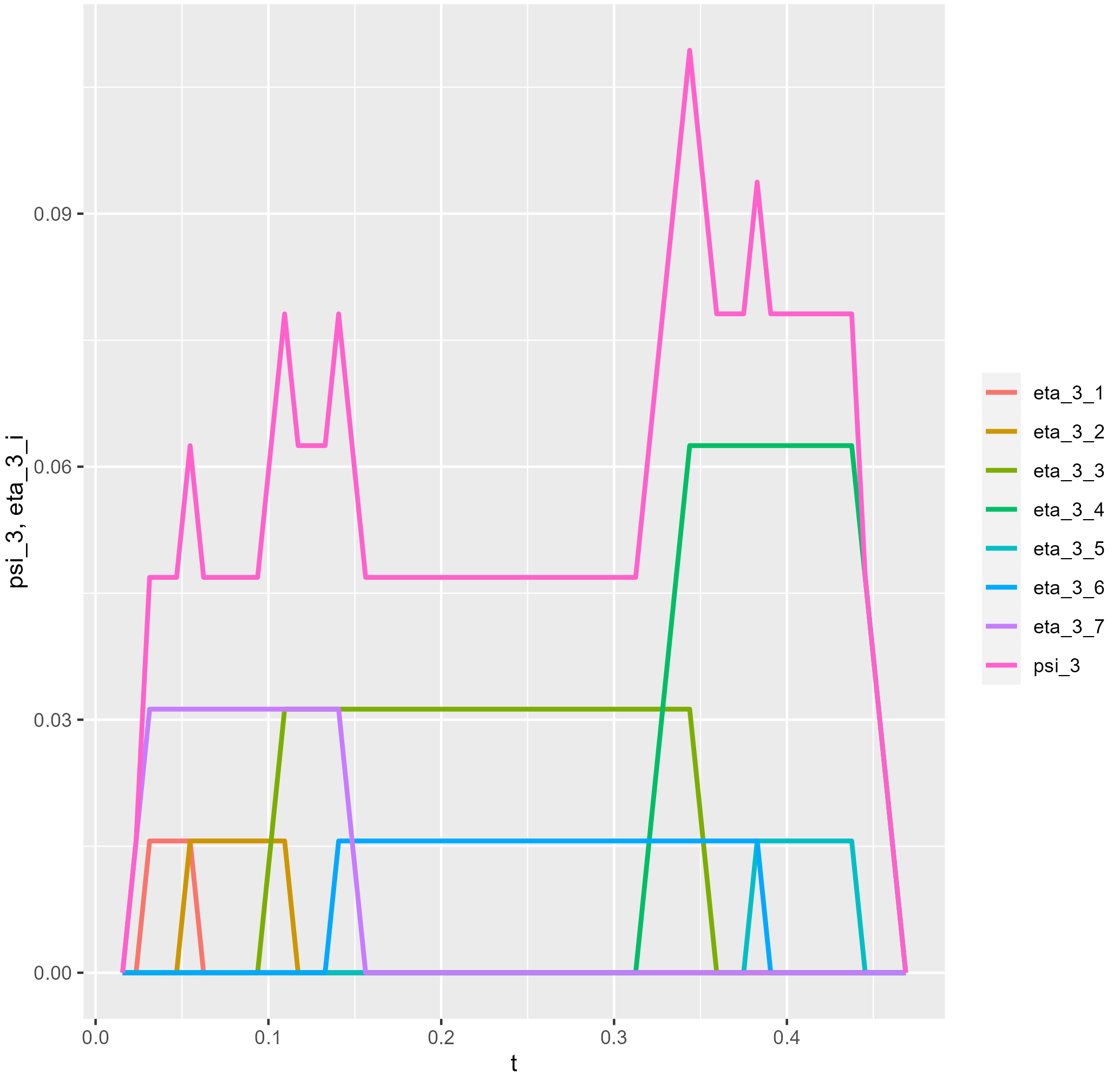}
\caption{For the sequence $S=\Big\{0,\frac{1}{64},\frac{1}{16},\frac{1}{8},\frac{1}{4},\frac{3}{4}\Big\}+\Z$ whose all points have radii 0, 
$\psi_3[S]$ has 5 local maxima.
}
\label{fig:set0_powers2_psi_3_eta}      
\end{figure}

\begin{figure}[H]
\includegraphics[width=\linewidth]{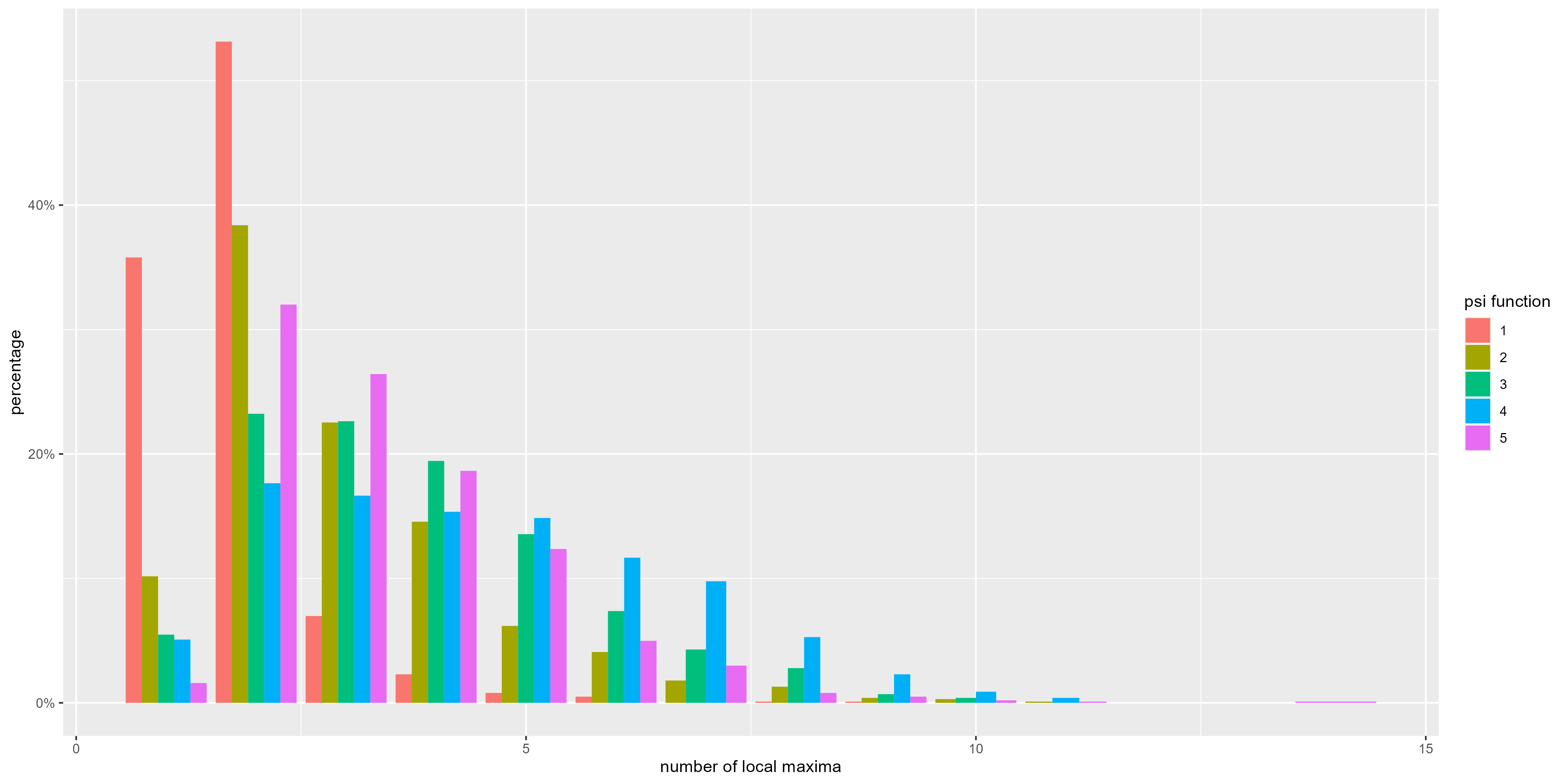}
\caption{Percentages of cases when the density functions $\psi_k(t)$, $k=1,\dots,5$ (shown in five different colors) have one or multiple local maxima for 1000 sequences of 10 points with zero radii, which are uniformly sampled in the internal $[0,1]$.}
\label{fig:histogram_maxima}      
\end{figure}

\begin{exa}[densities with multiple maxima]
\label{exa:multimax}
Fig.~\ref{fig:set0_0125_025_075_densities3} shows a simple 4-point sequence $S$ whose 2nd density $\psi_2[S]$ has two local maxima.
Figs.~\ref{fig:set0_powers3_psi_2_eta} and~\ref{fig:set0_powers2_psi_3_eta} show complicated sequences whose density functions have more than two maxima.
Fig.~\ref{fig:histogram_maxima} shows that two local maxima are more common than one maximum. 
\eexa
\end{exa}

\begin{pro}[density functions for 2D lattices]
\label{pro:lattices2D_density}
Find an analytic description of all density functions $\{\psi_k[\La]\}_{k=0}^+\infty$ for any lattices $\La\subset\R^2$ in terms of the root invariant $\RI(\La)$ from Definition~\ref{dfn:RI}.
\epro
\end{pro}

\bibliographystyle{plain}
\bibliography{Geometric-Data-Science-book}

%
%
%

\chapter{Pointwise Distance Distributions of periodic point sets in $\R^n$}
\label{chap:PDD-periodic} 

\abstract{
This chapter extends the Pointwise Distance Distribution (PDD) from the case of finite clouds of unordered points to arbitrary periodic point sets.
We prove that the PDD is Lipschitz continuous and generically complete for periodic point sets under isometry in $\R^n$.
The PDD is computable in a near-linear asymptotic time of key input sizes and detects numerous near-duplicates among about 2 million crystals in major materials databases within two hours on a modest desktop computer.
}

\section{Pointwise Distance Distributions for lattices and $l$-periodic sets
}
\label{sec:l-periodic}

This section follows papers \cite{widdowson2025pointwise,widdowson2022resolving,widdowson2025geographic} with minor updates.
Definition~\ref{dfn:l-periodic} extends 1-periodic point sets from Definition~\ref{dfn:1-periodic} and fully periodic sets in $\R^n$ from Definition~\ref{dfn:periodic} to the more general sets that are periodic in $l$ directions for $1\leq l\leq n$.

\begin{dfn}[$l$-periodic point set in $\R^n$] 
\label{dfn:l-periodic}
Let vectors $\vec v_1,\dots,\vec v_n\in\R^n$ form a basis of $\R^n$,  define the lattice $\La=\{\sum\limits_{i=1}^l c_i\vec  v_i \mid c_1,\dots,c_l\in\Z\}$
Fix any $1\leq l\leq n$.
The \emph{unit cell} defined by $\vec v_1,\dots,\vec v_n$ is $U=\{\sum\limits_{i=1}^n x_i\vec  v_i \mid x_1,\dots,x_l\in[0,1), x_{l+1},\dots,x_n\in\R\}\subset\R^n$.
\sskip

If $l=n$, then $U$ is an $n$-dimensional parallelepiped.
If $l<n$, then $U$ is an infinite slab over an $l$-dimensional parallelepiped on $\vec v_1,\dots,\vec v_l$.   
For any finite \emph{motif} of points $M\subset U$, the sum $S=M+\La=\{\vec p+\vec v \mid p\in M, v\in\La\}$ is an \emph{$l$-periodic point set}.
\edfn
\end{dfn}

\begin{figure}[h!]
\includegraphics[height=19.5mm]{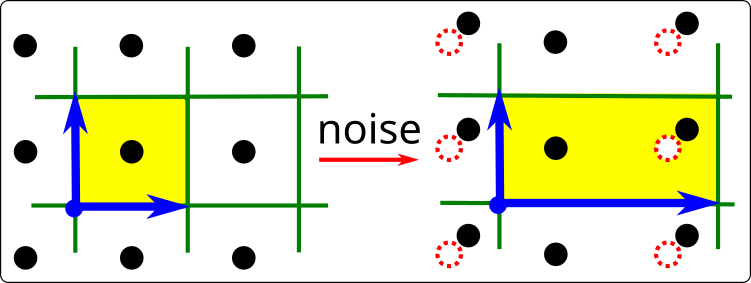}
\hspace*{0.5mm}
\includegraphics[height=19.5mm]{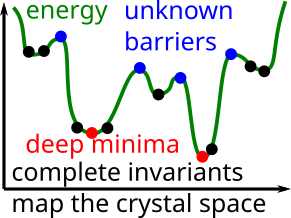}
\hspace*{0.5mm}
\includegraphics[height=19.5mm]{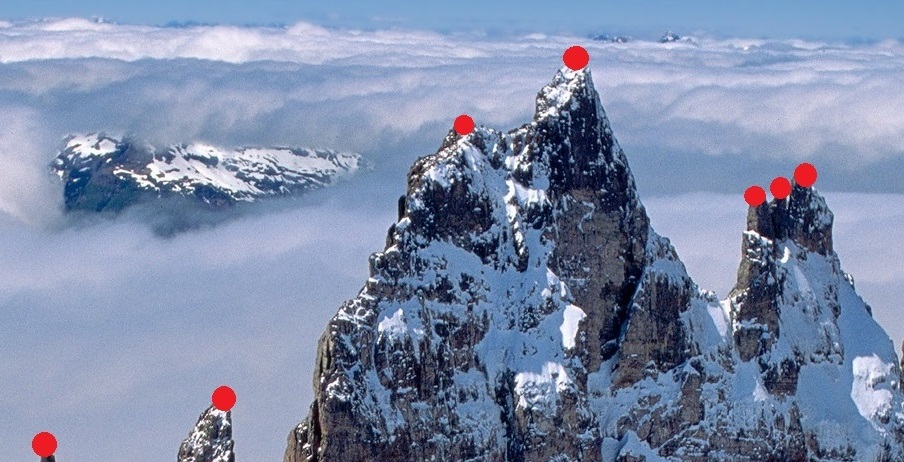}
\caption{
\textbf{Left}: the symmetry group and a reduced cell discontinuously change under tiny noise.
\textbf{Middle}: the most stable crystals are deep minima surrounded by high energy barriers in the crystal space.
\textbf{Right}: energy landscapes of crystals show optimised structures as isolated peaks of height$=-$energy.
To see beyond the `fog', we need an explicit \emph{geomap} parametrised by geocodes in Problem~\ref{pro:l-periodic}.}
\label{fig:noise_double_peaks}
\end{figure}

Fig.~\ref{fig:noise_double_peaks}~(left) illustrates the discontinuity of cell-based representations under noise.
Definition~\ref{dfn:crystal_spaces} of generic periodic sets and crystal spaces under equivalences extends to $l$-periodic point sets with motifs of $m$ points in $\R^n$.
We use similar notations and add $l$ as an extra parameter.
This chapter studies the \emph{Crystal Isometry Space} $\CIMS(\R^n;m,l)$ of $l$-periodic sets with motifs of up to $m$ points under isometry in $\R^n$.
\myskip

Problem~\ref{pro:l-periodic} extends
Problem~\ref{pro:periodic_density} to $l$-periodic point sets and add the extra condition of realisability in \ref{pro:l-periodic}(b).
This realisability condition was unrealistic for density functions from the previous chapter, at least in dimensions $n\geq 2$, but will be achieved for generic periodic sets with $l=n$ by Pointwise Distance Distributions.

\begin{pro}[isometry invariants of $l$-periodic point sets in $\R^{n}$]
\label{pro:l-periodic}
Design an invariant $I$ on the Crystal Isometry Space $\CIMS(\R^n;m,l)$
satisfying the following conditions.
\smallskip

\noindent
\tb{(a)} 
\emph{Generic completeness:}
let $S,Q$ be any \emph{generic} sets whose isometry classes are in a dense subspace of $\CIMS(\R^n;m,l)$, then $S,Q$ are isometric if and only if $I(S)=I(Q)$.
\myskip

\noindent
\tb{(b)} 
\emph{Reconstruction:}
any generic periodic point set $S\subset\R^n$ is reconstructable from its invariant 
$I(S)$, a lattice $\La$ of $S$ and the motif size $m$, uniquely under isometry in $\R^n$. 
\myskip

\noindent
\tb{(c)} 
\emph{Metric:} 
there is a distance $d$ on the Crystal Isometry Space $\CIMS(\R^n;m,l)$ satisfying all metric axioms in Definition~\ref{dfn:metrics}(a). 
\myskip

\noindent
\tb{(d)} 
\emph{Continuity:} 
there is a constant $\la>0$, such that, for all sufficiently small $\ep>0$, if a periodic point set $Q$ is obtained by perturbing every point of a periodic point set $S\subset\R^n$ up to Euclidean distance $\ep$, then $d(I(S),I(Q))\leq\la\ep$.
\myskip

\noindent
\tb{(e)} 
\emph{Computability:} 
for fixed $1\leq l\leq n$, the invariant $I(S)$, a reconstruction of $S\subset\R^{n}$ from $I(S)$, and the metric $d(I(S),I(Q))$ can be computed in times that depend polynomially on the maximum motif size of $l$-periodic point sets $S,Q$.
\epro
\end{pro}

Definition~\ref{dfn:PDD} introduces our main invariant $\PDD$ and its average $\AMD$. 

\index{$l$-periodic point set}

\begin{dfn}[$\PDD(S;k)$ and $\AMD(S;k)$ of any $l$-periodic set $S$]
\label{dfn:PDD}
Let $M=\{p_1,\dots,p_m\}$ be a motif of any $l$-periodic point set $S\subset\R^n$.
Fix an integer $k\geq 1$.
For every point $p\in M$, let $d_1(p)\leq\dots\leq d_k(p)$ be the distances from $p$ to its $k$ nearest neighbours within the full set $S$ (not restricted to $M$).
The matrix $D(S;k)$ has $m$ rows consisting of the distances $d_1(p_i),\dots,d_k(p_i)$ for $i=1,\dots,m$.
If any $l\geq 2$ rows coincide, we collapse them into a single row with the weight $l/m$.
The resulting unordered set (written as a matrix) of maximum $m$ rows and 
$k+1$ columns, including the extra column of weights, is the \emph{Pointwise Distance Distribution} $\PDD(S;k)$. 
The \emph{Average Minimum Distance} $\AMD_i$ is the weighted average of the $i$-th column in $\PDD(S;k)$ for each $i=1,\dots,k$.
Let $\AMD(S;k)$ denote the vector $(\AMD_1,\dots,\AMD_k)$. 
\edfn
\end{dfn}

Theorem~\ref{thm:PDD_periodic_invariance} shows that $\PDD$ and hence $\AMD$ are independent of a motif $M\subset S$, so there is no need to include a motif $M$ in the notation $\PDD(S;k)$.

\begin{thm}[isometry invariance of $\PDD$, {\cite[Theorem~3.3(b)]{widdowson2025pointwise}}]
\label{thm:PDD_periodic_invariance}
For any $l$-periodic point set $S\subset\R^n$, where $1\leq l\leq n$, 
$\PDD(S;k)$ and $\AMD(S;k)$ are invariants of $S$ (independent of a motif $M\subset S$) under isometry of $\R^n$ for $k\geq 1$.
\ethm
\end{thm}


\section{Asymptotic and time of Pointwise Distance Distributions}
\label{sec:PDD_asymptotic}

This section analyses the behaviour of $\PDD(S;k)$ as the index of neighbours $k\to+\infty$.
Definition~\ref{dfn:PPC} will help describe this asymptotic in Theorem~\ref{thm:asymptotic}.

\index{cell-periodic point set}
\index{Point Packing Coefficient}

\begin{dfn}[Point Packing Coefficient $\PPC$ of a cell-periodic set $S$]
\label{dfn:PPC}
\tb{(a)}
For $1\leq l\leq n$ and a basis $\vec v_1,\dots,\vec v_n\in\R^n$, consider the
\emph{lattice} $\La=\{\sum\limits_{i=1}^l c_i\vec  v_i \mid c_1,\dots,c_l\in\Z\}$ and the unit cell $U=\{\sum\limits_{i=1}^n x_i\vec  v_i \mid x_1,\dots,x_l\in[0,1), x_{l+1},\dots,x_n\in\R\}$.
A set $S\subset\R^n$ is \emph{cell-periodic} if $S$ has a fixed number $m$ points in every shifted cell $U+\vec v$ for all $\vec v\in\La$. 
\myskip

\nt
\tb{(b)}
If $l<n$, let $R^l\subset\R^n$ be the subspace spanned by $\vec v_1,\dots,\vec v_l$, then $U$ is an infinite slab based on the $l$-dimensional parallelepiped of volume $\vol[U\cap R^l]$ .
The volume of the unit ball in $\R^l$ is $V_l=\dfrac{\pi^{l/2}}{\Ga(\frac{l}{2}+1)}$, where Euler's Gamma function 
is 
$\Ga(m)=(m-1)!$ and $\Ga(\frac{m}{2}+1)=\sqrt{\pi}(m-\frac{1}{2})(m-\frac{3}{2})\cdots\frac{1}{2}$ for any integer $m\geq 1$.
Define the \emph{Point Packing Coefficient} of the cell-periodic set $S$ as $\PPC(S)=\sqrt[l]{\dfrac{\vol[U\cap R^l]}{mV_l}}$.
\end{dfn}

Any $l$-periodic set is cell-periodic, but all cell-periodic sets form a wider collection of Delone sets and model disordered solid materials that can have an underlying lattice with atoms at different positions in periodically translated cells $U+\vec v$, see Fig.~\ref{fig:lattice_periodic_set_hierarchy}.

\begin{thm}[asymptotic of $\PDD(S;k)$ as $k\to+\infty$, {\cite[Theorem~3.7]{widdowson2025pointwise}}]
\label{thm:asymptotic}
For a point $p$ in a cell-periodic set $S\subset\R^n$, let $d_k(S;p)$ be the distance from $p$ to its $k$-th nearest neighbour in $S$.
Then  
$\lim\limits_{k\to+\infty}\dfrac{d_k(S;p)}{\sqrt[l]{k}}=\PPC(S)$ and
$\lim\limits_{k\to+\infty}\dfrac{\AMD_k(S)}{\sqrt[l]{k}}=\PPC(S)$.
\ethm
\end{thm}

By Theorem~\ref{thm:asymptotic}, $\AMD_k(S)$ and all distances in the last column of $\PDD(S;k)$ asymptotically approach $\PPC(S)\sqrt[l]{k}$ as $k\to+\infty$ and hence are mainly determined by $\PPC(S)$ for large $k$.
That is why the most descriptive information is contained in $\PDD(S;k)$ for smaller values of $k$, e.g. we use $k=100$ atomic neighbours in most experiments on crystals.
To neutralise the asymptotic growth, we subtract and also normalise by the term $\PPC(S)\sqrt[l]{k}$ to get simpler invariants under uniform scaling.

\index{Average Deviation from Asymptotic}
\index{Pointwise Deviation from Asymptotic}
\index{Average Normalised Deviation}
\index{Pointwise Normalised Deviation}

\begin{dfn}[simplified invariants $\ADA$, $\PDA,\AND$, $\PND$]
\label{dfn:ADA}
Let $S\subset\R^n$ be any $l$-periodic set with an underlying lattice generated by $l$ vectors. 
\myskip

\nt
\tb{(a)}
The \emph{Average Deviation from Asymptotic} is $\ADA_k(S)=\AMD_k(S)-\PPC(S)\sqrt[l]{k}$, $k\geq 1$.
To get the \emph{Pointwise Deviation from Asymptotic} $\PDA(S;k)$ from $\PDD(S;k)$ subtracting $\PPC(S)\sqrt[l]{j}$ from each distance in a row $i$ and a column $j$ for $i\geq 1\leq j\leq k$.
\myskip

\nt
\tb{(b)}
The \emph{Average Normalised Deviation} is $\AND_k(S)=\ADA_k(S)/(\PPC(S)\sqrt[l]{k})$, $k\geq 1$.
The \emph{Pointwise Normalised Deviation} $\PND(S;k)$ obtained from $\PDA(S;k)$ by dividing every element in a row $i$ and a column $j$ by $\PPC(S)\sqrt[l]{j}$ for $i\geq 1\leq j\leq k$.
\edfn
\end{dfn} 

\begin{figure}[h!]
\caption{{The average invariants $\AMD_k$ and $\ADA_k$ from Definition}~\ref{dfn:ADA} {for $k=1,\dots,25$ and five simple crystals from the Materials Project.}}
\includegraphics[width=0.49\textwidth]{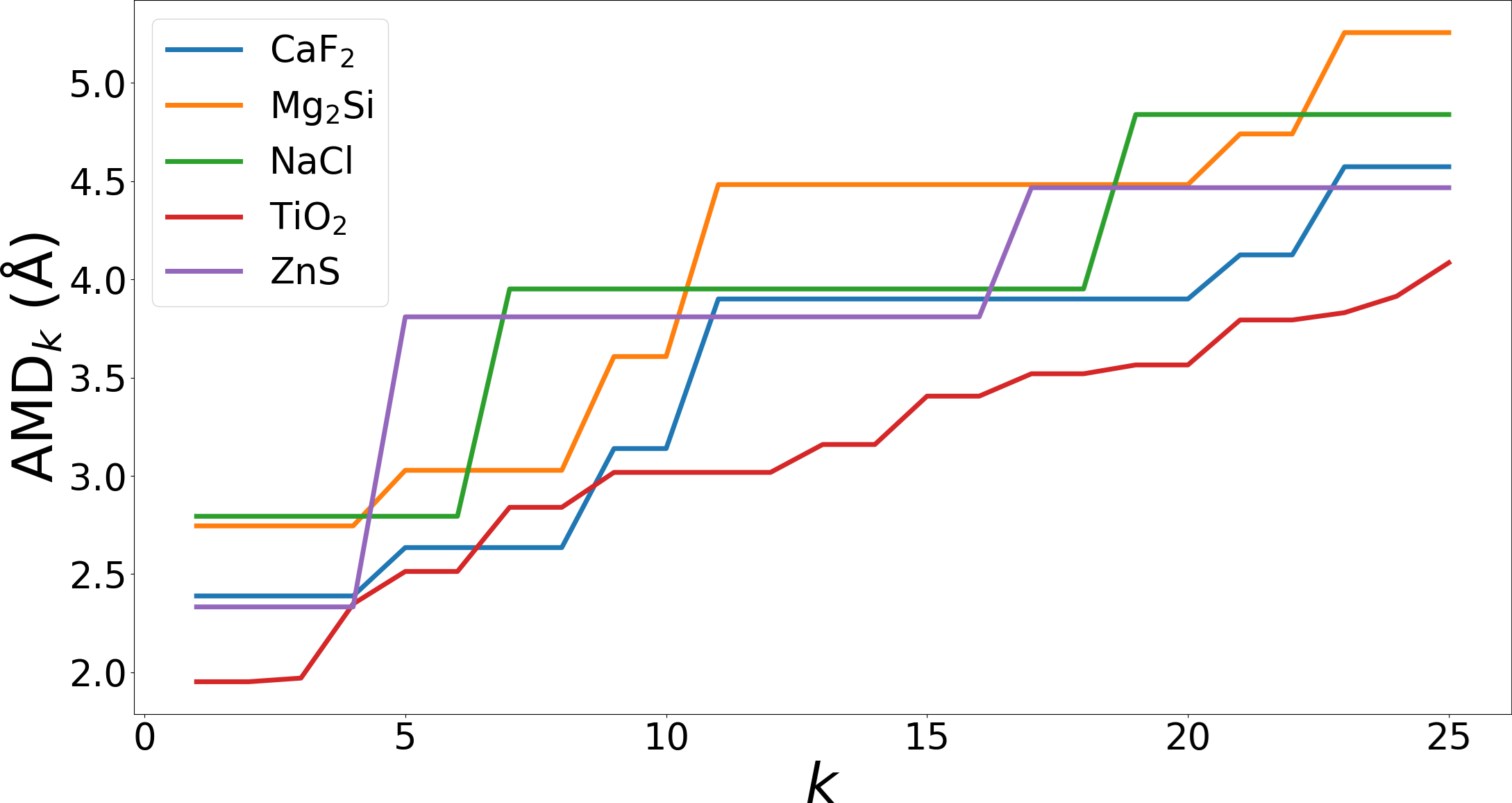}
\includegraphics[width=0.49\textwidth]{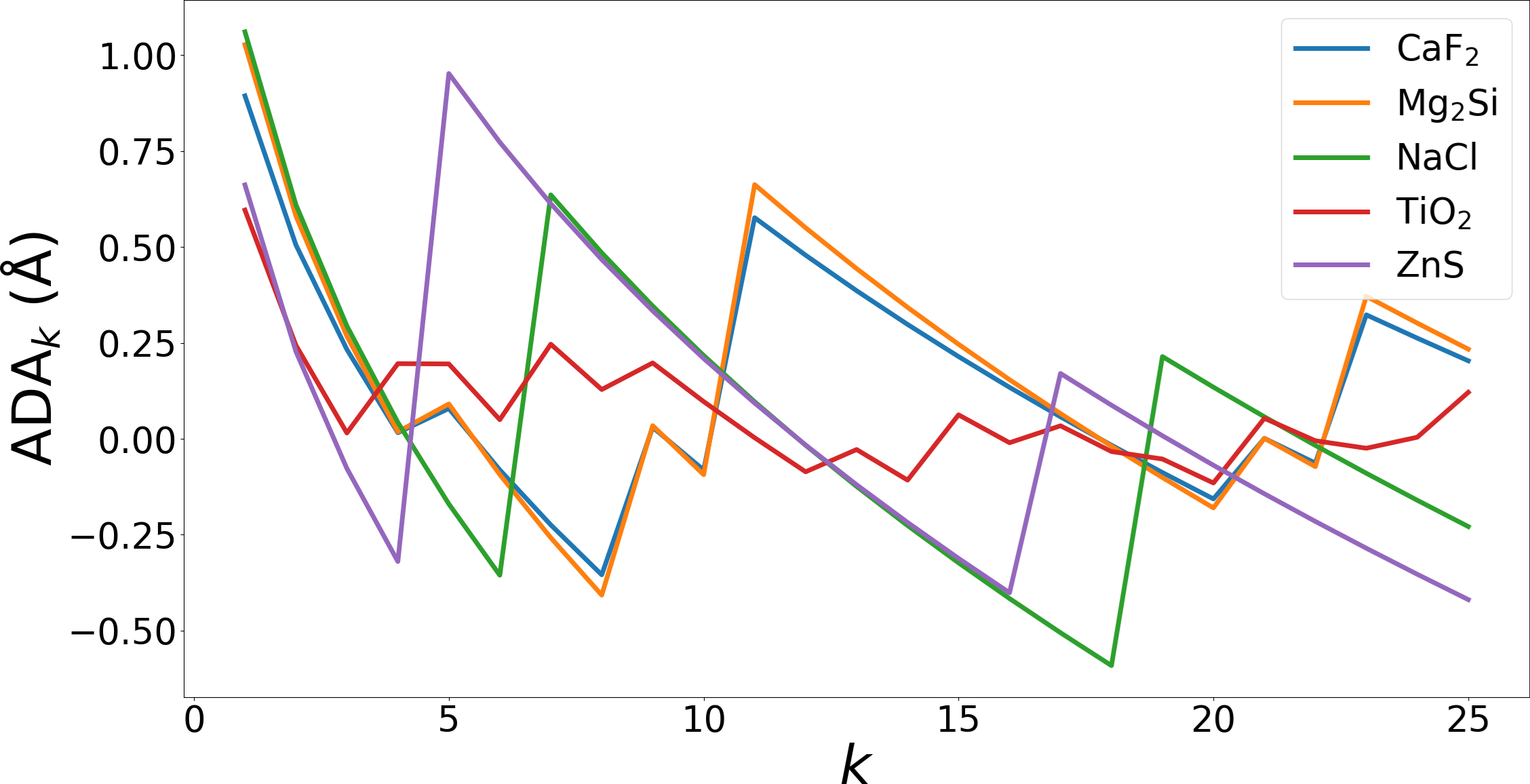}
\label{fig:AMD25ADA_5common_crystals}
\end{figure}

The invariants $\AMD_k$ and $\ADA_k$ form vectors of length $k$, e.g. set $\AMD(S;k)=(\AMD_1(S),\dots,\AMD_{k}(S))$ and $\ADA(S;k)=(\ADA_1(S),\dots,\ADA_{k}(S))$.
These vectors can be compared by many metrics.
The metric $L_\infty(u,v)=\max\limits_{i=1,\dots,k}|\vb*{u}_i-\vb*{v}_i|$ for any vectors $\vb*{u},\vb*{v}\in\R^k$ preserves the intuition of atomic displacements in the following sense. 
If $S$ is obtained from $Q$ by perturbing every point up to a small $\ep$, 
then $L_\infty(\AMD(S;k),\AMD(Q;k))\leq 2\ep$ by \cite[Theorem~9]{widdowson2022average}.
Other distances such as Euclidean can be considered but will accumulate a larger deviation depending on $k$.
\smallskip

All invariants above and metrics on them are measured in the same units as original coordinates, i.e. in Angstroms for crystals given by Crystallographic Information Files (CIFs). 
The Point Packing Coefficient $\PPC(S)$ was defined as the cube root of the cell volume per atom (of the same radius $1\angstrom$) and can be interpreted as an average radius of balls `packed' in a unit cell.
So $\PPC(S)$ is roughly inversely proportional to the physical density but they are exactly related only when materials have the same average atomic mass (total mass of atoms in a unit cell divided by the cell volume). 
\myskip

\begin{figure}[h!]
\caption{The averages of $\ADA_k$ and standard deviations (1 sigma shaded) vs $\sqrt[3]{k}$ for four databases. 
}
\includegraphics[width=\textwidth]{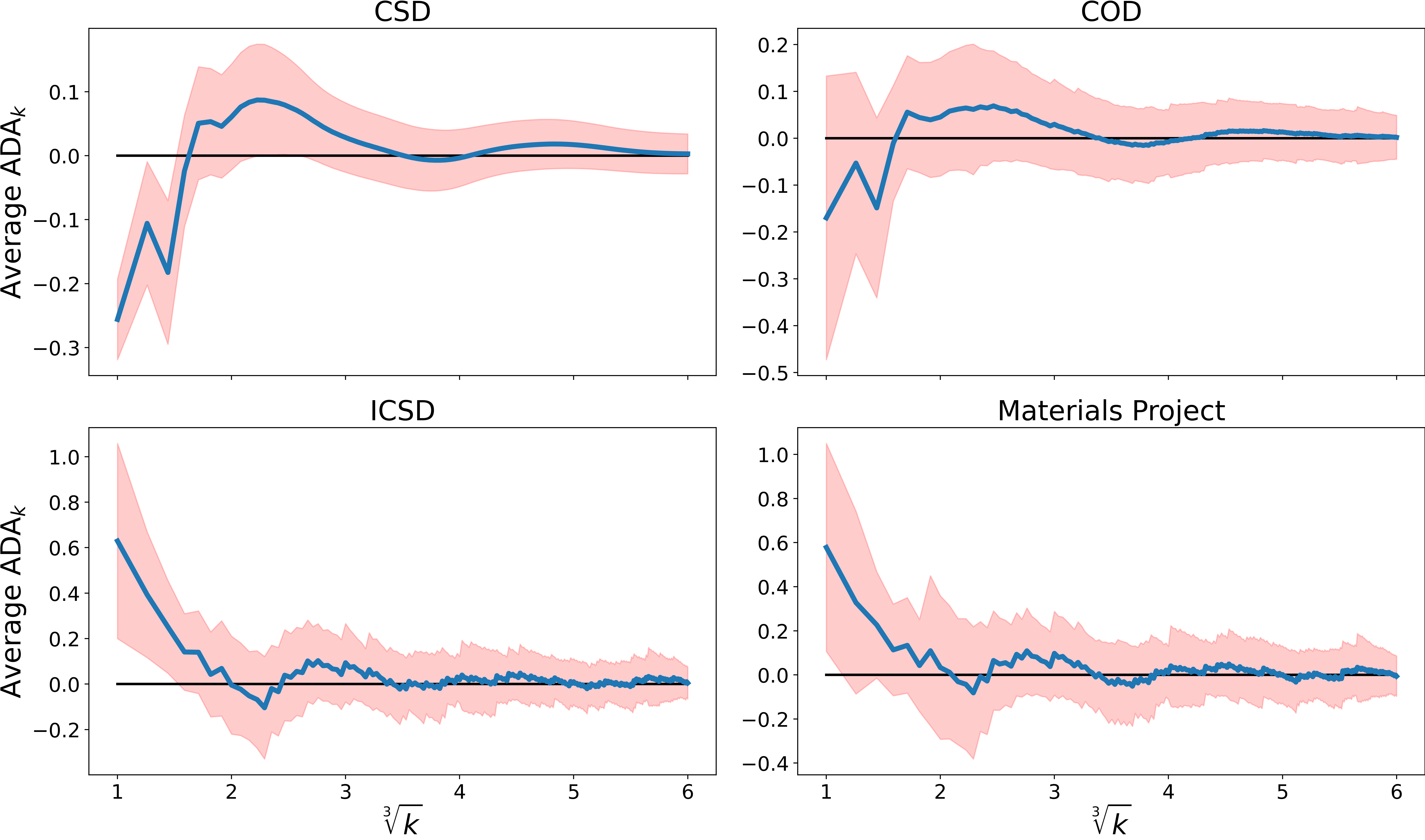}
\label{fig:ADA_averages}
\end{figure}

While $\AMD_k(S)$ monotonically increases in $k$, the invariants $\ADA_k(S)$ can be positive or negative as deviations around the asymptotic $\PPC(S)\sqrt[3]{k}$.
Fig.~\ref{fig:ADA_averages} reveals geometric differences between the mainly organic databases CSD and Crystallography Open Database (COD) \cite{gravzulis2009crystallography} versus the more inorganic collections ICSD and  MP.
\myskip

The first average of $\ADA_1\in[-0.25,-0.17]$ in the top images of Fig.~\ref{fig:ADA_averages} can be explained by the presence of many hydrogen atoms, which have distances smaller than $\PPC(S)$ to their first neighbour in most organic materials.
Indeed, hydrogens are usually bonded at distances less than $1.2\angstrom$, while $\PPC(S)$ is often larger than $1.2\angstrom$ because most chemical elements have van der Waals radii above $1.2\angstrom$ \cite{batsanov2001van}.
\myskip

For inorganic materials, metal atoms or ions have relatively large distances to their neighbours, so the average $\ADA_1$ is in $[0.58,0.62]$ in the bottom images of Fig.~\ref{fig:ADA_averages}.
\myskip

For all types of materials in Fig.~\ref{fig:ADA_averages}, the value of $\ADA_k$ experimentally converges to 0 on average, so there is no need to substantially increase $k$ because the important structural information emerges for smaller indices $k$ of neighbours.
\myskip
 
If we increase $k$, the matrix $\PDD(S;k)$ and hence the vector $\ADA(S;k)$ become longer by including distance data to further neighbours but all initial values remain the same.
Hence we consider $k$ not as a parameter that changes the output but as a degree of approximation similarly to the number of decimal places on a calculator.
\smallskip
 
The experimental convergence $\ADA_k\to 0$ as $k\to+\infty$ in Fig.~\ref{fig:ADA_averages}justifies computing the distance $L_\infty$ between $\ADA$ vectors up to a reasonable $k$.
We use $k=100$ because all $\ADA_k$ for $k>100$ are close to 0 (the range of 1 sigma between $\pm 0.2\angstrom$) in Fig.~\ref{fig:ADA_averages}.

\begin{cor}[invariance of $\AND,\PND$ under uniform scaling, {\cite[Corollary~3.9]{widdowson2025pointwise}}]
\label{cor:AND}
For any $l$-periodic set $S\subset\R^n$, $\AND_k(S)$ and $\PND(S;k)$ in Definition~\ref{dfn:ADA} are invariant under isometry and uniform scaling for $k\geq 1$.
Also, $\AND_k(S)\to 0$ as $k\to+\infty$.
\ecor
\end{cor}

The following conjecture is justified by Example~\ref{exa:ADA_asymptotic} and Fig~\ref{fig:ADA_averages}.

\begin{conj}[asymptotic of $\ADA(S;k)$ as $k\to+\infty$]
\label{conj:ADA_asymptotic}
For any periodic point set $S\subset\R^n$, we have $\ADA_k(S)\to 0$ as $k\to+\infty$.
\epro
\end{conj}

\begin{exa}[asymptotic $\ADA_k(S)\to 0$ as $k\to+\infty$ for the cubic lattice $S=\Z^n$]
\label{exa:ADA_asymptotic}
The survey \cite{ivic2004lattice} describes progress on the generalised Gauss circle problem expressing the number of points from the cubic lattice $\Z^n$ within a ball of a radius $r$ as $k=V_n r^n-O(r^{\al_n+\ep})$ for any $\ep>0$, where 
$\al_n<n-1$ for $n\geq 2$, e.g.
$\al_2\leq\frac{2}{3}$,  
$\al_3\leq\frac{3}{2}$, and 
$\al_n\leq n-2$ for any $n\geq 4$.
The cubic lattice $S=\Z^n$ has $\PPC(\Z^n)=1/\sqrt[n]{V_n}$.
Let $d_k$ denote the distance from the origin $0$ to its $k$-th neighbour in $\Z^n$.
Then 
$$k=V_n d_k^n-O(d_k^{\al_n+\ep}), \text{ so }
d_k = \sqrt[n]{\dfrac{k+O(d_k^{\al_n+\ep})}{V_n}}=\PPC(\Z^n)\sqrt[n]{k+O(d_k^{\al_n+\ep})},$$
$$\dfrac{\ADA_k(\Z^n)}{\PPC(\Z^n)}=\dfrac{d_k}{\PPC(\Z^n)}-\sqrt[n]{k}=\sqrt[n]{k+O(d_k^{\al_n+\ep})}-\sqrt[n]{k}
=\dfrac{O(d_k^{\al_n+\ep})}{P_n(\sqrt[n]{k+O(d_k^{\al_n+\ep})},\sqrt[n]{k})},$$ where $P_n$ is a homogeneous polynomial of degree $n-1$, e.g. $P_2(x,y)=x+y$, $P_3(x,y)=x^2+xy+y^2$. 
Since the numerator has the power $\al_n<n-1$ of $d_k=O(\sqrt[n]{k})$ for $n\geq 2$, the final expression of $\dfrac{\ADA_k(\Z^n)}{\PPC(\Z^n)}$ and hence $\ADA_k(\Z^n)$ have limit $0$ as $k\to+\infty$.
\eexa
\end{exa}

\begin{thm}[time of $\PDD$, {\cite[Theorem~3.10]{widdowson2025pointwise}}]
\label{thm:PDD_time}
Let $S\subset\R^n$ be any $l$-periodic set with a minimum inter-point distance $d_{\min}$ and a unit cell $U=P\times R^{n-l}$, where $P$ is a parallelepiped in the $l$-dimensional subspace $R^l$ with the orthogonal subspace $R^{n-l}$ in $\R^n$. 
Consider the \emph{width} $w=\sup\limits_{u,v\in P}|\vec u-\vec v|$ and the \emph{height} $h$ equal to the maximum distance between points in the orthogonal projection of $S$ to $R^{n-l}$.
If the motif $M=S\cap U$ consists of $m$ points, then $\PDD(S;k)$ can be computed for any $k\geq 1$ in time
$$O(km(2^{4n}\log k+\log m)+2^{12n}m\log^2 k +(2^{8n}/l)k\log k+2^{8n}a^l b k),$$
where $a=1+\dfrac{2.5w+2h}{\PPC(S)}$ and $b=\log(2\PPC(S)+3w+5h)-\log d_{\min}$.
The complexity of $\AMD(S;k)$ and invariants $\PDA(S;k),\PND(S;k)$ from Definition~\ref{dfn:ADA} is the same as for $\PDD(S;k)$, because the extra computations can be done in time $O(km)$.
\ethm
\end{thm}

The worst-case estimate in Theorem~\ref{thm:PDD_time} is conservative due to the upper bound $2^n$ for the expansion constants $c_{\min},c$ from \cite[Definition~1.4]{elkin2023new}.
We conjecture that this upper bound can be reduced to $2^l$ for any $l$-periodic point set $S\subset\R^n$. 
\smallskip

For any fixed dimensions $l\leq n$, if we ignore the parameters $a,b,d_{\min}$, and $\PPC(S)$, then the complexity in Theorem~\ref{thm:PDD_time} becomes $O(km(\log k+\log m))$, which is near-linear in both $k,m$. 
For the most practical dimensions $l=n=3$, experiments in section~\ref{sec:experiments} will report running times in minutes on a modest desktop computer for about 2 million real crystals from the world's largest materials databases.

\section{Lipschitz continuous metrics and local novelty distances}
\label{sec:PDD_continuous}

This section proves the Lipschitz continuity of the vectorial invariants $\AMD,\ADA,\AND$, matrix invariants $\PDD,\PDA,\PND$, and their averages in Theorem~\ref{thm:PDD_continuous}.
\myskip

The widely used Mercury software visually compares periodic structures \cite{chisholm2005compack} by minimizing the Root Mean Square Deviation (RMSD) of atomic positions from up to a given number $m$ (15 by default) of closest molecules in two structures.   
This comparison depends on many parameters (maximum number of matched molecules, thresholds for matched distances and angles), fails the triangle inequality in metric axioms, see Definition~\ref{dfn:metrics}, and is too slow for pairwise comparisons.
In fact, this RMSD tries to measure a maximum displacement of atoms only by using finite subsets of crystals.
\myskip

For full infinite crystals, the theoretically better alternative is the \emph{bottleneck distance} $\BD(S,Q)$ equal to the maximum Euclidean distance needed to perturb every point $p\in S$ to its unique match in $Q$, see Example~\ref{exa:metrics}(b).
\myskip

Example~\ref{exa:infinite_bottleneck}
shows that the periodic sequences have $d_B=+\infty$ for any $\de>0$.

\begin{exa}[infinite bottleneck]
\label{exa:infinite_bottleneck}
We show that $S=\Z$ and $Q=(1+\de)\Z$ for any $\de>0$ have $d_B(S,Q)=+\infty$.
Assuming that $d_B(S,Q)$ is finite, consider an interval $[-N,N]\subset\R$ containing $2N+1$ points of $S$.
If there is a bijection $g:S\to Q$ such that $|p-g(p)|\leq d_B$ for all points $p\in S$, the image of $2N+1$ points $S\cap [-N,N]$ under this bijection $g$ should be within the interval $[-N-d_B,N+d_B]$.
The last interval contains only $1+\frac{2(N+d_B)}{1+\de}$ points, which is smaller than $1+2N$ when $\frac{N+d_B}{1+\de}<N$ or $d_B<\de N$.
We get a contradiction by choosing a large $N>\frac{d_B}{\de}$.
\eexa
\end{exa}

If $S,Q$ are lattices of equal density (equal unit cell volume), they have a finite bottleneck distance $d_B$ by \cite[Theorem~1(iii)]{duneau1991bounded}. 
If we consider only periodic point sets $S,Q\subset\R^n$ with the same density (or unit cells of the same volume), $d_B(S,Q)$ becomes a well-defined \emph{wobbling} distance \cite{carstens1999geometrical}, which is discontinuous under perturbations below.

\begin{exa}[discontinuous wobbling]
\label{exa:dB_discontinuous}
Slightly perturb the basis $(1,0),(0,1)$ of the integer lattice $\Z^2$ to the basis vectors $(1,0),(\ep,1)$ of the new lattice $\La$.
We prove that $d_B(\La,\Z^2)\geq\frac{1}{3}$ for any $\ep>0$.
Map $\R^2$ by $\Z^2$-translations to the unit square $[0,1]^2$ with identified opposite sides (a torus).
Then the whole square lattice $\Z^2$ is mapped to the single point represented by the corners of the square $[0,1]^2$.
The perturbed lattice $\La$ maps to the sequence of points $\{k\ep\pmod{1}\}_{k=0}^{+\infty}\times\{0,1\}$ in the horizontal edges.
If $d_B(\La,\Z^2)=r<\frac{1}{3}$, then all above points should be covered by the closed disks of the radius $r$ centred at the corners of $[0,1]^2$.
For $0<\ep<\frac{1}{3}-r$, we can find $k\in\Z$ so that $k\ep$ is strictly between $r,1-r$, hence not covered by these disks, so $d_B(\La,\Z_2)\geq\frac{1}{3}$.
\eexa 
\end{exa}

We will use the Earth Mover's Distance ($\EMD$) from Definition~\ref{dfn:EMD}, which is well-defined for any normalised distributions of different sizes and makes sense for any matrix invariant $I(S)$ that is an unordered collection of row vectors $\vec R_i(S)$ with weights $w_i(S)\in(0,1]$ satisfying $\sum\limits_{i=1}^{m(S)} w_i(S)=1$.
Each row $\vec R_i(S)$ should have a size independent of $i$, e.g. the number $k$ of neighbours in $\PDD(S;k)$. 
For any vectors 
$\vec R_i=(r_{i1},\dots,r_{ik})$ and $\vec R_j=(r_{j1},\dots,r_{jk})$, we will use the Minkowski and Chebyshev distances from Example~\ref{exa:metrics}(b):  
$L_{q}(\vec R_i,\vec R_j)=\big(\sum\limits_{l=1}^k |r_{il}-r_{jl}|^q\big)^{1/q}$, 
$L_{\infty}(\vec R_i,\vec R_j)=\max\limits_{l=1,\dots,k}|r_{il}-r_{jl}|$.

\begin{exa}[$\EMD$ under noise]
\label{exa:EMD_noise}
We illustrate $\EMD$ for perturbations that scale up a unit cell as in Fig.~\ref{fig:noise_double_peaks}~(left).
The integer sequence $\Z$ has $\PDD(\Z;2)=(1;1,1)$, a single row of weight 1 and unit distances to 2 neighbours.
The periodic sequence $\Z_\ep=\{0,1+\ep,2-\ep\}+3\Z$ is obtained from $\Z$ by $\ep$-perturbations of points $1,2$ and all their translates with period 3.
Then $\PDD(\Z_\ep;2)=\left(\begin{array}{c|cccc} 
1/3 & 1+\ep & 1+\ep \\
2/3 & 1-2\ep & 1+\ep
\end{array} \right)$, where the 2nd row represents the shifted points $1+\ep,2-\ep$.
After splitting $\PDD(\Z;2)=(1;1,1)$ into two identical rows of weights $\frac{1}{3},\frac{2}{3}$ and using $L_\infty$ on vectors of two distances, a difference between $\PDD$s can be defined as the weighted average $\frac{1}{3}\ep+\frac{2}{3}2\ep=\frac{2}{3}\ep$. 
\eexa
\end{exa}

\begin{figure}[h!]
\caption{The Pointwise Distance Distribution $\PDD(S;4)$ of a unit square lattice $S$ changes continuously under noise, shown here in the case when a 1-point motif extends to a motif of 4-points that form an isosceles trapezium with parallel sides of lengths 0.8 (red) and 1.2 (blue).
}
\includegraphics[width=\textwidth]{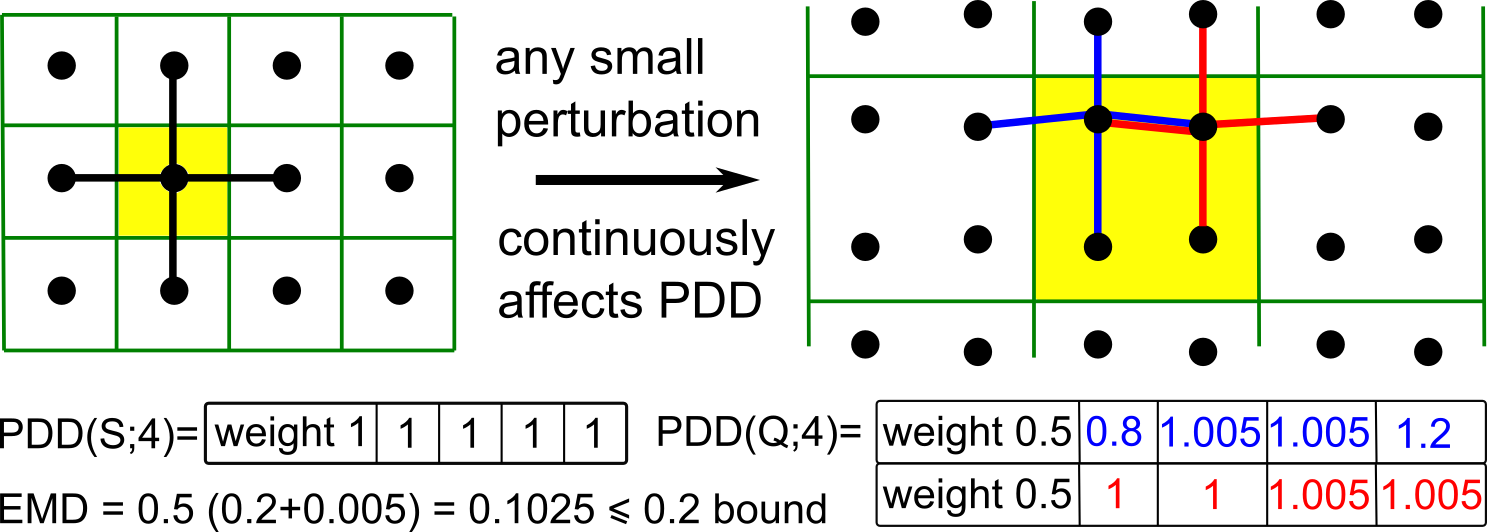}
\label{fig:PDD_under_noise_square4}
\end{figure}

Fig.~\ref{fig:PDD_under_noise_square4} illustrates the continuity of $\PDD(S;4)$ under a perturbation of a square lattice $S\subset\R^2$, which scales up an initial cell by a factor of $4$.

The Lipschitz continuity of invariants in $\EMD$ will use bounded perturbations of points up to Euclidean distance $\ep$ in $\R^n$.
Recall that the \emph{packing radius} $r(S)$ is the minimum half-distance between any points of $S$, see Definition~\ref{dfn:packing+covering}(a).

\begin{thm}[Lipschitz continuity of $\PDA$ and $\PND$, {\cite[Theorem~4.2]{widdowson2025pointwise}}]
\label{thm:PDD_continuous} 
Let $S,Q\subset\R^n$ be $l$-periodic point sets such that $Q$ is obtained from $S$ by perturbing every point of $S$ up to Euclidean distance $\ep$.
Fix any $q\in[1,+\infty]$ and an integer $k\geq 1$.
Interpret $\sqrt[q]{k}$ as 1 in the limit case $q=+\infty$.
If $\min\{r(S),r(Q)\}>\ep$, then $\PPC(S)=\PPC(Q)$, 
\sskip

\nt
\tb{(a)}
$\EMD_q(\PDA(S;k),\PDA(Q;k))\leq 2\ep\sqrt[q]{k}$, and
\sskip

\nt 
\tb{(b)}
$\EMD_q(\PND(S;k),\PND(Q;k))\leq 
\dfrac{2\ep\sqrt[q]{k}}{\PPC(S)}$.
\ethm
\end{thm}

All columns of $\PDD,\PDA,\PND$ are ordered by the index $k$ of neighbours.
Though their rows are unordered (as points of a motif $M$), all such matrices (even with different numbers of rows) can be compared by Earth Mover's Distance, or by any other metrics on weighted distributions, see Definition~\ref{dfn:EMD}. 
We can simplify any $\PDD$ into a fixed-size matrix, which can be flattened into a vector, while keeping the continuity and almost all invariant data. 
Any distribution of $m$ unordered values can be reconstructed from its $m$ moments in Definition~\ref{dfn:moments}.
When all weights $w_i$ are rational, as in our case, the distribution can be expanded to equal-weighted values $a_1,\dots,a_m$.
The $m$ moments can recover all $a_1,\dots,a_m$ as roots of a degree $m$ polynomial whose coefficients are expressed via the $m$ moments \cite{macdonald1998symmetric}, e.g. any $a,b\in\R$ can be found from $a+b,a^2+b^2$ as the roots of $x^2-(a+b)x+ab$, where $ab=\frac{1}{2}((a+b)^2-(a^2+b^2))$.
\myskip

For $t=1$, the $1\times k$ matrix $\mu^{(1)}[\PDD(S;k)]$ appeared in Definition~\ref{dfn:PDD} as the vector $\AMD(S;k)=(\AMD_1,\dots,\AMD_k)$ of column averages.
All rows and columns of $\mu^{(t)}[I(S)]$ are ordered, but this matrix is a bit weaker than $I(S)$ because each column can be reconstructed from its moments (for a large enough $t$) only up to permutation.
We can flatten any matrix $\mu^{(t)}[I(S)]$ 
to a vector for machine learning. 

\begin{thm}[lower bounds, {\cite[Theorem~4.4]{widdowson2025pointwise}}]
\label{thm:lower_bound}
For any $l$-periodic sets $S,Q\subset\R^n$, 
\myskip

\noindent
\tb{(a)} 
$\EMD_q(\PDD(S;k),\PDD(Q;k))\geq L_q(\AMD(S;k), \AMD(Q;k))$;
\myskip

\noindent
\tb{(b)} 
$\EMD_q(\PDA(S;k),\PDA(Q;k))\geq L_q(\ADA(S;k), \ADA(Q;k))$;
\myskip

\noindent
\tb{(c)} 
$\EMD_q(\PND(S;k),\PND(Q;k))\geq L_q(\AND(S;k), \AND(Q;k))$ for $q,k\geq 1$. 
\ethm
\end{thm}


\section{Generic completeness of PDDs for periodic point sets in $\R^n$}
\label{sec:PDD_gen_complete}

While the generic completeness of the $\PDD$ was easy for any finite clouds in $\R^n$, this section
extends Theorem~\ref{thm:PDD_finite_gen_complete}
to the much harder periodic case in Theorem~\ref{thm:PDD_periodic_gen_complete}.
\myskip

For a periodic point set $S\subset\R^n$, the generic completeness of $\PDD$ is not straightforward, because infinitely many distances between points of $S$ are repeated due to periodicity.
We introduce a few auxiliary concepts for distance-generic periodic sets in Definition~\ref{dfn:distance-generic_set}.
For any point $p$ in a lattice $\La\subset\R^n$, the \emph{open Voronoi domain} 
$$V(\La;p)=\{q\in\R^n \text{ such that } |q-p|<|q-p'| \mbox{ for any }p'\in\La-p\}$$ is the neighbourhood of all points $q\in\R^n$ that are strictly closer to $p$ than to all other points $p'$ of the lattice $\La$.
Definition~\ref{dfn:Voronoi_domain}(a) used the closed version $\bar V(\La;p)$. 
\myskip

Open Voronoi domains $V(\La;p)$ of different points $p\in\La$ are disjoint translation copies of each other and their closures tile $\R^n$, so $\cup_{p\in\La}\bar V(\La;p)=\R^n$.
For example, for a generic lattice $\La\subset\R^2$, the domain $V(\La;p)$ is a centrally symmetric hexagon.
\myskip

Points $p,p'\in\La$ are \emph{Voronoi neighbours} if their Voronoi domains share a boundary point, so $\bar V(\La;p)\cap\bar V(\La,p')\neq\emptyset$.
Below we always assume that any lattice $\La$ is shifted to contain the origin $0$, also any periodic point set $S=\La+M$ has a point at $0$.

\begin{dfn}[neighbour set $N(\La)$ and base distances]
\label{dfn:neighbour_set}
For any lattice $\La\subset\R^n$, the \emph{neighbour set} of the origin 0 is $N(\La)=\La\cap\bar B(0;r)\setminus\{0\}$ for a minimum radius $r$ such that $N(\La)$ is not contained in any affine $(n-1)$-dimensional subspace of $\R^n$, and $N(\La)$ includes all $n+1$ nearest neighbours (within $\La$) of any point $q\in V(\La;0)$.
\smallskip

Consider all sets of unordered points $p_1,\dots,p_n\in N(\La)$ that are \emph{linearly independent}, i.e. the vectors $\vec p_1,\dots,\vec p_n$ form a linear basis of $\R^n$.
For any point $q\in V(\La;0)$, a lexicographically smallest list of distances $d_1(q)\leq\dots\leq d_n(q)$ from $q$ to a set of linearly independent points $p_1,\dots,p_n\in N(\La)$ is the list of \emph{base distances} of $q$.
\edfn
\end{dfn}

The linear independence of vectors $\vec p_1,\dots,\vec p_n$ in Definition~\ref{dfn:neighbour_set} guarantees that any point $q$ is uniquely determined in $\R^n$ by its distances $|q|,d_1(q),\dots,d_n(q)$ to $n+1$ neighbours $0,p_1,\dots,p_n$, which are not in the same $(n-1)$-dimensional subspace.

\begin{exa}[neighbour sets]
\label{exa:neighbour_sets}
The vector $(2,0),(0,1)$ generate the rectangular lattice $\La\subset\R^2$.
The Voronoi domain $V(\La;0)$ is the rectangle $(-1,1)\times(-0.5,0,5)$.
The neighbour set $N(\La)\subset\La$ includes the 3rd neighbours $(0,\pm 2)$ of the points $(0,\pm 0.4)\in V(\La;0)$.
Indeed, if in Definition~\ref{dfn:neighbour_set} $\La$ has a radius $r<2$, then $\La\cap\bar B(0;r)\setminus\{0\}=\{(0,\pm 1)\}$ is in the 1-dimensional subspace ($y$-axis) of $\R^2$.
For $q=(0,0.4)$, considering all pairs $(\vec p_1,\vec p_2)$ that generate $\R^2$ among the four possibilities $((0,\pm 1),(\pm 2,0))$, we find the base distances $d_1(q)=0.6<d_2(q)=\sqrt{0.4^2+2^2}\approx 2.04$ for the 2nd and 3rd lattice neighbours $p_1=(0,1)$ and $p_2=(\pm 2,0)$ of $q$, respectively.
\eexa
\end{exa}

\begin{dfn}[a distance-generic set]
\label{dfn:distance-generic_set}
A periodic point set $S=M+\La\subset\R^n$ with the origin
$0\in\La\subset S$ is called \emph{distance-generic} if the following conditions hold.
\smallskip

\noindent
\tb{(a)}
For any points $p,q\in S\cap V(\La;0)$, the vectors $\vec p,\vec q$ are not orthogonal.
\myskip

\noindent
\tb{(b)}
For vectors $\vec u,\vec v$ between any two pairs of points in $S$, if $|\vec u|=l|\vec v|\leq 2R(\La)$ for $l=1,2$, then $\vec u=\pm l\vec v$ and $\vec v\in\La$.
\myskip

\noindent
\tb{(c)}
For any point $q\in S\cap V(\La;0)$, let $d_0=|q|$ be its distance to the closest neighbour $p_0=0$ in $\La$.
Take any linearly independent points $p_1,\dots,p_n\in N(\La)$ and any distances $d_1\leq\dots\leq d_n$ from $q$ to some points in $S\cap \bar B(0;2R(\La))$.
The $n+1$ spheres $\bd B(p_i;d_i)$ can meet at a single point of $S\cap V(\La;0)$ only if $d_1\leq\dots\leq d_n$ are the base distances of $q$ and only for two tuples $p_1,\dots,p_n\in N(\La)$ related by 
$\vec v\mapsto -\vec v$.
\edfn
\end{dfn}

Condition~\ref{dfn:distance-generic_set}(b) means that all inter-point distances are distinct apart from necessary exceptions due to periodicity.
Since any periodic set $S=M+\La\subset\R^n$ is invariant under translations along all vectors of $\La$, condition~\ref{dfn:distance-generic_set}(b) for $|\vec v|\leq 2R(\La)$ can be checked only for vectors from all points of $S$ in the original Voronoi domain $V(\La;0)$ to all points in the domain $3V(\La;0)$ extended by factor 3. 
Condition~\ref{dfn:distance-generic_set}(b) implies that $S$ has no points on the boundary $\bd V(\La;0)$, because any such point is equidistant to points $0,v\in\La$ and hence should belong to $\La$. 
Let a \emph{lattice distance} be the Euclidean distance from any $p\in M=S\cap V(\La;0)$ to its lattice translate $p+\vec v$ for all $\vec v\in\La$.
Condition~\ref{dfn:distance-generic_set}(a) guarantees that only a lattice distance $d$ appears together with $2d$ (and possibly with higher multiples) in a row of $\PDD(S;k)$.
Any such $d$ and its multiples are repeated twice in every row, because $\La$ is centrally symmetric. 

\begin{lem}[almost any periodic set is distance-generic, {\cite[Lemma~5.7]{widdowson2025pointwise}}]
\label{lem:distance-generic}
Let $S=M+\La\subset\R^n$ be any periodic point set. 
For any $\ep>0$, one can perturb coordinates of a basis of $\La$ and of points from $M$ up to $\ep$ such that the resulting perturbation $S'$ of $S$ is a distance-generic periodic point set in the sense of Definition~\ref{dfn:distance-generic_set}. 
\elem
\end{lem}

The size $m$ of a motif $M$ is an isometry invariant because any isometry maps $N$ to another hose motif of the same size.
In dimensions $n=2,3$, any lattice $\La$ can be reconstructed from its complete isometry invariants \cite{kurlin2024mathematics,kurlin2022complete}.
Theorem~\ref{thm:PDD_periodic_gen_complete} reconstructs a periodic point set $S=M+\La\subset\R^n$ in any dimension $n\geq 2$ from the invariant $I(S)$ consisting of $m$, $\PDD(S;k)$, and (complete invariants of) a lattice $\La$ to satisfy condition~\ref{pro:l-periodic}(a) for distance-generic periodic sets $S\subset\R^n$. 
Recall that the \emph{packing} radius $R(\La)$ is the smallest radius $R$ such that $\bigcup\limits_{p\in\La} \bar B(p;R)=\R^n$, see Definition~\ref{dfn:packing+covering}(b).
 
\begin{thm}[generic completeness of $\PDD$ for periodic sets, {\cite[Theorem~5.8]{widdowson2025pointwise}}]
\label{thm:PDD_periodic_gen_complete}
Let $S=M+\La\subset\R^n$ be any distance-generic periodic set. 
For any $k$ such that all distances in the last column of $\PDD(S;k)$ are larger than $2R(\La)$, the set $S$ can be reconstructed from $\La$, $\PDD(S;k)$ and the size $m$ of a motif of $S$, uniquely under isometry in $\R^n$.
\ethm
\end{thm}

\section{Detecting near-duplicates in the world's largest databases}
\label{sec:experiments}

This section reports thousands of previously unknown (near-)duplicates in the world's largest databases \cite{taylor2019million,gravzulis2009crystallography,zagorac2019recent,jain2013commentary,merchant2023millions}.
The sizes in Table~\ref{tab:databases} below are the numbers of all periodic crystals (with no disorder and full geometric data) in September 2024 (total number is 1,847,462, see Table~\ref{tab:COMPACK_times} and all experimental details in \cite[appendix SM1]{widdowson2025pointwise}. 
\myskip

\begin{table}[h!]
\caption{Links and verisons of the world's largest materials databases, see their sizes in Table~\ref{tab:COMPACK_times}.}
\label{tab:databases}
\begin{center}
\begin{tabular}{ll}
database and web address & version \\
\hline
CSD: Cambridge Structural Database,  http://ccdc.cam.ac.uk & version 6.00  \\
COD: Crystallography Open Database, crystallography.net/cod & July 30, 2024 \\
ICSD: Inorganic Crystal Structures, icsd.products.fiz-karlsruhe.de & Feb 25, 2025 \\
MP: Materials Project, http://next-gen.materialsproject.org & v2023.11.1 \\ 
GNoME: github.com/google-deepmind/materials\_discovery & Nov 29, 2023 \\
\end{tabular}
\end{center}
\end{table}

We first used the vector $\ADA(S;100)$ to find nearest neighbours across all databases by $k$-d trees \cite{gieseke2014buffer} up to $L_\infty\leq 0.01\angstrom$.
Since the smallest inter-atomic distances are about $1\angstrom=10^{-10}$m, atomic displacements up to $0.01\angstrom$ are 
considered experimental noise.
For the closest pairs found by $\ADA(S;100)$, the stronger $\PDA(S;100)$
can have only equal or larger $\EMD\geq L_\infty$ by Theorem~\ref{thm:lower_bound}.
The CSD, COD, ICSD should contain experimental structures.
MP is obtained from ICSD by extra optimisation.
\smallskip

Table~\ref{tab:EMD_Linf_PDA100_leq001A} shows that the well-curated 60-year-old CSD has 0.9\% near-duplicate crystals, while more than a third of the ICSD consists of near-duplicates that are geometrically almost identical so that all atoms can be matched by an average perturbation up to $0.01\angstrom$.
Table~1 in \cite[section~6]{anosova2024importance} reported many thousands of exact duplicates, where chemical elements were replaced while keeping all coordinates fixed.  
These replacements are physically impossible without more substantial perturbations. Five journals are investigating integrity \cite{chawla2023crystallography}.
\myskip

The bold numbers in Table~\ref{tab:EMD_Linf_PDA100_leq001A} count near-duplicates, and their percentages within each database, which should be filtered out, else the ground truth data becomes skewed.
Table~\ref{tab:near-duplicates_diff_cells} confirms that cell-based comparisons miss near-duplicates as in Fig.~\ref{fig:noise_double_peaks}.

\begin{table}[h!]
\begin{center}
\setlength{\tabcolsep}{3pt}
\caption{Count and percentage of all ideal periodic crystals in each database (left) found to have a near-duplicate in other databases (top) by the distance $\EMD\leq 0.01\angstrom$ on matrices $\PDA(S;100)$.}
\label{tab:EMD_Linf_PDA100_leq001A}
\begin{tabular}{l|cc|cc|cc|cc|cc}
duplicates & \multicolumn{2}{c|}{CSD} & \multicolumn{2}{c|}{COD} & \multicolumn{2}{c|}{ICSD} & \multicolumn{2}{c|}{MP} & \multicolumn{2}{c}{GNoME}  \\ 
in databases    & count  & \%     & count  & \%     & count & \%     & count & \%   & count & \%     \\ 
      \hline
CSD   & 8343 & 0.92 & 283000 & 31.19 & 26506 & 2.92 & 33 & 0.00 & 1 & 0.00\\
COD   & 286663 & 80.18 & 19568 & 5.47 & 47065 & 13.16 & 5231 & 1.46 & 2705 & 0.76\\
ICSD  & 26853 & 15.78 & 69948 & 41.10 & 51085 & 30.01 & 27194 & 15.98 & 15449 & 9.08\\
MP    & 73 & 0.05 & 11986 & 7.82 & 15188 & 9.91 & 19177 & 12.51 & 10681 & 6.97\\
GNoME & 2 & 0.00 & 1800 & 0.47 & 2614 & 0.68 & 3401 & 0.88 & 82859 & 21.53
\end{tabular}
\end{center}
\end{table}

\begin{table}[h!]
\setlength{\tabcolsep}{3pt}
\begin{center}
\caption{Near-duplicates from Table~\ref{tab:EMD_Linf_PDA100_leq001A} whose unit cells differ by $0.01\angstrom$.
Unit cells are compared by the Chebyshev metric $L_\infty$ between vectors of corresponding lengths of 3 edges and 3 face diagonals.}
\label{tab:near-duplicates_diff_cells}
\begin{tabular}{l|cc|cc|cc|cc|cc}
duplicates & \multicolumn{2}{c|}{CSD} & \multicolumn{2}{c|}{COD} & \multicolumn{2}{c|}{ICSD} & \multicolumn{2}{c|}{MP} & \multicolumn{2}{c}{GNoME}\\ 
in databases    & count  & \%     & count  & \%     & count & \%     & count & \%  & count & \%   \\ 
      \hline
CSD   & 776 & 0.09 & 419 & 0.05 & 210 & 0.02 & 29 & 0.00 & 1 & 0.00\\
COD   & 472 & 0.13 & 7263 & 2.03 & 8629 & 2.41 & 5059 & 1.42 & 2684 & 0.75\\
ICSD  & 462 & 0.27 & 28863 & 16.96 & 42946 & 25.23 & 26554 & 15.60 & 15360 & 9.02\\
MP    & 70 & 0.05 & 11790 & 7.69 & 14915 & 9.73 & 18582 & 12.13 & 10608 & 6.92\\
GNoME & 2 & 0.00 & 1786 & 0.46 & 2590 & 0.67 & 3346 & 0.87 & 60248 & 15.65
\end{tabular} 
\end{center}
\end{table}

\begin{figure}
\centering
\includegraphics[height=43mm]{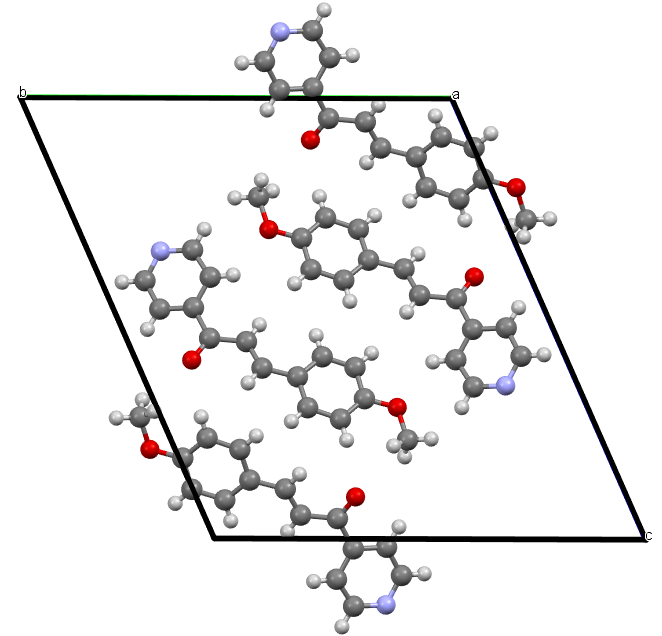}
\includegraphics[height=43mm]{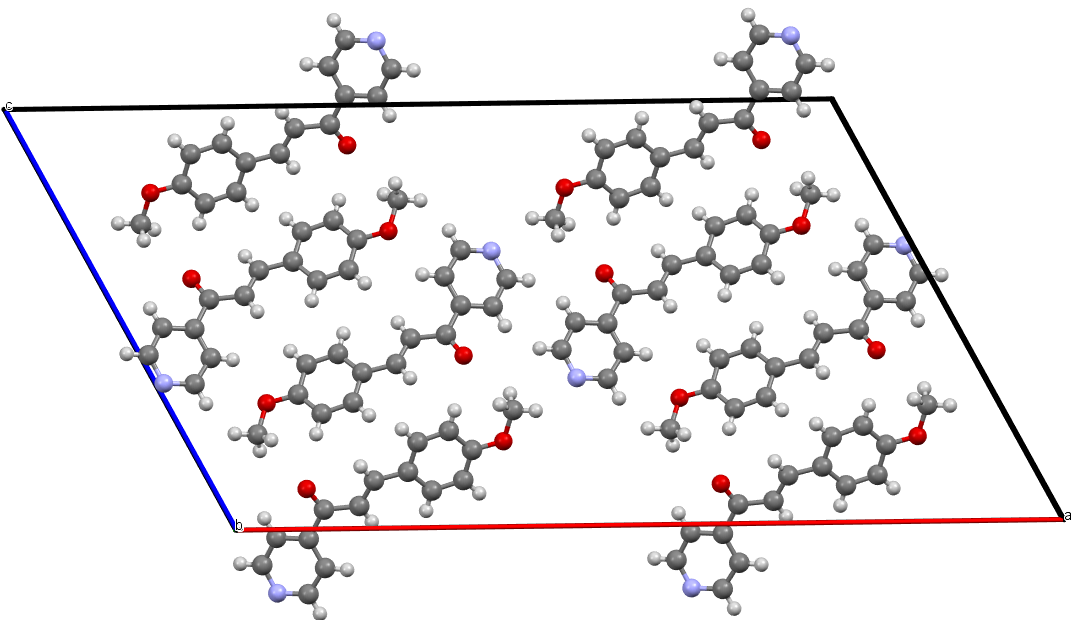}
\caption{In the CSD, near-duplicates PUBTEM (left) and PUBTEM01 (right) have a very small $\EMD= 0.00038\angstrom$ on invariants $\PDA(S;100)$, though their unit cells are rather different.}
\label{fig:CSD_near-duplicates_diff_cells}
\end{figure}

\begin{figure}
\centering
\includegraphics[width=\textwidth]{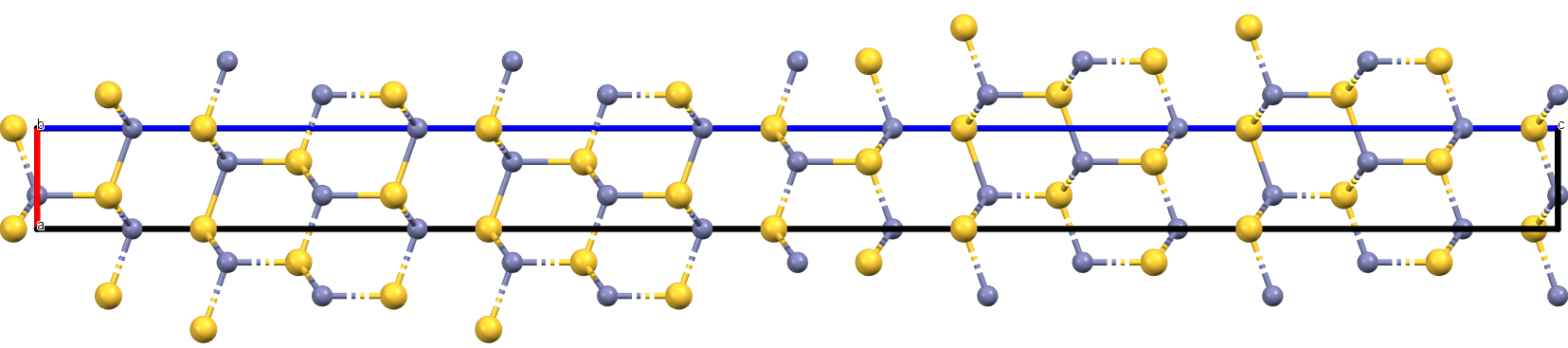}
\includegraphics[width=\textwidth]{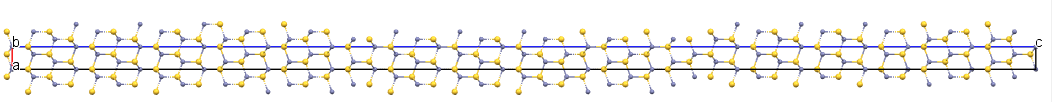}
\caption{In the COD, near-duplicates 2310812 (top) and 2310813 (bottom) have a very small $\EMD= 0.0008\angstrom$ on invariants $\PDA(S;100)$, though their unit  cells differ by a factor of about 3.}
\label{fig:COD_near-duplicates_diff_cells}
\end{figure}

\begin{figure}
\centering
\includegraphics[width=\textwidth]{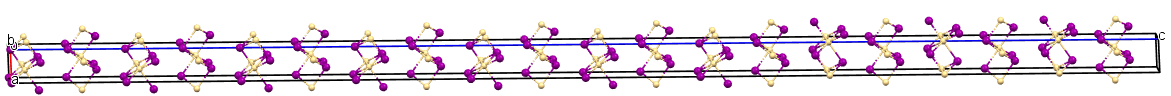}
\includegraphics[width=\textwidth]{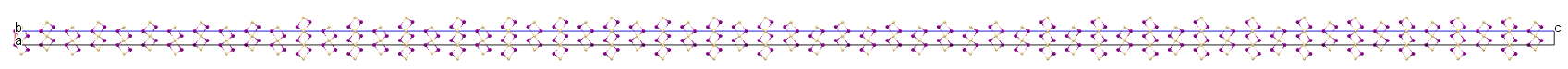}
\caption{In the ICSD, near-duplicates 42291 (top) and 42302 (bottom) have a very small $\EMD= 0.0024\angstrom$ on invariants $\PDA(S;100)$, though their unit  cells differ by a factor of about 3.}
\label{fig:ICSD_near-duplicates_diff_cells}
\end{figure}

\begin{figure}
\centering
\includegraphics[height=14.5mm]{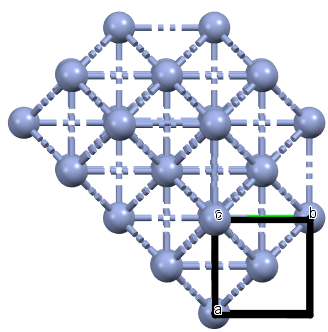}
\includegraphics[height=14.5mm]{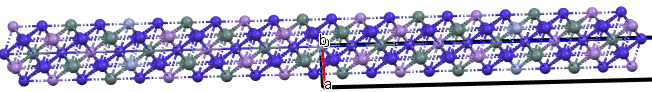}
\caption{In the Materials Project, near-duplicate entries mp-90 (left) and mp-1221808 (right) have a very small $\EMD= 0.0087\angstrom$ on invariants $\PDA(S;100)$, though their unit cells substantially differ.}
\label{fig:MP_near-duplicates_diff_cells}
\end{figure}

\begin{figure}
\centering
\includegraphics[height=42mm]{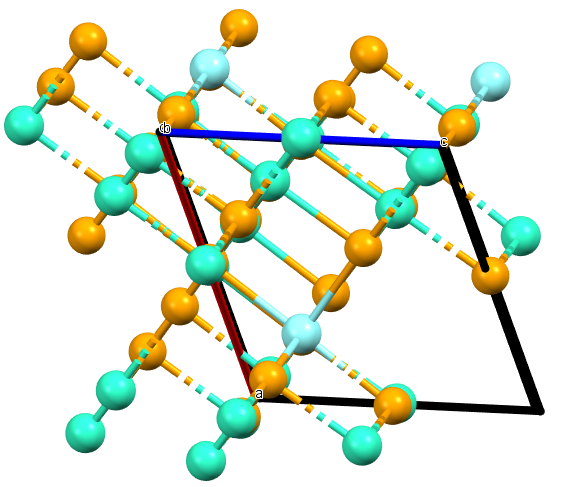}
\includegraphics[height=42mm]{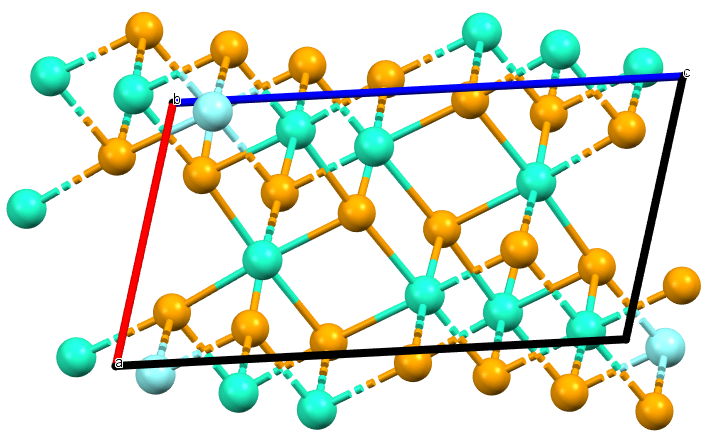}
\caption{In the GNoME, near-duplicates 4cb3b6ed9f (left) and 776c1b7570 (right) in the GNoME have $\EMD= 0.0079\angstrom$ on invariants $\PDA(S;100)$, though their unit cells are very different.}
\label{fig:GNoME_near-duplicates_diff_cells}
\end{figure}

Figures
\ref{fig:CSD_near-duplicates_diff_cells},
\ref{fig:COD_near-duplicates_diff_cells},
\ref{fig:ICSD_near-duplicates_diff_cells},
\ref{fig:MP_near-duplicates_diff_cells},
\ref{fig:GNoME_near-duplicates_diff_cells},
 show near-duplicates with very different cells, which were counted in Table~\ref{tab:near-duplicates_diff_cells}. 
In the past, the (near-)duplicates were impossible to detect at scale, because the traditional comparison through iterative alignment of 15 (by default) molecules by the COMPACK algorithm \cite{chisholm2005compack} is too slow for all-vs-all comparisons.
Tables~\ref{tab:times_duplicates_PDA100} and SM6 
compare the running times: \textbf{minutes} of $\PDA(S;100)$ vs \textbf{years} of RMSD, extrapolated for the same machine from the median time 117 milliseconds (582 ms on average) for 500 random pairs 
in the CSD. 
On the same 500 pairs, $\PDA(S;100)$ for two crystals and $\EMD$ together took only 7.48 ms on average.
All experiments were done on a modest desktop computer (AMD Ryzen 5 5600X 6-core, 32GB RAM). 

\begin{table}[h!]
\setlength{\tabcolsep}{4pt}
\caption{Times to compute $\PDA(S;100)$ and find all near-duplicates in Table~\ref{tab:EMD_Linf_PDA100_leq001A} with $\EMD\leq 0.01\angstrom$ across all major databases (seconds in the last 4 columns), compare with years in Table~\ref{tab:COMPACK_times}. }
\label{tab:times_duplicates_PDA100}
\begin{center}
\begin{tabular}{l|r|r|r|r|r|r|r}
database & $\PDA$, min:sec & $\EMD$, min:sec & CSD  & COD  & ICSD  & MP & GNoME  \\ 
\hline
CSD & 60:44  & 12:21  & 125.5 & 498.1 & 77.0 & 19.1 & 20.6 \\
COD & 30:16 & 16:29 & 524.5 & 122.0  & 235.1 & 79.6 & 27.0 \\
ICSD & 5:57  & 22:04 &  80.5  & 239.3  & 515.8  & 414.9 & 73.5 \\
MP & 1:40  & 13:31 & 28.2   & 82.9 & 413.8  & 222.8 & 63.0 \\
GNoME & 4:07  & 18:59 & 29.0   & 26.7 & 74.5  & 64.5 & 943.7
\end{tabular}
\end{center}
\end{table}

\begin{table}[h!]
\caption{These times for all comparisons by COMPACK \cite{chisholm2005compack} are extrapolated from the median time of 117 ms on 500 random pairs from the CSD 
on the same machine, which completed Table~\ref{tab:EMD_Linf_PDA100_leq001A} 
of near-duplicates across all the major databases within 2 hours.}
\label{tab:COMPACK_times}
\begin{center}
\begin{tabular}{l|c|r|r|r|r}
database & periodic crystals & unordered pairs & COMPACK time, sec & 
years \\
\hline
CSD & 907,246 & 411,547,198,635 & $4.81\times 10^{10}$ & 
1526 \\
COD & 357,510 & 63,906,521,295 & $7.48\times 10^{9}$ & 
237 \\
ICSD & 170,206 & 14,484,956,115 & $1.69\times 10^{9}$ & 
53 \\
MP & 153,235 & 11,740,405,995 & $1.37\times 10^{9}$ & 
43 \\
GNoME & 384,938 & 74,088,439,453 & $8.67\times 10^{9}$ & 274
\end{tabular}
\end{center}
\end{table}

Table~\ref{tab:periodic_invariant_comparisons} compares the proven properties of past and new descriptors.
All invariants based on cut-off atomic environments, such as MACE \cite{batatia2022mace}, discontinuously change under almost any perturbation  that arbitrarily scales up a primitive cell, as in Fig.~\ref{fig:noise_double_peaks}~(left). 
\myskip

The $\PDD$ remains continuous by taking into account only distances to neighbours rather than indices or relative positions of neighbours, which are discontinuous at cut-off boundaries.
Another exception is the complete isoset invariant in the next chapter.
\myskip

\begin{table}
\centering
\begin{tabular}{lccccc} 
	                 Descriptor                  & Invariant  & Continuity &  Complete   & Reconstruction & Time \\ 
	                 \hline
	             primitive cell                & $\times$ & $\times$ & $\times$  &   $\times$   & \checkmark \\
	                reduced cell                 & \checkmark & $\times$ & $\times$  &   $\times$   & \checkmark \\
	                space group                  & \checkmark & $\times$ & $\times$  &   $\times$   & \checkmark \\
	   PDF~\cite{terban2021structural}    & \checkmark & \checkmark & $\times$  &   $\times$   & \checkmark* \\
	      MACE~\cite{batatia2022mace}                       & \checkmark & $\times$ & \checkmark*  &   $\times$   & \checkmark* \\
	      \hline
	densities~\cite{edelsbrunner2021density} & \checkmark & \checkmark & \checkmark* &   $\times$   & \checkmark* \\
	       AMD \cite{widdowson2022average}      & \checkmark & \checkmark & $\times$  &   $\times$   & \checkmark \\
	       PDD \cite{widdowson2022resolving}     & \checkmark & \checkmark & \checkmark* &  \checkmark*   & \checkmark \\
	       	      isosets~\cite{anosova2021isometry,anosova2025recognition}       & \checkmark & \checkmark & \checkmark  &   \checkmark   & \checkmark* 
\end{tabular}
\caption{Comparison of crystal descriptors in the context of Problem~\ref{pro:l-periodic}. 
\checkmark* in the `Computable' column indicates that only an approximate algorithm exists for distances, and \checkmark* in the `Complete' and `Reconstruction' columns means that the condition holds only for generic periodic sets as in \ref{pro:l-periodic}(a).}
\label{tab:periodic_invariant_comparisons}
\end{table}

\section{Structural novelty of crystals and navigating materials space}
\label{sec:novelty}

This section leverages the strength of $\PDA$s to quickly and continuously  quantify the novelty of any periodic crystal relative to a given dataset in Definition~\ref{dfn:LND}.

\begin{dfn}[Local Novelty Distance $\LND(S;D)$]
\label{dfn:LND}
Let $D$ be a finite dataset of periodic point sets.
Fix an integer $k\geq 1$.
For any periodic point set $S$, the \emph{Local Novelty Distance} $\LND(S;D)=\min\limits_{Q\in D} \EMD(\PDA(S;k),\PDA(Q;k))$ is the shortest $L_\infty$-based $\EMD$ distance from $S$ to its nearest neighbour $Q$ in the given {crystal} dataset $D$. 
\edfn
\end{dfn}

If $S$ is already contained in the dataset $D$, then $\LND(S;D)=0$, so $S$ cannot be considered novel.
{Conversely, if $\LND(S;D)=0$ then $S$ highly likely belongs to $S$, because $\PDD(S;100)$ distinguished all non-duplicate periodic crystals in the CSD.}
\myskip

$\LND(S;D)$ is based on $\PDA$s instead of $\PDD$s because distances to $k$-th neighbours in $\PDD(S;k)$ asymptotically increase as $\PPC(S)\sqrt[3]{k}$ by Theorem~\ref{thm:asymptotic}.
If crystals $S,Q$ have $\PPC(S)\neq\PPC(Q)$, the Chebyshev distance $L_\infty$ between rows of $\PDD$s equals the largest {absolute difference of $i$-th distances, which likely happens for $i=k$. 
Hence, subtracting $\PPC(S)\sqrt[3]{k}$ in Definition~\ref{dfn:ADA} makes any metric on $\PDA$s more informative than on $\PDD$s.
If a newly synthesised periodic crystal $S$ is a near-duplicate of some known {$Q\in D$, then} $\LND(S;D)$ is small as justified below.
 
\begin{thm}[continuity of $\LND$, {\cite[Theorem~6]{widdowson2025geographic}}]
\label{thm:LND}
For periodic point sets $S,Q\subset\R^3$, if $S$ is obtained from $Q$ in a dataset $D$ by perturbing every point of $Q$ up to $\ep<r(Q)$, then $\LND(S;D)\leq2\ep$.
To get $S$ from a periodic point set $Q\in D$ with $\LND(S;D)<2r(Q)$, a point of $Q$ should be perturbed by at least $0.5\LND(S;D)$.
\ethm
\end{thm}

\cite[section~3]{widdowson2025geographic} describes how the 43 materials reported by Berkeley's A-lab \cite{peplow2023robot} can be automatically positioned relative to the ICSD and the MP within the full Crystal Isometry Space $\CIMS(\R^3)=\bigcup\limits_{m\geq 1}\CIMS(\R^3;m)$ in seconds, see
Table~\ref{tab:neighbour-times}.

\begin{table}[!h]
\begin{tabular}{l|ll}
Stage & ICSD (s) & MP (s) \\ \hline
Binary search on $\ADA(S;100)$ in the full database & 3.023 & 2.450 \\
$\PDA(Q;100)$ for 100 neighbours $Q$ of $S$ found by $\ADA$ & 5.272 & 5.990 \\
$\EMD$ on PDAs for 100 neighbours $Q$ found by $\ADA$ & 0.535 & 0.742 \\
Elemental Mover's Distance (ElMD) for 100 neighbours & 9.534 & 9.737
\end{tabular}
\caption{Time (seconds) to complete each stage of the process of finding nearest neighbours in the ICSD and Materials Project for 43 A-lab crystals \cite{peplow2023robot} on a modest desktop computer.
{The binary search used 6-cores for multiprocessing}.}
\label{tab:neighbour-times}
\end{table}

Two A-lab crystals were found to already exist in the ICSD with the same composition: \ce{KNaP6(PbO3)8} matched ICSD 182501 reported in 2011 \cite{azrour2011rietveld}, and \ce{MnAgO2} matched ICSD 670065 reported as a hypothetical structure in 2015 \cite{cerqueira2015identification}.
In particular, \ce{MnAgO2} was one of three crystals that the later rebuttal said was synthesized successfully \cite{leeman2024challenges}, and they go on to state that the material was first reported in 2021 \cite{griesemer2021high} (ICSD 139006), after the snapshot used to train the GNoME, and so was not included in the original training data and could be considered a success. 
Our findings show this crystal did in fact exist in the ICSD prior to the 2021 snapshot.
The pre-existing version of this crystal was not found by \cite{leeman2024challenges} using a unit cell search because the unit cell of ICSD 670065 significantly differs from that of the A-lab version or ICSD 139006, with the former listing its space group as A 2/m and the latter two having space group C 2/m, see Fig.~\ref{fig:MnAgO2}. 
Such cell-based search can always miss near-duplicates as in Fig.~\ref{fig:noise_double_peaks}~(left), while continuous invariants independent of a unit cell find near-duplicates despite disagreement on a space group, which breaks down under almost any noise.

\begin{figure}[h!]
\caption{\textbf{Left}: \ce{MnAgO2} synthesized by A-lab.
\textbf{Middle}: ICSD entry 670065 with the same composition and $\EMD=0.097\angstrom$ found by $\PDA(S;100)$ in \cite[Table~2]{widdowson2025geographic}, though its unit cell is very different from the cell of \ce{MnAgO2}.
\textbf{Right}: another ICSD entry 139006 from 2021 matched by \cite{leeman2024challenges} and found by unit cell search, but is more distant from \ce{MnAgO2} by $\EMD=0.368\angstrom$ on invariants $\PDA(S;100)$.}
\includegraphics[height=36mm]{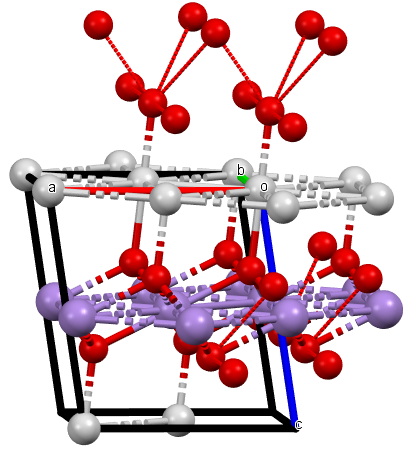}
\hspace*{1mm}
\includegraphics[height=36mm]{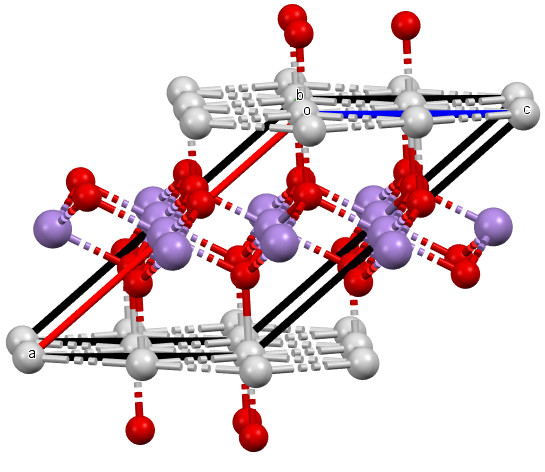}
\includegraphics[height=36mm]{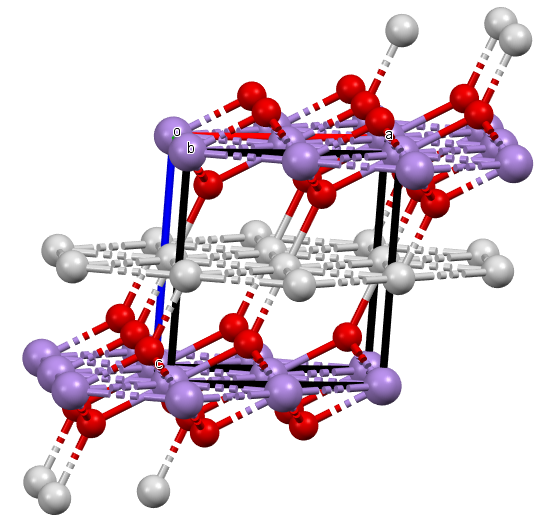}
\label{fig:MnAgO2}
\end{figure}

Aside from the two structures above, all other A-lab crystals were found to have a geometric near-duplicate in the ICSD with a different composition.
Many of these near-duplicates involve the substitution of only one atom, replacing a disordered site with a fully ordered one or adjusting the occupancy ratios of atoms at a site. 
\myskip

These structural analogues of A-lab’s claimed crystals are not surprising, as the GNoME  \cite{merchant2023scaling} used atomic substitution on existing crystals to generate potential new ones without substantially changing the atomic geometry.
The fact that pre-existing structures in the ICSD were missed by the later rebuttal \cite{leeman2024challenges} suggests that a more robust method is needed for comparing structures in the aid of materials discovery.
\myskip

In conclusion, crystals were classified for hundreds of years almost exclusively by discrete tools such as space groups or by using reduced cells, which are unique in theory.
Fig.~\ref{fig:noise_double_peaks}~(left) showed that any known crystal can be disguised by changing a unit cell, shifting atoms a bit, changing chemical elements, and then claiming them as `new'.
\myskip

Artificial near-duplicates threaten the integrity of experimental databases \cite{chawla2023crystallography}, which are skewed by previously undetectable near-duplicates.
These challenges motivated the stronger question (\emph{if different, by how much?}) that
was formalised by Lipschitz continuity in condition~\ref{pro:l-periodic}(d) and inspired the research leading to this book.
\myskip

Our future paper \cite{widdowson2025higher} extends the $\PDD$ to stronger higher-order invariants.
Since the ultra-fast $\PDD$ distinguished all non-duplicate crystals among all experimental materials in the periodic case, these invariants already parametrise the known `universe' of all existing crystals as `shiny stars', while all not yet discovered crystals remain hidden in empty spots on the same map.
Fig.~\ref{fig:heatmap_C_allotropes} shows an example for two invariant coordinates.

\begin{figure}
\centering
\includegraphics[width=\textwidth]{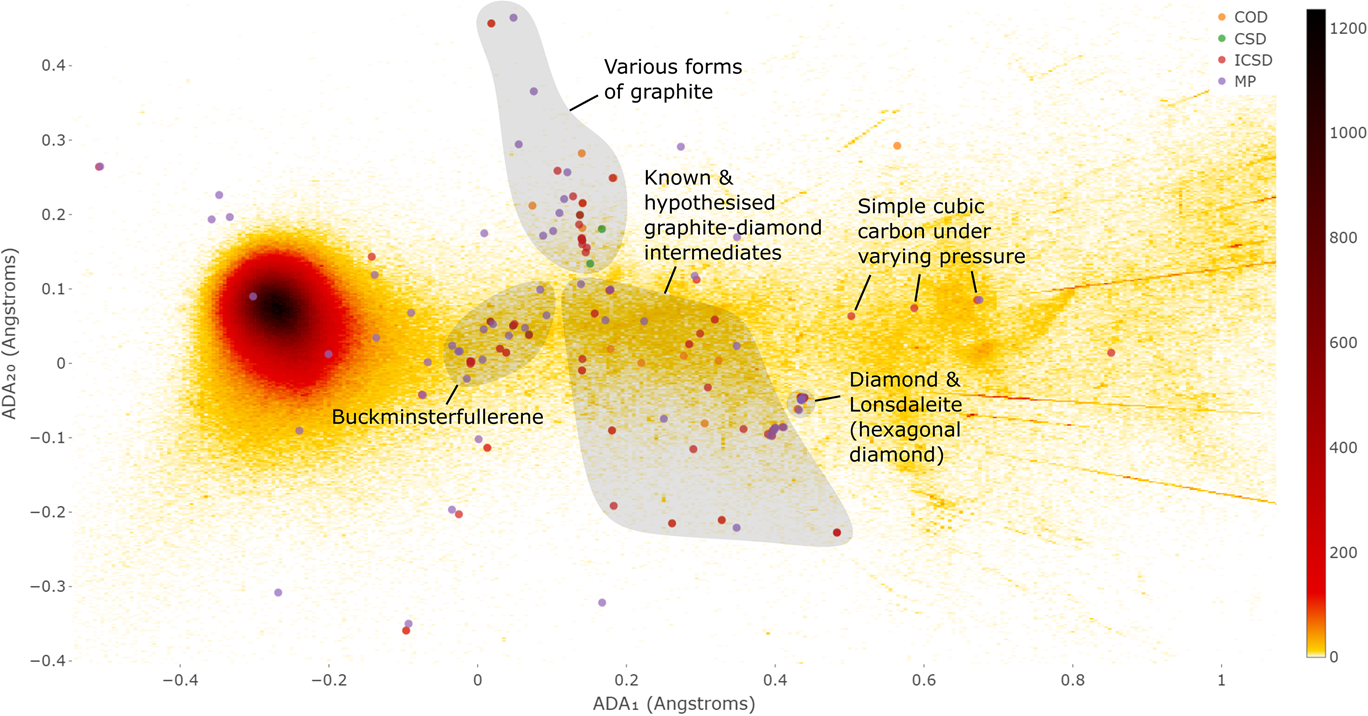}
\caption{All carbon allotropes (materials of pure carbon) on the heat map of four major databases.}
\label{fig:heatmap_C_allotropes}
\end{figure}

\bibliographystyle{plain}
\bibliography{Geometric-Data-Science-book}

%
%
%

\chapter{Complete and continuous isosets of all periodic point sets in $\R^n$}
\label{chap:isosets} 

\abstract{
This chapter adapts the general Geo-Mapping Problem to periodic point sets under rigid motion in $\R^n$.
We introduce a complete isoset invariant with a Lipschitz continuous metric.
For a fixed dimension $n$, the isoset is computable in polynomial time of the input size, while the distance metric is approximated with a constant factor in polynomial time.
The isoset distinguished all known homometric crystals that have identical diffraction patterns and detected several pairs of exact (but unlikely) mirror images in the Cambridge Structural Database of experimental materials.
}

\section{Geo-mapping for periodic point sets under rigid motion in $\R^n$}
\label{sec:periodic_rigid}

This chapter follows papers \cite{anosova2021isometry,anosova2025recognition}. 
In comparison with Problem~\ref{pro:l-periodic} in the previous chapter, Problem~\ref{pro:periodic_rigid} covers all (not only generic) periodic point sets and asks for completeness under rigid motion, which distinguishes mirror images.

\begin{pro}[invariants of periodic point sets under rigid motion in $\R^{n}$]
\label{pro:periodic_rigid}
Design an invariant $I$ on the Crystal Rigid Space $\CRIS(\R^n;m)$
satisfying the conditions below.
\smallskip

\noindent
\tb{(a)} 
\emph{Completeness:}
any periodic point sets $S,Q\subset\R^n$ are related by rigid motion in $\R^{n}$ if and only if $I(S)=I(Q)$.
\myskip

\noindent
\tb{(b)} 
\emph{Reconstruction:}
any periodic point set $S\subset\R^n$ is reconstructable from its invariant value 
$I(S)$, uniquely under rigid motion in $\R^n$. 
\myskip

\noindent
\tb{(c)} 
\emph{Metric:} 
there is a distance $d$ on the Crystal Rigid Space $\CRIS(\R^n;m)$ satisfying all metric axioms in Definition~\ref{dfn:metrics}(a). 
\myskip

\noindent
\tb{(d)} 
\emph{Continuity:} 
there is a constant $\la>0$, such that, for all sufficiently small $\ep>0$ and periodic point set $S,Q\subset\R^n$ if $Q$ is obtained by perturbing every point of $S$ up to Euclidean distance $\ep$, then $d(I(S),I(Q))\leq\la\ep$.
\myskip

\noindent
\tb{(e)} 
\emph{Computability:} 
for a fixed dimension $n$, the invariant $I(S)$ of any periodic point set $S$ and a reconstruction of $S\subset\R^{n}$ from $I(S)$ can be computed in times that depend polynomially on the motif size $m$ of $S$, while the metric $d(I(S),I(Q))$ can be approximated in polynomial time of the maximum motif size of $S,Q$.
\epro
\end{pro}

One remaining limitation in Problem~\ref{pro:periodic_rigid} is an approximate algorithm for a Lipschitz continuous metric, which will be improved to an exact one in future work.
\myskip

Fig.~\ref{fig:CRISP-infinitely-many-layers} visualises the practical consequences of discontinuous cell-based representations.
Starting from any periodic crystal with a unit cell of (say) $m$ atoms, one can extend this cell to get an arbitrarily large motif of $2m$ or $3m$ or $sm$ atoms for any integer scale $s\geq 2$.
If $s$ is not prime, then any extension can be done in different geometric ways.
For example, an extension of unit cell $U(\vec v_1,\vec v_2,\vec v_3)\subset\R^3$ by factor $s=8$, can be by factor $8$ in the direction of $\vec v_i$ for $i=1,2,3$, or by factor $4$ in the direction of $\vec v_1$ and by factor 2 in the direction of $\vec v_2$ or $\vec v_3$, or by factor in each of the directions $\vec v_2, \vec v_2, \vec v_3$.

\begin{figure}[h!]
\centering
\includegraphics[width=\textwidth]{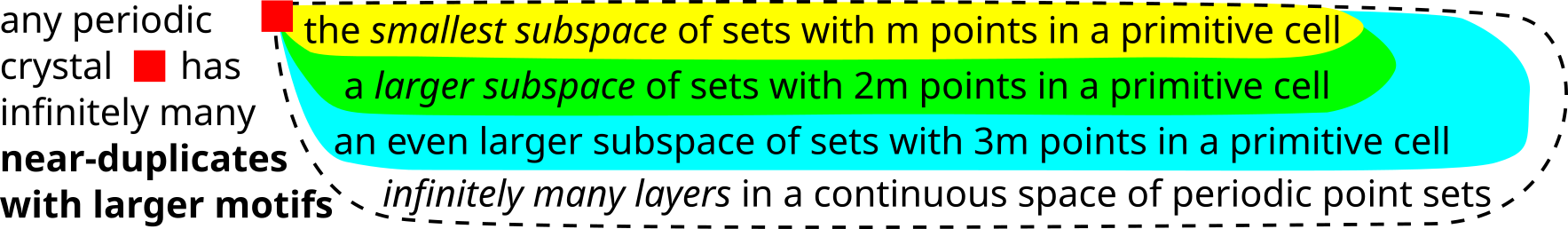}
\caption{The full Crystal Rigid Space $\CRIS(\R^3)$ is `totally singular' in the sense that every periodic structure $S$ with $m$ points in a minimal cell is infinitesimally close to infinitely many subspaces of periodic structures with $km$ points in a minimal cell for any integer factor $k=1,2,3,\dots$}
\label{fig:CRISP-infinitely-many-layers}
\end{figure}

To solve Problem~\ref{pro:periodic_rigid}, we introduce an isoset $I(S)$ consisting of local clusters around points $p$ in a motif of $S$.
Each cluster of a point $p$ is considered under rotations from the group $\SO(\R^n)$.
The completeness of isosets under rigid motion in $\R^n$ essentially reduces Problem~\ref{pro:periodic_rigid} to several subproblems for finite clouds with fixed centres under rotations from $\SO(\R^n)$.
Then we define a continuous metric on isosets in several steps.
\myskip

The first step introduces a boundary tolerant metric $\BT$ on local clusters around points of a periodic set $S$, which continuously changes when points cross a cluster boundary.
This discontinuity at the boundary can be formally resolved by an extra factor, which smoothly goes down to 0 depending on an extra parameter. 
Without using extra parameters, the metric $\BT$ will be expressed in terms of simpler distances. 
\medskip

The second step uses the Earth Mover's Distance from Definition~\ref{dfn:EMD} to extend the boundary tolerant metric $\BT$ to complete invariants \cite{anosova2021isometry} that are weighted distributions of local clusters under rotations.
The resulting metric on periodic sets in $\R^n$ is approximated with a factor $\eta$, e.g. $\eta\approx 4$ in $\R^3$, in a time depending polynomially on the input size.
\medskip

The third step proves the metric axioms and continuity $d(S,Q)\leq 2\ep$, which also has practical importance.
Indeed, if $d(S,Q)$ is approximated by a value $d$ with a factor $\eta$, we get the lower bound $\ep\geq\frac{d}{2\eta}$ for the maximum displacement $\ep$ of points.
Such a lower bound is impossible to guarantee by analysing only finite subsets, which can be very different in identical periodic sets, see Fig.~\ref{fig:crystal_finite_subsets}.

\begin{figure}[h!]
\centering
\includegraphics[width=\textwidth]{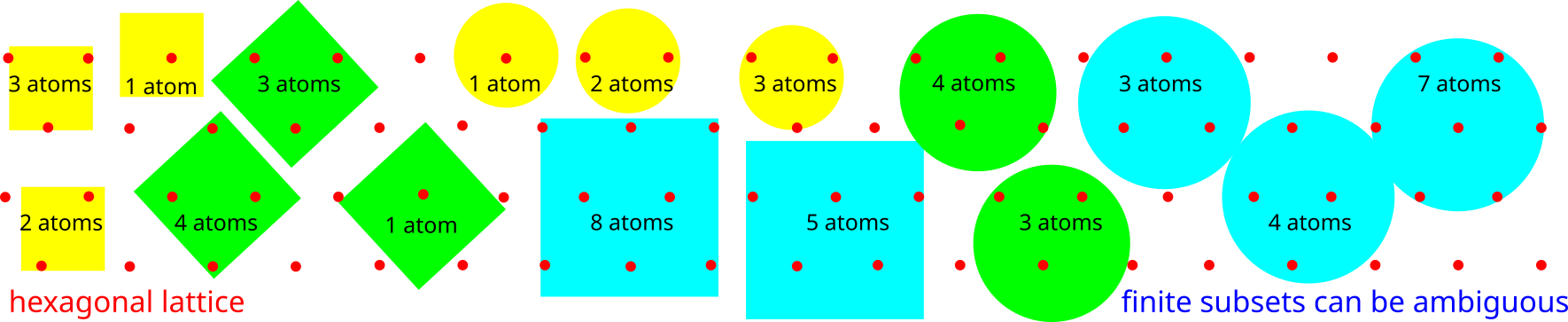}
\caption{Comparing periodic crystals by finite subsets is hard to justify because any periodic point set 
has many non-isometric subsets of different numbers of points within a box or a ball of a fixed size.}
\label{fig:crystal_finite_subsets}
\end{figure}

\section{Isotree of local clusters of points in a periodic set in $\R^n$}
\label{sec:local_clusters}

This section defines the complete invariant \cite{anosova2021isometry}
based on local clusters and their symmetry groups, which were previously studied in \cite{delone1976local,dolbilin1998multiregular}.

\index{$m$-regular periodic point set}

\begin{dfn}[\emph{global} clusters and \emph{$m$-regular} periodic sets]
\label{dfn:m-regular}
\tb{(a)}
For any point $p$ in a periodic set $S\subset\R^n$, the \emph{global cluster} is $C(S,p)=\{\vec q - \vec p \vl q\in S\}$.
For any $p,q\in\R^n$, let the set $\Or(\R^n;p,q)$ consist of all isometries of $\R^n$ that map $p$ to $q$.
\myskip

\nt
\tb{(b)}
Global clusters $C(S,p)$ and $C(S,q)$ are called \emph{isometric} if there is $f\in\Or(\R^n;p,q)$ such that $f(S)=S$. 
A periodic point set $S\subset\R^n$ is called \emph{$m$-regular} if all global clusters of $S$ form exactly $m\geq 1$ isometry classes.
\edfn
\end{dfn}

For any point $p\in S$, its global cluster is a view of $S$ from the position of $p$.
We view all astronomical stars in the universe $S$ from our Earth at $p$.
Any lattice is 1-regular since all its global clusters are related by translations.
Though global clusters $C(S,p),C(S,q)$ at any different points $p,q\in S$ contain the same set $S$, they may not match under the translation shifting $p$ to $q$.
The global clusters are infinite, hence distinguishing them under isometry is not easier than the original periodic sets.
However, the $m$-regularity of a periodic set can be checked in terms of local $\al$-clusters below.

\index{local cluster}

\begin{dfn}[\emph{local $\al$-clusters} $C(S,p;\al)$ and \emph{symmetry groups} $\sym(S,p;\al)$]
\label{dfn:local_cluster}
For a point $p$ in a periodic point set $S\subset\R^n$ and any $\al\geq 0$, the local \emph{$\al$-cluster} $C(S,p;\al)$ is the set of all vectors $\vec q - \vec p$ such that $q\in S$ and $|\vec q - \vec p|\leq\al$.
Let the group $\Or(\R^n;p)$ consist of all isometries that fix $p$.
If $p=0$ is the origin, $\Or(\R^n;0)$ is the usual orthogonal group.
The \emph{symmetry} group $\sym(S,p;\al)$ consists of all isometries $f\in\Or(\R^n;p)$ that map $C(S,p;\al)$ to itself so that $f(p)=p$.
\edfn
\end{dfn}

Fig.~\ref{fig:alpha-clusters}~(left) shows the 1-regular periodic set $S_1\subset\R^2$ whose all points (close to vertices of square cells) have isometric global clusters related by translations and rotations through $90^\circ,180^\circ,270^\circ$. 
The 2-regular periodic set $S_2$ has extra points at the centers of all square cells.
The local $\al$-clusters around these centers are not isometric to $\al$-clusters around the points close to cell vertices for any $\al\geq3\sqrt{2}$. 
\medskip

The 1-regular periodic point set $S_1$ in Fig.~\ref{fig:alpha-clusters} for any $p\in S_1$ has the symmetry group $\sym(S_1,p;\al)=\Or(\R^2)$ for $\al\in[0,4)$. 
Then $\sym(S_1,p;\al)$ stabilizes as $\Z^2$ with one reflection for $\al\geq 4$ as soon as $C(S_1,p;\al)$ includes one more point.

\begin{figure}[h!]
\includegraphics[height=55mm]{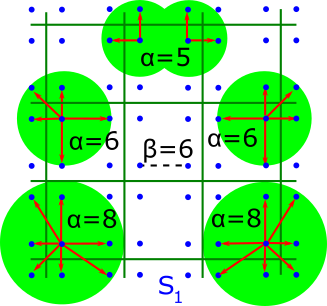}
\hspace*{2mm}
\includegraphics[height=55mm]{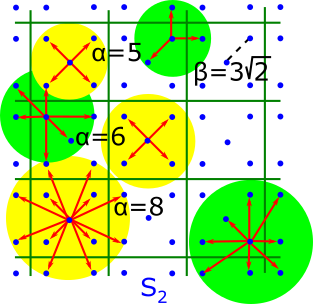}
\caption{
\textbf{Left}: 
in $\R^2$, the periodic point set $S_1$ has the square unit cell $[0,10)^2$ containing the four points $(2,2),(2,8),(8,2),(8,8)$, so $S_1$ isn't a lattice, but is 1-regular by Definition~\ref{dfn:m-regular}, and $\be(S_1)=6$.
All local $\al$-clusters of $S_1$ are isometric, shown by red arrows for $\al=5,6,8$, see Definition~\ref{dfn:local_cluster}.
\textbf{Right}: 
$S_2$ has the extra point $(5,5)$ in the center  of the cell $[0,10)^2$ and is 2-regular with $\be(S_2)=3\sqrt{2}$, so $S_2$ has green and yellow isometry types of $\al$-clusters. }
\label{fig:alpha-clusters}
\end{figure}

For any periodic set $S$, if $\al$ is smaller than the minimum distance between all points of $S$, then any $\al$-cluster $C(S,p;\al)$ is one point $\{p\}$.
Its symmetry group  consists of all isometries fixing the centre $p$, so $\sym(S,p;\al)=\Or(\R^n;p)$.
When $\al$ is increasing, the $\al$-clusters $C(S,p;\al)$ become larger and there can be fewer (not more) isometries $f\in\Or(\R^n;p)$ that bijectively map $C(S,p;\al)$ to itself.
\myskip

So the group $\sym(S,p;\al)$ can become smaller (not larger) and eventually stabilises (stops changing), which is formalised in Definition~\ref{dfn:stable_radius}.
This stabilization uses the bridge length extending the idea of a longest edge in a Minimum Spanning Tree to a periodic set $S$, which does not easily reduce to the finite case \cite{mcmanus2025computing}.

\index{bridge length}

\begin{dfn}[\emph{bridge length} $\be(S)$]
\label{dfn:bridge_length}
For a periodic point set $S\subset\R^n$, the \emph{bridge length} is a minimum distance $\be(S)>0$ such that any $p,q\in S$ can be connected by a sequence of points $p_0=p,p_1,\dots,p_k=q$ such that 
any two successive points $p_{i-1},p_{i}$ are close so that $|\vec p_{i-1} - \vec p_{i}|\leq\be(S)$ for $i=1,\dots,k$.
\edfn
\end{dfn}

The seminal result in \cite[Theorem~1.3]{dolbilin1998multiregular} described how a family of clusters determines a periodic point set under isometry.
These results motivated the \emph{isotree}, \emph{stable} radius, and \emph{isoset} in Definitions~\ref{dfn:isotree}, \ref{dfn:stable_radius},~\ref{dfn:isoset}, respectively, leading to the isometry classification of periodic point sets via isosets in Theorem~\ref{thm:isoset_complete}.
\myskip

The \emph{isotree} in Definition~\ref{dfn:isotree} is inspired by a clustering dendrogram because points of $S$ split into isometry classes of $\al$-clusters at variable radii $\al$, not at a fixed $\al$.
\myskip

Fig.~\ref{fig:1-regular_set_isotree} 
illustrates the isotrees of the periodic sets $S_1,S_2$ in Fig.~\ref{fig:alpha-clusters}, as defined below.

\begin{figure}[h!]
\includegraphics[width=\textwidth]{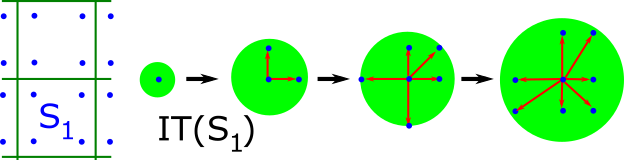}
\myskip

\includegraphics[width=\textwidth]{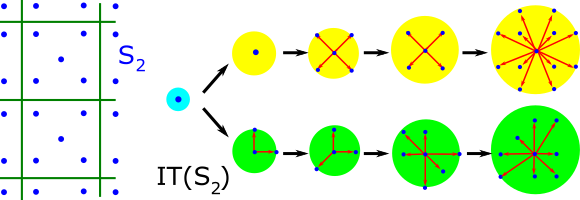}
\caption{
\textbf{Top}: the isotree $\IT(S_1)$ from Definition~\ref{dfn:isotree} of the 1-regular set $S_1$ in Fig.~\ref{fig:alpha-clusters} for any $\al\geq 0$ has one isometry class of $\al$-clusters under rotations.
\textbf{Bottom}: the isotree $\IT(S_2)$ of the 2-regular set $S_2$ in Fig.~\ref{fig:alpha-clusters} stabilizes with two non-isometric classes of $\al$-clusters for $\al\geq 4$.}
\label{fig:1-regular_set_isotree}
\end{figure}

\index{isotree}

\begin{dfn}[\emph{isotree} $\IT(S)$ of \emph{$\al$-partitions}]
\label{dfn:isotree}
\tb{(a)}
Fix a periodic point set $S\subset\R^n$.
Points $p,q\in S$ are \emph{$\al$-equivalent} if their $\al$-clusters $C(S,p;\al)$ and $C(S,q;\al)$ can be related by an isometry that matches their centres.
The \emph{isometry class} $[C(S,p;\al)]$ consists of all $\al$-clusters isometric to $C(S,p;\al)$.
The \emph{$\al$-partition} $P(S;\al)$ is the splitting of $S$ into $\al$-equivalence classes of points.
Call a value $\al$ \emph{singular} if $P(S;\al)\neq P(S;\al-\ep)$ for any small enough $\ep>0$.
\myskip

\nt
\tb{(b)}
Represent each $\al$-equivalence class by a vertex of the \emph{isotree} $\IT(S)$.
The top vertex of $\IT(S)$ represents the $0$-equivalence class coinciding with $S$.
For any successive singular values $\al<\al'$, connect the vertices representing any classes $A\in P(S;\al)$ and $A'\in P(S;\al')$ such that $A'\subset A$ by an edge of the length $\al'-\al$ in $\IT(S)$.
\edfn
\end{dfn}

For any periodic point set $S\subset\R^n$, the root vertex of $\IT(S)$ at $\al=0$ is the single class $S$, because any 0-cluster $C(S,p;0)$ of a point $p\in S$ consists only of its centre $p$.
When $\al$ is increasing, $\al$-clusters $C(S,p;\al)$ include more points and hence may not be isometric.
In other words, any $\al$-equivalence class from $P(S;\al)$ may split into two or more classes, which cannot merge at any larger $\al'$.
Branched vertices of $\IT(S)$ correspond to the values of $\al$ when an $\al$-equivalence class is split into subclasses for $\al'$ slightly larger than $\al$.
So the number $|P(S;\al)|$ of $\al$-equivalence is non-decreasing in $\al$, see Fig.~\ref{fig:4-regular_set_isotree}.
The isotree $\IT(S)$ is continuously parametrised by $\al\geq 0$ and is visualised as a tree of $\al$-equivalence classes in Fig.~\ref{fig:4-regular_set_isotree}.
due to Lemma~\ref{lem:isotree}.
\myskip

\begin{figure}[h]
\includegraphics[width=\textwidth]{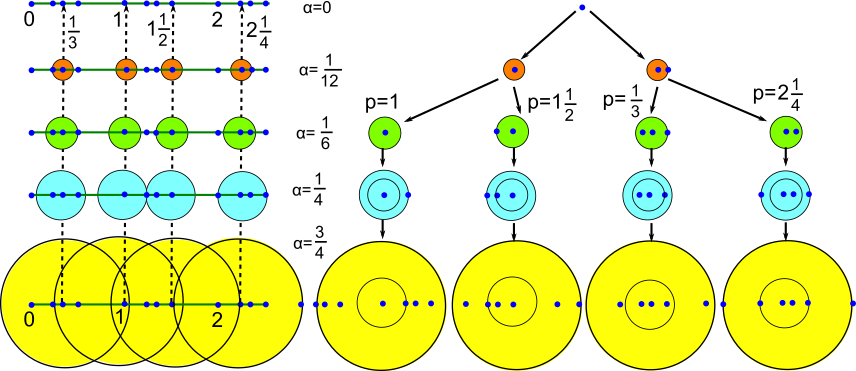}
\caption{\textbf{Left}: the 1-dimensional set $S_4=\{0,\frac{1}{4},\frac{1}{3},\frac{1}{2}\}+\Z$ has four points in the unit cell $[0,1)$ and is 4-regular by Definition~\ref{dfn:m-regular}.
\textbf{Right}: the colored disks show $\al$-clusters in the line $\R$ with radii $\al=0,\frac{1}{12},\frac{1}{6},\frac{1}{4},\frac{3}{4}$ and represent points in the isotree $\IT(S_4)$ from Definition~\ref{dfn:isotree}.
}
\label{fig:4-regular_set_isotree}
\end{figure}

Our proofs of Lemmas~\ref{lem:isotree}, \ref{lem:local_extension},~\ref{lem:global_extension}, and \ref{lem:stable_radius} were peer-reviewed but were published only in appendices of \cite{anosova2025recognition} online, so we include them for completeness.

\begin{lem}[properties of $\al$-partitions and isotrees, {\cite[Lemma 4.1]{anosova2025recognition}}]
\label{lem:isotree} 
The isotree $\IT(S)$ from Definition~\ref{dfn:isotree} 
has the following properties.
\medskip

\noindent
\tb{(a)}
for $\al=0$, the $\al$-partition $P(S;0)$ consists of one class.
\medskip

\noindent
\tb{(b)}
if $\al<\al'$, then $\sym(S,p;\al')\subseteq\sym(S,p;\al)$ for any point $p\in S$.
\medskip

\noindent
\tb{(c)}
If $\al<\al'$, the $\al'$-partition $P(S;\al')$ \emph{refines} $P(S;\al)$, i.e. any $\al'$-equivalence class from $P(S;\al')$ is included into an $\al$-equivalence class from the partition $P(S;\al)$.
So the cluster count $|P(S;\al)|$ is a non-strictly increasing integer-valued function of 
$\al$, i.e. $|\iso(S;\al)|\leq|\iso(S;\al')|$ for $\al<\al'$.
\elem
\end{lem}
\begin{proof}
\tb{(a)}
Let $\al\geq 0$ be smaller than the minimum distance $2r(S)$ betweens any points of $S$.
Then any cluster $C(S,p;\al)$ is the single-point set $\{p\}$.
All these 1-point clusters are isometric to each other.
So $|P(S;\al)|=1$ for $\al<2r(S)$.
\medskip

\noindent
\tb{(b)}
For any $p\in S$, the inclusion of clusters $C(S,p;\al)\subseteq C(S,p;\al')$ implies that any isometry $f\in\Or(\R^n;p)$ that isometrically maps the larger cluster $C(S,p;\al')$ to itself also maps the smaller cluster $C(S,p;\al)$ to itself.
Hence any element of $\sym(S,p;\al')\subseteq\Or(\R^n;p)$ belongs to $\sym(S,p;\al)$. 
\medskip

\noindent
\tb{(c)}
If points $p,q\in S$ are $\al'$-equivalent at the larger radius $\al'$, i.e. the clusters $C(S,p;\al')$ and $C(S,q;\al')$ are related by an isometry from $\Or(\R^n;p,q)$, then $p,q$ are $\al$-equivalent at the smaller radius $\al$.
Hence any $\al'$-equivalence class of points in $S$ is a subset of an $\al$-equivalence class in $S$.
\end{proof}

The $\al$-clusters of the periodic sequence $S_4\subset\R$ in Fig.~\ref{fig:4-regular_set_isotree} are intervals in $\R$, shown as disks for better visibility.
In Fig.~\ref{fig:4-regular_set_isotree}, the initial 1-point class persists until $\al=\frac{1}{12}$, when all points $p\in S_4$ are split into two classes: one represented by 1-point cluster $\{p\}$ for $p\in\{0,\frac{1}{2}\}+\Z$, and another represented by 2-point clusters $\{p,p+\frac{1}{12}\}$, $p\in\{\frac{1}{4},\frac{1}{3}\}+\Z$.
\medskip

The sequence $S_4$ has four $\al$-equivalence classes for any radius $\al\geq\frac{1}{6}$.
For any point $p\in\Z\subset S_4$, the symmetry group $\sym(S_4,p;\al)=\Z_2$ is generated by the reflection in $p$ for $\al\in[0,\frac{1}{4})$.
For all $p\in S_4$, the symmetry group $\sym(S_4,p;\al)$ is trivial for $\al\geq\frac{1}{4}$.
\medskip

Lemmas~\ref{lem:local_extension} and~\ref{lem:global_extension} are key steps towards a complete classification of periodic point sets under isometry and rigid motion in Theorem~\ref{thm:isoset_complete}. 

\begin{lem}[local extension]
\label{lem:local_extension}
Let $S,Q\subset\R^n$ be periodic point sets and $\sym(S,p;\al-\be)=\sym(S,p;\al)$ for some point $p\in S$ and $\al>\be$.
Assume that there is an isometry $g\in\Or(\R^n;p,q)$ such that $g(C(S,p;\al))=C(Q,q;\al)$.
Let $f\in\Or(\R^n;p,q)$ be any isometry such that $f(C(S,p;\al-\be))=C(Q,q;\al-\be)$.
Then $f$ isometrically maps the larger clusters: $f(C(S,p;\al))=C(Q,q;\al)$.
\bs
\end{lem}
\begin{proof}
The composition $h=f^{-1}\circ g$ fixes $p$ and isometrically maps $C(S,p;\al-\be)$ to itself, so $h\in\sym(S,p;\al-\be)$.
The condition $\sym(S,p;\al-\be)=\sym(S,p;\al)$ implies that $h\in\sym(S,p;\al)$, so the isometry $h\in\Or(\R^n;p)$ isometrically maps the larger cluster $C(S,p;\al)$ to itself.
Then the given isometry $f=g\circ h^{-1}$ isometrically maps $C(S,p;\al)$ to $f(C(S,p;\al))=g(C(S,p;\al))=C(Q,q;\al)$. 
\end{proof}

\begin{lem}[global extension]
\label{lem:global_extension}
Let periodic point sets $S,Q\subset\R^n$ have a common stable radius  $\al$ satisfying Definition~\ref{dfn:stable_radius} for an upper bound $\be\geq\be(S),\be(Q)$.
Let $I(S;\al)=I(Q;\al)$ and $p\in S$, $q\in Q$ be any points with an isometry $f\in\Or(\R^n;p,q)$ such that $f(C(S,p;\al))=C(Q,q;\al)$.
Then $f(S)=Q$.
\elem
\end{lem}
\begin{proof}
To show that $f(S)\subset Q$, it suffices to check that the image $f(a)$ of any point $a\in S$ belongs to $Q$.
By Definition~\ref{dfn:bridge_length} the points $p,a\in S$ are connected by a sequence of points $p=a_0,a_1,\dots,a_k=a\in S$ such that the distances $|a_{i-1}- a_{i}|$ between any successive points have the upper bound $\be$ for $i=1,\dots,k$.
\medskip

We will prove that $f(C(S,a_k;\al))=C(Q,f(a_k);\al)$ by induction on $k$, where the base $k=0$ is given.
The induction step below goes from $i$ to $i+1$. 
\medskip

The ball $\bar B(a_i;\al)$ contains the smaller ball $\bar B(a_{i+1};\al-\be)$ around the closely located center $a_{i+1}$.
Indeed, since $|a_{i+1}-a_i|\leq\be$, the triangle inequality for the Euclidean distance implies that any point $a'_{i+1}\in\bar B(a_{i+1};\al-\be)$ with $|a'_{i+1}-a_i|\leq\al-\be$ satisfies 
$$|a'_{i+1}-a_i|\leq |a'_{i+1}-a_{i+1}|+|a_{i+1}-a_i|\leq(\al-\be)+\be=\al, \text{ so } 
\bar B(a_{i+1};\al-\be)\subset\bar B(a_i;\al).$$
Then the inductive assumption 
$f(C(S,a_i;\al))=C(Q,f(a_i);\al)$ gives 
$$\begin{array}{l}
f(C(S,a_{i+1};\al-\be))=f(C(S,a_i;\al))\cap f(\bar B(a_{i+1};\al-\be))= \\ 
C(Q,f(a_i);\al)\cap \bar B(f(a_{i+1});\al-\be)=C(Q,f(a_{i+1});\al-\be)
\end{array}.$$

Due to $I(S;\al)=I(Q;\al)$, the isometry class of $C(S,a_{i+1};\al)$ equals an isometry class of $C(Q,b_{i+1};\al)$ for some point $b_{i+1}\in Q$, i.e. there is an isometry $g\in\Or(\R^n;a_{i+1},b_{i+1})$ such that $g(C(S,a_{i+1};\al))=C(Q,b_{i+1};\al)$.
\myskip

Since $f\circ g^{-1}\in\Or(\R^n;b_{i+1})$ isometrically maps $C(Q,b_{i+1};\al-\be)$ to $C(Q,f(a_{i+1});\al-\be)$, the points $b_{i+1},f(a_{i+1})\in Q$  are in the same $(\al-\be)$-equivalence class of $Q$.
\medskip

By condition~(\ref{dfn:stable_radius}a), the splitting of the periodic point set $Q\subset\R^n$ into $\al$-equivalence classes coincides with its splitting into $(\al-\be)$-equivalence classes. 
Hence, the points $b_{i+1},f(a_{i+1})\in Q$ are in the same $\al$-equivalence class of $Q$.
Then $C(Q,f(a_{i+1});\al)$ is isometric to $C(Q,b_{i+1};\al)=g(C(S,a_{i+1};\al))$.\medskip

Now we can apply Lemma~\ref{lem:local_extension} for $p=a_{i+1},q=f(a_{i+1})$ and conclude that the given isometry $f$, which satisfies $f(C(S,a_{i+1};\al-\be))=C(Q,f(a_{i+1});\al-\be)$, isometrically maps the larger clusters: $f(C(S,a_{i+1};\al))=C(Q,f(a_{i+1});\al)$.
\myskip

The induction step is finished.
The inclusion $f^{-1}(Q)\subset S$ is proved similarly.
\end{proof}

\section{A complete isoset of periodic point sets under rigid motion in $\R^n$}
\label{sec:isosets}

This section shows that, for any periodic point set $S\subset\R^n$, the $\al$-partition of $S$ stabilises, which allows us to form a complete invariant at a stable radius defined below.

\index{stable radius}

\begin{dfn}[the \emph{minimum stable radius} $\al(S)$]
\label{dfn:stable_radius}
Let $S\subset\R^n$ be a periodic point, $\be\geq\be(S)$ be an upper bound of the bridge length $\be(S)$ from Definition~\ref{dfn:bridge_length}.
A radius $\al\geq\be$ is called \emph{stable} if the following conditions hold:
\medskip

\noindent
\tb{(a)} 
the $\al$-partition $P(S;\al)$ equals the $(\al-\be)$-partition $P(S;\al-\be)$;
\medskip

\noindent
\tb{(b)} 
the groups stabilise so that
$\sym(S,p;\al)=\sym(S,p;\al-\be)$ for any $p\in S$, i.e.
any isometry $f\in\sym(S,p;\al-\be)$ preserves the larger cluster $C(S,p;\al)$.
\medskip

\noindent
A minimum value of a stable radius $\al$ satisfying \ref{dfn:stable_radius}(a,b) for $\be=\be(S)$ from Definition~\ref{dfn:bridge_length} is called \emph{the minimum stable radius} and denoted by $\al(S)$. 
\edfn
\end{dfn}

Due to the upper bounds in Lemma~\ref{lem:upper_bounds}(b,c), the minimum stable radius $\al(S)\geq 0$ exists and is achieved because $P(S;\al)$ and $\sym(S,p;\al)$ are continuous on the right (unchanged when $\al$ increases by a sufficiently small value). 
\medskip

Any $m$-regular periodic point set $S\subset\R^n$ has at most $m$ $\al$-equivalence classes, so the isotree $\IT(S)$ stabilises with maximum $m$ branches.
Though \ref{dfn:stable_radius}(b) is stated for all points $p\in S$ for simplicity, it suffices to check condition~\ref{dfn:stable_radius}(b) for points only from a finite motif $M$ of $S$ due to periodicity.
\medskip

\begin{lem}[all stable radii $\al\geq\al(S)$]
\label{lem:stable_radius}
If $\al$ is a stable radius of a periodic point set $S\subset\R^n$, then so is any  larger radius $\al'>\al$.
Then all stable radii form the interval $[\al(S),+\infty)$, where $\al(S)$ is the minimum stable radius of $S$.
\bs
\end{lem}
\begin{proof}
Due to Lemma~(\ref{lem:isotree}bc), conditions (\ref{dfn:stable_radius}ab) imply that the $\al'$-partition $P(S;\al')$ and the symmetry groups $\sym(S,p;\al')$ remain the same for all $\al'\in[\al-\be(S),\al]$, where $\be(S)$ is the bridle length. 
We need to show that they remain the same for any $\al'>\al$ and will apply Lemma~\ref{lem:global_extension} for $S=Q$ and $\be=\be(S)$.
\myskip

Let points $p,q\in S$ be $\al$-equivalent, i.e. there is an isometry $f\in\Or(\R^n;p,q)$ such that $f(C(S,p;\al))=C(S,q;\al)$.
By Lemma~\ref{lem:global_extension}, $f$ isometrically maps the full set $S$ to itself.
Then all larger $\al'$-clusters of $p,q$ are matched by $f$, so $p,q$ are $\al'$-equivalent  and $P(S;\al)=P(S,\al')$.
Similarly, any isometry $f\in\sym(S,p;\al)$ by Lemma~\ref{lem:global_extension} for $S=Q$ and $p=q$, isometrically maps the full set $S$ to itself.
Then $\sym(S,p;\al')$ coincides with $\sym(S,p;\al)$ for any $\al'>\al$. 
\end{proof}

All stable radii of $S$ form the interval $[\al(S),+\infty)$ by
Lemma~\ref{lem:stable_radius}.
The periodic sequence $S_4$ in Fig.~\ref{fig:4-regular_set_isotree} has $\be(S_4)=\frac{1}{2}$ and $\al(S)=\frac{3}{4}$ since the $\al$-partition and symmetry groups $\sym(S_4,p;\al)$ are stable for $\frac{1}{4}\leq\al\leq\frac{3}{4}$.
\myskip

Condition \ref{dfn:stable_radius}(b) doesn't follow from condition \ref{dfn:stable_radius}(a) due to the following example. 
Let $\La$ be the 2D lattice with the basis $(1,0)$ and $(0,\be)$ for $\be>1$.
Then $\be$ is the bridge length of $\La$. 
Condition \ref{dfn:stable_radius}(a) is satisfied for any $\al\geq 0$, because all points of any lattice are equivalent under translations.
\myskip

However, condition \ref{dfn:stable_radius}(b) fails for any $\al<\be+1$.
Indeed, the $\al$-cluster of the origin $(0,0)$ contains five points $(0,0),(\pm 1,0),(0,\pm\be)$, whose symmetries are generated by the two reflections in the axes $x,y$, but the $(\al-\be)$-cluster of the origin $(0,0)$ consists of its centre and has the symmetry group $\Or(\R^2)$. 
\medskip

It is possible that condition \ref{dfn:stable_radius}(b) might imply \ref{dfn:stable_radius}(a), but in practice it makes sense to verify \ref{dfn:stable_radius}(b) only after checking much simpler condition \ref{dfn:stable_radius}(a). Both conditions are essentially used in the proof of Isometry Classification Theorem~\ref{thm:isoset_complete}. 
\medskip

Conditions \ref{dfn:stable_radius}(ab) appeared in \cite{dolbilin1998multiregular} with different notations $\rho,\rho+t$.
Since many applied papers use $\rho$ for the physical density and have many types of bond distances, we replaced $t$ and $\rho+t$ with the bridge length $\be$ and radius $\al$, respectively, as for growing $\al$-shapes in Topological Data Analysis \cite{smith2024generic}. 
\medskip

Recall that the \emph{covering radius} $R(S)$ of a periodic point set $S\subset\R^n$ is the minimum radius $R>0$ such that $\bigcup\limits_{p\in S}\bar B(S;R)=\R^n$, or the largest radius of an open ball in the complement $\R^n\setminus S$.
For $m$-regular point sets in $\R^n$, an upper bound of $\al(S)$ can be extracted from \cite[Theorem 1.3]{dolbilin1998multiregular} whose proof motivated a stronger bound in Lemma~\ref{lem:upper_bounds}(c), see comparisons in Example~\ref{exa:upper_bounds}(c). 
\medskip

A periodic point set $S$ is \emph{locally antipodal} if the local cluster $C(S,p;2R(S))$ is centrally symmetric for any point $p\in S$, i.e. bijectively maps to itself under $\vec q\mapsto 2\vec p-\vec q$, $q\in\R^n$.
The important result in \cite[Theorem~1]{dolbilin2016uniqueness} says that all locally antipodal Delone sets, hence all periodic sets $S$, are globally antipodal, i.e. $S$ is preserved under the isometry $\vec q\mapsto 2\vec p-\vec q$ for any fixed $p\in S$, e.g. any lattice is antipodal.

\begin{lem}[upper bounds for a stable radius $\al(S)$, {\cite[Lemma 3.6]{anosova2025recognition}}]
\label{lem:upper_bounds}
\textbf{(a)}
Let $S\subset\R^n$ be a periodic point set with a unit cell $U$, which has the longest edge $b$ and longest diagonal $d$.
Set $r(U)=\max\{b,\frac{d}{2}\}$.
Then the bridge length $\be(S)$ from Definition~\ref{dfn:bridge_length} 
has the upper bound $\min\{2R(S),r(U)\}\geq\be(S)$.
\medskip

\noindent
\textbf{(b)}
For any antipodal periodic set $S\subset\R^n$ whose covering radius is $R(S)$, the minimum stable radius has the upper bound $2R(S)+\be(S)>\al(S)$.  
\medskip

\noindent
\textbf{(c)}
Let $S\subset\R^n$ be any periodic point set with the bridge length $\be$.
For any point $p\in S$ and a radius $\al_0\geq 2R(S)$, the order $|\sym(S,p;\al_0)|$ of the group $\sym(S,p;\al_0)$ should be finite.
Let $p_1,\dots,p_m\in S$ be all points of an asymmetric unit of $S$.
Set $L=\left[\sum\limits_{i=1}^{m}\big(\log_2|\sym(S,p_i;\al_0)|-\log_2|\sym(S,p_i)|\big)\right]$. 
Then the minimum stable radius $\al(S)$ from Definition~\ref{dfn:stable_radius} 
has the upper bound $\al_0+(L+m)\be\geq\al(S)$.
If $\al_0=2R(S)$, then $(L+m+1)2R(S)\geq\al(S)$. 
\elem
\end{lem}

The upper bound in Lemma~\ref{lem:upper_bounds}(a) holds for any unit cell of $S$.
If a cell is non-reduced and too long, its reduced form can have smaller bounds for $\be(S)$.   

\begin{exa}[upper bounds for $\al(S)$ and $\be(S)$]
\label{exa:upper_bounds}
Let $\La(b)\subset\R^n$ be a lattice whose unit cell is a rectangular box with the longest edge $b\geq 1$.
\smallskip

\noindent
\textbf{(a)}
In Lemma~\ref{lem:upper_bounds}(a), the upper bound $b\geq\be(S)$ is tight because $\be(\La(b))=b$.
\smallskip

\noindent
\textbf{(b)}
In Lemma~\ref{lem:upper_bounds}(b), the ratio $(2R(S)+\be(S))/\al(S)\geq 1$ tends to $1$ as $b\to+\infty$ for any fixed $n$.
Indeed, a cluster $C(\La(b),0;\al)$ is $n$-dimensional only for $\al\geq b$, so the group $\sym(\La(b),0;\al)$ stabilises at $\al=b$, hence $\al(S)=b+\be(\La(b))=2b$ is the minimum stable radius. 
The covering radius $R(\La(b))$ is half of the longest diagonal of the rectangular cell $U$.
If $b\to+\infty$ and all other sizes of $U$ remain fixed, the ratio $(2R(\La(b))+\be(\La(b)))/\al(S)$ tends to $1$ for any fixed $n$.
\medskip

\noindent
\textbf{(c)}
Lemma~\ref{lem:upper_bounds}(c) was motivated by \cite[Theorem~1.3]{dolbilin1998multiregular}, which implies the upper bound $\be(S)+2m(n^2+1)\log_2(2+R(S)/r(S))>\al(S)$ for $m$-regular point sets.
Let $\La\subset\R^2$ be a lattice whose unit cell is a rhombus with sides 1.
Then $m=1$, $n=2$, $r(\La)=0.5$, $\be(\La)=1$, and $\al(\La)=2$.
If $\La$ deforms from a square lattice to a hexagonal lattice, the covering radius $R(\La)$ varies in the range $[\frac{1}{\sqrt{3}},\frac{1}{\sqrt{2}}]$.
The past bound above gives the estimate $1+2(2^2+1)\log_2(2+\frac{2}{\sqrt{3}})\approx 17.6>\al(\La)=2$.  
For any lattice $\La$ in this family, the symmetry group $\sym(\La,0)=\sym(\La,0;1)$ stabilises at $\al_0=1$. 
\myskip

Lemma~\ref{lem:upper_bounds}(c) for $\al_0=1$ gives $L=\log_2(2)-\log_2(2)=0$, so the upper bound $\al_0+(L+m)\be(S)\geq\al(S)$ is tight: $2\geq\al(\La)$. 
In practice, if $L$ is large because some local clusters $C(S;p;\al_0)$ have too many symmetries, one can increase the radius $\al_0$ to reduce $L$ for a better bound of $\al(S)$. 
\eexa
\end{exa}

Definition~\ref{dfn:isoset} reminds of the \emph{isoset}, which was initially introduced in \cite[Definition~9]{anosova2021isometry}. 
We also cover the case of rigid motion 
and prove Completeness Theorem~\ref{thm:isoset_complete} in the appendix in more detail than in \cite[Theorem~9]{anosova2021isometry}.
  
\index{isoset}
  
\begin{dfn}[isoset $I(S;\al)$ at a radius $\al\geq 0$]
\label{dfn:isoset}
Let a periodic point set $S\subset\R^n$ have a motif $M$ of $m$ points.
Split all points $p\in M$ into $\al$-equivalence classes.
Each $\al$-equivalence class of (say) $k$ points in $M$ can be associated with the \emph{isometry class} $\si=[C(S,p;\al)]$ of an $\al$-cluster centreed at some $p\in M$.
The \emph{weight} of $\si$ is $w=k/m$.
Then the \emph{isoset} $I(S;\al)$ is the unordered set of all isometry classes $(\si;w)$ with weights $w$ for all points $p$ in the motif $M$.
If we replace isometry with rigid motion, we get the \emph{oriented} isoset $I^o(S;\al)$.
\edfn
\end{dfn}

All points $p$ of a lattice $\La\subset\R^n$ from one $\al$-equivalence class for any radius $\al\geq 0$ because all $\al$-clusters $C(\La,p;\al)$ are isometrically equivalent to each other by translations. 
Hence the isoset $I(\La;\al)$ is one isometry class of weight 1 for $\al\geq 0$, see examples in Fig.~\ref{fig:square_vs_hexagon}.
All isometry classes $\si$ in $I(S;\al)$ are in a 1-1 correspondence with all $\al$-equivalence classes in the $\al$-partition $P(S;\al)$ from Definition~\ref{dfn:isotree}.
\myskip

Hence, the isoset $I(S;\al)$ without weights can be viewed as a set of points in the isotree $\IT(S)$ at the radius $\al$.
The size of the isoset $I(S;\al)$ equals the number $|P(S;\al)|$ of $\al$-equivalence classes in the $\al$-partition. 
Formally, $I(S;\al)$ depends on $\al$ because $\al$-clusters grow in $\al$.
To distinguish any $S,Q\subset\R^n$ under isometry, we will compare their isosets at a maximum stable radius of $S,Q$.

\begin{exa}[isosets of simple lattices]
\label{exa:isosets_lattices}
\textbf{(a)}
Any lattice $\La\subset\R^n$ is 1-regular by Definition~\ref{dfn:m-regular} and can be assumed to contain the origin $0$ of $\R^n$.
Then the isoset $I(\La;\al)$  consists of a single isometry class of a cluster $C(\La,0;\al)$.
So the isotree $\IT(\La)$ is a linear path, which is horizontally drawn for the hexagonal and square lattices $\La_6,\La_4$ in Fig.~\ref{fig:lattice_isotree}.
If both $\La_6,\La_4$ have a minimum inter-point distance 1, then the bridge length from Definition~\ref{dfn:bridge_length} is $\be=1$.
\medskip

\begin{figure}[h!]
\includegraphics[width=\textwidth]{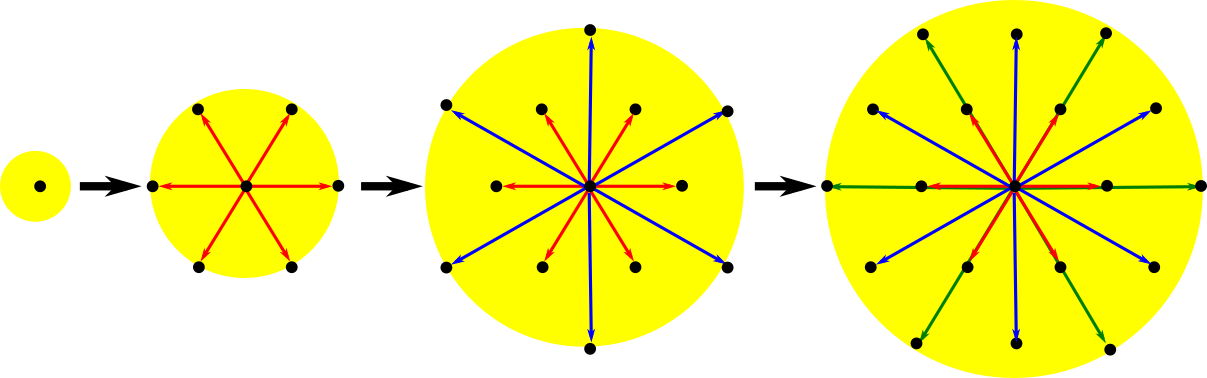}
\sskip

\includegraphics[width=\textwidth]{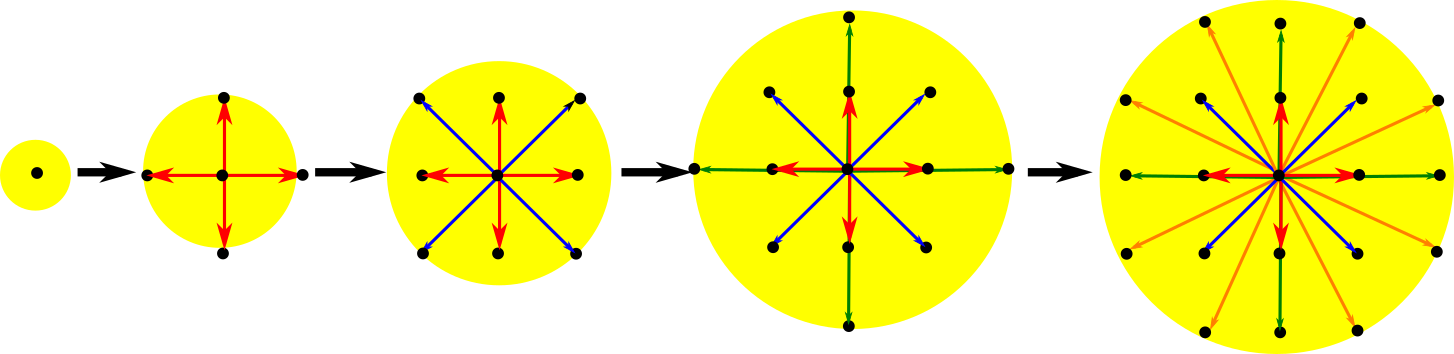}
\caption{
The isotree of any lattice $\La$ is $[0,+\infty)$ is a line $\R$ parametrised by the radius $\al$.
\textbf{Top}: the isotree of the hexagonal lattice $\La_6$.
\textbf{Bottom}: the isotree of the square lattice $\La_4$.
}
\label{fig:lattice_isotree}
\end{figure}

\noindent
\textbf{(b)}
For the hexagonal lattice $\La_6\subset\R^2$, 
$C(\La_6,(0,0);\al)$ includes points $p\neq(0,0)$ only for $\al\geq 1$.
The cluster $C(\La_6,(0,0);1)=\{(0,0),(\pm 1,0),(\pm\frac{1}{2},\pm\frac{\sqrt{3}}{2})\}$ appears in the 2nd step of  Fig.~\ref{fig:lattice_isotree}~(left).
The symmetry group $\sym(\La_6,(0,0);\al)$ becomes the dihedral group $D_6$ (all symmetries of a regular hexagon) for $\al\geq 1$.
Hence any $\al\geq\be+1=2$ is stable.
The isoset $I(\La_6;1)$ is the isometry class of the cluster $C(\La_6,(0,0);1)$ of six vertices of the regular hexagon and its centre.
\medskip

\noindent
\textbf{(c)}
For the square lattice $\La_4\subset\R^2$, $C(\La_4,(0,0);\al)$ has points $p\neq(0,0)$ only for $\al\geq 1$.
 $C(\La_4,(0,0);2)=\{(0,0),(\pm 1,0),(0,\pm 1),(\pm\sqrt{2},\pm\sqrt{2}),(\pm 2,0),(0,\pm 2)\}$ includes the origin $(0,0)$ with its 12 neighbors in the 4th step of Fig.~\ref{fig:lattice_isotree}~(right).
The group $\sym(\La_4,(0,0);\al)$ becomes the dihedral group $D_4$ (all symmetries of a square) for $\al\geq 1$.
So any $\al\geq\be+1=2$ is stable.
The isoset $I(\La_4;1)$ is the isometry class of $C(\La_4,(0,0);1)$ of four vertices of the square and its centre.
\eexa
\end{exa}

An equality $\si=\xi$ between isometry classes of clusters means that some (hence any) clusters $C(S,p;\al)$ and $C(Q,q;\al)$ representing $\si,\xi$, respectively, are related by $f\in\Or(\R^n;p,q)$, which will be algorithmically tested in Corollary~\ref{cor:compare_isosets}.

\begin{thm}[isometry classification of periodic point sets, {\cite[Theorem 3.10]{anosova2025recognition}}]
\label{thm:isoset_complete}
For any periodic point sets $S,Q\subset\R^n$, let $\al$ be a common stable radius satisfying Definition~\ref{dfn:stable_radius} for an upper bound $\be\geq\be(S),\be(Q)$.
Then $S,Q$ are isometric (related by rigid motion, respectively) if and only if there is a bijection $\ph:I(S;\al)\to I(Q;\al)$ (between oriented isosets, respectively) that preserves all their weights.
\ethm
\end{thm}

Theorem~\ref{thm:isoset_complete} was inspired by the seminal result in \cite[Theorem~1.3]{dolbilin1998multiregular} saying that, for a multi-regular point set $X$, ``the only Delone sets $Y$ all of whose $\rho$-stars are isometric to $\rho$-stars of $X$ are sets globally isometric to $X$''.
After renaming $\rho$-stars as $\al$-clusters, we collected their isometry classes (with weights) into the \emph{isoset} to rephrase \cite[Theorem~1.3]{dolbilin1998multiregular} as a classification of all periodic point sets by isosets. 
\myskip

The $\al$-equivalence and isoset in Definition~\ref{dfn:isoset} can be refined by labels such as chemical elements, which keeps Theorem~\ref{thm:isoset_complete} valid for labelled points.
\smallskip

When comparing sets from a finite database, it suffices to build their isosets only up to a common upper bound of a stable radius $\al$ in Lemma~\ref{lem:upper_bounds}(c).

\section{Continuous metrics on isometry classes of periodic sets in $\R^n$}
\label{sec:isosets_metrics}

This section proves the continuity of the isoset $I(S;\al)$ in Theorem~\ref{thm:isoset_continuous} by using the Earth Mover's Distance (EMD) from Definition~\ref{dfn:EMD_isosets}.
\myskip

For a point $p\in\R^n$ and a radius $\ep$, 
the {\em closed} ball $\bar B(p;\ep)=\{q\in\R^n \vl |\vec q - \vec p|\leq\ep\}$ has as its the boundary $(n-1)$-dimensional sphere $\bd\bar B(p;\ep)\subset\R^{n}$.
The \emph{$\ep$-offset} of any set $C\subset\R^n$ is the Minkowski sum $C+\bar B(0;\ep)=\{\vec p+\vec q \vl p\in C, q\in \bar B(0;\ep)\}$.
\medskip

Then the directed \emph{Hausdorff} distance from Example~\ref{exa:metrics}(b) $d_H(C,D)$ 
is the minimum radius $\ep\geq 0$ such that $C\subseteq D+\bar B(0;\ep)$.
Definition~\ref{dfn:tolerant_metric} introduces the crucial new metric, which will be explicitly computed in Lemma~\ref{lem:max-min_formula}. 

\begin{dfn}[boundary tolerant metric $\BT$ on isometry classes of clusters] 
\label{dfn:tolerant_metric}
For a radius $\al$ and periodic point sets $S,Q\subset\R^n$, let clusters $C(S,p;\al),C(Q,q;\al)$ represent isometry classes $\si\in I(S;\al),\xi\in I(Q;\al)$, respectively.
The \emph{boundary tolerant} metric $\BT(\si,\xi)$ 
 is defined as the minimum $\ep\geq 0$ such that 
\smallskip

\noindent
(\ref{dfn:tolerant_metric}a) 
$C(Q,q;\al-\ep)\subseteq f(C(S,p;\al))+\bar B(0;\ep)$ for some $f\in\Or(\R^n;p,q)$, and 
\smallskip

\noindent
(\ref{dfn:tolerant_metric}b) 
$C(S,p;\al-\ep)\subseteq g(C(Q,q;\al))+\bar B(0;\ep)$ for some $g\in\Or(\R^n;q,p)$.
\edfn
\end{dfn}

In Definition~{\ref{dfn:tolerant_metric}}, if one cluster consists of only its centre, e.g. $C(S,p;\al)=\{p\}$, then the boundary tolerant metric is $\BT=\max\{|\vec s-\vec q| \mid s\in C(Q,q;\al)\}$.
\cite[Lemma~4.2]{anosova2025recognition} proves that $\BT$ is independent of cluster representatives and satisfies all metric axioms from Definition~\ref{dfn:metrics}(a). 

\begin{exa}[square lattice vs hexagonal]
\label{exa:square_vs_hexagon}
The isoset $I(\La;\al)$ of any lattice $\La\subset\R^n$ containing the origin $0$ consists of a single isometry class $[C(\La,0;\al)]$, see Example~\ref{exa:isosets_lattices}.
For the square (hexagonal) lattice with minimum inter-point distance 1 in Fig.~\ref{fig:square_vs_hexagon}, the cluster $C(\La,0;\al)$ consists of only 0 for $\al<1$ and includes four (six) nearest neighbors of 0 for $\al\geq 1$.
Hence $\sym(\La,0;\al)$ stabilises as the symmetry group of the square (regular hexagon) for $\al\geq 1$.
The lattices have the minimum stable radius $\al(\La)=2$ and $\be(\La)=1$ by Example~\ref{exa:upper_bounds}(c).

\begin{figure}[h!]
\includegraphics[width=\linewidth]{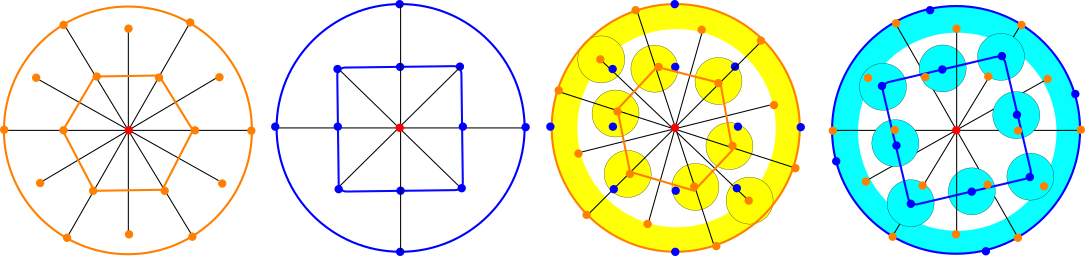}
\caption{
Example~\ref{exa:square_vs_hexagon} computes the metric $\BT$ from Definition~\ref{dfn:tolerant_metric} for the isometry classes of the $2$-clusters
in the square and hexagonal lattices $\La_4,\La_6$.
\textbf{1st}: the 2-cluster $C(\La_6,0;2)$ with its boundary circle $\bd\bar B(0;2)$;
\textbf{2nd}: the 2-cluster $C(\La_4,0;2)$ with its boundary circle $\bd\bar B(0;2)$;
\textbf{3rd}: for $\ep=\sqrt{2}-1\approx 0.41$, the cluster $C(\La_4,0;2)$ is covered by the yellow $\ep$-offset of $C(\La_6,0;2)\cup\bd\bar B(0;2)$ rotated through $15^\circ$ clockwise.
\textbf{4th}: $C(\La_6,0;2)$ is covered by the blue $\ep$-offset of $C(\La_4,0;2)\cup\bd\bar B(0;2)$ rotated through $15^{\circ}$ anticlockwise, so $\BT=\sqrt{2}-1$.}
\label{fig:square_vs_hexagon}
\end{figure}

Fig.~\ref{fig:square_vs_hexagon} shows the stable 2-clusters $C(\La_4,0;2)$ and $C(\La_6,0;2)$ of the square ($\La_4$) and hexagonal ($\La_6$) lattices.
Without rotations, the 1st picture of Fig.~\ref{fig:square_vs_hexagon} shows the directed Hausdorff distance $d_H=\sqrt{(1-\frac{\sqrt{3}}{2})^2+(\frac{1}{2})^2}=\sqrt{2-\sqrt{3}}\approx 0.52$ between clusters with the added boundary circle $\bd B(0;2)$.
Due to high symmetry, it suffices to consider rotations 
of the square vertex $(1,1)$ 
for angles $\ga\in[45^{\circ},60^{\circ}]$ because all other ranges can be isometrically mapped to this range for another vertex of the square.
\myskip

We find the squared distances $s_1(\ga)$ and $s_2(\ga)$ from the vertex $(\sqrt{2}\cos\ga,\sqrt{2}\sin\ga)$ rotated from $(1,1)$ at $\ga=45^\circ$ through the angle $\ga-45^{\circ}$ to its closest neighbors $(\frac{1}{2},\frac{\sqrt{3}}{2})$ and $(\frac{3}{2},\frac{\sqrt{3}}{2})$ in $C(\La_6,0;2)$.
\medskip

$$\begin{array}{l}
s_1(\gamma)=\left|(\sqrt{2}\cos\ga,\sqrt{2}\sin\ga)-\Big(\frac{1}{2},\frac{\sqrt{3}}{2}\Big)\right|^2= \\
\left(\sqrt{2}\cos\ga-\frac{1}{2}\right)^2+\left(\sqrt{2}\sin\ga-\frac{\sqrt{3}}{2}\right)^2=\\
3-\sqrt{2}\cos\ga-\sqrt{6}\sin\ga,
\end{array}$$ 

$$\frac{ds_1}{d\ga}=\sqrt{2}\sin\ga-\sqrt{6}\cos\ga=0,\;
\tan\ga=\sqrt{3}, \; 
\ga=60^{\circ}, 
s_1=(\sqrt{2}-1)^2$$ 
is minimal for the points in 
the line $y=\sqrt{3}x$ at distances $1,\sqrt{2}$ from $0$.

$$\begin{array}{l}
s_2(\gamma)=\left|(\sqrt{2}\cos\ga,\sqrt{2}\sin\ga)-\Big(\frac{3}{2},\frac{\sqrt{3}}{2}\Big)\right|^2=\\
\left(\sqrt{2}\cos\ga-\frac{3}{2}\right)^2+\left(\sqrt{2}\sin\ga-\frac{\sqrt{3}}{2}\right)^2=\\
5-3\sqrt{2}\cos\ga-\sqrt{6}\sin\ga,
\end{array}$$

$$\frac{ds_2}{d\ga}=3\sqrt{2}\sin\ga-\sqrt{6}\cos\ga=0,
\ga=30^{\circ}, 
s_2=(\sqrt{3}-\sqrt{2})^2$$ 
is minimal for the points in the line $y=\frac{x}{\sqrt{3}}$ at distances $\sqrt{2},\sqrt{3}$ from $0$. 
\medskip

It might look that the second minimum is smaller.
However, for the angle $\ga=30^{\circ}$, another vertex $(-1,1)$ rotated through $\ga-45^\circ=-15^{\circ}$ has distance $\sqrt{2}-1$ to its closest neighbor $(-\frac{1}{2},\frac{\sqrt{3}}{2})\in C(\La_6,0;2)$.
For any angle $\ga\in[45^\circ,60^\circ]$, the second function has the minimum $s_2(45^{\circ})=2-\sqrt{3}=d_H^2$ in the 1st picture of Fig.~\ref{fig:square_vs_hexagon}.
\myskip

Hence, the vertex $(1,1)$ has the minimum distance $\sqrt{2}-1\approx 0.41<\sqrt{2-\sqrt{3}}\approx 0.52$ in the 3rd picture of Fig.~\ref{fig:square_vs_hexagon}.
All other points of the square cluster $C(\La_4,0;2)$ are even closer to their neighbors in $C(\La_6,0;2)$.
For example, the point $(1,0)$ rotated by $15^\circ$ has the distance to $(1,0)$ equal to $\sqrt{(\cos 15^\circ-1)^2+\sin^2 15^\circ}\approx 0.26$.
\myskip

The final picture in Fig.~\ref{fig:square_vs_hexagon} confirms that all points of the hexagonal cluster $C(\La_6,0;2)$ are covered by the $(\sqrt{2}-1)$-offset of 
$C(\La_4,0;2)$ and the boundary circle.
So 
$\BT=\sqrt{2}-1\approx 0.41$.
\eexa
\end{exa}

Non-isometric periodic sets $S,Q$ can have isosets of different numbers of isometry classes.
The distance between these weighted distributions of different sizes can be measured by the Earth Mover's Distance  below.

\index{Earth Mover's Distance}

\begin{dfn}[Earth Mover's Distance on isosets]
\label{dfn:EMD_isosets}
Let periodic point sets $S,Q\subset\R^n$ have a common stable radius $\al$ and isosets $I(S;\al)=\{(\si_i,w_i)\}$ and $I(Q;\al)=\{(\xi_j,v_j)\}$, where $i=1,\dots,m(S)$ and $j=1,\dots,m(Q)$.
The \emph{Earth Mover's Distance} is 
$$\EMD( I(S;\al), I(Q;\al) )=\sum\limits_{i=1}^{m(S)} \sum\limits_{j=1}^{m(Q)} f_{ij} \BT(\si_i,\xi_j)$$ minimised over \emph{flows} 
$f_{ij}\in[0,1]$ subject to 
the conditions 
$\sum\limits_{j=1}^{m(Q)} f_{ij}\leq w_i$ for $i=1,\dots,m(S)$,
$\sum\limits_{i=1}^{m(S)} f_{ij}\leq v_j$ for $j=1,\dots,m(Q)$, and
$\sum\limits_{i=1}^{m(S)}\sum\limits_{j=1}^{m(Q)} f_{ij}=1$.
\edfn
\end{dfn}

\begin{exa}[EMD for lattices with $\BD=+\infty$]
\label{exa:EMD_cluster}
\cite[Example~2.1]{widdowson2022resolving} 
showed that the lattices $S=\Z$ and $Q=(1+\de)\Z$ have the bottleneck distance $\BD(S,Q)=+\infty$ for any $\de>0$.
We show that $S,Q$ have Earth Mover's Distance $\EMD=2\de$ at their common stable radius $\al=2+2\de$. 
The bridge lengths are $\be(S)=1$ and $\be(Q)=1+\de$.
The $\al$-cluster $C(S,0;\al)$ contains non-zero points for $\al\geq 1$, e.g. $C(S,0;1)=\{0,\pm 1\}$.
\myskip

The symmetry group $\sym(S,0;\al)=\Z_2$ includes a non-trivial reflection with respect to 0 for all $\al\geq 1$, so the stable radius of $S$ is any $\al\geq \be+1=2$.
Similarly, $Q$ has $\be(Q)=1+\de$ and stable radii $\al\geq 2(1+\de)$.
The Earth Mover's Distance between $I(S;\al)$ and $I(Q;\al)$ at the common stable radius $\al=2+2\de$ equals the metric $\BT$ between the only $\al$-clusters $C(S,0;\al)=\{0,\pm 1,\pm 2\}$ and $C(Q,0;\al)=\{0,\pm (1+\de),\pm2(1+\de)\}$. 
\medskip

By Definition~\ref{dfn:tolerant_metric} we look for a minimum $\ep>0$ such that the cluster $C(S,0;\al-\ep)$ is covered by $\ep$-offsets of $\pm (1+\de),\pm2(1+\de)$ and vice versa.
If $\ep<2\de$, the points $\pm 2\in C(S,0;\al-\ep)$ cannot be $\ep$-close to $\pm (1+\de),\pm1(+\de)$, but $\ep=2\de$ is large enough.
The cluster $C(Q,0;\al-2\de)=\{0,\pm (1+\de)\}$ is covered by the $2\de$-offset of $C(S,0;\al)=\{0,\pm 1,\pm 2\}$, so $\EMD(I(S;\al),I(Q;\al))
=2\de$.
\eexa 
\end{exa}

For rigid motion instead of general isometry, Definition~\ref{dfn:tolerant_metric} of a boundary tolerant metric $\BT$ is updated to $\BT^o$ by considering only orientation-preserving isometries from $\SO(\R^n;p,q)$, which also makes the continuity below valid for oriented isosets $I^o(S;\al)$ under $\EMD$ using $\BT^o$ instead of $\BT$ in Definition~\ref{dfn:EMD_isosets}. 

\begin{thm}[continuity of isosets under perturbations, {\cite[Theorem 4.9]{anosova2025recognition}}]
\label{thm:isoset_continuous}
Let periodic point sets $S,Q\subset\R^n$ have a bottleneck distance $\BD(S,Q)<r(Q)$, where $r(Q)$ is the packing radius in Definition~\ref{dfn:packing+covering}(a). 
Then the isosets $I(S;\al),I(Q;\al)$ are close in the Earth Mover's Distance: 
$\EMD( I(S;\al), I(Q;\al) )\leq 2\BD(S,Q)$ for $\al\geq 0$. 
\ethm
\end{thm}

Corollary~\ref{cor:isosets_metrics}(a) justifies that the EMD satisfies all metric axioms for periodic point sets that have a stable radius $\al$.
Corollary~\ref{cor:isosets_metrics}(b) avoids this dependence on $\al$ and scales any periodic point set $S$ to the minimum stable radius $\al(S)=1$.

\begin{cor}[metric on periodic point sets, {\cite[Corollary 4.10]{anosova2025recognition}}]
\label{cor:isosets_metrics}
\textbf{(a)}
For $\al>0$, $\EMD( I(S;\al), I(Q;\al) )$ is a metric on the space of isometry classes of all periodic point sets with a stable radius $\al$ in $\R^n$.
\medskip

\noindent
\textbf{(b)}
For a periodic point set $S\subset\R^n$, let $S/r(S)\subset\R^n$ denote $S$ after uniformly dividing all vectors by the packing radius $r(S)$.
Then $$|r(S)-r(Q)|+\EMD\big( I(S/r(S);1), I(Q/r(Q);1) \big)$$ is a metric on all periodic point sets. 
\ecor
\end{cor}

\section{Algorithms to compute isosets and their approximate metrics}
\label{sec:isosets_algorithms}

This section describes time complexities for computing the complete invariant isoset (Theorem~\ref{thm:compute_isoset}), comparing isosets (Corollary~\ref{cor:compare_isosets}), approximating 
the boundary tolerant metric $\BT$ and 
Earth Mover's Distance on isosets (Corollary~\ref{cor:approximate_EMD}).  
\myskip

All time estimates will use the geometric complexity $\GC(S)$ defined below.

\index{geometric complexity}

\begin{dfn}[geometric complexity $\GC$]
\label{dfn:geom_complexity}
Let a periodic point set $S\subset\R^n$ have an asymmetric unit of $m$ points in a  cell $U$ of volume $\vol[U]$.
Let $L$ be the symmetry characteristic for $\al_0=2R(S)$ in Lemma~\ref{lem:upper_bounds}(c), where $R(S)$ is the covering radius.
The \emph{geometric complexity} is 
$\GC(S)=\dfrac{(10(L+m+2)R(S)/n)^n}{2\vol[U]}$.
\edfn
\end{dfn}

Let $V_n=\dfrac{\pi^{n/2}}{\Ga(\frac{n}{2}+1)}$ be the volume of the unit ball in $\R^n$, where the Gamma function $\Ga$ is defined as $\Ga(k)=(k-1)!$ and $\Ga(\frac{k}{2}+1)=\sqrt{\pi}(k-\frac{1}{2})(k-\frac{3}{2})\cdots\frac{1}{2}$ for any integer $k\geq 1$.
Set $\nu(U,\al,n)=\dfrac{(\al+d)^n V_n}{\vol[U]}$, where 
$d=\sup\limits_{p,q\in U}|\vec p-\vec q|$ is a longest diagonal of $U$.
\medskip

The main input size of a periodic set is the number $m$ of motif points because the length of a standard Crystallographic Information File (CIF) is linear in $m$.
\myskip

For a fixed dimension $n$, the big $O$ notation $O(m^n)$ in all complexities means a function $t(m)$ such that $t(m)\leq C m^n$ for a fixed constant $C$ independent of $m$.
We will include all other parameters depending on a periodic point set $S$.

\index{isoset}

\begin{thm}[time of an isoset, {\cite[Theorem 5.3]{anosova2025recognition}}]
\label{thm:compute_isoset}
For any periodic point set $S\subset\R^n$ given by a motif $M$ of $m$ points in a unit cell $U$, the isoset $I(S;\al)$ at a stable radius $\al$ can be found in time $O(m^2 k^{\lceil n/3\rceil}\log k)$, where $k=\nu m$ 
for $\nu\leq \GC(S)$.
\ethm
\end{thm}

\index{isoset}
\index{directed distance}

\begin{cor}[comparing isosets, {\cite[Corollary 5.4]{anosova2025recognition}}]
\label{cor:compare_isosets}
There is an algorithm to check if any periodic point sets $S,Q\subset\R^n$ with motifs of at most $m$ points are isometric in total time $O(m^2 k^{\lceil n/3\rceil}\log k)$, where $k=\nu m$ for $\nu\leq\max\{\GC(S),\GC(Q)\}$.
\ecor
\end{cor}

\begin{dfn}[directed distances $d_R$ and $d_M$]
\label{dfn:directed_distances}
\textbf{(a)}
For any sets $C,D\subset\R^n$, the directed \emph{rotationally invariant} distance  $d_R(C,D)=\min\limits_{f\in\Or(\R^n)}d_H(C,f(D))$ is minimised over all maps $f\in\Or(\R^n;0)$, which fix the origin $0\in\R^n$. 
\smallskip

\noindent
\textbf{(b)}
For any finite sets $C,D\subset\R^n$, order all points $p_1\dots,p_k\in C$ by increasing distance to the origin $0$.
The \emph{radius} of $C$ is $R(C)= \max\limits_{p\in C}|p|$.
Define the directed \emph{max-min distance} as $d_M(C,D)=\max\limits_{i=1,\dots,k}\min\{\;\al-|p_i|,\; d_R(\{p_1,\dots,p_i\}, D) \;\}$.
\edfn
\end{dfn}

If $C'\subset C$, then $d_R(C',D)\leq d_R(C,D)$.
Let $C,D\subset\bar B (0; \al)$ be finite sets including the origin $0$.
If $C=\{0\}$, then $d_R(C,D)=0$  because $C\subset D$, but $d_R(D,C)=R(D)$ is the radius of $D$ because $D\subset \{0\}+\bar B(0;\ep)$ only for $\ep\geq R(D)$.
\myskip

Definition~\ref{dfn:directed_distances}, Lemma~\ref{lem:max-min_formula} and hence all further results work for rigid motion by restricting all maps to the special orthogonal group $\SO(\R^n;0)$.

\begin{lem}[max-min formula for $d_R$, {\cite[Lemma 5.6]{anosova2025recognition}}]
\label{lem:max-min_formula}
For any finite sets $C,D\subset\R^n$, if $\al\geq R(C)$, then $d_R(C\cup\bd\bar B(0;\al),D\cup\bd\bar B(0;\al))$ equals $d_M(C,D)$.
\elem
\end{lem}

\begin{exa}[max-min formula]
\label{exa:max-min_formula}
Consider the subcluster $C\subset C(\La_4,0;2)$ of the points $p_1=(1,0)$, $p_2=(1,1)$, $p_3=(1,-1)$, $p_4=(2,0)$ from the square lattice $\La_4$ in Fig.~{\ref{fig:square_vs_hexagon}}.
Let $\al=2$ and $D=C(\La_6,0;2)$ be the 2-cluster of the hexagonal lattice $\La_6$.
Then $d_R(p_1,D)=0$ because $p_1$ coincides with $(1,0)\in D$.
Then $d_R(\{p_1,p_2\},D)=\sqrt{2}-1$, because the cloud $D$ after the clockwise rotation through $15^\circ$ has the points $(\cos 15^\circ,-\sin 15^\circ)$ and $(\frac{1}{\sqrt{2}},\frac{1}{\sqrt{2}})$ at distances $\sqrt{(\cos 15^\circ -1)^2+\sin^2 15^\circ}\approx 0.26$, $\sqrt{2}-1\approx 0.41$ to $p_1,p_2$, respectively.
Then $d_R(\{p_1,p_2,p_3\},D)=\sqrt{2}-1$ because the same rotated image of $D$ has $(\sqrt{\frac{3}{2}}, -\sqrt{\frac{3}{2}})$ at the distance $\sqrt{3}-\sqrt{2}\approx 0.32$ to $p_3$. 
\myskip

For $i=1$, $\min\{\al-|p_1|,d_R(p_1,D)\}=\min\{2-1,0\}=0$.
For $i=2,3$, 
$$\begin{array}{l}
\min\{\al-|p_2|,d_R(\{p_1,p_2\},D)\}=\\
\min\{\al-|p_3|,d_R(\{p_1,p_2,p_3\},D)\}=\\
\min\{2-\sqrt{2},\sqrt{2}-1\}=\sqrt{2}-1.
\end{array}$$

For $i=4$, $\min\{\al-|p_4|,d_R(C,D)\}=0$ since $\al=2=|p_4|$.
\myskip

The maximum value is $\sqrt{2}-1$, so Example~\ref{exa:square_vs_hexagon} fits Lemma~\ref{lem:max-min_formula}.
\eexa
\end{exa}

Lemma~\ref{lem:rot-inv_distance} extends \cite[section 2.3]{goodrich1999approximate} from $n=3$ to any dimension $n>1$. 

\begin{lem}[approximating $d_R$, {\cite[Lemma 5.8]{anosova2025recognition}}]
\label{lem:rot-inv_distance}
Let a cloud $C\subset\R^n$ consist of $k=|C|$ points ordered by distances $|p_1|\leq\dots\leq|p_k|$ from the origin 
and $\ar{C}$ denote the number of different vectors $\vec p/|\vec p|$ for $p\in C$. 
For each $j=1,\dots,k$, consider the subcloud $C_j=\{p_1,\dots,p_{j}\}$.
For any cloud $D\subset\R^n$ of $|D|$ points, all distances $d_j=d_R(C_j,D)$ from Definition~{\ref{dfn:directed_distances}} for $j=1,\dots,k$ can be approximated by some $d'_j$ in time $O(|C|\ar{C}^{n-1}|D|)$ so that $d_j\leq d'_j\leq \om d_j$, $\omega=1+\frac{1}{2}n(n-1)$.
\elem
\end{lem}

The proof of Lemma~\ref{lem:rot-inv_distance} uses only orientation-preserving isometries from $\SO(\R^n,0)$.
Hence the upper bounds from Lemma~\ref{lem:rot-inv_distance}, Theorem~\ref{thm:approximate_BT}, and Corollary~\ref{cor:approximate_EMD} work for both cases of rigid motion and general isometry in $\R^n$.

\begin{thm}[approximating $\BT$, {\cite[Theorem 5.9]{anosova2025recognition}}]
\label{thm:approximate_BT}
Let periodic point sets $S,Q\subset\R^n$ have isometry classes $\si,\xi$ represented by clusters $C,D$ of a radius $\al$, respectively.
In the notations of Lemma~{\ref{lem:rot-inv_distance}}, 
$\BT(\si,\xi)$ from Definition~{\ref{dfn:tolerant_metric}} can be approximated with the factor $\omega=1+\frac{1}{2}n(n-1)$ in time $O(|C|(\ar{C}^{n-1}+\ar{D}^{n-1})|D|)$.  
\ethm
\end{thm}

\begin{cor}[approximating EMD on isosets, {\cite[Corollary 5.10]{anosova2025recognition}}]
\label{cor:approximate_EMD}
Let $S,Q\subset\R^n$ be periodic point sets whose motifs have at most $m$ points $p$ and $\chi$ different vectors $\vec p/|\vec p|$.
For any $\al>0$, the metric $\EMD(I(S;\al),I(Q;\al))$ can be approximated with the factor $\omega=1+\frac{1}{2}n(n-1)$ in time $O(\nu^2 m^4 \chi^{n-1})$, where $\nu\leq\max\{\GC(S),\GC(Q)\}$.
\ecor
\end{cor}

Counting directions $\vec p/|p|$ as points ($\chi\leq m$), for dimension $n=3$, the rough bounds for the isoset and its approximate $\EMD'$ in Theorem~\ref{thm:compute_isoset} and Corollary~\ref{cor:approximate_EMD} are $O(m^3\log m)$ and $O(m^6)$, respectively.
Algorithms 1-2 in the appendix describe pseudocodes for
Lemma~\ref{lem:rot-inv_distance}, Theorem~\ref{thm:approximate_BT}, and Corollary~\ref{cor:approximate_EMD}.

\section{Comparisons of isosets with past invariants and experiments}
\label{sec:isosets_experiments}

This section justifies that the isoset can be efficiently used with the faster $\PDD$ due to a lower bound in Theorem~\ref{thm:isosets_lower_bound}.

\begin{thm}[lower bound for $\EMD$, {\cite[Theorem 6.5]{anosova2025recognition}}]
\label{thm:isosets_lower_bound}
Let $S,Q\subset\R^n$ be periodic sets with a common stable radius $\al$.
Let $\ep=\EMD(I(S;\al),I(Q;\al))$ and $k$ be the maximum number of points of $S,Q$ in their $(\al-\ep)$-clusters.
If $\ep$ is less than the half-distance between any points of $S,Q$, then $\EMD(\PDD(S;k),\PDD(Q;k))\leq\ep$.
\ethm
\end{thm}

Hence, $\PDD(S;k)$ can be used for a fast filtering of distant crystals so that the isoset is computed only for near-duplicates that are hard to distinguish.
\myskip

In 1930, future Nobel laureate Linus Pauling noticed the ambiguity of crystal structures obtained by diffraction \cite{pauling1930crystal}.
Such \emph{homometric} crystals with identical diffraction patterns were only manually distinguished until now because even the generically complete $\PDD$s coincide for the Pauling periodic sets $P(\pm u)$ for all $u\in(0,0.25)$, see the real overlaid crystals for $u=0.03$ in Fig.~\ref{fig:Pauling}~(left).
\myskip

\newcommand{\mh}{45mm}
\begin{figure}[h!]
\includegraphics[height=\mh]{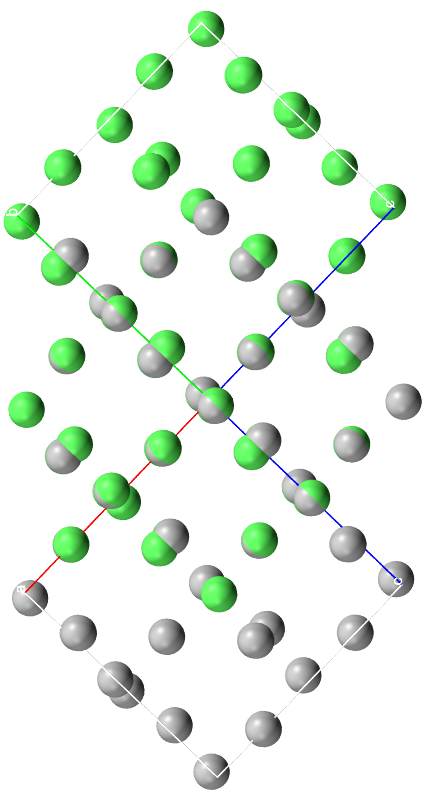}
\hspace*{0mm}
\includegraphics[height=\mh]{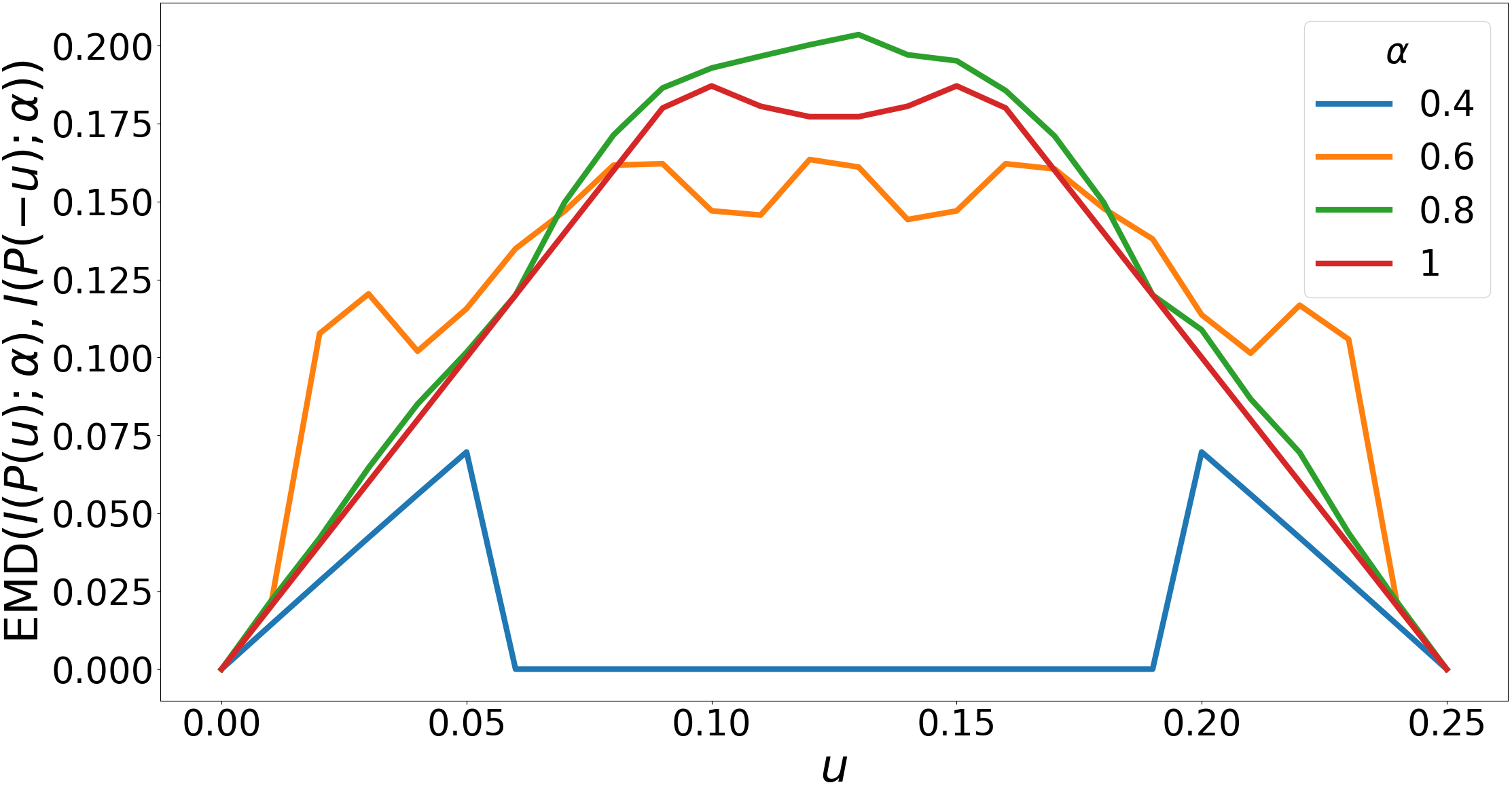}
\caption{\textbf{Left}: a comparison of Pauling's homometric crystals $P(\pm u)$ for $u=0.03$ \cite{pauling1930crystal}, by COMPACK \cite{chisholm2005compack} 
aligning subsets of 48 atoms and outputs RMSD, which fails the triangle inequality. 
The atoms from different $P(\pm 0.03)$ are shown in green and grey.
\textbf{Right}: the pairs of $P(\pm u)$ have $\EMD'>0$ for all $u\in(0,0.25)$ and $\al>0.4$ (running time 50 ms for $u=0.03$ and $\al=0.5$). }
\label{fig:Pauling} 
\end{figure}

The strongest past invariant PDD is based on distances and cannot distinguish mirror images. 
In the CSD, we found four pairs that have identical PDDs but are mirror images shown in Fig.~{\ref{fig:Pauling}}~(right), distinguished by isosets with $\al\geq 1.5\angstrom$ in Fig.~{\ref{fig:mirror_images}}~(left).
For WODLOS vs XAWGAE and $\al=2$, the total time including isosets and EMD is about 4.3 seconds.
All experiments were run on CPU AMD Ryzen 5 5600X, 32GB RAM. 
\myskip

\begin{figure}[h!]
\includegraphics[width=\textwidth]{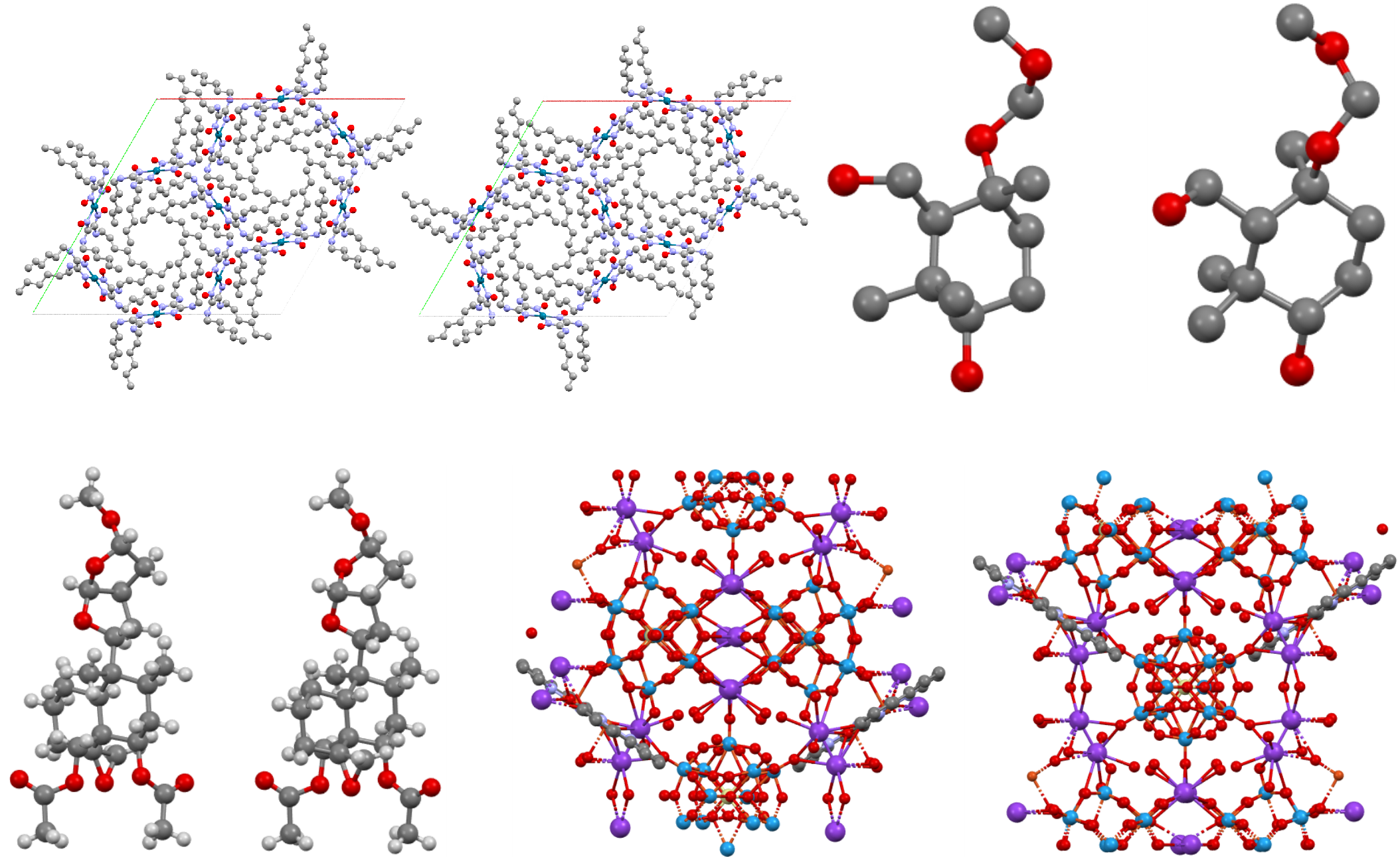}
\caption{Four pairs of mirror images in the CSD are indistinguishable by all past invariants but have approximate $\EMD'>0$ for all radii $\al>1.5\angstrom$ in Fig.~\ref{fig:mirror_images}~(left). 
}
\label{fig:CSD_mirror_images} 
\end{figure}

The limitations of the $\EMD$ metric on isosets in Definition~\ref{dfn:EMD_isosets} are a slower running time than for $\AMD,\PDD$ and the approximate (not exact) algorithm in Corollary~\ref{cor:approximate_EMD}, which are outweighed by the following crucial advantages.
\medskip

\begin{figure}[h!]
\includegraphics[width=\textwidth]{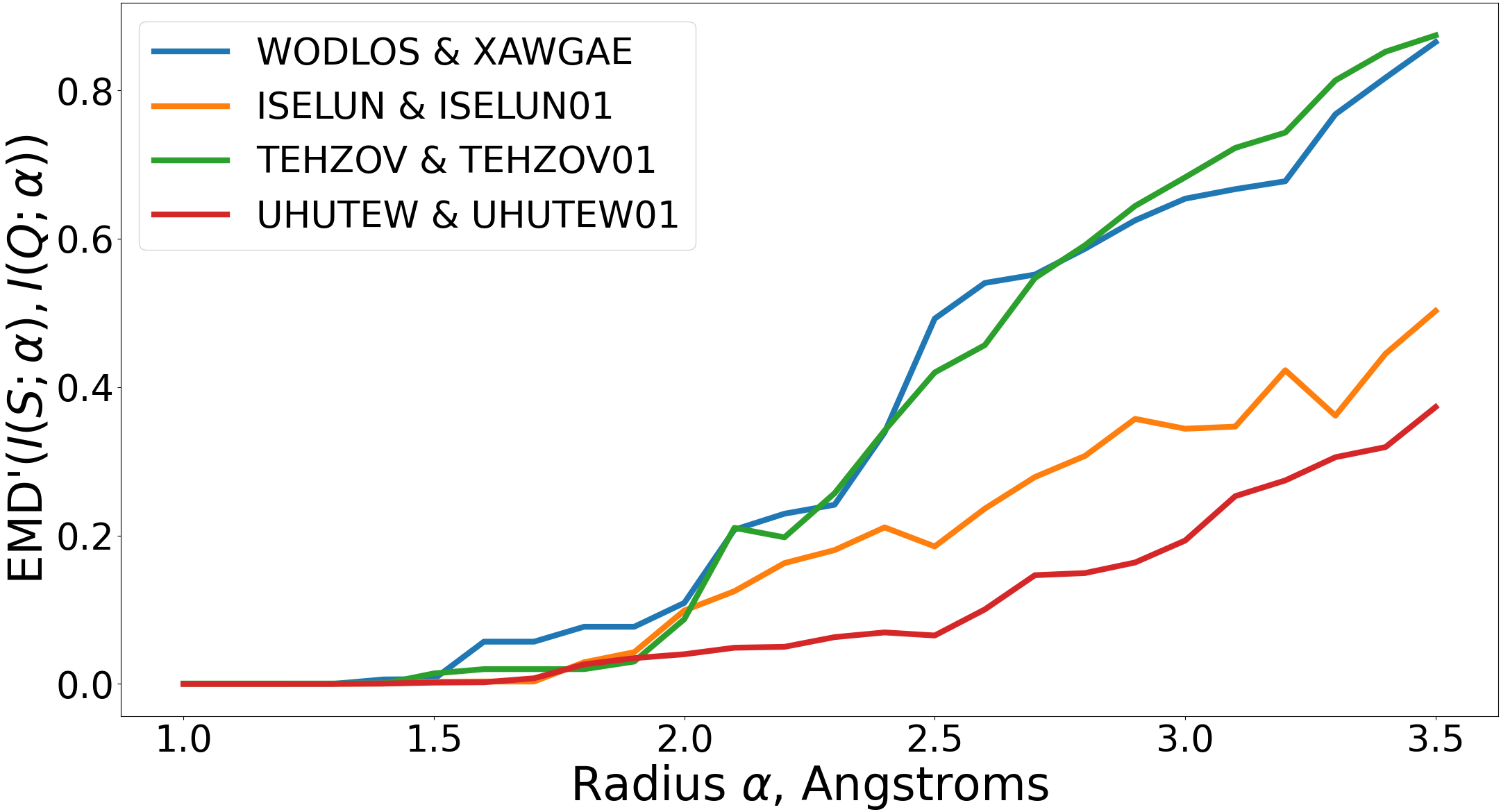}
\myskip

\includegraphics[width=\textwidth]{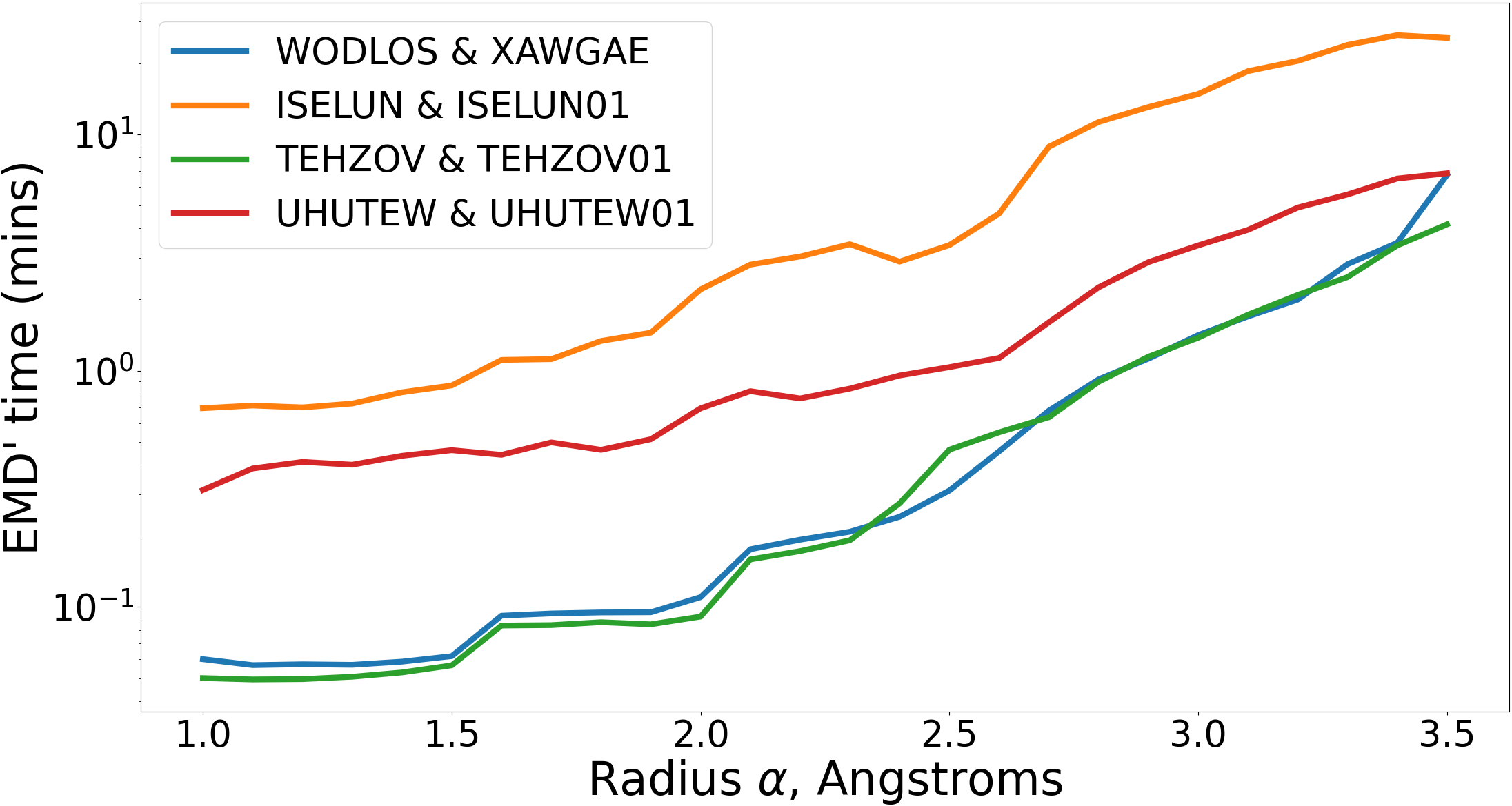}
\caption{
The isosets distinguish all four pairs of mirror images given by their codes in the CSD.
\textbf{Top}: approximate $\EMD'$ for different radii $\al$.
\textbf{Bottom}: running times on a modst desktop.
}
\label{fig:mirror_images} 
\end{figure}

First, all past invariants could not distinguish infinitely many periodic sets (including all mirror images) under rigid motion, e.g. the real crystals in Fig.~\ref{fig:Pauling}.
The new continuous $\EMD$ fully solved  Problem~\ref{pro:periodic_rigid}, which remained open since 1965 \cite{lawton1965reduced}.  
\myskip

Second, because the proved error factor in the practical dimension $n=3$ is close to $4$, any near-duplicate crystals that differ by atomic deviations of up to $\ep$ 
have an exact distance $\EMD\leq 2\ep$ by main Theorem~\ref{thm:isoset_continuous} and hence an approximate distance up to about $8\ep$ by Corollary~\ref{cor:approximate_EMD}.
\smallskip

Any crystals that can be matched under rigid motion are recognisable since our approximation of $\EMD=0$ is also $0$.
Any approximate value $\de$ of $\EMD$ for real crystals $S,Q$ implies that all atoms of $S$ should be perturbed by at least $\de/8$ on average for a complete match with $Q$.
\smallskip

Future work can use the $\EMD$ to continuously quantify changes in material properties under perturbations of atoms and extend Problem~\ref{pro:periodic_rigid} to metrics on finite or periodic sets of points under affine and projective transformations.
\smallskip

In conclusion, sections~\ref{sec:isosets} and \ref{sec:isosets_metrics} prepared the complexity results in section~\ref{sec:isosets_algorithms}: algorithms for computing and comparing isosets (Theorems~\ref{thm:compute_isoset}, Corollary~\ref{cor:compare_isosets}),  and approximating the new boundary tolerant metric $\BT$ (Theorem~\ref{thm:approximate_BT}), and EMD on isosets (Corollary~\ref{cor:approximate_EMD}). 
The proofs expressed polynomial bounds in terms of the motif \emph{size} $m=|S|$ of a periodic set $S$ because the input size of a Crystallographic Information File is linear in $m$, e.g. any lattice has $m=1$. 
\smallskip

The factors depending on the dimension and geometric complexity $\GC(S)$ are inevitable due to the curse of dimensionality and the infinite nature of crystals.
In practice, crystal symmetries reduce a motif to a smaller asymmetric part, which usually has fewer than 20 atoms, even for large molecules in the CSD.
The lower bound via faster PDD invariants in Theorem~\ref{thm:isosets_lower_bound} justifies applying the algorithm of Corollary~\ref{cor:approximate_EMD} only for a final confirmation of near-duplicates. 
So the isosets finalised the hierarchy of the faster but incomplete invariants.
\smallskip

The main novelty is the boundary-tolerant metric in Definition~\ref{dfn:tolerant_metric} that makes the complete invariant isoset Lipschitz continuous (Theorem~\ref{thm:isoset_continuous}) without extra parameters that are needed to smooth past descriptors, such as powder diffraction patterns and atomic environments with fixed cut-off radii.
Since the isoset is the only Lipschitz continuous invariant whose completeness under isometry was proved for all periodic point sets in $\R^n$, the isoset was used to confirm near-duplicates in the CSD (Table~\ref{tab:CSD}) and GNoME (Table~\ref{tab:GNoME}). 

\begin{table}[H]
\label{tab:CSD}
\caption{The first pair consists of rigidly different mirror images from Fig.~{\ref{fig:mirror_images}}~(right).
All others are geometric near-duplicates from (surprisingly) different families in the CSD, confirmed by tiny values of the EMD metric on isosets.
The distance units are in \emph{attometers}: 1 am $=10^{-8}\angstrom=10^{-18}$ meter.
The run times in milliseconds (ms) depend on the cluster size (maximum number of atoms in $\al$-clusters) according to Theorem~\ref{thm:compute_isoset} and Corollary~\ref{cor:approximate_EMD}.}

\begin{tabular}{l|l|r|r|r|r}
CSD id1  & CSD id2  & EMD, am & isosets time, ms & EMD time, ms & cluster size \\
WODLOS   & XAWGAE   & 85856.22         & 129.619         & 1204.58               & 7            \\
TAFQIA   & VAVQIS   & 952.96           & 1690.949        & 321603.86             & 20           \\
FIJKIU   & IPEQUR   & 728.43           & 407.579         & 77455.47              & 16           \\
JIZMIR01 & JIZNAK   & 496.08           & 40.454          & 634.73                & 5            \\
HIYVUG01 & MASPIF   & 334.62           & 35.518          & 543.45                & 7            \\
KIVXEW10 & KIWCEC   & 125.03           & 22.456          & 32.32                 & 5            \\
XAYZOP   & ZEMDAZ   & 89.47            & 301.217         & 1697.26               & 4            \\
KIVXEW10 & KIWCEC28 & 83.07            & 21.701          & 32.49                 & 5            \\
AFIBOH   & NENCUF   & 31.67            & 126.582         & 1160.58               & 5            \\
KIVXEW07 & KIWCEC09 & 31.11            & 22.287          & 36.95                 & 5            \\
KIVXEW07 & KIWCEC11 & 31.11            & 22.434          & 36.75                 & 5            \\
KIVXEW11 & KIWCEC26 & 26.11            & 21.646          & 32.33                 & 5            \\
SERKIL   & SERKOR   & 23.78            & 2444.885        & 18485.57              & 6            \\
ADESAG   & REWPOB   & 5.81             & 54.675          & 5689.27               & 15           \\
GEQRAX   & IFOQOL   & 0.05             & 265.4           & 2090.42               & 6            \\
BUKYEN   & UYOCES   & 0.03             & 398.129         & 15739.11              & 11           \\
GOHYOT   & VIHCEY   & 0.01             & 100.031         & 940.16                & 5            \\
JUMCUP   & QAHBOT   & 0.01             & 179.367         & 4234.6                & 5            \\
CALMOV   & CALNAI   & 0.01             & 128.437         & 3913.32               & 4            \\
NABKOT   & ZIVSEF   & 0.01             & 75.401          & 796.14                & 5            \\
LIBGAE   & VESJUY   & 0.01             & 41.535          & 403.87                & 3            \\
AMEVEV   & OLERON   & 0                & 70.172          & 558.78                & 4            \\
SIHFIZ   & TEZBUV   & 0                & 207.984         & 1761.39               & 5            \\
XATCAA   & ZAQMEN   & 0                & 60.254          & 394.74                & 4            \\
PIDREA   & XIZNOL   & 0                & 94.135          & 243.46                & 5              
\end{tabular}
\end{table}

\begin{table}[H]
\label{tab:GNoME}
\caption{{After excluding 3248 exact numerical duplicates from} \cite[Table~1]{anosova2024importance}, {the next 25 pairs of closest near-duplicates in the GNoME database are confirmed by tiny values of EMD on isosets, see the first pair in Fig.~\ref{fig:exact_duplicates}.
The distance units are \emph{attometers}: 1 am $=10^{-8}\angstrom=10^{-18}$ meter.
The run times are in milliseconds (ms).
The cluster size is the maximum number of atoms in $\al$-clusters.}
}
\begin{tabular}{l|l|l|r|r|r}
GNoME id1  & GNoME id2  & EMD, am & isosets time, ms & EMD time, ms & cluster size \\
\hline
1547d30046 & ddc216e80c & 1                & 1.659            & 434.362               & 14           \\
b4065a4798 & e78d3559e6 & 1.7              & 3.034            & 13.271                & 6            \\
98ab164895 & df1252bc44 & 2                & 1.002            & 419.142               & 14           \\
0de9d25713 & b1733941a7 & 2.7              & 1.971            & 49.816                & 6            \\
0e79f7c053 & 6cf951ac6f & 3                & 1.035            & 429.487               & 14           \\
07ece241f0 & 45cacc8d45 & 3.2              & 0.618            & 14.374                & 6            \\
a58dc74a92 & c16bf63220 & 4.1              & 2.532            & 641.086               & 14           \\
5023e3a4b8 & 8f7ffb4d4a & 4.6              & 2.776            & 10.02                 & 6            \\
3198d1a3ea & 35f67abe6d & 5                & 1.031            & 403.398               & 14           \\
6826b81efb & 76ee112799 & 5                & 0.985            & 407.618               & 14           \\
6826b81efb & e9be17f0ee & 5                & 1.008            & 404.306               & 14           \\
2cff5f2fa0 & f470a5f6fa & 5.3              & 0.635            & 169.911               & 17           \\
2ce912f039 & 9de239ee0c & 5.5              & 0.632            & 3.456                 & 2            \\
c9f5a7a51b & fd9f40e0e1 & 6                & 1.14             & 195.261               & 10           \\
18078e002b & aca2a892a5 & 6                & 1.028            & 421.009               & 14           \\
18078e002b & b9722429b1 & 6                & 1.18             & 445.453               & 14           \\
18078e002b & b702e73db3 & 6                & 1.035            & 414.325               & 14           \\
34b4204eee & adee17535b & 6                & 1.017            & 396.855               & 14           \\
506b8b5646 & 60d266db80 & 6                & 1.174            & 413.254               & 14           \\
506b8b5646 & ec7b789cb3 & 6                & 1.174            & 403.014               & 14           \\
780741962f & a19688f106 & 6.5              & 1.804            & 870.629               & 15           \\
780741962f & c6af1fc763 & 6.5              & 2.731            & 921.496               & 15           \\
780741962f & c64c3e245c & 6.5              & 1.792            & 829.979               & 15           \\
b06353561c & b6d2341d32 & 6.6              & 2.387            & 259.359               & 12           \\
ebb33e044c & ebc9a4db61 & 6.8              & 1.232            & 450.351               & 14           \end{tabular}
\end{table}

\begin{figure}[h!]
\includegraphics[height=100mm]{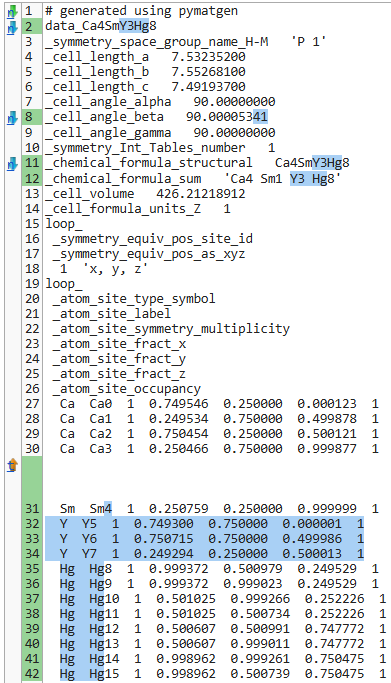}
\includegraphics[height=100mm]{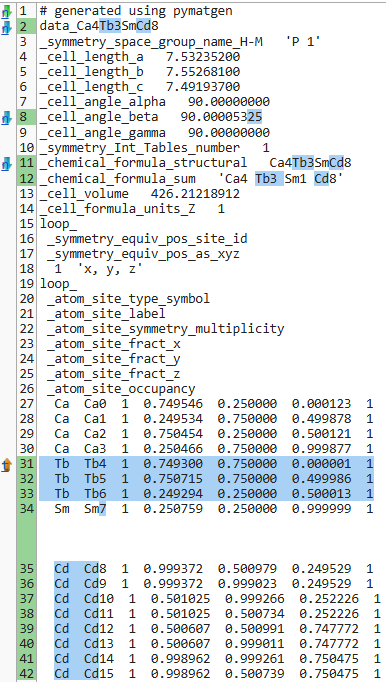}
\caption{The GNoME crystals 1547d30046 and ddc216e80c in the first row of Table~A.2 are compared as texts by https://text-compare.com.
All differences are highlighted in blue.
}
\label{fig:exact_duplicates}
\end{figure}

Fig.~\ref{fig:repeated_atoms} shows the most striking pair of exact duplicates in the GNoME is cdc06a1a2a and 0e2d8f26d6, whose CIFs are identical symbol by symbol in addition to two pairs of atoms at the same positions (Na1=Na2 and Na3=Na4).

\begin{figure}[h!]
\includegraphics[height=30mm]{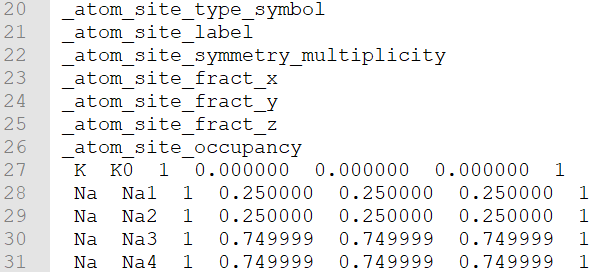}
\includegraphics[height=30mm]{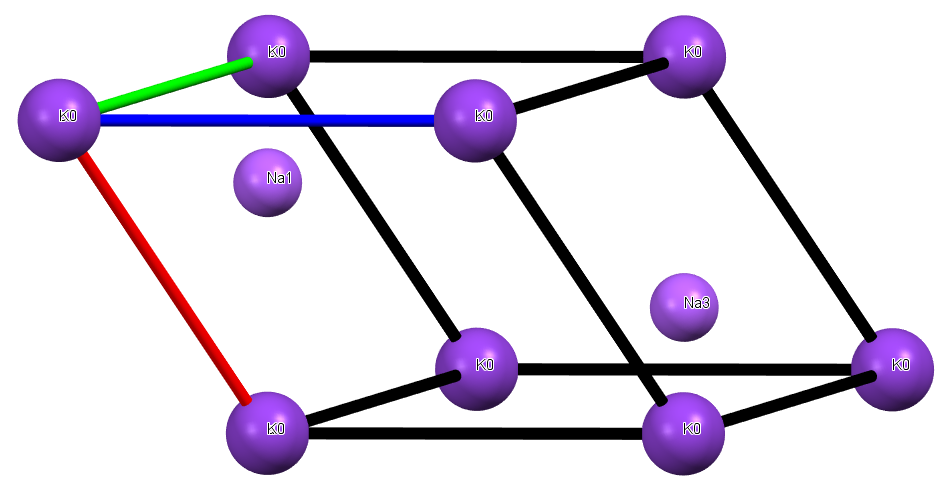}
\caption{Different entries cdc06a1a2a and 0e2d8f26d6 in the GNoME database are not only identical symbol by symbol but also contain two pairs of atoms (Na1=Na2 and Na3=Na4) at the same positions.
\textbf{Left}: a screenshot from the CIF.
\textbf{Right}: Mercury visualisation can show only one atom in each pair of coinciding atoms, e.g. only Na1 and not Na2 from the CIF.}
\label{fig:repeated_atoms}
\end{figure}

\bibliographystyle{plain}
\bibliography{Geometric-Data-Science-book}

%
%
%

\chapter{The most significant results, new concepts, and open problems}
\label{chap:conclusions} 

\abstract{
This chapter concludes with the most important theorems from each of the previous chapters.
The resulting hierarchies of Lipschitz continuous invariants from the ultra-fast to complete ones (under rigid motion in Euclidean space $\R^n$) allowed us to distinguish all non-duplicate objects in major databases of molecules and crystals.
Experimental validation justified new concepts of geometric structures, such as a crystal structure defined as an equivalence class of periodic sets of only atomic centres without chemical elements, under rigid motion in $\R^3$.
The resulting Crystal Isometry Principle uniquely identifies any real periodic material in the continuous moduli space of all periodic point sets.
The book finishes by highlighting open problems and future work.
}

\section{The most important results of Geometric Data Science}
\label{sec:results}

This section briefly summarises the most significant results from Chapters 1-11.
\myskip

\nt
\textbf{Chapter~\ref{chap:intro}} motivated and then formalised the three practical questions 
\sskip

\nt
\emph{(1) Same or different? 
(2) If different, by how much? 
(3) Where do all real objects live?} 
\sskip

\nt
in Geo-Mapping Problem~\ref{pro:geocodes} for all types of geometric data, which can be studied under practically important equivalence relations, such as rigid motion in Euclidean space $\R^n$.
If all conditions of Problem~\ref{pro:geocodes} are satisfied, the resulting moduli spaces of data objects can be explored by \emph{geocodes}, similar to geographic coordinates on Earth.
\myskip

\nt
\textbf{Chapter~\ref{chap:proteins}} described the Backbone Rigid Invariant ($\bri$) in Definition~\ref{dfn:bri}.
Theorems~\ref{thm:BRI_complete}, \ref{thm:BRI_continuous}, \ref{thm:BRI_inverse} proved that the $\bri$ is a geocode (a geographic-style invariant) on the moduli space of non-degenerate protein backbones 
in $\R^3$.
Within a few hours on a modest desktop computer, this invariant detected thousands of exact duplicate backbones and many more near-duplicate chains in the Protein Data Bank (PDB). 
\myskip

\nt
\textbf{Chapter~\ref{chap:directional}} presented in Theorem~\ref{thm:WMI_completeness} complete polynomial-time invariants for finite clouds of unordered points under rigid motion in $\R^n$.
These invariants were outperformed by faster invariants in later chapters, but influenced the latest developments towards a full solution of Geo-Mapping Problem~\ref{pro:geocodes} for unordered points in $\R^n$.
\myskip

\nt
\textbf{Chapter~\ref{chap:PDD-finite}} extended earlier invariants of finite clouds of unordered points (incomplete under isometry in $\R^n$) to the stronger and continuous Pointwise Distance Distribution ($\PDD$) in Definition~\ref{dfn:PDD_finite}.
In addition to being generically complete under isometry in any $\R^n$, the $\PDD$ was proved to be fully complete in Theorem~\ref{thm:PDD_complete_m<5} for any 4 unordered points under isometry in $\R^n$.
It seems the first faster-than-brute-force extension of the side-side-side theorem, known in Euclidean geometry for 2000+ years.
\myskip

\nt
\textbf{Chapter~\ref{chap:SDD}} generalised the $\PDD$ to a stronger Simplexwise Distance Distribution ($\SDD$) in Definition~\ref{dfn:SDD} for an arbitrary metric space.
Theorems~\ref{thm:SDD_metrics} and \ref{thm:SDD_continuity} proved the Lipschitz continuity and polynomial-time computability of the $\SDD$, which is simple enough
to distinguish all (infinitely many) known pairs of non-isometric clouds in $\R^3$ with the same $\PDD$ through manual computations in Examples~\ref{exa:5-point_sets}, \ref{exa:7-point_sets}, and~\ref{exa:6-point_sets}. 
\myskip

\nt
\textbf{Chapter~\ref{chap:SCD}} improved the $\SDD$ to the faster Simplexwide Centred Distribution ($\SCD$) in Definition~\ref{dfn:SCD} for 
clouds of unordered points in $\R^n$.
Theorems~\ref{thm:SCD_complete}, ~\ref{thm:SCD_metrics}, and \ref{thm:OSD+SCD_continuous} for guaranteed completeness, Lipschitz continuity, and polynomial-time computability of the $\SCD$ for any $n$-dimensional cloud of unordered points under rigid motion in $\R^n$.
The key ingredient was the concept of the strength of a simplex in Definition~\ref{dfn:strength_simplex}.
\myskip

\nt
\textbf{Chapter~\ref{chap:1-periodic}}
initiated a continuous approach to point sets that are periodic in one direction in a high-dimensional space $\R\times\R^{n-1}$.
Theorem~\ref{thm:Euclidean_ordered} improved the distance matrix 
to a Lipschitz continuous invariant that is complete under rigid motion, distinguishing all mirror images in $\R^n$.
Theorem~\ref{thm:1-periodic} developed a complete invariant with Lipschitz continuous and polynomial-time metrics for 1-periodic point sets in $\R\times\R^{n-1}$. 
\myskip

\nt
\textbf{Chapter~\ref{chap:lattices2D}} 
fully solved Geo-Mapping Problem~\ref{pro:geocodes} (restated as Problem~\ref{pro:lattices2D}) for all 2-dimensional lattices under four equivalences in $\R^2$ in Theorem~\ref{thm:lattices2D_classification} and Corollary~\ref{cor:lattices2D_classification}.
The case of rigid motion remained discontinuous since the time of Lagrange \cite{lagrange1773recherches}, who classified 2D lattices in terms of quadratic forms, not distinguishing mirror images.
\myskip

\nt
\textbf{Chapter~\ref{chap:densities}} 
studied density functions, which are Lipschitz continuous isometry invariants of arbitrary periodic point sets in $\R^n$.
Theorem~\ref{dfn:densities} proved the generic completeness of density functions under isometry in $\R^3$.
Theorems~\ref{thm:0-th_density}, \ref{thm:1st_density}, and \ref{thm:k-th_density} analytically described density functions for periodic sequences of intervals in $\R$.
\myskip

\nt
\textbf{Chapter~\ref{chap:PDD-periodic}}
extended the Pointwise Distance Distribution ($\PDD$) 
to infinite point sets that are periodic in $l$ directions.
Theorems~\ref{thm:PDD_continuous} and  \ref{thm:PDD_periodic_gen_complete} proved Lipschitz continuity and generic completeness of the $\PDD$ for periodic point sets under isometry in any $\R^n$.
\myskip

\nt
\textbf{Chapter~\ref{chap:isosets}} developed the \emph{isoset} invariant for any periodic point sets in $\R^n$.
Theorem~\ref{thm:isoset_complete} finalised the full completeness of the isoset under isometry and rigid motion in $\R^n$. 
Definition~\ref{dfn:EMD_isosets} introduced the Earth Mover's Distance ($\EMD$) on isosets.
Theorem~\ref{thm:isoset_continuous} and Corollary~\ref{cor:approximate_EMD} proved Lipschitz continuity and approximate polynomial-time algorithms for the $\EMD$ on isosets in $\R^n$, for a fixed dimension $n$.
\myskip

The next section highlights the most important conclusions from applications.

\section{New definitions of geometric structures and verified principles}
\label{sec:principles}

Mathematics allows other sciences to progress by developing new concepts that formalise practical questions.
Hence, it was crucial to state the challenges of ambiguity and discontinuity in terms of complete invariants with Lipschitz continuous metrics.
\myskip

After continuity was stated for metrics on lattices \cite{mosca2020voronoi} in January 2020, the subsequent work \cite{widdowson2022resolving,widdowson2023recognizing,kurlin2024mathematics}
 gradually added more conditions to Geo-Mapping Problem~\ref{pro:geocodes}. 
\myskip

The resulting hierarchies of continuous invariants (from the simplest and ultra-fast to slower but complete) allowed us to open the `black boxes' of major databases, which keep their data in ambiguous photograph-style forms.
For example, crystals are usually stored  as Crystallographic Information Files (CIFs), while molecular geometries are represented by xyz files listing atomic positions in an arbitrary coordinate system.
\myskip

The continuous invariant-based approach revealed thousands of exact geometric duplicates and many more near-duplicates in the PDB, CSD, ICSD, and other databases.
Some exact duplicates could have been found by dataset creators through direct comparisons of digital representations.
For instance, \cite[Table 1]{anosova2024importance} counts thousands of symbol-by-symbol duplicate CIFs in Google's GNoME \cite{google2023}.
These experimental validations led us to new concepts of geometric structures and principles for crystals and molecules.
\myskip

Our \cite[Definition~6]{anosova2024importance} introduced the \emph{crystal structure} as a class of all periodic sets of atoms (as in a CIF) that can be \emph{exactly} matched with each other under rigid motion.
\myskip

This definition emphasises the importance of \emph{exact} matching.
Indeed, ignoring noise up to any tiny threshold $\ep>0$ leads to the \emph{sorites paradox}  \cite{hyde2011sorites} and a trivial classification of all objects within a continuous space due to the transitivity axiom in Definition~\ref{dfn:equivalence}.

\begin{figure}[h!]
\centering
\includegraphics[width=\textwidth]{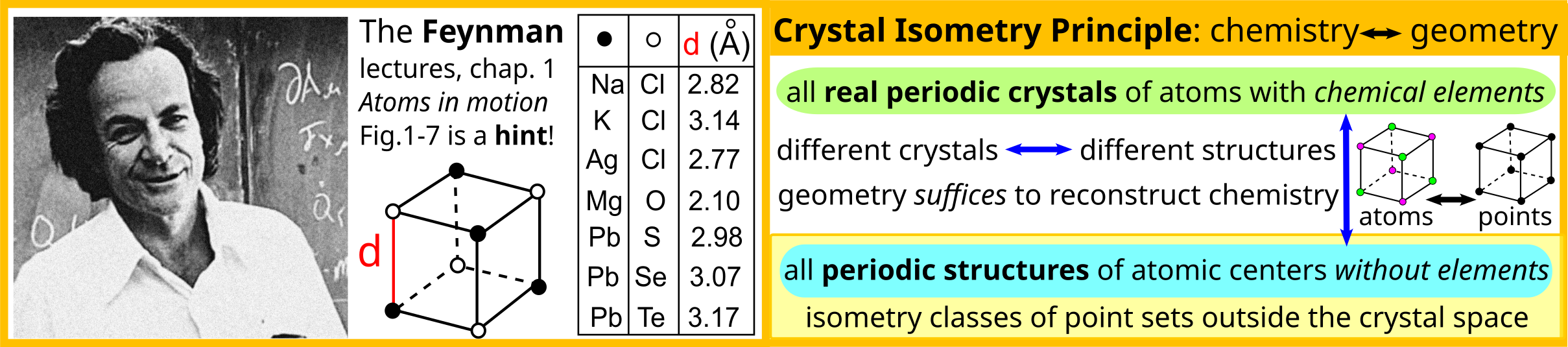}
\caption{
\textbf{Left}: Feynman's table in his first lecture ``Atoms in motion'' \cite{feynman2011lectures}
showed that 7 cubic crystals differ by the smallest interatomic distance $d$, which was the Eureka moment for the second author in May 2021 to realise that all real periodic crystals can be distinguished as periodic sets of atomic centres without chemical elements.
\textbf{Right}: our experiments on the world's largest databases of real materials confirmed the original intuition, now stated as the \emph{Crystal Isometry Principle}: any real periodic material (under fixed ambient conditions) is uniquely identified by a precise enough atomic geometry within a common continuous space of all periodic structures independent of their chemistry and symmetries.
}
\label{fig:Feynman+CRISP}
\end{figure}

Fig.~\ref{fig:Feynman+CRISP}~(left) shows Feynman's table distinguishing 7 cubic crystals by their single geometric invariant $d$ that is the smallest interatomic distance. 
Our much stronger invariants $\AMD,\PDD$, and isosets distinguished all non-duplicate structures among in major materials databases.
Some of the found duplicates had different chemical compositions.
In the first striking example, the pair of CSD entries HIFCAB and JEPLIA has all numbers in their CIFs identical almost to the last digit, but one atom (Mn) is replaced with a different one (Cd).
Since these elements have very different atomic masses 25 and 48, these coincidences of all coordinates seem physically impossible.
\myskip

Indeed, any atomic replacement should change inter-atomic interactions and hence distances to neighbours, which is immediately detected by the $\PDD$ invariant.
Our colleagues at the Cambridge Crystallographic Data Centre, who curate the CSD, checked that the raw diffraction data (structure factors) were also identical in this case.
Several more pairs of duplicates are under investigation by five journals for data integrity.
\myskip

The much more important consequence of the ability to distinguish all periodic crystals by geometry is the \emph{Crystal Isometry Principle} (CRISP) in Fig.~\ref{fig:Feynman+CRISP}~(right).
\medskip

A mathematical formulation of the CRISP says that all real periodic crystals (independent of their chemistry and symmetries) live in a common moduli space $\CRIS(\R^3)$, which is now called the \emph{Crystal Rigid Space}, see its projections in Fig.~\ref{fig:heatmap_C_allotropes} and \cite{widdowson2024continuous}.
Since the first invariants $\AMD,\PDD$ were developed for isometry, not for rigid motion, the initial name was the \emph{Crystal Isometry Space} \cite{widdowson2022resolving}, now denoted by $\CIMS(\R^3)$.
\myskip

Of course, not any periodic set of points can be realised as a periodic crystals.
For instance, distances between neighbouring atoms are usually in small ranges, especially for fixed elements.
In the geographic analogy, not every location on Earth is suitable for humans to live.
However, the knowledge of a full geographic map certainly helped to find all hospitable places and not to waste time on exploring many hostile regions.  
\medskip 

Chemistry has substantially benefited from the periodic table, though it was initially half-empty, as a map of all known chemical elements.
Indeed, \emph{Mendeleev's geocode} consisting of the period and group number provided a complete invariant and guided the search for new elements.
More than 150 years after Mendeleev's breakthrough in 1869, coordination chemistry \cite{bernhardt2025introduction} can progress from studying isolated shapes of of atomic environments to continuous maps parametrised by complete invariants \cite{widdowson2023recognizing}.
\medskip

\begin{figure}[h!]
\centering
\includegraphics[width=\textwidth]{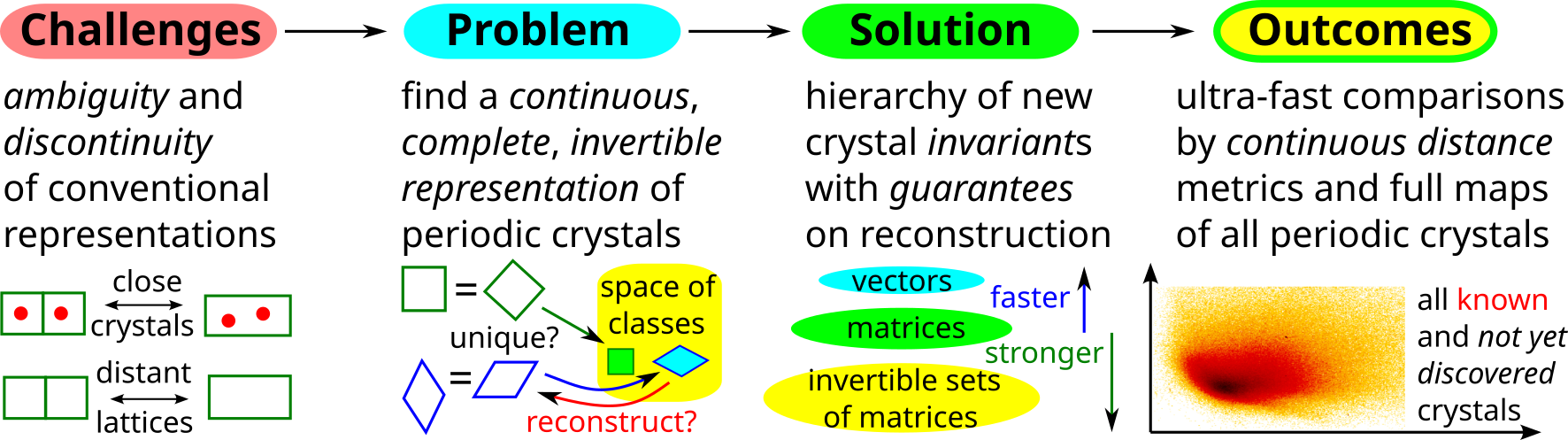}
\caption{
This book exemplified how mathematics can advance physical sciences through rigorous definitions and theorem-based principles verified on all available experimental data, as in Fig.~\ref{fig:Feynman+CRISP}~(right). 
}
\label{fig:challenges-problem-solution-outcomes}
\end{figure}

The complete hierarchy of point cloud invariants \cite{widdowson2023recognizing} helped us to geometrically compare all molecules as clouds of unordered atomic centres, even without covalent bonds and chemical elements.
Large-scale experiments on molecular databases confirmed that a precise enough atomic geometry determines any real \emph{molecular structure}, now defined as an equivalence class of only atomic centres under rigid motion in $\R^3$.
\myskip

The book included several conjectures and also postponed for future work some conditions of Geo-Mapping Problem~\ref{pro:geocodes}, such as inverse continuity and Euclidean embeddabilty  for finite clouds of unordered points in $\R^n$.
We highlight Conjecture~\ref{conj:PDD_complete_n=2} about the completeness of the Pointwise Distance Distribution for any $m\geq 5$ unordered points under isometry in $\R^2$, which can be accessible even to school children.  
\medskip

Fig.~\ref{fig:GDS-theory-applications} illustrates theoretical sources and applications of Geometric Data Science.

\begin{figure}[h!]
\centering
\includegraphics[width=\textwidth]{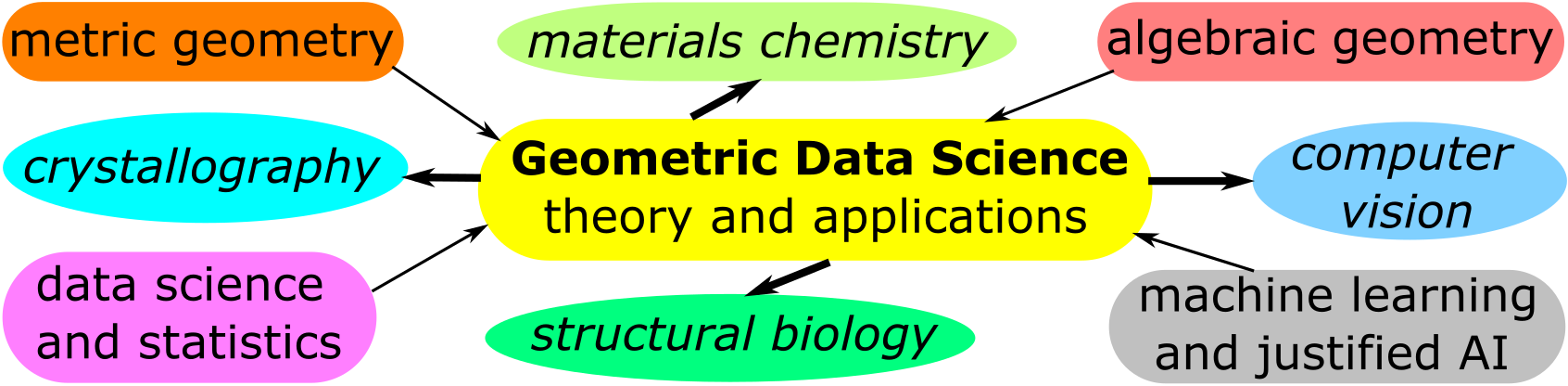}
\caption{Geometric Data Science develops methods of metric and algebraic geometry to enrich statistics and machine learning for applications in computer vision, chemistry, and structural biology.
}
\label{fig:GDS-theory-applications}
\end{figure}

In conclusion, Geometric Data Science `connected the dots' in practical challenges by unifying all important requirements to rigorously answer the basic questions 
(\emph{Same or different?}
\emph{If different, by how much?}\emph{Where do all objects live?}) 
into a list of verifiable conditions in Geo-Mapping Problem~\ref{pro:geocodes} as a guide for future developments.

\bibliographystyle{plain}
\bibliography{Geometric-Data-Science-book}

\begin{thebibliography}{10}

\bibitem{alfakih2018euclidean}
Abdo~Y Alfakih et~al.
\newblock {\em Euclidean distance matrices and their applications in rigidity theory}.
\newblock Springer, 2018.

\bibitem{anosova2025complete}
Olga Anosova et~al.
\newblock A complete and bi-continuous invariant of protein backbones under rigid motion.
\newblock {\em MATCH Comm. Math. Comp. Chemistry}, 94:97--134, 2025.

\bibitem{anosova2025recognition}
Olga Anosova et~al.
\newblock Recognition of near-duplicate periodic patterns by continuous metrics with approximation guarantees.
\newblock {\em Pattern Recognition}, 171:112108, 2025.

\bibitem{arnold2013age}
David Arnold.
\newblock {\em The age of discovery, 1400-1600}.
\newblock Routledge, 2013.

\bibitem{arnold2013real}
V~I Arnold.
\newblock {\em Real algebraic geometry}, volume~66.
\newblock Springer Science \& Business Media, 2013.

\bibitem{berisha2021digital}
Visar Berisha et~al.
\newblock Digital medicine and the curse of dimensionality.
\newblock {\em NPJ digital medicine}, 4(1):153, 2021.

\bibitem{delone1934mathematical}
B.N.Delone, N.~Padurov, and A.~Aleksandrov.
\newblock {\em Mathematical Foundations of Structural Analysis of Crystals}.
\newblock State Technical-Theoretical Press, USSR, 1934.

\bibitem{bochnak2013real}
Jacek Bochnak, Michel Coste, and Marie-Fran{\c{c}}oise Roy.
\newblock {\em Real algebraic geometry}, volume~36.
\newblock Springer Science \& Business Media, 2013.

\bibitem{breiding2024metric}
Paul Breiding, Kathl{\'e}n Kohn, and Bernd Sturmfels.
\newblock {\em Metric algebraic geometry}.
\newblock Springer, 2024.

\bibitem{bright2023geographic}
Matthew Bright et~al.
\newblock Geographic-style maps for 2{D} lattices.
\newblock {\em Acta Cryst A}, 79:1--13, 2023.

\bibitem{brock2021change}
C~Brock.
\newblock Change to the definition of "crystal" in the {I}{U}{C}r {O}nline {D}ictionary of {C}rystallography.
\newblock www.iucr.org/news/newsletter/etc/articles?issue=151351\&result\_138339\_result\_page=17, 2021.

\bibitem{bronstein2017geometric}
Michael Bronstein, Joan Bruna, Yann LeCun, Arthur Szlam, and Pierre Vandergheynst.
\newblock Geometric deep learning: going beyond {E}uclidean data.
\newblock {\em IEEE Signal Processing Magazine}, 34:18--42, 2017.

\bibitem{bronstein2021geometric}
Michael~M Bronstein, Joan Bruna, Taco Cohen, and Petar Veli{\v{c}}kovi{\'c}.
\newblock Geometric deep learning: Grids, groups, graphs, geodesics, and gauges.
\newblock {\em arXiv:2104.13478}, 2021.

\bibitem{carlsson2021topological}
Gunnar Carlsson and Mikael Vejdemo-Johansson.
\newblock {\em Topological data analysis with applications}.
\newblock Cambridge University Press, 2021.

\bibitem{conway1992low}
J~Conway and N~Sloane.
\newblock Low-dimensional lattices. {V}{I}. {V}oronoi reduction of three-dimensional lattices.
\newblock {\em Proceedings Royal Society A}, 436(1896):55--68, 1992.

\bibitem{lagrange1773recherches}
Joseph~Louis De~Lagrange.
\newblock Recherches d’arithm{\'e}tique.
\newblock {\em Nouveaux M{\'e}moires de l’Acad{\'e}mie de Berlin}, 1773.

\bibitem{derksen2001computational}
Harm Derksen and Gregor Kemper.
\newblock {\em Computational invariant theory}.
\newblock Springer, 2001.

\bibitem{dey2022computational}
Tamal~Krishna Dey and Yusu Wang.
\newblock {\em Computational topology for data analysis}.
\newblock Cambridge University Press, 2022.

\bibitem{deza2009encyclopedia}
Michel~Marie Deza and Elena Deza.
\newblock {\em Encyclopedia of distances}.
\newblock Springer, 2009.

\bibitem{dolbilin1998multiregular}
N.~Dolbilin, J~Lagarias, and M.~Senechal.
\newblock Multiregular point systems.
\newblock {\em Discrete \& Computational Geometry}, 20(4):477--498, 1998.

\bibitem{dolgachev2003lectures}
Igor Dolgachev.
\newblock {\em Lectures on invariant theory}.
\newblock Number 296. Cambridge University Press, 2003.

\bibitem{edelsbrunner2010computational}
Herbert Edelsbrunner and John Harer.
\newblock {\em Computational topology: an introduction}.
\newblock AMS, 2010.

\bibitem{feynman2011lectures}
Richard Feynman.
\newblock {\em The Feynman lectures on physics: the new millennium edition}, volume~1.
\newblock 2011.

\bibitem{gordon1992isospectral}
Carolyn Gordon, David Webb, and Scott Wolpert.
\newblock Isospectral plane domains and surfaces via riemannian orbifolds.
\newblock {\em Inventiones mathematicae}, 110(1):1--22, 1992.

\bibitem{gordon1992one}
Carolyn Gordon, David~L Webb, and Scott Wolpert.
\newblock One cannot hear the shape of a drum.
\newblock {\em Bulletin of the American Mathematical Society}, 27(1):134--138, 1992.

\bibitem{hassett2018stable}
Brendan Hassett, Alena Pirutka, and Yuri Tschinkel.
\newblock Stable rationality of quadric surface bundles over surfaces.
\newblock 220(2):341--365, 2018.

\bibitem{hyde2011sorites}
Dominic Hyde.
\newblock The sorites paradox.
\newblock In {\em Vagueness: A guide}, pages 1--17. Springer, 2011.

\bibitem{joharinad2023mathematical}
Parvaneh Joharinad and J{\"u}rgen Jost.
\newblock {\em Mathematical principles of topological and geometric data analysis}.
\newblock Springer, 2023.

\bibitem{kac1966can}
Mark Kac.
\newblock Can one hear the shape of a drum?
\newblock {\em Amer. Math. Monthly}, 73(4P2):1--23, 1966.

\bibitem{kendall2009shape}
David Kendall, Dennis Barden, Thomas Carne, and Huiling Le.
\newblock {\em Shape and shape theory}.
\newblock Wiley \& Sons, 2009.

\bibitem{kraft2000classical}
H~Kraft and C~Procesi.
\newblock {\em Classical Invariant Theory. A Primer}.
\newblock 2000.

\bibitem{kurlin2024mathematics}
Vitaliy Kurlin.
\newblock Mathematics of 2{D} lattices.
\newblock {\em Found. Comp. Mathematics}, 24:805–863, 2024.

\bibitem{landweber2016fiber}
Peter Landweber, Emanuel Lazar, and Neel Patel.
\newblock On fiber diameters of continuous maps.
\newblock {\em The American Mathematical Monthly}, 123(4):392--397, 2016.

\bibitem{liberti2017euclidean}
Leo Liberti and Carlile Lavor.
\newblock {\em Euclidean distance geometry}.
\newblock Springer, 2017.

\bibitem{marron2021object}
James~S Marron and Ian~L Dryden.
\newblock {\em Object oriented data analysis}.
\newblock Chapman and Hall/CRC, 2021.

\bibitem{mumford1994geometric}
David Mumford, John Fogarty, and Frances Kirwan.
\newblock {\em Geometric invariant theory}, volume~34.
\newblock Springer Science \& Business Media, 1994.

\bibitem{olver1999classical}
Peter~J Olver.
\newblock {\em Classical invariant theory}.
\newblock Number~44. Cambridge University Press, 1999.

\bibitem{penner2012decorated}
R~Penner.
\newblock {\em Decorated {T}eichm{\"u}ller theory}, volume~1.
\newblock European Mathematical Society, 2012.

\bibitem{putinar2008emerging}
Mihai Putinar and Seth Sullivant.
\newblock Emerging applications of algebraic geometry.
\newblock 2008.

\bibitem{rass2024metricizing}
Stefan Rass et~al.
\newblock Metricizing the euclidean space towards desired distance relations in point clouds.
\newblock {\em IEEE Trans. Information Forensics and Security}, 19:7304--7319, 2024.

\bibitem{sacchi2020same}
Pietro Sacchi et~al.
\newblock Same or different -- that is the question: identification of crystal forms from crystal structure data.
\newblock {\em Cryst Eng Comm}, 22(43):7170--7185, 2020.

\bibitem{scheiderer2024course}
Claus Scheiderer.
\newblock {\em A course in real algebraic geometry}.
\newblock Springer, 2024.

\bibitem{smith2024generic}
Philip Smith and Vitaliy Kurlin.
\newblock Generic families of finite metric spaces with identical or trivial 1-dimensional persistence.
\newblock {\em Journal of Applied and Computational Topology}, 8:839–855, 2024.

\bibitem{theobald2024real}
Thorsten Theobald.
\newblock {\em Real algebraic geometry and optimization}, volume 241.
\newblock AMS, 2024.

\bibitem{tschinkel2006geometry}
Yuri Tschinkel.
\newblock Geometry over nonclosed fields.
\newblock In {\em International Congress of Mathematicians}, volume~2, pages 637--651, 2006.

\bibitem{widdowson2023recognizing}
Daniel~E Widdowson and Vitaliy~A Kurlin.
\newblock Recognizing rigid patterns of unlabeled point clouds by complete and continuous isometry invariants with no false negatives and no false positives.
\newblock In {\em Computer Vision and Pattern Recognition}, pages 1275--1284, 2023.

\end{thebibliography}


\begin{thebibliography}{10}

\bibitem{anfinsen1973principles}
Christian~B Anfinsen.
\newblock Principles that govern the folding of protein chains.
\newblock {\em Science}, 181(4096):223--230, 1973.

\bibitem{anosova2025complete}
Olga Anosova et~al.
\newblock A complete and bi-continuous invariant of protein backbones under rigid motion.
\newblock {\em MATCH Comm. Math. Comp. Chemistry}, 94:97--134, 2025.

\bibitem{baek2021accurate}
Minkyung Baek et~al.
\newblock Accurate prediction of protein structures and interactions using a three-track neural network.
\newblock {\em Science}, 373(6557):871--876, 2021.

\bibitem{burley2017protein}
Stephen~K Burley, Helen~M Berman, Gerard~J Kleywegt, John~L Markley, Haruki Nakamura, and Sameer Velankar.
\newblock Protein {D}ata {B}ank ({P}{D}{B}): the single global macromolecular structure archive.
\newblock {\em Protein crystallography: methods and protocols}, pages 627--641, 2017.

\bibitem{dekster1987edge}
Boris~V Dekster and John~B Wilker.
\newblock Edge lengths guaranteed to form a simplex.
\newblock {\em Archiv der Mathematik}, 49(4):351--366, 1987.

\bibitem{heifetz2003effect}
A~Heifetz and Miriam Eisenstein.
\newblock Effect of local shape modifications of molecular surfaces on rigid-body protein--protein docking.
\newblock {\em Protein Engineering}, 16(3):179--185, 2003.

\bibitem{hekkelman2024pdbredo}
Maarten Hekkelman, Anastassis Perrakis, and Robbie Joosten.
\newblock P{D}{B}-redo: updated and optimised crystallographic structures.
\newblock \url{http://https://pdb-redo.eu}.

\bibitem{holm2024dali}
Liisa Holm.
\newblock Dali: Protein structure comparison server.
\newblock \url{http://ekhidna2.biocenter.helsinki.fi/dali}.

\bibitem{horn2012matrix}
Roger~A Horn and Charles~R Johnson.
\newblock {\em Matrix analysis}.
\newblock Cambridge University Press, 2012.

\bibitem{jones2022impact}
David~T Jones and Janet~M Thornton.
\newblock The impact of {A}lpha{F}old2 one year on.
\newblock {\em Nature methods}, 19(1):15--20, 2022.

\bibitem{jumper2021highly}
John Jumper et~al.
\newblock Highly accurate protein structure prediction with alphafold.
\newblock {\em Nature}, 596(7873):583--589, 2021.

\bibitem{kabsch1983dictionary}
Wolfgang Kabsch and Christian Sander.
\newblock Dictionary of protein secondary structure: pattern recognition of hydrogen-bonded and geometrical features.
\newblock {\em Biopolymers}, 22(12):2577--2637, 1983.

\bibitem{kendall1977diffusion}
David Kendall.
\newblock The diffusion of shape.
\newblock {\em Advances in Applied Probability}, 9(3):428--430, 1977.

\bibitem{kendall2009shape}
David Kendall, Dennis Barden, Thomas Carne, and Huiling Le.
\newblock {\em Shape and shape theory}.
\newblock Wiley \& Sons, 2009.

\bibitem{kruskal1978multidimensional}
Joseph~B Kruskal and Myron Wish.
\newblock {\em Multidimensional scaling}.
\newblock Number~11. Sage, 1978.

\bibitem{liberti2017euclidean}
Leo Liberti and Carlile Lavor.
\newblock {\em Euclidean distance geometry}.
\newblock Springer, 2017.

\bibitem{linderstrom1952lane}
Kaj~Ulrik Linderstr{\o}m-Lang.
\newblock {\em Lane medical lectures: proteins and enzymes}, volume~6.
\newblock Stanford University Press, 1952.

\bibitem{mariani2013lddt}
Valerio Mariani, Marco Biasini, Alessandro Barbato, and Torsten Schwede.
\newblock lddt: a local superposition-free score for comparing protein structures and models using distance difference tests.
\newblock {\em Bioinformatics}, 29(21):2722--2728, 2013.

\bibitem{murshudov2011refmac5}
Garib Murshudov et~al.
\newblock Refmac5 for the refinement of macromolecular crystal structures.
\newblock {\em Acta Cryst D}, 67(4):355--367, 2011.

\bibitem{press2007numerical}
William~H Press, Saul~A Teukolsky, William~T Vetterling, and Brian~P Flannery.
\newblock {\em Numerical recipes: the art of scientific computing}.
\newblock Cambridge University Press, 2007.

\bibitem{rass2024metricizing}
Stefan Rass et~al.
\newblock Metricizing the euclidean space towards desired distance relations in point clouds.
\newblock {\em IEEE Trans. Information Forensics and Security}, 19:7304--7319, 2024.

\bibitem{rossmann1990molecular}
Michael~G Rossmann.
\newblock The molecular replacement method.
\newblock {\em Acta Cryst A}, 46(2):73--82, 1990.

\bibitem{schoenberg1935remarks}
Isaac Schoenberg.
\newblock Remarks to {M}aurice {F}rechet's article ``{S}ur la definition axiomatique d'une classe d'espace distances vectoriellement applicable sur l'espace de {H}ilbert.
\newblock {\em Annals of Mathematics}, pages 724--732, 1935.

\bibitem{scott2017mathematical}
L~Ridgway Scott and Ariel Fern{\'a}ndez.
\newblock {\em A {M}athematical {A}pproach to {P}rotein {B}iophysics}.
\newblock Springer, 2017.

\bibitem{weyl1946classical}
Hermann Weyl.
\newblock {\em The classical groups: their invariants and representations}.
\newblock Princeton Univ., 1946.

\bibitem{wlodawer2025duplicate}
A~Wlodawer et~al.
\newblock Duplicate entries in the protein data bank: how to detect and handle them.
\newblock {\em Acta Cryst D}, 81:170--180, 2025.

\bibitem{zhang2004scoring}
Yang Zhang and Jeffrey Skolnick.
\newblock Scoring function for automated assessment of protein structure template quality.
\newblock {\em Proteins: Structure, Function, and Bioinform.}, 57(4):702--710, 2004.

\end{thebibliography}


\begin{thebibliography}{10}

\bibitem{abdi2010principal}
Herv{\'e} Abdi and Lynne~J Williams.
\newblock Principal component analysis.
\newblock {\em Wiley interdisciplinary reviews: computational statistics}, 2(4):433--459, 2010.

\bibitem{boutin2004reconstructing}
Mireille Boutin and Gregor Kemper.
\newblock On reconstructing n-point configurations from the distribution of distances or areas.
\newblock {\em Advances in Applied Mathematics}, 32(4):709--735, 2004.

\bibitem{caelli1979generating}
Terry Caelli.
\newblock On generating spatial configurations with identical interpoint distance distributions.
\newblock In {\em Proc. 7th Australian Conference on Combinatorial Mathematics}, pages 69--75, 1979.

\bibitem{chew1997geometric}
Paul Chew, Michael Goodrich, Daniel Huttenlocher, Klara Kedem, Jon Kleinberg, and Dina Kravets.
\newblock Geometric pattern matching under {E}uclidean motion.
\newblock {\em Comp. Geom.}, 7(1-2):113--124, 1997.

\bibitem{fan2018eigenvector}
Jianqing Fan, Weichen Wang, and Yiqiao Zhong.
\newblock An eigenvector perturbation bound and its application to robust covariance estimation.
\newblock {\em Journal of Machine Learning Research}, 18(207):1--42, 2018.

\bibitem{jonker1987shortest}
Roy Jonker and Anton Volgenant.
\newblock A shortest augmenting path algorithm for dense and sparse linear assignment problems.
\newblock {\em Computing}, 38(4):325--340, 1987.

\bibitem{kantorovich1960mathematical}
Leonid Kantorovich.
\newblock Mathematical methods of organizing and planning production.
\newblock {\em Management science}, 6(4):366--422, 1960.

\bibitem{kroto1985c60}
Harold~W Kroto, James~R Heath, Sean~C O’Brien, Robert~F Curl, and Richard~E Smalley.
\newblock C60: Buckminsterfullerene.
\newblock {\em Nature}, 318(6042):162--163, 1985.

\bibitem{kurlin2024polynomial}
Vitaliy Kurlin.
\newblock Polynomial-time algorithms for continuous metrics on atomic clouds of unordered points.
\newblock {\em MATCH Comm. Math. Comp. Chemistry}, 91:79--108, 2024.

\bibitem{rubner2000earth}
Yossi Rubner, Carlo Tomasi, and Leonidas Guibas.
\newblock The {E}arth {M}over's {D}istance as a metric for image retrieval.
\newblock {\em International Journal of Computer Vision}, 40(2):99--121, 2000.

\bibitem{vaserstein1969markov}
Leonid Vaserstein.
\newblock Markov processes over denumerable products of spaces, describing large systems of automata.
\newblock {\em Probl. Per. Inf.}, 5(3):64--72, 1969.

\bibitem{widdowson2022average}
Daniel Widdowson, Marco~M Mosca, Angeles Pulido, Andrew~I Cooper, and Vitaliy Kurlin.
\newblock Average minimum distances of periodic point sets - foundational invariants for mapping all periodic crystals.
\newblock {\em MATCH Commun. Math. Comput. Chem.}, 87:529--559, 2022.

\end{thebibliography}


\begin{thebibliography}{10}

\bibitem{boutin2004reconstructing}
Mireille Boutin and Gregor Kemper.
\newblock On reconstructing n-point configurations from the distribution of distances or areas.
\newblock {\em Advances in Applied Mathematics}, 32(4):709--735, 2004.

\bibitem{caelli1979generating}
Terry Caelli.
\newblock On generating spatial configurations with identical interpoint distance distributions.
\newblock In {\em Proc. 7th Australian Conference on Combinatorial Mathematics}, pages 69--75, 1979.

\bibitem{elkin2022counterexamples}
Yury Elkin and Vitaliy Kurlin.
\newblock Counterexamples expose gaps in the proof of time complexity for cover trees introduced in 2006.
\newblock In {\em Topological Data Analysis and Visualization}, pages 9--17, 2022.

\bibitem{elkin2023new}
Yury Elkin and Vitaliy Kurlin.
\newblock A new near-linear time algorithm for k-nearest neighbor search using a compressed cover tree.
\newblock In {\em Int. Conf. Machine Learning}, pages 9267--9311, 2023.

\bibitem{memoli2011gromov}
Facundo M{\'e}moli.
\newblock Gromov--{W}asserstein distances and the metric approach to object matching.
\newblock {\em Foundations of Computational Mathematics}, 11(4):417--487, 2011.

\bibitem{rubner2000earth}
Yossi Rubner, Carlo Tomasi, and Leonidas Guibas.
\newblock The {E}arth {M}over's {D}istance as a metric for image retrieval.
\newblock {\em International Journal of Computer Vision}, 40(2):99--121, 2000.

\bibitem{schoenberg1935remarks}
Isaac Schoenberg.
\newblock Remarks to {M}aurice {F}rechet's article ``{S}ur la definition axiomatique d'une classe d'espace distances vectoriellement applicable sur l'espace de {H}ilbert.
\newblock {\em Annals of Mathematics}, pages 724--732, 1935.

\bibitem{shirdhonkar2008approximate}
S~Shirdhonkar and D~Jacobs.
\newblock Approximate earth mover’s distance in linear time.
\newblock In {\em Conference on Computer Vision and Pattern Recognition}, pages 1--8, 2008.

\bibitem{vance1982minimum}
Irvin~E Vance.
\newblock Minimum conditions for congruence of quadrilaterals.
\newblock {\em School Science and Mathematics}, 82(5):403--15, 1982.

\bibitem{weyl1946classical}
Hermann Weyl.
\newblock {\em The classical groups: their invariants and representations}.
\newblock Princeton Univ., 1946.

\bibitem{widdowson2022resolving}
Daniel Widdowson and Vitaliy Kurlin.
\newblock Resolving the data ambiguity for periodic crystals.
\newblock {\em Advances in Neural Information Processing Systems}, 35:24625--24638, 2022.

\bibitem{widdowson2025pointwise}
Daniel Widdowson and Vitaliy Kurlin.
\newblock Pointwise distance distributions for detecting near-duplicates in large materials databases.
\newblock {\em SIAM Journal on Applied Mathematics (doi:10.1137/25M1736657, arxiv:2108.04798)}, 2025.

\end{thebibliography}


\begin{thebibliography}{1}

\bibitem{boutin2004reconstructing}
Mireille Boutin and Gregor Kemper.
\newblock On reconstructing n-point configurations from the distribution of distances or areas.
\newblock {\em Advances in Applied Mathematics}, 32(4):709--735, 2004.

\bibitem{efrat2001geometry}
Alon Efrat, Alon Itai, and Matthew~J Katz.
\newblock Geometry helps in bottleneck matching and related problems.
\newblock {\em Algorithmica}, 31(1):1--28, 2001.

\bibitem{pozdnyakov2020incompleteness}
S~Pozdnyakov et~al.
\newblock Incompleteness of atomic structure representations.
\newblock {\em Phys. Rev. Lett.}, 125:166001, 2020.

\bibitem{fredman1987fibonacci}
Michael~L Fredman and Robert~Endre Tarjan.
\newblock Fibonacci heaps and their uses in improved network optimization algorithms.
\newblock {\em Journal of the ACM}, 34(3):596--615, 1987.

\bibitem{keeping1995introduction}
Ernest~Sydney Keeping.
\newblock {\em Introduction to statistical inference}.
\newblock Courier Corporation, 1995.

\bibitem{kurlin2023simplexwise}
Vitaliy Kurlin.
\newblock Simplexwise {D}istance {D}istributions for finite spaces with metrics and measures.
\newblock {\em arXiv:2303.14161 (latest version at http://kurlin.org/projects/cloud-isometry-spaces/SDD.pdf)}, 2023.

\bibitem{memoli2011gromov}
Facundo M{\'e}moli.
\newblock Gromov--{W}asserstein distances and the metric approach to object matching.
\newblock {\em Foundations of Computational Mathematics}, 11(4):417--487, 2011.

\bibitem{widdowson2023recognizing}
Daniel~E Widdowson and Vitaliy~A Kurlin.
\newblock Recognizing rigid patterns of unlabeled point clouds by complete and continuous isometry invariants with no false negatives and no false positives.
\newblock In {\em Computer Vision and Pattern Recognition}, pages 1275--1284, 2023.

\end{thebibliography}


\begin{thebibliography}{1}

\bibitem{chew1997geometric}
Paul Chew, Michael Goodrich, Daniel Huttenlocher, Klara Kedem, Jon Kleinberg, and Dina Kravets.
\newblock Geometric pattern matching under {E}uclidean motion.
\newblock {\em Comp. Geom.}, 7(1-2):113--124, 1997.

\bibitem{kurlin2023strength}
Vitaliy Kurlin.
\newblock The strength of a simplex is the key to a continuous isometry classification of {E}uclidean clouds of unlabelled points.
\newblock {\em arXiv:2303.13486 (latest version at http://kurlin.org/projects/cloud-isometry-spaces/SCD.pdf)}, 2023.

\bibitem{sippl1986cayley}
Manfred Sippl and Harold Scheraga.
\newblock Cayley-{M}enger coordinates.
\newblock {\em PNAS}, 83:2283--2287, 1986.

\bibitem{widdowson2023recognizing}
Daniel~E Widdowson and Vitaliy~A Kurlin.
\newblock Recognizing rigid patterns of unlabeled point clouds by complete and continuous isometry invariants with no false negatives and no false positives.
\newblock In {\em Computer Vision and Pattern Recognition}, pages 1275--1284, 2023.

\end{thebibliography}


\begin{thebibliography}{1}

\bibitem{deza2009encyclopedia}
Michel~Marie Deza and Elena Deza.
\newblock {\em Encyclopedia of distances}.
\newblock Springer, 2009.

\bibitem{kurlin2025complete}
Vitaliy Kurlin.
\newblock Complete and continuous invariants of 1-periodic sequences in polynomial time.
\newblock {\em SIAM Journal on Mathematics of Data Science}, 7:1643--1663, 2025.

\bibitem{pozdnyakov2022incompleteness}
Sergey~N Pozdnyakov and Michele Ceriotti.
\newblock Incompleteness of graph neural networks for points clouds in three dimensions.
\newblock {\em Machine Learning: Science and Technology}, 3(4):045020, 2022.

\end{thebibliography}


\begin{thebibliography}{10}

\bibitem{norms}
Relations between norms, 2021.

\bibitem{aroyo2011crystallography}
Mois~I Aroyo, JM~Perez-Mato, Danel Orobengoa, EMRE Tasci, Gemma de~la Flor, and Asel Kirov.
\newblock Crystallography online: Bilbao crystallographic server.
\newblock {\em Bulg. Chem. Commun}, 43(2):183--197, 2011.

\bibitem{aroyo2013international}
Mois~Ilia Aroyo and H~Wondratschek.
\newblock {\em International Tables for Crystallography}.
\newblock Wiley Online Library, 2013.

\bibitem{delone1934mathematical}
B.N.Delone, N.~Padurov, and A.~Aleksandrov.
\newblock {\em Mathematical Foundations of Structural Analysis of Crystals}.
\newblock State Technical-Theoretical Press, USSR, 1934.

\bibitem{bright2021complete}
Matthew Bright, Andrew~I Cooper, and Vitaliy Kurlin.
\newblock A complete and continuous map of the lattice isometry space for all 3-dimensional lattices.
\newblock {\em arXiv:2109.11538}, 2021.

\bibitem{bright2023geographic}
Matthew Bright et~al.
\newblock Geographic-style maps for 2{D} lattices.
\newblock {\em Acta Cryst A}, 79:1--13, 2023.

\bibitem{bright2023continuous}
Matthew~J Bright, Andrew~I Cooper, and Vitaliy~A Kurlin.
\newblock Continuous chiral distances for 2-dimensional lattices.
\newblock {\em Chirality}, 35:920--936, 2023.

\bibitem{conway1992low}
J~Conway and N~Sloane.
\newblock Low-dimensional lattices. {V}{I}. {V}oronoi reduction of three-dimensional lattices.
\newblock {\em Proceedings Royal Society A}, 436(1896):55--68, 1992.

\bibitem{lagrange1773recherches}
Joseph~Louis De~Lagrange.
\newblock Recherches d’arithm{\'e}tique.
\newblock {\em Nouveaux M{\'e}moires de l’Acad{\'e}mie de Berlin}, 1773.

\bibitem{delone1938geometry}
B.~N. Delone.
\newblock Geometry of positive quadratic forms. part ii (in russian).
\newblock {\em Uspekhi Matematicheskikh Nauk}, (4):102--164, 1938.

\bibitem{delone1975bravais}
BN~Delone, RV~Galiulin, and MI~Shtogrin.
\newblock On the {B}ravais types of lattices.
\newblock {\em Journal of Soviet Mathematics}, 4(1):79--156, 1975.

\bibitem{engel2004lattice}
Peter Engel, Louis Michel, and Marjorie S{\'e}n{\'e}chal.
\newblock Lattice geometry.
\newblock Technical Report IHES-P-2004-45, 2004.

\bibitem{gruber1989reduced}
B~Gruber.
\newblock Reduced cells based on extremal principles.
\newblock {\em Acta Cryst A}, 45(1):123--131, 1989.

\bibitem{jones1987complex}
Gareth~A Jones and David Singerman.
\newblock {\em Complex functions: an algebraic and geometric viewpoint}.
\newblock Cambridge University press, 1987.

\bibitem{jost2013compact}
J{\"u}rgen Jost.
\newblock {\em Compact Riemann surfaces: an introduction to contemporary mathematics}.
\newblock Springer Science \& Business Media, 2013.

\bibitem{kurlin2022complete}
V~Kurlin.
\newblock A complete isometry classification of 3{D} lattices.
\newblock {\em arxiv:2201.10543}, 2022.

\bibitem{kurlin2024mathematics}
Vitaliy Kurlin.
\newblock Mathematics of 2{D} lattices.
\newblock {\em Found. Comp. Mathematics}, 24:805–863, 2024.

\bibitem{lawton1965reduced}
Stephen Lawton and Robert Jacobson.
\newblock The reduced cell and its crystallographic applications.
\newblock Technical report, Ames Lab, Iowa State University, 1965.

\bibitem{minkowski1891ueber}
Hermann Minkowski.
\newblock Ueber die positiven quadratischen formen und {\"u}ber kettenbruch{\"a}hnliche algorithmen.
\newblock {\em Journal f{\"u}r die reine und angewandte Mathematik (Crelles Journal)}, (107):278--297, 1891.

\bibitem{niggli1928krystallographische}
P.~Niggli.
\newblock {\em Krystallographische und strukturtheoretische Grundbegriffe}, volume~1.
\newblock Akademische verlagsgesellschaft mbh, 1928.

\bibitem{selling1874ueber}
Eduard Selling.
\newblock Ueber die bin{\"a}ren und tern{\"a}ren quadratischen formen.
\newblock {\em Journal f{\"u}r die reine und angewandte Mathematik}, 77:143--229, 1874.

\bibitem{voronoi1908nouvelles}
Georges Voronoi.
\newblock Nouvelles applications des param{\`e}tres continus {\`a} la th{\'e}orie des formes quadratiques.
\newblock {\em J. Reine Angew. Math}, (133):97--178, 1908.

\bibitem{widdowson2022average}
Daniel Widdowson, Marco~M Mosca, Angeles Pulido, Andrew~I Cooper, and Vitaliy Kurlin.
\newblock Average minimum distances of periodic point sets - foundational invariants for mapping all periodic crystals.
\newblock {\em MATCH Commun. Math. Comput. Chem.}, 87:529--559, 2022.

\bibitem{zhilinskii2016introduction}
B.~Zhilinskii.
\newblock {\em Introduction to lattice geometry through group action}.
\newblock EDP sciences, 2016.

\end{thebibliography}


\begin{thebibliography}{10}

\bibitem{anosova2022density}
Olga Anosova and Vitaliy Kurlin.
\newblock Density functions of periodic sequences.
\newblock In {\em Lecture Notes in Computer Science (Proceedings of DGMM)}, volume 13493, pages 395--408, 2022.

\bibitem{anosova2023density}
Olga Anosova and Vitaliy Kurlin.
\newblock Density functions of periodic sequences of continuous events.
\newblock {\em Journal of Mathematical Imaging and Vision}, 65:689–701, 2023.

\bibitem{anosova2023R}
Olgay Anosova.
\newblock R code for density functions of periodic sequences, 2023.

\bibitem{edelsbrunner2021density}
Herbert Edelsbrunner, Teresa Heiss, Vitaliy Kurlin, Philip Smith, and Mathijs Wintraecken.
\newblock The density fingerprint of a periodic point set.
\newblock In {\em SoCG}, volume 189, pages 32:1--32:16, 2021.

\bibitem{edelsbrunner2023simple}
Herbert Edelsbrunner and Georg Osang.
\newblock A simple algorithm for higher-order delaunay mosaics and alpha shapes.
\newblock {\em Algorithmica}, 85(1):277--295, 2023.

\bibitem{edelsbrunner1986voronoi}
Herbert Edelsbrunner and Raimund Seidel.
\newblock Voronoi diagrams and arrangements.
\newblock {\em Discrete \& Computational Geometry}, 1(1):25--44, 1986.

\bibitem{grunbaum1995use}
Gr{\"u}nbaum and Moore.
\newblock The use of higher-order invariants in the determination of generalized {P}atterson cyclotomic sets.
\newblock {\em Acta Cryst A}, 51:310--323, 1995.

\bibitem{nguyen2009low}
Phong~Q Nguyen and Damien Stehl{\'e}.
\newblock Low-dimensional lattice basis reduction revisited.
\newblock {\em ACM Transactions on algorithms (TALG)}, 5(4):1--48, 2009.

\bibitem{pulido2017functional}
A~Pulido et~al.
\newblock Functional materials discovery using energy--structure maps.
\newblock {\em Nature}, 543:657--664, 2017.

\bibitem{smith2022practical}
Phil Smith and Vitaliy Kurlin.
\newblock A practical algorithm for degree-k voronoi domains of three-dimensional periodic point sets.
\newblock In {\em Lecture Notes in Computer Science (Proceedings of ISVC)}, volume 13599, pages 377--391, 2022.

\bibitem{widdowson2022resolving}
Daniel Widdowson and Vitaliy Kurlin.
\newblock Resolving the data ambiguity for periodic crystals.
\newblock {\em Advances in Neural Information Processing Systems}, 35:24625--24638, 2022.

\bibitem{widdowson2022average}
Daniel Widdowson, Marco~M Mosca, Angeles Pulido, Andrew~I Cooper, and Vitaliy Kurlin.
\newblock Average minimum distances of periodic point sets - foundational invariants for mapping all periodic crystals.
\newblock {\em MATCH Commun. Math. Comput. Chem.}, 87:529--559, 2022.

\end{thebibliography}


\begin{thebibliography}{10}

\bibitem{anosova2021isometry}
O~Anosova and V~Kurlin.
\newblock An isometry classification of periodic point sets.
\newblock In {\em Proceedings of Discrete Geometry and Mathematical Morphology}, pages 229--241, 2021.

\bibitem{anosova2025recognition}
Olga Anosova et~al.
\newblock Recognition of near-duplicate periodic patterns by continuous metrics with approximation guarantees.
\newblock {\em Pattern Recognition}, 171:112108, 2025.

\bibitem{anosova2024importance}
Olga Anosova, Vitaliy Kurlin, and Marjorie Senechal.
\newblock The importance of definitions in crystallography.
\newblock {\em International Union of Crystallography Journal}, 11:453--463, 2024.

\bibitem{azrour2011rietveld}
M~Azrour et~al.
\newblock Rietveld refinements and vibrational spectroscopic studies of \ce{Na_{1- x}K_x Pb_4 (PO4)_3} lacunar apatites ($0\leq x\leq$ 1).
\newblock {\em Journal of Physics and Chemistry of Solids}, 72(11):1199--1205, 2011.

\bibitem{batatia2022mace}
Ilyes Batatia, David~P Kovacs, Gregor Simm, Christoph Ortner, and G{\'a}bor Cs{\'a}nyi.
\newblock Mace: Higher order equivariant message passing neural networks for fast and accurate force fields.
\newblock {\em Advances in Neural Information Processing Systems}, 35:11423--11436, 2022.

\bibitem{batsanov2001van}
Stepan Batsanov.
\newblock Van der {W}aals radii of elements.
\newblock {\em Inorganic mat.}, 37:871--885, 2001.

\bibitem{carstens1999geometrical}
Hans-Georg Carstens et~al.
\newblock Geometrical bijections in discrete lattices.
\newblock {\em Combinatorics, Probability and Computing}, 8(1-2):109--129, 1999.

\bibitem{cerqueira2015identification}
T~Cerqueira et~al.
\newblock Identification of novel {C}u, {A}g, and {A}u ternary oxides from global structural prediction.
\newblock {\em Chemistry of Materials}, 27(13):4562--4573, 2015.

\bibitem{chawla2023crystallography}
Dalmeet~Singh Chawla.
\newblock Crystallography databases hunt for fraudulent structures.
\newblock {\em ACS Central Science}, 9:1853–1855, 2023.

\bibitem{chisholm2005compack}
J.~Chisholm and S.~Motherwell.
\newblock Compack: a program for identifying crystal structure similarity using distances.
\newblock {\em J. Applied Crystal.}, 38:228--231, 2005.

\bibitem{duneau1991bounded}
Michel Duneau and Christophe Oguey.
\newblock Bounded interpolations between lattices.
\newblock {\em Journal of Physics A: Mathematical and General}, 24(2):461, 1991.

\bibitem{edelsbrunner2021density}
Herbert Edelsbrunner, Teresa Heiss, Vitaliy Kurlin, Philip Smith, and Mathijs Wintraecken.
\newblock The density fingerprint of a periodic point set.
\newblock In {\em SoCG}, volume 189, pages 32:1--32:16, 2021.

\bibitem{elkin2023new}
Yury Elkin and Vitaliy Kurlin.
\newblock A new near-linear time algorithm for k-nearest neighbor search using a compressed cover tree.
\newblock In {\em Int. Conf. Machine Learning}, pages 9267--9311, 2023.

\bibitem{gieseke2014buffer}
Fabian Gieseke, Justin Heinermann, Cosmin Oancea, and Christian Igel.
\newblock Buffer kd trees: processing massive nearest neighbor queries on {G}{P}{U}s.
\newblock In {\em Intern. Conf. Machine Learning}, pages 172--180, 2014.

\bibitem{gravzulis2009crystallography}
Saulius Gra{\v{z}}ulis et~al.
\newblock Crystallography open database--an open-access collection of crystal structures.
\newblock {\em J Appl. Crystallography}, 42(4):726--729, 2009.

\bibitem{griesemer2021high}
Sean~D Griesemer, Logan Ward, and Chris Wolverton.
\newblock High-throughput crystal structure solution using prototypes.
\newblock {\em Physical Review Materials}, 5(10):105003, 2021.

\bibitem{ivic2004lattice}
A~Ivic, E~Kr{\"a}tzel, M~K{\"u}hleitner, and WG~Nowak.
\newblock Lattice points in large regions and related arithmetic functions: recent developments in a very classic topic.
\newblock {\em Publications of the Scientific Society at the Johann Wolfgang Goethe University}, pages 89--128, 2006.

\bibitem{jain2013commentary}
Anubhav Jain et~al.
\newblock Commentary: The materials project: A materials genome approach to accelerating materials innovation.
\newblock {\em APL materials}, 1(1), 2013.

\bibitem{kurlin2022complete}
V~Kurlin.
\newblock A complete isometry classification of 3{D} lattices.
\newblock {\em arxiv:2201.10543}, 2022.

\bibitem{kurlin2024mathematics}
Vitaliy Kurlin.
\newblock Mathematics of 2{D} lattices.
\newblock {\em Found. Comp. Mathematics}, 24:805–863, 2024.

\bibitem{leeman2024challenges}
Josh Leeman, Yuhan Liu, Joseph Stiles, Scott~B Lee, Prajna Bhatt, Leslie~M Schoop, and Robert~G Palgrave.
\newblock Challenges in high-throughput inorganic materials prediction and autonomous synthesis.
\newblock {\em PRX Energy}, 3(1):011002, 2024.

\bibitem{macdonald1998symmetric}
Ian~Grant Macdonald.
\newblock {\em Symmetric functions and Hall polynomials}.
\newblock Oxford University Press, 1998.

\bibitem{merchant2023scaling}
Amil Merchant, Simon Batzner, Samuel~S Schoenholz, Muratahan Aykol, Gowoon Cheon, and Ekin~Dogus Cubuk.
\newblock Scaling deep learning for materials discovery.
\newblock {\em Nature}, pages 80--85, 2023.

\bibitem{merchant2023millions}
Amil Merchant and Ekin~Dogus Cubuk.
\newblock Millions of new materials discovered with deep learning.
\newblock \url{https://deepmind.google/discover/blog/millions-of-new-materials-discovered-with-deep-learning/}, 2023.

\bibitem{peplow2023robot}
Mark Peplow.
\newblock Robot chemist sparks row with claim it created new materials.
\newblock \url{}https://www.nature.com/articles/d41586-023-03956-w, 2023.

\bibitem{taylor2019million}
Robin Taylor and Peter~A Wood.
\newblock A million crystal structures: The whole is greater than the sum of its parts.
\newblock {\em Chemical reviews}, 119(16):9427--9477, 2019.

\bibitem{terban2021structural}
Maxwell~W Terban and Simon~JL Billinge.
\newblock Structural analysis of molecular materials using the pair distribution function.
\newblock {\em Chemical Reviews}, 122(1):1208--1272, 2021.

\bibitem{widdowson2022resolving}
Daniel Widdowson and Vitaliy Kurlin.
\newblock Resolving the data ambiguity for periodic crystals.
\newblock {\em Advances in Neural Information Processing Systems}, 35:24625--24638, 2022.

\bibitem{widdowson2025geographic}
Daniel Widdowson and Vitaliy Kurlin.
\newblock Geographic-style maps with a local novelty distance help navigate in the materials space.
\newblock {\em Scientific Reports}, 15:27588, 2025.

\bibitem{widdowson2025higher}
Daniel Widdowson and Vitaliy Kurlin.
\newblock Higher-order, generically complete, continuous, and polynomial-time isometry invariants of periodic sets.
\newblock {\em arXiv:2509.15088}, 2025.

\bibitem{widdowson2025pointwise}
Daniel Widdowson and Vitaliy Kurlin.
\newblock Pointwise distance distributions for detecting near-duplicates in large materials databases.
\newblock {\em SIAM Journal on Applied Mathematics (doi:10.1137/25M1736657, arxiv:2108.04798)}, 2025.

\bibitem{widdowson2022average}
Daniel Widdowson, Marco~M Mosca, Angeles Pulido, Andrew~I Cooper, and Vitaliy Kurlin.
\newblock Average minimum distances of periodic point sets - foundational invariants for mapping all periodic crystals.
\newblock {\em MATCH Commun. Math. Comput. Chem.}, 87:529--559, 2022.

\bibitem{zagorac2019recent}
Dejan Zagorac, H~M{\"u}ller, S~Ruehl, J~Zagorac, and Silke Rehme.
\newblock Recent developments in the inorganic crystal structure database.
\newblock {\em J Applied Crystallography}, 52(5):918--925, 2019.

\end{thebibliography}


\begin{thebibliography}{10}

\bibitem{google2023}
Google {D}eep{M}ind's {G}{N}o{M}{E} database of 384,398 {C}{I}{F}s.
\newblock \url{https://deepmind.google/discover/blog/millions-of-new-materials-discovered-with-deep-learning}, 2023.

\bibitem{anosova2024importance}
Olga Anosova, Vitaliy Kurlin, and Marjorie Senechal.
\newblock The importance of definitions in crystallography.
\newblock {\em International Union of Crystallography Journal}, 11:453--463, 2024.

\bibitem{bernhardt2025introduction}
Paul~V Bernhardt and Geoffrey~A Lawrance.
\newblock {\em Introduction to coordination chemistry}.
\newblock John Wiley \& Sons, 2025.

\bibitem{lagrange1773recherches}
Joseph~Louis De~Lagrange.
\newblock Recherches d’arithm{\'e}tique.
\newblock {\em Nouveaux M{\'e}moires de l’Acad{\'e}mie de Berlin}, 1773.

\bibitem{feynman2011lectures}
Richard Feynman.
\newblock {\em The Feynman lectures on physics: the new millennium edition}, volume~1.
\newblock 2011.

\bibitem{hyde2011sorites}
Dominic Hyde.
\newblock The sorites paradox.
\newblock In {\em Vagueness: A guide}, pages 1--17. Springer, 2011.

\bibitem{kurlin2024mathematics}
Vitaliy Kurlin.
\newblock Mathematics of 2{D} lattices.
\newblock {\em Found. Comp. Mathematics}, 24:805–863, 2024.

\bibitem{mosca2020voronoi}
M~Mosca and V~Kurlin.
\newblock Voronoi-based similarity distances between arbitrary crystal lattices.
\newblock {\em Crystal Research and Technology}, 55(5):1900197, 2020.

\bibitem{widdowson2022resolving}
Daniel Widdowson and Vitaliy Kurlin.
\newblock Resolving the data ambiguity for periodic crystals.
\newblock {\em Advances in Neural Information Processing Systems}, 35:24625--24638, 2022.

\bibitem{widdowson2024continuous}
Daniel Widdowson and Vitaliy Kurlin.
\newblock Continuous invariant-based maps of the cambridge structural database.
\newblock {\em Crystal Growth and Design}, 24:5627–5636, 2024.

\bibitem{widdowson2023recognizing}
Daniel~E Widdowson and Vitaliy~A Kurlin.
\newblock Recognizing rigid patterns of unlabeled point clouds by complete and continuous isometry invariants with no false negatives and no false positives.
\newblock In {\em Computer Vision and Pattern Recognition}, pages 1275--1284, 2023.

\end{thebibliography}

\backmatter

%
%

\extrachap{Acronyms}



\begin{description}[ABCDE]
\item[ADA]{Average Deviation from Asymptotic}
\item[ADD]{Average Distance Distribution}
\item[AMD]{Average Minimum Distances}
\item[AND]{Average Normalised Distances}
\item[ASD]{Average Simplexwise Distribution}
\item[ASM]{Average Simplexwise Moments}
\item[BD]{Bottleneck Distance}
\item[BRI]{Backbone Rigid Invariant}
\item[BRIS]{Backbone Rigid Invariant Space}
\item[BRAIN]{Backbone Rigid Average Invariant}
\item[BT]{Boundary Tolerant metric}
\item[CDM]{Cyclic Distance Matrix}
\item[CDS]{Cyclic Distances with Signs}
\item[CIF]{Crystallographic Information File}
\item[CIS]{Cloud Isometry Space}
\item[CIM]{Cyclic Isometry Metric}
\item[CIMS]{Crystal Isometry Space}
\item[COD]{Crystallography Open Database}
\item[CR]{Cyclic Rigid invariant}
\item[CRM]{Cyclic Rigid Metric}
\item[CRS]{Cloud Rigid Space}
\item[CRIS]{Crystal Rigid Space}
\item[CRISP]{Crystal Isometry Principle}
\item[CSD]{Cambridge Structural Database}
\item[DIM]{Dihedral Isometry Metric}
\item[DR]{Dihedral Rigid invariant}
\item[DRM]{Dihedral Rigid Metric}
\item[EMD]{Earth Mover's Distance}
\item[GC]{Geometric Complexity}
\item[GDS]{Geometric Data Science}
\item[GNoME]{Graph Network Materials Exploration dataset}
\item[HD]{Hausdorff Distance}
\item[ICSD]{Inorganic Crystal Structural Database}
\item[LAC]{Linear Assignment Cost}
\item[LDS]{Lattice Dilation Space}
\item[LHS]{Lattice Homothety Space}
\item[LIS]{Lattice Isometry Space}
\item[LRS]{Lattice Rigid Space}
\item[MSD]{Measured Simplexwise Distribution}
\item[MP]{Materials Project}
\item[OCD]{Oriented Centred Distribution}
\item[ORD]{Oriented Relative Distribution}
\item[OSD]{Oriented Simplexwise Distribution}
\item[PCA]{Principal Components Analysis}
\item[PCI]{Principal Coordinates Invariant}
\item[PCM]{Principal Coordinates Matrix}
\item[PDA]{Pointwise Deviation from Asymptotic}
\item[PDD]{Pointwise Distance Distribution}
\item[PND]{Pointwise Normalised Distribution}
\item[PPC]{Point Packing Coefficient}
\item[PI]{Projected Invariant}
\item[QT]{Quotient Triangle}
\item[RDD]{Relative Distance Distribution}
\item[RI]{Root Invariant}
\item[SBD]{Superbases under Dilation}
\item[SBH]{Superbases under Homothety}
\item[SBI]{Superbases under Isometry}
\item[SBR]{Superbases under Rigid Motion}
\item[SCD]{Simplexwise Centred Distribution}
\item[SDD]{Simplexwise Distance Distribution}
\item[SDM]{Superbase Dilation Metric}
\item[SHM]{Superbase Homothety Metric}
\item[SIM]{Superbase Isometry Metric}
\item[SLM]{Spherical Lattice Map}
\item[SM]{Symmetrised Metric}
\item[SPD]{Sorted Pairwise Distances}
\item[SRD]{Sorted Radial Distances}
\item[SRM]{Superbase Rigid Metric}
\item[TC]{Triangular Cone}
\item[TDA]{Topological Data Analysis}
\item[TRIN]{Triangular Invariant}
\item[WMI]{Weighted Matrices Invariant}
\item[WSD]{Weighted Simplexwise Distribution}

\end{description}

\printindex


\end{document}